%% file: BZSVpaper.tex
\documentclass[11pt]{amsart} 
\usepackage{float}
\usepackage{amsmath}
  \usepackage{amsfonts}
\usepackage{amsthm}
\usepackage{epsfig}

 \usepackage{amssymb} 
 \usepackage{imakeidx}
 \usepackage{mathabx}
 \usepackage{textcomp}
 \usepackage{subcaption}
 \usepackage{caption}
 \usepackage{color}
\usepackage{tikz-cd}
 \usepackage[backref=page]{hyperref}
 \usepackage{xypic}
\usepackage{mathrsfs}
\newcommand{\Sph}{\mathcal{H}}
\newcommand{\SHV}{\mathrm{SHV}}

\usepackage{stmaryrd}
\usepackage{latexsym}
\usepackage[OT2,T1]{fontenc}
\usepackage{enumerate}
\DeclareSymbolFont{cyrletters}{OT2}{wncyr}{m}{n}
\DeclareFontFamily{OT1}{rsfs}{}

\usepackage{verbatim}
\usepackage{graphicx}
\usepackage{nomencl}
\DeclareFontEncoding{OT2}{}{}  

\usepackage[OT2,T1]{fontenc}
\DeclareSymbolFont{cyrletters}{OT2}{wncyr}{m}{n}
\DeclareFontFamily{OT1}{rsfs}{}
\DeclareFontEncoding{OT2}{}{}  
  
     \DeclareFontShape{OT1}{rsfs}{n}{it}{<-> rsfs10}{}
\DeclareMathAlphabet{\mathscr}{OT1}{rsfs}{n}{it}

 \newcommand{\scc}{\mathrm{sc}}
\newcommand{\Spec}{\mathrm{Spec}}
\newcommand{\Z}{\mathbb{Z}}

\newcommand{\Ext}{\mathrm{Ext}}
\newcommand{\spec}{\mathrm{spec}}
\newcommand{\Ad}{\mathrm{Ad}}
\newcommand{\uHom}{\underline{\mathrm{Hom}}}
 \newcommand{\Hom}{\mathrm{Hom}}
 \newcommand{\renormalized}{ind-finite }
 \newcommand{\unrenormalized}{ind-safe }
 
 \newcommand{\inv}{{}^{-1}}
 \newcommand{\Planch}{{\mathbb P \mathbb L}}
  \newcommand{\RTF}{{\mathbb R \mathbb T \mathbb F}}
\newcommand{\Pl}{\Planch}
 \newcommand{\Sym}{\mathrm{Sym}}
 \newcommand{\Gm}{{\mathbb{G}_m}}
 \newcommand{\GGm}{\mathbb{G}_{gr}}
 \newcommand{\Ga}{{\mathbb{G}_a}}
 
  \newcommand{\ga}{{\mathfrak{g}_a}}
  
    \newcommand{\gax}{{\mathfrak{g}_a^\ast}}
    \newcommand{\Gr}{\mathrm{Gr}}
    \newcommand{\SL}{\mathrm{SL}}
\renewcommand{\sl}{\mathfrak{sl}}
\newcommand{\Lcal}{\mathcal{L}}
\newcommand{\Pcal}{\mathcal{P}}

\newcommand{\Frob}{\mathrm{Fr}}
\newcommand{\Fr}{\Frob}  
\newcommand{\QC}{\mathrm{QC}}
\newcommand{\Shv}{\mathrm{Shv}}
\newcommand{\Diag}{\varpi}

\newcommand{\SO}{\mathrm{SO}}
\newcommand{\tw}{\mathsf{tw}}

\newcommand{\Aut}{\mathrm{Aut}}
\newcommand{\AUT}{\mathrm{AUT}}

\renewcommand{\AA}{\mathbb{A}}
\newcommand{\CC}{\mathbb{C}}
\newcommand{\PP}{\mathbb{P}}

\newcommand{\cc}{\mathfrak{c}}

\newcommand{\cg}{\check{\mathfrak{g}}}
\newcommand{\Spin}{\mathrm{Spin}}

\newcommand{\B}{\mathcal{B}}
 \renewcommand{\dim}{\mathrm{dim}}

\newcommand{\Perf}{\mathrm{Perf}}
\newcommand{\Q}{\mathbb{Q}}
\newcommand{\QQ}{\mathbb{Q}}
\newcommand{\GSp}{\mathrm{GSp}}
\newcommand{\aff}{{\operatorname{aff}}}
\newcommand{\can}{{\operatorname{can}}}
\newcommand{\val}{{\operatorname{val }}}

\newcommand{\cA}{\mathcal{A}}

\newcommand{\rG}{{\mathrm G}}
\newcommand{\IndCoh}{\mathrm{QC}^{!}}
\newcommand{\geom}{\mathrm{geom}}

\newcommand{\cY}{{\mathcal Y}}

\newcommand{\cDD}{\mathcal{D}}
\newcommand{\cNN}{\mathcal{N}}
\newcommand{\cLL}{\mathcal{L}}

\newcommand{\cHH}{\mathcal{H}}
\newcommand{\bL}{\mathbf{L}}

\newcommand{\C}{\mathbb{C}}
\newcommand{\adele}{\mathbb{A}}
\newcommand{\OO}{\mathcal{O}}
\newcommand{\et}{et}
\newcommand{\R}{\mathbb{R}}

\newcommand{\kk}{k}
\newcommand{\Zrep}{\mathcal{T}}
\newcommand{\Vect}{\mathrm{Vect}}
\newcommand{\degsheaf}{\underline{\mathrm{deg}}}

\newcommand{\PGL}{\mathrm{PGL}}

\newcommand{\Sp}{\mathrm{Sp}}
\newcommand{\autshv}{\mathrm{Shv}}
\newcommand{\autshvspec}{\mathrm{Aut}}
\newcommand{\Autshv}{\mathrm{Shv}}
\newcommand{\bigautshvspec}{\mathrm{AUT}}
\newcommand{\bigautshv}{\mathrm{SHV}}
\newcommand{\G}{\mathcal{G}}

 \DeclareFontShape{OT1}{rsfs}{n}{it}{<-> rsfs10}{}
\DeclareMathAlphabet{\mathscr}{OT1}{rsfs}{n}{it}
\newcommand{\GL}{\mathrm{GL}}
\newcommand{\vol}{\mathrm{vol}}
\newcommand{\Perv}{\mathrm{Perv}}

\newcommand{\Loc}{\mathrm{Loc}}
\newcommand{\LOC}{\mathrm{LOC}}
\newcommand{\tr}{\mathrm{tr}}

\newcommand{\dR}{\mathrm{dR}}

\newcommand{\F}{\mathbb{F}}

\newcommand{\Hecke}{\mathcal{H}}
\newcommand{\HECKE}{{\overline{\mathcal{H}}}}

\newcommand{\wt}{\widetilde}
\newcommand{\cF}{\mathcal F}
\newcommand{\cW}{\mathcal W}
\newcommand{\cX}{\mathcal X}
\newcommand{\Conf}{\mathrm{Conf}}
\newcommand{\cK}{\mathcal K}

\newcommand{\D}{\mathcal{D}}
\newcommand{\FF}{\mathbb F}
\newcommand{\zBun}{{}_z\Bun}

\newcommand{\actson}{\circlearrowright}
\newcommand{\onacts}{\circlearrowleft}
\newcommand{\GvX}{{{\check G}_X}}
\newcommand{\cG}{\mathcal G}
\newcommand{\cZ}{\mathcal Z}
 
  \newcommand{\bO}{\mathbb O}
 \newcommand{\hyperspherical}{hyperspherical }
  
   \newcommand{\bH}{\mathbb H}
\newcommand{\comod}{\mhyphen\mathrm{comod}}
\newcommand{\module}{\mhyphen\mathrm{mod}}

\newcommand{\mhyphen}{{\hbox{-}}}
\def\fgxv{ \mathfrak{\check g}^\star}
\def\fhxv{ \mathfrak{\check h}^\star}
\def\fgv{ \mathfrak{\check g}}

\newcommand{\Rep}{\text{Rep}}
\newcommand{\Bun}{\mathrm{Bun}}
\newcommand{\ad}{\mathrm{ad}}
\newcommand{\Ind}{\operatorname{Ind}}
\newcommand{\fc}{\mathfrak{c}}
\newcommand{\ff}{\mathfrak{f}}
\newcommand{\Eis}{\mathrm{Eis}}
\newcommand{\fo}{\mathfrak{o}}
\newcommand{\Map}{\mathrm{Map}}
\newcommand{\End}{\mathrm{End}}
\newcommand{\cO}{\mathcal O}
\newcommand{\cM}{\mathcal M}
\newcommand{\bs}{\backslash}

\def\ul{\underline}
\newcommand{\ot}{\otimes}
\def\Gv{{\check G}}

\newcommand{\GIT}{{/\! /}}
\newcommand{\lGIT}{\bs \! \bs}

\newcommand{\fh}{\mathfrak{h}}

\newcommand{\fg}{\mathfrak{g}}
\newcommand{\fu}{\mathfrak{u}}
\newcommand{\fgx}{\mathfrak{g}^\star}
\newcommand{\invertQM}{{\text{?`}}}
\newcommand{\cC}{\mathcal{C}}
\newcommand{\Ran}{\mathrm{Ran}}
\newcommand{\cB}{\mathcal{B}}
\newcommand{\cP}{\mathcal{P}}

\newcommand{\cS}{\mathcal{S}}

\newcommand{\Ccirc}{{{\cC}^\circ}}
\newcommand{\Coo}{{{\cC}^{\circ\circ}}}

\newcommand{\Bv}{{\check B}}
\newcommand{\Hv}{{\check H}}

\newcommand{\Mv}{{\check M}}
\newcommand{\Xv}{{\check X}}
\newcommand{\Yv}{{\check Y}}

\newcommand{\pv}{{\check p}}
\newcommand{\qv}{{\check q}}

\newcommand{\Id}{\operatorname{Id}}

\newcommand{\shear}{{\mathbin{\mkern-6mu\fatslash}}}

\newcommand{\unshear}{{\mathbin{\mkern-6mu\fatbslash}}}
\newcommand{\QCshear}{{QC^{\shear}}}
\newcommand{\Coh}{\mathrm{Coh}}

\newcommand{\LL}{{\mathbb L}}
\newcommand{\Ll}{{\mathcal L}}

\newcommand{\diff}{\mathfrak{d}}

\newcommand{\norm}{\mathrm{norm}}
\newcommand{\slice}{\mathsf{slice}}

\newcommand{\pt}{\mathrm{pt}}
\newcommand{\aut}{\mathrm{aut}}
\newcommand{\Pic}{\mathrm{Pic}}
\newcommand{\super}{{\mathrm{super}}}
\newcommand{\la}{\langle}
\newcommand{\ra}{\rangle}

\newcommand{\ssslash}{\mathbin{
  \mathchoice{/\mkern-10mu/\mkern-10mu/} 
    {/\mkern-10mu/\mkern-10mu/} 
    {/\mkern-8mu/\mkern-8mu/} 
    {/\mkern-8mu/\mkern-8mu/}}}

\newtheorem{theorem}[subsubsection]{Theorem}
\newtheorem{prop}[subsubsection]{Proposition}
\newtheorem{conj}[subsubsection]{Conjecture}
\newtheorem{conjecture}[subsubsection]{Conjecture}
\newtheorem{proposition}[subsubsection]{Proposition}
\newtheorem{expectation}[subsubsection]{Expectation}
\newtheorem{corollary}[subsubsection]{Corollary}
\newtheorem{definition}[subsubsection]{Definition}

\newtheorem{lemma}[subsubsection]{Lemma}

\theoremstyle{definition}
\newtheorem{example}[subsubsection]{Example}
\newtheorem{terminology}[subsubsection]{Terminology}

 \newtheorem{remark}[subsubsection]{Remark}
 \newtheorem{problem}[subsubsection]{Problem}

\counterwithin{figure}{subsection}
\counterwithin{table}{subsection}
\numberwithin{equation}{section}

\indexsetup{othercode=\label{index}}

\makeindex

\begin{document}

\author{David Ben-Zvi, Yiannis Sakellaridis and Akshay Venkatesh}
\title{Relative Langlands Duality} 
 \begin{abstract}
      We propose a duality in the relative Langlands program.
  This duality
pairs a Hamiltonian space for a group $G$ with a Hamiltonian space under its dual group $\check{G}$,
and recovers 
  at a numerical level the relationship between a period on $G$ and an $L$-function attached to $\check{G}$;
  it is an  arithmetic analog of the electric-magnetic duality of boundary conditions in
  four-dimensional supersymmetric Yang--Mills theory.
  \end{abstract}

 \maketitle

%{\color{red} This is an updated version. The section numbering coincides with the original version posted online in July 2023, except that the subsections of \S 2 have been renumbered and some subsections added.  If we have failed to attribute or properly reference a work it is most likely due to either ignorance or forgetfulness - please tell us.}

\newpage

 {\tiny 
 \tableofcontents
 }

\newpage

\input{setup}

\newpage
\part{Structure theory} \label{part:structure}

\input{spherical}

\input{rationality}

\part{Local  theory}
\input{shearingPart2}

\input{local-geometric}

\input{local-geometric-Plancherel}

\input{local-numerical}

\newpage

\part{Global theory}

\input{global-geometric}

\input{global-geometric-2}

 \input{global-geometric-P1}

\input{global-numerical}

\newpage
\part{Local-to-global aspects}

 \input{automorphic-factorization}

 \input{spectral-factorization}

\newpage
\part{Appendices}
\appendix
\input{shearing}

\input{sheaf-theory}

 \input{geometric-Langlands}

\input{fieldtheory}

\input{MiscellaneousComputations}

\bibliographystyle{alphaurl}
\bibliography{BZSVpaper}

\printindex
\end{document}

%% file: setup.tex
\newcommand{\iX}{(i)\textsubscript{X}}
\newcommand{\iiX}{(ii)\textsubscript{X}}
\newcommand{\iiiX}{(iii)\textsubscript{X}}
\newcommand{\ivX}{(iv)\textsubscript{X}}

\newpage

\section{Introduction}  

 One of the fundamental properties of automorphic forms is that, when integrated against certain distinguished cycles or distributions, they give special values of $L$-functions.
 The study of these integrals, or ``periods,'' 
 has a long history starting at least with the 1937 work of Hecke \cite{Hecke}.
Now, the Langlands program posits that automorphic forms correspond to Galois representations, and Hecke's formulas
 and their sequels can be expressed as a commutative diagram:
\begin{equation} \label{match diagram}
 \xymatrix{
 \mbox{automorphic forms} \ar[rr] \ar[rd]^{\text{period}}  & &  \ar[ld]^{\text{L-function}}\mbox{Galois representations} \\
 & \mbox{complex numbers}.
 }
 \end{equation}

 That is to say, ``periods'' and ``$L$-functions'' are specific ways to extract numerical invariants
 from the two sides of the Langlands program; and in interesting cases, they
 match with one another. 
 We are going to propose that:
 
 \begin{itemize}
 \item[(i)] We should index  both periods and $L$-functions by suitable {\em Hamiltonian spaces} -- these are, in particular, symplectic algebraic varieties
 with a group action.  
 \item[(ii)] The passage from Hamiltonian space to period or $L$-function can be considered as (an incarnation of) quantization.
 \item[(iii)]  ``Relative Langlands Duality:'' when viewed from this point of view, the relationship between $L$-functions and periods becomes symmetric.
   \item[(iv)]  Similar structures exists at all ``tiers'' of the Langlands program (global, local, geometric, arithmetic, etc.). In the local tier, point (ii) is familiar: it is 
   the philosophy that one can construct representations of Lie groups by quantization. 
   \end{itemize}

 For example, two of the earliest and best-studied periods
 are the Godement--Jacquet integral  and the Rankin--Selberg integral, both for the groups $G=\GL_n \times \GL_n$.
 They are switched under the duality of (iii); let us   sketch what this means in the global context.
  Firstly, to these periods will be associated certain Hamiltonian spaces
$M_{GJ}$ and $M_{RS}$ for $G$, explicitly
 $$ M=T^*\left[G \times^{\Delta \GL_n} \mathbb{A}^n \right]  \mbox{ and } M_{RS} = T^*\left[  n \times n \mbox{ matrices}.\right].$$

Now, the group $G$ being self-dual,
it plays a symmetric role on the automorphic and spectral (Galois) sides of the Langlands program.  
We will see that, when considered on the automorphic side, the data $M_{GJ}$ or $M_{RS}$ index  periods -- distributions on the space of automorphic forms;
 and when placed on the Galois side they index certain $L$-functions. For $\phi$ an automorphic form on $G$ we will then have
\begin{equation} \label{GJeq1} \mbox{$M_{GJ}$-period of  $\phi$}= \mbox{ $L$-function for  $\phi$ indexed by $M_{RS}$.}\end{equation}
\begin{equation} \label{GJeq2}  \mbox{$M_{RS}$-period of  $\phi$}= \mbox{$L$-function for  $\phi$ indexed by $M_{GJ}$.}\end{equation} 
   
How to turn translate these words into the more standard
language of periods  will be explained in  Part 3 of this paper. When thus translated, 
\eqref{GJeq1} and\eqref{GJeq2} are exactly the results proved in \cite{GJ, JS1, JS2};
however, what is made clear by this phrasing is that there is a duality between the two results. 
Thus, we can regard the ``$L$-function for $\phi$ indexed by $M_{RS}$'' 
as the ``spectral period for $M_{RS}$.''

The duality is perhaps more visible when the same statements are formulated in the context of geometric Langlands. 
Writing $M_{GJ} =T^* X_{GJ}$ and $M_{RS} = T^* X_{RS}$, 
we can consider two constructions:

\begin{itemize}
\item Automorphically, we consider $G$-bundles on a curve together with a section of the
associated $X_{GJ}$-bundle. We push forward the constant sheaf on this space
to moduli of $G$-bundles.

\item Spectrally, we consider $G$-local systems on a curve together with a 
flat section of the associated $X_{RS}$-bundle. We push forward the dualizing sheaf
on this space to moduli of $G$-local systems.
\end{itemize}

The geometric analogue of \eqref{GJeq1} is now that the constructions of (i) and (ii)
should match with respect to the conjectural geometric Langlands equivalence;
\eqref{GJeq2} corresponds to switching the roles of $X_{RS}, X_{GJ}$ above.

 A simpler but conceptually significant example is the duality of Whittaker and trivial periods. Manifestations
 of this duality account for the central role of the Whittaker model in many contexts; from our point of view
 it should not be distinguished from other dual pairs of periods. A more interesting example
 is the duality between the Gan--Gross--Prasad period \cite{GGP} and
 the $\theta$-correspondence between equal rank orthogonal and symplectic groups \cite{Rallis}.
 Other examples are given in \S \ref{IntroEx}.
 
 Although much of the ultimate payoff may be  
investigating periods whose duals are {\em not} currently known,
our focus here is to formulate carefully the duality at least in a certain very well-behaved setting, spell out what it predicts,
and see how it unifies a large class of phenomena in the Langlands program.

{\em Overview of the introduction:}  The complexity of our situation  warrants a longer explanation 
 of what we are trying to do, and how we think about the situation. 
 To that end,  we will first review in \S \ref{whatis} various  objects that enter into the relative Langlands program.
  These objects
  are the analogues of periods/$L$-functions in other ``tiers'' of the Langlands program. 
  The discussion of \S \ref{whatis} organizes the situations of the paper
  into a tesseract, rendered in figures \ref{autdiagram} and \ref{specdiagram}. 
These figures invoke a great variety of different mathematical structures
and therefore may seem very confusing at first;  
fortunately, the ideas of quantum field theory give
rise to a very appealing  metaphor with which to organize them, as we explain in  \S \ref{TFTintro} and \S \ref{AFT2} -- skip forward to   Table \ref{tab:output} to get a sense of what this discussion aims at. Finally, in \S \ref{aims}, we will spell out what we actually accomplish in this paper. 

\subsection{What is a ``period?''} \label{whatis}
 As mentioned above, the current subsection \S \ref{whatis} spells out some of the data on the
 automorphic  and spectral sides of the relative Langlands program.
On both sides, the data will be organized into diagrams -- Figures \ref{autdiagram} and \ref{specdiagram} -- and the proposal,
  is, of course, that these diagrams should match.  
More organized ways of thinking about these diagrams will be given in subsequent subsections.
 
\subsubsection{The automorphic side of the Langlands program.} \label{whatisAut} 

  Take $G$ a reductive group --- split, for simplicity. 

Some key objects of study in the  Langlands program are:

 \begin{itemize}
 \item the trace formula for $G$ (output a complex number);
 \item automorphic functions for $G$ (a vector space); 
 \item representations of $G(\R)$ or $G(\Q_p)$ (a category).
 \end{itemize}
 
 The geometric Langlands program, in the setting of curves over algebraically closed fields, 
 adds to these objects ``categorified'' analogs:
  \begin{itemize}
 \item automorphic sheaves for $G$ -- that is, sheaves on the moduli $\Bun_G$ of principal $G$-bundles on a curve (a category);
 \item categorical representations of the loop group $LG$ (a 2-category). 
 \end{itemize}

 There is a rich web of interconnections between these objects,
and we will examine how to organize all this information in  \S \ref{TFTintro} and \S \ref{AFT2} (see in particular Table \ref{tab:output}).  
 We preview the discussion of \S \ref{AFT2} in particular by noting that the ``categorical complexity'' of the output is tied to the ``arithmetic dimension'' in a way reminiscent of topological field theory: global fields (which should be considered as $3$-dimensional objects) are assigned vector spaces, arithmetic local fields and geometric global fields (of dimension $2$) are assigned categories, and geometric local fields (of dimension $1$) are assigned 2-categories.

 A core goal of the Langlands program is to  give a ``dual description'' of all 
 these objects in terms of Galois theory and the dual group $\check G$,  compatible with
 these various connections. In
 all cases the key player is the space of representations of a Galois group into $\check G$. 
 
 \subsubsection{Relative Langlands: Automorphic Side} \label{whatisRelAut}
Now let $X$ be an algebraic variety with $G$-action. It produces objects of study in the various ``tiers'' of the Langlands program, whose
 unified study forms the topic of the relative Langlands program. Some of the most familiar are:
 \begin{itemize} 
  \item the $X$-theta series, or $X$-period functional (a vector in the space of automorphic functions), and

\item functions $L^2(X(F))$, for $F$ a local field (a representation of $G(F)$, and a classical object of harmonic analysis) 
  \end{itemize}

These objects are often studied through their ``squares'' or self-pairings, after the insertion of Hecke operators:
\begin{itemize}
  \item the relative trace formula for $X$, the self-pairing of the $X$-theta series (output in $\C$, or generally functional on the global Hecke algebra);
  \item the Plancherel formula for spherical functions on $X$, obtained from the self-pairing of the basic spherical vector with Hecke modifications (output a functional on the spherical Hecke algebra).
\end{itemize}

As we detail in this paper, the geometric Langlands program also admits categorified counterparts of these objects, which have not been studied before in a uniform fashion:    
\begin{itemize}
   \item the $X$-period sheaf (an object in the category of automorphic sheaves);   \item the category of sheaves on $LX$ (a categorical representation of $LG$);
   \end{itemize}
as well as ``squares'' or self-pairings: 
\begin{itemize}
\item the {\em RTF algebra:} an associative algebra in the global Hecke category (see \S \ref{L-observables}) which conjecturally encodes maps between Hecke functors applied to the period sheaf;
   \item the {\em Plancherel algebra:} a commutative algebra\footnote{As we discuss at length, the Plancherel algebra is a derived object which is commutative only on the level of cohomology, and its structure is closely related to that of {\em little 3-disc algebras} in topology.} in the spherical Hecke category, which encodes all maps between Hecke functors applied to the basic spherical sheaf on $LX$. 
\end{itemize}

 These objects and some of their relations are captured schematically in the diagrams of Figure \ref{autdiagram} below. Here, as elsewhere in the paper, we focus on the {\em unramified} settings  (arithmetic and geometric) over function fields -- see Remark \ref{why no nf}.  In particular, we replace the $G(F)$-representation $L^2(X(F))$ by its $G(O)$-invariants, a module for the spherical Hecke algebra.
 Our conventions for these diagrams are:
 \begin{itemize}
 \item the top row lives in geometric Langlands, the bottom row lives in arithmetic Langlands
 over a finite field; the left hand column is local, and the right column is global.
 \item vertical arrows come from ``trace of Frobenius.'' 
Dashed  horizontal arrows    suggest that the right-hand object
has the nature of an ``Euler product''  or ``integrated version'' of the left-hand object.\footnote{In each case, there are precise statements
formalizing this, but we don't need this level of precision here.}
\item The diagrams come in pairs, which have been termed as ``states'' and ``observables.'' In each case,
the ``observable'' diagram arises from self-pairings or endomorphisms of the ``states'' diagram, informally:
$$ \mbox{observables} = \langle \mbox{states}, \mbox{states} \rangle,$$
 The nomenclature arises from the analogy between
these concepts and states and observables in quantum mechanics (or between {\em geometric} and {\em deformation} quantization), a recurrent theme in the paper. 
See also Remark \ref{square of relative Langlands}. \index{states} \index{observables}
 
 \end{itemize}

  \begin{figure}
  \centering
\begin{subfigure}{0.5 \textwidth}

\centering
\caption{{\em Automorphic States:}}

 {\small
 \xymatrix@C=1cm{
		&	 \txt{Local 	}								
		& \txt{Global} 					\\
{\tiny \text{Geometric}}	& *++[F-,]\txt{Sheaves on $X_F/G_O$ \\  \S \ref{section-unramified-local} }  \ar[d]  			
& *++[F-,]\txt{$X$-period sheaf \\  \S \ref{section:global-geometric} esp. \S \ref{periodX} }  \ar[d]		 			
			\\
 {\tiny \text{Arithmetic}}
 			&*++[F-,]\txt{Functions on $X_F/G_O$\\ \S  \ref{Hecke module structure}}   	
 			& *++[F-,]\txt{$X$-period function\\   \S \ref{section:global-geometric} esp. \S \ref{periodX}}		
			\\
}
}
 \end{subfigure}

\vspace{1cm}
\begin{subfigure}{0.5 \textwidth}
\caption{ {\em Automorphic Observables:}}
 {\small
 \xymatrix@C=1cm{
		&	 \txt{Local 	}								
		& \txt{Global} 					\\
{\tiny \text{Geometric}}	
			&  *++[F-,]\txt{Plancherel algebra, \\ \S \ref{PlancherelCoulomb}}    \ar[d]  \ar@{.}[r]^{\S \ref{automorphic-factorization}}
			& *++[F-,]\txt{$X$-RTF algebra \\   \S \ref{automorphic-factorization},  \S \ref{factorizable Theta}}  \ar[d]		 			
			\\
 {\tiny \text{Arithmetic}}
 			&  *++[F-,]\txt{Plancherel measure \\ \S \ref{section:numerical Plancherel} }\ar@{.}[r]^{\textrm{\cite{SV}}}
 			& *++[F-,]\txt{ (a part of) $X$-RTF\\ \S \ref{RTF discussion} 
 }		
			\\
}
}
\end{subfigure}
\caption{Automorphic states and observables}
 \label{autdiagram}
\end{figure}

\begin{remark}[Number fields and ramification]  \label{why no nf}
These diagrams cover many interesting tiers of the Langlands program, but also
omit many important ones. For example, we do not discuss the case of number fields,
local fields of characteristic zero,  or ramification.

However, we wish to emphasize here that we anticipate
 {\em the general picture  should apply equally to these cases}. 
 That is, as we see it, the story we tell here is not one specific to
 geometric Langlands or to the case of function fields, but a general feature
 of the paradigm of Langlands duality.  Incorporating number fields and ramification would
 require the development of the local theory over an archimedean local field, and for ramified representations over nonarchimedean local fields; there should be many similarities between these two cases. 
To examine these questions is therefore a very interesting open problem. See also Remark \ref{remark-global-ramified}.

 \end{remark} 

An important principle underlying this work is the {\em microlocal} nature of periods. One manifestation of this principle is that many of the above structures can be thought of as associated not to the $G$-variety $X$, but to the {\em Hamiltonian $G$-variety} $M=T^*X$. 
For example, the question of reconstructing $L^2(X(F))$ from $T^* X$
has a standard name: it is the question of {\em geometrically quantizing} the symplectic space $T^*X$. This passage from $G$-varieties to their cotangents as Hamiltonian $G$-spaces encodes important symmetries: for example, in Tate's thesis (independently due to Iwasawa: \cite{Iwasawa, Tate}) the functional equation of abelian $L$-functions derives from the Fourier transform, corresponding to a symmetry of the Hamiltonian $G=\Gm$-space $T^*\AA^1$.

Another manifestation of the microlocal nature is that there are many important examples of periods (such as  the Whittaker periods and the $\Theta$-correspondence) that carry structures that resemble the above, but arise not from a $G$-space $X$ but rather from a Hamiltonian $G$-space $M$ which is not a cotangent bundle. We shall think broadly of the passage from $M$ to the list of data above as a form of quantization.
 To avoid getting into it at the moment, however, 
  let us continue to work with the ``polarized case'' where one has a $G$-variety  $X$ (or its cotangent bundle $T^*X$) rather than a general Hamiltonian action.   

\subsubsection{Relative Langlands: Spectral Side} \label{whatisRelAutSpec}
We now wish to propose a spectral, or Galois, side for the relative Langlands program that matches the structures enumerated above.
The experience in the study of periods is that the collection of $X$ for which
one has a satisfactory theory has rather large overlap with the
 set of {\em spherical varieties} for $G$.  
We shall eventually extend our setting to include
a slightly larger class of examples, the ``hyperspherical varieties'' of \S \ref{sphvar}. In any case the main goal of this paper is to propose that 
\begin{quote}
For favorable $G$-spaces $X$ (or Hamiltonian $G$-spaces $M$), the various structures of the $X$-relative Langlands program are simultaneously encoded, on the dual side, by 
a Hamiltonian $\check G$-variety $\Mv$. \end{quote}

We have already asserted that the passage from $M$ to its associated automorphic data
 should be seen as a quantization. Similarly, the passage from $\check{M}$
 to the associated spectral data should be seen as a quantization.
 If we want to distinguish the two, we will say ``automorphic quantization''
 and ``spectral quantization''. Thus the philosophy pervading the paper, sometimes
 implicitly and sometimes explicitly, can be summarized as saying 
 \begin{quote}
The automorphic quantization of a favorable Hamiltonian $G$-variety is Langlands dual to the spectral quantization of a dual Hamiltonian $\Gv$-variety $\Mv$. \end{quote}

The notion of spectral quantization incorporates a new geometric perspective on the notion of $L$-functions of Galois representations, mirror to the microlocal perspective on periods: we propose to think of an $L$-function $L(\rho,V)$ as associated to the Hamiltonian $\Gv$-variety $T^* V$ by a form of geometric quantization.  This is not completely achieved in this paper, but ideally speaking the situation is as follows: 
\begin{itemize}
\item We first attach $L$-functions to $\Gv$-varieties $\Xv$, rather than linear representations\footnote{A motivation for this passage, from a classical point of view, is that the local factor of an $L$-function has the form of the trace of Frobenius on a symmetric algebra -- i.e., on the ring of functions on a vector space. This suggests, at leasts, that what one needs is not the vector space but only its function ring.} of $\Gv$--- the familiar $L$-functions arise by linearization at Galois fixed points on $\Xv$ as
in e.g.\ \eqref{PXshriekprediction}. 
\item We then view the $L$-functions as attached not to $\Xv$ but to the Hamiltonian $\Gv$-variety $T^*\Xv$. This encodes additional symmetries, such as the functional equation.
\item Finally we seek to attach $L$-functions to Hamiltonian $\Gv$-varieties that are not necessarily polarized -- in particular, to extract square roots of $L$-functions associated to symplectic representations (see, e.g., \eqref{L^2 normalized equation}, where we leave the signs unspecified, and compare with \cite{AV}). 
\end{itemize}

As with automorphic quantization, spectral quantization manifests itself in every tier of the Langlands program, which we predict match the corresponding automorphic objects under Langlands duality.
 An informal flavor of what is studied in this paper is captured by matching the two automorphic diagrams above with two spectral counterparts
 in Figure  \ref{specdiagram} below,
 which will follow the same conventions as for the previous diagrams. \footnote{It might be better to picture two matching cubes, or a single tesseract. We leave it as an open problem to develop better visualizations of this dictionary.} 
 To avoid getting into details of quantization we restrict to the case when $\check{M}$
 is polarized, i.e., $\Mv=T^*\Xv$.

\medskip

\begin{figure}[h] 
\centering

\begin{subfigure}{0.5\textwidth}
 \caption{\em{Spectral states}}
 {\small
 \xymatrix@C=1cm{
		&	 \txt{Local }				
			& \txt{Global} 					\\
{\tiny \text{Geom.}}	& *++[F-,]\txt{coh. sheaves on \\  $\Mv/\Gv$, \S \ref{section-unramified-local} }  \ar[d]   
			& *++[F-,]\txt{$L$-sheaf for $\check{X}$, \\ \S \ref{Lsheaf}}  \ar[d]		 			
			\\
 {\tiny \text{Arith.}}
 			&*++[F-,]\txt{functions on Frobenius\\fixed points of $\check{M}$, \S \ref{Hecke module structure}} 
 			& *++[F-,]\txt{$L$-function for $\check{X}$ \\  \S \ref{Lsheaf} }		
			\\
}
}
\end{subfigure}
 
 \vspace{1cm}

  \begin{subfigure}{0.5 \textwidth}
\caption{{\em Spectral Observables}}
 {\small
 \xymatrix@C=1cm{
		&	 \txt{Local }								
		& \txt{Global} 					\\
{\tiny \text{Geom.}}
			&  *++[F-,]\txt{regular functions \\ on $\Mv$ }    \ar[d]  \ar@{.}[r]^{\S \ref{spectral-factorization}} 
			& *++[F-,]\txt{$L$-observable algebra, \\ \S \ref{case5}, \S \ref{L-observables}}  \ar[d]		 			
			\\
 {\tiny \text{Arith.}}	
 			&  *++[F-,]\txt{twisted character  \\ of $\check{G}$ on $\Mv$ \S \ref{section:numerical Plancherel} }\ar@{.}[r] 
 			& *++[F-,]\txt{$L$-function for $\Mv$ \\  \S \ref{L2numerical} esp. \eqref{L^2 normalized equation}    }		
			\\
}
}
 
 \end{subfigure}
 \caption{Spectral states and observables}
 \label{specdiagram}

\end{figure}

New in this diagram, and
of particular importance to this paper, are the {\em $L$-sheaves}, new counterparts of $L$-functions in the setting of geometric Langlands. These are objects of the categories of sheaves on moduli of Langlands parameters given by considering algebraic volume forms on fixed point spaces of Galois representations on $\Xv$. 
They have not been systematically studied in the literature, but in specific cases have been considered, most notably the ``Whittaker sheaf'' and  work of Lysenko (e.g.\ \cite{Lysenko1, Lysenko2}).

  \subsubsection{Matching automorphic and spectral} \label{whatisRelAutMatching}
  
   In Figures \ref{autdiagram} and \ref{specdiagram} we have summarized some data on
   the automorphic and spectral side of the relative Langlands program. The proposal, of course, 
is that each automorphic square should match with the corresponding spectral square!
  
  Precise conjectures
 to this effect are contained in the text (most importantly in \S \ref{section-unramified-local}, \S \ref{GGC}, \S \ref{L2numerical} and \S \ref{local-global}).  
In the arithmetic setting -- i.e., in the bottom row of both squares --
many of these conjectures, for specific choices of $X, \check{X}$,  are theorems in the relative Langlands program.
 In the geometric setting they are mainly conjectural; we 
regard the lower row as being evidence for them on the level of Frobenius traces.

\begin{remark} [Squares and square-roots]\label{square of relative Langlands}
Note that we see some appearances of $\check{X}$ and some of $\check{M}$.
 This is related to a more general phenomenon, which will occur throughout this paper:
frequently the ``square'' expressions (the observables) on both sides
  can be formulated with reference only to $M$ or $\check{M}$, but to extract their ``square-roots'' (the states) requires extra structure.
Correspondingly, the  ``squared'' version of the
conjectures can be formulated without reference to extra choices such as polarization.
This reflects the fact that the geometric quantization $V$ 
of a symplectic manifold often requires choosing some extra structure (such as a polarization),
but its ``square'', i.e., $V \otimes \check{V}$, is much more directly related to the symplectic manifold itself and its {\em deformation} quantization.

For example,  in relation
to the square of a global period, such ideas will be a running thread in this paper, discussed in
 \S \ref{Lindep}, 
 \S \ref{case5}, \S \ref{L2numerical}, 
 and \S \ref{spectral-factorization}.
\end{remark}

\subsection{Electric-Magnetic Duality and Topological Field Theory} \label{TFTintro}
 The picture that we just described has been very much influenced by ideas of topological field theory. We will now
 explain this point of view in more detail,  before returning to our main goals in \S \ref{aims}. 
 At the moment (\S \ref{TFTintro}) we will use the language of topological field theory, and then
 in \S \ref{AFT2} we will relate this language to our arithmetic setting.
 In order to understand the utility of the language of TFT for representation theory, we strongly recommend the toy model provided by finite group gauge theory, a synopsis of which is found in Appendix \S \ref{finite group TFT}.

The seminal work of Kapustin and Witten~\cite{KapustinWitten} and subsequent developments demonstrated that the geometric Langlands correspondence can be profitably seen through the lens of four-dimensional topological
quantum field theory (TFT),  \index{TFT} specifically as an aspect of electric-magnetic duality in gauge theory. 
 The idea of TFT, although formally inapplicable, 
remains a powerful metaphor to structure discussion of the arithmetic Langlands program, as we shall see in \S \ref{AFT2},
and for this reason we will review it now.

Recall that (in the mathematical viewpoint)
an extended $4d$ TFT is a representation of the higher category of bordisms of manifolds of dimension $\leq 4$, with its symmetric monoidal structure given by disjoint union.
It assigns data to manifolds of dimension $\leq 4$, as follows: 
{\small
\begin{itemize}
\item $4$-manifolds are assigned numbers;
\item $3$-manifolds are assigned vector spaces; \\ \qquad a bounding $4$-manifold produces a vector in the space.
\item $2$-manifolds are assigned categories; \\  \qquad a bounding $3$-manifold produces an object in the category.
\item $1$-manifolds are assigned $2$-categories; \\
\qquad a bounding $2$-manifold produces an object of the $2$-category. 
\end{itemize}}

In practice only the invariants associated to manifolds of dimension less than 4 can be made to fit this rigid algebraic formalism; defining {\em partition functions}, numerical invariants of 4-manifolds (possibly with boundaries or corners) requires analytic renormalization procedures, which are closely analogous to the issues one encounters  in handling the trace formula. 

Before proceeding any further let us note the analogy with the data on the automorphic and spectral sides of the Langlands program:

Kapustin and Witten study a specific pair of 4d TFTs
 sometimes called the 4d A-model and B-model, associated to a compact Lie group (or its complexification, a complex reductive group $G$). A
 special case of a fundamental conjecture in gauge theory, ``electric-magnetic S-duality of $\mathcal{N}=4$ super-Yang-Mills'',   implies an equivalence of 4d TFTs $$\cA_G\simeq \cB_{\Gv}$$ associated to Langlands dual groups.\footnote{This can be considered as a 4-dimensional counterpart to the celebrated mirror symmetry identifying the 2d A-model associated to a real symplectic manifold and the 2d B-model of a holomorphic mirror manifold.} When evaluated on a 2-manifold $\Sigma$ we obtain an equivalence of categories $\cA_G(\Sigma)\simeq \cB_{\Gv}(\Sigma)$, which Kapustin and Witten then interpret as identifying the automorphic and spectral sides of the geometric Langlands correspondence. This gauge theoretic perspective on geometric Langlands has been extensively developed in subsequent physics papers, including~\cite{GukovWitten,FrenkelWitten,Wittenwild, Witten:2009mh,Witten6, GaiottoSduality,FrenkelGaiotto,Wittenmore} (see the reviews~\cite{WittenKWreview,Frenkelgauge}) and, to a lesser extent, in the mathematics literature, see in particular ~\cite{Elliott:2015rja, BettiLanglands,ElliottYoosingular,elliottgwilliamwilliams}.
 
The structure of a topological field theory imposes strong relations between its outputs, and can be used to constrain  -- and sometimes characterize --  the values on higher dimensional manifolds in terms of those of lower dimension.  A crucial part of this package is the study of {\em defects} --- roughly speaking, field theories living on embedded submanifolds of spacetime. Defects of dimensions 2,1, and 0, known respectively as surface, line and local operators, account in the physical interpretation for much of the rich structure of the geometric Langlands program (ramification~\cite{GukovWitten}, Hecke operators~\cite{KapustinWitten} and singular support~\cite{ElliottYoosingular}, respectively). 
The duality $\cA_G\simeq \cB_{\Gv}$ of field theories implies a duality for defects of all dimensions, which in the case of line operators was interpreted by Kapustin and Witten as recovering the geometric Satake correspondence. 

The most subtle and interesting defects are those of dimension 3 (i.e., of codimension one) -- these form {\em boundary theories} or more generally {\em domain walls} or interfaces, the natural class of morphisms between TFTs. Indeed a crude paraphrase of the Cobordism Hypothesis~\cite{BaezDolan,jacobTFT} (see also~\cite{Freed:2012hx,Ayala:2017wcr}) is that the collection of domain walls completely determines a fully extended topological field theory.
 Among the many structures induced by a boundary theory, we may now formally view an $n+1$-manifold $N$ ($n<4$) with marked boundary as a new closed manifold, and every $n$-manifold $M$ as bounding a new $n+1$-manifold ($M\times [0,1]$ with $M\times \{0\}$ marked by the boundary theory), thus defining a distinguished element of the invariant of $M$.  Therefore, a TFT together with a choice of boundary theory induces the following data: 
{\small
\begin{itemize}
\item a $4$-manifold with marked boundary is assigned a number;
\item A closed $3$-manifold acquires a vector in its vector space; \\ \qquad a $3$ manifold with marked boundary
is now assigned a vector space. 
\item A closed $2$-manifold acquires an object in its category; \\ \qquad
a $2$-manifold with marked boundary is now assigned a category, etc.
\end{itemize}}

 The study of boundary theories in the physical setup for the geometric Langlands correspondence, and in particular the effect of electric-magnetic (Langlands) duality on such boundary theories, was pioneered by Gaiotto and Witten~\cite{GaiottoWittenboundary,GaiottoWittenSduality}, explored further in~\cite{GaiottoSduality,FrenkelGaiotto}, and studied mathematically recently  by Hilburn and Yoo~\cite{HilburnYoo}.

The study of boundary theories is a natural generalization of geometric quantization -- by which a symplectic manifold $M$ defines a quantum mechanical theory, i.e., 1-dimensional quantum field theory -- and the orbit method -- by which a Hamiltonian $G$-action on $M$ upgrades this quantum mechanics to a representation of $G$, i.e., a boundary theory for a suitable 2d QFT (a gauge theory with gauge group $G$). In the higher-dimensional setting of Gaiotto-Witten, the quantization of a {\em hyperk\"ahler} manifold $M$ defines a 3d QFT (an $\mathcal{N}=4$ supersymmetric sigma model), while an action of $G$ by isometries upgrades this 3d QFT to a boundary theory for a 4d $G$-gauge theory (4d $\mathcal{N}=4$ super-Yang-Mills). If we are interested in the underlying {\em topological} field theories we only need $M$ as a holomorphic symplectic manifold, and the $G$-action as a holomorphic Hamiltonian action. (On the B-side the 3d TFTs are the well-studied Rozansky-Witten theory, while the A-side theory -- the 3d analogs of the A-model in mirror symmetry -- are much less understood.)

Thus the work of Gaiotto and Witten suggests that
\begin{itemize}
\item one should study a higher form of geometric quantization for holomorphic Hamiltonian $G$- and $\Gv$-spaces which outputs boundary conditions for the topological field theories $\cA_G$ and $\cB_\Gv$, respectively; and
\item the electric-magnetic duality $\cA_G\simeq \cB_{\Gv}$ implies an identification of boundary conditions for the two theories, which we might expect to lead to a correspondence between Hamiltonian actions of the two dual groups.
\end{itemize}
Not all boundary theories come from Hamiltonian actions. However, any boundary theory can be approximated by one coming from a Hamiltonian action -- the ``Higgs branch'' of its moduli space of vacua -- and it also gives rise directly to a Hamiltonian action of the dual group -- the ``Coulomb branch'', mathematically described by the ``relative'' variant~\cite{BFNring} of the construction of Braverman-Finkelberg-Nakajima~\cite{BFN}.

Thus the electric-magnetic duality for boundary theories of Gaiotto and Witten suggests a correspondence (or partially defined duality) between Hamiltonian actions of Langlands dual groups. 
We propose in this paper that the dual of a boundary condition coming from a {\em hyperspherical} action is again of the same type, so that we obtain a duality between hyperspherical varieties for Langlands dual groups.   {\em In fact, we are proposing that the hyperspherical actions have the feature 
that the ``approximation'' alluded to above is exact -- the boundary theory is entirely determined by the Hamiltonian action.  } Moreover we explain how the structure of boundary theories for the field theories $\cA_G$ and $\cB_\Gv$ provides both a geometric counterpart to the relative Langlands program and a new perspective on the theory of periods and $L$-functions.
 
\begin{remark}[Gaiotto-Witten data, Nahm poles and Whittaker induction]\label{GaiottoWitten intro}
In Part 1 of this work, we prove a general structure theory for hyperspherical varieties. 
The data entering into this structure theorem, on the one hand, aligns
with invariants studied in the Langlands program (cf. Remark \ref{WhittArth0}, \S \ref{NTmov}),
but is  also closely related to the indexing of Gaiotto and Witten, as we will explain 
(in reading the following, it will be helpful to be familiar with the notation of \S \ref{sphvar}). 

\index{Nahm pole}

 Gaiotto and Witten describe $1/2$ BPS boundary conditions as associated to a triple $(\rho,H,\cZ)$ of data, where $\rho$ is an $\mathfrak sl_2$-triple in $G$, $H\subset G$ is a subgroup and $\cZ$ is a 3d superconformal field theory with $H$-symmetry. These data line up precisely with the description we give for hyperspherical varieties. First, Gaiotto-Witten used the datum $\rho$ alone to specify a new family of {\em Nahm pole} boundary conditions. These correspond to the Whittaker period $T^*G\GIT_\psi U$ we associate to $\rho$ -- more precisely, our construction corresponds to the topologically twisted form (of either $A$ of $B$-type) of this boundary condition, when considered on two-manifolds times $\mathbb R$ using the $R$-symmetry provided by the $\GGm$-action. (See also Remark ~\ref{domain wall remark}.)  
   Next, the special case when $\cZ$ is a theory of free hypermultiplets -- the $\sigma$-model into a vector space, which is a quaternionic representation of $H$ -- corresponds in our language to the periods specified by the underlying complex symplectic representation. Finally the full data $(\rho,H,\cZ)$ labels what we study (in this vectorial case) as the Whittaker induction of $\cZ$. In other words, the boundary conditions we study -- Whittaker-inductions $M$ from $H$ to $G$ of special symplectic representations -- are natural special cases of the Gaiotto-Witten classification. 

A key contribution of our work, as paraphrased in~\cite{gaiottosphere}, is
to isolate within this broader class the {\em hyperspherical} varieties.  These appear to form a subclass of boundary conditions
that are exceptionally well behaved, see the subsequent remark.
Our prediction is that, when the assembled Hamiltonian space $M$ is {\em hyperspherical}
and anomaly-free, then the dual boundary condition is  also hyperspherical and anomaly--free,  i.e., of the same special form, and no exotic CFTs appear on either side of the duality. 
\end{remark}
 
\begin{remark}  \label{WhyHyp} Why hyperspherical?

As just noted, a key feature of this work is the notion of hyperspherical variety. 
 We arrived at this by abstracting the role of ``multiplicity one'' in the Langlands program. 
 
 Namely,
 a key guiding principle in the study of automorphic periods is the following:
 for $X$ a $G$-space, the $X$-period of an automorphic form is closely 
 related to $L$-functions when the following multiplicity one property
 holds:\footnote{We have phrased this informally, ignoring details of functional analysis -- i.e.
 what exact space of functions on $X_F$ to choose.}
\begin{quote}
for a local field $F$, every representation of $G_F$ occurs at most once inside functions on $X_F$.
\end{quote}
 
 Abstracting this leads to at least part of the definition of hyperspherical variety (although
 we emphasize this is a simplified picture:
the multiplicity one property certainly does not hold for all $X$ such that $T^*X $ is hyperspherical.)

\end{remark}

\subsection{Arithmetic Field Theory} \label{AFT2}
We now explain that the zoo of data from \S \ref{whatis} 
 looks very structured   
 when viewed through the lens of topological field theory.
 
The central metaphor of arithmetic topology (as suggested by Mazur and developed by Kapranov, Reznikov, Morishita and others~\cite{Mazur, Reznikov, MorishitaAnalogies}) is that 
objects of number theory are analogous to manifolds of a suitable dimension. 
In our case, the ``manifolds'' of interest are 
global\footnote{More precisely, instead of talking of a global field,
we should talk about its ring of integers or a curve with this function field, but for brevity we will not do so.} and local fields, or their geometric analogues, i.e., the function fields of curves over an algebraically closed field $\FF$,
or the Laurent series field $\FF((t))$.  
In the arithmetic topology metaphor, then, global fields are akin to $3$-manifolds,
local fields and geometric global fields to $2$-manifolds, and geometric local fields to $1$-manifolds. 

 Kim~\cite{KimChernSimons} suggested the notion of 3d {\em arithmetic quantum field theory} -- namely, that one consider these ``new'' manifolds as inputs to topological quantum field theories, assigning numbers to arithmetic $3$-manifolds (global fields) and vector spaces to arithmetic $2$-manifolds (local fields). Specifically Kim and collaborators have studied arithmetic analogs of 3d TFTs including Chern-Simons theory and Dijkgraaf-Witten theories with varied applications~\cite{KimCS1,KimCS2,Kimlinking,KimBF,KimLfunctions}, see also~\cite{pappas}. 

 We are interested in applying this notion, but one dimension higher:
Observe that the various items in our discussion of the Langlands program (\S \ref{whatisAut}) loosely resemble the structure of a  {\em four dimensional }topological field theory evaluated on manifolds of different dimensions \footnote{An apparently different analogy between the Langlands correspondence and topological quantum field theory was proposed by Kapranov in~\cite{KapranovTQFT}.}
, while the items in the relative Langlands program (\S \ref{whatisRelAut}) 
 resemble the data provided by a bounding manifold -- or more precisely, the extra data of a boundary theory. 

Let us then informally say that a ``4-dimensional arithmetic quantum field theory''
is a mechanism that associates to arithmetic $j$-manifolds
(for $1 \leq j \leq 3$) vector spaces, categories or $2$-categories as appropriate,
satisfying various natural compatibilities with reference to 
the arithmetic analogue of one manifold bounding another.  Moreover, we expect such a theory to carry the rich structure of defects - local, line and surface operators as well as boundary theories.

The duality posited by the Langlands correspondence can be described
as an equivalence of two 4d arithmetic quantum field theories associated to Langlands dual reductive groups $G$, $\Gv$:
$$ \mbox{automorphic theory $\cA_G$} \simeq \mbox{spectral theory $\cB_\Gv$}.$$
Thus we think of the arithmetic correspondence as the same duality as in the Kapustin-Witten interpretation of geometric Langlands, but where we've extended the possible inputs into the arithmetic regime. 
 
The descriptions of the automorphic and spectral theory are ``mirrors'' of each other:
the automorphic theory studies the {\it topology} of moduli spaces of {\it algebraic} $G$-bundles (or arithmetic locally symmetric spaces), while the spectral theory studies the  {\it algebraic geometry} of moduli or deformation spaces of {\it topological} $\Gv$-bundles (local systems or Galois representations). 
For example, in the everywhere unramified geometric setting, the automorphic theory attaches to a smooth projective curve
the constructible sheaves on the space of $G$-bundles of $\Sigma$, and the spectral theory
attaches the coherent sheaves on the space of $G$-local systems on $\Sigma$; the basic
structure of both sides is the same, but different topologies have been used in defining sheaf theory, 
the notion of bundle, and on $G$ itself.  This parallel is much more visible in the geometric than the arithmetic
aspects of the Langlands program. Indeed one of the great advantages of the physical setting of electric-magnetic duality is that it provides complete symmetry between the two sides.  In this paper, we have tried to systematically take this perspective in studying the spectral counterpart of the theory of periods.  

 Our proposal is that the theory of periods attached to a spherical $G$-variety $X$ (or its cotangent $M$, as Hamiltonian $G$-variety) 
can profitably be viewed as defining a {\em boundary theory}
$\Theta_M$ for the automorphic field theory $\cA_G$. Informally this means for any arithmetic $j$-manifold $N$ we have a new kind of arithmetic $j+1$-manifold bounding $N$ (the product of $N$ by an interval with one end labelled by $\Theta_M$), producing an object in $\cA_G(N)$. Moreover we explain that various structures in the theory of periods fit naturally into this framework.

In Table \ref{tab:output} we explain what the automorphic theory for $G$
attaches to various arithmetic manifolds, both without and with ``boundary.''

\begin{remark}[The fourth dimension]
Note that the 4th dimension in this picture -- ``time'' -- plays a formal role and is not carrying any arithmetic structure. There are no genuine $4$-manifolds in the arithmetic analogy, but nonetheless we can build  
objects from the available arithmetic $3$-manifolds that behave like $4$-manifolds insofar
as our arithmetic TQFT goes. For example, 
  the relative trace formula and the period functionals on automorphic forms arise as partition functions of the automorphic theory on what can be considered ``arithmetic $4$-manifolds.''

  For example, the action of Hecke operators on automorphic forms associated to an arithmetic 3-manifold $M$ 
 should be considered as the invariant associated to $M\times I$ with a line defect inserted along a knot in $M$. 
The trace formula arises as the invariant associated to $M\times S^1$ (with the insertion of line defects / Hecke operators), while the 
relative trace formula for $G$-spaces $X,Y$ arises as the invariant associated to a 3-manifold times an interval, with the two ends marked by $X$ and $Y$ respectively (see also Remark \ref{square of relative Langlands}). 
\end{remark}

  {\small
\begin{table} 
\centering
\caption{ $F$ a global field, $X,Y$ spherical varieties, $F_v$ a local field; $\bar{F}$ a geometric global field e.g.\ $\overline{\mathbb{F}_q}(t)$,
$\bar{F}_v$ a geometric local field e.g.\ $\overline{\mathbb{F}_q}((t))$.}
\label{tab:output}
 \begin{tabular}{|c|c|c|c|}
\hline
  Dim. 	& ``Manifold''							&``boundary''						&  	Field theory	 				\\ 
\hline
 4  			& global   $F$						&	 sph. var. $X_i$							&    Relative trace formula($X_1, X_2$) $\in\mathbb{C}$	 \\ 
  4  			& global  $F$  						&	 sph. var. $X$ 							&   $X$-period functional	$\in \mathcal{H}$	 	 		 \\
 3  			&  global  $F$						& 								&  $\mathcal{H}$=Hilbert space of automorphic functions 	 	 \\
\hline
 3 			& local  $F_v$ 						&  $F$							&  aut. forms with ramification at $v\in \mathcal{C}$  	 \\
 3 		& local $F_v$ 								&   $X$							&   functions on $X(F_v) \in \mathcal{C}$			\\
 2			&  local  $F_v$								& 								&	category $\mathcal{C}$ of $G(F)$-representations 		\\
\hline
 
 3 			&  geom. global $\bar{F}$					&	$X$							& period sheaf $\in \mathcal{A}$						\\	
 2 			&  geom. global $\bar{F}$					&								& $\mathcal{A}=$  category of automorphic sheaves 			\\
\hline
 2			& geom.  local  $\bar{F}_v$			&	$X$							&cat. of aut. sheaves w/ ramification	$\in \mathfrak{C}$ \\
 2			& geom. local  $\bar{F}_v$		&	$\Sigma$							&cat. of sheaves on $X(\bar{F}_v)$ $\in \mathfrak{C}$	\\
 1 			&  geom. local  $\bar{F}_v$			&	 							& $2$-category $\mathfrak{C}$ of $G(\bar{F}_v)$-categories	\\
\hline
 \end{tabular}
 \end{table}
}

Likewise, a Hamiltonian $\Gv$-space $\Mv$ indexes a boundary theory $\Ll_\Mv$ for the spectral theory $\cB_\Gv$, which (in the polarized case $\Mv=T^*\Xv$) encodes the  structure of Galois fixed points on $\Xv$ and associated $L$-functions. 
 Our proposed duality (*) can be understood in terms of a meta-conjecture (formulated more precisely, in the function field case,  in Conjecture~\ref{meta-conjecture}):
 
\begin{conjecture}[Meta-conjecture]
Under the conjectural Langlands correspondence of arithmetic quantum field theories $\cA_G\simeq \cB_\Gv$, we have an identification of boundary theories
$$ \Theta_M\in \cA_G \longleftrightarrow \Ll_{\Mv}\in \cB_\Gv$$ associated to dual hyperspherical varieties $G\actson M \leftrightarrow \Gv\actson \Mv$. 
\end{conjecture}

This meta-conjecture encodes, in particular, all the matching of data between
Figures \ref{autdiagram} and Figures \ref{specdiagram}.
 
 \subsection{Aims and outline of the current paper} \label{aims}
 Informally, the aim of this paper is
 \begin{quote} {\em to put the meta-conjecture just described on a rigorous footing, at least
 in a certain subset of the phenomena it covers.}
 \end{quote}

That is: we shall try to formulate precise conjectures in both arithmetic  Langlands (everywhere unramified over a function field) and geometric Langlands in a reasonable level of generality, 
clarify their relation to existing (proved or conjectural) numerical statements in the relative Langlands program, 
 and provide where possible other supporting evidence for them. The primary contribution of this paper, then, is in finding
what we hope are  appropriate definitions and formulations. 

 To describe things more formally we need to introduce some more precise
notation.
 
 \label{subsection-evenHamiltonian}
   Let $G$ be a split connected reductive algebraic group over a local or global field $F$,
  with dual group $\check G$, which we will regard as a split reductive group over a coefficient field $k$ in characteristic zero (for example, $k=\overline{\mathbb Q_\ell}$).
 
In this paper, by a graded Hamiltonian $G$-space we will mean a smooth, symplectic variety $M$ over $F$,
equipped with an action of $G$ and a commuting ``grading'' action of $\Gm$,\   such that
\begin{enumerate}
 \item 
there is a $G$-moment map 
\[\mu: M\to \mathfrak g^*\]
(i.e., a map such that the vector field induced by $Z\in \mathfrak g$ is the Hamiltonian vector field associated to the function $M\ni x \mapsto \left<Z,\mu(x)\right>$);
 \item the action of $\Gm$ scales the symplectic form by the \emph{square} of the tautological character, i.e., 
 $$a^* \omega = a^2 \omega.$$ 
\end{enumerate}
In particular, the moment map is $\Gm$-equivariant when $\Gm$ acts by the square of the tautological character on $\mathfrak g^*$.  
The spaces $M$ that we consider will satisfy certain other conditions, involving ``parity'' and ``anomaly'', 
which we will discuss later e.g \S \ref{analyticarithmetic}, \S \ref{parity}, \S \ref{remark-quantization}.
  The $\Gm$-factor in the condition above will be a component of the ``extended reductive group'' that we will define in \ref{subsection-extended-group}. This is a group that appears naturally on the spectral side of the Langlands correspondence, and, in order to distinguish it from other instances of $\Gm$ we will denote it throughout this paper by $\GGm$, and call it the ``grading group.'' Its action on the spectral side is related to a cohomological grading in a geometric setting, 
   and to the Frobenius in an arithmetic setting.

 The paper is organized as follows. (We give here a high-level description; each
 section begins with a detailing of its contents.)
 
 \medskip
 
$\bullet$ {\bf Part 1: Structure  theory (\S \ref{sphvar} -- \S \ref{goodhypersphericalpairs}):}
In \S \ref{sphvar} we introduce the {\em hyperspherical varieties}, a convenient class of graded Hamiltonian $G$-spaces
that is closely related to the class of cotangent bundles of spherical varieties.
Formally, they are are affine Hamiltonian varieties with a suitably normalized commuting $\GGm$-action whose invariant functions Poisson-commute (for example, cotangents to affine spherical varieties), which also satisfy mild conditions on the moment map and generic stabilizers.
 
 We establish a rigid structure theorem for hyperspherical varieties, Theorem~\ref{thm:structure}
-- they are all given by Whittaker-twisted symplectic induction
from symplectic representations of reductive subgroups $H \leqslant G$.\footnote{Note that, although
all hyperspherical varieties arise thus, the converse is not the case - it is, in general, not obvious
when a Whittaker-twisted symplectic induction will be hyperspherical.}

We then introduce the notion of ``distinguished polarization''
for a hyperspherical variety, and show  that it is essentially unique when it exists.
Here the notion of distinguished polarization is slightly
weaker than requiring that $M$ be equivariantly a cotangent bundle,
and in particular it includes the Whittaker case.

In \S\ref{dualofX} we attach a dual Hamiltonian $\Gv$-space to a polarized hyperspherical variety. This dual is built explicitly as a Whittaker-twisted symplectic induction  using three main ingredients:
\begin{itemize}
\item the dual subgroup $\Gv_X\subset \Gv$ associated to a spherical $G$-variety~\cite{GaitsgoryNadler,SV,KnSch};
\item the commuting $SL_2\to \Gv$ of~\cite{SV,KnSch}, and
\item a symplectic representation of $\Gv_X$, denoted $S_X$. 
\end{itemize}

Finally, \S \ref{goodhypersphericalpairs} discusses 
our expectations concerning the exact domain of the duality $M \leftrightarrow \check{M}$
(see in particular Expectation \ref{anomaly expectation})
and also tentatively discusses 
  issues of rationality, i.e., what the ``split'' form of a hyperspherical
variety is when the ground field is not algebraically closed. 

In summary,  this part of the paper produces an explicit and readily computable class of pairs $(M, \check{M})$
that are candidate pairs for relative Langlands duality.
They are all ``Whittaker induced'' from symplectic representations -- which allows us to reduce many questions
to the case where $M$ or $\Mv$ is in fact  a symplectic representation.

\medskip
 
$\bullet$ {\bf Part 2: Local   Theory. (\S \ref{shearingsec0} - \S \ref{section:numerical Plancherel})} \label{part2intropage}
In these sections we formulate and study 
our general story in the local unramified setting. 
That is to say,  given  $(M, \check{M})$ as above, we will spell out the matching between
``automorphic'' and ``spectral'' data corresponding to $M$ and $\check{M}$ 
under certain restrictions -- in particular, the existence of a polarization on $M$ --
 starting from a finite characteristic local field such as $F=\F_q((t))$,
or its ``geometric'' analogue $F=\C((t))$. 
As we have explained in \S \ref{AFT2}, this matching takes the form of an equivalence of categories in the geometric case,
and an isomorphism of Hilbert spaces in the $\F_q((t))$ case.

\begin{itemize}
\item 
In \S \ref{shearingsec0} we discuss the notion of {\em shearing}, or shifting cohomologically graded vector spaces according to their weights. This basic operation is familiar from Koszul duality (or from considerations of Tate twists), as reviewed in this section, but thanks to the machinery of higher category theory~\cite{1affine} it can be carried out in very general settings. The widespread systematic use of shearing is an underlying current in our work, especially in describing the spectral sides of our conjectures;
in particular it will be used in \S \ref{section-unramified-local}.

\item 
In \S \ref{section-unramified-local}  we attach a category of ``spherical sheaves on a spherical variety'' - the automorphic quantization of $M=T^*X$ over $F$ a geometric local field. 
We formulate a spherical counterpart to the geometric Satake correspondence, identifying this category with the spectral quantization (a category of coherent sheaves) associated to the dual hyperspherical variety $\Mv$. 
This conjecture encounters many technical complications related to sheaf theory on infinite dimensional spaces; 
these are   important issues for further study. 
 We also spell out various constraints this equivalence is expected to satisfy.
 
 \item 
In \S \ref{PlancherelCoulomb} we explicate 
a part of the local conjecture from \S \ref{section-unramified-local} that
avoids most of the intricacies of sheaf theory from the previous section. 
By studying the internal endomorphisms of a certain basic object -- the so-called {\em Plancherel algebra}  -- we deduce a version of the conjecture
that is simultaneously 
  related to the study of the Coulomb branch by Braverman, Finkelberg and Nakajima \cite{BFN},
  and, as we see in the next section, very closely related to the study of the Plancherel measure. 

\item 
In \S \ref{section:numerical Plancherel} we establish the conjecture of \S \ref{PlancherelCoulomb} at the level of Frobenius traces, in many cases, by relating it to the known description of the Plancherel measure for spherical functions (as computed in \cite{SaSph,SaWang}). This also accounts for our name ``Plancherel algebra,''
which we have used in preference to ``Coulomb branch'' for our purposes because it is more evocative of its arithmetic role. 
\end{itemize}
\medskip

$\bullet$ {\bf Part 3: Global Theory (\S  \ref{section:global-geometric} - \S \ref{L2numerical})}
\label{part3intropage}
We study the global story, by which we mean
the story over a curve over either a finite field or the complex numbers,
with function field $F$. 
Just as in Part 2, we will start with $(M, \check{M})$ as in Part 1,
and use $M$   to construct additional data
 on the automorphic side of the global Langlands program  (as in  \S \ref{whatisRelAut})
 and we use $\check{M}$ to construct corresponding data on the spectral side.  
 The central proposal is that they match;
 for example, in the ``arithmetic'' case, this amounts 
  to the diagram \eqref{match diagram}.  
 
\begin{itemize}

 \item \S \ref{section:global-geometric}
 constructs an automorphic quantization (here, a
 automorphic sheaf or an automorphic function, depending on context) when
 $M$ is polarized and $F$ is a global field or a geometric global field. 
 These are ``periods'' or ``period sheaves.'' 
  
 \item \S \ref{Lsheaf} constructs a spectral quantization (here, an object of the spectral category -- a sheaf on the moduli space of local systems) when $\check{M}$ is polarized
 and $F$ is a geometric global field.
 These are termed $L$-sheaves and we will explain
 how their stalks recover $L$-functions. 
 
 \item \S \ref{GGC} formulates the geometric matching statement:  ``period sheaf matches $L$-sheaf''
 in the settings specified in the prior two sections.  Interesting subtleties arise here, for the sheaves
 do not always live in the most natural categories for the Langlands correspondence, and have to be forced into them. We 
then study various formal properties of the conjecture as well as some examples,  and
formulate a conjecture in some unpolarized cases.
  
\item \S \ref{P1} performs a useful reality-check by studying the matching of period and $L$-sheaves in the case of the projective line $\mathbb{P}^1$.

\item \S \ref{L2numerical}
passes to the arithmetic  setting of everywhere unramified forms over a global field of finite characteristic. 
We formulate the matching of automorphic and spectral quantizations in this context --
this amounts to \eqref{match diagram}, that is, to the matching of a period and an $L$-value. 
This  is part of the ``classical'' setting for relative Langlands duality and so can
be compared with numerical predictions. 
Our analysis gives evidence for the Conjecture of  \S \ref{GGC}.
 We study several phenomena of independent interest in passing,
in particular   the star-period.

\end{itemize}

\medskip
$\bullet$ {\bf Part 4: Local-to-Global and Factorization Aspects (\S \ref{local-global} - \S \ref{spectral geometric quantization})} \label{part4intropage}

In Part 4, we  work in a geometric setting and relate the local theory from Part 2 with the global theory from Part 3. This gives rise to rich algebraic structures of interest in their own right.
A key role here is played by {\em factorization}.
 
\begin{itemize}

\item \S \ref{local-global} introduces the (one-point, unramified) $\Theta$-series functors on the automorphic and spectral sides, and formulates the local-global compatibility conjecture: the $\Theta$-series intertwine the local and global period conjectures of \S \ref{section-unramified-local} and \S \ref{GGC}.

\item \S \ref{automorphic-factorization} describes the role of factorization structures on the automorphic side. In particular we explain how the Plancherel algebra (categorifying the Plancherel measure for spherical functions) extends to a factorization associative algebra on any curve,
 closely related to the relative trace formula, and also to a factorizable form of the $\Theta$-series construction.

\item \S \ref{spectral-factorization} describes factorization structures on the spectral side with much greater precision but at the cost of restricting to the Betti setting. We formulate the notion of a ``spectral deformation quantization,'' producing the spectral counterpart to the factorizable Plancherel algebra, and construct one when $\Mv$ is a (possibly twisted) cotangent bundle.

\item \S \ref{spectral geometric quantization} studies the {\em algebra of $L$-observables} and its action on the $L$-sheaf, a geometric counterpart of the $L$-function of $\Mv$ and its square root given by the $L$-function of $\Xv$, via coherent sheaf forms of microlocalization and quantization. 
We also explain a spectral construction of 
``geometric Arthur parameters'' using the tools of shearing and $L$-sheaves.
\end{itemize}

We refer here to page
\pageref{part4longintropage} and the section introductions for lengthier and more motivated discussions of the contents of this part. 
  
 \medskip
 
$\bullet$ {\bf Part 5: Appendices}   
In Appendix \S \ref{shearingsec} we collect basic properties of Koszul duality that appear repeatedly in the text. We also introduce a simple gadget, the ``spectral exponential sheaf'', which underlies our approach to Whittaker-type constructions on the spectral side.

\S \ref{sheaf theory} we survey the menagerie of sheaf theories that we use throughout. 

\S \ref{geometric Langlands} gathers background regarding the geometric Langlands correspondence in its various forms.

 \S \ref{fieldtheorysec} contains a discussion of structures coming from topological field theory. After reviewing factorization and $E_n$-algebras we sketch a formal approach to the definition of algebraic quantum field theory on curves, capturing part of the rich structure predicted by the metaphor of arithmetic quantum field theory. We explain how this formalism conjecturally houses many features of the Langlands program and its relative version. 
 
 \S \ref{Eisappendix}  is a garbage can full of miscellaneous computations, which we could neither bring ourselves
 to discard entirely, nor to leave in the main part of the file.

\subsection{Some examples} \label{IntroEx}

 The following tables give a list of sample examples of hyperspherical dual pairs $(M, \check M)$.
   Many of these examples are discussed in some more detail in the paper, and the notation explained in more detail. 
This list is very far from comprehensive, and are chosen to some extent to reflect examples discussed in the paper. 
   Not all examples are, to our knowledge, studied in the automorphic literature;
   but some cases have been so extensively studied that  is impossible to even begin to summarize the work. 
We have at least  tried to give representative citations to   papers discussing the associated global period,
 from which precise details  can be extracted.

{ \begin{table} 
\centering
\caption{ Examples of hyperspherical dual pairs}
\label{tab:examples}
\tiny
\begin{tabular}{|c|c|c|c|}\hline
attribution/name &  $G/H$ or $(G,X)$  &  $\check G/\check H$ or $(\check G, \check X)$ 
 &  attribution/name  \\ \hline \hline
{\tiny Hecke} & $\PGL_2/A$ & $(\SL_2, \mathbb{A}^2)$ & {\tiny normalized}  \\
 		&		&						& {\tiny Eisenstein } \\
		\hline
{\tiny Iwasawa--Tate} & $(\mathbb{G}_m, \mathbb{A}^1)$ & $(\mathbb{G}_m, \mathbb{A}^1)$ & {\tiny Iwasawa--Tate} \\ \hline
\cite{LM} &  $H \times H/\Delta H$ & $H \times H/\Delta'(H)$  &  \cite{LM} \\  \hline
 \cite{JacquetRallis} &  $\GL_{2n+1}/\GL_n \times \GL_{n+1} $ & $\GL_{2n+1}/\mathrm{Sp}_{2n}$  & \\
  \hline
  \cite{JacquetRallis}  &  $\GL_{2n}/\GL_n \times \GL_n$ & $\GL_{2n} \times_{\mathrm{Sp}_{2n}} \mathrm{std}$  & \\
  \hline
  \cite{GJ}  &  $\GL_{n} \times \GL_n, M_n$ &$\check{G} \times_{\GL_n} \mathbb{A}^n$ &  \cite{JPSRS} \\
  \hline
\cite{WanThesis}   & $\mathrm{SO}_{4n+1}/\GL_{2n}$ & $\check{G} \times_{\mathrm{Sp}_{2n}^2} \mathrm{std}^2$  & \\

\hline
\hline
   point & $G/G$ &   $G/(U, \psi)$  &  Whittaker, \cite{LM} \\  \hline
\cite{JacquetShalika} &  $ \GL_{2n}/ \GL_n. (M_n, \psi)$   & $\GL_{2n}/\mathrm{Sp}_{2n}$  &  \cite{JRsymplectic} \\
 \hline
 ``std. $L$-function'' & $(\GL_n, \mathrm{Std})$. & 		$(\GL_n,  \GL_n/\Gm\cdot (U, \psi))$ &  \cite{JS1, JS2}  \\
 \hline
 {\tiny ``sph. harmonic''} &  $\mathrm{SO}_{2n}/\mathrm{SO}_{2n-1}$ & $\SO_{2n}/  \SO_3.(U, \psi)$  &  Bessel \\
    \hline
   & $\mathrm{SO}_{4n}/\GL_{2n}$ & $\mathrm{SO}_{4n} \times_{\mathrm{Sp}_{2n}(U. \psi)} \mathrm{std}^2$  & \\ 
      \hline
   \cite{BumpGinzburgSpin}    & $\GSp_6\times_{\GL_2.(U, \psi)} \mathbb{A}^2$& $(\mathrm{GSpin}_7, \mathrm{spin} \times \mathbb{G}_m)$ & \\ 
\hline
\cite{GinzburgE6} & $GE_6 \times_{\GL_3 .(U, \psi)} \mathbb{A}^3$ & $(GE_6,  \mathrm{std}_{27} \times \mathbb{G}_m) $ & \\ 
     \hline\hline
 attribution/name &  $G/H$ or $(G,X)$  &  $\check{M}$, {\bf not} $\check{X}$
 &  attribution/name  \\ \hline \hline
    \cite{GGP}&  $\mathrm{SO}_{2n} \times \mathrm{SO}_{2n+1}/\mathrm{SO}_{2n}$ & $\mathrm{SO}_{2n} \times \mathrm{Sp}_{2n}, \mathrm{std} \otimes \mathrm{std}$  & \cite{RallisIP} \\ 
    Gan--Gross--Prasad  & & &   Rallis' inner product  \\ \hline
\cite{GinzburgRallis}& $(\PGL_6, (\GL_2)^3 \cdot (V, \psi)/\Gm)$ & $(\SL_6, \check{M} = \wedge^3)$ & \\
  \hline 
\cite{WanZhang} &   $\GSp_6 \times \GSp_4/ (\GSp_4 \times \GSp_2)^0$ & $(\check{G}, \check{M} = \mathrm{spin}_7 \otimes \mathrm{spin}_5)$ & \\  \hline
\cite{WanZhang}  &   $E_7/ \PGL_2 \cdot (U, \psi)$ & $(\check{G}, \check{M} = \omega_7)$ & \\  \hline
\cite{WanZhang}  &   $\mathrm{GSO}_{12}/ \GL_2 \cdot (U, \psi)$ & $(\check{G}, \mathrm{half-spin} \times T^*\mathbb{G}_m)$ & \\  \hline

  \hline

\end{tabular}
\end{table}
}
 We have separated them into three general classes. 
The first are the simplest class, where both $M$ and $\check{M}$ are polarized by $G$-varieties.
The second are also polarized but now allowing ``twisted'' polarizations, as indicated by the presence of $\Psi$. 
Finally, in the third class, one side does not admit a polarization; in the examples we present,
these cases have vectorial $\check{M}$ with $\GGm$ action by scaling, which corresponds to an automorphic period
which squares to the central value of an $L$-function. 
We also warn that, although we often listed $X$ and not $M$, it is really the  latter and not the former that is intrinsic, cf. Examples \ref{ex:SX3}, \ref{ex:SX3p}.

 The generality we propose above is certainly far from  the end of the story; we have many interesting examples of
 duality-type phenomena when $M$ is either:
 \begin{itemize}
 \item a graded Hamiltonian space that is hyperspherical {\em except for} the connectedness condition \eqref{condtorsion} in the definition therein,  or
 \item not affine, or 
 \item not smooth, or 
 \item not even a variety (e.g., a stack or derived scheme),
 \end{itemize}
and indeed some of the most interesting consequences for number theory may reside in such instances.  

A simple but already important example is the case of $M = T^*(U \backslash G)$, 
  considered as a space under $T \times G$; its putative dual should
  be $\check{M} = T^*(\check{U}^- \backslash \check{G})$ considered as a space under $\check{T} \times \check{G}$.
This example is related to the theory of Eisenstein series, but 
does not fall in our general framework because $M$ is not affine.
See in particular Example \ref{Eisgeom}. 
\footnote{In this connection it is interesting to observe
that the Whittaker model, which is related to a bundle over $U \backslash G$, is very well behaved; and
this is precisely because the associated $M$, which is a {\em twisted} cotangent bundle of $T^*(U \backslash G)$,
in fact {\em is} affine.}

 \subsection{Some open questions}
Put charitably, this paper   leaves open many more questions than it answers.
  Some of these
are formulated as conjectures throughout the paper; but we also take the opportunity to draw
attention to some that are perhaps less clearly formulated but seem important. 

\begin{itemize}
\item 
Can we  extend the theory of spherical
varieties to a theory of hyperspherical Hamiltonian spaces?
Many of the Knop's techniques for studying spherical $G$-varieties
$X$ already go through the cotangent bundle $T^*X$;
many of his results have been generalized to Hamiltonian spaces
by Losev \cite{Losev}. The ideas of Knop and Losev are used heavily in \S \ref{sphvar}. 
It would be particularly desirable if the duality $(G, M) \leftrightarrow (\check{G}, \check{M})$
could be read off from matching combinatorial invariants.   
 
\item From the point of view of supersymmetric quantum field theory, one should consider not just symplectic varieties with Hamiltonian $G$-action, but hyperk{\"a}hler manifolds with isometric actions of the compact real form of $G$. So it is natural to ask if complex hyperspherical varieties admit such metrics. This is known in the abelian case (the theory of hypertoric varieties) as well as for cotangent bundles of flag varieties 
and symmetric spaces.

\item What is the precise role of (hyper-)sphericity?  (cf. Remark
\ref{L-observables}). In the local setting, our conjectures (for $M=T^*X$) predict strong rigidity for the categories of spherical sheaves
$\Hecke^X$ of spherical sheaves on $LX$. In particular, we predict they are generated by the basic sheaf, admit a ``graded'' lift in the spirit of Koszul duality~\cite{BGS}, and admit a locally constant factorization structure (an algebraic analog of a braided monoidal structure in topology), see  Problem \ref{localconstancyHeckeX}.  These properties appear to be strongly tied to sphericity, and it would be very useful to have a formal statement to this effect.

\item   The singularities present in  local loop spaces pose
a number of interesting questions that are important for a better understanding of the local conjecture. 
This is discussed at varoius points in \S \ref{section-unramified-local}, see 
for example Remark \ref{! vs * conjecture}  and 
\S \ref{Sheafdesiderata} 
for some specific problems. 

\item The papers \cite{SaSph}, \cite{SaWang} of the second-named author, the second in collaboration with
J.\ Wang,   provide
a computation of spherical Plancherel measure for a large class of spherical varieties. 
This is the numerical form of the local conjecture (cf. \S \ref{section:numerical Plancherel}); so 
is it possible to categorify the proof of \cite{SaWang} to
give a parallel proof of the local conjecture ?

\item What is the dual version of the theory of asymptotics  and boundary degenerations of spherical varieties, as in~\cite{SV}?

 \item 
 A crucial question for applications is to extend the applicability of the theory past the split case
 (and indeed to study more carefully the properties of the split form of a Hamiltonian space)
 cf.  \S   \ref{GMdesiderata}.
 
 \item  Extend the global story to number fields, and to ramified situations. For example, to do so, one should -- starting from $(G, M)$ -- specify
for each local field $F$, a  representation $\Omega_M$ of $G(F)$; and for each global field $F$ with ring of adeles $\mathbb A$, a specified morphism  $\Theta$
from the restricted tensor product of the $\Omega_{M,v}$  to the space of automorphic functions, 
 and (the hard part) give  compatible descriptions of both data on the spectral side in terms of $\check{M}$.
 
 \item  The star-period, discussed in \S \ref{starperiod}, is obviously also an important numerical object,
and deserves further examination.

 \item  Clarify the global quantization of non-polarized $M$:
 \begin{itemize}
 \item On the automorphic side, the theory of the Weil representation does this, 
 but one needs to systematically understand the splitting in order to have a fully satisfactory theory
 (see \S~\ref{remark:unpol} for a suggested answer). 
 \item On the spectral side, one needs to construct a spectral analogue of the Weil representation,
 which would be sufficient to spectrally quantize general pairs $(\check{G}, \check{M})$.
 We have carried this out in the Betti case and will present it separately. 
  \end{itemize}

\item Develop spectral quantization (\S\ref{spectral-factorization}, \S\ref{spectral geometric quantization}) in the positive characteristic setting. What is the counterpart of $E_2$-algebras and braided tensor categories as the algebraic structure underlying locally constant factorization algebras? What takes the place of shifted differential operators and relative flat connections in the description of the global Hecke category and $L$-observables?

\item  An extremely intriguing direction is to find new ``exotic'' examples
of dual pairs $(G, M)$ and $(\check{G}, \check{M})$ that are outside the framework of this paper. Already many
Rankin-Selberg integrals seem to require that we allow $M$ or $\check{M}$ to be singular, or stacks, to properly
fit into this framework. In some cases these Rankin-Selberg integrals are 
still related to spherical varieties, but which fail some of our assumptions -- for example they are singular \cite{SaRS}, 
or they possess roots of type $N$.  See also \cite{ChenVenkatesh}
 for some preliminary work in this direction.

\item The field theory perspective suggests a higher categorical structure to the collection of periods for $G$. In particular it is natural to consider morphisms of periods, given by (quantizations of) Lagrangian correspondences of Hamiltonian $G$-spaces. Examples of this appear in the theory of unfolding, and it would be interesting to place these examples in a more structured framework, in particular to study duality on morphisms of periods.  
 
\item Periods for product groups $G\times H$ are closely related to instances of Langlands functoriality. The field theory setup suggests that it may be interesting to consider an expanded notion of functoriality: to study linear maps from 
the vector space of automorphic functions for $H$ to the vector space of automorphic functions for $G$ which arise
 from Hamiltonian actions of $\Gv\times \Hv$ rather than only homomorphisms $\Hv\to \Gv$.

\end{itemize}

 \subsection{Recent developments} 
 
 There have been substantial pieces of work around the topics of interest in this paper,
 during the time that the paper was being prepared.
 We try to summarize a few of the relevant works known to us here. Note that our description of these papers are often incomplete;
 we note only the results that they contain that are most directly related to this work. 
 
 \begin{itemize}
 
 \item {\it For recent results related to geometric local duality, see \S \ref{OtherExamples}.  }
 
 \end{itemize}

 In the setting
 of local fields:
 \begin{itemize}
 	\item 
 	Wee Teck Gan and Bryan Wang Peng Jun have provided evidence \cite{GanJun1}, arising from $\theta$-correspondence,
 	for the role of hyperspherical duality in {\em ramified} local Langlands.
 \end{itemize}

 In the global setting:
 \begin{itemize}
 	\item Eric Chen has studied \cite{Chen2, ChenVenkatesh}, partly in collaboration with the third author,
 	several examples of numerical global duality involving {\em singular} spaces.
 	\item 
 	Tony Feng and Jonathan Wang have examined \cite{FengWang}
 	the geometric conjecture in the Hecke case, and proven it up to an issue of identifying two extensions. 
 	
 	\item Zhengyu Mao, Chen Wan, Lei Zhang  
 	\cite{MWZ} have proposed a relative trace formula comparison to reduce a version of the global numerical conjecture (even allowing ramification!) to the ``strongly tempered'' case; this is based on related conjectures ({\em not} following from hyperspherical duality) about degenerate Whittaker models. They have also \cite{MWZ2} provided a classification of a large class of strongly tempered hyperspherical varieties.
 	
 	\item Gan and Jun have also checked in \cite{GanJun2},
 	in the same set of examples as \cite{GanJun1},  the validity of the
 	results of our \S \ref{section:numerical Plancherel} as well as the global conjectures of \S 14 (allowing ramification, although
 	they do not examine the constant).

 \end{itemize} 
 
 On the general geometry of hyperspherical varieties:
 \begin{itemize}
 	
 	\item 
 	Finkelberg, Ginzburg, Travkin  have proposed \cite{FGT} the statement
 	that, given dual hyperspherical varieties $M, M^{\vee}$, there is a close
 	relation between ``symplectic Borel orbits on $M$'' (by which we really mean Borel orbits on $X$ if $M=T^* X$,
 	which can be formulated purely in terms of $M$), and symplectic Borel orbits on $M^{\vee}$.
 	
 	\item   In the context of his study of endoscopy for symmetric varieties, Leslie \cite{Leslie} 
 		has both clarified some of the constructions
 		of our Part 1, in particular explicating the symplectic representation that enters the definition of the dual hyperspherical variety and confirming our Conjectures \ref{conjsymplectic}, \ref{conjGaloisaction} in that case. 
 		Much of Leslie's work applies to general symmetric varieties, not only those with hyperspherical cotangent. 
 	
 	\item  Work of Jiajun Ma, Congling Qiu, Zhiwei Yun, Jialiang Zou in progress provides a categorification
 	of the work of \cite{FGT} and proof of this statement in some cases related to $\theta$ correspondence. 
 	
 \end{itemize}

 \subsection{Acknowledgments}
 This project was initiated at the Institute of Advanced Study during the special year program of 2017--2018 
 on locally symmetric spaces. We thank the IAS for its hospitality and for providing the right environment
 in which to launch our collaboration. During the course of more than 5 years that this work took, D.B.-Z.\ was supported by a fellowship and a Visiting Professorship at MSRI from the Simons Foundation and NSF individual grants DMS-1705110 and DMS-2001398 as well as NSF Grants No. 1440140 and PHY-1066293 while in residence at the MSRI  and the Aspen Center for Physics (respectively), Y.S.\ was supported by NSF grants DMS-1801429, DMS-1939672, and DMS-2101700, and Simons Fellowship 616018, and A.V.\
 was supported by NSF grants DMS-1401622 and DMS-1931087. 
   
We would like to acknowledge Tony Feng and Jonathan Wang for many helpful conversations. In particular, Feng and Wang studied independently \cite{FengWang} the Hecke case,
 and their ideas there were very helpful for orientation.
Joint computations 
 of one of the authors (A.V.) with Feng and Wang helped to support
 some of the finer features of the global conjecture and influenced our thinking.
We are indebted to Sam Raskin for numerous consultations and providing ``urgent care'', in particular on sheaf theory (Appendix~\ref{sheaf theory}), factorization (\S \ref{automorphic-factorization}), and restricted geometric Langlands. 
We would also like to acknowledge Sam Gunningham for helpful discussions surrounding shearing and geometric Satake,  John Francis for suggesting the parallel between $L$-functions and factorization homology, David Nadler for discussions of the Langlands program as a TQFT,  Pavel Safronov for teaching us his point of view on quantization for observables and states, Justin Hilburn, Tudor Dimofte and Andy Neitzke for discussions of S-duality for boundary conditions, and Simon Riche for discussions about modular Satake correspondence. 
Finally, we would like to thank Alexander Braverman, Dario Beraldo, Eric Chen,  Minhyong Kim,  David Hansen, Tony Feng, Sam Raskin, Pavel Safronov and Harold Williams for 
  feedback on,  and corrections to,  an earlier version of this paper.

  \section{Notation and conventions}
\label{subsection-notation} 

We are going to summarize here some of the notation and conventions used throughout the paper. Since we will endeavor to also define notation where it is used, the reader should refer to this section
only as needed.

\subsection{The coefficient fields $\FF$ and $\kk$} \label{coefficients}

There will be two base fields used in the paper.  

The field $\FF$ will usually be the base field on the automorphic side. It will (most of the time) either be  the algebraic closure of a finite field, or the field of complex numbers, 
according to context.

By $k$ we will denote a (usually algebraically closed) coefficient field of characteristic zero, which will be the base field on the spectral side.
When $\mathbb F = \overline{\FF_q}$ , this coefficient field will be taken to be $\overline{\mathbb Q_l}$, for some $l$ different from the characteristic of $\FF$.
Moreover, in this case, we will fix once and for all a square root $\sqrt{q}$ of the order of $\mathbb{F}_q$ inside $k$,
which will permit us to think about half-Tate twists. In particular, we assume the existence of such a square root, when half-integral Tate twists appear. Our conjectures,
however, will not depend on this choice, and we   discuss at various points 
how to formulate statements
in an invariant fashion, see in particular \S \ref{extended-group appendix}.  

We also fix, once and for all, an additive character \index{$\psi=$ additive character of $\FF_q$} $\psi: \FF_q\to k^\times$. In \S~\ref{Affcase} we will, correspondingly, fix an Artin--Schreier sheaf on which all ``Whittaker-type'' constructions will depend.

\subsection{Curves and their fundamental/Galois groups} \label{curvenotn} \index{$\Sigma=$ curve} \index{$\mathbb F$}
We fix a smooth projective curve $\Sigma$ over the field $\mathbb F$.  The symbol $F$ will we will be used to denote a variety of fields, depending on the setting,
loosely related to functions on $\Sigma$. This includes: 
\begin{itemize}
 \item  the function field of $\Sigma$;
 \item the completion of the above at a point of the curve;
 \item when $\Sigma$ is equipped with a model over $\FF_q$, the function field of $\Sigma_{\FF_q}$;
 \item occasionally, we will also use $F$ to denote number fields and their completions.
\end{itemize}

\subsubsection{Global setting} \label{curvenotnglobal}
When $F$ is a global field,  i.e.\ the function field of a curve over $\FF_q$ or a number field, we denote by $\mathbb A$ its ring of adeles, and for a linear algebraic group $G$ over $F$ we set $[G]=G(F)\backslash G(\mathbb A)$. \index{$[G]$}

There is an adelic norm
\begin{equation} \label{idelenorm}  \mathbb{A}^{\times}/F^{\times} \rightarrow k^{\times}, x \mapsto |x|\end{equation}
 which sends each uniformizer in the adeles to the inverse of the associated residue field cardinality.

 Let 
$\Gamma=\Gamma_F$ be the Weil group of $F$ when $F$ is a function field of a curve over $\FF_q$.
This comes with the cyclotomic character: \index{cyclotomic character} 
\begin{equation} \label{cyclotomicdef} \varpi: \Gamma_F \longrightarrow \Q^{\times}\end{equation}
 which sends a(n arithmetic) Frobenius element to $q \in \Q^{\times}$. 
 The chosen square root of $q$  in $k$ defines in particular a square root $\varpi^{1/2}: \Gamma_F \rightarrow k^{\times}$.

  \subsubsection{Local setting} \label{localfieldnotn}
    
   Let $F = \FF((t))$, with integer ring $O = \FF[[t]]$. 
   
   In the case when $\FF = \overline{\FF_q}$ and we need to distinguish, we will use the notation
   $\ff, \fo$ for the same objects defined over the finite field $\FF_q$:
   $$ \ff = \FF_q((t))\ \mbox{ and } \fo = \FF_q[[t]].$$
   
   We will also use the following notation for arc and loop spaces:
   For $X$ a scheme,  we shall write $X_O$ for the formal arc space, representing the functor $$R\mapsto \mbox{$R[[t]]$-points of $X$}, $$
    a scheme over $\FF$; 
     and $X_F$ for the formal loop space, representing 
     $$R\mapsto \mbox{$R((t))$-points of $X$}, $$
     which for $X$ affine is represented by an ind-scheme over $\FF$, see \S \ref{arc loop space}. 

       \subsection{Reductive group notation} \label{redgroupnotn}

   \subsubsection{}  \label{reductive group notation}
   Let $G$ be a split reductive group over $\FF$. (Our reductive groups will be understood to be split by default, unless otherwise stated). 
  We will usually use  
  $$ U \subset B \subset G$$
  to denote a maximal unipotent subgroup and  Borel subgroup of the reductive group $G$. 
    A {\em pinning} of $G$ is as usual a choice of $T \subset B$
  and an isomorphism $\mathbb{G}_a \simeq U_{\alpha}$ of each root space. 
  In this situation there is a distinguished character $B \rightarrow \Gm$  
  (descending to the torus quotient), given by the sum of all positive roots.
  We will denote it by $e^{2 \rho}$
  or by $2 \rho$ at the Lie algebra level. 

  We denote by $\check{G}$ the Langlands dual group to $G$, which we take as a split pinned
  group over $\kk$.

 The {\em exponents} $G$ of a reductive group are, by definition,  \index{exponents}
 the dimension of homogeneous polynomials generating the polynomial ring of $G$-invariant regular functions on $\mathfrak g$. For example, the exponents of $\SL_n$ are $2,3,\dots, n$. 

  \subsubsection{The duality involution} \label{dualityinvolution} \index{duality involution}\index{$d=$ duality involution}
   There are two closely related involutions defined for a pinned reductive group, which we review now. 
The first version is the {\em Chevalley involution}, and we will denote it by $c$.  It is uniquely characterized
by the fact that it preserves the pinning and acts on  the torus according to   $t \mapsto w(t^{-1})$.

The second is what D.\ Prasad has called the {\em duality involution} \cite{Prasad-contragredient}, and we will denote it by $d$. It is the composite of the Chevalley involution 
with conjugation by $e^{\rho}(-1)$, i.e., it is uniquely 
that it {\em negates} the pinning and acts on  the torus according to   $t \mapsto w(t^{-1})$.

For example, for $\SL_n$, the duality involution is given by
$$g \mapsto \mathrm{Ad}(w) (g^t)^{-1},$$
where $w$ is the matrix with all entries $1$ on the anti-diagonal, and other entries zero.

We will use a superscript $d$ for various involutions induced by the duality involution. For example,  \index{$X^d$}
  for $X$ a $G$-variety we denote by $X^d$ the same variety
but with $G$-action twisted by  means of $d$.

\subsubsection{The notation $H$ and $\check{G}_X$; the Arthur $\SL_2$s}\label{dual to H}

Let $G, \check{G}$ be a pair of dual reductive groups as above. 
We will often deal with a Hamiltonian $G$-space $M$ and Hamiltonian $\check{G}$-space $\check{M}$
that will be, in a suitable sense, in duality with each other. Because of Theorem \ref{thm:structure}, 
we will be particularly interested  in a cases where
both sides are determined by linear-algebraic data:  an $\SL_2$
in the group, a commuting reductive subgroup, and a symplectic representation of that subgroup. 
The notation we will use for this data will be, however, somewhat asymmetric,
reflecting notation used in the automorphic literature, where
the symmetry of the two sides is not apparent:

\begin{itemize}
\item[(a)] On  the $G$ side,   we will use the notation
\begin{equation} \label{basic0.1} (H \subset G, \mathfrak{sl}_2 \subset \mathfrak{g}, S \textrm{ a symplectic $H$-representation)}\end{equation}
for the data defining $M$.

\item[(b)] 
On the spectral side, in the case when $\check{M}$ is in duality with $M=T^*X$, we will rather use the notation
\begin{equation} \label{basic0.2} (\check{G}_X \subset \check{G},  
\mathfrak{sl}_2 \subset \check{\mathfrak{g}}, \mbox{$S_X$ a symplectic $\check{G}_X$-representation})\end{equation}
for the corresponding data defining $\check{M}$.
\end{itemize}

\begin{remark} \label{WhittArth0} (Whittaker and Arthur role of $\SL_2$):
This linear-algebraic data for $M$ and $\check{M}$
play an important role in  shaping the automorphic and spectral quantizations
in the sense of \S \ref{whatis}. The role of $\sl_2$ is quite different on the two sides:
informally, automorphically it measures the involvement of ``Whittaker characters,''
whereas spectrally it relates to the ``Arthur'' $\SL_2$ that quantifies the failure of temperedness.

\end{remark}

\subsection{Navigating the assumptions on hyperspherical spaces} \label{navigating} 

In Section \ref{sphvar} we will introduce the central objects in the relative Langlands duality that we propose in this paper, the \emph{hyperspherical} Hamiltonian spaces. Since the definitions are given at first over algebraically closed fields in characteristic zero, while these spaces are used later over non-algebraically closed fields in arbitrary characteristic, we would like to point the reader to the places where this leap is explained. 

Central in our use of hyperspherical spaces over arbitrary fields will be the Structure Theorem \ref{thm:structure}, which states that hyperspherical $G$-spaces over algebraically fields in characteristic zero are obtained by a process of ``Whittaker Hamiltonian induction'' from symplectic representations of reductive subgroups of $G$. It is such a structure that we \emph{assume} over different fields, in order to talk about hyperspherical spaces there, as outlined in the preliminary discussion of rationality issues in \S~\ref{hdprings}. 

Moreover, at many points in this paper we will need our spaces to admit a \emph{distinguished polarization}, in the sense introduced in \eqref{dp}, over algebraically closed fields in characteristic zero, and in \S~\ref{hdprings} over more general rings. Such spaces are closely related to spherical varieties, with additional assumptions (such as smoothness) inherited from the conditions defining ``hyperspherical,'' see Proposition \ref{isspherical}.

In Section \ref{dualofX}, which deals mostly with such polarized hyperspherical spaces, proposing a construction for their Hamiltonian dual, the effect of working over a non-algebraically closed field is discussed in \S~\ref{CheckMRat}. In particular, when the field of definition of is not algebraically closed, the dual side should come with additional structure, which extends the definition of the $L$-group of $G$; we define the $L$-group of a spherical variety (Definition \ref{definition-LGX}), and discuss, but do not quite define, the action of the Galois group on the dual hyperspherical variety (see Conjecture \ref{conjGaloisaction}). 

Finally, in \S~\ref{GMdesiderata} we postulate that there should be a distinguished ``split'' form of a hyperspherical variety over general rings, and provide our working definition for split forms, Definition \ref{dhpFqdef}, which is used throughout the paper. This definition requires some split form of the data that go into the structure theorem to descend in a unique way from $\mathbb Z$ localized at a finite number of places, an assumption that does not always hold. However, it works in large enough characteristics and over large enough finite fields (Proposition \ref{rigidity}), allowing us to talk about a distinguished split form in those settings in the absence of an abstract theory for those. 

For the remainder of the paper, whenever not specified, a hyperspherical (possibly polarized) hyperspherical variety is one constructed as in the structure theorem over $\mathbb C$ (as in \S~\ref{hdprings}), and it is called split whenever it satisfies our working definition of \S~\ref{GMdesiderata}.

  \subsection{Shifting, super-vector spaces and Frobenius traces} \label{anglebracketnotation}

  For $V$ a vector space over $\kk$ the dual of $V$ will be denoted by $V^{\vee}$, pronounced ``vee-vee,''
  and the symmetric algebra on $V$ will be denoted by $\Sym^* V$ or simply $\Sym V$ to keep typography simple. 
We denote by
  $$ \det(V) \textrm{ or } \|V\| = \wedge^{\dim V} V$$
  the top exterior power of $V$. 
  If $V$ is equipped with an action of the Galois group of a finite field $\mathbb{F}_q$, we define
  $$  [V] = \textrm{trace of {\em geometric} Frobenius on } V$$
  where the geometric Frobenius is inverse to $x \mapsto x^{q}$.
  The geometric Frobenius will be denoted by $\Frob$.  The notation $[\dots]$
  will also be used for shifts (see below), but we hope this will not cause too much trouble.
  
  \subsubsection{Sheaf-function correspondence}
  We record the normalization of sheaf-function correspondence to try to clarify signs. 
  
  For $X$ a variety over a finite field $\mathbb{F}_q$
  and a {\'e}tale sheaf $\mathcal{F}$ over $X$ with coefficients in $\kk$, 
  the associated trace function is, by definition, that function  
  $$X(\mathbb{F}_q) \rightarrow \kk$$
whose value at $x$ is
given by $[\mathcal{F}_{\bar{x}}]$,
with   the bracketed vector space the stalk $\mathcal{F}$ at a geometric point $\Spec \mathbb{F}_q \rightarrow  X$ above $x$, 
and $\mathrm{Gal}(\overline{\mathbb{F}_q}/\mathbb{F}_q)$ is acting by pullback (note that this is a left action,
because the association from fields to their spectra is contravariant). 

  \subsubsection{Shifting}
 For a complex of vector spaces $V$ with a Frobenius action,  we denote by
 $$ V(a,b] = V[b] \otimes_{k} k(a),$$
 the combination of a cohomological shift by $b$ and a Tate twist by $a$. 
 The latter means that we twist the Frobenius action 
 by $q^{-a}$; note this is well-defined for $a \in \frac{1}{2} \Z$ because
 we have fixed $\sqrt{q} \in \kk$.

 In particular, we have
\begin{equation} \label{BTT} \left[  V(b, b/2] \right] =  (-1)^b q^{-b/2} [V]. \end{equation}
 The sign $(-1)^b$ that appears here comes from the fact
 that trace on a complex is an alternating construction.

  \subsubsection{Multiplication by $\sqrt{q}$ is super} \label{Sqrtqsuper} \index{Supervector spaces}
  
  Throughout this paper we will have to deal with supervector spaces, and it is important to clarify why this is. 
  
     There are many contexts in number theory over $\mathbb{F}_q$
  where one wants to multiply by $\sqrt{q}$ or $q^{-1/2}$.  {\em This requires
  supervector spaces to properly geometrize!} 
Indeed, equation \eqref{BTT} can be considered as a first attempt to do this categorically, while preserving
  purity (thus, the cohomological shift).  
  However, this has the unfortunate feature of introducing the sign $(-1)^b$. 
  To have a way of multiplying trace by $q^{-b/2}$, without any sign, we want to include
  in our twist a {\em change of parity}. 
  
  More formally,  following \cite{BD}, we work with 
{\tiny supercallifragalisticexpialodocious}-vector spaces over $\kk$, which will be abridged
 for typographical reasons 
  to ``supervector spaces.'' That means we consider $\Z/2$-graded vector spaces, but with the symmetric monoidal 
 structure given by the Koszul rule of signs. As a formal result of the theory of traces, the trace of an endomorphism (such as Frobenius) in this context becomes the supertrace,  and we will use the extra freedom to fix the signs.     
 \index{angle bracket twists $\langle d \rangle$}

Let us now introduce a special notation for the pure twisting  
 with the right sign: 
\begin{equation} \label{ultimate shearing} V \langle b \rangle = \Pi^b V [b, b/2)\end{equation}
  where $\Pi$ is the shift of parity functor.
  We then have
\begin{equation} \label{BTT2} [V\langle b \rangle ] = q^{-b/2} [V]\end{equation} 
and $[ \| V \langle b \rangle \|] = q^{-b \dim(V)/2} [\|V\|]$. 
We will also use the notation $V \langle b \rangle$ in contexts where there is no Frobenius action to mean the super-analog $\Pi^bV[b]$ of the shift.

For obvious practical reasons, we will not refer to the super-grading every time that such shifts are introduced, i.e., we will often write $V \langle b \rangle = V [b, b/2)$, recalling the super-grading when it is relevant.

 \begin{remark}[Super-vector spaces and Galois descent]\label{super descent}
 The role of super-vector spaces in this paper can be nicely encapsulated in the idea~\cite{theosuper} that $\Vect^{super}$ is a $\Z/2$-Galois extension (in fact algebraic closure) of $\Vect$ in the world of tensor categories. Various constructions, most notably shearing (\S \ref{shearing notation}, \S \ref{commutativity}), are naturally defined as operations on super vector spaces which we descend to less evident operations involving ordinary vector spaces. 
 \end{remark}

\subsubsection{Shifting of complexes and super signs} \label{shift}
To avoid any confusion, and to illustrate the role of super-vector spaces, we are going to describe carefully various
signs that occur in shifting of complexes, deriving from the ``super'' symmetric monoidal structure on complexes (the Koszul rule of signs). 
The reader may want to skip this discussion on a first reading.

Let $k[a]$ be the complex consisting of a copy of $k$ in degree $-a$. 
There are ``evident'' identifications $k[a] \otimes k[b] \simeq k[a+b]$, which are not naively compatible with the symmetric monoidal structure, in that
the composite $$k[a+b]\to k[b]\otimes k[a] \to k[a]\otimes k[b] \to k[a+b]$$ is not the identity.  

For a complex of $k$-vector spaces $M$, we denote by $M[a]$ the 
shift, which is convenient to regard as the tensor product
$$ M[a] =  k[a] \otimes M$$
with respect to the symmetric monoidal structure on complexes. 
With this convention the differential is multiplied by $(-1)^a$. 
We may then identify
\begin{equation} \label{Sms2} M[a] \otimes N[b] \simeq (M \otimes N)[a+b]\end{equation}
using the symmetric monoidal structure and the noted identification $k[a] \otimes k[b] \simeq k[a+b]$. 
The identification is associative, i.e., the two resulting identifications
$ M[a] \otimes N[b] \otimes O[c] \rightarrow (M \otimes N \otimes O) [a+b+c]$ coincide. However as before there is a sign that enters into the commutativity.
Similarly, we may identify
$$ \Hom(M[a], N[b]) \simeq \Hom(M, N)[b-a],$$
in a fashion that is compatible with composition.  Here one must take care with
sign; these can be handled abstractly using the fact that complexes form a closed symmetric monoidal category, and are therefore
enriched over themselves, and then considering the convolution action of the invertible objects $k[a]$.

On the other hand, if we consider instead the {\em super}lines $\{k\langle a \rangle\}_{a\in \Z}$ (with $k\langle a \rangle$ in super-parity $(-1)^a$ and in cohomological degree $-a$) then these sign issues disappear, 
and the composite of natural maps $$k\langle a+b\rangle \to k\langle b\rangle \otimes k\langle a \rangle \to k\langle a\rangle
\otimes k\langle b\rangle \to k\langle a+b\rangle$$ is the identity.

   \begin{remark}
A formal way to encode these signs can be expressed as follows (see~\cite{Kapranovsuper, dugger} for more comprehensive accounts of the underlying issues). The tensor-invertible graded vector spaces -- the full subcategory $\{k[a]\}_{a\in \Z}$ -- form the Picard groupoid  (commutative group object in groupoids) of $\Z$-graded lines, which is a {\em nontrivial} extension of $\Z$ by $B\Gm$, even though on the level of monoidal categories (and underlying group objects) this extension does split. In other words in a suitable higher-categorical sense the action of $\Z$ by shifts is only a projective action, but lifts to a genuine action of a $B\Gm$ (or in fact $B\Z/2$) extension.
By contrast the  Picard groupoids formed by the coreresponding super-lines form the trivial extension $\Z\times B\Gm$, and we have a genuine action of $\Z$ by shifts. 
 \end{remark} 
 Thus if we denote $$ M\langle a \rangle =  k\la a\ra \otimes M,$$ we again have identifications 
$$ M\la a\ra \otimes N\la b\ra \simeq (M \otimes N)\la a+b\ra$$ and $$ \Hom(M\la a\ra, N\la b\ra) \simeq \Hom(M, N)\la b-a\ra,$$ but now without any signs modifying the commutativity.  
  
\subsubsection{Shearing} \label{shearing notation} 
   
 Given a $\Gm$-equivariant complex $N\in \Rep(\Gm)$ of $k$-vector spaces, we define its {\em shear} as a $\Gm$-equivariant complex of {\em super} $k$-vector spaces by combining the cohomological and $\Gm$-grading of $N$, and modifying the parity accordingly:
 \[
N = \bigoplus_i N_i \mapsto N^{\shear} := \bigoplus_{i \in \Z} N_i\la i\ra \stackrel{\eqref{ultimate shearing}}{=} \bigoplus_{i \in \Z} \Pi^i N_i [i] \left( \frac{i}{2}\right)
\]
Here $N_i$ is the component upon which $\Gm$ acts by $\lambda \mapsto \lambda^i$. 
Here and below, as remarked after \eqref{BTT2}, our convention is that the Tate twist  $(i/2)$ embedded in the definition of $\la i \ra$ is  
to be ignored for the moment since we have no Frobenius action; but we will continue to write it because it is very helpful  to keep in mind for settings
where a Frobenius will be present. 

Note that even shifts don't involve parity shifts,  so that even iterates of the shearing operation 
don't require super vector spaces.
In general when the symmetric monoidal structure is not needed we can apply the forgetful functor from super 
to ordinary vector spaces, and abuse notation to write  \[
N  \mapsto N^{\shear} \]
for the resulting endofunctor of $\Rep(\Gm)$.
 
The process of shearing will be discussed in more detail in Section \ref{shearingsec0}.

\subsubsection{Shearing in the presence of a Frobenius action}  \label{shearing Frobenius notation}
    In the case that $N$ has an action of Frobenius (or other Galois or Weil groups) the Tate twist included in the shear modifies the Frobenius action.  
 
In particular,  observe that:
 
\begin{itemize}
\item If $N$ is pure, the action of Frobenius on the sheared vector space $N^{\shear}$
is obtained by introducing half-integral Tate twist to preserve purity; 
\item 
In particular, for $N$ in degree zero with trivial Frobenius action, 
  the action of Frobenius on $N^{\shear}$ will be very easy to remember and reconstruct from the shearing:
 vector spaces in cohomological degree $(-n)$ will be twisted by $k(\frac{n}{2})$. In particular, for odd cohomological degrees, this will depend on the choice of a square root $k(\frac{1}{2})$, i.e., on the fixed choice of square root $q^\frac{1}{2}\in k$. 
 \end{itemize}

    \begin{example} Take $G$ to be a split reductive group over a finite field. Then
    there is an isomorphism in the derived category of $\Q_{\ell}$-vector spaces
    with Frobenius action:
     $$H^*(BG, \Q_{\ell}) \simeq \left( \mbox{$G$-invariant regular functions on Lie algebra $\mathfrak{g}_{\Q_{\ell}}$}\right)^{\shear},$$
    where on the left we have geometric etale cohomology,
    and on the right  we regard  the  functions on the Lie algebra to have trivial Frobenius action
    and grade  according to the squaring action of $\Gm$ on the Lie algebra.
   Thus, for example, quadratic functions lie in degree $-4$, and are sheared to cohomological degree $+4$
    and with a Tate twist $-2$, so the statement says:
    $$H^4(BG, \Q_{\ell}) \simeq \mbox{(invariant quadratic functions on $\mathfrak{g}$)}(-2).$$
        \end{example}

   \subsubsection{Shifting of vector bundles} \label{shiftingofschemes}
$V$ continues to be a vector space over $\kk$, which we shall consider now as an affine scheme over $\kk$. 
The ``shifted vector space $V[-1]$'' is understood to be the derived scheme with ring of functions $\mathcal{O} =\Sym V^*[1]$,
which is abstractly an exterior algebra (note that there is no parity shift on the right). 

By contrast, when we attempt to apply positive shifts to $V$, we encounter  coordinate rings  such as $\Sym V^*[-2]$ which are {\em cohomological} (i.e., coconnective rather than connective or homological) graded rings, and thus are no longer affine objects in derived algebraic geometry -- we will use them only as convenient placeholders for their rings of functions and the corresponding categories of modules.
We warn the reader that another possible interpretation of $V[2]$ is the coaffine stack represented by this coconnective ring,
which has the same ring of functions, and this is {\em never} what we refer to: its category of quasicoherent sheaves is different from the category of modules
for the ring. 

For a variety $Y$, we denote by $TY$ and $T^* Y$ its usual tangent and cotangent bundle (or (co)tangent {\em complex} for $V$ singular).
We will denote by $T[-1]Y$ the shifted tangent bundle, i.e.
if $\mathcal{T}$ is the tangent sheaf, then $T[-1]Y $
is the relative spectrum of the symmetric algebra of $\mathcal{T}^*[1]$. 
Similarly we define other shifts; particularly important for us will be $T^*[2]Y$, the formal placeholder object whose ring of functions is 
the symmetric algebra of $\mathcal{T}[-2]$ and whose quasicoherent sheaves are modules over this ring. 
Note that $T[-1] Y$ is a derived scheme, whereas $T^*[2]Y$ is graded in cohomological degrees,
and in both cases there are no parity shifts unless explicitly stated otherwise. 

As an example of our notation, suppose that $\GGm$ acts on the vector space $V$ by scaling.
Then: 
\begin{itemize}
\item We do not use the object $T^*[2] V$ directly, but will
allow ourselves to refer to its ring of functions $\mathcal{O}_{T^*[2] V}$: this is the symmetric algebra on $V^* \oplus V[-2]$. 
\item The sheared ring
$\mathcal{O}_{T^*[2] V}^{\shear}$
is (therefore) the symmetric algebra on $V^*[-1] \oplus V[-1]$
where {\em both factors are put in odd parity}; thus
this ``symmetric algebra'' is infinite-dimensional.
\end{itemize}

 \subsection{Inner products of functions and sheaves}  
 
 The following lemma relates sheaf homomorphisms
 to inner product of trace functions. It will be used at various
 points in the text to connect  
 sheaf-theoretic and function-theoretic considerations. 
 
 \begin{lemma} \label{Homlemma} 
 Given two Weil sheaves $\mathcal{F}, \mathcal{G}$ on an $\mathbb{F}_q$-variety $X$, 
 let $f$ and $\check{g}$ be the trace functions associated to, respectively, $\mathcal{F}$
 and $D \mathcal{G}$; then 
  $$ \sum_{x \in X(\mathbb{F}_q)} f(x) \check{g}(x) = \mbox{(geometric) Frobenius trace on }\,  \Hom(\mathcal{F}, \mathcal{G})^\vee.$$
 Here, and in what is written below, we will keep denoting by $\Hom(\mathcal F, \mathcal G)$ the (derived) homomorphisms over the base change of schemes and stacks to $\FF = \overline{\FF_q}$  -- \textbf{not} on the corresponding stacks over $\mathbb F_q$, even if the sheaves descend there.

  The same holds true replacing $X$ by a quotient stack $X/G$ of $X$ by a connected linear
  algebraic group (hence, $\mathcal{F}, \mathcal{G}$ are Weil sheaves in the 
  equivariant derived category),
  where, on the left,    we understand
  each point $\FF_q$-point $x$ of $X/G$ to be counted with weight equal to $\frac{1}{\# G_x(\mathbb{F}_q)}$. 
 \end{lemma}
  
  Note that if $\mathcal{G}$ is pure of weight zero, then the left hand side
  is  (up to a power of $q$) the usual inner product $\sum_{x} f(x) \overline{g(x)}$. 
  Then the statement says precisely  that ``inner product'' corresponds to sheaf Hom
  under the function-sheaf dictionary.   
Taking $\mathcal{F} = \mathcal{G} = \underline{\kk}$ and $X$ smooth, the equality of the Lemma is the assertion
$$q^{-\dim X} |X(\mathbb F_q)| = \mbox{(geometric) Frobenius trace on } H^*(X, \kk)^{\vee}   $$
which follows readily from the Grothendieck--Lefschetz trace formula and Poincar{\'e} duality.

 \proof 
 We first take the case where $G$ is trivial. 
 Writing $\langle - \rangle$ for ``trace of Frobenius on'', and writing underline for the Hom-sheaf
 rather than its global sections, the right hand side above equals 
 \begin{multline*}     \langle D \Hom(\mathcal{F}, \mathcal{G}) \rangle   =  \langle  \mathbb{H}^*_c D \underline{\Hom}(\mathcal{F}, \mathcal{G}) \rangle \\=   \sum_{x \in X(\mathbb{F}_q)} \langle i_x^* D \underline{\Hom}(\mathcal{F}, \mathcal{G})  \rangle   =   \sum_{x \in X(\mathbb{F}_q)} \langle i_x^* (\mathcal{F}\otimes D\mathcal{G})\rangle,
 \end{multline*}
where we used the canonical isomorphism $D \underline{\Hom}(\mathcal{F}, \mathcal{G}) \simeq \mathcal{F}\otimes D\mathcal{G}$.  
 
 In the case of nontrivial $G$   we observe
 that the dual of $\mathcal{G}$ as a usual sheaf on $X$
 is obtained by shifting via $\langle +2 \dim(G) \rangle$ its dual as an equivariant sheaf (i.e., a sheaf on $X/G$);
 thus, the trace function $\check{g}_{\mathrm{eqvt}}$ associated to $D \mathcal{G}$ 
 dualized as an equivariant sheaf differs from the same function computed as a regular sheaf:
 $$ \check{g}_{\mathrm{eqvt}} = q^{\dim(G)} \check{g}$$
 On the other hand we have a spectral sequence 
 $$ H^i(BG; k) \otimes \Ext_{X}^j(\mathcal{F},   \mathcal{G}))  \langle -2 \dim G \rangle  \implies \Ext^{i+j}_{X/G}(\mathcal{F},    \mathcal{G}),$$
from which we see that the right-hand side \eqref{Homlemma} is altered from
the corresponding non-equivariant computation by the factor
$$ \mbox{Frobenius trace on $H^*(BG, \kk)^{\vee}$} =  \frac{q^{\dim(G)}}{\# G(\mathbb{F}_q)},$$
the equality verified by passing from $G$ to its reductive quotient (since the cohomology of connected unipotent groups is trivial), where one can explicitly compute in a standard way.  On the other hand, by the connectedness of $G$ (Lang's theorem), the $\FF_q$-points of $X/G$ coincide with $X(\FF_q)/G(\FF_q)$.
\qed

\begin{remark}\label{remarkconvergence}
  Note that, in the equivariant setting, the graded vector space $\underline{\Hom}(\mathcal{F}, \mathcal{G})$ is, in general, infinite dimensional. The lemma above entails the assertion that the corresponding series of Frobenius traces converges.
\end{remark}

\subsection{Analytic versus arithmetic normalization. Parity} \label{analyticarithmetic}
\index{analytic normalization} \index{arithmetic normalization}

In many areas of mathematics we see arising half-twists -- square roots of various types -- when we pass from functions
to half-forms. The latter manifestly have a unitary structure and the former
are often easier to manipulate because of the absence of a twist.    In our context we will  often use the words
\begin{quote} {\em analytic} and {\em arithmetic} \end{quote}
or occasionally
\begin{quote} {\em normalized} and {\em un-normalized} \end{quote}
to denote, respectively, the ``unitary picture'' and the ``untwisted picture.''

The origin of these terms  is the theory of $L$-functions, where
 are two standard conventions as to how to index an $L$-function.
One, which is often called ``analytic normalization'', has the property that the center of symmetry of the $L$-function lies at $s=1/2$,
and is commonly used in analytic number theory. 
The other, which differs from it from a half-shift, and could be called ``motivic'' or ``arithmetic' normalization, has the property
that the points of arithmetic interest are always integral $s \in \mathbb{Z}$. 
 Mathematically, these conventions are entirely equivalent. 
  
Correspondingly, all the principal conjectures of this paper can be expressed in two forms:
 
 \begin{itemize}
 \item Analytic normalization -- better suited to $L^2$-theory, tends to be naturally self-dual, but involves taking square roots
 (e.g.\ $\sqrt{q}$ is chosen). 
 \item Arithmetic normalization -- better suited to arithmetic situations, no need to choose $\sqrt{q}$. 
  \end{itemize}

We have chosen the analytic normalization as our primary way to formulate statements, although we have also
included in the text discussions of the arithmetic reformulations.   Both pictures have benefits:
\begin{itemize}
\item The analytic picture corresponds to a geometrically natural $\GGm$ grading on the hyperspherical spaces. 
\item The arithmetic picture interacts better with parity considerations, as discussed below. 
\end{itemize}

In any case, however we choose to 
express things, there are various half-twists
of various sorts embedded in the story. This brings us to the concept of parity.  \index{parity} Informally,
when we refer to ``parity'' conditions, they all have the following feature:
\begin{equation} \label{Parity Philosophy} \mbox{Parity is  extra structure
that causes all half-twists to cancel out.} 
\end{equation}
Our conjectures  {\em a priori} depend on square roots (square roots of $q$, square roots of canonical bundles...).
The validity of the conjectures is in fact independent of these choices;
this is not a formality, but  follows from certain constraints -- which we call ``parity'' constraints -- on the
  data entering into them. Such parity constraints exist, of course, in either the analytic or arithmetic version,
but they look more transparent in the latter.

Here are some places in the paper where the reader will find versions of the analytic/arithmetic dilemma:

\begin{itemize}
\item the arithmetic version of geometric Satake is
discussed in  \S \ref{sss:Satake-shearing} and \S \ref{sssderivedSatake};
more generally, the distinction between analytic and arithmetic is discussed
extensively in \S \ref{shearSatake} and \S \ref{Satake finite field}. 
\item 
the arithmetic version of the local conjecture is stated in Remark \ref{unnormalized local}.
\item In the discussion of the global conjecture, ``normalized'' periods correspond to the 
analytic picture, and unnormalized periods to the algebraic picture; 
\item
 the arithmetic formulation of the global conjecture is given in \S \ref{mitch}. 
 \end{itemize}

 \subsubsection{The grading $\GGm$} \label{Grading Gm}
  \index{$\GGm$}

 Throughout this paper, we will frequently encounter objects
 that have both a $G$-action and a commuting $\Gm$-action.
 We will often use the word ``graded'' to connote or remind
 the existence of this $\Gm$-action. 
 
 Now, the $\Gm$ has nothing to do with $G$; it is auxiliary, and we will 
use a special name for it:
 $$\mbox{$\GGm$ = an ``auxiliary'' copy of $\Gm$  often used to shear or regrade.}
 $$
 From a formal point of view, then, $\GGm$ {\em means exactly the same thing as $\Gm$, i.e.,  the multiplicative group $\GL_1$.}
 However, the notation is meant to hint to the reader that this plays a different role 
 to the reductive group $G$ or $\check{G}$.

Frequently, the parity considerations and the arithmetic/analytic dilemma mentioned above 
will be encoded by this $\GGm$. For example, parity will be related to the action of $-1 \in \GGm$, 
 and passage from analytic to arithmetic will often be effected by modifying the $\GGm$ action through a central co-character of $G$ or $\check{G}$.

 \subsubsection{Parity and superspaces} \label{supersloppy}
 Analytic formulations involve half-twists and $\sqrt{q}$,
 and, as we discussed   in \S \ref{Sqrtqsuper}, the
 appearance of super-vector spaces is essential to the algebraic formulation of such
 half-twists.

The super-structure can, however,{\em often be forgotten -- it needs to be remembered
only for certain specific purposes.}
   
 Here is a simple example. The shearing operation
 $$E \mapsto E^{\shear}$$
 briefly discussed in 
\S \ref{shearing notation}  and examined at more length in \S \ref{shearingsec0} 
can be used to define two closely related equivalences of categories:
\begin{itemize}
\item[(a)] a monoidal functor from the dg category of
complexes of $\Gm$-representations to itself, or

 \item[(b)] a {\em symmetric} monoidal functor
 from the same category, to the category $\Rep^{super}_{\epsilon}\Gm$ of $\Gm$-representations on super vector spaces
with the property that their parity coincides with the action of $-1 \in \Gm$ (\S \ref{commutativity}).
 \end{itemize}
 \index{$\Rep^{super}_{\epsilon}\Gm$}

We obtain version (a) from version (b) simply by forgetting the super-structure. 
Version (b) is appropriate if one is interested in issues involving the symmetric monoidal structure;
for us, the most important such issue is involving traces,
for shearing preserves traces only in the sense (b), as we already noticed.
However, if one is not interested in such computations, it is perfectly 
fine to work with version (a).

 In our text, all our equivalences of categories
can be formulated in the general form of (b). However, 
(b) is admittedly something of a mouthful,
and to avoid weighing down our discussion with super-vector spaces
 we will often formulate the equivalences in the form (a),
and then describe the  relevant information for form (b) 
in subsequent remarks.

 \subsection{Extended dual group} \label{subsection-extended-group}  \index{extended dual group}
 
 The material here is related to the considerations of parity mentioned above in \S \ref{analyticarithmetic}. 
We will  sometimes use an extended  version of the Langlands dual group -- called the \emph{$C$-group} in \cite{BuzzardGee}, and the \emph{Langlands-data group} in \cite{Bernsteinhidden}.
 See \S \ref{extended-group appendix} for a more thorough discussion of the role of the extended group in defining a cleaner form of the Langlands correspondence. 
  
  Let ${}^CG$ be the quotient of $G \times \GGm$ by the central element 
  $(e^{2 \check\rho}(-1), -1)$. 
 We remind that $\GGm$ denotes the group $\Gm$, but we use different notation in order to refer to this distinguished instance of the group, and sometimes refer to it as the ``grading group''. The canonical map ${}^CG \rightarrow \mathbb{G}_m$
  descending from $(g, t) \mapsto t^2$ will also be called the cyclotomic character.  
  For example, if $G=\SL_2$, then ${}^CG \stackrel{\sim}{\rightarrow} \GL_2 $ 
  via $(g, \lambda) \mapsto \lambda g$, 
  and the ``cyclotomic character'' just defined is the determinant. 
With this notation we have an exact sequence
  $$ 1\rightarrow G \rightarrow {}^C G \rightarrow \Gm\rightarrow 1.$$ 
  We apply the same definitions to the Langlands dual group $\check G$ of $G$ to define the {\em $C$-group} ${}^C \Gv$, though note that the $C$-group is {\em not} the Langlands dual of the group ${}^CG$.
  
  \index{$C$-group}
  \index{${}^CG$} \index{${}^C\Gv$}
  
For a $G \times \Gm$-space (such as the graded Hamiltonian
spaces that are the main concern of our paper), the condition that the $G \times \Gm$-action on $M$ descend
to ${}^C G$  is the condition
that $e^{2\check\rho}(-1) \in G$ acts the same way as $-1 \in \Gm$. 
This condition,
or variants of it, will arise often as parity conditions in the sense of \S \ref{analyticarithmetic}. 
To keep track of these variants, we will want to use
to other central elements besides $e^{2\check\rho}(-1)$ and accordingly   given a central involution $z \in G$,
we will occasionally use the following notation:
\begin{equation} \label{CGz def}   {}^CG_z:= \mbox{quotient of $(G \times \GGm)$ by the central element $(z, -1)$. }\end{equation}
  The notation is regrettably heavy, but will be used only in a very few points.

  \index{${}^CG_z$}

\subsection{Langlands parameters, extended Langlands parameters, and their $L$-functions} \label{Ldef}
  In what follows, we restrict $k$ to be an algebraic closure of $\mathbb{Q}_{\ell}$ and we fix an isomorphism $k \simeq\C$,
 and $F$ will be a global function field. 

A Langlands parameter $\phi$ is a Frobenius-semisimple morphism $\Gamma_F \rightarrow \check{G}(k)$. 
A Langlands parameter into $\check{G}=\GL_n$ gives rise to  an $L$-function
$L(s, \phi)$, defined by the 
analytic continuation of the usual Euler product with factors 
the inverse characteristic polynomial of geometric Frobenius on inertial invariants.
For example if $\phi$ is trivial the resulting function is the zeta-function 
of the field $F$ and has poles at $s=0$ and $s=1$.

\index{extended Langlands parameter}
We will  usually use the  term {\em extended Langlands parameter} to mean 
\begin{itemize}
\item
a Langlands parameter with the role of $\check{G}$ replaced by $\check{G} \times \Gm$,  projecting to the positive square root $\varpi^\frac{1}{2}$ of the cyclotomic character in $\Gm$ 
(``positive'' with reference to an isomorphism $\kk \simeq \C$).
\end{itemize}

Obviously, there is a bijection between Langlands parameters and extended Langlands parameters, once $\kk\simeq \C$ (or simply the square root of $q$) has been fixed. 
In practice, all that matters will be the projection of this parameter to a
$C$-group ${}^C\Gv_z\simeq \check{G} \times \Gm/\mu_2$
of \S \ref{subsection-extended-group}, i.e., a homomorphism from $\Gamma_F$
to ${}^C \Gv_z$  projecting
to the cyclotomic character under $\Gm/\mu_2 \simeq \Gm$,
and we will also allow ourselves occasionally to use ``extended Langlands parameter''
in this context.  In the arithmetic context, this type of parameter is more canonical, while the other versions are obtained from it by noncanonical choices of $\varpi^\frac{1}{2}$, $\kk\simeq\CC$, etc.
 However, since the relevant parity element $z$ will vary, it is more convenient
 for us to work with the less intrinsic form above.

Our convention will be that:
\begin{itemize}
\item $\phi_E, \psi_{E} \dots$ will denote an {\em extended} Langlands parameter, and 
\item  $\phi_L, \psi_L, \dots$ will denote {\em usual} Langlands parameters, and
\item $\phi_A, \psi_A, \dots$ will denote {\em Arthur} parameters,
\end{itemize}
with the indices omitted when the type of parameter used is clear from the context. 

Moreover, 
\begin{itemize}
\item $f_{\phi}$ will denote an automorphic form attached to a parameter $\phi_E/\phi_L/\phi_A$.
\end{itemize}
When $\phi$ is fixed or clear from the situation we will abridge this simply to $f$.

Given an extended Langlands parameter 
$$ \phi: \Gamma \rightarrow \check{G} \times \Gm$$
and a representation $\check{G} \times \Gm \rightarrow \GL(V)$,
we define the associated $L$-function to be the
composite Langlands parameter into $\GL(V)$.  We denote this by\footnote{We will use the notation $L(V^{\shear}, s)$ in situations where
there is no ambiguity as to what $\phi$ could be.}
\begin{multline*} L(\phi, V^{\shear}, s) \mbox{ or } L(V^{\shear}, \phi, s) \mbox{ or } L(V^{\shear}, s) := \\ \mbox{$L$-function of the composite $\Gamma_F \rightarrow \check{G} \times \Gm \rightarrow V$.}
\end{multline*}
 The $\shear$ in the notation reminds us that the $\Gm$ action is relevant here. 
 More explicitly, the $\Gm$ action on $V$ grades it as $V = \bigoplus V_k$, and then
\begin{equation} \label{Lsheardef} L(V^{\shear}, s) = \prod_{k} L(s+\frac{k}{2},V_k).\end{equation}
 
 Indeed, we can think of  $L(V^{\shear}, s)$ as either:
 \begin{itemize}
 \item  an Euler product, whose local factors come from characteristic polynomial of Frobenius on the sheared local system $V^{\shear}$, or
 \item the characteristic polynomial of Frobenius on $(H^* V^{\shear}) = (H^* V)^{\shear}$, where $H^* V$ denotes the {\'e}tale cohomology of the curve
 associated to $F$ with coefficients in the local system defined by $V$, and $(H^* V)^{\shear}$ denotes its shearing with respect to the grading
 induced by that on $V$.  
 \end{itemize}
In both cases it is important that we take account of the super-structure on shearing in considering determinants, e.g.\ 
 the characteristic polynomial of Frobenius on $V \langle 1 \rangle$ equals $\det(1-q^{-1/2} \Fr|V)$ and {\em not} its inverse.

  Other notation on $L$ and $\epsilon$-factors is set up in \S \ref{Leps}.

\subsection{Function spaces, left and right actions} \label{leftrightconventions}

We will want the freedom to use both left and right actions for $G \times \GGm$ acting on $X$ or $M$.  \index{left versus right actions}
To pass between the two  we shall use  equivalence of categories 
\begin{equation} \label{leftrightequiv} \mbox{left $G \times\GGm$ spaces} \leftrightarrow \mbox{right $G \times \GGm$ spaces}\end{equation}
 wherein we invert the $G$-action {\em but not the $\GGm$-action}.   
  
In Parts 2 and 3 of the paper (local and global theory), our
convention about these will be
 \begin{equation} \label{leftrightequation }\mbox{ Automorphic actions on right, spectral actions on left. } \end{equation}
 That is to say, when considering the role of $(G, M)$ on the automorphic side,
 we will use right actions, and when considering the role of $(\check{G}, \check{M})$ on the spectral side, we will use left actions. 
 The switch is somewhat unfortunate, but the former is more in line with conventions about automorphic periods, and the latter is more in line
 with conventions about $L$-functions.

In any case, for the action of $G$  and $\GGm$ on functions, sheaves, forms etc. derived from this we will {\em always use left actions,}
derived in the standard way from the left   or right  actions of $G$ on $M$, and the action of $\GGm$. 

For a a symplectic vector space $(V,\omega)$, the identification $\iota: V\simeq V^*$ is defined by contraction in the first variable, i.e., $\left<\iota(v),w\right> = \omega(v,w)$. 
Similarly, for a symplectic manifold $M$ the isomorphism $T_x M \simeq T_x^*M$ is defined by contraction in the first variable, so that the Hamiltonian vector field $X_H$ associated to a function $H$ satisfies $dH(Y) = \omega(X_H, Y)$. 

  \index{moment map convention}
Given $G$ acting by symplectomorphisms on $M$, a moment map for the $G$-action will then be a map $\mu: M \rightarrow \mathfrak{g}^*$
with the property that for any $X\in \mathfrak g$, the vector field $\underline{X}$ defined by the infinitesimal action of $X$  
satisfies 
\begin{equation} \label{momentmapdef} \omega(\underline{X}, Y) = d\mu_X(Y),\end{equation} 
where $d\mu_X$ is the pairing of $\mu$ with $X$ and $Y \in TM$ is any tangent vector. Note that this implies, in particular, that when we translate between left and right actions as above, the moment map gets \emph{reversed}.
\index{symplectic space convention}
 In particular for right actions:

\begin{itemize}
\item[-]
 For $M = T^* X$ such a moment map is given simply 
by dualizing the orbit map $\mathfrak{g} \rightarrow T_x X$. To make this valid,  
we adopt the following normalization for the symplectic form: writing $\theta$ for the tautological $1$-form on $T^*X$,
which pairs tangent and cotangent directions, we take the symplectic form to be
$$ \omega_{T^*X} := -d\theta.$$

\item[-]
For a symplectic vector space $(M,\omega)$, the moment map
 paired with $Z \in \mathfrak{sp}_M$ is the function 
 \begin{equation} M\ni m\mapsto \label{smm} \frac{1}{2} \omega(Zm, m).\end{equation}

 \item[-]
The cotangent bundle $T^*G$ has a right action of $G \times G$ which arises from the action $(x,y): g \mapsto x^{-1} g y$
of $G \times G$ on $G$.  Let us label the two copies of $G$ acting as $G_l \times G_r$
to distinguish ``left'' and ``right.''
If we identify $TG = \mathfrak{g} \times G$ according to the rule where
\begin{equation} \label{TGidentity} (X \in \mathfrak{g}, g \in G)\end{equation} is sent to the derivative of the curve $e^{tX} \cdot g$ at $t=0$,
and correspondingly identify $T^* G = \mathfrak{g}^* \times G$, then the moment map for
the $G_l \times G_r$ action on $T^*G$ just mentioned is 
$$(\xi \in \mathfrak{g}^*, g) \mapsto (-\xi, \Ad(g) \xi) \in \mathfrak{g}_l \times \mathfrak{g}_r$$
and the action of $G_l \times G_r$ on $T^*G$ is given by
$$(g_l, g_r) \cdot (\xi, g) =  (\Ad(g_l^{-1}) \xi, g_l^{-1} g g_r).$$

 Here, $\mathrm{Ad}$ denotes the \emph{left} (co)adjoint action of $G$, which we will also denote by $\Ad(g) Z = {^gZ}$. We will denote the corresponding right action by $Z^g :=\mathrm{Ad}(g^{-1}) Z$. 
\end{itemize}

\subsubsection{Inner products}

We will abuse the terminology and notation for inner products, and use them for the corresponding \emph{bilinear} forms instead. That is, for measurable functions $f, g$ on a measure space $E, \mu$, we set
$$ \langle f, g \rangle := \int_{E} f(e) g(e) d\mu_e.$$
This will make it easier to pass between algebraic considerations and analytic ones. 

\subsection{Categorical background} \label{cat0}

We will make use of higher categories throughout the paper.
We will give details in \S \ref{HigherCatAppendix}; at the moment we just point out some key features. 

 First of all, some words to orient the reader  (particularly those whose background is the arithmetic Langlands program) as to why we are using such language.
 Loosely speaking, in the geometric statements that we study, categories play the role of function spaces; and then: 
 \begin{itemize}
 \item 
The fact that there are several different options for categories of sheaves on, e.g., the space of $G$-bundles should be thought
of as related to the fact that there are many reasonable topological vector spaces incarnating the space of functions
on an ad\`elic quotient $G_F \backslash G_{\mathbb{A}}$. As in the latter case, {\bf this type of detail is important
to make a mathematically precise statement, but probably should be ignored at a first reading, for it does not carry
the essential content of the conjecture. }
 
\item  Why not work with triangulated categories?  Unfortunately it is a well-known problem that
it is not easy to perform natural categorical operations on triangulated categories,
so for internal arguments it is extremely convenient to work with dg enhancements. Again, the reader
can ignore this at a first reading, for the distinction between dg and triangulated categories
is of a nature orthogonal to the essential content of our conjectures.   
\end{itemize}

\index{category over $\kk$}
Unless otherwise specified, {\em category} will always mean a differential graded (dg) category over a field $\kk$ of characteristic zero. Moreover our dg categories will always be stable (i.e., pre-triangulated) -- such categories are equivalently described as stable $\kk$-linear $\infty$-categories, and are closed under finite (homotopy) limits and colimits. Moreover, dg categories are considered up to quasi-equivalence, i.e., we work in a suitable $\infty$-category of dg categories (again see \S \ref{HigherCatAppendix} for details).

\subsection{Categories of sheaves} \label{intro sheaf cat}
  \index{$\Shv$} \index{category of sheaves}
Categories of sheaves will play an important role. 
As just mentioned, all our categories of (coherent or constructible) sheaves on various spaces will be, by definition and without extra mention, \emph{$k$-linear dg-categories}. 
These come in a variety of flavors, but a few general words about the notation:

 Our notation for ``coherent'' sheaves is straightforward: For $X$ a (derived) stack over $\kk$,  $QC(X)$ and $QC^{!}(X)$ denote, respectively, quasi-coherent and ind-coherent sheaves on a scheme or stack. These are large versions of the categories $\Perf(X)$ and $\Coh(X)$ of perfect complexes and bounded coherent complexes. In between lie the categories $QC^!_\Lambda(X)$ (and $\Coh_\Lambda(X)$) in which we specify singular support in the sense of~\cite{ArinkinGaitsgory}.  For details see \S \ref{coherent sheaf theories}.
 
Constructible sheaves will come with a depressing number of variants.   These are described in more detail in \S \ref{constructible
sheaf theories} and \S \ref{sheaves on stacks} and we will summarize some important points.  For $Y$ a (derived) stack over $\FF$,
$$ \Shv^{?}_{\invertQM}(X) \mbox{ or } \SHV^{?}_{\invertQM}(X)$$
will denote a LARGE (=presentable, see Appendix~\ref{HigherCatAppendix}) or small \footnote{In some situations where a category of sheaves is denoted by calligraphic font, e.g.\ $\mathcal{H}$, we use the notation $\overline{\mathcal{H}}$ for the corresponding large category.} dg-category of ``constructible'' sheaves on $X$ with $\kk$-coefficients,  and $?, \invertQM$ will be various adornments modifying the category.
The  options are:
\begin{itemize}
\item $? = B, dR, et$: records whether we are in the Betti, de Rham, or {\'e}tale settings:
\begin{itemize}
\item Betti, see \S \ref{CST:B}: only applies when $\FF=\C$, but with $\kk$ arbitrary. Sheaf theory is built from sheaves of $\kk$-vector spaces on $Y(\C)$ for the complex analytic topology, and in particular contains constructible complexes.
\item de Rham, see \S \ref{CST:DR}:  only applies when $\FF=\C, \kk=\C$. Sheaf theory is built from  $D$-modules on $Y$. 
\item et, see \S \ref{CST:ET} applies for any $\FF$, with coefficients  
 $k=\overline{\mathbb Q_l}$. Sheaf theory is built from  constructible $l$-adic \'etale sheaves.
\end{itemize}
 \item  $\invertQM = s$ denotes ``safety,'' and its absence
 denotes that  we work with the ind-finite category,   see \S \ref{renormalization section}. 
\item $\invertQM=\Lambda$ records that we consider sheaves with fixed {\em singular support} in the sense of~\cite{KashiwaraSchapira}. 
\end{itemize}

 \subsubsection{Notation for automorphic sheaves}
 To avoid having to keep track of this bewildering array of notation,
 we will specify in Appendix
\ref{geometric Langlands}
 a standardized set of options for ``automorphic sheaves.''
 Thus, when we write e.g.\ 
 $$ \Aut(\Bun_G)$$
 in the body of the paper, 
 we actually have implicitly chosen various adornments $?, ?$, which
 depend on the context (Betti versus de Rham versus finite) in which we are working; 
 the reader should refer to  Appendix \ref{geometric Langlands} as necessary to recall.

 \subsection{Basic notation} \label{miniindex} 
The Index on p.~\pageref{index}  contains many of the notations we use; some of the more persistent ones include the following:
 
\begin{description}

 \item[$\mathcal P_X, \mathcal L_X$] Period and $L$-sheaves (\S \ref{periodX}, \ref{LsheafX}). 
 \item[$P_X, L_X$] Period and $L$-functions (\S \ref{periodX}, \ref{locXex}).  \index{$P_X$} \index{$L_X$}
 \item[$k$] Coefficient field (\S \ref{coefficients}).
 \item[$\mathbb F$] (Algebraically closed) field of definition of the curve (\S \ref{coefficients}).
 \item[$\Sigma$] Curve. 
 \item[$\ssslash$] Hamiltonian reduction (\S \ref{symplecticinduction}). 
 \item[$\sslash$] GIT quotient.
 \item[$A \times^H B$] Contracted product, i.e., then $A$ is a space with a right action of a group $H$, and $B$ is a space with a left $H$-action, this is the quotient of $A\times B$ by the equivalence relation $(ah,b) \sim (a, hb)$ (for $a\in A, \, b\in B, \, h\in H$). Equivalently, it is the quotient of $A\times B$ by the diagonal right $H$-action $(a,b)\cdot h = (ah, h^{-1} b)$. Often, $A$ and $B$ will be given with right actions (in which case it is the quotient by the diagonal $H$-action), and typically they will be schemes, in which case the quotient is understood as a stack-theoretic quotient (although, in almost all cases, the action will be free and the quotient is again a scheme).
 \item[$Alg(\mathcal{C})$]  Algebra objects {\em in the category $\mathcal{C}$}.
\index{Alg}
\end{description}
\index{$A\times^H B$}
\index{$k$}\index{$\mathbb F$} \index{$\sslash$} \index{$\ssslash$}

%% file: spherical.tex
In the next couple of sections, we discuss reductive groups and Hamiltonian spaces defined over an \emph{algebraically closed field $\FF$  in characteristic zero}.

However, at the end of each section, we will discuss various issues related to the general case where the field of definition is not algebraically closed, in particular \S \ref{hdprings}, \S \ref{CheckMRat}, and
 \S \ref{dhpFq}. We are particularly interested in the case of $\FF=\overline{\FF_q}$, with the varieties defined over $\FF_q$. \footnote{Note that by general spreading out arguments,
 one can pass results from characteristic zero to ``large enough'' characteristic: the spaces that we consider will admit models over the $S$-integers of an algebraic number field, where $S$ is a finite number of places, and can then be reduced to finite fields; but we will be more precise where possible. }

\section{Hyperspherical Hamiltonian spaces}  \label{sphvar} 
 
\subsection{Introduction} \label{sphvarintro}
Let $M$ be a Hamiltonian $G$-space over an algebraically closed field $\FF$  in characteristic zero, by which, in this paper, we will always mean a \emph{smooth, symplectic} variety with a $G$-equivariant moment map $\mu:M\to \mathfrak g^*$. 
Sign conventions related to moment maps and symplectic spaces have been discussed in \S \ref{miniindex}.  
 Our Hamiltonian spaces will almost always be {\em graded:} equipped with a commuting action of the grading group $\GGm$
(cf. \S \ref{Grading Gm}) compatible with the action on $\mathfrak g^*$ and on the symplectic form \emph{by the square character}. 

\begin{example}
 If $X$ is a smooth $G$-variety, we can grade its cotangent bundle $T^*X$ by letting $\GGm$ act on the fibers by the square character.
\end{example}

The goal of this section is to {\em explicate a class of such graded Hamiltonian $G$-spaces that is well-adapted to relative Langlands duality.}  We do not claim this is the correct generality for the story,
only that it seems to be a context where the relative Langlands duality plays out nicely. It is {\em plausible} that relative Langlands duality 
in fact gives an exact duality on a slight restriction of this class, and we formulate our preliminary expectations
on this issue in \S \ref{hdpC}.  

In more detail:  \begin{itemize}
\item  \S~\ref{model examples} describes a class of ``model examples'' to motivate the later reasoning.
\item   \S~\ref{symplecticinduction} gives background on the processes of Hamiltonian (also known as symplectic) reduction and induction. 
\item   \S~\ref{Whittaker induction}
describes the process of Whittaker induction, which constructs a Hamiltonian $G$-space from a datum
 \begin{equation} \label{basic0} (H \subset G, \;  \mathfrak{sl}_2 \to \mathfrak{g}^H, \; S: \textrm{ a Hamiltonian $H$-space}).\end{equation}

Whittaker induction carries vector spaces $S$ to 
   vector bundles over $H \backslash G$
 as is explained in \S~\ref{Slodowy}.
 
  The class of Whittaker inductions, and even
  the subclass of spaces Whittaker-induced from symplectic representations of $H$, contains all spaces of interest to us in this paper, but is too big; the
  remainder of the section addresses this.
  
 \item   \S~\ref{ssconditions}
 axiomatically describes a subclass of Hamiltonian actions, the class of ``\hyperspherical'' Hamiltonian spaces, which is most relevant for our story. 

\item\S~\ref{ss:structuretheorem} proves that all \hyperspherical spaces arise as Whittaker inductions. 
 \item \S~\ref{dp} examines when a \hyperspherical  space can be polarized. 
\item \S~\ref{eigencharacter} discusses volume forms on a polarized \hyperspherical space. 
\item \S~\ref{hdprings} discusses rationality issues.

  \end{itemize}
 \index{parity}

 \begin{remark}
  There is one important feature of graded Hamiltonian spaces whose discussion is deferred to the next section
 (\S~\ref{parity}), and that is the issue of {\em parity}, cf. \eqref{Parity Philosophy},
 \S \ref{analyticarithmetic}, and the discussion above \eqref{CGz def}: 
 It will be important in our later examples  that there is a central involution $z \in G$ which acts on $M$ as the element $-1 \in \GGm$. But we will not impose this explicitly as part of the definition.
\end{remark}

 \begin{remark} \label{fancy induction} \index{shifted symplectic geometry}
Many of our constructions, both here and in later sections, are conveniently expressed in terms of ``the category of Hamiltonian spaces and Lagrangian correspondences.''
Unfortunately, for reasons of transversality, this does not form a category in the classical setting; 
see~\cite{Weinstein} for discussion. 
This problem can be resolved by passing to derived geometry:
the language of shifted symplectic geometry~\cite{PTVV,Safronov}, specifically the higher category of Lagrangian correspondences of shifted symplectic stacks constructed in~\cite{haugseng} (see~\cite{calaquesymplectic}),
 give a convenient setting in which to work. 
In this language, $\mathfrak{g}^*/G=T^*[1]\text{pt}/G$ has a 1-shifted symplectic structure and Hamiltonian $G$-spaces $M$ are identified with Lagrangians in this space -- 
more precisely with the structure of shifted Lagrangian (see~\cite{calaquelagrangian}) on the equivariant moment map
$$M/G \rightarrow \mathfrak{g}^*/G.$$
We will not explicitly use this language in the current section, but it is often helpful later in the paper
when we come to consider quantizations. 
  \end{remark}

\subsection{Some motivating examples} \label{model examples}
  
 Here are three important 
 examples of graded Hamiltonian $G$-spaces:
  \begin{itemize}

 \item[(a)] {\em Spherical case.}
Recall that a normal $G$-variety $X$ is {\em spherical} if a Borel subgroup of $G$ acts with a Zariski open orbit of $G$. 
 Attached to a smooth spherical variety $X$, we have the Hamiltonian $G$-space
 $$ M =T^* X$$
 with $\GGm$-action squaring along the fibers.

 \item[(b)]  {\em Whittaker-type cases.}
 
  The basic
 example here is obtained by twisting $T^*(U\backslash G)$ (where $U$ is a maximal unipotent subgroup) by an additive character $\psi:  U \rightarrow \mathbb{G}_a$:
 we take $M$ to be the Hamiltonian reduction of 
 $T^*G$ by the character $\psi$, which is to say, we consider
 the preimage of $d\psi \in \mathfrak{u}^*$ under the moment map
 $T^* G \rightarrow \mathfrak{u}^*$, and take its quotient by $U$.
The action of $\lambda \in \GGm$ composes left translation
by $\lambda^{2 \check{\rho}}$ on $U\backslash G$ with the squaring action along fibers.

 \item[(c)] The {\em vectorial case}; here 
 \begin{equation} \label{stdvect} (M,\omega) \end{equation} is a symplectic vector space, $G\subset \Sp_M$, and the moment map $M \rightarrow \mathfrak{g}^*$  
factors through $M \rightarrow \mathfrak{sp}_{M}^*$ sending $m \in M, X \in \mathfrak{sp}_M$
to $\frac{1}{2} \langle Xm, m \rangle$.  Here, $\GGm$ acts on $M$ by linear scaling.

  \end{itemize}
  
 The earlier work of Y.S.\ and A.V.\ \cite{SV} focused entirely on  the first 2 cases.
 However, this class is certainly not closed under the conjectural duality that we want to introduce,
 and it does not include many important examples in the theory of automorphic forms, such as the theory of the theta correspondence; as was observed in \cite{SaHowe}, the conjectures of \cite{SV} extend to that case. 
 
 \subsubsection{Twisted cotangent bundles}\label{twisted cotangent bundle section} 
 \index{twisted cotangent}
 \index{$\Psi$} \index{$\Ga$-torsor} \index{twisted polarization}
 It will be useful to reformulate case (b) in  a setting that is closer to that of (a),
 by considering   {\em twisted cotangent bundles} associated to affine bundles. Although this will be a special case of a more general construction, we single it out here for its recurrent appearance in this paper, and point the reader to \S~\ref{symplecticinduction}, \ref{Whittaker induction} for more general constructions, and a recollection of notions such as ``Hamiltonian reduction.''
 
 Let
$$ \Psi \rightarrow X$$
be an equivariant $\Ga$-torsor over a $G \times \GGm$ variety $X$, where
$\GGm$ acts on $\Ga$ by the character $x \mapsto x^{2}$. 
To spell out, this means:
 \begin{itemize}
 \item $\Psi$ is a $\Ga$-torsor over $X$, which is to say
that it is equipped
 with an action of $\Ga$ over $X$ and is {\'e}tale locally (on $X$) isomorphic
 to $\Ga \times X$;  
 \item The action of $ \GGm$
 lifts to an action of $\Ga \rtimes  \GGm$ on $\Psi$, both commuting with $G$. 
  \end{itemize}

  The $\GGm$-action on $\Psi$ induces an action on $T^* \Psi$; this
  is equivariant for the moment map where $\GGm$ acts on $\mathfrak{g}_a^*$ by the character $\lambda\mapsto\lambda^{-2}$. 
We modify this $\GGm$-action by composing it with the commuting $\GGm$-action wherein $\lambda \in \GGm$
scales fibers of $T^* \Psi$ by $\lambda^2$. The resulting $\GGm$-action is equivariant for the 
  moment map
  $$T^* \Psi \rightarrow \mathfrak{g}_a^*$$
  wherein $\GGm$ acts trivially on the target. 
 We now define the graded Hamiltonian $G$-space, the ``twisted cotangent bundle'' of $(X,\Psi)$, \index{twisted cotangent}
 $$ M := T^*(X, \Psi):= T^*_\Psi X:=  T^*\Psi\ssslash_1 \Ga$$ 
 to be the Hamiltonian reduction of the cotangent bundle of the total space $\Psi\to X$ of $\Psi$  
 under the action of $\Ga$, at the moment map value $1 \in \mathfrak g_a^*$.

 \begin{example}
 Here is how the Whittaker case fits into this framework. 
 We may take $X = U \backslash G$ (where $U$ is a maximal unipotent subgroup), and $\Psi = U_0 \backslash G$ where $U_0$ is the kernel
 of a generic additive character $U \rightarrow \Ga$; and the $\GGm$-action is given
 by left multiplication by $2\check{\rho}$.  The $\Ga$-action arises from 
 left multiplication by $U/U_0 \simeq \Ga$. 
 \end{example}

\begin{example}The twisted cotangent bundle $M$ above is a torsor under the usual cotangent bundle $T^*X$;
the transition functions are obtained from differentials of the transition functions for the $\Ga$-bundle $\Psi$. 
A situation more commonly encountered in representation theory is the twist of a cotangent bundle
associated to a {\em line} bundle, wherein the transition functions are obtained by $d\log$ of the transition
functions of a $\Gm$-bundle.
\end{example}
  
 \subsection{Hamiltonian reduction and induction}   \label{symplecticinduction}
We review the operations of reduction and induction of Hamiltonian spaces, usually called \emph{symplectic reduction/induction} -- but we will use the term \emph{Hamiltonian reduction/induction}, to emphasize the dependence on the Hamiltonian structure, i.e., on the moment map. 
 \index{Hamiltonian induction}

The Hamiltonian reduction of the Hamiltonian $G$-space $M=T^*X$ given as the cotangent bundle to a $G$-space $X$ is the cotangent bundle of the quotient $X/G$ (assuming the quotient exists in the desired category; for the purposes of this paper, it is enough to consider quotients which are  schemes). Modelling on this example, we define Hamiltonian reduction in general as  \index{$\ssslash$}
\[ M\ssslash G = M \times_{\mathfrak g^*}^G \{0\},\]
the quotient of the fiber of $0$ under the moment map by $G$. More generally, we may reduce at a different element $f\in \mathfrak g^*$, by the formula 
\[ M\ssslash_f G = M \times_{\mathfrak g^*}^G \mathcal O_f,\]
where $\mathcal O_f$ is the $G$-orbit of $f$. (In particular, this operation depends only on the coadjoint orbit of $f$.) In general, the Hamiltonian reduction is a derived symplectic stack, but we will only use it in cases where the action is free, hence is a symplectic variety.

Similarly, suppose $H \subset G$ is an inclusion\footnote{
 We may also replace  the inclusion $H \subset G$ by an arbitrary morphism but will have no need of this.
  } of algebraic groups and $S$ a Hamiltonian $H$-space,
the {\em Hamiltonian induction} of $S$ to $G$ is the semiclassical version of induction of unitary representations.  In the case when $S$ is the cotangent bundle $T^*Y$ of a smooth $H$-variety then the Hamiltonian induction is simply the cotangent bundle of the $G$-space $Y\times^H G$ induced from $Y$.

Modelling on this example, we define the ``Hamiltonian induction''
\[\textrm{h-ind}_H^G(S) := (S\times T^*G)\ssslash H,\]
the Hamiltonian reduction of $S \times T^*G$ under $H$.
Here, $T^* G$ is considered as a right Hamiltonian $H$-space, with the action induced from the action $h: g \mapsto h^{-1} g$
  of $H$ on $G$.

   Equivalently, we can consider the \emph{left} action of $H$ by left multiplication on $G$, $h\cdot g = hg$, and the induced Hamiltonian \emph{left} structure (note that this changes the $H$-moment map by a factor of $(-1)$), and then the right-hand side is given by the contracted fiber product 
 \begin{equation}
 \label{cfp} 
 M = S \times_{\mathfrak{h}^*}^H T^*G, \end{equation}
 that is, the fiber product over $\mathfrak h^*$, divided by the equivalence relation $(sh, \xi)\sim (s, h\xi)$ (for $s\in S,\, \xi\in T^*G,\, h \in H$, such that $sh$ has the same moment image as the left-$H$-action moment image of $\xi$). 
  Taking account of our identifications of $T^*G$, see \eqref{TGidentity}, we 
see that this can also be written as 
\begin{equation} \label{cfp2} M \simeq (S \times_{\mathfrak{h}^*} \mathfrak{g}^*)\times^H G,\end{equation}
where $H$ acts on $\mathfrak{g}^*$ by the right coadjoint action $h: \xi \mapsto \mathrm{Ad}(h^{-1}) \xi$, and the moment map for $M$ is induced by the right coadjoint map $\mathfrak{g}^*\times^H G\to \mathfrak g^*$.
 
 In particular, the Hamiltonian induction comes with the structure of a fiber bundle $\text{h-ind}_H^G(S) \rightarrow  H \backslash G$.
 If $S$ carries the structure of graded Hamiltonian $H$-space (\S \ref{sphvarintro}), then so does the Hamiltonian induction, using the diagonal $\Gm$-action on $S\times T^*G$, which commutes with $H$ and the moment map (when $\GGm$ acts on the fibers of $T^*G$ by the square character).  
 \index{symplectic normal bundle}
 A relevant notion to Hamiltonian induction is that of the \emph{symplectic normal bundle} to a $G$-orbit $\mathcal O$ in a symplectic manifold $M$: It is a vector bundle over $\mathcal O$, whose fiber over $x\in \mathcal O$ is equal to the space $S=T_x \mathcal O^\perp/(T_x \mathcal O^\perp\cap T_x\mathcal O)$. This is a symplectic vector space, equipped with an action of the stabilizer $H=G_x$ of $x$, hence with a quadratic moment map $S \to \mathfrak h^*$. When $M = \textrm{h-ind}_H^G(S)$, where $S$ is a symplectic $H$-vector space with the quadratic moment map, the symplectic normal bundle construction at the orbit $\mathcal O = \{0 \} \times^H G\subset M$ recovers $S$. See Remark~\ref{SS2} for the closely related operation of {\em Hamiltonian restriction}.

 \subsubsection{Recognizing Hamiltonian induction} \label{recognizinginduction}
 We spell out a property, analogous to Frobenius reciprocity, that will be used later to recognize a Hamiltonian induction. 
Write $L = S \times_{\mathfrak{h}^*} \mathfrak{g}^*$; then, by \eqref{cfp2},
$L$ is embedded in $M=\text{h-ind}_H^G(S)$
as the fiber above the identity of $M \rightarrow H \backslash G$. 
The symplectic form restricted to $L$ is pulled back from $S$; 
said differently, $L$ defines a Lagrangian correspondence 
\begin{equation} \label{Lc0}
  M^\circ \leftarrow L  \to  S
\end{equation}
where $M^\circ$ is the space $M$ with the opposite symplectic form.
  
Now take {\em any} Hamiltonian $G$-space $M$
equipped with such an $H$-stable Lagrangian correspondence $L$ as in \eqref{Lc0}, and with the property that the composites
\[L \rightarrow M\to \mathfrak{g}^*  \rightarrow \mathfrak{h}^*, \,\, L \rightarrow S \rightarrow \mathfrak{h}^*\]
coincide.

 Then, there is an induced Lagrangian correspondence  
\begin{equation} \label{inducedLL} M^\circ \leftarrow L\times^H G \rightarrow \textrm{h-ind}_H^G(S),\end{equation} 
 compatible with the moment maps of $M$ and $\text{h-ind}_H^G(S)$. Indeed, consider first the case $H=\{1\}$, and define the map as
 \[ L \times G \ni (l, g) \mapsto (m(l)g, s(l), \mu(l), g) \in M^\circ\times S \times \mathfrak g^* \times G.\]
 Here, $m\times s$ is the given correspondence, and $\mu$ is the pullback of the moment map from $M$. It is easily verified that this is a Lagrangian. In the presence of an $H$-action, this map is equivariant for the diagonal $H$-action on $L\times G$ (by the given action on $L$, and left multiplication on $G$), and the $H$-action on $M\times S \times T^*G$ which is trivial on $M$, the given one on $S$, and induced by left multiplication on $T^*G$; moreover, by our assumption on the maps to $\mathfrak h^*$, it lives over the kernel of the moment map for $H$, and the existence of the Lagrangian correspondence \eqref{Lc0} follows from Hamiltonian reduction by $H$ and dimension counting. In favorable circumstances, one can argue (e.g., by studying tangent spaces) that this correspondence comes from an isomorphism between the symplectic spaces.

  \begin{remark}[Hamiltonian restriction] \label{SS2}\index{shifted symplectic geometry}
In the language of Remark \ref{fancy induction}, 
 a Hamiltonian space $S$ is encoded by the shifted Lagrangian $S/H \rightarrow \mathfrak{h}^*/H$ (i.e., a Lagrangian correspondence from a point to $\fh^*/H$),
 and Hamiltonian induction amounts to composing this map with the Lagrangian correspondence
 $ \fg^*/G \leftarrow \fg^*/H \rightarrow \fh^*/H$. Composing with this Lagrangian in the opposite direction gives an adjoint operation of Hamiltonian restriction, through which one can formulate the relevant Frobenius reciprocity above. 
\end{remark}

  \subsection{Whittaker induction}\label{Whittaker induction}
   In this section we introduce the operation of {\em Whittaker induction} of Hamiltonian spaces, associated to a homomorphism $H\times SL_2\to G$
 and a Hamiltonian $H$-space:
  $$ \mbox{(graded) Hamiltonian $H$-spaces} \longrightarrow \mbox{(graded) Hamiltonian $G$-spaces},$$
   reducing to ordinary Hamiltonian induction when the $SL_2$-homomorphism is trivial. 
   The class of Whittaker inductions of vectorial  representations of subgroups $H\subset G$ will
subsume the examples of \S~\ref{model examples}, and have many favorable properties;
 as $G$-spaces, they are simply vector bundles over homogeneous $G$-varieties -- see further
 Examples \ref{ttt}, \ref{ttt2}.

 \subsubsection{The notion of an $\mathfrak{sl}_2$-pair} \label{sl2pair}
 \index{$\mathfrak{sl}_2$ pair}
 \index{$\mathfrak{sl}_2$ triple}
 
For what follows, one can think of a fixed invariant identification 
\begin{equation}\label{ggstar}\mathfrak{g} \simeq \mathfrak{g}^*,
\end{equation}
 and a morphism 
\[\mathfrak{sl}_2 \to \mathfrak{g},\]
expressed via a triple $(h,e,f)$ of elements in $\mathfrak g$. 

More canonically,  the construction that we are about to describe depends only on an element $f \in \mathfrak g^*$, and a cocharacter \index{$\Diag$}  $\varpi: \Gm \to [G,G]$, such that $\varpi$ normalizes $f$, and the pair $(h=d\varpi(1), f)$ belongs to an $\mathfrak{sl}_2$-triple $(h, e, f)$ under some (equivalently, any) invariant identification \eqref{ggstar}. Such a pair $(\varpi, f)$ will be called an ``$\sl_2$-pair.''  Notice that the centralizer of the triple $(h,e,f)$ does not depend on the identification \eqref{ggstar}, thus, it is a subgroup (or sub-Lie algebra) associated with the $\sl_2$ pair. Throughout the discussion that follows, the subgroup $H\subset G$ is a subgroup of the centralizer of such an $\sl_2$-triple.
 
\subsubsection{The Hamiltonian space $(\mathfrak u/\mathfrak u_+)_f$} \label{uuplus}
 
 Let  $(\varpi, f)$ be an $\sl_2$-pair as above. When there is no danger of notational clash with a subgroup $H$, we will also denote the cocharacter $\varpi$ by
 $\lambda \mapsto \lambda^h$.

 Decompose
\begin{equation}\label{hdecomp} {\mathfrak g} = \mathfrak j \oplus  \overline{\mathfrak{u}} \oplus \mathfrak{u}^{0} \oplus \mathfrak{u},
\end{equation}
where $\mathfrak j$ is the centralizer of $\mathfrak{sl}_2$ and $\overline{\mathfrak{u}} \oplus \mathfrak{u}^{0} \oplus \mathfrak{u}$ is the sum of all nontrivial $\mathfrak{sl}_2$-subrepresentations, decomposed into the sum of negative, zero, and positive weight spaces for the left adjoint action of $h$;
thus, $f \in \bar{\mathfrak{u}}$. We denote by $\bar{U},U$ the associated unipotent subgroups. Observe that the $\Gm$-action 
 $\lambda^h$ normalizes $U$, i.e., we consider $\fu$ as a {\em graded} Lie algebra.  \index{graded Lie algebra}

 \index{$\mathfrak{u}_+$}
 Let   \begin{equation} \label{uplusdef} \mathfrak{u}_+ \subset \mathfrak{u}\end{equation}
be the sum of all $h$-eigenspaces of weight $\geq 2$.
Then $\mathfrak{u} = \mathfrak{u}_+$ exactly when all the weights of the adjoint $\SL_2$-action
on the Lie algebra are even, i.e., when $-1 \in \SL_2$ is central in $G$. 
This situation is somewhat simpler, and we  {\em recommend that the reader assume at first reading 
 that   all weights of $\mathfrak{u}$ are even, i.e., $\mathfrak{u}=\mathfrak{u}_+$.}

\begin{example} \label{ttt}
Suppose that $S$ is trivial and there are no odd weights. 
In this case, 
$$ \mbox{Whittaker induction of trivial space from $H$ to $G$} = T^*(HU\backslash G, \Psi),$$
  with the notation as in \S \ref{twisted cotangent bundle section}.
 Here,  the element $f$ defines an additive character $HU \rightarrow \Ga$, which gives rise to
a $\Ga$-bundle $\Psi$ over $X=HU \backslash G$. 
Equivalently, it is the Hamiltonian induction from $\pt_f$, the trivial $HU$-symplectic space with moment map $\text{pt} \to f \in (\mathfrak h+\mathfrak u)^*$. 
\end{example}
\index{$\pt_f$}

 Let   $\mathfrak{u}_+$ be as in \eqref{uplusdef}  and let  $U_+$ be the associated
unipotent group. 
We now treat $f$, restricted to $\mathfrak u_+$, as a Lie algebra homomorphism $\mathfrak u_{+} \to \mathfrak g_a$.
 The quotient $\mathfrak{u}/\mathfrak{u}_+ $ carries an $H$-invariant symplectic form
\begin{equation} \label{uquot symplectic}
(x,y) \in \mathfrak{u} \times \mathfrak{u} \mapsto \langle f, [x,y] \rangle.\end{equation}
 Indeed,  it descends to $\mathfrak{u}/\mathfrak{u}_+$ in both factors
 by weight arguments;
 also, the right-hand side equals $\langle [f,x], y \rangle$,
 and  if $x\ne 0$ has weight $1$, then $[f,x]$ is necessarily nonzero of weight $-1$ and therefore there is a $y\in \mathfrak u$ with $\left< [f,x], y \right>\ne 0$.  Now define
 $$ (\mathfrak{u}/\mathfrak{u}_+)_f = \mbox{ 
 $(\mathfrak{u}/\mathfrak{u}_+)$ considered as a 
 Hamiltonian $HU$-space }$$
  where
   \begin{itemize}
  \item 
the $H$-action is through the adjoint action;
\item  $U$ acts by translation 
  via $U/U_+ \simeq \fu/\fu_+$;
  \item the moment map on the $H$ factor arises from the structure of symplectic representation, i.e., through \eqref{stdvect};
  \item on $U$  
 we use the $f$-shifted moment map
  $$\xymatrix{\mathfrak{u}/\mathfrak{u}_+ \ar[r]^-{\sim}_{\eqref{uquot symplectic}}& (\mathfrak{u}/\mathfrak{u}_+)^*\ar[rr]^-{X \mapsto X+f}&& \mathfrak{u}^*},$$
  \end{itemize}

 \subsubsection{Whittaker induction} \label{ss:Whitinduction}
 For $S$ a Hamiltonian $H$-space we put
\begin{equation} \label{tildeSdef} \tilde{S} = S \times (\mathfrak{u}/\mathfrak{u}_+)_f, \end{equation}
which we consider as a Hamiltonian $HU$-space (with $U$ acting trivially on $S$). 
The Hamiltonian $G$-space giving the Whittaker induction will be, by definition,
the Hamiltonian induction
\[\textrm{h-ind}_{HU}^G \tilde{S}\]
  of $\tilde{S}$ from $HU$ to $G$ defined as in \S~\ref{symplecticinduction}.
Explicitly, by \eqref{cfp},
\begin{equation}\label{Whittinduction} \textrm{h-ind}_{HU}^G \tilde{S} = \Big(S \times (\mathfrak{u}/\mathfrak{u}_+)_f \Big) \times^{HU}_{(\mathfrak h + \mathfrak u)^*} (\mathfrak g^*\times G),\end{equation}
where, as before, $\mathfrak g^*\times G$ stands for $T^*G$ via \eqref{TGidentity}.

Note that the Whittaker induction from a symplectic $H$-vector space $S$ (with its natural, quadratic, moment map to $\mathfrak h^*$) contains a canonical basepoint, the point $(0, 0) \times (f, \mathrm{id}_G)$,
in the above presentation, whose moment image is $f$. 
Note also that, when we choose an invariant identification $\mathfrak g^*\simeq \mathfrak g$, identifying $\mathfrak u^*\simeq \bar{\mathfrak u}$ (opposite nilpotent subalgebra), the moment image of $(\mathfrak{u}/\mathfrak{u}_+)_f$ is $f+\mathfrak u_{-1}$.

\subsubsection{The $\GGm$-action on a Whittaker induction} \label{ggmm} 
  
 For $S$ a graded Hamiltonian $H$-space the Whittaker induction inherits a natural grading.
  But the definition of this grading is slightly more complicated than for symplectic induction of a graded space; the space $(\mathfrak{u}/\mathfrak{u}_+)_f$
is {\em not} naturally graded, because of the $f$-shift. 
 We will write a formula below but the reader might want to skip to the 
more conceptual reformulation outlined in \S~\ref{symplecticshear} and use the formula only in case of emergency.

Before we proceed to the formula we give an example:

 \begin{example} \label{ttt2}
As described in Example \ref{ttt} the Whittaker induction in the case
of $S$ trivial and even weights is
 $ M = T^*(X, \Psi)$
 arising from an affine bundle $\Psi \rightarrow X = HU \backslash G$.

 The  action of $\Gm$ via left multiplication by $\varpi$, 
  on $G$ descends to an action on 
 $\Psi \rightarrow X$, and thereby to an action of $\Gm$
 on $T^* \Psi$; but this action does not descend to $M$. 
 We modify 
   this action by multiplying it by the action on $T^* \Psi$ 
which scales by the squaring character along fiber. The resulting
$\Gm$ action on $T^* \Psi$ is now equivariant for the moment map
$T^* \Psi \rightarrow \mathfrak{g}_a^*$ and descends 
to $T^*(X, \Psi)$; this is the desired grading on $T^*(X, \Psi)$. 

Explicitly, the fiber over $T^*(X, \Psi)$ over the identity coset in $X$
is identified with the elements of $\mathfrak{g}^*$ that
restrict to $f \in \mathfrak{h}^* \oplus \mathfrak{u}^*$, i.e.
$f + (\mathfrak{h} \mathfrak{u})^{\perp}$. 
The action of $G$ then defines an isomorphism
\begin{equation} \label{TXPexplicit}  T^*(X, \Psi) \simeq 
HU \backslash  \left ( f + (\mathfrak{h} \mathfrak{u})^{\perp})  \times  G\right)\end{equation}
and the $\GGm$ action is given by left multiplication by $\varpi$
on the $G$ factor, and  by the composition of the left adjoint action of $\varpi$ and squaring on the first factor. 
 
\end{example}

In the general case, referring to 
 \eqref{Whittinduction}, we let $\GGm$ act as follows:
\begin{enumerate}
 \item by the given grading on $S$;
 \item by the tautological character (scalar action) on the symplectic vector space $(\mathfrak{u}/\mathfrak{u}_+)_f$;
 \item by the square character composed with the \emph{left} coadjoint action of $\varpi$ on $\mathfrak g^*$;
 \item by left multiplication by the left action of the cocharacter $\varpi$ on $G$ (i.e., $\lambda: g \mapsto \varpi(\lambda) g$). 
\end{enumerate}
In the next subsection we explain a more conceptual viewpoint on this action, which in particular implies that it scales the symplectic form appropriately.

  \subsubsection{Shearing} \label{symplecticshear}
 Suppose that $\varpi: \GGm \rightarrow \mathrm{Aut}(G)$
  is an action of $\GGm\simeq \Gm$ on $G$ by automorphisms, i.e., $G$ is a ``graded group'' -- this terminology seems most natural
  when $G$ is e.g., unipotent, but we will be applying it to more general affine groups, thinking of the grading of their coordinate rings. \footnote{In our applications the grading will be either on a unipotent group or {\em inner}, in fact: by left conjugation composed with a cocharacter $\varpi:\GGm\to G$.}
  It is most convenient to denote this here as a right action. We let $\GGm$ act on $\mathfrak g^*$ by the induced action, composed with dilation by the square character.
  Then:
  \begin{quote}
A {\em sheared} Hamiltonian $G$-space $M$ (relative to the grading on $G$)
will be a Hamiltonian $G$-space with $\GGm$-action compatible with the grading on $G$ and $\fgx$.
\end{quote}

In other words, to give a Hamiltonian $G$-space $M$ (with moment map $\mu$) the structure of a sheared space
means that we should give a $\GGm$-action on $M$ with the following properties for $x\in M, g \in G, \lambda \in \GGm$:
\begin{equation} \label{GGU} x \cdot g \cdot \lambda= x \cdot  \lambda \cdot g^{\varpi(\lambda)} \mbox{ and }
\mu(x \cdot \lambda) = \lambda^2 \mu(x)^{\varpi(\lambda)}. \end{equation} 
To avoid confusion, we emphasize that a \emph{graded} Hamiltonian $G$-space, defined in \S~\ref{sphvarintro}, is the same as sheared Hamiltonian $G$-space for the \emph{trivial} grading (i.e., $\GGm$-action) on $G$.

   Here are some examples:
   \begin{itemize}
   \item[(i)]  If $\varpi$ is trivial, this recovers the notion of a graded Hamiltonian space. 
   \item[(ii)]  If $M$ is a graded Hamiltonian $G$-space, and $\varpi: \Gm \rightarrow G$ is a cocharacter,
we can alter the $\GGm$-action by composing it with the (right) action of $\varpi$ on $M$, and this  gives a sheared Hamiltonian space, where $G$ is graded through the right inner action of $\varpi$. 
\item[(iii)]  
The point as a space under $\mathbb{G}_a$ (as a graded group with $\GGm$-action by the square character), but with moment map
image $1 \in \mathfrak{g}_a^*$, is a sheared Hamiltonian space. 
\item[(iv)]  Let $W$ be a symplectic vector space 
over $\FF$; then $W$ is a sheared Hamiltonian space under the Heisenberg group $W \ltimes \mathbb{G}_a$, 
which acts on $W$ through translation by the quotient $W$;  the $\GGm$-action is by scaling on  $W$, by  square scaling
on the center of the Heisenberg group, and the moment map is given by
$W \mapsto W^* \oplus \mathfrak{g}_a^*$ where the first coordinate is the identification $w\mapsto \omega(w,\bullet)$ and the second
coordinate is $1$.   
\item[(v)]  The Hamiltonian space $(\mathfrak{u}/\mathfrak{u}_+)_f$ defined in \S~\ref{uuplus} is a sheared Hamiltonian space for $U$, when $\GGm$ acts by left conjugacy on $U$ via the cocharacter $\varpi$ associated to the $\sl_2$-triple (= right conjugacy through $\varpi^{-1}$). The $\GGm$-action on $\mathfrak u/\mathfrak u_+$ is the induced left adjoint action, namely,  
scaling by the tautological character. Combined with the square action on $\mathfrak u^*$, this causes the $f$-shifted moment map defined in \S~\ref{uuplus} to be equivariant under the $\GGm$-action. Finally, this extends to the structure of sheared Hamiltonian space on $HU$, where $H$ is
graded trivially.
\end{itemize}

  Using this terminology,  we can describe the Whittaker induction process as follows: Fix the ``$\sl_2$-pair'' $(\varpi,f)$ or, if desired, an isomorphism $\mathfrak g \simeq \mathfrak g^*$ and an $\sl_2$-triple, and fix a subgroup $H\subset G$ of the centralizer of the $\sl_2$-triple. We use the action induced by right conjugation action via $\varpi$ to define shearing below. Whittaker induction is the process of assigning to a graded Hamiltonian $H$-space a (non-sheared) graded Hamiltonian $G$-space, via
\begin{multline} \label{Whit induction in stages}
  \mbox{graded Hamiltonian $H$-spaces} \stackrel{\times (\mathfrak{u}/\mathfrak{u}_+)_f}{\longrightarrow}
\mbox{sheared Hamiltonian $HU$-spaces } \\
 \stackrel{\text{h-ind}}{\rightarrow} \mbox{sheared Hamiltonian $G$-spaces } 
\rightarrow \mbox{graded Hamiltonian $G$-spaces}.
\end{multline}
 For the last arrow,  we use the fact  -- inverting (ii) above -- that for the $\varpi$-grading on $G$,   
  any sheared Hamiltonian space arises from a usual one by twisting the $\GGm$-action through $\varpi$. 
 
 This point of view on the $\GGm$-action is more conceptual, and parallel constructions will be very useful later on in describing $L$-sheaves and spectral quantizations of Whittaker inductions.

  \subsubsection{The vector bundle structure of a Whittaker induction}
  \label{Slodowy}
We shall now prove that, ignoring the symplectic structure,  
any Whittaker-induced \emph{linear} space  (i.e., for $S$ a symplectic $H$-representation, equipped with the scaling action of $\GGm$) is simply a vector bundle over the homogeneous space $H \backslash G$.
More precisely, as a $G$-space, we may identify
\begin{equation} \label{fund_id} M \simeq  V \times^{H} G, \ \ V := \left[ S \oplus (\mathfrak h^\perp \cap \mathfrak g^{*,e}) \right].\end{equation}
Here $\mathfrak{g}^{*,e}$ is the  kernel of the action of $e$ on $\mathfrak{g}^*$, which under an isomorphism $\mathfrak g^*\simeq \mathfrak g$ can be identified with the centralizer Lie algebra $\mathfrak g_e$.  \index{$\mathfrak{g}_e$} 

Moreover,
the resulting isomorphism \eqref{fund_id} respects $\GGm$ actions,
where $\GGm$ acts on $G$ through left multiplication by $\varpi$, and where 
the action on $V$ is as follows:
\begin{itemize}
\item $\GGm\simeq \Gm$ acts by linear scaling on $S$,  and
\item  it acts with weight $2+t$ on the weight-$t$ component of $\mathfrak{g}_e$ under the left adjoint action of $\varpi$ (equivalently, the right adjoint action of $\varpi^{-1}$).
\end{itemize} 

This isomorphism is in some ways a little artificial -- $M$ is more canonically an affine bundle over $H \backslash G$,
and some choices are required to make the above identification --
but in any case it is very convenient for our purposes. It is in the form \eqref{fund_id}
that the space $M$ previously appeared in the theory of automorphic forms.  Moreover, the isomorphism \eqref{fund_id}
has the following perhaps surprising consequence: 
\begin{lemma}\label{lemma:affine} 
In the setting of \S~\ref{Slodowy}, if $H$ is reductive, then $M$ is affine.
\end{lemma}

\begin{proof} 
 It is enough to show that the natural map from the stack quotient of $V\times G$ by $H$ to the invariant-theoretic quotient $(V\times G)\sslash H$ is an isomorphism. Both live over $H\backslash G$ (which is affine since $H$ is reductive), and by reductivity the restriction map $F[V\times G]^H\to F[V\times H]^H$ is surjective. In other words, the map is an isomorphism over the fibers of $H\backslash G$, and by homogeneity it is an isomorphism everywhere.
\end{proof}

From \eqref{cfp2} the Whittaker induction is identified with 
\begin{equation} \label{cfp3}  \tilde{S} \times^{HU}_{(\mathfrak{h}+ \mathfrak{u})^*}  (\mathfrak{g}^* \times G),
\mbox{ with } \tilde{S} := S \times (\mathfrak{u}/\mathfrak{u}_+)_f .
\end{equation}

Projecting to the $S$- and $\mathfrak{g}^*$-coordinates
  defines an $HU$-equivariant isomorphism
\begin{equation} \label{cfp4}
   S \times (\mathfrak{u}/\mathfrak{u}_+)_f  \times_{(\mathfrak{h}+ \mathfrak{u})^*} \mathfrak{g}^*
   \simeq \{s \in S,  t\in f+\mathfrak{u}_+^{\perp}: \mu(s) = t|_{\mathfrak{h}} \} .\end{equation}
   Next, by the theory of Slodowy slices \cite[Lemma 2.1]{Gan-Ginzburg}, the action of $U$ on $f+\mathfrak u_+^\perp$ is free, and admits a transversal section equal to $f+ \mathfrak g_e$ (considered as an affine subspace of $\mathfrak{g}^*$ under an identification $\mathfrak g^*\simeq \mathfrak g$ -- this subspace only depends on the $\sl_2$ pair $(\varpi,f)$). Note that this section is invariant under the action of the group $H$, since this is contained in the centralizer of $\mathfrak{sl}_2$,   and the action of $U$ does not affect the projection $f+\mathfrak u_+^\perp \to \mathfrak h^*$, by
 a simple weight argument. 
 Correspondingly, the map $(s,z) \mapsto (s, f+z)$,
 together with the identification $\mathfrak{g} \simeq \mathfrak{g}^*$, gives rise to an
 $H$-equivariant identification
 $$ \left[ S \times_{\mathfrak{h}^*} \mathfrak{g}_e \right] \times U \rightarrow  \{s \in S, t\in f+\mathfrak{u}_+^{\perp}: \mu(s) = t|_{\mathfrak{h}} \}$$

In particular, \eqref{cfp3} is identified with 
\begin{equation}\label{cfp5}\left[ S \times_{\mathfrak{h}^*} \mathfrak{g}_e\right] \times^H G.
\end{equation}
More canonically, one should replace $\mathfrak{g}_e$ by $\mathfrak g^{*,e}$, and then the isomorphism \eqref{cfp5} depends only on the data of $H$ and the $\sl_2$-pair $(\varpi, f)$. 
Now, fixing
an $H \times \SL_2$-equivariant splitting of $\mathfrak{g}^* \rightarrow \mathfrak{h}^*$, we get an identification
of $S \times_{\mathfrak{h}^*} \mathfrak{g}_e$ with the vector space $V$
appearing in \eqref{fund_id}.   Our description of the $\GGm$ action follows from 
  \S~\ref{ggmm} (see in particular after \eqref{TXPexplicit}). 
  
  Note, in particular, that when we embed \begin{equation} \label{Gmprime} \Gm \stackrel{\varpi^{-1} \times \Id}\hookrightarrow \check G \times \GGm\end{equation}    this $\Gm$ fixes the identity coset of  $H\backslash G$, and provides a grading on $V$; we will sometimes denote this copy of $\Gm$ by $\GGm'$ \index{$\GGm'$}.

\begin{example}\label{Slodowy slice example}
Consider the case when $H$ and $S$ are both trivial, and all weights of the $\SL_2$ on the Lie algebra are even.  Then the Whittaker induction  $M= T^*(U\backslash G,\Psi)$ is the generalized-Whittaker twisted cotangent bundle,
and our above considerations reduce to the isomorphism
 \[M\simeq \fg_e \times G.\]

 Moreover, if we choose $H$ commuting with $\SL_2$, this $M$ can also be identified 
 with the Whittaker induction from $H$ of $T^* H$, and correspondingly acquires an 
$H$-action and moment map. The $H$-action is, explicitly, the action arising from
 left multiplication on $U \backslash G$, and the moment map   is identified with the projection $\fg_e\to \fh^*$.

 In other words (a point of view that will be useful later)  we may identify $M/(H\times G) \simeq \fg_e/H$ as spaces over $\fh^*/H \times pt/G$. 
\end{example}

 \subsection{Hyperspherical Hamiltonian spaces} \label{ssconditions} 
 In this section we shall describe a class of graded Hamiltonian spaces, which we will call ``hyperspherical.''
 This class contains the  
  cotangent spaces of \emph{smooth, affine} spherical varieties
  satisfying a certain connectedness condition on stabilizers (see Proposition \ref{hypersphericalspherical}), and seems suitable for our conjectural duality. The most important property is the coisotropic property, which in representation theory is closely related to the ``multiplicity one'' property, and plays an important role in the theory of automorphic forms.
 
  As a matter of notation, we will often refer to a ``hyperspherical $G$-space $M$.''
In other words, we do not explicitly include the Hamiltonian structure and grading  in the notation,
even though it is understood to be part of the structure.
  
 \subsubsection{The conditions} \label{condALL}
 Consider graded  irreducible (and smooth, by definition) Hamiltonian $G$-varieties $M$ satisfying the following conditions (which we discuss in detail in the following sections); these spaces will be referred to as ``hyperspherical'' in this paper: 
\begin{enumerate}
\item \label{condaffine}  $M$  is affine;
 \item \label{condcoisotropic} the field $\FF(M)^G$ of 
$G$-invariant rational functions on $M$ is commutative with respect to the Poisson bracket (i.e., $M$ is ``coisotropic'');
 \item \label{condzero} the moment map image has nonempty intersection with the nilcone in $\mathfrak g^*$;
 \item \label{condtorsion} the stabilizer (in $G$) of a generic point of $M$ is connected; 
\item \label{condGm} the   $\GGm$-action is ``neutral'' (to be defined in \S~\ref{Sneutral}). \index{neutral $\GGm$}

 \end{enumerate}

Neutrality will be introduced after we have seen some consequences of the other conditions: the rest of the conditions imply the existence of a unique closed $G \times \GGm$-orbit $M_0 \subset M$,
whose moment map image is a nilpotent orbit, and neutrality
essentially means that the $\GGm$-action near $M_0$
is determined by an $\mathfrak{sl}_2$-triple associated to that nilpotent orbit. It is the most subtle condition, and perhaps least satisfying. 
The neutral $\GGm$-action is geometrically natural and quite rigid. It has one main drawback, 
namely, it does not not always satisfy the even parity condition
  described in \S~\ref{aims}, i.e.,  that the action of $\GGm \times G$ factors through the extended dual group ${}^CG$. 
  For this reason, we will sometimes modify it for specific purposes, but it is very convenient as a general definition. 

Any Whittaker-induced linear symplectic space will satisfy \eqref{condaffine}
  \eqref{condzero}, and \eqref{condGm}.
What we will show below implies, more precisely, that if we add \eqref{condcoisotropic}
the converse is also true: any such space must in fact arise from the construction of \S~\ref{ss:Whitinduction}.

\subsubsection{Some consequences of the conditions on $M$}
 \index{$\mathfrak c^*$}
 \index{$\mathfrak c_M^*$}
 
  Passing to the invariant-theoretic quotient $$\mathfrak g^*\to \mathfrak c^*:= \mathfrak g^*\sslash G = \mathfrak a^*\sslash W,$$  we obtain the \emph{invariant moment map} $\mu_G:M\to \mathfrak c^*$. 
  Following Knop, Losev \cite{Losev} introduces a Stein factorization of the invariant moment map
  \begin{equation}\label{Stein} M \xrightarrow{\tilde\mu_G} \mathfrak c_M^* \to \mathfrak c^*,
  \end{equation}
  such that the first map is dominant with connected generic fiber, and the second map is finite. The space $\mathfrak c_M^*$ is defined as the spectrum of the integral closure (normalization) of the image of $\FF[\mathfrak c^*]$ inside of the function field $\FF(M)$.

The condition of $M$ being \emph{coisotropic} is defined by either of the equivalent criteria of the following proposition:
\begin{proposition}\label{coisotropic}
 The following are equivalent:
  \begin{enumerate}
   \item[(i)]  \label{invariantscommute} the field $\FF(M)^G$ is commutative with respect to the Poisson bracket; 
   \item[(ii)] \label{coisotropicorbit} the generic $G$-orbit on $M$ is coisotropic;
   \item[(iii)]\ \label{opengenericfiber} the generic fiber of $\tilde \mu_G$ contains an open $G$-orbit.
  \end{enumerate}
\end{proposition}

\begin{proof}
The equivalence (i) $\iff$ (ii) is essentially \cite[II.3, Proposition 5]{Vinberg-commutative}, except that there it was stated in the differentiable setting. We repeat the argument, with the details necessary for the algebraic setting: 

It is known that, for any action of an algebraic group on an irreducible variety $M$, a finite number $f_1, \dots, f_r$ of elements of the field $\FF(M)^G$ of rational invariants separate generic orbits \cite[Theorem 2]{Rosenlicht}.  
In particular, there is an open dense $G$-stable subset $M'$ where the $f_i$'s are defined, and the fibers of the resulting morphism $M' \to \mathbb A^r$ are $G$-orbits (if nonempty). If $\mathfrak c'\subset \mathbb A^r$ is the spectrum of the subalgebra spanned by the $f_i$'s, we may, by further restricting $M'$ and $\mathfrak c'$, assume that the morphism $M'\to \mathfrak c'$ is smooth and surjective. In particular, the differentials of the $f_i$'s span the orthogonal complement to the image $\mathfrak g_x$ of $\mathfrak g^*\to T_x M'$, at every point $x\in M'$. It follows easily from the definitions, now, that these functions Poisson-commute iff $\mathfrak g_x$ is a coisotropic subspace of $T_xM'$, for all $x\in M'$. The equivalence of the first two statements follows. 

To prove their equivalence with the third statement (iii), we use \cite[Theorem 1.2.4]{Losev}, which states that $\FF(M)^G$ is the fraction field of $\FF[M]^G$; equivalently, the generic fiber of $M\to M\sslash G$ contains a dense $G$-orbit. Moreover, \cite[Proposition 5.9.1]{Losev} identifies regular functions in the image of $\tilde\mu_G$ as the intersection of $\FF[M]^G$ with the Poisson center of $\FF(M)^G$. Hence, if $\FF(M)^G$ is Poisson-commutative, we have a dense embedding $M\sslash G\hookrightarrow \mathfrak c_M^*$, and the generic fiber over $\mathfrak c_M^*$ contains an open $G$-orbit. Vice versa, if the generic fiber of $\tilde\mu_G$ contains an open $G$-orbit, the fact that the elements of $\FF[\mathfrak c_M^*]$ Poisson-commute in $\FF(M)$ implies, by the same argument as before, that the generic $G$-orbit is coisotropic. 
\end{proof}

 To emphasize some corollaries of the preceding proof, under the equivalent conditions of  Proposition \ref{coisotropic}, \cite[Proposition 5.9.1]{Losev} states that the invariant-theoretic quotient $M\sslash G$ is equal to the image of $\tilde\mu_G$, and  \cite[Theorem 1.2.4]{Losev} identifies $\FF(M)^G$ with the quotient field of $\FF[M]^G = \FF[\text{Im}\tilde\mu_G]$.

  \medskip

 Next, condition \eqref{condzero} that the moment map image meet the nilcone is equivalent to asserting that $0 \in \mathfrak a^*\sslash W$ is in the image of the invariant moment map $\mu_G$.
This has several important consequences: 
\begin{itemize}
\item[(i)]  $\mathfrak{c}_M^*$ contains a
 unique point (also to be denoted by $0$) over the point $0\in \mathfrak c^*$.
 Indeed, the  $\GGm$-action lifts to $\mathfrak{c}_M^*$ by functoriality, and the points over $0$ will be the closed $\GGm$-orbits. Hence, they are the points in the spectrum of the $0$-th graded piece of $\FF[\mathfrak c_M^*]$. By finiteness of the morphism $\mathfrak c_M^*\to \mathfrak c$, the latter is an integral extension of $\FF=$ the $0$-th graded piece of $\FF[\mathfrak c^*]$, and therefore equal to $\FF$. 
\item[(ii)] The image of $\tilde{\mu}_G$ is all of $\mathfrak{c}_M^*$.   Indeed this image is open by \cite[Theorem 1.2.2]{Losev}, so  its complement is a closed set that is stable by $\GGm$; since it cannot contain the unique point over $0\in \mathfrak c^*$, it must be empty.   
   \item[(iii)]
      There is a unique closed $G\times \GGm$-orbit $M_0\subset M$. Indeed, this is equivalent to $\FF[M]^{G\times \GGm} = \FF$, which follows from $\FF[M]^{G\times \GGm}  = (\FF[M]^G)^{\GGm} = \FF[\mathfrak c_M^*]^{\GGm} = \FF$.  

      Note a special case of this is when  the $\GGm$-action on $M$ is contracting, in which case the image of $M_0$ is the origin of $\mathfrak{g}^*$.

\item[(iv)] The closed $G \times \GGm$-orbit $M_0$ is  in fact a single $G$-orbit. This is because all of them map to $0$ in $\mathfrak{c}_M^* = M\sslash G$, whose preimage contains a unique closed $G$-orbit; but, by $\GGm$-transitivity, if one of those orbits is closed, all of them are. 
     \end{itemize}
    
 \medskip
 
 Finally, condition \ref{condtorsion} seems to be of technical nature, but is essential for the form of the conjectures that we present here. For example, in the spherical case $M=T^*X$, it is essentially equivalent to the statement that the Langlands dual of the ``universal Cartan'' of $X$ embeds into the Langlands dual of the Cartan of $G$, see \S~\ref{sec:invariants}. Without this condition, we do not fully understand, even conjecturally, the Langlands dual picture, but any such duality would involve different constructions, such as covering groups like the ones introduced in \cite[\S~3.2]{SV}, stacks,  or, in more complicated cases,  quantum groups such as in \cite{MT}, see \S~\ref{OtherExamples}. 
  
 \subsubsection{Neutrality} \label{Sneutral}
    Choose $x \in M_0$ (the closed $G\times \GGm$-orbit), and let $f$ be its (nilpotent) image in $\mathfrak g^*$.  
 (It is easy to see that the definition of neutrality that follows will be independent of choice of $x$.)
Since $M_0$ is affine, the stabilizer $H:=G_{x}$ of $x$ inside $G$ is reductive. 

\begin{remark}\label{Hconnected}
 A priori, the group $H$ may be disconnected. We do not know of such examples, and we expect that $H$ will always be connected, but we only have a proof in the polarized case, see Proposition \ref{isspherical}.
\end{remark}

 Now, $H$ fixes $f$ under the coadjoint action,
and the moment map on $M_0 \simeq H \backslash G$ is given by $Hg \mapsto f^g$. 
The $\GGm$-action on $M_0$ commutes with $G$ and therefore
is given by left multiplication by a cocharacter $ \varpi: \Gm \rightarrow N(H)/H$, such that $f^{\varpi(\lambda)} = \lambda^2 f$.

This cocharacter can also be thought of as a cocharacter 
\[ \varpi: \Gm \to Z_G(H)/Z(H) \]
(where $Z_G(H)$ is the centralizer of $H$ in $G$, and $Z(H)$ is the center of $H$). Indeed, both $H$ and $N(H)$ are reductive groups, hence at the level of Lie algebras we have $\mathfrak n(\mathfrak h) = \mathfrak z_{\mathfrak g}(\mathfrak h) \times^{\mathfrak z(\mathfrak h)} \mathfrak h$; therefore, the embedding $Z_G(H)/Z(H) \to N(H)/H$ is an isomorphism on identity components.

\begin{definition} \label{condneutral}
With notation as above, the $\GGm$-action on $M$ will be called {\em neutral}
when both of the following conditions are satisfied:

\begin{quote}
\begin{itemize}
\item[(i)] The pair $(\varpi, f)$ lifts to an $\sl_2$-pair for $G$ (\S~\ref{sl2pair}), that is, under an invariant identification $\mathfrak g \simeq \mathfrak g^*$, the cocharacter $\varpi$ lifts to a cocharacter of the form $\lambda\mapsto \lambda^h$, for an $\sl_2$-triple $(h,e,f)$. 
\item[(ii)]  (i) implies that the action of $(\lambda^{-h},\lambda)\in G\times \GGm$ stabilizes $x$; 
write \begin{equation}\label{varpixdef}
\varpi_{x}: \lambda \mapsto (\lambda^{-h}, \lambda) \in G \times \GGm \end{equation} for this one-parameter subgroup. Then 
 $\varpi_{x}$ acts by the identity cocharacter (i.e., by simple scaling) on the fiber $S$ of the   symplectic normal bundle 
(cf. \S~\ref{symplecticinduction})
 to the orbit $M_0\subset M$.\footnote{Since the moment map for the fibers of the symplectic normal bundle is quadratic, this requirement is compatible with our condition that $\GGm$ act on $\mathfrak g^*$ by squares. }

\end{itemize}
\end{quote}

\end{definition}

Let us observe that:
\begin{itemize} 
\item  The lift $h$ of $\varpi$ has image in $Z_G(H)$. Indeed, otherwise, by the theory of $\sl_2$-modules, $\ad(f)(\mathfrak h)\ne 0$; however, the intersection of $\ad(f)(\mathfrak g)$ with the centralizer of $f$ is normal in the latter, and nilpotent \cite[Theorem 3.6]{KostantTDS}, while $H$ does not contain a normal unipotent subgroup, a contradition.
\item Since $H$ centralizes both $f$ and $h$, it commutes with the corresponding $\mathfrak{sl}_2$-triple (obtained by an identification $\mathfrak g \simeq \mathfrak g^*$) in full. 
\item The $\mathfrak{sl}_2$-triple occurring here is {\em unique}.  Indeed, for a given $f$,
any other choice necessarily has the form $h' =  h + z$ 
where (by definition) $z$ lies in the center of $H$.
By  \cite[Theorem 3.6]{KostantTDS}, 
however, $z$ necessarily has negative weight under $h$,
contradicting the fact that $h$ centralizes $H$. 
\item The unique $\mathfrak{sl}_2$ will sometimes be called the {\em Arthur-$\mathfrak{sl}_2$}
attached to $M$, for reasons motivated by the role it will play when $M$ is placed on the  spectral side of the Langlands correspondence. 
\end{itemize}

In other words, we have extracted from $M$ a commuting pair $H \times \mathfrak{sl}_2$ in $G$,
as well as a symplectic $H$-vector space $S$; we will prove in Theorem \ref{thm:structure}
that, in fact, $M$ can be identified with  the Whittaker-induction of $S$. 
  
\subsection{The structure theorem} \label{ss:structuretheorem}
Let $M$ be a Hamiltonian $G \times \GGm$-space satisfying the conditions of \S~\ref{condALL}.   
Recall that it admits a unique $G \times \GGm$-closed orbit $M_0$. 
Fix a point $x\in M_0$ as above, with stabilizer $G_x=H$ and image $f$ under the moment map, and
recall from  the discussion after Definition
\ref{condneutral}
that its image $f\in \mathfrak g^*$ (under the moment map) belongs to a unique $\sl_2$-pair $(\varpi, f)$ (or, after fixing $\mathfrak g \simeq \mathfrak g^*$, an $\sl_2$-triple $(h,e,f)$),  commuting with $H$ and describing the action of $\GGm$ near $M_0$. 

\begin{theorem}
\label{thm:structure}
 Let $M$ be a Hamiltonian $G$-space satisfying the conditions of \S~\ref{condALL}. 
  Let $x\in M_0$ be a point in the closed $G\times \GGm$-orbit $M_0$,  with stabilizer $H$, and let $S$ be the fiber of the symplectic normal bundle to $M_0$ at $x$. Then, there is a unique  $G \times \GGm$-equivariant isomorphism of Hamiltonian $G$-spaces
 \begin{equation}\label{eq:slice} M\simeq 
 \textrm{Whittaker induction of $S$ from $(H, \mathfrak{sl}_2)$}
 \end{equation}
 which carries $x$ to the basepoint of the Whittaker induction (see discussion after \eqref{Whittinduction}) and 
induces there the identity on symplectic normal bundles.
 \end{theorem}

\begin{remark}
 The inducing space $\tilde S$ of \eqref{tildeSdef} is coisotropic for $H$, and in fact also for a smaller subgroup -- see Proposition \ref{Scoisotropic} below.
\end{remark}

The idea of the proof is as follows. Recall, from the discussion of \S~\ref{recognizinginduction}, that one can ``recognize'' $M$ as a Hamiltonian induction from $HU$ by producing
an $HU$-stable Lagrangian correspondence between  $M$
and a Hamiltonian $HU$-space.   Fix $x \in M_0$
and let $\Gm$ act by \eqref{varpixdef}; the Hamiltonian $HU$-space $\tilde S$ of \eqref{tildeSdef}
will be the weight-one subspace 
of $T_x M_0$; the Lagrangian  correspondence will have image 
  the set of points that  contract to $x$ under 
\eqref{varpixdef}.

\begin{proof}
We use the $\mathfrak{sl}_2$-pair 
to decompose
$\mathfrak{g}$ as in \eqref{hdecomp} and use notation as described there. 
Let $\Gm$ act by the cocharacter $$\varpi_{x}: (\lambda^{-h},\lambda)\in G\times\GGm$$ on $M$, so that it stabilizes $x$, by condition \eqref{condGm}.

Let $M_+\subset M$ be the subscheme of points that this $\Gm$-action contracts to $x$; as a functor, for a test scheme $T$, $M_+(T)$ is the set of $\Gm$-equivariant maps
$m: \mathbb A^1 \times T \to M$
such that $m(\{0\}\times T)= x$. In particular, there is a morphism
\[m: \mathbb{G}_a \times M_+ \rightarrow M,\]
classified by the identity morphism of $M_+$. 
 By the theorem of Bialynicki-Birula \cite[Theorem 4.1]{BB}, $M_+$ is a smooth scheme, whose tangent space at $x$ is identified with the sum of positive weight spaces for $\Gm$ on $T_x M$. Moreover, $M_+$ is (non-canonically) $\Gm$-equivariantly isomorphic to this tangent space. 
 $M_+$ is fixed under $HU$ because
 $H$ fixes $x$, and $\Gm$ contracts $U$; that is to say, for $z \in M, \lambda \in \Gm, u \in U$ we have
\begin{equation} \label{roofius} z \cdot u \cdot \varpi_{x}(\lambda) = z \cdot \varpi_{x}(\lambda) \cdot [\lambda^h u \lambda^{-h}].\end{equation}

(For example,  in the setting of Example \ref{ttt2},
if we take   $x$  to be the basepoint $f \times \mathrm{id}_G$, then
$M_+$ is identified with a single cotangent fiber of $T^*(X, \Psi)$.)

  Repeating the same considerations for the coadjoint representation $\mathfrak g^*$,
  where $\Gm$ is now acting through the product of the scaling square action
  and the left coadjoint action of $\lambda^h$, 
   we see that the map $m$ lives over a correspondingly defined 
\[ \Ga \times (f+ \mathfrak g^*_{\ge -1}) \to (f+ \mathfrak g_{\ge -1}^*).\]
Here, $\mathfrak g_{\ge -1}^*$ is the sum of weight spaces with weights $\ge -1$ for the left coadjoint action of $h$, so that $(f+ \mathfrak g_{\ge -1}^*)$ is the subset of points contracting to $f$ under the $\Gm$-action. 

Let us study the tangent space $T_xM$. It contains $T_x M_0$,
which is the tangent space to the $G$-orbit through $x$ and so identified with $\mathfrak{g}/\mathfrak{h}$. The orthogonal complement to this 
orbit is the kernel of the derivative of the moment map at $x$. Recalling the convention $d\left<Z,\mu\right> = \omega(Z,\bullet)$ for the moment map, 
we see that $T_x M_0 \cap (T_x M_0)^\perp = \mathfrak{g}_f/\mathfrak{h}$
(with $\mathfrak{g}_f$ the stabilizer of $f \in \mathfrak{g}^*$). We get
a filtration 
\begin{equation} \label{V1V2} 0 \subset \underbrace{ T_x M_0  \cap (T_x M_0)^\perp}_{V_1 = \mathfrak{g}_f/\mathfrak{h}}  \subset  \underbrace{ \mathfrak{g}/\mathfrak{h} }_{T_x M_0} \subset \underbrace{ T_x M_0 + (T_x M_0)^\perp}_{V_2=V_1^{\perp}} \subset T_x M, \end{equation}
about which we know the following:
\begin{enumerate}
 \item The orbit map gives an injection $\mathfrak{g}/\mathfrak{h} \hookrightarrow T_x M$,
 and the restriction of the symplectic form to $\mathfrak g/\mathfrak h$ is given by 
 \[ \omega(Z_1, Z_2) = \left< \text{ad}^*(Z_1)(f) , Z_2\right> = - \left< f,[Z_1,Z_2]\right>.\]  
 
  \item The quotient $ V_2/V_1$ is symplectic,  
and $\mathfrak{g}/\mathfrak{g}_f$, endowed with the symplectic form above, injects into it.
Let $S$ be the orthogonal complement of $\mathfrak{g}/\mathfrak{g}_f$
inside $V_2/V_1$; hence, $S$ is a symplectic vector space, and  
 \[ V_2/V_1 = S \oplus \mathfrak g/\mathfrak g_f,\]
an isomorphism of symplectic vector spaces. 

 \item  The $\Gm$-action via $\varpi_{x}$ on 
  the embedded $\mathfrak{g}/\mathfrak{h}$
  is given by $\lambda \mapsto \mathrm{Ad}(\lambda^h)$ (left adjoint action),
  and the restriction of this to $V_1=\mathfrak{g}_f/\mathfrak{h}$
  has weights $\leq 0$. 
  Since $\GGm$ acts on the symplectic form through squaring,
  the weights of $\varpi_{x}$ on $V_1$ are $\le 0$, and its weights on its dual, $ T_x M/V_2$, are $\ge 2$. 
 
 \item The weights of $\Gm$ on $T_x M/T_x M_0$ are, by the discussion above, all $\geq 1$,
 since they are all either weights on $S$ or on $T_X M/V_2$.

\end{enumerate}

It follows that the weight-$1$ subspace of $T_x M_+$ can be identified with the weight-$1$
subspace of $V_2/V_1$, which in turn is identified with the symplectic space $$\tilde{S} := S\oplus  \mathfrak{u}/\mathfrak{u}_+$$

  \begin{quote}
 {\em Claim:} There is a unique $HU\rtimes \Gm$-equivariant morphism  \begin{equation}\label{eq:mapMplus} \Lambda: M_+ \to  \tilde{S} \end{equation}
(where $\Gm$ is embedded in $G\times \GGm$ by $\varpi_x$), sending $x$ to $0$ and such that the differential at $x$ induces the natural projection $T_x M\to \tilde S$
to the weight one subspace. It has the properties that the restriction of the symplectic form to $M_+$ is obtained by pullback, and 
is compatible with moment maps, in the sense that the diagram 
\begin{equation} \label{rodacet}
\xymatrix{
M_+ \ar[r]^\Lambda   \ar[d]^{\mu} &  \tilde{S} \ar[d] \\
\mathfrak{g}^*  \ar[r] &  (\mathfrak{h} +   \mathfrak{u})^*
}
\end{equation}
commutes; the right vertical arrow here is the moment map for $HU$ acting on $\tilde{S}$,
defined as in \eqref{tildeSdef}. 
The induced map
\begin{equation} \label{mpp} M_+ \rightarrow \tilde{S} \times_{(\mathfrak h + \mathfrak u)^*} \mathfrak g^*
\end{equation}
is an isomorphism. 

 \end{quote}

{\em Proof of claim:} 
Set $\Lambda: M_+ \rightarrow T_x M_+$ to be the partial differential of  $m$ at $0\in \Ga$,
\[ d_{\Ga}m: \mathfrak g_a \times M_+ \to T_x M_+\]
evaluated at $1\in \mathfrak g_a$. Clearly, $\Lambda$ takes image in the weight-1 space $\tilde S:= S\oplus \mathfrak u_1$. Then the differential $d\Lambda: TM_+ \rightarrow \tilde S \subset T_x M_+$ can also be computed as the limit
$ \lim_{t \to 0} \lambda(t)$
where $\lambda(t): TM_+ \rightarrow TM_+$  is induced by the $\varpi_{x}$-action on $M_+$ multiplied by the \emph{inverse} of the scaling action on the fibers -- this can easily be inferred from the existence of a (noncanonical) $\Gm$-equivariant isomorphism $M_+ \simeq T_x M_+$, guaranteed by the theorem of Bialynicki-Birula. 
 Clearly, $\lambda(t)$ preserves symplectic forms, so $\Lambda$ also preserves symplectic pairings. 
Uniqueness of $\Lambda$ follows from the fact that the coordinate ring of $M_+$ is graded by the $\Gm$-action, with the $1$-graded piece identified with the dual vector space to $\tilde S$.

The map $\Lambda$ is $H$-equivariant because   the $H$-action on $M_+$ commutes with the $\Gm$-action.
 It is also $U$-equivariant with reference to the natural action of $U$
on $\tilde{S}$ -- trivial on $S$, and translation on $\mathfrak{u}_1$ via $U \rightarrow \mathfrak{u}/\mathfrak{u}_+ \simeq \mathfrak{u}_1$; 
this follows from \eqref{roofius}. Since any two moment maps for the action of $HU$ differ by translation, to verify that the diagram \eqref{rodacet} commutes it is enough to show that they agree at a single point. Taking $x$ to be that point, we see that both routes evaluate to the image of $f$ in $(\mathfrak{h}+ \mathfrak{u})^*$. 
 
The induced map \eqref{mpp} is an isomorphism because its differential at $x$ induces an isomorphism
 (because there are  compatible
  $\Gm$-actions on both sides with unique fixed point
$x$ and its image).  In more detail, given the fact that the moment map $S\to \mathfrak h^*$ is quadratic, the tangent space of $\tilde{S} \times_{(\mathfrak h + \mathfrak u)^*} \mathfrak g^*$ at the image of $x$ is naturally identified with $\tilde{S} \oplus (\mathfrak{h} \mathfrak{u})^{\perp}$; hence, the differential at $x$ is a map
\[T_x M_+ \rightarrow \tilde{S} \oplus (\mathfrak{h} \mathfrak{u})^{\perp}.\]
The projection to $\tilde S$ is the natural projection to the weight-$1$ subspace, whose kernel 
is the sum $(T_x M_+)^{\geq 2}$ of weight spaces on $T_x M_+$
with weight $\geq 2$, or equivalently the  corresponding definition $(T_x M)^{\geq 2}$ for $M$. 
The differential of the moment map gives $T_x M \rightarrow \mathfrak{g}^*$;
the image of this map is $\mathfrak{h}^{\perp}$ and its 
kernel is $T_x M_0^{\perp}$.
The weight $\geq 2$ subspace of this kernel is however trivial,
because it is a subspace of the weight $\geq 2$ subspace of
of $(T_x M/ T_x M_0)^{*}$.

 Therefore, $(T_x M_+)^{\geq 2}$ is 
 identified  by means of the moment map with the weight $\geq 2$ subspace of $\mathfrak{h}^{\perp}$,
which is the same as $(\mathfrak{h} + \mathfrak{u})^{\perp}$.  This concludes the proof of the claim. 
 
 \medskip

We are now in the situation discussed in \eqref{Lc0} and the subsequent discussion, namely, we have a Lagrangian
correspondence as in \eqref{Lc0}:
\begin{equation} \label{Lc1}
  \xymatrix{
 M_+  \ar[d] \ar[r]& M \\
 \tilde{S}^\circ  
  }
\end{equation}
where as before the superscript $\circ$ means that the symplectic form has been negated. 
Then, as in \eqref{inducedLL},  $M_+ \times^{HU} G$
gives a Lagrangian correspondence between $M$
and the Hamiltonian induction of $\tilde{S}$, 
and the induced map to this Hamiltonian induction
is an isomorphism, by what we just proved in \eqref{mpp}. 
That is to say, we get a $G$-equivariant
map
$$ \varpi: \textrm{h-ind}_{HU}^G \tilde{S}  (= M_+ \times^{HU} G) \rightarrow M$$
preserving symplectic forms. Since -- examining \eqref{V1V2} --  both 
sides have the same dimension, 
 $\varpi$ is an isomorphism 
on tangent spaces everywhere.   Therefore $\varpi$
is an unramified morphisms between smooth varieties
of the same dimension; it is then automatically {\'e}tale. 
 
 We will now argue that $\varpi$ is an isomorphism.
To do so it is enough to argue that $\varpi$ is {\em finite};
it is then a finite {\'e}tale cover, and its degree is constant
on the target; then use the fact that the preimages of points in $M_0$
under $\varpi$ have size $1$, which  will follow from inspection of what is happening on the closed orbit.

 To verify the desired finiteness, we invoke
a lemma of Luna \cite[Lemme, p89]{LunaSlice}
deduced from Zariski's main theorem.
Set $M'$ to be the source $\textrm{h-ind}_{HU}^G \tilde{S}  $ of $\varpi$; so we have a morphism
$$\varpi: M' \rightarrow M$$
of $G \times \Gm$-varieties. Luna's lemma shows that finiteness is assured 
if    both varieties are affine and $\varpi$ has finite fibers,   carries closed orbits to closed orbits,
and induces a finite map on invariant theoretic quotients. 
 
 That $M'$ is affine follows from the structure theorem \eqref{fund_id}
for vectorial Whittaker induction.   The only $G \times \Gm$-invariant
functions in $\FF(M)$ are constants. Therefore, $\FF(M')$ being algebraic
over $\FF(M)$, the same is true for $M'$. Therefore, $\FF[M']^{G \times \Gm} =\FF$;
so $M'$ has a unique $G \times \Gm$-closed orbit, which by construction
is carried to the closed $G \times \Gm$-orbit on $M$, and moreover
the invariant theoretic quotients are both points. By inspection, 
the map $\varpi: M' \rightarrow M$ is an isomorphism on these closed orbits. 
This concludes
our proof that $\varpi$ is finite, and therefore, as argued above, that it is an isomorphism.

Regarding uniqueness, take any $G\times\GGm$-equivariant 
automorphism of $M$ as a Hamiltonian space preserving $x$
and acting trivially on the symplectic normal bundle $S$. 
This induces an automorphism of $M_+$ that 
is necessarily trivial by the uniqueness statement of the {\em Claim} above. 
Thus the original automorphism of $M$ is also trivial on $G \cdot M_+$,
which contains an open subset of $M$ (e.g., by consideration of tangent spaces at $x$),  and so the automorphism is trivial. 
 \end{proof}

The following proposition provides sufficient and necessary conditions for a Whittaker induction as in the Structure Theorem \ref{thm:structure} to be coisotropic.

\begin{proposition}\label{Scoisotropic}
 For a Whittaker-induced space $M$ as in \eqref{Whittinduction}, if the subgroup $H$ is reductive, the quotient $Y=HU\backslash G$ is quasiaffine. The Whittaker induction \eqref{Whittinduction} is coisotropic if and only if $Y$ is spherical, and the inducing symplectic space $\tilde S$ (in the notation of \ref{tildeSdef}) is coisotropic for the generic stabilizer of $G$ on $T^*Y$. It is hyperspherical under the $\GGm$-action described in \S~\ref{ggmm} iff, in addition, it satisfies \eqref{condtorsion}.
\end{proposition}

\begin{proof}

To prove that $Y$ is quasiaffine, it is enough to show that $H$ is contained in the kernel of a $P$-regular dominant character of $L$, where $P=LU$ is the parabolic defined by the nonnegative eigenspaces for the element $h$ of the $\sl_2$-triple. (We call a character $P$-regular if it does not extend to a larger Levi subgroup.) 
Such a character is provided by the element $h$ itself, which becomes the differential of a character by means of an invariant isomorphism $\mathfrak g^*\simeq \mathfrak g$ (restricted to the Lie algebra of $L$). This character is trivial on $H$ because $h$ belongs to an $\sl_2\hookrightarrow \mathfrak g$ that commutes with $\mathfrak h$.  This proves the quasiaffine property.

By Proposition \ref{coisotropic}, $M$ is coisotropic iff the general nonempty fiber of the invariant moment map $M\to \mathfrak c^*$ contains an open (not necessarily dense) orbit. In our arguments, we will repeatedly use partial ``Springer--Grothendieck resolutions'' of $\mathfrak g^*$, determined by classes $\mathcal Q$ of parabolics, and defined as
\[ \widetilde{\mathfrak g}^{\mathcal Q*} = \{(Z,Q)| Z\in \mathfrak u_Q^\perp\subset \mathfrak g^*\},\]
that is, pairs consisting of a coadjoint vector and a parabolic in the given class whose nilradical is orthogonal to the given vector. We recall that the forgetful map $\widetilde{\mathfrak g}^{\mathcal Q*}\to \mathfrak g^*$ is surjective and proper, and that it is finite over the set of \emph{$Q$-regular elements} -- this can be taken as a definition of $Q$-regular, but under $\mathfrak g^*\simeq \mathfrak g$ this is also the set of all $Z$ that belong to a finite number of parabolic Lie algebras $\operatorname{Lie}(Q)$,  $Q\in \mathcal Q$. Choosing a representative $Q \in \mathcal Q$, the space $\widetilde{\mathfrak g}^{\mathcal Q*}$ can also be written $\mathfrak u_Q^\perp\times^Q G$. It lives over $\mathfrak c_{\mathcal Q}^*=$ the analog of $\mathfrak c^*$ for the Levi quotient of any representative $Q\in \mathcal Q$ (as they are all canonically isomorphic).

Therefore, if the Hamiltonian space $M$ is such that its moment image consists, generically, of $\mathcal Q$-regular elements, it is coisotropic iff the general nonempty fiber of the base-change map $\tilde M^{\mathcal Q}\to \mathfrak c_{\mathcal Q}^*$ contains an open orbit, where 
\[ \tilde M^{\mathcal Q} = M\times_{\mathfrak g^*} \widetilde{\mathfrak g}^{\mathcal Q*} = M^{\mathfrak u_Q^\perp}\times^Q G.\]

Closely related to this construction is the \emph{polarized} cotangent bundle, defined as
 \begin{equation}\label{Mpolarized}
  \hat M = M \times_{\mathfrak c^*} \mathfrak a^*.
 \end{equation}

By \eqref{Whittinduction}, the space $M$ is fibered over $Y$. 
By \cite[Satz 2.3]{KnWeyl}, there is a parabolic $Q$, a Levi decomposition $Q=RU_Q$, and an open $Q$-stable subvariety $Y_Q$ with  
\begin{equation}\label{XQisom} Y_Q \simeq R_0\backslash Q \times V
\end{equation}
 as a $Q$-variety, where $R_0\supset [R,R]$, and $V$ is a variety with trivial $Q$-action. We will use the presentation \eqref{Whittinduction} for $M$, and, without loss of generality, the coset of $1$ in $Y=HU\backslash G$ belongs to $Y_Q$, which means that $U_Q$ is ``opposite'' to $U$, i.e., $Y_Q$ lives over the open $Q$-orbit on $P\backslash G$. Note that $Y$ is spherical iff $V$ is $0$-dimensional (hence, by connectedness, a point). 
 
 Set $A_Y = R/R_0$.  
 From \eqref{XQisom} it is clear that the generic stabilizer for the $Q$-action on $Y$ is conjugate to $R_0$. It is known from a construction of Knop  that this is also the generic stabilizer for the action of $G$ on $T^*Y$. More precisely, Knop  \cite[\S~3]{KnMotion} produces polarized cotangent vectors on $Y$ by using the $Q$-moment map on $Y_Q$. Namely, the isomorphism \eqref{XQisom} implies that $T^*Y_Q^{\mathfrak u_Q^\perp} = T^*A_Y \times T^*V \times U_Q$, in such a way that the diagram 
 \[ \xymatrix{ T^*Y_Q^{\mathfrak u_Q^\perp} \ar[r]\ar[d]& T^*Y \ar[d] \\ 
 \mathfrak a_Y^* \subset \mathfrak a^* \ar[r] &\mathfrak c^*}
 \]
 commutes, where the projection to $\mathfrak a_Y^*$ is trivial on $T^*V$, and the natural projection on $T^*A_Y$. Since $Y$ is quasiaffine, \cite[Lemma 3.1]{KnMotion} implies that the generic element of $\mathfrak a_Y^*$ is $Q$-regular, and this implies that the map 
 \[  \widetilde{T^*Y}^{\mathcal Q} \supset T^*Y_Q^{\mathfrak u_Q^\perp} \times^Q G \to T^*Y\]
(where the embedding on the left is open) is generically finite and that the $G$-stabilizer of a generic element in $T^*Y$ is conjugate to the $Q$-stabilizer of a generic element in $\widetilde{T^*Y}^{\mathcal Q}$, hence conjugate to $R_0$.  
 
We can adapt Knop's construction to the Whittaker-induced space $M$, as follows: Recall, first, the structure of $M$ as a Whittaker induction \eqref{Whittinduction},  
\[ M= \tilde{S} \times^{HU}_{(\mathfrak h + \mathfrak u)^*} T^*G.\]
 It maps naturally to $Y$, and we can restrict it to the open set $Y_Q$. If we denote by $G_Q$ the preimage of $Y_Q$ under the action map, 
 we have $G_Q= UH \times^{R_0} Q$, equivariantly under the left-$HU$- and right-$Q$-action, hence 
\[T^*G_Q = T^*(UH) \times^{R_0}_{\mathfrak r_0^*} T^*Q \times T^*V,
\]
compatibly with the left moment map to $(\mathfrak h+\mathfrak u)^*$ and the right moment map to $\mathfrak q^*$ (where the $T^*V$-factor does not affect the moment map). 
Hence, 
\[ M|_{Y_Q} = \tilde{S} \times^{R_0}_{\mathfrak r_0^*} T^*Q \times T^*V,\]
compatibly with the moment map to $\mathfrak q^*$. Hence, the open subset $(\tilde M^{\mathcal Q})'= {M|_{Y_Q}}^{\mathfrak u_Q^\perp}\times^Q G$ of the $\mathcal Q$-cover $\tilde M^{\mathcal Q}$
is equal to 
\begin{equation}\label{MYQ}
(\tilde M^{\mathcal Q})' = \tilde{S} \times^{R_0}_{\mathfrak r_0^*} T^*R \times T^*V \times^R G.
\end{equation}

We have a natural action of $\mathfrak a_Y^*$ on the affine space $\mathfrak c_R^*$ (descending from its action on $\mathfrak a^*$), with quotient $\mathfrak c_{R_0}^*$.
The presentation \eqref{MYQ} shows that the invariant moment image of $(\tilde M^{\mathcal Q})'$ in $\mathfrak c_R^*$ is equal to the preimage of the invariant moment image of $\tilde{S}$ (considered as an $R_0$-space) in $\mathfrak c_{R_0}^*$. In particular, by the $Q$-regularity of the generic element of $\mathfrak a_Y^*$, mentioned above, the invariant moment image consists generically of $Q$-regular elements, and $M$ is coisotropic iff the general nonempty fiber of $(\tilde M^{\mathcal Q})'\to \mathfrak c_R^*$ contains an open $G$-orbit. Clearly, again by \eqref{MYQ}, this is equivalent to requiring that $\dim T^*V = 0$, and the generic nonempty fiber of $\tilde S\to \mathfrak c_{R_0}^*$ contain an open $R_0$-orbit; that is, that $Y$ be spherical and $\tilde S$ be a coisotropic $R_0$-space.

If, in addition, we assume \eqref{condtorsion}, then all the defining conditions for ``hyperspherical'' are satisfied (the affine property by Lemma \ref{lemma:affine}). 
\end{proof}
 
\begin{remark} \label{rem:uniqueness}
A neutral $\Gm$-action on such a Hamiltonian $G$-space $M$ is clearly quite rigid and it is likely
that it is in fact unique; it would be nice to prove this. 
 \end{remark}

 \subsection{Polarization by twisted cotangent bundles} \label{dp} 
We now discuss the question of polarizing $M$ in a fashion that is compatible with the $G$-action
and our later needs.  When such a polarization
exists there is a distinguished class of them --   ``up to the question of polarizing a symplectic vector space.''

 We continue to consider a hyperspherical Hamiltonian $G$-space, i.e., one that  
satisfies the conditions of
\S~\ref{condALL}.  As in the discussion after Definition \ref{condneutral}, 
 any point $x\in M_0$ (the closed $G\times\GGm$-orbit)  gives rise to an $\sl_2$-pair commuting with the stabilizer $H=G_x$,
 and via Theorem \ref{thm:structure} we get a map $M\to HU\backslash G$.  Let $S$ be, as before, the symplectic normal bundle to $M_0$. 
 
\begin{definition} \label{distinguishedpolarizationdef}
 In the setting above, we say that $M$ admits a {\em distinguished polarization} if
  the weight-1 component  $\mathfrak u_1\subset \mathfrak{u}$ vanishes, and 
 there is a  Lagrangian $H$-stable decomposition
\[S = S^+ \oplus S^-,\]
 In this case we will also call the choice of such a decomposition a {\em distinguished polarization} of $S$.
 \end{definition}
  
  \begin{remark} \label{dp-generalized} 
  The assumption that the weight-1 component of $\mathfrak{u}$ vanishes can be replaced,
  for many purposes of this paper, by the assumption that there exists an $H$-stable Lagrangian space in $\mathfrak u_1$.
We have not done so, largely because this would add the additional complication of verifying
independence of Lagrangians at various stages. We do not foresee any difficulty in doing this.  
\end{remark}
  
Note that, by the theory of $\sl_2$-modules, the vanishing of $\mathfrak u_1$ is equivalent to \emph{evenness}, i.e., $h$ acts on  $\mathfrak g$ with even weights only.

In this setting, notice the following:
\begin{lemma}
 \label{lemma-generic}
Assume that $h$ acts on $\mathfrak g$ with even weights. Then, the additive character $f$ of $\mathfrak u$ is \emph{generic}; that is, it lies in the open $L$-orbit on $\Hom_{\text{Lie}}(\mathfrak u, \mathfrak g_a)$, where $L$ is the Levi quotient of the parabolic whose unipotent radical is $U$.  
\end{lemma}

\begin{proof}
By the theory of $\sl_2$-modules, the abelianization $\mathfrak u /[\mathfrak u,\mathfrak u]$ is isomorphic to the weight-2 eigenspace $\mathfrak u_2$ under the canonical projection. Indeed, all eigenspaces of $\mathfrak u$ of weight $>2$ are generated by the $2$-eigenspace under the adjoint action of $e$, and, vice versa, all weight $2$ vectors commute up to vectors of weight $>2$. 

The dual $\mathfrak u_2^*$ decomposes as an $L$-module into a direct sum $\sum_{\alpha \in \Delta\smallsetminus \Delta_L} \mathfrak u^*_{2,\alpha}$, where $\Delta$, $\Delta_L$ denote the simple roots of $G$ and $L$, respectively, and the center  of $L$ acts on $\mathfrak u_{2,\alpha}^*$ by the character $-\alpha$ \cite{ABS}. 

We claim that 
the character $f$ has nontrivial projections to all summands $\mathfrak u^*_{2,\alpha}$, i.e., nontrivial restrictions to all the irreducible $L$-submodules $\mathfrak u_{2,\alpha}$ of $\mathfrak u_2$. Indeed, the adjoint action of $f$ sends $ \mathfrak u_{2,\alpha}$ injectively into $\mathfrak l$, therefore for every nonzero $x\in \mathfrak u_{2,\alpha}$ there is a $y\in \mathfrak l$ such that (under the invariant symmetric bilinear form identifying $\mathfrak g$ with $\mathfrak g^*$)
\[ 0\ne \langle [f, x], y \rangle  = \langle f , [x,y] \rangle,\]
hence $f$ does not vanish on the $L$-module generated by $x$.

By the equivalent conditions of \cite[Lemma 2.6.1]{SV}, the nonvanishing of $f$ on all $\mathfrak u_{2,\alpha}$ means that it is generic.
\end{proof}

In this case, $M$ has the structure of twisted cotangent bundle:
Take
\begin{equation} \label{polarizationXstructure} X = S^+ \times^{HU} G, \ \ \Psi = S^+ \times^{HU'} G,\end{equation}
with $U'$ the kernel of $U \rightarrow \mathbb{G}_a$. 
 We may identify $M = T^*(X, \Psi),$ in the notation of \S~\ref{twisted cotangent bundle section}. We can also write this as 
 \begin{equation}\label{Whittakerinduction}
  (X,\Psi) = \Ind_P^G (X_L,\Psi_L),
 \end{equation}
 where $P$ is the parabolic with unipotent radical $U$ (and Levi quotient $L$, $\Ind_P^G Y$ denotes the variety $Y\times^P G$, and $(X_L, \Psi_L)$ denote the pair of $P$-varieties $(S^+\times^{HU} P, S^+\times^{HU'} P')$.  Since the additive character $U/U'\xrightarrow\sim \Ga$ is generic by Lemma \ref{lemma-generic}, we are in the setting of ``Whittaker-type reduction,'' as introduced in \cite[\S~2.6]{SV}.

The neutral $\GGm$-action on $M$, described in \S~\ref{ggmm}, is induced by one on $X$ that is covered by an action on $\Psi$; namely, 
the action of $\lambda \in \GGm$ on $X$ is given by $\lambda\cdot (s^+,   g) =(\lambda s^+,   \lambda^h g)$,
and the same formula on $\Psi$. The $\Ga$- and $\GGm$-actions on $\Psi$ combine to a left action of $\Ga \rtimes \GGm$, where $\GGm$ acts on $\Ga$ by $\lambda \cdot x  \cdot \lambda^{-1} = \lambda^2 x$.

We now examine what conditions on $X$ are forced by the conditions of \S~\ref{ssconditions}. 
\begin{proposition} \label{hypersphericalspherical} \label{isspherical}
When a hyperspherical variety $M$ admits a distinguished polarization $M=T^*(X,\Psi)$,  then
\begin{itemize}
\item[(a)] $X$ is a spherical $G$-variety, and
\item[(b)] the $B$-stabilizers of points in the open $B$-orbit on $X$ are connected.
\end{itemize}
In particular, $X$ has no roots of type $N$ (for this  terminology, see
 \S \ref{subsubsec:colors}). Moreover, in that case the subgroup $H$ of $G$ as above is connected. Vice versa, a Whittaker-induced Hamiltonian space $M$ which admits a polarization of the form \eqref{polarizationXstructure} (with $H$ reductive) is hyperspherical if it satisfies these conditions.
\end{proposition}
 We recall that a (normal) $G$-variety is {\em spherical} if a Borel subgroup of $G$
 acts with a Zariski open orbit.  See also Remark \ref{notypeN}
 for an interpretation of condition (b) in terms of the combinatorial data attached to a spherical variety. 
 
\begin{proof}
The equivalence of (a)--(b) with the hyperspherical property follows from Proposition \ref{Scoisotropic}, if we prove that the connectedness of generic stabilizers on $M$ (Condition \ref{condtorsion}) is equivalent to Condition (b). Note that, again by Proposition \ref{Scoisotropic}, $X$ is quasiaffine. Let $B$ be a Borel subgroup, and $P(X)\supset B$ the parabolic stabilizing the open Borel orbit. As we recall later in \S~\ref{sec:invariants}, the generic stabilizer of $B$ on $X$ is equal to the intersection of $B$ with the generic stabilizer of $P(X)$ on $X$, which is a subgroup contained in a Levi of $P(X)$, and containing its derived group. In particular, the generic stabilizer in $B$ is connected iff the generic stabilizer in $P(X)$ is connected. The latter, by \cite[Korollar 8.2]{KnWeyl} (whose proof also applies to the twisted case), is conjugate to the generic stabilizer of $G$ on $T^*X$. This proves the desired equivalence.

The statement about roots of type $N$ is a consequence of (b), see Remark \ref{notypeN}.

If $H$ were not connected, we would have a finite \'etale map of $G$-spherical varieties $X' = S^+ \times^{H^\circ U} G\to X = S^+ \times^{HU} G$, and in particular the stabilizers in the open $B$-orbit on $X'$ would be of finite index $>1$ in the stabilizers of the open $B$-orbit on $X$; this contradicts (b).

\end{proof}

  \subsection{Eigenmeasures}
 \label{eigencharacter}
\index{eigenmeasure} \index{eigencharacter}

Given a distinguished polarization $M=T^*(X, \Psi)$ 
there is a further condition on $X$ that is important:

 We may ask that $X$ admits a {\em nowhere vanishing eigen-volume form}, which
 we sometimes just call {\em eigenmeasure}:  a nowhere vanishing algebraic differential form
of top degree, with the property that, up to scaling, it is preserved by $G$. 
Such a form is then also automatically preserved up to scaling by $\GGm$;
for its translate by $\GGm$ gives another form $\omega'$ with the same $G$-eigencharacter,
and then $\omega'/\omega$ is a nowhere vanishing $G$-invariant function on $X$,
so constant (since $X$ admits an open $G$-orbit, by Proposition \ref{isspherical}).

 \index{$\eta$=eigencharacter of $G$} \index{$\gamma=$ eigencharacter of $\Gm$}
 In particular, having fixed such a form $\omega$, it determines a character $\eta: G \rightarrow \Gm$ 
 and an integer $\gamma \in \mathbb{Z}$ with the property that
\begin{equation} \label{gammaXdef} (g, \lambda)^* \omega = \eta(g)  \lambda^{\gamma} \cdot \omega.\end{equation}
 
 To appreciate the relevance of this condition, we note that the conjectures that we are about to formulate in this paper are best ``calibrated'' by working with half-densities on $X$, rather than functions. However, half-densities are a little awkward geometrically and arithmetically, and we would like to have the possibility to translate to functions. The existence of such an eigenform allows us to do so.

\subsubsection{Translation in terms of $H$}  \label{Htrans}
 
Although the above definition makes sense for an arbitrary $G \times\GGm$-space $X$
with open $G$-orbit and eigenform,  
 let us specialize to the case of \eqref{polarizationXstructure} and explicate the situation. 
 
 Let $\eta$ be the character by which $H$ acts on the top exterior power of the tangent space at the point $(0,1) \in  S^+\times^{HU}G=X$.
 Then $X$ admits a nowhere vanishing eigen-volume form $\omega$ if and only if
 $\eta$ extends to a character    of $G$ (to be denoted by the same letter);
 for any such extension, there is a unique, up to scaling, eigen-volume form on $X$ with eigencharacter $\eta$.   
 (To help with signs, it may be helpful to observe that that the action of $H$ on the tangent space is defined as a \emph{right} action, via pushforward $v \mapsto h_* v$ of tangent vectors, while the action on functions and differential forms is defined as a \emph{left} action.)
 In this situation we readily compute
that the quantity $\gamma$ of \eqref{gammaXdef} is given by \begin{equation} \label{gammadef} \gamma = \dim(S^+) - \langle 2 \rho, \varpi \rangle,\end{equation}
where $\varpi$ is the character associated to the $\SL_2$ for $(G, M)$,
and $2\rho$ is the sum of roots on $U$.  Later on we will frequently encounter
the following expression, which we also compute for later reference
 \begin{equation} \label{sigmaXdef}
\beta_X:= \dim(G) + \gamma - \dim(X)  = \dim(HU)  - \langle 2 \rho, \varpi \rangle.\end{equation}

 \subsubsection{Why we allow ourselves to often assume that $X$ has an eigenmeasure} \label{ssseigencharinocuous} 
The assumption that $X$ admits an eigen-volume form is innocuous for us. Most of the
issues we consider in this paper can be reduced to that case, although we have
not at present written this out in all cases; and, moreover,
the choice of eigen-volume form makes no difference. 

First of all, even if $X$ does not have such an eigen-volume form, there is a $\Gm$-cover of it that does: Indeed, replace $H$ by the kernel $H_0$ of the character $\eta$ above.\footnote{
Note that $\eta$ is necessarily nontrivial, for otherwise $X$ would have an eigen-volume form, and so $H_0$ is a proper subgroup of $H$. }
Consider the variety $\tilde X = S^+\times^{H_0U}G$, as a $\tilde G = H/H_0 \times G$-space. 
By this technique, many issues studied in this paper can be reduced to the  
 case when $X$ has an eigenmeasure.

\begin{example}
Take $G=\GL_3$ and  $X=\Gm \cdot U\backslash G$, 
where $\Gm$ is embedded in the $(1,1)$ entry and $U$ is $\small{\begin{pmatrix} 1 & 0 & * \\0 & 1 & *  \\ 0 & 0 &1 \end{pmatrix}}$;
we equip $X$ with the $\Ga$-torsor defined by the homomorphism $U \rightarrow \Ga$ defined by the $(2,3)$ entry. 
Then $X$ does not have an eigenmeasure; however,
$$\tilde{X} = U\backslash G \mbox{ as a $G \times \Gm$-space}$$
does have an eigenmeasure, where $\Gm$ is again acting through the $(1,1)$-entry, and 
$\tilde{X}/\Gm = X$ as a $G$-space. 
 The space $\tilde{X}$ represents the standard $L$-function of $G \times \Gm$.
 \end{example}

 Secondly, an eigenmeasure is not unique. What
 actually plays a bigger role for us is
  the character $\eta$, and 
  it too is therefore not uniquely defined.  However, the resulting ambiguity is again essentially irrelevant for us.
  From the point of view of our later duality theory, two eigenmeasures
 $\omega$ and $\omega'$ induce characters $\eta, \eta'$ with the following property:
 the dual of $\eta'/\eta$, considered now as a central cocharacter $\Gm \longrightarrow \check{G}$,
 will act {\em trivially} on the Hamiltonian space dual to $X$.  
Our local and global conjectures will apparently use the choice of $\eta$, but this
fact means that the choice does not matter.

\subsection{Hyperspherical varieties over general fields}  \label{hdprings} 
 
Throughout this chapter we have worked with hyperspherical varieties $(G, M)$
over an algebraically closed field  $\FF$ in characteristic zero.

  It is, of course, desirable to have a theory over general fields or  rings. While it would be nice to develop such a notion from a list of properties such as the conditions used to define hyperspherical varieties over $\C$ in \S~\ref{condALL}, we will not do so in this paper. Rather, we will use the Structure Theorem \ref{thm:structure}, and \emph{mandate} that the varieties that we will call hyperspherical over more general rings have this structure. 
  
  Our working definition of ``hyperspherical variety''  for this paper will be Definition \ref{def:hyperspherical}. 
This definition is far from satisfactory,  and it is not clear whether it produces the right objects in arbitrary characteristic,  but it will be complemented by the definition of a ``distinguished split form,'' Definition \ref{dhpFqdef}, which, when available, serves as a distinguished base point for the dualities that we propose later in this paper.

 \subsubsection{Forms of a hyperspherical datum over a ring} \label{hdr-ring}
  
How should we define the notion of ``hyperspherical variety $(G, M)$ over a ring $R$?''
 
The theory of reductive group schemes, mentioned again below,
gives a satisfactory notion of a form of $G$ over $R$. 
For 
  $M$, the most optimistic and pleasing interpretation would be a smooth affine $R$-scheme equipped with an action
of $(G \times \GGm)_{/R}$ and a Poisson bracket of degree $-2$
relative to the $\GGm$-action, which arise in fact from a symplectic
structure on each geometric fiber.    
   Ideally speaking, we would then formulate the definition from \S \ref{condALL} in a way that made sense relative to $\mathrm{Spec} \ R$, and   
  then develop our structure theory in that context. 
 
We have not proceeded in such a systematic way. 
  Rather we will content ourselves with a simple way of constructing
  useful examples over rings.

  Our Structure Theorem \ref{thm:structure}
  asserts that, over $\C$, a hyperspherical $G$-variety is defined  (up to isomorphism) by a 
        linear-algebraic datum $\mathcal{D}$:
    \begin{equation} \label{hypersphericaldatadef} \mathcal{D} = \mbox{  $\iota: H \rightarrow G$,
    commuting $\mathfrak{sl}_2$-pair $(\varpi,f)$, and $\rho: H \rightarrow \mathrm{Sp}(S)$, }
\end{equation}
where $S$ denotes a symplectic vector space. 
We recall (\S~\ref{sl2pair}) that $\varpi$ is a cocharacter $\Gm \rightarrow G$. We will call $\mathcal D$ a \emph{hyperspherical datum}. \index{hyperspherical datum} \index{$\mathcal D$}

A hyperspherical $G$-variety with a distinguished polarization (\S~\ref{dp}) $(X, \Psi)$  is also determined (up to isomorphism)
by linear-algebraic data, namely   \begin{equation} \label{dplusdef} \mathcal{D}^+ = \mbox{$\iota: H \rightarrow G$,
commuting \emph{even} $\mathfrak{sl}_2$-pair $(\varpi,f)$, and
  $\rho^+: H \rightarrow \GL(S^+)$. }
\end{equation}
 where $S^+$ is a vector space, and we recall that ``even'' means that $\mathfrak g$ is a sum of odd-dimensional $\SL_2$-representations. The associated datum $\mathcal{D}$ is obtained
by setting $S=S^+\oplus (S^+)^*$ and replacing $\rho^+$ by its composition $\rho$  with the standard inclusion $\GL(S^+) \rightarrow \Sp(S)$. We will call $\mathcal D^+$ a \emph{polarized hyperspherical datum}. \index{polarized hyperspherical datum} \index{$\mathcal D^+$}
 
 \begin{remark} \label{Xdatumuniqueness}
For later use, we observe that, all other data of $\mathcal D^+$ being equal, any two choices $(\rho_1^+,S_1^+)$ and $(\rho_2^+,S_2^+)$ giving the same hyperspherical datum $\mathcal{D}$ over $\C$ are related as follows: 
There are decompositions of $H$-representations, $S_i^+ = \bigoplus_i S_{i,j}^+$, such that, for every $j$, $S_{1,j}^+$ is equal to $S_{2,j}^+$ or to its dual. 
  \end{remark}

It is very easy to define the notion of a ``form of $\mathcal{D}$ over a ring,''
by using the existing satisfactory theory of reductive group schemes developed in \cite{SGA3} (see also \cite[\S 3, \S 5]{ConradSGA3}); 
in particular, the geometric fibers of such a group scheme are reductive, 
the isomorphism class of the associated root data is locally constant on the base,
and there is a notion of a split reductive group scheme, obtained by base change from the Chevalley group scheme over $\mathbb Z$. (See also the discussion surrounding \cite[Definition 5.1.1]{ConradSGA3}.)

We note that, following the usual convention, ``reductive group schemes'' are required to have connected geometric fibers.
Therefore, the following definition only models hyperspherical varieties where $H$ is connected -- see Remark \ref{Hconnected}.

  \begin{definition}  \label{ringdatum} 
 Let $R$ be a subring of $\C$, $G$ a reductive group scheme over $R$, and $\mathbb F$ a field contained in the algebraic closure of $\mathbb F_p$, for some prime $p$. 
  
 \begin{itemize}
 \item[(i)] A $\mathfrak{sl}_2$-pair over $R$ is a pair
 $$(\varpi: \Gm \rightarrow G, f \in \mathfrak{g}^*(R))$$
arising from a homomorphism $\rho: \SL_2 \rightarrow G$ of group schemes over $R$,
where $\varpi$ comes by restriction to the maximal torus, 
and $f$ arises from $d\rho {\small \begin{pmatrix} 0 & 0 \\ 1 & 0 \end{pmatrix}}$
by means of a $G$-equivariant  isomorphism $\mathfrak{g}_R \rightarrow \mathfrak{g}_R^*$.
(See below for discussion.)

 \item[(ii)]
 Given a hyperspherical datum $\mathcal{D}$ 
 as in \eqref{hypersphericaldatadef},
an $R$-form $\mathcal D_R$ of $\mathcal{D}$ 
consists of a triple $(H_R, G_R, \Sp_R)$ of reductive group schemes over $R$, together with an injective morphism $\iota: H_R \rightarrow G_R$ and a morphism
$\rho: H_R \rightarrow \mathrm{Sp}_R$, and an $\mathfrak{sl}_2$-pair $(\varpi, f)$ over $R$,
centralized by $H_R$,  which recover (the isomorphism class of) $\mathcal{D}$ after base change via $R \rightarrow \mathbb{C}$. \index{hyperspherical datum, over rings}

\item[(iii)]  A hyperspherical datum over $\FF$ is a similar collection of data over $\mathbb F$, which is obtained by base change via a homomorphism $R\to \FF$ from a hyperspherical datum over a subring $R$ of $\CC$. 

\item[(iv)]
We say that $\mathcal{D}_{R}$ (resp.\ $\mathcal D_{\FF}$) is \emph{split} if $H_R, G_R$ are split, i.e., they are the Chevalley
forms over $R$. 
 \end{itemize}
 
 Similar language will be used for a ``polarized hyperspherical datum''
 as in \eqref{dplusdef}. 
 \end{definition}

   \begin{remark} (The bilinear forms in the definition of $\mathfrak{sl}_2$-pair):  
Our use of a $G$-equivariant  isomorphism $\mathfrak{g}_R \rightarrow \mathfrak{g}_R^*$  in (i) (equivalently, a nondegenerate, $G$-invariant bilinear form on $\mathfrak g_R$) is, certainly,  a little crass. 
We leave to future work  a more intrinsic formulation. 
Such a bilinear form as above always exists if \footnote{
 When $G$ is the Chevalley form of a semisimple group the existence of such a form is proved in  \cite[Lemma 4.2.3]{Riche}; see also \cite[\S 4]{SteinbergSpringerConjugacy}.
 The general case reduces to this one, as follows: 
  We have 
an isogeny 
$Z \times G^{\scc} \rightarrow G$
where $Z$ is the connected center of $G$,
and $G^{\scc}$ is the simply connected cover of the derived group of $G$.
The kernel of this isogeny is isomorphism to a subgroup of the center 
of  $G^{\scc}$, and in particular has order divisible only by 
``bad'' primes excluded above.
It follows that we have an isomorphism of Lie algebras
\begin{equation} \label{splitLie} \mathfrak{z}_{R} \oplus \mathfrak{g}^{\scc}_{R} \rightarrow \mathfrak{g}_{R}.\end{equation}
which permits us to extend an invariant bilinear form  from $\mathfrak{g}^{\scc}_R$ to $\mathfrak{g}_R$.}
  \index{$N_G$} \index{bad primes}
\begin{equation} \label{NGdef} N_G = \prod_{p \in 2 \cup B} p\end{equation}
is invertible in $R$, where $B$ ranges over the set of all primes that are not very good for $G$. 
 Explicitly, the primes  in $B$ are
divisors of $n+1$ for type $A_n$,
 $2$ for other classical types, $\{2,3\}$ for all exceptional types except $E_8$,
 and $\{2,3,5\}$ for $E_8$.
We will only use the above definition when $N_G$ is invertible in $R$.

    \end{remark}

\begin{definition} \label{def:hyperspherical}
\index{hyperspherical variety, over rings} Let $R$ (or $\mathbb F$) be as in Definition \ref{ringdatum}. A \emph{hyperspherical scheme} $M$ over $R$ (resp.\ $\mathbb F$) is a $G$-space that is obtained from a  hyperspherical datum by the Whittaker induction process of \eqref{Whittinduction}, 
\begin{equation} \label{YMdef} M= \Big(S \times (\mathfrak{u}/\mathfrak{u}_+)_f \Big) \times^{HU}_{(\mathfrak h + \mathfrak u)^*} (\mathfrak g^*\times G).\end{equation}
We say that $M$ is split if the datum is split. 

Similarly, a \emph{polarization} of (the isomorphism class of) $M$ is (the isomorphism class of) the pair of $G$-spaces $(X,\Psi)$ obtained from a polarized hyperspherical datum over $R$ (or $\FF$) by \eqref{polarizationXstructure}.
\end{definition}

\begin{remark} \label{ringdatumremark}    
 It is not clear that every datum obtained by the process above can be part of our conjectural dualities and deserves to be called ``hyperspherical;'' we have not examined the peculiarities that could arise over general rings, and definitely a robust structure theory of hyperspherical schemes needs to be developed. Nonetheless, the introduction of a ``distinguished split form'' that follows introduces at least one form (over sufficiently large finite fields) that should be part, and in some sense the ``distinguished base point,'' of our dualities.

 One of the issues that might arise has to do with the quotient  by $HU$ implicit in \eqref{YMdef}: the resulting $M$,  {\em a priori} a stack, may not be
  represented by an affine scheme.
 {Let's say it's a scheme, if we have confirmed this.} Our analysis used
 the theory of Slodowy slices (\S \ref{Slodowy}) and to use this argument one must assume that a certain list of ``bad'' primes, relative to the $\mathfrak{sl}_2$-pair,
  is invertible.
 See \cite[Theorem 4.3.3]{Riche} for the  case of the principal $\mathfrak{sl}_2$;
 it should be routine to formulate an explicit list of bad primes. 
  It is entirely possible, however, that the quotient of \eqref{Whittinduction} is in fact represented by an affine 
  scheme under weaker conditions. This is an interesting and important question to investigate, and
 seems to be the case in some simple examples we examined. 
\end{remark}

\subsubsection{Distinguished split form}

 The following Proposition will be used to show that there is a distinguished split form of a hyperspherical pair $(G, M)$ over finite fields, in some cases.

\begin{proposition} \label{rigidity} 
  Let $(G \times \GGm, M)$ be a hyperspherical
 variety over $\C$, in the sense of  \S \ref{ssconditions},
 defined by a datum $\mathcal{D}_{\C}$ as in Definition \ref{ringdatum}. 
Let $N_G$ be as in \eqref{NGdef}
and let $\Z'=\Z[\frac{1}{N_G}]$. There are integers $p_0, N$ such that the following holds.
 
 For $p \geq p_0$ and for any field $\FF$  algebraic over  $\mathbb{F}_{p^N}$, 
there is at most one, up to isomorphism, 
   datum $\mathcal{D}_{\FF}$ as in Definition \ref{ringdatum}
   which fits
   into a diagram
\begin{equation} \label{FZC} \mathcal{D}_{\FF} \leftarrow \mathcal{D}_{\Z'} \rightarrow \mathcal{D}_{\C},\end{equation}
 where  $\mathcal{D}_{\Z'}$ is a hyperspherical datum over $\Z'$ and the arrows are obtained by base change.
Moreover, if the automorphism group of $\mathcal{D}_{\C}$ is connected\footnote{By this we mean
that the centralizer of $H \times \SL_2$ in $G$ and the centralizer of $H \rightarrow \mathrm{Sp}_{2g}$ are both connected.}, then 
we   may take $N=1$, i.e., the above statement is valid for all fields $\FF$
   of sufficiently large finite characteristic.  
 
\end{proposition}
\proof
Suppose first that $\FF$ is the algebraic closure
of a finite field.  One first checks that any two choices of
data $\mathcal{D}_{\FF}$ which fit into a diagram \eqref{FZC} 
are conjugate so long as $p$ is sufficiently large. 
This follows for $(H \times \Gm \rightarrow G, H \rightarrow \mathrm{Sp}_{2g})$
via standard constructibility arguments, using ``rigidity'' of homomorphisms
between reductive groups (see
 Appendix \ref{dumblemmaproof} for a spelling out). 
Then (again, in large enough characteristic) any two choices for $f \in \mathfrak{g}_{\FF}^*$ are conjugate under
the centralizer of $H \times \Gm$ by the same reasoning
as in \cite[Theorem 4.2]{KostantTDS}.

Consequently,  once we fix one choice
of $\mathcal{D}_{\FF_p}$, 
the other choices
 are indexed by a Galois cohomology group $H^1(\mathrm{Gal}_{\FF_p}, Z)$, 
where $Z$ is the automorphism group of $\mathcal{D}$, considered
as an algebraic group over $\FF_p$.

 We now use the following fact: the restriction map  \begin{equation} \label{vanFF} H^1(\mathbb{F}_p, Z) \rightarrow H^1(\mathbb{F}_{p^k}, Z)\end{equation} 
 vanishes identically whenever $k$ is divisible by the product of $\# \pi_0(Z)$ and $\# \mathrm{Aut}(\pi_0 Z)$,
 where $\pi_0(Z)$ denotes the component group of $Z_{\overline{\mathbb{F}_p}}$. 
  In fact by Lang's theorem the map from $H^1(\mathbb{F}_{q}, Z)$ to $H^1(\mathbb{F}_q, \pi_0 Z)$ is injective
 which permits us to replace $Z$ by $\pi_0 Z$, considered as a finite {\'e}tale group scheme over $\mathbb{F}_p$.  Let us do so. Then a $1$-cocycle for the Galois group of $\mathbb{F}_p$  
 is determined by the image  of Frobenius $g \in   Z(\overline{\mathbb{F}_p})$;
 the image of the $k$th power of Frobenius is then given by the product 
 $$ g \cdot g^{F} \cdot g^{F^2} \dots g^{F^k},$$
 where the superscript $F$ denotes the Frobenius action on $Z(\overline{\mathbb{F}_p})$;
 and this product is automatically trivial when $k$ is divisible by the product of the order of $Z$
 and the order of the automorphism group of $Z$. 

Given this statement about  vanishing of the restriction morphism \eqref{vanFF},  the desired conclusion follows
because the set of possibilities for $\pi_0(Z)$ is finite -- this is deduced by means of general constructibility arguments, see e.g.\ \cite[9.7.9]{EGAIV}. 
Similarly, if the automorphism group of $\mathcal{D}_{\C}$ is connected, we see by general constructibility arguments that the same
is true for $Z$ so long as the characteristic $p$ is large enough. 
 \qed

Based on Proposition \ref{rigidity}, 
we offer  the following working definition
of a distinguished class of split  forms; it does not capture all forms that should be considered split forms,
but where it works, it should give the correct form.  However, first, a

\begin{quote}  {\bf Suggestion to the reader:}   rather than using the definition below,
take the more practical attitude that,  in most examples, the 
split form is either obvious or can be determined from its expected properties
and the conjectures of this paper
by a small amount of experimentation.
\end{quote}

 \begin{definition} 
\label{dhpFqdef}\index{distinguished split form}  
  Let $(G \times \GGm, M)$ be a hyperspherical
 variety over $\C$, in the sense of  \S \ref{ssconditions},
 defined by a datum $\mathcal{D}_{\C}$ as in Definition \ref{ringdatum}. 
Let $N_G$ be as in \eqref{NGdef}, 
 let $\Z'=\Z[\frac{1}{N_G}]$ and $\FF$ any field of characteristic not dividing $N_G$. 
   \begin{itemize}
 \item[(a)]
 A \emph{distinguished split form}   $(G\times \GGm, M)_{\FF}$
 over $\FF$ is one satisfying the following two conditions: 
 \begin{itemize}
 \item[(i)]  It arises as in Definition \ref{def:hyperspherical} and Definition \ref{ringdatum} from a hyperspherical datum $\mathcal{D}_{\FF}$  over $\FF$ which fits into a diagram
 $\mathcal{D}_{\FF} \leftarrow \mathcal{D}_{\Z'} \rightarrow \mathcal{D}_{\C}$;
   \item[(ii)]
 $\mathcal{D}_{\FF}$ is the unique (up to isomorphism) datum
 that fits into such a diagram.
 \end{itemize}

\item[(b)]
A \emph{distinguished split form}  of a twisted polarization $(X, \Psi)$ of $M$ 
over $\FF$  is
one arising from a polarized datum $\mathcal{D}^+_{\FF}$ (see \eqref{dplusdef})
polarizing some $\mathcal{D}_{\FF}$ that satisfies the conditions of (a).

  \end{itemize}
 \end{definition}
 
 Let us summarize what Proposition \ref{rigidity} says about the existence of such forms: 
 \begin{itemize}
 \item 
In (a), assuming that $\mathcal{D}_{\C}$ admits a lift to $\mathcal{D}_{\Z'}$,  
 such a form  exists and is unique when  $\FF=\overline{\FF_p}$
  for all sufficiently large $p$, or when
  $\FF = \FF_{p^k}$ for all sufficiently large $p$ and all sufficiently divisible $k$. \footnote{Both ``sufficiently large'' and ``sufficiently divisible''
  arise from the general constructibility arguments
in the proof of Proposition \ref{dumblemma}; with enough work, this could probably be made explicit.}
In the favorable case when the automorphism group of $\mathcal{D}_{\C}$ is connected, such forms
exist for all finite fields of sufficiently large characteristic. 

\item 
In (b), 
  such split forms exist under the same conditions as in (a),  and are unique up 
   the ambiguity in polarizing a symplectic vector space, as specified in
  Remark \ref{Xdatumuniqueness}.
  
  This does not follow directly from Proposition \ref{rigidity}: 
 we need to show
 that  a datum $\mathcal{D}_{\FF}$ arising in (a)
 can in fact always be polarized. The datum $\mathcal{D}_{\Z'}$
involves a homomorphism $\rho: H \rightarrow \mathrm{Sp}_{2g}$
 over $\Z'$ and, here, $\rho_{\C}$ can be polarized, by assumption. 
Extension of scalars gives a bijection between irreducible representations of $H_{/\Q}$
and $H_{/\C}$, and this preserves both symplectic self-duality and orthogonal self-duality.
 From this, we deduce that $\rho_{\Q}$, too, can be polarized,
and choosing an arbitrary integral lattice we see that $\rho_{\mathbb{F}_p}$
can be polarized for all sufficiently large $p$. Uniqueness here 
 follows as in Remark \ref{Xdatumuniqueness} so long 
 as we assume the characteristic is sufficiently large
that everything is semisimple, in particular so that there are no $\Ext$s between simple factors of the representation underlying $\rho$
 (which is readily seen to be valid in large enough characteristic; for instance it follows from \cite[Part II, 6.17]{JantzenAlgebraicGroups}). 
 \end{itemize}

  We presume this definition will be rendered obsolete by a more sophisticated study of rationality questions.

We finish this section with a useful lemma about the application of the duality involution (\S~\ref{dualityinvolution}) on a distinguished split form of a hyperspherical variety. Let $M^d$ denote the space $M$ with $G$-action and moment map twisted by the duality involution, and $\overline{M}$ the space $M$ where we negate the symplectic form and moment map. 

\begin{lemma} \label{dualitystableX} 
  Let $(G, M)_{\F_q}$ be a  distinguished split form of a hyperspherical variety
over $\FF_q$ in the sense of Definition \ref{dhpFqdef}. Then
  $\bar{M} \simeq M^d$ as $G \times \GGm$-spaces over $\FF_q$.  \end{lemma}

Before we give the proof, we emphasize that this is less interesting than it appears-- the notion
of ``distinguished split form'' in Definition \ref{dhpFqdef} is very restrictive, and allows us to avoid
serious subtleties.  But, as we will explain in \S~\ref{dhpFq}, there should be a distinguished split form in a more general setting, and the statement of the Lemma may be a good desideratum for what the correct notion of such a form should be.

To prove that the statement, we will use some facts about the theory of Cartan involutions on real groups. 
This is not because of any special role of the real numbers,  but rather
this is a context where involutions closely related to the duality involution have been studied.
 \proof

Let us pick data $\mathcal{D}_{\C}$  defining $M$  over $\C$ as in \eqref{hypersphericaldatadef}. 
Thus $\mathcal{D}_{\C}$ consists in  $\iota: H \rightarrow G$, a
    commuting $\mathfrak{sl}_2$-pair $(h, f)$, and $\rho: H \rightarrow \mathrm{Sp}_{2g}$. 
    The negated manifold $\bar{M}$ is then defined by negating $f$ and 
    replacing $\rho$ by $\bar{\rho}$ (i.e., conjugating through an element of $\mathrm{GSp}_{2g}$
    that negates the form).

It will be enough to show that $\overline{\mathcal{D}}_{\C} \simeq \mathcal{D}_{\C}^d$ where a superscript
$d$ means that we twist the datum through $d: G \rightarrow G$, and the bar means that we negate
the datum in the sense just described.
\footnote{Thus, if $\mathcal{D}$ is 
defined by $(\iota, h, f, \rho)$, we write $\overline{\mathcal{D}} = (\iota, h, -f, \bar{\rho})$,
and $\mathcal{D}^d = (d \circ \iota, d(h), d(f), \rho)$. }
In fact, over $\C$, $\bar{\rho}$ and $\rho$ are automatically conjugate, but we will
distinguish them for conceptual clarity.

  Indeed 
assuming  $\overline{\mathcal{D}}_{\C} \simeq \mathcal{D}_{\C}^d$ is valid, any diagram
 $\mathcal{D}_{\mathbb{F}_q} \leftarrow  \mathcal{D}_{\Z'} \rightarrow  \mathcal{D}_{\C}$
 also gives a diagram
  $\overline{\mathcal{D}}_{\mathbb{F}_q}^d \leftarrow  \overline{\mathcal{D}_{\Z'}}^d \rightarrow  \mathcal{D}_{\C}$.
Since part of the very definition of the distinguished split form is that any diagram of this form 
 results in the same $\mathbb{F}_q$-form, we find $ \mathcal{D}_{\mathbb{F}_q}^d \simeq \overline{\mathcal{D}_{\mathbb{F}_q}}$.

Let us recall that a {\em Cartan involution} of a real reductive group
 $J$  is an  involution $\theta_J$
with the property that the associated antiholomorphic involution $g \mapsto \overline{\theta(g)}$
of $J_{\C}$  has compact fixed set in $J$. (A general reference in essentially this context is  \cite{AdamsTaibi}). Moreover:
\begin{itemize}
\item  All such involutions are conjugate by $J(\mathbb{R})$,
which follows from \cite[Theorem 3.1]{MostowSAG}, see e.g
\cite[Theorem 3.12]{AdamsTaibi}. 

\item If $\theta_J$ is a Cartan involution, and $J_0 \subset J$ is a real reductive subgroup
stable under $\theta_J$, then the restriction of $\theta_J$ to $J_0$ is again a Cartan involution,
for it evidently has the same compactness property. 
 
\item
 If $J$ is split, there exists such an involution  which, taken with respect to 
  a Chevalley basis inverts the split torus and   sends the basis element $X_{\alpha}$
to $-X_{-\alpha}$, see e.g.\ \S 2.2 of \cite{DFdG}.  In particular, the complex-linear  extension $\theta_J^{\C}$ of the Cartan involution
lies in the inner class of the duality involution.
\item  
  Mostow  \cite[Theorem 5.1]{MostowSAG} has proven that for an embedding $G_1 \subset G_2 \subset \GL_{n}$
  of reductive groups over $\R$ there   exists a quadratic form on $\R^n$ with the property that the inverse-transpose
$\theta:g \mapsto (g^t)^{-1}$ fixes both $G_1$ and $G_2$. 
Then $\theta$ is   a Cartan involution for $\GL_n$ and simultaneously induces one on $G_1, G_2$. 
 \end{itemize}

Returning now to our context, 
we may, by definition of the distinguished split form,
assume that the datum $\mathcal{D}_{\C}$
defining $M$ is in fact defined over $\R$.
Let us fix an invariant bilinear form 
on $\mathfrak{g}_{\R}$, such that $(h,f)$ arises from a morphism of real algebraic groups
$\SL_2 \rightarrow G$, whose image commutes with $H$.  In what follows, we will work with real algebraic groups.

  We fix a Cartan involution $\theta_G$ of $G$ that fixes  $H$;
call that restriction $\theta_H$.  By direct computation, $\theta_G$ fixes the bilinear form on $\mathfrak{g}$. 
$\theta_G$ also induces 
a Cartan involution on the connected centralizer $Z$ of $H$;
by conjugacy of such involutions, we can further conjugate by an element of $Z(\R)$
so that $\theta$ preserves also the image of $\SL_2$ in $Z$, therefore
inducing a Cartan involution on it, which we can further assume to be the conjugation action of the standard Weyl element $w ={\small \left[\begin{array}{cc} 0 & 1\\ -1 & 0 \end{array}\right]}$.
Write $w^* = w \cdot e^{2 \rho}(\sqrt{-1})$.  Then
\[ \overline{\mathcal{D}} =  (\iota, h, -f, \bar{\rho}), \mbox{ and } \]
\begin{multline*}
\mathcal{D}^d \simeq (\theta_G \iota, \theta_G   h, \theta_G(f), \rho) =  (\iota \theta_H,-h, -e, \rho) = \\ \mathrm{Ad}(w^*)  (\iota \theta_H,  h, -f, \rho) 
\simeq (\iota, h, -f, \rho \theta_H).
\end{multline*}
 
  Finally, we readily verify that $\rho \theta_H$ and $\bar{\rho}$ are conjugate inside
 $\mathrm{Sp}_{2g}$. 
  \qed

\section{The dual Hamiltonian space to a polarized hyperspherical variety}
\label{dualofX}\label{spherical}

\subsection{Outline and motivation} \label{omspherical}
 In the previous section \S \ref{sphvar}, whose notation will be used here,  we introduced a certain class of ``hyperspherical'' Hamiltonian $G$-spaces (see \S~\ref{ssconditions}) over an algebraically closed field $\FF$ in characteristic zero.
Throughout this section, we continue to assume that $M$ satisfies the assumptions of \S~\ref{ssconditions}, but we assume, in addition, that it admits a distinguished polarization, in the sense of \S~\ref{dp}:
$$ M =T^*X \mbox{ or } M =T^*(X, \Psi).$$

In this polarized case, the theory is quite a lot better developed,  and, in the current section,
we will construct  an explicit candidate for its Hamiltonian dual $\check M$ in terms
of the geometry of spherical varieties.\footnote{{Several of the invariants associated to spherical varieties and more general $G$-spaces, that we are using, have been generalized by Losev \cite{Losev} to general Hamiltonian $G$-spaces. However, several essential elements of our construction of $\check M$, such as the dual group of a spherical variety, are still missing in the general case.}}  
This $\check{M}$ will be defined over an algebraically closed field $k$
of characteristic zero (although in  \S~\ref{CheckMRat} we will also discuss the case where $k$ is not algebraically closed).
The construction will depend on a conjecture (\ref{conjsymplectic}), which we have confirmed on all examples that we checked. 
We anticipate that the dual will not depend on the choice of polarization of $M$ (i.e., on the choice of polarization of the symplectic vector space $S$, in the notation of \S~\ref{dp}).

In this way, we will have constructed a class
of pairs 
\[(G, M) \mbox{ and } (\check{G}, \check{M}),\]
with   $M$ a  polarizable hyperspherical Hamiltonian space, and $\check{M}$
satisfying at least conditions \eqref{condaffine}, \eqref{condzero}, \eqref{condGm} of \S~\ref{condALL}.
\begin{itemize}
\item[-] We expect $\check{M}$ to be a hyperspherical Hamiltonian space, i.e., to satisfy all the conditions of \S~\ref{condALL}, but we cannot prove this.  (The local
conjecture of \S~\ref{Slocalconjecture} implies that it satisfies condition \eqref{condcoisotropic}.)

 \item[-]
We anticipate that much of the discussion that follows can
be generalized to avoid the ``polarizability'' condition.  Indeed,
the most ideal state of affairs would be that the duality
$$(G, M) \leftrightarrow (\check{G}, \check{M})$$
can be constructed on all hyperspherical Hamiltonian spaces,
satisfying suitable auxiliary conditions. 
A more precise expectation is formulated in Expectation \ref{anomaly expectation}. 
 \end{itemize}

\begin{remark}
For most of this paper, we will use $\check M$ on the ``spectral'' side of Langlands dualities. According to the convention \eqref{leftrightequation }, the group $\check G$ will act on the left on $\check M$. However, in this section, we use right actions for both $M$ and $\check M$; the translation to left actions is immediate based on the conventions of \S~\ref{leftrightconventions}.
\end{remark}
 
\label{sphericalclass}

\subsubsection{The hyperspherical data for $M$ and $\check{M}$} \label{hyp data checkM}
We now explain what we do in slightly more detail.
We fix a distinguished polarization, thus realizing $M$ as a twisted cotangent bundle $$M = T^*(X,\Psi).$$ In particular, by Proposition
\ref{hypersphericalspherical} $X$ is spherical and Borel stabilizers in the open orbit are connected. 
We also assume that $X$ admits a nowhere vanishing eigen-volume form with eigencharacter $\eta$, which we fix 
(see \S~\ref{eigencharacter});  the general situation will be reduced to that case (Lemma \ref{freecolors} (a), just as 
  in   \S~\ref{ssseigencharinocuous}).  
 
\textbf{For notational simplicity we will use $X$ as a symbol for the space together with the bundle $\Psi$ in what follows.} All invariants to be associated to $X$, such as its dual group, are, in general, \emph{different} from the invariants that would be associated to the space $X$ without the extra data.

Let us recall from Theorem \ref{thm:structure} that our general class of hyperspherical spaces all arise from triples (fixing $\mathfrak g \simeq \mathfrak g^*$, for notational simplicity),
\begin{equation} \label{basic1} (H \subset G, \mathfrak{sl}_2 \subset \mathfrak{g}, S \textrm{ a symplectic $H$-representation.})
\end{equation}
where the $\sl_2$ and $H$ commute.

 In the polarized case, as explained in \S~\ref{dp}, $S = S^+ \oplus S^-$ as $H$-representation and $M$
 is a twisted cotangent bundle over the space $X=S^+ \times^{HU} G$. 
 In more detail, the nilpotent element $f$
of the $\mathfrak{sl}_2$-triple defines a character of the Lie algebra $\mathfrak u$, normalized by $H$, 
and  integrates to an additive character $U\to \Ga$ whose kernel we will denote by $U'$. Thus,  we have a $\Ga$-bundle (or trivial bundle, if $f=0$) $\Psi := S^+ \times^{HU'} G$ over $X$;
then, $M$ is the  twisted cotangent bundle associated to $\Psi$.

We are going to construct the dual $\check{M}$ of $X$ from a corresponding datum on the side of $\check{G}$: \begin{equation} \label{basic2} (\check{G}_X \subset \check{G},  
\mathfrak{sl}_2 \to \check{\mathfrak{g}}, \mbox{$S_X$ a self-dual representation  of $\check{G}_X$})\end{equation}
where the $\sl_2$ and $\check{G}_X$ commute, which should be
understood as dual to \eqref{basic1}. We conjecture (Conjecture \ref{conjsymplectic}) that $S_X$ admits
a $\check G_X$-invariant symplectic form. This is valid in all cases we have checked, and we can prove it in the {\em strongly tempered} case when $\check G=\check G_X$ (Lemma \ref{vrhopair}).   Assuming the conjecture, we define the Hamiltonian $G \times \GGm$-space dual to $M$ to be
\begin{equation}\label{defM} \check{M} = \mbox{Whittaker induction (\S~\ref{Whittaker induction}) of $S_X$ from $(\check G_X, \sl_2)$}.
\end{equation}

\index{strongly tempered}
\begin{remark} 
\begin{itemize}
\item[(a)] (Notational warning): 
As remarked in \S \ref{dual to H}, although $\check{G}_X$ plays a role dual to $H$, 
we will \emph{not} denote that group by $\check H$, to avoid confusion with the Langlands dual group of $H$. 
\item[(b)] In discussing the twisted case it is often convenient to write in the notation of \ref{polarizationXstructure}:
\begin{equation} \label{Whittwisted notation} X = X_L\times^P G, \textrm{ with } X_L = S^+\times^H L\end{equation}
 where $L \subset P$ are the Levi and parabolic determined
 by the $\mathfrak{sl}_2$-triple,   and $X_L$ is an affine spherical $L$-variety. Moreover, the unipotent radical $U$ comes equipped with an additive character, which by Lemma \ref{lemma-generic} is generic.  
\end{itemize}
\end{remark}

 \subsubsection{Contents}
With reference to the general outline above, the contents of the remainder of the section are as follows.

 \begin{itemize}
  \item[-] \S~\ref{sec:invariants} constructs
 $\check{G}_X$ and $\sl_2$.
 \item[-] We will construct  the $\check{G}_X$-representation $S_X$ 
in \S~\ref{SXdef} (in a special case) as well as \S \ref{SX2sec} (general case).  
The reader might want to look directly at Definitions \ref{def:SX} and \ref{def:SX2} to get a sense
of what ingredients go into the pot: the highest weights are explicitly determined
in terms  of $B$-stable divisors on $X$. 

\item[-]   \S \ref{NTmov} is not used elsewhere in this section. It  discusses a certain $\check{G}_X$-representation
$V_X$ derived from $S_X$. It is $V_X$, rather than $S_X$,  which
is most visible in number theory. It
helps motivate how we found the formula for $S_X$. 
\item[-] 
  \S~\ref{parity}  examines the issue of ``parity'' (cf. \S \ref{analyticarithmetic}). 

\item[-]   \S~\ref{regnilp}  considers the image of the moment map for $\check{M}$; 
this is useful in our study of rationality issues. 
\item[-]  
  \S~\ref{CheckMRat}  studies 
 certain issues of rationality, both for $M$ and $\check{M}$: 
  \begin{itemize}
  \item[-] When $M$ is defined over a finite field, we will, in favorable cases, endow $\check{M}$ with a Frobenius action. 
  \item[-] When $M=T^* X$  (and $\Psi$ is trivial) we will  construct a preferred ``split'' form of $\check{M}$ even over a field $\kk$ that is not algebraically closed. 
  \end{itemize} 
 \end{itemize}
 
 The section has a quite liberal sprinkling of examples which, we hope, will help the reader digest the general constructions.

\subsection{The dual group of a spherical variety} \label{sec:invariants}
 
We will start by recalling some notions from the structure theory of spherical varieties. In particular,
we shall begin by describing the dual group, assuming at first that $f=0$ (i.e., no Whittaker induction).
$X$ will therefore be a spherical $G$-variety over the algebraically closed field $\FF$
whose generic $B$-stabilizers are connected (cf. 
Proposition \ref{hypersphericalspherical})
and admitting a $G$-eigenmeasure (\S \ref{eigencharacter});
cf. \S \ref{ssseigencharinocuous} as to why this assumption is inocuous.
 
If $B\subset G$ is a Borel subgroup, the spherical variety $X$ has an open $B$-orbit $ X^\circ$.  Let $P(X)\supset B$ \index{$P(X)$} be the stabilizer of $ X^\circ$, and $U(X)$ its unipotent radical. It is known \cite{KnMotion} that $U(X)$ acts freely on $ X^\circ$, and that the Levi quotient $L(X)$ acts on $ X^\circ/U(X)$ through a faithful action of a torus quotient $L(X)\twoheadrightarrow A_X$; this torus is the \emph{universal Cartan} of the spherical variety $X$. 
More precisely, there is a choice of Levi subgroup $L(X)\hookrightarrow P(X)$ such that we have an isomorphism of $P(X)$-spaces
\begin{equation}\label{openorbit}
 X^\circ\simeq T_X\times U(X),
\end{equation}
where $T_X$ is a torsor for $A_X$. 

The existence of an eigenmeasure with eigencharacter $\eta$ means that 
\begin{equation}\label{rhorho1}
\eta + 2 \rho  - 2 \rho_{L(X)} \in X^*(A_X),
\end{equation} 
where we use additive notation for the group of characters of $A$, $2\rho$ denotes the sum of positive roots of $G$, and similarly for $2\rho_{L(X)}$ for the Levi $L(X)$. Indeed, the difference  $2 \rho  - 2 \rho_{L(X)}$ is the modular character of the parabolic $P(X)$, so the condition above is the condition for the existence of $P(X)$-eigenmeasure with eigencharacter $\eta$ on the open orbit \eqref{openorbit}. 
  To see this we recall from \S \ref{Htrans} that the restriction of $\eta$ to a  point stabilizer  
  is simply the determinant of its right adjoint action on the tangent space at that point,
  which gives the {\em inverse} of the modular character.

For any character $\chi: A_X\to \mathbb G_m$, there is a unique up to scalar $(P(X),\chi)$-eigenfunction $f_\chi$ on $ X^\circ$, whose logarithmic differential $\frac{d \log f_\chi}{f_\chi}$ defines a section of the cotangent bundle $T^*X$ that is independent of the choice of $f_\chi$. Taking linear combinations of those (over $F$), we obtain a map 
\[ \mathfrak a_X^* \times  X^\circ \to T^* X^\circ,\]
where $\mathfrak a_X^*$ is the character group of $A_X$ tensored with $F$. 
Allowing now the parabolic $P(X)$ to vary in its conjugacy class $\mathcal P$, we obtain a $G$-equivariant map
\[ \mathfrak a_X^* \times (X \times \mathcal P)^\circ \to T^*X,\]
where $(X \times \mathcal P)^\circ$ denotes pairs $(x,P)$ with $x$ in the open $P$-orbit.

Knop shows that the image of this map is dense, and descends to an isomorphism 
\[ T^*X\sslash G \xrightarrow\sim \mathfrak c_X^*,\]
where $\mathfrak c_X^* = \mathfrak a_X^*\sslash W_X$ for a reflection group $W_X$, the \emph{little Weyl group} of $X$, acting on $\mathfrak a_X^*$. 
This $\mathfrak c_X^*$ is the spectrum of the integral closure of $F[\mathfrak g^*]^G$ in $F(T^*X)$, mentioned in \S~\ref{ssconditions} and denoted there by $\mathfrak c_M^*$. 

Moreover, Knop shows that this Weyl group $W_X$ is the same as the one constructed by Brion \cite{Brion-generalisation-symetriques}, a fundamental chamber of which on $\mathfrak a_{X,\mathbb R}:= X_*(A_X)\otimes_{\mathbb Z} \mathbb R$ is the \emph{cone of $G$-invariant valuations} $\mathcal V_X\subset \mathfrak a_{X,\mathbb R}$, where a $G$-invariant valuation on the function field $F(X)$ is considered as an element in the dual of the character group $X^*(A_X)$ by restriction to the Borel-eigenfunctions $F(X)^{(B)}$. (It is known that this restriction completely identifies a $G$-invariant valuation, see \cite{KnLV}.)

The data above play a role in the construction of the \emph{dual group} of $X$, which we will think here as 
a reductive subgroup $\check G_X \subset \check G$, with a canonical maximal torus $\check A_X\subset \check A$ dual to $A\to A_X$. The cone $\mathcal V_X$ of invariant valuations turns out to be the \emph{negative} Weyl chamber for $\check G_X$, i.e., the one corresponding to a Borel subgroup \emph{opposite} to $\check G_X\cap \check B$, where $\check B$ is the distinguished Borel subgroup of $\check G$. 
The group $\check G_X$  was constructed by Knop and Schalke \cite{KnSch}, and we expect it to be the same as the one constructed via the Tannakian formalism by Gaitsgory and Nadler \cite{GaitsgoryNadler}.  At this point, we will consider $\check G_X$ as a subgroup of $\check G$, unique up to conjugation by the maximal torus $\check A$; later, in \S~\ref{simple-M}, we will describe a precise choice of conjugate, compatible with the pinning of $\check G$, in order to define the \emph{$L$-group} of the spherical variety. 

\begin{remark}\label{notypeN}
For clarity, we emphasize that condition \eqref{condtorsion} of \ref{condALL} ensures that the group denoted by $\check G_X$ in \cite{SV} is a subgroup of $\check G$, and does not differ from the one of Gaitsgory--Nadler. Moreover, the same condition ensures that $X$ does not have any ``spherical roots of type $N$,'' in the language of \cite[\S~3.1]{SV}
which will be recalled in \S \ref{subsubsec:colors}. 
Indeed, if $X$ has spherical roots of type $N$, which by definition means that there is a simple root $\alpha$ such that the corresponding $\PGL_2$-variety $X^\circ P_\alpha/\mathcal R(P_\alpha)$ is isomorphic to $\mathcal N(\Gm)\backslash\PGL_2$ (see also \S~\ref{SXdef}), then generic stabilizers for the Borel orbit are disconnected, contradicting Proposition
\ref{isspherical}.
Spherical varieties with roots of type $N$ (a standard example being $\SO_n\backslash\SL_n$ for $n\ge 3$) are {\em not} contained in our conjectural duality of Hamiltonian spaces, as we expect their Hamiltonian dual to be a stack, rather than a smooth variety.   This is an issue that needs to be understood in future work.
\end{remark}

As was shown in \cite[Proposition 9.10]{KnSch} (see also \cite[\S~3.6]{SV}), the embedding $\check G_X\hookrightarrow \check G$ commutes with a principal $\SL_2$ into the standard Levi $\check L(X)$ \index{$\check L(X)$} dual to $P(X)$:
\begin{equation}\label{ArthurSL2}
 \iota: \check G_X \times \SL_2 \to \check G.
\end{equation}

In particular, we observe for later use that the ``$h$'' of the associated $\mathfrak{sl}_2 \rightarrow \check{\mathfrak{g}}$
is given by
\begin{equation} \label{2rhoreminder} h = 2 \rho_{L(X)},\end{equation}
the sum of positive roots for $L(X)$, considered as a coweight for $\check{G}$;
in the notation of \S \ref{uuplus}, this is the differential of the cocharacter denoted $\varpi$. 

\subsubsection{The dual group in the case of nontrivial $\Psi$} \label{dual-Whittaker}
In the case of Whittaker induction, it was explained in \cite[\S~2.6]{SV} how to attach a dual group, that \emph{differs} from the dual group of the space $X$, considered without the Whittaker character. 
Namely, if $(X,\Psi)$ is as in \eqref{Whittakerinduction}, we have the map \eqref{ArthurSL2} for the dual group of $X_L$, 
\[  \iota_L: \check G_{X_L} \times \SL_2 \to \check L.\]

We will consider the abstract based root system of $L$ as a subset of the abstract based root system of $G$ via the \emph{opposite} of the parabolic $P$ from which it is Whittaker-induced,
let $\check G_{X,\Psi}$ be the reductive subgroup of $\check G$ generated by $\check G_{X_L}$ and all the simple root spaces corresponding to $\pm \alpha$, for $\alpha \in \Delta\smallsetminus \Delta_L$.

\begin{proposition}\label{prop:psiroots}
	The subgroup of $\check G$ generated by $\check G_{X_L}$ and all the simple root spaces corresponding to $\pm \alpha$, for $\alpha \in \Delta\smallsetminus \Delta_L$ is a reductive subgroup $\check G_{X,\Psi}$ with the same Cartan as $\check G_{X_L}$,  and set of simple roots the union of the simple roots of $\check G_{X_L}$ and the set $\Delta\smallsetminus\Delta_L$. The subgroup $\check G_{X,\Psi}$ centralizes the image of $\SL_2$ under $\iota_L$.
\end{proposition}

\begin{proof}
	Let $\Delta_{X_L}$ be the set of simple roots of $X_L$, and $\Delta_{X,\Psi}=\Delta_{X_L}\cup \Delta\smallsetminus\Delta_L$. We will use a check ($\check \Delta$) to denote corresponding sets of coroots. It suffices to show that $\Delta_{X,\Psi}$, together with the root lattice of $A_X=A_{X_L}$ gives rise to a root datum, and that the corresponding reductive algebraic group embeds into the centralizer of $\SL_2$ in $\check G$, extending the embedding of $\check G_{X_L}$.
	
	First of all, we notice that any $\alpha\in \Delta\smallsetminus\Delta_L$ belongs to the character group of $A_X$. The proof is the same as in \cite[Proposition 2.6.2]{SV}; we outline the argument: First of all, by construction and \cite[Lemma 2.6.1]{SV}, the (additive) character by which the subgroup $U$ acts on the fiber of $\Psi\to X$ is nontrivial on each simple root subgroup $U_{-\alpha}$. Now, the Cartan $A$ of $G$ acts on the Lie algebra of $U_{-\alpha}$ via the character $-\alpha$, and the kernel of the map $A\to A_X$ stabilizes the aforementioned additive character; therefore, $\alpha$ has to be trivial on the kernel, hence factors through $A_X$. 
	
	The rest of the argument follows the construction of the dual group of a spherical variety by Knop and Schalke \cite{KnSch}. Knop and Schalke construct the dual group $\check G_{X_L}\subset \check L$ by a certain process of ``folding'' on the root datum of a full-rank subgroup $\hat G_{X_L}\subset\check L$ (i.e., $\hat G_{X_L}$ contains the full dual Cartan $\hat A$). If $\check\Delta_{X_L,\rm{as}}$ denotes the set of simple roots of $\hat G_{X_L}$ (the index stands for ``associated'' roots), the set $\check\Delta_{X,\rm{as}}=\check\Delta_{X_L,\rm{as}}\cup \check\Delta\smallsetminus\check\Delta_L$ forms the basis of an ``additively closed'' root subsystem by the criterion of \cite[Lemma 3.3]{KnSch}, hence, as in Theorem 7.3 of \emph{op.cit.}, corresponds to a full-rank subgroup $\hat G_{X,\Psi}$ of $\check G$. Then, as in Lemma 7.6 of \emph{op.cit.}, there is a ``folding'' involution $s$ which corresponds to the desired subgroup $\check G_{X,\Psi}$: It is obtained by trivially extending the folding involution of $\check\Delta_{X_L,\rm{as}}$ to $\check\Delta_{X,\rm{as}}$; that is, the involution fixes all elements of $\check\Delta\smallsetminus\check\Delta_L$. The verification of the ``folding'' property follows as in \emph{op.cit.}, namely, the only nontrivial property to check is that, for all $\beta\in \Delta\smallsetminus\Delta_L$, and all $\check\alpha\in \check\Delta_{X_L,\rm{as}}$, we have $\langle \check\alpha - {^s\check\alpha}, \beta\rangle = 0$, and this follows from Lemma 6.4 of \emph{op.cit}. As in Theorem 7.7 of \emph{op.cit.}, this implies that the embedding of $\check G_{X_L}$ into $\check L$ extends to an embedding of $\check G_{X,\Psi}$ into $\check G$. 
	
		By the ``folding'' construction of Knop--Schalke, then, the simple roots of $\check{G}_{X, \Psi}$ are exactly the simple roots of $\check G_{X_L}$, together with $\check{\Delta}\smallsetminus\check{\Delta}_L$; in particular, as a subgroup of $\check G$, $\check G_{X,\Psi}$ is generated by $\check G_{X_L}$ and the simple root spaces corresponding to $\pm \alpha$, for $\alpha \in \Delta\smallsetminus \Delta_L$. There remains to show that this subgroup commutes with $\SL_2$, which, by \cite[Proposition 9.10]{KnSch}, boils down to showing that these simple root spaces centralize a certain subgroup that is denoted by $L^\wedge_{\mathcal S}$. As in the proof of Theorem 9.7 of \emph{op.cit.}, this boils down to Lemma 6.6, and eventually to the study of the rank-1 ``spherical variety'' corresponding to each $\alpha \in \Delta\smallsetminus \Delta_L$ -- but the situation in this case is particularly simple, because $\alpha$ is a simple root of $G$ and, at the same time, orthogonal to the simple roots of $P(X_L)=P(X)$ (since, by what we just proved, it belongs to the character group of $A_X$). Therefore, the image of the map $\SL_2\to\check G$ corresponding to $\alpha$ commutes with the Levi $\check L(X)$, and therefore with its subgroup $L^\wedge_{\mathcal S}$ and with the principal $\SL_2$ which appears in $\iota_L$.
\end{proof}

Under a slightly restrictive condition called the ``wavefront'' property of $X_L$, it was proven in \cite[Proposition 2.6.2]{SV} that the 
Weyl group of $\check G_{X,\Psi}$ is the one attached by Knop \cite{KnMotion} to the $\Ga$-bundle $\Psi$, viewed as a (nonspherical) $G$-space -- see also \cite[\S~5.4]{SaSpc}. It would be desirable to extend this result to the general case; if it does not hold, our definition of $\check G_X$ should be changed, for it to have the Weyl group defined by Knop. However, in this paper we will proceed with the definition above. 

 As remarked in \S \ref{hyp data checkM}, in this situation we will be using $X$ to denote not just the space but also the data of this $\Ga$-bundle, hence will be denoting this dual group and its Weyl group simply by $\check G_X$, $W_X$, etc., keeping in mind that this is different from the dual group of $X$ without the $\Ga$-torsor.

\subsection{The $\check{G}_X$-representation $S_X$ in the case of affine closures} \label{SXdef}

In the following subsections,  we give an ad hoc description of the $\check G_X$-representation $S_X$, which in all examples that we have considered matches results of \cite{SaSph,SaWang} in a sense to be described in \S \ref{ssPlancherel},  and is conjecturally symplectic (Conjecture \ref{conjsymplectic}).

It is likely that we are ``working too hard'' here -- for example, in the examples we have examined, all the weights of $S_X$ are minuscule, which greatly simplifies things. This might always be the case, for $X$ affine, spherical, and smooth; in principle, this could be checked ``by hand,'' based on the classification of such varieties by Knop and Van Steirteghem \cite{KnVS}. Until such simplifications are established, or the reader should take the general definition with a grain of salt, as a working hypothesis.

\subsubsection{The canonical affine open within $X$} 

\begin{lemma}\label{affineopen}
 If $X$ is an affine spherical $G$-variety, and $X^\bullet$ denotes the open $G$-orbit, the canonical map $\overline{X^\bullet}^\aff \to X$ is an open embedding; in particular, if $X$ is smooth, so is $\overline{X^\bullet}^\aff$.  
\end{lemma}
Here, $\overline{X^\bullet}^\aff := \Spec \mathbb F[X^\bullet]$ is the affine completion  of the homogeneous part of $X$, which we will call the \emph{canonical affine open} subset of $X$, and denote by $X^\can$. (The spherical condition implies that the coordinate rings $\mathbb F[X]$ is finitely generated, since $B$-eigenspaces are $1$-dimensional and the monoid of $B$-eigencharacters is finitely generated.)  Here, we assume that $X$ is untwisted, but we will extend the definition to the twisted case below.

\begin{proof} 
 Let $X_1 \subset X$ be the complement of all $G$-stable divisors. It is an affine, open, and normal (actually smooth, since $X$ is smooth) subvariety, where the complement of the open $G$-orbit is of codimension $\ge 2$. By normality, every regular function on $X^\bullet$ extends to $X_1$, identifying it with its affine closure. 
\end{proof}

In the twisted case, i.e., when $X$ is endowed with a $\Ga$-torsor $\Psi$, by \eqref{Whittakerinduction} we can write
$(X,\Psi) = \Ind_P^G (X_L,\Psi_L),$ with $X_L$ a smooth, affine, spherical $L$-variety, where $L$ is the Levi quotient of a parabolic $P$, and $\Psi_L$ is a $P$-equivariant $\Ga$-torsor defined by a generic character of the unipotent radical. In this case, we will define the canonical open subset of $X$ to be 
\begin{equation}\label{Xcantwisted} X^\can = \Ind_P^G (\overline{X_L^\bullet}^\aff)\subset X,
\end{equation}
 endowed with $\Psi^\can:=$ the pullback of the torsor $\Psi$ to it.

In the current section, we will give the definition of $S_X$
in a special case, namely, 
 when  $(X, \Psi)=(X^\can,\Psi^\can)$ and satisfies a certain additional  combinatorial condition,
namely, {\em freeness of the color monoid} (to be explained). 
As it will turn out, this condition is automatically satisfied under our smoothness assumptions, but we will see this when discussing the general case, in \S \ref{SX2sec}. For notational simplicity, we will only be referring to $X$, but the discussion will apply verbatim to the twisted case. 
 
\subsubsection{Colors}  \label{subsubsec:colors}
When $X= X^\can$, the $\check A_X$-space $S_X$ will be determined (up to isomorphism) by the \emph{colors} of $X$.

By definition, colors on a spherical variety $X$ are the irreducible $B$-stable divisors\footnote{If $\FF$ were not algebraically closed, these should
be considered geometrically, i.e., over the algebraic closure} that are not $G$-stable; when $X=X^\can$, these are all the $B$-stable divisors in $X^\bullet$. (The choice of a Borel subgroup is immaterial, since $B$-orbits on $X$ are the same as $G$-orbits on $X\times\mathcal B$, where $\mathcal B$ is the flag variety.)

Let $\Delta$ denote the set of simple roots of $G$, and $\Delta_{(X)}$ the subset of those belonging  to a Levi subalgebra of the parabolic $P(X)$. There is a crucial diagram:
$$ \xymatrix{  \mbox{colors}  \ar[d]^{\check v}  \ar[r]^{r \qquad \qquad} & \textrm{subsets of }\Delta\smallsetminus \Delta_{(X) } \\ 
  X_*(A_X)} $$
 The horizontal map associates to each color $D$ the set of $\alpha \in \Delta$
for which the associated parabolic $P_\alpha$ of semisimple rank one satisfies $DP_{\alpha} \supset  X^\circ$. 
The vertical map $\check v$ takes a color $D$ to the corresponding valuation, restricted to rational $B$-eigenfunctions $f_\chi$ on $X$; 
this gives as homomorphism $X^*(A_X) \rightarrow \Z$, 
or, what is the same a cocharacter 
$ \check v_D \in X_*(A_X).$ 
We will often abuse language and talk of the $\check v_D$'s as the ``colors'', considered as a multiset in $X_*(A_X)$.
 
For each $\alpha \in \Delta \smallsetminus \Delta_{(X)}$, the preimage has either size $1$ or $2$. 
More precisely, for every such $\alpha$,  taking the (geometric) quotient of $X^\circ P_\alpha$ by the radical $\mathcal R(P_\alpha)$ of $P_\alpha$, we obtain a homogeneous spherical variety $X_\alpha$ for $\PGL_2$, which, under our current assumptions can only be isomorphic (over the algebraic closure) to
one of the following:
\begin{itemize}
\item[(type U)]
 $RU\backslash \PGL_2$, where, $U\simeq \Ga$ is a unipotent subgroup, and $R$ is a finite subgroup in its normalizer.\footnote{
the case where $R\simeq \Gm$ is excluded, because $X^\bullet$ is quasiaffine.}

   If $X$ is equipped with a nontrivial affine bundle $\Psi$, this splits into two types: The $G$-equivariant $\Ga$-torsor $\Psi\to X$ 
    also admits a similar reduction $\Psi_\alpha \rightarrow X_{\alpha}$
    and -- following the terminology of \cite{SaSph}) -- we say that
    the root has type $U$ or type $(U, \psi)$ according to whether $\Psi_{\alpha}$ is trivial or not.  
   
 \item[(type N)]
 $\mathcal N(\Gm)\backslash \PGL_2$. This is  excluded by our assumptions -- these are the ``roots of type $N$,'' see Remark \ref{notypeN}.
 \item[(type $T$)] 
 $\Gm\backslash \PGL_2$
 \item[(type G)]
 $\PGL_2\backslash\PGL_2$.
 \end{itemize}
We will say that a simple root $\alpha$ is of a certain type, \index{type of a root} if the $\PGL_2$-variety above is of that type. 

We will now identify a distinguished subset of the colors, which we will call \emph{of even sphere type}  (for reasons that we will explain). The main representatives of those are the colors belonging to simple roots $\alpha$ of type $T$. For those roots, there are precisely 2 colors $D, D'$ contained in $X^\circ P_\alpha$
(these will also be called ``colors of type $T$.'') 
 Their valuations satisfy
\begin{equation}\label{coloreq}
 \check v_{D'} = -w_\alpha \check v_D,\,\,\, \mbox{ and } \check v_D + \check v_{D'} = \check\alpha,
\end{equation}
and in particular $\langle \check v_D,\alpha\rangle = \langle \check v_{D'},\alpha\rangle  = 1$; see \cite[\S~1.4]{Luna-typeA}. Here, in the first equality, $w_{\alpha}$
is regarded as an automorphism of $X_*(A_X)$ because $\alpha$ is itself
a spherical root, see {\em op.\ cit.} 
 These colors play (whose set is denoted by $\mathbf A_X$ in  \cite{Luna-typeA}) play an important role in the classification theory of spherical varieties. It can be shown that if a color $D$ is of type $T$, then \emph{every} simple root in its image under the horizontal map $r$ above is of type $T$. (This is an easy exercise, based on \cite[\S~1.3]{Luna-typeA}, and is left to the reader.)

Simple roots of type $T$ can be considered as a special case of the following phenomenon: a parabolic $P\subset G$, such that $X^\circ P/\mathcal R(P)$ is isomorphic to $\SO_{2n}\backslash \SO_{2n+1}$.
 \emph{We want to include here the case of the spherical variety $\SL_3\backslash G_2$}, which is isomorphic as a variety to $\SO_6\backslash\SO_7 = \Spin_6\backslash\Spin_7$ (via the embedding $G_2\hookrightarrow\Spin_7$) so we adopt the following formal definition: a standard parabolic $P$ is of ``even spherical type''  when the associated spherical variety
  $(P/R(P),X^{\circ}P/R(P))$ is isomorphic either to $ (\SO(2n+1), \SO(2n)\backslash \SO(2n+1))$
   or to $(G_2, \SL_3\backslash G_2).$ When $n\ge 2$, the geometry of colors is different (in fact, there is a unique color meeting $X^\circ P$ in those cases), but those varieties also need to be considered alongside roots of type $T$, because they contribute to the definition of the space $S_X$; in terms of number theory, their periods contribute an $L$-value at $\frac{1}{2}$.  
\begin{definition}\label{def-spheretype}
 \index{colors of even sphere type} The colors which meet $X^\circ P$, for some parabolic $P$   of even sphere type as above, will be called \emph{colors of even sphere type}. Sometimes, we will identify them with their valuations, i.e., we will say ``colors of even sphere type'' for the corresponding (multi)set of elements of $X_*(A_X)$. The corresponding spherical root (see Remark \ref{remark-spheretype} below) will be called a spherical root of even sphere type.
\end{definition}

\begin{remark}
The spherical varieties occurring in the definition are, by inspection,  the spherical varieties of rank $1$ whose ``associated $L$-value,'' according to the recipe of \S \ref{ssPlancherel} for the unramified Plancherel formula, contains a factor evaluated at $\frac{1}{2}$. We will recall more details about this recipe in \S \ref{ssPlancherel}. In a future paper, we will give a more conceptual interpretation of the special role that these colors play, at least for those of type $T$ when $\check G_X=\check G$; see \cite[Theorem 6.2.2]{SaICM} for a statement.\footnote{This was previously observed, in examples, by V.~Lafforgue (unpublished).}
\end{remark}

\begin{remark}\label{remark-spheretype}
 Except for the ``type $T$'' case ($\SO_2\backslash\SO_3$),  which was already discussed around \eqref{coloreq},  in all other cases there is a unique color $D$ in $X^\circ P$ -- see \cite[6.14]{SaSph}, and a similar diagram can be calculated when the action is restricted to $G_2$. Moreover, except in type $T$, that color satisfies 
 \begin{equation}\label{colorSeq}
\check v_D = \frac{\check\gamma}{2},
 \end{equation}
 where $\gamma$ is the spherical root and $\check\gamma$ is the associated coroot. The spherical root, in those cases, is the short root which is the sum of all simple roots in the case of $\SO_{2n+1}$, and the sum of the long simple root with twice the short simple root in the case of $G_2$. Those are orthogonal to all simple roots but one (the one furthest, in the Dynkin diagram, from the short root in $\SO_{2n+1}$; the short root, in the case of $G_2$), and all those roots that are orthogonal to $\gamma$ will be contained in the Levi of $P(X)$. See Example \ref{ex:5sphere} below.
 
 However, we hasten to clarify that the unique associated color in this case is just a placeholder for a $B$-orbit of larger codimension, or rather for a ``formal difference'' of such $B$-orbits -- see the discussion of Galois actions preceding Definition \ref{simpleaction}. This is a pedantic detail that the reader can safely ignore in most cases, but it is necessary in order to get the Galois actions right. 
\end{remark}

Let us denote by $\mathcal C_X$ the set of colors of $X$ of even sphere type, considered as a multiset of elements of $X_*(A_X)$.
Let us first consider the case where the valuations $\check v_D \in \mathcal C_X$ freely generate a direct summand of $X^*(\check A_X)$. This need not be the case; indeed, 
it is not necessary that the images of distinct colors under $\check v_D$ are distinct, as the example of $\Gm\backslash\PGL_2$ shows, and  is also not necessary that these valuations generate the subgroup of $X_*(A_X)$ which lies in their $\mathbb Q$-span, as the example of $\Gm\backslash (\Gm\times\PGL_2)$ shows.  
But we will explain how to reduce the general case to this one in \S \ref{SX2sec}.

 Let $\mathcal D_X\subset X_*(A_X)$ denote the \emph{dominant} $W_X$-translates of elements of $\mathcal C_X$ (just as a subset, without multiplicity), and let $\mathcal D_X^{\text{max}}$ denote the subset of maximal elements of $\mathcal D_X$ with respect to the standard coroot ordering: $\check v_1\ge \check v_2 \iff \check v_1-\check v_2$ can be written as a linear combination of positive coroots of $\check G_X$ with \emph{integral}, non-negative coefficients.

\begin{definition}\label{def:SX}
 When $(X,\Psi) = (X^\can,\Psi^\can)$ and the elements of $\mathcal C_X$ freely generate a direct summand of $X_*(A_X)$, we let $S_X$ be the representation of $\check G_X$ with highest weights $\mathcal D_X^{\text{max}}$.
\end{definition}

\begin{remark}\label{remark-dmax}

\begin{enumerate}
 \item We do not know if $\mathcal D_X^{\text{max}}$ is actually ever strictly smaller than $\mathcal D_X$. In fact, it may well be that in the smooth affine case all elements of $\mathcal D_X$ are minuscule. In that case, the weights of the representation $S_X$ are \emph{precisely} the $W_X$-translates of colors. We do not currently know an example with non-minuscule weights, but we also do not know how to prove that they should always be minuscule. In \cite[Corollary 7.3.4]{SaWang}, it was asserted that this is the case when $X = H\backslash G$ is affine, but the claimed proof is incomplete and leads to a weaker conclusion.

\item In Lemma \ref{freecolors} we will see that there is always a finite cover $Y$ of $X$ which satisfies the condition about the colors. This means that we could omit it from our definitions (at the expense of working with \emph{multisets}, not sets of valuations), if we knew that $\mathcal D_X^{\text{max}}=\mathcal D_X$ -- in principle, this can change when passing from $X$ to $Y$. 

Hence, we may be ``working too hard:'' it is entirely possible that, when $X = \overline{X^\bullet}^\aff$ is smooth, $S_X$ can be defined simply as the representation with highest weights the $W_X$-translates of the colors -- but note that, when colors give the same valuations, they have to be included with the corresponding multiplicity, as Example \ref{ex:Hecke} below shows.
\end{enumerate}

\end{remark}

\begin{example} \label{ex:Hecke}
We consider the case of $X = \Gm \backslash \GL_2$
where $\Gm$ is embedded in $\GL_2$ via $e_1^{\vee}$ in standard notation.  
We have $$X^{\circ}/U \simeq \Gm^2, \ \begin{pmatrix} a&b \\ c & d \end{pmatrix}
\mapsto (c, a^{-1} \det )$$ and the corresponding map $B \rightarrow \Gm^2$
is given by $(e_1, e_2)$. In particular,
the function $f_{\chi}$ of \S \ref{sec:invariants} corresponding to $\chi = x e_1 + ye_2 \in X^*(B)$
is given by $c^{x} a^{-y} (\det)^y$. 
 The two $B$-stable divisors are defined by $c=0$ and
$a=0$,  and the valuations correspondingly are given by
$(x,y) \mapsto x$ and $(x,y) \mapsto -y$, i.e., 
 $ \mathcal{C}_X =  \{  e_1^{\vee},  -e_2^{\vee} \}$
It follows from this that
$$S_X = \mathrm{std} \oplus \mathrm{std}^*$$
as a representation of the dual $\GL_2$. 

Note that, if we replace $\GL_2$ by $\PGL_2$, both colors give valuation $\frac{\check\alpha}{2}$, failing to satisfy the ``free generation'' condition of Definition \ref{def:SX}. This case will be treated in the next subsection by considering the cover $\Gm \backslash \GL_2 \to \Gm \backslash \PGL_2$, resulting in essentially the same answer for $S_X$ (but, now, as a representation of the dual group $\SL_2$ of $\Gm\backslash\PGL_2$). 
\end{example}

 \begin{example}
A particularly trivial case is the Whittaker case $X=U \backslash G$, together with the affine
bundle $\Psi$ arising from a nondenerate (generic) homomorphism $U \rightarrow \mathbb{G}_a$. Here 
the set $\mathcal{D}_X$ is empty, and $S_X$ is trivial.  Accordingly
we in fact have $\check{M}$ trivial. 
\end{example}

\begin{example} \label{ex:SX3} 
 Consider ``Hecke periods for $\GL_n$.''  Let $V_n = V_1\oplus V_1'\oplus V_{n-2}$ be an $n$-dimensional vector space with a decomposition of the indicated dimensions, let $U\subset G:= \GL(V_n)$ be the unipotent radical of the  parabolic stabilizing $V_{n-2}$, and let $\psi$ be a generic character of $U$ stabilized by the subgroup $\Gm\simeq \GL(V_1)$ of $G$. This  gives some Whittaker-type induction, in the sense of \cite[\S~2.6]{SV}, of the ``Hecke period'' $X_L =\GL_1\backslash\GL_2$, for which the dual Hamiltonian space is $\check M_L = S_{X_L}=$ the cotangent space of the standard representation of $\GL_2$. For general $n\ge 2$, the corresponding period is known to represent the standard representation of $\check G$ -- hence, we would like to say that $\check M = S_X=$ the cotangent space of the standard representation of $\GL_n$.

As a warning, this $X$ is {\em not} a distinguished polarization of $M=T^*X$;
we will see the distinguished polarization in Example \ref{ex:SX3p} below.
From the automorphic point of view, both represent the same period (cf. \cite[(4.1.1)]{JPSSGL3} and following discussion). 
 
\end{example}

\begin{example} \label{ex:SX3p}  (``How to construct Rankin--Selberg integrals''): We consider the case $X=\mathbb{A}^n$ under the right action $G=\GL_n$. 

The only $B$-stable divisor is defined by $x_1=0$ and its stabilizer
is the parabolic $P(X) = P_{1,n-1}$ (in usual notation: upper triangular with $1$ and $n-1$ blocks). 
The only root $\alpha \in \Delta \backslash \Delta_{(X)}$ is the
simple root $e_1-e_2$, and it has type $U$. 
Therefore $S_X$ is the trivial representation of $\check{G}_X = \Gm$. 
The associated $\SL_2$ into $\check{G}$ is principal
for the Levi of type $(1, n-1)$; 
the morphism 
$\Gm \times \SL_2 \rightarrow \check{G} = \GL_n$
factors through $\GL_1 \times \GL_{n-1}$
in the evident way. 

For $n$ even the dual space is
$$ \check{M} = T^* (\check{X}, \check{\Psi}), \check{X} =  \GL_n /  \Gm (U, \psi),$$
for suitable $(U, \psi)$ derived from the $\SL_2$. 
(For $n$ odd it is does not admit a distinguished polarization, but it can be polarized as in Remark \ref{dp-generalized}.)
In the automorphic situation, this $(\check{X}, \check{\Psi})$ indexes a 
(a slight variation -- see below) of the standard integral representation of the standard $L$-function for $\GL_n$,
as in the example just studied; thus
we have {\em derived} this integral representation from our general recipe.
Our general proposal is that all Rankin--Selberg integrals representing
$V$ {\em for which $T^* V$ is hyperspherical} can be derived from first principles this way. 
Note, however, that many examples of interesting Rankin--Selberg integrals fall outside this class,
and we will discuss them elsewhere (see e.g.\ \cite{ChenVenkatesh}).

To clarify the relation between this and \cite{JPSSGL3}  
consider the example for $n=4$: $\check{X}$ is the quotient of $\GL_4$ by the subgroup
 \begin{equation} \label{MM3} \begin{pmatrix} * & 0 & 0 & * \\ * & 1 & x & * \\ 0 & 0 & 1 & y \\ 0 & 0 & 0 & 1  \end{pmatrix}  \end{equation}
 and the character to $\mathbb{G}_a$ defining $\Psi$ is given by $x+y$. 
 The integral of \cite{JPSSGL3} corresponds to $\check{X}'$ wherein the $(2,1)$-entry
 has been replaced by the $(1,3)$. The corresponding spaces are not $\GL_4$-equivalent,
 but their cotangent bundles are equivalent, as follows from the fact
 that the two spaces are both inductions from  a certain Heisenberg group. 
 In other words, from the point of view of this paper, the integral of \cite{JPSSGL3} 
 is not the ``standard'' one that one would derive from a hyperspherical presentation,
 but rather something equivalent to it.

   \end{example}

\begin{example} \label{ex:5sphere}
 (See \cite[\S~6.9]{SaSph}.) We consider the case of $X = \SO_4\backslash \SO_5$. Then $\check G_X = \SL_2$ with spherical root $\gamma = \alpha+\beta =$ the sum of the two simple roots of $G$, while $P(X)$ is the parabolic that has the short root in its Levi;
 correspondingly, 
 $$ \check{G}_X \times \SL_2 \rightarrow \check{G} = \mathrm{Sp}_4$$
 is the direct sum of standard representations on either factor.
 
If we present $X$ as the sphere $q(x_1, \dots, x_5)=1$
inside the $5$-dimensional quadratic space with form given in an orthogonal basis $\nu_i$ by
$q(\sum x_i \nu_i) = x_1x_5+x_2x_4+x_3^2$, 
and we take the Borel subgroup to be the stabilizer of the isotropic flag
 $\nu_1 \subset \langle \nu_1, \nu_2 \rangle \subset \langle \nu_1, \nu_2, \nu_3 \rangle \subset \langle \nu_1, \nu_2, \nu_3, \nu_4 \rangle \subset V$, 
 with roots $e_1, e_2$ corresponding to the action on $\nu_5, \nu_4$ respectively. 
 The positive simple roots are $e_1-e_2$ and $e_2$. 
 Then $x_5$ gives a $U$-invariant function on $X$ and indeed 
$$X^{\circ}/U \simeq \Gm, (x_i) \mapsto x_5.$$
The function  $f_{\chi}$ of \S \ref{sec:invariants} corresponding to $\chi = a e_1 \in \langle e_1 \rangle = X^*(A_X) \subset X^*(B)$
is given by $x_5^a$ and its order of vanishing along the color $x_5=0$
is given by $ae_1 \mapsto a$.  This gives an element of $X_*(A_X)$
whose pairing with the spherical root is equal to $1$. Considered
inside  $ X^*(A_X^{\vee})$ this is a weight for $\check{G}_X=\SL_2$
whose pairing with the coroot is $1$, i.e., the highest weight of the standard representation. 
Correspondingly we find
$$S_X = \mbox{standard repesentation of $\check{G}_X=\SL_2$.} $$

\end{example}

\subsubsection{Symplecticity of $S_X$}

For the rest of this subsection, we assume that the assumptions of Definition \ref{def:SX} are satisfied. Note that $S_X$ is self-dual, by construction, since its multiset of weights is invariant under $\{\pm 1\}$; this follows from \eqref{coloreq}. 
If it admits a symplectic structure,  that structure is unique up to isomorphism, in the sense that any other is obtained by applying a $\check G_X$-automorphism of $S_X$. 
This is a general fact about representations of reductive groups over algebraically closed fields, which we record for clarity:  

\begin{lemma} \label{uniqueness of symplectic form} 
Let $W$ be an irreducible representation of the reductive group $\Xi$
over the algebraically closed field $\kk$, such that $W$ is abstractly isomorphic to its dual representation $W^*$. Any two $\Xi$-invariant symplectic forms on $\Xi$
differ by an element of $\mathrm{Aut}_{\Xi}(W)$. 
\end{lemma}

\proof
We reduce to the following basic cases: $W = (\sigma \oplus \sigma^*) \otimes E$,
where $\sigma$ is irreducible and not self-dual, and $E$ has trivial $\Xi$ action;
and $W = \sigma \otimes F$, where $\sigma$ is irreducible self-dual and $F$ has trivial $\Xi$-action. 
In the latter case, symplectic forms on $W$ correspond to nondegenerate symmetric or skew-symmetric
pairings on $F$, and these are all conjugate under $\GL(F) \subset \Aut_{\Xi}(W)$. 
In the former case, symplectic forms on $W$ correspond to (not necessarily symmetric) perfect pairings
$E \times E \rightarrow \kk$, all of which are conjugate under $\GL(E) \times \GL(E) \subset \Aut_{\Xi}(W)$. 
\qed

Therefore, for our purposes, the only question regards existence. Let us denote by $\mathfrak B_X$ the multiset of weights of $S_X$.

\begin{conjecture}  \label{conjsymplectic}
The $\check G_X$-representation $S_X$ is symplectic. Moreover, the multiset $\mathfrak B_X$ has the following properties: 
\begin{enumerate}
 \item The valuations $\check v_D$ associated to colors of even sphere type appear with multiplicity one\footnote{For the multiplicity statement, we are relying on the assumptions of Definition \ref{def:SX}; when different colors induce the same valuation, as in Remark \ref{remark-dmax}.(2), the multiplicities should be adjusted accordingly.} in $\mathfrak B_X$; hence, we can consider $\mathcal C_X$ as a subset of $\mathcal B_X$.
 \item There is a decomposition $\mathfrak B_X = \mathfrak B_X^+ \sqcup \mathfrak B_X^-$ such that the weights of $\mathfrak B_X^+$ lie in $X_*(A_X)^D$ and the weights of $\mathfrak B_X^-$ lie in $-X_*(A_X)^D$. Here $X_*(A_X)^D\subset X_*(A_X)$ is the monoid generated by \footnote{There is redundancy in this generating set, as by \eqref{coloreq} we could omit the $\check\gamma$'s of even sphere type.} 
  valuations $\check v_D$ attached to colors $D$ of even sphere type and the  simple roots $\check\gamma$ of $\check G_X$. 
 \item For every simple root $\check\gamma$ of $\check G_X$, the negative root space $\check{\mathfrak g}_{X,- \check\gamma}$  maps the weight spaces with weights in $\mathfrak B_X^+$ to each other, except when $\gamma$ is a spherical root of even sphere type and $\mathcal C_X^\gamma\subset \mathfrak B_X^+$ is the set of valuations of the corresponding colors of even sphere type, in which case the $\mathcal C_X^\gamma$-weight spaces are mapped onto the $(-\mathcal C_X^\gamma)$-weight spaces by $\check{\mathfrak g}_{X,- \check\gamma}$:
 \begin{equation}\label{switchcolors}
  \check{\mathfrak g}_{X,- \check\gamma} (S_X)_{\mathcal C_X^\gamma} = (S_X)_{-\mathcal C_X^\gamma}.
 \end{equation}
\end{enumerate}
\end{conjecture}

\begin{lemma} \label{vrhopair}
If $\check G=\check G_X$ then Conjecture \ref{conjsymplectic} holds, except perhaps for the multiplicity statement when $\mathcal D_X^{\text{max}}\ne \mathcal D_X$.
\end{lemma}

\begin{proof}
 Since $S_X$ self-dual, for the symplectic property it is enough to show that any \emph{irreducible} self-dual subrepresentation of $S_X$ is symplectic. This, in turn, is equivalent to the condition that the central element $(-1)^{2\rho}\in \check G$ acts by $-1$ on it. It is enough to show that $\langle \check v_D,2\rho\rangle$ is odd for all colors $D$ whose valuation appears in the subrepresentation under consideration, and this follows from the existence of a $G$-eigenmeasure on $X$, as follows: 
  
  When $\check G_X=\check G$, for any $\alpha\in \Delta$ the open $P_\alpha$-orbit $X^\circ P_\alpha$ is isomorphic to $\Gm\backslash P_\alpha$, where $\Gm$ is embedded into $P_\alpha$ (up to conjugacy) via the cocharacter $\check v_D$ corresponding to a color in $X^\circ P_\alpha$. (Which color is chosen does not matter, because of \eqref{coloreq}.) For there to exist an eigen-volume form on $\Gm\backslash P_\alpha$ with eigencharacter $\eta$, arguing as in \eqref{rhorho1} gives 
  \[\langle \check v_D,  \eta + 2\rho_{P_\alpha}\rangle = 0,\] where $2\rho_{P_\alpha}$ is the sum of roots in the unipotent radical of $P_\alpha$ (i.e., its modular character). In other words, 
  \begin{equation} \label{modularchararg} 
      \langle \check v_D, \eta + 2\rho\rangle = \langle \check v_D ,\alpha\rangle =1, 
  \end{equation}
  by \eqref{coloreq}. 
  
 On the other hand, since $\eta$ is a character of $G$, $\langle \check v_D,\eta\rangle$ must be zero, if $\check v_D$ is to be the weight of an \emph{irreducible} self-dual representation.  Thus, $\langle  \check v_D, 2\rho \rangle = 1$ for all such colors.

 When $\mathcal D_X^{\text{max}} = \mathcal D_X$, the statement on the multiplicities of elements of $\mathcal C_X$ in $\mathfrak B_X$ is obvious by construction; namely, they appear with multiplicity one, since they are $W_X$-translates of the dominant weights, which also appear with multiplicity one (under our current assumptions). 
 
 Ignoring multiplicities, our multiset $\mathfrak B_X$ is the same as the multiset of weights of a crystal (in the sense of Kashiwara), denoted by the same symbol in \cite[7.1.4]{SaWang}. The claims about its decomposition into $\mathfrak B_X^+ \sqcup \mathfrak B_X^-$ (and the action of $  \check{\mathfrak g}_{X,- \check\gamma}$) are then contained in \cite[Theorem 7.1.9]{SaWang}. Note that $\mathfrak B_X^+$ can be directly characterized by \eqref{modularchararg} by the property that its elements $\theta$ satisfy $      \langle \theta, \eta + 2\rho\rangle >0$. 

\end{proof}

We mention here that all smooth affine spherical varieties have been classified (``modulo center'') by Knop and Van Steirteghem \cite{KnVS}; thus, it is possible to check Conjecture \ref{conjsymplectic} ``by hand'' -- but we haven't done so. Our
local conjecture gives a conceptual reason to expect symplecticity, as we discuss in \S \ref{Poisson local conj}, \S \ref{Poisson to loop} and more at length in \S \ref{spectral-factorization}.

\subsection{The $\check{G}_X$-representation $S_X$ in the general case } \label{SX2sec}

We are now going to  define $S_X$ in the general case. We will start with the untwisted case, $M = T^*X$, where $S_X$ is described in Definition \ref{def:SX2}, and will treat the twisted (Whittaker-induced) case, by a straightforward generalization of this definition, in \S~\ref{SX3sec}.

We  first use  the following statement:

\begin{lemma}\label{freecolors}
 Let $X^{\bullet}$ be  the open $G$-orbit on the smooth spherical affine $G$-variety $X$. There is a central extension $G'\to G$, whose kernel is a torus $T$, and a homogeneous $G'$-variety $Y$ which is (equivariantly) a $T$-torsor $Y\to X$, with the following properties:
\begin{itemize}
\item  the valuations $\check v_D$ associated to colors of $Y$ are distinct, and they freely generate a direct summand of the group of cocharacters of $A_Y$.  
\item $Y$ admits a $G'$-eigenmeasure. 
\end{itemize}
\end{lemma}

\proof 
This is  \cite[Lemma 5.3.3]{SaWang} except for the statement concerning eigenmeasure;
that follows, taking a further central extension if necessary, as in \S \ref{ssseigencharinocuous}.
  \qed 

\begin{remark}\label{Yaffsmooth} 
 One can show that the map $\overline{Y}^{\text{aff}}\to \overline{X^\bullet}^{\text{aff}}$ between affine completions is also a $T$-torsor; in particular, $\overline{Y}^{\text{aff}}$ is also smooth. This is not stictly needed for what follows, but since we have so far been associating invariants to affine closures, it is reassuring to know that we are not leaving the world of smooth varieties.
\end{remark}

Following the notation of the Lemma, the colors of $Y$ are in bijection with colors of $X$ under the quotient map. 
 Also, the dual group of $X$ is a subgroup of the dual group of $Y$ (with quotient equal to the dual of the torus $T$ acting on $Y\to X^\bullet$), and therefore, assuming Conjecture \ref{conjsymplectic}, we can will consider the representation\footnote{Note that the definition of $S_X$ in the previous paragraph only used the colors in $X^\bullet$, and thus makes sense for $Y$, as well, even if it is not affine.} $S_Y$ as a representation of $\check G_X$. To make sure that it is independent of $Y$ (which could, in principle, happen if valuations whose dominant $W_X$-translates are comparable for one cover, but not for another),\footnote{We have no examples where this happens.} we note that for a pair $Y_1, Y_2$ of such covers (for groups $G_1', G_2'$), the space $Y_1\times_{X^\bullet} Y_2$ (for the group $G_1'\times_G G_2'$) is also such a cover, so we can, and will, choose $Y$ large enough so that any two color valuations whose dominant translates become incomparable in some cover are incomparable in $Y$.

\begin{definition}  \label{def:SX2}
 Let $\mathcal D^G(X)$ be the set of 
valuations associated with $G$-stable prime divisors in $X$, considered, by restriction to $F(X)^{(B)}$, as elements of $X_*(A_X)$.
Let $Y$ and $S_Y$ be as above. 
We define 
 \begin{equation}\label{eq:nonaffineclosure}S_X = S_Y \oplus \bigoplus_{\check\lambda \in \mathcal D^G(X)} T^*V_{\check\lambda},
 \end{equation}
 where $V_{\check\lambda}$ denotes the irreducible representation of $\check G_X$ of \emph{lowest} weight $\check\lambda$. 
 \end{definition} 
 
 \begin{example} 
A fairly trivial example of this situation is provided by taking $X=\mathbb{A}^r$
to be affine $r$-space considered
as a variety under $G =  \mathbb{G}_m^r$. Then  $X=Y$;
there are no colors of even sphere type, but there are 
  $r$ $G$-stable divisors, the coordinate planes in $X=Y$, 
corresponding to the standard basis 
of co-characters for the dual $\check{G} = \mathbb{G}_m^r$. Correspondingly,
$$S_X = T^* (\mathbb{A}^r).$$
\end{example}
 
 \begin{example}\label{ex:GJ}
  Consider the Godement--Jacquet case, $X = \text{Mat}_n$ under the (right) action of $G= \GL_n\times\GL_n$, that is to say, $A \cdot (g_1, g_2) = g_1^{-1} A g_2$.  We have $Y = \GL_n$ and  
  \[\check G_X= \check G_Y = \GL_n \overset{C\times I}\hookrightarrow \check G = \GL_n\times\GL_n,\] where $C=$ the Chevalley involution. The colors are precisely the simple positive coroots of $Y$;
  none of them have even sphere type. 
  The set $\mathcal D^G(X)$ that appears above
 arises from the valuation induced   by the divisor $X\smallsetminus Y$ of singular matrices.

  This valuation, identified with an element of $X_*(A_X)$, is the lowest weight of the standard representation (the standard coweight $\check\epsilon_n$ into the last diagonal entry of the upper triangular Borel).   
 To see this, we may compute as follows:  taking the reference Borel of $G$ as the lower triangular Borel in the first $\GL_n$
 and the upper triangular entry in the second, the identity matrix lies
 in the open Borel orbit $X^{\circ}$.
 Write $m_j$ for the function on $\text{Mat}_n$ given by the determinant of the upper $n \times n$ block;
 then  $m_{j}/m_{j-1}$ transforms under the Borel character $(-e_j, e_j) \in X^*(A_X) \subset X^*(A)$.
In particular, writign $\chi = (-\sum x_j e_j, \sum x_j e_j)$, we have
in the notation of \S \ref{sec:invariants} that $f_{\chi} = m_1^{x_1-x_2} \dots m_n^{x_n}$.
The valuation along the divisor $X \smallsetminus Y$ precisely
extracts the coefficient $x_n$, i.e.,  corresponds to the cocharacter $\check{e}_n$. 

From this we deduce
 $$S_Y=\mbox{trivial}, S_X=\mbox{cotangent bundle of standard representation,}$$
 with grading $1$.  
\end{example}

 \subsubsection{Whittaker induction} \label{SX3sec} 
 
  Thus far we have considered the case $M=T^*(X, \Psi)$ where $\Psi$ was trivial.

 Now we examine the case of twisted polarizations, that is to say, when $\Psi$ is nontrivial;
 in fact, Definition \ref{def:SX2} applies without change, but we need to clarify the nature of the elements that comprise it. 
 
 Write $(X,\Psi)$ as a parabolic
 induction of $(X_L,\Psi_L)$, as in \eqref{Whittakerinduction}. Recall that we have already defined the analog of ``affine closure of the open $G$-orbit'' in the twisted case by \eqref{Xcantwisted}. We now apply Lemma \ref{freecolors} to $X_L$, obtaining a homogeneous $L$-space $Y_L$, and we let $\Psi_L^Y$ denote the pullback of the $\Ga$-torsor to it. We define $Y$ by the analogous induction from $Y_L$; it comes equipped with an induced $\Ga$-torsor $\Psi^Y$. To the pair $(Y,\Psi^Y)$, we have associated by Definition \ref{def:SX} a $\check G_Y$-representation $S_Y$, which we restrict again to the subgroup $\check G_X$. Note that $\check G_Y$, $\check G_X$, here, denote the dual groups associated to the $\Ga$-bundles, not just to the spaces $Y, X$.
 
 The space $S_Y$ is the first ingredient of the definition of $S_X$, and the rest will come from $G$-stable divisors.  
 We observe that there is a bijection between $G$-stable divisors on $X$ and $L$-stable divisors on $X_L$, that is to say
we have
 $\mathcal D^G(X) = \mathcal D^L(X_L)$,
 where the equality is not merely of sets of divisors but as subsets of $X_*(A_X)$. 
 
 \begin{lemma}
In the setting above, the elements of $\mathcal D^G(X)$ are antidominant as weights of $\check G_X$.  
 \end{lemma}

\begin{proof}
A priori, these $A_X$-coweights are nonpositive on the spherical roots of $X_L$,  so, by the definition of $\check G_X$ in \S~\ref{dual-Whittaker}, it suffices to check that they are also nonpositive on the simple roots that lie in the opposite of the Lie algebra $\mathfrak u$ of $U$, or, equivalently, nonnegative on the simple roots in $\mathfrak u$.  We will, in fact, argue that they vanish on those simple roots.

Indeed, $X_L$ is a smooth affine spherical variety, hence a vector bundle over some affine homogeneous space $H\backslash L$. Let $H_0\backslash L$ denote the open $L$-orbit on $X_L$, with $H_0\subset H$. We have $A_X=A_{X_L}= (H_0\cap B_L)\backslash B_L/N_L$, where $B_L\supset N_L$ are a Borel subgroup of $L$ and its unipotent radical, assumed in general position with respect to $H_0$ (i.e., $H_0B_L$ is open in $L$). The valuation $\check v_D$ associated to an $L$-stable divisor $D$ has the property that $\check v_D(\Gm)\subset A_{X_L}$ is the image of $H_1\cap B_L$, where $H_1\subset H$ is the stabilizer of a point of $D$ in general position (in the fiber over $H1\in H\backslash L$). In particular, $\check v_D$ has image in the kernel of the canonical map $A_{X_L} \to A_{H\backslash L}$. 

We claim that all roots $\alpha\in \Delta\smallsetminus \Delta_L$ vanish on the kernel of this map. Indeed, $H$ normalizes the additive character $U\to \Ga$, and these roots are spherical roots for the Whittaker-induction of this character from the variety $H\backslash L$ to $G$, hence elements of $X^*(A_{H\backslash L})$ (as in the proof of Proposition \ref{prop:psiroots}).

\end{proof}

 We may now define the space $S_X$ by the same formula \eqref{eq:nonaffineclosure}.

 \subsubsection{} 
This concludes our definition of $S_X$ in the case at hand. However, 
we would like to clarify the relationship of this definition to the 
  results of \cite{BNS} (for certain reductive monoids), and \cite{SaWang} (when $\check G_X = \check G$). This will be accomplished by the discussion below and   Proposition \ref{Gvaluations}.

Namely, we consider the  $A_X$-toric variety $X\sslash N$; it corresponds to the saturated submonoid $\mathfrak c_X \subset X_*(A_X)$ of all cocharacters $\check\mu$ such that $\lim_{t\to 0} \mu(t)$ exists in $X\sslash N$. The valuations of colors generate a submonoid $\mathfrak c_X^D\subset \mathfrak c_X$, and the intersection of $\mathfrak c_X$ with the antidominant cone of $\check G_X$ is another submonoid $\mathfrak c_X^-$. We define a set of antidominant weights $\mathcal D_{\text{sat}}^G(X)$ (following the notation of \cite{SaWang}, where ``sat'' stands for ``saturation''), as consisting of those nonzero elements of $\mathfrak c_X^-$ which 
\begin{itemize}
 \item are primitive in $\mathfrak c_X^-$, i.e., cannot be written as sums of two nonzero elements of if;
 \item cannot be written as $\check\theta + \check v$, with $\check\theta\in \mathfrak c_X^-$ (possibly zero) and $\check v$ a nonzero element of $\mathfrak c_X^D$. 
\end{itemize}

Note that the monoid that is generated by $\mathfrak c_X^D$ and $\mathcal D_{\text{sat}}^G(X)$ contains $\mathfrak c_X^-$, although it may not be equal to $\mathfrak c_X$, for a general spherical variety.\footnote{The latter was incorrectly asserted in the introduction of \cite{SaWang}; the example of $\Gm\backslash(\Gm\times\PGL_2)$ shows that it doesn't have to be true. The definition of the set $\mathcal D_{\text{sat}}^G(X)$ suffered a lot in that paper, with an incorrect description given in \cite[Corollary 5.1.4]{SaWang}; however, the correct definition appears before Corollary 5.6.5, and is the one that is being used consistently in all proofs.}  

\begin{proposition}\label{Gvaluations}
 Assume that $X$ is a smooth, affine spherical variety. Then $\mathcal D_{\text{sat}}^G(X)=\mathcal D^G(X)$, that is to say, $\mathcal D_{\text{sat}}^G(X)$ is the set of 
valuations associated with $G$-stable prime divisors in $X$, considered, by restriction to $F(X)^{(B)}$, as elements of $X_*(A_X)$. 
\end{proposition}

It is the set $\mathcal{D}_{\text{sat}}^G(X)$, instead of  $\mathcal D^G(X)$, that appears
in \cite{SaWang}; but this proposition shows that the two coincide, for $X$ smooth affine. 
Note that, unlike the case of colors, for valuations induced by $G$-stable divisors (and for $G$-invariant valuations, more generally), the map 
\[ \mbox{valuation } \mapsto \mbox{ its restriction to } F(X)^{(B)}\]
is injective; see \cite[Proposition 7.4]{LV}, \cite[Corollary 1.8]{KnLV}. 

\begin{proof}
Recall that a spherical variety is called \emph{simple} if it has a unique closed $G$-orbit; an affine spherical variety is always simple. By \cite[\S~5.1]{Brion-varietesspheriques}, a simple spherical variety is locally factorial iff the set of valuations of all $B$-stable divisors which contain the closed $G$-orbit (restricted to $F(X)^{(B)}$) forms part of a basis of $X^*(A_X)$. Let us, for the purposes of this proof, denote the set of those valuations by $\Delta_X$. Let us also denote by $\Delta_X'$ the valuations corresponding to colors that don't contain the closed $G$-orbit (in their closure). The rational cone in $X_*(A_X)\otimes \mathbb Q$ spanned by $\mathfrak c_X$ is the same cone as that spanned by $\Delta_X\cup \Delta_X'$. 

If, in particular, $X$ is a vector space with a linear $G$-action, then every $B$-stable divisor contains the origin, therefore $\mathcal D_{\text{sat}}^G(X)$ consists of the elements of $\Delta_X$ that are do not come from colors, i.e., come from $G$-stable divisors. This proves the proposition in the case of vector spaces. 

For the general smooth affine case, where $X$ is a vector bundle over an affine homogeneous space $H\backslash G$, we let $x_0\in H\backslash G$ be the point $H1$, and denote by $L\subset H$ the stabilizer of a generic point on $\mathfrak h^\perp\subset \mathfrak g^*$. (This is defined up to conjugacy in $H$, and its conjugacy class in $G$ is a generic stabilizer for the $G$-action on $T^*X$.) By \cite[Lemma 5.2]{KnVS}, if $B$ is a Borel subgroup such that $x_0B$ is open in $H\backslash G$, then $B\cap H$ is a Borel subgroup of a representative for $L$. Now, the variety $X$ is spherical iff $B_L:= B\cap H$ acts with an open orbit on the fiber $V$ of $X\to H\backslash G$ over $x_0$; that is, iff $V$ is $L$-spherical. Moreover, the closed $G$-orbit in $X$ being $H\backslash G$, we have a clear bijection, by inclusion:
\[ \{ B_L\mbox{-stable divisors in }V\} \leftrightarrow \{ B\mbox{-stable divisors in $X$ that contain the closed orbit}\}.\]

Moreover, by \cite[Theorem 6.7]{KnLV}, $X$ being affine implies that there is a $\chi\in X^*(A_X)$ such that $\chi$ is strictly negative on $\Delta_X'$ (the valuations of colors that \emph{don't} contain the closed $G$-orbit), and zero on $\mathfrak c_X^-$. Therefore, in the definition of $\mathcal D_{\text{sat}}^G(X)$ we could have replaced $\mathfrak c_X^D$ by the submonoid generated by the colors that belong to $\Delta_X$ (i.e., contain the closed $G$-orbit). This reduces us to the case of the spherical $L$-module $V$, where the claim has already been proven. 
\end{proof}

\subsection{The space $V_X$; how we arrived at the formula for $S_X$} \label{NTmov}

In arithmetic applications, what plays a more important role than $S_X$ is the space $V$ of \eqref{fund_id}, which here we will denote by 
\begin{equation}  \label{vx} V_X = S_X \oplus [\mathfrak{g}_X^{\perp} \cap \check{\mathfrak{g}}_e].\end{equation} It is a  $\check{G}_X$-representation,  
which is naturally graded by the action of a group $\GGm'\simeq\Gm$ with $S_X$ in degree $1$ and the remainder
graded by $2+$ the weight under the action of $h$ of the $\mathfrak{sl}_2$-triple. 
It is self-dual, but not necessarily symplectically self-dual.

There is, by \eqref{fund_id}, an identification 
of $\check{G}$-spaces   
\[\check M = V_X \times^{\check G_X} \check G.\]
With reference to this identification, the $\GGm'$ corresponds to 
 the image of $\Gm\xrightarrow{(-2\rho_{L(X)},\Id)} \check G\times \GGm$, 
 see \eqref{Gmprime}.

 Moreover, as was remarked after \eqref{cfp5}, while the isomorphism above depends on the choice of a splitting of the canonical map $\check{\mathfrak g}^*\to \check{\mathfrak g}_X^*$, the map $\check M\to \check G_X\backslash \check G$ is intrinsic, i.e., determined by the data $(\check G_X, \sl_2, S_X)$ of the Whittaker induction.

  The space $V_X$, rather than its subspace $S_X$, is what 
 appears naturally in the theory of automorphic forms:
 it is the $\check{G}_X$-representation that appears in the  local Plancherel formula ( \S~\ref{section-unramified-local})
 as well as the theory of global periods (e.g., as in  \eqref{XperiodSV}).
   These interpretations  lead to a candidate combinatorial description of $V_X$ by
extrapolating from results of \cite{SaSph} (when $X= H\backslash G$ is affine homogeneous), and their extension to non-homogeneous affine varieties by \cite{BNS} (for certain reductive monoids) and \cite{SaWang} (when $\check G_X = \check G$). However, such a description is quite involved, and not very enlightening. 
Nonetheless, by studying these candidate descriptions, we arrived at the
conjectural description of $S_X$ that has been presented in this chapter. 
  
   \begin{example}
  Consider the Shalika model of $\GL_{2n}$, which is by definition the Whittaker induction of the variety $X_L = \GL_n^{\text{diag}} \backslash\GL_n^2$ along
  the homomorphism $$ \GL_n \times \SL_2 \rightarrow \GL_{2n}$$
  arising from the tensor product of the standard representations. 
   We have $\check G_{X_L} = \GL_n$ (embedded Chevalley-diagonally as in Example \ref{ex:GJ}), and $\check G_X = \Sp_{2n}$. Moreover, 
   $S_{X_L}$ and $S_X$ are both trivial; and $V_{X_L} = \mathfrak{gl}_n = \mathfrak g_a \oplus \mathfrak{pgl}_n$, and $V_X = \wedge^2$, the exterior square of the standard representation (which includes the trivial representation $\mathfrak g_a$ as a direct summand). 
 \end{example}

\subsection{Parity} \label{parity} \index{parity element}

Now we discuss the parity of $\check{M}$; see \S \ref{analyticarithmetic} for a general discussion that this plays in the paper. 
 We continue to assume the existence of a nowhere vanishing eigen-volume form on $X$ with eigencharacter $\eta$, as in \eqref{eigencharacter}. We will be using additive notation for characters, and exponential notation when we want to think of their images in $\Gm$.

We define the following central element of $\check{G}$,
\begin{equation}\label{zXdef} z_X := (-1)^{\eta + 2 \rho}.
\end{equation}

\begin{proposition} \label{zXparity}
 The action on $\Mv$ of the central involution $z_X$ may be identified with that of the involution $(-1)\in \GGm$.  
\end{proposition}
For example, if we take the case $X=\mathbb A^1$ of Tate's thesis ($G=\Gm$), we have $z_X = -1 \in \check G = \Gm$ 
and indeed the action of $-1 \in \check{G}$ on the dual space $T^* \mathbb{A}^1$ coincides
with the scaling action of $-1 \in \GGm$.

 There is a reformulation of the Proposition using the following
 \begin{definition} \label{arithmetic Gm action}
 Suppose that $(G, M)$ and $(\check{G}, \check{M})$ is as above. 
 The {\em arithmetic $\GGm$-action} on $\check{M}$ is the product of the neutral action,
 and the action through $\GGm \stackrel{\eta}{\rightarrow} \check{G}$
 dual to the eigencharacter $\eta$ of the eigen-volume form on $X$. 
 \end{definition}
 Then the Proposition says that the {\em arithmetic} action of $G \times \GGm$
 factors through the extended group ${}^CG_z$ of \eqref{CGz def}.
 \index{${}^CG_z$}

\begin{proof}  (of the Proposition)
First, we will confirm that $z_X$ and $(-1)\in\GGm$ act the same way on $\check G_X\backslash\check G$. 

Let us consider the $\Gm$ action from \eqref{varpixdef}, that is to say, 
$\Gm$ acting diagonally via the embedding $$(\mathrm{id}, -2\rho_{L(X)}): \Gm \rightarrow \GGm \times  \check{G},$$
where we recall that the character $2\rho_{L(X)}$ is the datum $\varpi$ of Whittaker induction in the definition of $\check M$
(see \eqref{2rhoreminder}). 
This action preserves the coset of $1$ in $\check G_X\backslash \check G$; in particular, this is the case for the element $-1\in \Gm$, which is embedded as $(-1, (-1)^{-2\rho_{L(X)}})$ in $\GGm\times\check G$. 

Now, from \eqref{rhorho1}, we have that $\eta + 2 \rho - 2\rho_{L(X)}$ is a cocharacter into $\check G_X$; thus, the action of $(-1, z_X)\in \GGm\times\check G$ also preserves the coset $\check G_X \cdot 1$; since this element is central, it acts trivially on $\check G_X\backslash \check G$. 

Next, we will show that $(-1, z_X)$ acts trivially on the space $V_X$ of \eqref{vx}. The action on the summand $\check{\mathfrak g}_X^\perp\cap \check{\mathfrak g_e}$ is clearly trivial, since it is induced (see \S~\ref{ggmm}) by a combination of the square action of $\GGm$ (which is trivial on $-1$) and the coadjoint action of $\check G$ on $\check{\mathfrak g}^*$ (which is trivial for the central element $z_X$). Finally, we are left with showing that the action of $(-1)^{\eta + 2 \rho - 2\rho_{L(X)}}\in \check G_X$ is odd on $S_X$. 

Note that the element $(-1)^{\eta + 2 \rho - 2\rho_{L(X)}}$ is central in $\check G_X$, since $z_X$ is central in $\check G$ and $2\rho_{L(X)}$ commutes with $\check G_X$. Therefore, it suffices to show that it is odd on a set of representatives for the $W_X$-orbits of weights on $S_X$. By the construction of $S_X$ in \S~\ref{SXdef}, such representatives consist of 
\begin{itemize}
 \item valuations of colors of even sphere type; 
 \item the weights in $\mathcal D^G(X_L)$, where $X_L$ is the Whittaker-inducing variety $S^+\times^H L$. (Here we recall that $L$ is the Levi centralizing the element $h$ of the $\sl_2$-triple defining $X$.)
\end{itemize}

We first deal with the colors of type $T$ (a subset of the colors of even sphere type). We claim that, for every color $D$ of type $T$, we have 
\begin{equation}\label{colorparity}
\left<\check v_D,  \eta + 2\rho - 2\rho_{L(X)}\right> = 1. 
\end{equation}

The argument is very similar to that we have already given in Lemma \ref{vrhopair}. 
Indeed, let $\alpha$ be a simple root such that $DP_\alpha$ contains the open Borel orbit $X^\circ$, then  the $P_\alpha$-stabilizer of a point in $D\cap X^\circ P_\alpha$ is of the form $T\cdot N_{L(X)}$, where $N_{L(X)}$ is the unipotent radical of $L(X)\cap P_\alpha$, and $T$ is a torus such that the image of $T\to A$ is the preimage of $\Gm\overset{\check v_D} \hookrightarrow A_X$ under the quotient map $A\to A_X$. \footnote{This follows from the description of the open $P(X)$-orbit in \eqref{openorbit}, and the isomorphism $X^\circ P_\alpha/\mathcal R(P_\alpha)\simeq \Gm\backslash \PGL_2$.} The existence of an eigen-volume form with eigencharacter $\eta$, now, restricted to the open $P_\alpha$-orbit, implies that  
\[\left<\check v_D,  \eta + 2\rho -\alpha - 2\rho_{L(X)}\right> = 0 \iff \left<\check v_D,  \eta + 2\rho - 2\rho_{L(X)}\right> = \left<\check v_D,  \alpha \right>,\] 
which, by \eqref{coloreq}, is $1$. 

By a similar argument, we can show that the remaining valuations of colors associated to spherical roots $\gamma$ of even sphere type have odd pairing with $\eta + 2 \rho - 2\rho_{L(X)}$. Namely, if $P$ is the parabolic whose Levi has simple roots the simple roots in the support of $\gamma$, so that the adjoint group of the Levi is $\SO_{2n+1}$ with $n \ge 2$ or $G_2$, as we noted in Remark \ref{remark-spheretype} there is a unique corresponding color with valuation $\check v_D= \frac{\check\gamma}{2}$, and the Levi of $P(X)$ contains all but one of the simple roots of the Levi of $P$ (namely, those that are orthogonal to $\gamma$).  We can now calculate $\left<\check v_D,  \eta + 2\rho - 2\rho_{L(X)}\right> $ as follows: Let $P_1$ be the parabolic generated by $P(X)$ and $P$. Then $2\rho - 2\rho_{L(X)} = 2\rho_{P_1} + 2\rho_{L \cap P(X)}$, where $L\cap P(X)$ is the intersection of $P(X)$ with the Levi of $P$. Since $\check\gamma$ maps into the derived group of that Levi, its pairing with $\eta$ and $2\rho_{P_1}$ is zero, and we are left with computing $\left<\frac{\check \gamma}{2},  2\rho_{L\cap P(X)}\right>$, which is equal to $2n-1$ in the case of $\SO_{2n+1}$ and $5$ for the case of $G_2$.

Finally, we show that the pairing of $\eta + 2 \rho - 2\rho_{L(X)}$ with the elements of $\mathcal D^G(X_L)$ is odd.
Let $X_B$ be the complement of the union of colors. It is stable under the parabolic $P(X)$. By the local structure theorem of \cite{BLV} (see also \cite[Theorem 2.3]{KnMotion}), the isomorphism \eqref{openorbit} extends to an isomorphism $X_B\simeq \overline{T_X}\times U_{P(X)}$, where $\overline{T_X}$ is a smooth toric embedding of $T_X$. In particular, $G$-stable divisors in $X$ are in bijection with $T_X$-stable divisors in $\overline{T_X}$. The non-vanishing eigen-volume form on $X$, now, restricts to a volume form $\omega_{\overline{T_X}}\otimes \omega_{U(X)}$ on $X_B$, where $\omega_{U(X)}$ is a Haar volume form on $U(X)$. The eigencharacter for the action of $T_X$ on $\omega_{\overline{T_X}}$ is $\eta+2\rho - 2\rho_{L(X)}$. For that to be the eigencharacter of a nonvanishing volume form in the neighborhood of a $T_X$-stable divisor $D$ on the toric variety $\overline{T_X}$, the corresponding valuation must satisfy, again, 
\[\left<\check v_D,  \eta + 2\rho - 2\rho_{L(X)}\right> = 1.\]

This completes the proof of the proposition. 
\end{proof}

\subsection{Regular nilpotent elements in the image of the moment map} \label{regnilp}

The goal of this subsection is to prove the following:

\begin{proposition}\label{regnilpotent}
Let $X$ be an untwisted spherical $G$-variety satisfying our standard assumptions. Assume Conjecture \ref{conjsymplectic}, and let $\mu:\check M \to \check{\mathfrak g}^*$ be the dual Hamiltonian space with its moment map. The image of $\mu$ contains a regular nilpotent element in $\check{\mathfrak g}^*$. 
\end{proposition}

This statement will be useful in our discussion of rationality. It also plays an important philosophical role, which  
will be discussed in a sequel to this paper (some related remarks here:  Proposition \ref{invariantregular}, Example \ref{JHJG}).

\begin{proof}

The definition \eqref{defM} of $\check M$ as Whittaker induction of the space $S_X$ (\S~\ref{SXdef}) with respect to the $\sl_2$-pair $(\varpi, f)$ defined by the parabolic $P(X)$ means that we have a map $\check M \to U\check G_X\backslash \check G$, where $U$ is the unipotent radical of the parabolic associated to the cocharacter $\varpi$ (as in the definition of Whittaker induction, \S~\ref{Whittaker induction}); moreover, the fiber of this map over the identity coset in $U\check G_X\backslash\check G$ is (in the notation of \eqref{tildeSdef}) 
\[ (S_X \times (\mathfrak{u}/\mathfrak{u}_+)_f) \times_{(\check{\mathfrak g}_X + \mathfrak u)^*} \mathfrak g^*.\]
In particular, taking only the zero point of the vector space $\mathfrak{u}/\mathfrak{u}_+$, the image of the moment map contains 
\begin{equation}\label{subspaceM}
(f+\mu(S_X)) \times_{(\check{\mathfrak g}_X + \mathfrak u)^*} \check{\mathfrak g}^*,
\end{equation}
where $\mu:S_X\to (\check{\mathfrak g}_X +\mathfrak u)^*$ denotes, here, the moment map for $S_X$, considered orthogonal to $\mathfrak u$.

Let us first explain how to deal with the case when no simple root of $G$ is a spherical root, i.e., there are no colors of type $T$ (see \S~\ref{SXdef}). In most cases, this means that $S_X=0$, but, in any case, the space \eqref{subspaceM} contains
\[ f+(\check{\mathfrak g}_X + \mathfrak u)^\perp.\]
Choose also an invariant identification $\check{\mathfrak g} = \check{\mathfrak g}^*$, for ease of description; then
 $f$ -- which, we recall, arises from a principal $\SL_2$ into $L(X)$ -- can be written as $\sum_{\alpha \in \Delta_{L(X)}} u_{-\check\alpha}$, where $\Delta_{L(X)}$ denotes the set of simple roots for the
 Levi $L(X)$ of $P(X)$, and the $u_{-\check\alpha}$ are basis vectors in the opposite root subspaces in $\check{\mathfrak g}$. (We refer, here, to the root decomposition with respect to the standard maximal torus $\check A$ and Borel of $\check G$.) The regular nilpotent element of $\check M$ that we will construct will belong to the space above and have the form $m = \sum_{\alpha \in \Delta_G} u_{-\check\alpha}$, a similar sum but over the entire set of simple roots of $G$. The issue is to show that we can choose the remaining basis vectors $u_{-\check\alpha}$, for $\alpha\in \Delta_G\smallsetminus\Delta_{L(X)}$, so that their sum is orthogonal to $\check{\mathfrak g}_X + \mathfrak u$. 

Regarding $\mathfrak u$, this is automatic for \emph{any} choice of basis vectors in these negative root spaces. Indeed, since $\varpi = 2\rho_{L(X)}$ is a sum of positive roots of $L(X)$, we have $\langle \varpi, \check\alpha\rangle \le 0$ for every $\alpha\in \Delta_G\smallsetminus\Delta_{L(X)}$, and that means that the simple root spaces $\check{\mathfrak g}_{-\check\alpha}$ belong to $\mathfrak u$; under the identification $\check{\mathfrak g} \simeq \check{\mathfrak g}^*$, they are orthogonal to $\mathfrak u$. 

Thus, to finish with the case when there are no spherical roots which are simple roots of $G$, there remains to show that the sum 
\[\sum_{\alpha \in \Delta\smallsetminus\Delta_{L(X)}} u_{-\check\alpha}\]
can be taken to be orthogonal to $\check{\mathfrak g}_X$. 
For this we will recall some of the details of the Knop--Schalke construction.
  
Recall from \cite[\S 6]{KnSch} that we have inclusions $\check G_X \subset \hat G_X \subset \check G$, for some intermediate reductive subgroup $\hat G_X$ that contains $\check A$, and that $\check G_X$ is obtained by a process of ``folding'' of the roots of $\hat G_X$. This means that the simple coroots of $\check G_X$ are either coroots of $\hat G_X$ (and, therefore, roots of $G$), or sums $\gamma=\alpha+\beta$ of two simple coroots of $\hat G_X$. Moreover, by inspection of \cite[Table (6.1)]{KnSch}, in the latter case either both or none of $\alpha$ and $\beta$ are simple roots of $G$. Finally, such roots $\alpha$, $\beta$ are associated, in this sense, to a unique simple coroot $\gamma$ of $\check G_X$; this follows from \cite[Lemma 6.4]{KnSch} which shows that, given a simple root $\alpha$ of $G$ which is associated to a simple spherical root $\gamma$, that spherical root $\gamma$ is characterized, among simple coroots of $\check G_X$, by the property that $\langle \gamma, \check\alpha\rangle \ge 0$.

Now, if no simple root of $G$ is a spherical root, then each $\alpha\in \Delta_G\smallsetminus\Delta_{L(X)}$ either is not a root of $\hat G_X$, in which case the entire root space $\check{\mathfrak g}_{-\check\alpha}$ will be orthogonal to $\check{\mathfrak g}_X$; or, there is another simple root $\beta$ of $G$, such that $\alpha$, $\beta$ are associated to a spherical root $\gamma$, in which case we should choose $u_{-\check\alpha}$, $u_{-\check\beta}$ so that it is orthogonal to ``the'' simple root space $\check{\mathfrak g}_{X,\check\gamma}$ of $\check{\mathfrak g}_X$. (Of course, here, we have chosen an embedding $\check{\mathfrak g}_X \hookrightarrow \check{\mathfrak g}$, together with a commuting $\sl_2\to \check{\mathfrak g}$; in the Knop--Schalke construction, these data are determined up to $\check A$-conjugacy.)

Finally, we deal with the case of simple spherical roots of type $T$. It is enough to assume that $X$ is the affine closure of its open $G$-orbit; indeed, in the general case, the space $S_X$, and hence also the image of the moment map, only gets larger, by \eqref{eq:nonaffineclosure}. 
Using the identification $\check{\mathfrak g} \simeq \check{\mathfrak g}^*$, which restricts to $\check{\mathfrak g}_X \simeq \check{\mathfrak g}_X^*$, let us write the space \eqref{subspaceM} as 
\[ f+\mu(S_X)+ (\check{\mathfrak g}_X + \mathfrak u)^\perp.\]
Now, we can still choose basis vectors $u_{\check\alpha}$ as above for the simple root spaces $\check{\mathfrak g}_{-\check\alpha}$, for those $\alpha \in \Delta_G\smallsetminus\Delta_{L(X)}$ that are not spherical roots.  Let $\Delta_X^T\subset \Delta_G$ be the set of simple roots of $G$ that are also spherical roots. We will show that there is an element $s\in S_X$ with 
$$\mu (s) \perp ( \check{\mathfrak b}_X^-) $$
 where $\check{\mathfrak b}_X^- \subset \check{\mathfrak g}_X$ denotes the opposite of the standard Borel subalgebra   and moreover, for all simple roots $\check\gamma$ of $\check G_X$,
 \[\langle \check{\mathfrak g}_{X,\check\gamma}, \mu(s)\rangle \neq 0 \iff \gamma \in \Delta_X^T.\]  The resulting element \[f+\mu(s)+\sum_{\alpha \in \Delta_G\smallsetminus(\Delta_{L(X)}\cup \Delta_X^T)} u_{-\check\alpha}\] will then be the desired regular nilpotent element in the image of the moment map. 

Let $S_X= S_X^+ \oplus S_X^-$ be the decomposition into a sum of Lagrangians, corresponding to the decomposition $\mathfrak B_X = \mathfrak B_X^+ \cup \mathfrak B_X^-$ of the weight multisets of Conjecture \ref{conjsymplectic}. Recall that $\mathfrak B_X^+$ contains $\mathfrak B_X^C:=$ the set of valuations associated with colors of type $T$.  
Choose basis vectors $s_b$, $b\in \mathfrak B_X^C$, for the corresponding weight spaces, and set $s = \sum_{b\in \mathfrak B_X^C} s_b \in S_X^+$. 

Recall that the moment map $\mu: S_X \to \check{\mathfrak g}_X^*$ for a symplectic representation is defined by $\langle Z, \mu(v) \rangle = \frac{1}{2} \omega(v\cdot Z , v)$ for $Z\in \check{\mathfrak g}_X$ (see \eqref{smm}). 
By construction, the space $S_X^+$ is stable under the action of the negative Borel subalgebra $\check{\mathfrak b}_X^-$ of $\check{\mathfrak g}_X$. (Indeed, recall that $\check G_X$ is acting on the right, so the action of the $(-\check\gamma)$-root space adds $\check\gamma$ to the weight.) Moreover, Conjecture \ref{conjsymplectic} implies that, for a simple root $\check\gamma$ of $\check G_X$ and a basis element $e_{\check\gamma}\in \check{\mathfrak g}_{X,\check\gamma}$, we have 
\[ s\cdot e_{\check\gamma} \in \begin{cases} S_X^+,& \mbox{ if } \check\gamma\notin \Delta_X^T; \\ 
                                c_1 s_{-\check v_1} + c_2 {s_{-\check v_2}} + S_X^+, & \mbox{ if }\check\gamma\in \Delta_X^T,
                               \end{cases}\]
where, in the latter case, $\check v_1$, $\check v_2$ are the valuations of the two colors associated to $\check\gamma$, and $c_1, c_2$ are two nonzero constants
 and $s_{-\check v_1}, s_{-\check v_2}$ have weights $-\check{v}_1, -\check{v}_2$ respectively.

The  symplectic pairing $\omega(s, s\cdot e_{\check\gamma})$ is then nonzero,
as follows from the nondegeneracy of the symplectic form on $S_X$ and the fact that the weight spaces for $\pm \check v_i$ are each one-dimensional, again
by Conjecture \ref{conjsymplectic}. 
This shows that
\[ \langle \check{\mathfrak g}_{X,\check\gamma}, \mu(s)\rangle  \begin{cases} =0,& \mbox{ if } \check\gamma\notin \Delta_X^T; \\ 
                                \ne 0, & \mbox{ if }\check\gamma\in \Delta_X^T,
                               \end{cases}\]
as desired.

\end{proof}

\begin{remark} \label{twisted} 
In the twisted case, Proposition \ref{regnilpotent} is generally (likely always) false. For example, if the dual $\check M$ is polarized, $\check M = T^*\check Y$, and our duality is involutive, as predicted, then the Arthur-$\sl_2$ defining the twisting for $M$ is dual to the parabolic $P(\check X)$, which is therefore not minimal. By \cite[Satz 5.4]{KnWeyl}, the moment image of $T^*\check X$ only contains vectors which are perpendicular to $[\mathfrak p, \mathfrak p]$ for some parabolic $P$ in the conjugacy class of $P(\check X)$, and the largest nilpotent orbit in this set is the \emph{Richardson orbit} associated to $P(\check X)$, which is not regular. It would be interesting to examine if that Richardson orbit is in the image of the moment map for $\check M$ -- but, already, extracting the parabolic $P(\check X)$ out of our description of $\check M$ does not seem straightforward, and we will not attempt to make any progress on these questions here.
\end{remark}

\subsection{ Rational and Frobenius structures on $\check M$} \label{CheckMRat}

For $M=T^*X$ a polarized {\em untwisted} hyperspherical space over $\mathbb F$, we have constructed a dual Hamiltonian space $\check M$ over $k$; both $\mathbb F$ and $k$ have been assumed to be algebraically closed of characteristic zero. In this section, we will discuss the following two issues: 

\begin{description}
 \item[Rationality] Is there a ``distinguished'' form of the dual $\check M$, if $k$ is not algebraically closed?
 \item[Galois action] If $M$ (along with its polarization) is defined  over a subfield $\FF_0\subset \FF$, so that $\FF$ is the algebraic closure of $\FF_0$, is there a natural action of $\Gamma =  \mathrm{Gal}(\FF/\FF_0)$ on $\check{M}$?
\end{description}

 In this paragraph we will describe a ``simple'' action of the Galois group on $\check M$ as well as a distinguished class of $\kk$-rational structures,
where $\kk$ is an arbitrary field {\em in which $2$ is invertible} (see footnote \ref{bad2footnote});
the $\kk$-rational structure, which we will call the ``distinguished split form,'' will be uniquely specified in the absence of even sphere roots not of type $T$, see Remark \ref{dependence-egamma}. 
The Galois action is compatible with the standard, ``analytic'' Galois action on $\check G$
(\S \ref{sss:Satake-shearing-arithmetic}) according to which the action of $\Gamma$ preserves a pinning, corresponding to a choice of basis vectors $e_{\check\alpha} \in \check{\mathfrak g}_{\check\alpha}$, for all simple roots $\check\alpha$ of $\check G$, or dual basis vectors $f_{\check\alpha} \in (\check{\mathfrak g}^*)_{-\check\alpha}$; set $f = \sum_{\check \alpha} f_{\check\alpha}$. 
 Also, this simple action will be used
later on (\S \ref{Mshear1})  to produce Frobenius actions on the sheared coordinate ring of $\check{M}$. 
 
The principle guiding many of our constructions is that we should seek {\em pinned hyperspherical varieties}. Recall from Proposition \ref{regnilpotent} that, when $\check M$ is the dual of an untwisted spherical variety, there is a regular nilpotent element $f$ in the image of the moment map $\check M\to \check{\mathfrak g}^*$. Such a nilpotent element was actually described, in the proof of that proposition, in terms of the structure
\begin{equation} \label{cMstructurerecall}
 \check M = (S_X \times (\mathfrak{u}/\mathfrak{u}_+)_f) \times^{U\check G_X}_{(\check{\mathfrak g}_X + \mathfrak u)^*} T^*G.\end{equation}
We will now take this regular nilpotent element $f$ to be the canonical one coming from the pinning on $\check G$, and we postulate:  \begin{quote}
Guiding principle: \emph{When $\check M$ is the dual of an untwisted spherical variety $X$ the ``simple'' action of $\Gamma$ on $\check M$ should preserve an element $m\in \check M$ with image $f$ under the moment map. When $k$ is not algebraically closed, there should be a distinguished rational form of $\check M$ such that $m\in \check M(k)$.}
\end{quote}
We will refer to the pair $(\check{M},m)$ as a {\em pinned hyperspherical $G$-space}, relative to the pinned quasisplit reductive group $G$. 
Certainly, not all
hyperspherical spaces admit a pinned form; indeed, the moment map need not even meet the regular locus,
and  it is not clear to us that a pinned form is always unique up to unique isomorphism.  
Nonethelss, the above
definition is useful particularly for duals to (untwisted) polarized hyperspherical space.

\begin{itemize}
\item In \S \ref{simple-M} and Example \ref{unstableBC} we discuss some motivation coming from the theory of automorphic forms, in particular, the consideration of ``stable'' versus ``unstable'' base change.
\item In \S \ref{constructionGaloissplit} we construct the distinguished split form and Galois action in the case when $X$ is untwisted, i.e., no $\Ga$-torsor.  
The definition in the case $S_X=0$ is completed by
Definition \ref{definition-LGX}; the symplectic form in $S_X$
is constructed in Lemma \ref{lemmasymplform}, essentially
by specifying it on some distinguished basis elements, and
then we use the same basis elements to pin down the Galois action thereafter.
\item In 
Proposition \ref{invariantregular},   we verify that the construction just given does in fact have the property
quoted in the ``guiding principle'' above. 

\end{itemize}

\subsubsection{Motivation from the theory of automorphic forms} \label{simple-M}
The question we are discussing here is related to one
that has been studied in the theory of automorphic forms. 
Namely, 
it is related to  the question of extending the dual group $\check G_X$ of a spherical variety to an $L$-group of the form $\check G_X\rtimes \Gamma$, together with an embedding to the $L$-group of $G$ over $\Gamma$. This has been addressed in the literature \cite[\S~10]{KnSch}, but does not have a definite answer yet;
presumably, one could recover the correct Galois action by extending the work of Gaitsgory--Nadler \cite{GaitsgoryNadler} to a curve over a non-algebraically closed field.

To begin with, whenever $X$ is defined over a field $\mathbb F$ that is not algebraically closed, the set of its simple spherical roots (i.e., simple coroots of $\check G_X$) admits an action of the Galois group $\Gamma$. The easiest way to see this is to notice that the abstract Cartan $A_X$, introduced in \S~\ref{sec:invariants}, is defined canonically up to unique isomorphism over $\FF$ (as is the abstract Cartan $A$ of $G$ -- even if $G$ is not quasisplit!), giving rise to an action of the Galois group on $\mathfrak a_{X,\mathbb R}  = \Hom_{\overline{\FF}}(\Gm, A_X)\otimes \mathbb R$, which preserves the cone $\mathcal V_X$ of invariant valuations. This suggests a naive action of $\Gamma$ on $\check G_X$ -- namely that there is a pinning on the triple $\check A_X \subset \check G_X\cap \check B\subset \check G_X$ which is preserved by the action of $\Gamma$.
However, this proposal does not interact
well with the role of $\check G_X$ in the Langlands program,  as the following example shows. 

\begin{example}\label{unstableBC}
 Let $\mathbb E/\mathbb F$ be a quadratic extension, $G=$ the Weil restriction of scalars to $\FF$ of $\GL_{n, \mathbb E}$, $H=\GL_n$ over $\FF$, and $X = H\backslash G$. Then, the ``correct'' $L$-group of $X$ is known to be the subgroup
 \[ {^LG_X} \subset {^LG}\]
 described as follows: First, recall that ${^LG} \simeq (\GL_n\times\GL_n)\rtimes\Gamma$, with Galois conjugation for $\mathbb E/\FF$ switching the two copies of $\GL_n$ (and preserving their pinning, which we will take to be the standard pinning defined by upper triangular matrices with $1$'s on the $(i,i+1)$-entries). Now, ${^LG_X}\simeq \GL_n\rtimes \Gamma$, where $\check G_X\simeq \GL_n \hookrightarrow \GL_n^2$ is embedded as $g\mapsto (g, g^d)$, where the exponent $d$ denotes the \emph{unpinned} duality involution introduced in \S~\ref{dualityinvolution}. Explicitly, 
 \[ g^d = w \cdot {^tg^{-1}} w,\]
 where $w$ is the matrix with $1$'s on the antidiagonal and $0$'s elsewhere.
 
 The reason that this is known to be the ``correct'' $L$-group has to do with poles of the Asai $L$-function, which are detected by the $H$-period in $G$, see \cite{Flicker}. When $n$ is odd, this $L$-group is $\check G$-conjugate to the one that one would obtain from the pinned embedding $g\mapsto (g, g^c)$ of $\GL_n$ into $\GL_n^2$, therefore the difference does not matter for the purposes of Langlands duality; thinking of the latter as the $L$-group of the unitary group $U_n$, this embedding is the standard, so-called ``stable'' embedding ${^LU_n}\hookrightarrow {^LG}$. However, when $n$ is even, although ${^LG_X}$ is still abstractly isomorphic to ${^LU_n}$, its embedding into ${^LG}$ is not $\check G$-conjugate to the standard one; it is customary to call it the ``unstable'' embedding of ${^LU_n}$ (or, more often, to talk of ``unstable base change'' of automorphic representations from $U_n$ to $G$. 
 
Our interpretation of the example above, and other examples that we will discuss below, is that the appropriate Galois-fixed pinning should not be on the subgroup $\check G_X$ itself, but on the Hamiltonian space $\check M$, in the sense we have discussed after \eqref{cMstructurerecall}. 
\end{example}

 \subsubsection{Construction of the Galois action and rational structure, untwisted case} \label{constructionGaloissplit}
 The construction of this action goes through the following steps:

\begin{enumerate}
\item 
Recall that in \S~\ref{sec:invariants} we considered the dual group $\check G_X$ as a subgroup of $\check G$, unique up to $\check A(\bar k)$-conjugacy. We will start here by describing a precise $\check A$-conjugate $\check G_X\subset \check G$, with a commuting $\sl_2 \to \check{\mathfrak g}$, largely following \cite[\S~10]{KnSch}. Once we have done this, the rational structure and Galois action on $\check G$ will then give rise to a compatible rational structure and $\Gamma$-action on $\check G_X$. We will call this the ``simple'' or ``analytic'' action of $\Gamma$ on $\check G_X$.

\item 
The map $\sl_2\to \check{\mathfrak g}$ will map the standard nilpotent $e\in \sl_2$ to $\sum_{\alpha\in \Delta_{L(X)}} e_{\check\alpha}$. We will call this the ``pinned'' $\sl_2$. 

\item 
The subgroup $\check G_X$ is determined by its simple root spaces $\check{\mathfrak g}_{X,\check\gamma}$. Each spherical root $\gamma$ is either a root of $G$, or the sum of two strongly orthogonal roots, $\gamma = \alpha+ \beta$. (This includes the twisted case, see \S~\ref{dual-Whittaker}.) In the first case, the simple root space $\check{\mathfrak g}_{X,\check\gamma}$ will map isomorphically onto the corresponding simple root space for $\check{\mathfrak g}$. In the second case, there is a choice to be made. However, \cite[Lemma 10.4]{KnSch} shows that there is a \emph{unique} $1$-dimensional subspace of $\check{\mathfrak g}_{\check\alpha}+\check{\mathfrak g}_{\check\beta}$ that commutes with the pinned $\sl_2$, \emph{unless} $\alpha$ and $\beta$ are simple (as roots of $G$). If $\check\alpha,\check\beta$ are simple, our desideratum that the moment image of the fiber of $\check M\to U\check G_X\backslash \check G$ over the coset of $1$ contains $f$ implies that $\check{\mathfrak g}_{X,\check\gamma} \hookrightarrow \check{\mathfrak g}_{\check\alpha}+\check{\mathfrak g}_{\check\beta}$ is \emph{antidiagonal} with respect to the pinning of $\check{\mathfrak g}$, that is, its image is spanned by $e_{\check\alpha} - e_{\check\beta}$.

\begin{definition}\label{definition-LGX}
 Let $\check G_X\hookrightarrow \check G$ be embedded as described above. This embedding being stable under the action of $\Gamma = \text{Gal}(\overline{\mathbb F}/\mathbb F)$ on $\check G$, we define the \emph{$L$-group of $X$} as the semidirect product 
 \begin{equation}
{^L G_X}:= \check G_X \rtimes \Gamma.  
 \end{equation}
\end{definition}

 \item Finally, to complete the definition of the rational structure and action of $\Gamma$ on $\check M$, we need to define these on the symplectic vector space $S_X$, compatibly with the rational structure and $\Gamma$-action on $\check G_X$. 
 
Denote by $\mathfrak B_X^S$ the set of colors of even sphere type,
  and by $\mathfrak B_X^D$ its union with the set of $G$-invariant divisors, or equivalently (by Proposition \ref{Gvaluations}): $\mathfrak B_X^D = \mathfrak B_X^S \sqcup  \mathcal D^G(X).$
  The symplectic representation $S_X$ is determined by these colors (\S~\ref{SXdef}, \ref{SX2sec}), but, when $k$ is not algebraically closed, the symplectic form on $S_X$ is not necessarily unique up to $k$-isomorphism. For example, scaling $S_X$ scales the symplectic form by squares, which leaves us with the problem of choosing a square class for it.
  
  In order to do so, note that for every simple root $\alpha$ of $G$ which is also a spherical root of type $T$, the pinning of $\check G$ gives us generators $e_{\check\alpha}$ of the associated root spaces in $\check{\mathfrak g}_X$. We would also like to pick generators of the root spaces associated to the rest of the simple spherical roots $\gamma$ of even sphere type (``of type $\SO_{2n}\backslash\SO_{2n+1}$'' with $n\ge 2$). By \cite[Lemma 10.4]{KnSch}, one can choose such generators $e_{\check\gamma}$ so that they are preserved by the action of the Galois group of Definition \ref{definition-LGX}. Let us fix such a choice.
  
 \begin{remark}\label{dependence-egamma}
  The choice of $e_{\check\gamma}$  introduces  an ambiguity into the rational structure on $S_X$ that we are about to describe.  
 This ambiguity exists, of course, only when there are roots of even sphere type, and moreover
 the choice of $e_{\check\gamma}$ matters only up to square in $k^{\times}$. One could prescribe the generators $e_{\check\gamma}$ more precisely, starting from the pinning of $\check G$, but we have no compelling evidence in order to make this choice. 
 \end{remark}

\begin{lemma} \label{lemmasymplform}
There is a unique (up to isomorphism) form of the symplectic $\check G_X$-vector space $(S_X,\omega)$ over $k$, with the following property: 

There is a  ``subbasis'' (i.e., a linearly independent subset) $(s_b)_{b\in \mathfrak B_X^S}$, with $s_b$ in the weight space corresponding to the color $b$ (cf. Conjecture   \ref{conjsymplectic} and Lemma \ref{vrhopair}) such that 
\begin{itemize}
 \item 
for every simple root $\alpha$ such that $(X,\alpha)$ is of type $T$, with associated colors $b_1, b_2\in \mathfrak B_X^C$, the symplectic form satisfies\footnote{Recall that we are considering \emph{right} actions here, so the highest weight vectors of a representation are annihilated by ``$f$-elements'' rather than ``$e$-elements.'' } 
 \begin{equation}\label{symplrankone}
 \omega(s_{b_i}, e_{\check\alpha} s_{b_j}) = \begin{cases} 1, & \mbox{if } i\neq j; \\ 0, & \mbox{ otherwise;} \end{cases}  
 \end{equation}  
 \item for every other other color $b\in \mathfrak B_X^S$ of even sphere type, corresponding to a spherical root $\gamma$, we have\footnote{\label{bad2footnote} Note that this equation implies that the moment map evaluated at $s_b$ involves $\frac{1}{2}$;
 this is why we needed $\frac{1}{2} \in k$ here.}
 \begin{equation} \label{symplrankonetwo}
 \omega(s_{b}, e_{\check\gamma} s_{b}) = 1.
 \end{equation} 
\end{itemize}

\end{lemma}

Note that the span of $s_{b_1}, s_{b_2}$ is necessarily isotropic, since they are both highest-weight vectors for the copy of $\sl_2$ associated to the root $\check\alpha$. 

\begin{definition}\label{standardsubbasis}
A subbasis as in the lemma will be called \emph{standard}. The same term will be used for any enlargement of it to a subbasis indexed by the set $\mathfrak B_X^D$. 
\end{definition}

\begin{proof}

In the notation of \eqref{eq:nonaffineclosure}, we have an orthogonal decomposition
\[ S_X = S_Y \oplus \bigoplus_{\check\lambda \in \mathcal D^G(X)} T^*V_{\check\lambda},\]
the subspace $S_Y$ is generated under the $\check G_X$-action by the subspaces indexed by $\mathfrak B_X^S$, and the sum of isotropic subspaces $V_{\check\lambda}$ is generated by the rest of the elements of $\mathfrak B_X^D$. The $k$-isomorphism class, as symplectic vector spaces with a $\check G_X$-action, of all summands of the form $T^*V_{\check\lambda}$ is uniquely determined. There remains to determine the $k$-isomorphism class of $S_Y$, hence we will now assume that $S_X=S_Y$.  

Start with \emph{any} $\check G_X$-invariant symplectic form $\omega$ on $S_X$. We will now use the same construction as in the proof of existence of an element of $\check M$ with regular nilpotent image, Proposition \ref{regnilpotent}. Namely, choose a subbasis $(s_b')_{b\in \mathfrak B_X^S}$, and note that the basis elements corresponding to two colors $b_1, b_2$ associated to a simple root $\alpha$ of type $T$ satisfy relations analogous to \eqref{symplrankone}, but with nonzero constants $c_\alpha$ rather than $1$,  
\[ \omega(s_{b_i}', e_{\check\alpha} s_{b_j}') = \begin{cases} c_\alpha, & \mbox{if } i\ne j; \\ 0, & \mbox{ otherwise.} \end{cases}  \]
A similar relation, with a constant $c_{\check\gamma}$, holds for the rest of the colors of even sphere type, each associated to a unique spherical root $\gamma$.

Since the $\check\alpha$'s and $\check\gamma$'s above are all simple roots of $\check G_X$, there is an inner automorphism $\iota: \check G_X \to \check G_X$, defined over $k$, which acts on simple root spaces $\check{\mathfrak g}_{\check\alpha}$ (resp., $\check{\mathfrak g}_{\check\gamma}$) as above by multiplication by $c_\alpha^{-1}$ (resp.\ $c_\gamma^{-1}$). The representation $\iota^*S_X$, now, satisfies \eqref{symplrankone}. Uniqueness is clear, since the subspaces associated to elements of $\mathfrak B_X^S$ generate $S_X=S_Y$.
\end{proof}

Now we discuss the action of the Galois group $\Gamma = \text{Gal}(\overline{\FF}/\FF)$ on $S_X$. Note that $\Gamma$ acts on the sets $\mathfrak B_X^S$ and $\mathcal D^G(X)$, hence on their valuations, compatibly with its action on the cocharacter group $X_*(A_X)$.\footnote{We do not need the Galois group to fix a Borel subgroup here, although this, of course will be automatic over a finite field. In general, we can think of colors as $G$-invariant divisors on $(\mathcal B \times X^\bullet)_{\overline{\mathbb F}}$, where $\mathcal B$ is the flag variety of Borel subgroups of $G$, and this description makes clear that the Galois group acts on them.} We would roughly like to say that the Galois group should act on the aforementioned basis elements $s_b$ as it acts on their indices; however, this is not quite right for varieties such as $\SO_{2n}\backslash \SO_{2n+1}$, $n\ge 2$, with $\SO_{2n}$ nonsplit. (We know that this is not right by comparison with the results of \cite{SaSph}; see \S~\ref{ssPlancherel}.) The ``reason'' is that the unique associated color $b$, in that case, is really just a placeholder for a ``formal difference of $B$-orbits of higher codimension.'' More precisely, for every spherical root $\gamma$ of even sphere type, there are, over the algebraic closure, \emph{two} $B$-orbits $b_+^\gamma$, $b_-^\gamma$ of minimal dimension in $X^\circ P$, where $P$ is the associated parabolic $P$ (so that $X^\circ P/\mathfrak R(P)\simeq \SO_{2n}\backslash\SO_{2n+1}$). If $n=1$, these $B$-orbits are colors, but if $n\ge 2$, they are of codimension $>1$. Choose an arbitrary labeling $b^\gamma_+, b^\gamma_-$ of these two orbits, for each spherical root $\gamma$ of even sphere type that is not of type $T$ (i.e., such that $n>1$), and let ${\mathfrak B_X^D}'$ denote the disjoint union of: (a) the set of colors of type $T$, (b) the disjoint union of the sets $\{b^\gamma_+, b^\gamma_-\}$, over all other spherical roots of even sphere type, and (c) the $G$-invariant divisors on $X$. The Galois group $\Gamma$ acts naturally on the set ${\mathfrak B_X^D}'$, and hence on the free $k$-module on its elements. Denoting by $s_{b'}$ the basis element of this free module corresponding to an element $b'$, we now define $s_b:= s_{b^\gamma_+}- s_{b^\gamma_-}$, when $b$ is the unique color of even sphere type associated to a spherical root $\gamma$ as above. We let $S_X^D$ be the $k$-submodule spanned by those elements $s_b$, as well as by the elements $s_{b'}$ associated to colors $b'$ of type $T$ and $G$-stable divisors. The Galois group, then, acts on $S_X^D$; moreover, we can identify $S_X^D$ as a subspace of $S_X$ via a standard subbasis of Definition \ref{standardsubbasis}.

\begin{conjecture}\label{conjGaloisaction}
 There is an action of $\Gamma$ on $S_X$ by $k$-rational symplectomorphisms, compatible with its action on $\check G_X$ (Definition \ref{definition-LGX}), which extends its natural action on the subspace $S_X^D$. 
\end{conjecture}

Although we have not been able to prove this conjecture directly, using the combinatorics of spherical varieties, there are geometric reasons to believe that it is true. The conjecture is trivial, of course, when for all spherical roots of even sphere type the associated subquotients are of the form $\SO_{2n}\backslash \SO_{2n+1}$ with $\SO_{2n}$ and $\SO_{2n+1}$ split  (i.e., when the $\Gamma$-action on $S_X^D$ is trivial). 

\begin{definition}\label{simpleaction}
 A symplectic action of $\Gamma$ on $S_X$ as in Conjecture \ref{conjGaloisaction} will be called a \emph{simple} action. The $\Gamma$-action on $\check M$ induced from its analytic action on the pair $\check G_X\subset \check G$ and its simple action on $S_X$ will also be called \emph{simple}.
\end{definition}

It is immediate that -- having chosen $e_{\gamma}$ for even sphere roots as in Remark \ref{dependence-egamma} -- the resulting symplectic $\check G_X\rtimes\Gamma$-representation $S_X$, and the resulting ${^LG}$-Hamiltonian space $\check M$ are unique up to isomorphism.

\end{enumerate}

This completes the description of the $k$-rational structure and the simple action of $\Gamma$ on $\check M$, conditionally on Conjecture \ref{conjGaloisaction} -- which we will from now on assume. The resulting $k$-rational form of $\check M$ will be called its \emph{distinguished split form}.
We conclude by verifying that this indeed is a ``pinned hyperspherical variety,'' a notion
described at the start of this subsection \S~\ref{CheckMRat}.

\begin{proposition} \label{invariantregular}
 Assume Conjecture \ref{conjGaloisaction}. Then there is an $m\in \check M(k)$ which is stable under the $\Gamma$-action, and whose moment image is the distinguished element $f\in \mathfrak g^*$ of the pinning. 
\end{proposition}
\begin{proof} 
As in the proof of Proposition \ref{regnilpotent}, we start by choosing an element $m' \in \check M(k)$, of the form (after choosing an invariant identification $\check{\mathfrak g}\simeq \check{\mathfrak g}^*$, which here we also need to require to be $\Gamma$-invariant) 
\[m'= \sum_{b\in \mathfrak B_X^C} s_b + \sum_{\gamma \in \Delta_X \smallsetminus \Delta_G} f_{\check\alpha_\gamma + \check\beta_\gamma} + \sum_{\alpha \in \Delta_{L(X)}} f_{\check\alpha},\]
where the first sum ranges over all the colors of type $T$, the second sum ranges over all simple spherical roots $\gamma$ which are not simple roots of $G$, with associated roots $\alpha_\gamma$ and $\beta_\gamma$, and the third sum ranges over the simple roots of the Levi of $P(X)$. 

The element above is evidently defined over $k$, and fixed under the action of $\Gamma$. Moreover, it was shown in Proposition \ref{regnilpotent} that is has regular nilpotent image under the moment map, of the form $f + f'$, where $f'$ belongs to the simple root spaces indexed by roots $(-\check\alpha)$, with $\alpha$ positive but not simple. In particular, $f+f' \in \Ad^*(\check N^-) f$, where $\check N^-$ denotes the ``negative'' maximal unipotent subgroup of $\check G$. Since both $f$ and $f+f'$ are $\Gamma$-stable and defined over $k$, and $\check N^-$ is unipotent, it follows easily that there is an $n\in (\check N^-(k))^\Gamma$ with $f+f' =  f\cdot n$ (under the right coadjoint action). Then, the element $m = m'\cdot n^{-1}$ is $k$-rational, $\Gamma$-fixed, and with moment image equal to $f$.

 \end{proof}

%% file: rationality.tex
\section{Towards hyperspherical duality} \label{sec:rationality}  \label{goodhypersphericalpairs} 

 The current section, which is tentative or speculative
 at several parts, discusses how
 the results thus far obtained 
 may fit into a future ``ideal picture.''
 
 \begin{itemize}
  
  \item  \S~\ref{quantization} introduces an important mod $2$ characteristic class associated to a Hamiltonian space, the {\em anomaly}, whose vanishing should remove all ``metaplectic obstructions''
  to quantization. 
Our  definition here is a provisional one; further study of examples is required. 

 \item  In \S \ref{hdpC} we 
 describe what an ideal statement of ``hyperspherical duality''
  $$(G, M)  \leftrightarrow (\check{G}, \check{M})$$
would be, and recall how the previous sections
give partial results in this direction. We also describe what some consequences of this ideal statement would be.

\item 

In \S \ref{GMdesiderata} we propose that   hyperspherical pairs $(G, M)$ without anomaly
should admit a distinguished ``split'' form over $\mathbb{Z}$, and in particular over any ring $R$, 
   where the form $G_{/R}$ is the Chevalley split form of $G$. 
We will formulate in \S \ref{usf} a working definition of a certain
class of hyperspherical dual pairs that is well-suited for Langlands
   program -- for motivation, see prior discussion in   \S \ref{hdprings}, \S \ref{CheckMRat}.
    
   \end{itemize}

We emphasize again that  this section should be regarded as providing starting or working definitions, 
  which we leave open to being revised as further computations are carried out -- the main goal being to motivate further research into these questions.

 \subsection{Automorphic quantization and anomaly}
 \label{quantization}
 \label{remark-quantization}

  As mentioned in \S \ref{TFTintro}, we should like to equip each $(G, M)$  with an  ``automorphic quantization'' and a ``spectral quantization.'' 
As is usual in the theory of quantization, there can exist an ``anomaly'' (and here we follow  language that is used in the physics literature)
which means that the automorphic quantization does not exist without the passage to a metaplectic cover of $G$.
For example, the automorphic quantization of $(G=\mathrm{Sp}_{2n}, M=\mathbb{A}^{2n})$
does not exist over a local field without passing from $G$ to its metaplectic cover. 
There is a similar phenomenon on the spectral side that we do not discuss in this paper:
spectral quantization requires, in general, the consideration of {\em twisted} sheaves.

 The notion of ``metaplectic cover''
 is not an algebraic one: e.g.\ the metaplectic cover
 of $\mathrm{Sp}_{2g}(\R)$ is not the real points of an algebraic group.   Since our framework of hyperspherical varieties $(G, M)$
 is entirely algebraic, we should like to find a purely algebraic condition that eliminates the appearance
of metaplectic covers and twisted sheaves.  Again we stress that our chosen condition of Definition \ref{anomalyfree} below is 
\emph{highly provisional}: roughly, we are confident it is the correct condition ``up to isogeny issues,''  
but we have not carefully analyzed these issues; see also \S \ref{whichanomaly}.

We will make free use of Chern classes of representations.
Namely, to any complex representation $V$ of a topological group $H$ we can associate
a Chern class $c_{i}(V) \in H^{2i}(BH, \Z)$.
This is, by definition, the Chern class of the vector bundle over the classifying space $BH$ defined by $V$. We are using topological (Betti) cohomology here, but if $H$ is an algebraic group defined over a field $F$, and the representation $V$ is defined over $F$, we can replace Betti cohomology by \'etale cohomology, by substituting the coefficients with $\Z_l(i)$. (The case $l=2$ is particularly relevant for our discussion here.)

In what follows, we will define an ``anomaly'' criterion of algebraic nature for hyperspherical spaces defined over a separably closed field (actually, we will work with $\CC$), but we motivate it by the following proposition, which gives an algebraic criterion for the splitting of metaplectic covers over local fields:

 \begin{proposition} \label{anomalousautomorphic} Suppose that $F$ is a local
 field. Let  $H \leqslant \Sp(V)$ be an algebraic $F$-subgroup
 of a symplectic group over $F$ 
 such that  there exists a character $\theta: H \rightarrow \Gm$  with 
\begin{equation} \label{Vtheta} c_2(V) =c_1(\theta)^2 \mbox{ in } H^4_{\et}(BH, \Z/2),\end{equation}
where the right hand side denotes the {\'e}tale cohomology
of the classifying stack $BH$ considered as an algebraic
stack over $\Spec F$.  \footnote{We emphasize that this is absolute {\'e}tale cohomology, and not geometric {\'e}tale cohomology.}
   If  $F$ is nonarchimedean and doesn't have
residue characteristic $2$,\footnote{There is an analogous statement without this
restriction, if one first pushes out the metaplectic cover to $S^1$.}
   the metaplectic cover of $\Sp_{2r}(F)$ splits over $H(F)$.
 \end{proposition} 
 
 This proposition will not be used in any significant way so we give the proof in an appendix, \S~\ref{prooflemmametaplectic}.
We use it only as a potential motivation for the following definition:

Let $G$ be a reductive
group over $\CC$, and $M$ a symplectic $G$-variety. The invariant we shall use 
   is
 the $G$-equivariant second Chern class of (the tangent bundle of) $M$, considered mod $2$:  
  \[ c_2 \in H^4_G(M, \Z) \otimes \Z/2.\]
\begin{definition}\label{anomalyfree}
  We shall say that $M$ is:
  \begin{itemize}
  \item {\em strongly anomaly free}, if  
  $c_2= 0$ (mod $2$), and
  \item {\em anomaly free}, if  there exists  $\beta \in H^2_G(M, \Z)$
  such that $c_2 \equiv \beta^2$ (mod $2$). (Note that $\beta$ is an integral cohomology class, not just mod $2$. The analogy of this condition with \eqref{Vtheta} will become clear in the proof of Proposition \ref{calculateanomaly} below.) \end{itemize}
\end{definition}

  Our expectation is that  if $M$ is anomaly-free, it will admit an automorphic and spectral quantization. 
  Partial justification for this expectation comes from Proposition \ref{anomalousautomorphic}
on the automorphic side and Remark \ref{anomalyspectral} on the spectral side. 
 Our definition was motivated by these facts, and also a rather loose parallel
with the idea of a spin-c structure.

\begin{remark}  (Rationality issues in defining anomaly:) The above definition of ``anomaly-free'' is over $\CC$, and one can formulate it over any separably closed field $\FF$, replacing Betti cohomology with \'etale cohomology with $\Z_2$-coefficients. 
It does not automatically imply that the corresponding
statements hold in absolute {\'e}tale cohomology when $(G, M)$ are defined over a subfield $F\subset \FF$.
It is the latter statement that is more directly related to metaplectic splitting, as in  Proposition \ref{anomalousautomorphic}. 
Sometimes one can deduce results over $F$ from those over $\bar{F}$;
see e.g.\   \cite[1.10]{Deligne}  for the semisimple, simply connected case
or Lemma \ref{anomaly rationality} for an explicit result when the situation is defined over $\Z$. 
As with other issues concerning rationality, we do not
understand the situation well,  and simply wish to point
to this as an important question for further study. 
\end{remark}

 \begin{remark} \label{whyrationality}

  Besides Proposition \ref{anomalousautomorphic} there 
  are other reasons to expect that this $c_2$ plays an important role.  
The reduction of $c_2$ mod $2$ arises in the work of the third-named author with A. Abdurrahman
on a closely related topic \cite{AV}; and in physics, the reduction of $c_2$ mod $2$
is again known in certain contexts to be an obstruction to quantization (see in particular~\cite{WittenAnomaly}).

 At the same time, we leave open the possibility that 
 our definition is not the optimal one.  
 For example \cite{RaskinNisyros} formulate
 a different but related condition,  which we will discuss in \S~\ref{whichanomaly}  below.  There are relatively
 few examples to check for hyperspherical varieties, and
 a thorough study of them should reveal the best definition;  we are presently not aware of an example where the two conditions differ.
  
\end{remark}

  The condition of being anomaly-free admits
a readily computable reformulation in
the case of hyperspherical varieties.
Let $(G, M)$ be a hyperspherical variety over $\CC$, 
associated to a datum $H \times \SL_2 \rightarrow G$ and
$H \rightarrow \mathrm{Sp}(S)$, as in Theorem \ref{thm:structure}.
Recall that, in addition to $S$, there is a second symplectic $H$-representation
of interest, namely, $\mathfrak{u}/\mathfrak{u}_+$ in the notation
of \S \ref{uuplus}, that is,   the weight one space of $\mathbb{G}_m \subset \SL_2$
acting on the Lie algebra of $G$.

  \begin{proposition} \label{calculateanomaly}
With notation as described,  let
$T$ be a maximal torus of $H$, let 
$V$ be the $H$-representation $\mathfrak{u}/\mathfrak{u}_+ \oplus S$, 
and let 
 $$\Xi = \mbox{the nonzero weights of $T$ acting on  $V$}.$$
 In what follows, $c_2(V)$ refers to the second Chern class
 of $V$ in $H^4(BH, \Z)$, and $X^*$ to the character lattice for $T$.

 \begin{description}
 \item[(a)] $M$ is strongly anomaly free if 
 $c_2(V) =0$ modulo $2$,   
 equivalently, 
  \begin{equation} \label{qqq}  \sum_{\chi \in \Xi/ \{\pm 1\}} \chi \in 2 X^*\end{equation}

 \item[(b)] $M$ is anomaly free if and only if there is a character
 $\theta: H \rightarrow \mathbb{G}_m$
 such that $c_2(V) = c_1(\theta)^2$ modulo $2$, 
 equivalently,  
   \begin{equation} \label{qqq2}  \sum_{\chi \in \Xi/ \{\pm 1\}} \chi \in (X^*)^W + 2 X^*\end{equation}
i.e., the sum of \eqref{qqq} is congruent
 modulo $2 X^*$ to a character that extends to $H$. 
 \end{description}
 
\end{proposition}

Here the notation in \eqref{qqq} means that we sum over
  an arbitrary set of representatives for $\Xi$ modulo $\pm 1$, noting that $\Xi=-\Xi$ because the representation in question is symplectic; we 
observe that the resulting class in $X^*/2 X^*$ does not depend on the choice of representatives.

\proof 
   The inclusion of $G/H$ into $M$ (as the closed $G\times\GGm$-orbit) is a homotopy equivalence and correspondingly 
the $G$-equivariant cohomology of $M$ is identified with the cohomology of $BH$.
Under this identification, the Chern class of the tangent bundle 
is carried to the   Chern class of the $H$-representation
\begin{equation} \label{obsidian}
(\mathfrak{g}/\mathfrak{h}) \oplus (\mathfrak{g}/\mathfrak{h})^{e} \oplus S\end{equation}
as an $H$-representation; we used \eqref{fund_id} to identify the tangent space. Decompose $(\mathfrak{g}/\mathfrak{h})$ as an $\SL_2 \times H$-representation as $\bigoplus  [m] \otimes W_{m}$
with $[m]$ the $m$-dimensional irreducible representation. Then, as an $H$-representation, 
\eqref{obsidian} is identified with 
$ \bigoplus_{m} W_m^{\oplus m+1} \oplus S$, and this has the same $H$-equivariant second Chern class (mod $2$) as
$ \bigoplus_{m \in 2\Z} W_m \oplus S \simeq \mathfrak{u}/\mathfrak{u}_+ \oplus S$, i.e.
the $V$ defined in the statement of the proposition. 

This proves the first statement of (a). For the first statement of (b)
we note that $H^2_G(M, \Z)$ is similarly identified with $H^2(BH, \Z)$,
which is identified by the Chern class map with $\Hom(H, \Gm)$. Consequently, the condition that the second Chern class
of the tangent bundle is a square of a class in $H^2_G$
is equivalent to the condition that $c_2(V)$ is the square
of $c_1(\theta)$ for some character $\theta: H \rightarrow \Gm$.

To prove the second statements of (a) and (b), that
is to say, the numerical criteria \eqref{qqq} and \eqref{qqq2}, 
we use the following Proposition \ref{Liegroupcohomology}
  to reduce to computing in the maximal torus,
where the computation is straightforward and left to the reader. 
 \qed
    
\begin{proposition} \label{Liegroupcohomology}
Let $H$ be a reductive group over $\C$ with maximal torus $T$.
\begin{itemize}
\item[(a)] 
The restriction map 
$$H^4(BH, \Z) \rightarrow H^4(BT, \Z) $$
identifies the source with the Weyl-fixed part of the target.
\item[(b)] The maps
$$ \Sym^2 X^*(T) \rightarrow \Sym^2 H^2(BT, \Z) \rightarrow H^4(BT, \Z)$$
are isomorphisms. The first map is induced from $X^*(T) \rightarrow H^2(BT, \Z)$
that attaches to a $T$-representation its first Chern class. The second
map is the cup product.  
 \item[(c)]
With reference to these isomorphisms, the second Chern class
of a symplectic $H$-representation with nonzero weights $\Xi \subset X^*(T)$, is given by
$$ \sum_{\chi \in \Xi/\{\pm 1\}} (-\chi^2) \in \Sym^2 X^*(T),$$
\end{itemize}
\end{proposition}
 
 \proof
 (a) is proved in \cite[Theorem 6]{Henriques}.
  (b) reduces to the case $T=\mathbb{G}_m$ by the K\"unneth formula
where it is standard. (c) arises from computing the total Chern class of the restriction of the representation to $T$ as $\prod_{\chi \in \Xi} (1+t\chi)$ (with $t$ the grading variable), and noticing that 
the coefficient of $t^2$ is $\sum_{\chi \in \Xi/\{\pm 1\}} (-\chi^2)$.  \qed

\begin{example}
 If $M$ admits a distinguished polarization, it is anomaly-free:
 In the notation established prior to Proposition \ref{calculateanomaly} the representation $\mathfrak{u}/\mathfrak{u}_+$
 vanishes and  $S=W \oplus W^*$ as $H$-representation.  Then
 the sum of \eqref{qqq} coincides with the 
 character of the determinant of $W$. 
 \end{example}

\begin{example}
If $H$ has the property that its distinguished central element
$(-1)^{2 \check{\rho_H}} \in H$
acts trivially on 
the symplectic $H$-representation $W$, 
 then automatically the condition \eqref{qqq2} is satisfied:
 any such representation  has no symplectic irreducible factor, 
 and so is isomorphic 
to the sum of an {\em even} number of orthogonal $H$-representations,
and a number of representations of the form $E \oplus E^*$. 
 \end{example}

 The example $G=\mathrm{Sp}(V)$ and $M=V$ 
is {\em not} anomaly-free, but 
 it can happen that the restriction of the $\Sp(V)$ action on $M=V$
to a subgroup of $\mathrm{Sp}(V)$ is anomaly-free. 
By Proposition \ref{calculateanomaly} this question can be readily computed in terms of weights.

 \begin{example} 
Some interesting anomaly-free hyperspherical examples
where $M$ is a vector space, 
taken from the table of \cite{KnopMFSR}, are the following:
\begin{equation} \label{MFSR1} \SO_{2n} \times \Sp_{2m} \rightarrow \mathrm{Sp}_{4nm}, E_7 \rightarrow \mathrm{Sp}_{56},  \mathrm{SL}_3 \rightarrow \Sp_{20},\end{equation} 
given, respectively, by the
tensor product of defining representations, the $56$-dimensional fundamental representation 
and the exterior cube of the standard representation; similarly,  also,   the spin or half-spin representations 
\begin{equation}  \label{MFSR2} 
\mathrm{Spin}_{10} \rightarrow \mathrm{Sp}(16), \mathrm{Spin}_{11} \rightarrow \mathrm{Sp}_{32},
\mathrm{Spin}_{12} \rightarrow \mathrm{Sp}_{32}\end{equation}
are nonanomalous. 
\footnote{Note that the fact these examples 
are hyperspherical requires some verification.
The quoted reference shows these examples are coisotropic; but we must also verify that
 the stabilizer of a generic point is connected, which appears to be known in all cases: see \cite[p 81]{SatoKimura}, \cite[Table 1]{GuralnickGaribaldi},
and the splendid \cite[Table 3]{GuralnickLawther}.   }
 In all these cases except the final one the dual $\check{M}$ is known
at least in an isogenous example
 (some are listed in \S \ref{IntroEx}, see also the tables of C. Wan and L. Zhang in \cite{WanZhang}). 
 \end{example}
 
 \begin{example}
 Several other examples in Table 1.1 of \cite{KnopMFSR} are, however, anomalous,
 for example the standard representation of $\SO_{2n+1} \times \Sp_{2m}$,
 and the $14$-dimensional representation of $\SL_2 \times G_2$. An interesting non-anomalous example that {\em fails} to be hyperspherical,
 because it fails the connectedness criterion, is the action of $G = \SL(2)$ on $M = \mathrm{Sym}^3 (\mathrm{std})$. 

 \end{example}

\subsubsection{Relation with the anomaly condition of \cite{RaskinNisyros}.} \label{whichanomaly}

Proposition \ref{calculateanomaly} explicates the anomaly condition in terms of the symplectic representation $V$ of the subgroup $H\subset G$ (of the structure theorem). The paper \cite{RaskinNisyros} introduces a different ``anomaly vanishing'' condition for a symplectic $H$-representation $V$; without going into much detail, we reformulate their condition as follows:

\begin{definition}\label{anomalyfree-Coulomb}
 $M$ satisfies the anomaly-vanishing condition of \cite{RaskinNisyros} if the pullback of $V$ to the simply-connected cover $H_{sc}$ of the derived group of $H$ is anomaly-free in the sense of Proposition \ref{calculateanomaly}, that is to say, 
 $c_2(V)$ vanishes mod $2$ in the cohomology of $H_{sc}$.  \end{definition}

It might as well be that this is the ``correct'' anomaly-vanishing condition; we presently do not know any examples of hyperspherical varieties where the two conditions differ.  
 In any case, we can directly confirm that the dual of a 
 tempered spherical variety is anomaly-free in the sense of \eqref{anomalyfree-Coulomb}:
 
\begin{proposition} 
\label{dual is Coulomb nonanomalous} Assuming Conjecture \ref{conjsymplectic}, if $P(X)=B$, the dual Hamiltonian space of the spherical variety $X$ satisfies the anomaly-vanishing condition of Definition \ref{anomalyfree-Coulomb}. 
\end{proposition}

\begin{proof}
In the notation of Proposition \eqref{calculateanomaly}, but with $\check G_X$ in place of $H$, it is enough to show that the sum $\sum_{\chi\in \Xi/\{\pm 1\}} \chi$ (mod 2) vanishes on the coroots of $\check G_X$. Indeed, its pullback to the simply connected cover of the derived group of $\check G_X$ will then be  trivial (mod 2). 

Since we have assumed that $P(X)=B$, the symplectic representation $V$ coincides with the representation $S_X$ of Definition \ref{def:SX2}. The only summand of \eqref{eq:nonaffineclosure} for which this could fail is the summand $S_Y$. In the notation of Conjecture \ref{conjsymplectic} (but replacing $X$ by $Y$), the sub-sum of $\Xi$ that corresponds to $S_Y$ can be written as $\chi_Y = \sum_{b\in \mathfrak B_Y^+} \textrm{wt}(b)$. Finally, we have, for every simple coroot $\gamma$ of $\check G_X$, with $w_\gamma$ the corresponding Weyl reflection, that  
\[  \chi_Y - w_\gamma \chi_Y = \langle \chi_Y, \gamma \rangle \check\gamma.\]
By Conjecture \ref{conjsymplectic}, $w_\gamma\chi_Y = \chi_Y$, unless $(Y,\gamma)$ is of type $T$ with associated valuations $\check v_1, \check v_2$, in which case $w_\gamma\chi_Y = \chi_Y - 2 (\check v_1+ \check v_2) = \chi_Y - 2 \check\gamma$. In particular, 
$\langle \chi_Y, \gamma \rangle  \in 2\mathbb Z$. 
\end{proof}

  \subsection{Hyperspherical dual pairs over $\C$}  \label{hdpC}
  
The most ideal form of hyperspherical duality would be:  
 \begin{expectation} \label{anomaly expectation}
  There exists a  bijection
\[ (G,M) \leftrightarrow (\check{G}, \check{M}),\]
 between isomorphism classes of anomaly-free hyperspherical $(G, M)$ over $\C$
with anomaly-free hyperspherical $(\check{G}, \check{M})$ over $\C$, with the following properties:
 \begin{quote}
\item If $M$ admits a distinguished polarization $M=T^*(X, \Psi)$, as in Definition \ref{distinguishedpolarizationdef}, then $\check{M}$ arises from  $(X, \Psi)$ via the procedure of \S \ref{dualofX},  and vice versa.
\end{quote}
\end{expectation}

We must admit at the moment that this expectation remains somewhat tentative: we would not formulate it as a ``conjecture,''
  but we believe that something like it should be true, perhaps after slight modifications of the
  definitions of ``anomaly-free'' or ``hyperspherical.'' 

\index{hyperspherical dual pair}
In the remainder of this paper, we will use the phrase {\em hyperspherical dual pair} (over $\C$)
to mean either a pair as $(G, M), (\check{G}, \check{M})$ which
arises via the construction of \S \ref{dualofX}, or its reverse. 
We anticipate, however, that 
all the statements of the paper will apply to the
class of dual pairs of Expectation \ref{anomaly expectation}.

\subsubsection{Some consequences of Expectation \ref{anomaly expectation}} \label{anomexpcon}
We will explicitly note several consequences of the expectation.  They are  are valid in examples that we have checked,  but we have no general proof:
\begin{itemize}  

\item[(i)] The construction of \S \ref{dualofX} is independent of distinguished polarization.  
It is likely this is provable using ideas from the theory of spherical varieties. 

\item[(ii)] Any $(\check{G}, \check{M})$ arising via the procedure of \S \ref{dualofX} from $(G, M=T^*(X, \Psi))$ is non-anomalous.  
For what we have proved in this direction, see Proposition \ref{dual is Coulomb nonanomalous}. 

\item[(iii)]  If $(G , M)$ and $(\check{G}, \check{M})$ 
 are a distinguished hyperspherical pair, where
 both sides admit distinguished polarizations, and $(\check{G}, \check{M})$
 arises from $(G, M)$ via the procedure of \S \ref{dualofX}, then
 the reverse is also true: $(G, M)$ arises from $(\check{G}, \check{M})$
 from the procedure of  \S \ref{dualofX}.  
 \index{$\eta_{\aut}$ automorphic-side eigenmeasure}
 
 \item[(iv)] In the setting of  $M=T^*(X, \Psi)$ and $\check{M}$ as in \S \ref{dualofX},  we have proved in  Proposition \ref{zXparity} \index{parity element}
 \begin{equation}  \label{parity1} \mbox{parity: } e^{2 \check{\rho}} \eta(-1) \in \check G \mbox{ acts trivially on $\check M$.} \end{equation} 
where $\eta: G \rightarrow \Gm$ is the character of a $G$-eigenmeasure on $X$, identified with a dual central cocharacter in $\Gv$.

Expectation \ref{anomaly expectation}   implies that the  parity condition \eqref{parity1}
 also holds ``in reverse'' for a   hyperspherical dual pair if $\check{M}$ also admits a distinguished polarization:
\begin{equation}  \label{P2} \mbox{dual parity: }  e^{2 \check{\rho}} \check\eta(-1) \in G  \mbox{ acts trivially on $M$.} \end{equation} 
where $\check \eta: \check G \rightarrow \Gm$ is the character of a $\check G$-eigenmeasure on $\check X$, again identified to a cocharacter for $G$.\footnote{This statement
doesn't depend on the choice of eigenmeasure, cf.   \S \ref{ssseigencharinocuous} .}
\index{$\eta_{spec}$ spectral-side eigenmeasure}

\end{itemize}

 \begin{remark}  \label{noeigenmeasure} 
Although
we have often used them as a crutch, no special role in the duality should be played by the existence
of eigenforms as in \S \ref{eigencharacter}. 
A typical example of a dual pair, where one side has no eigenmeasure, is given by
 $$ (G, (X, \Psi))= (\PGL_3, \Gm \cdot (\Ga^2, \psi)\backslash \PGL_3) \mbox {and }  (\check{G}, \check{X})  = (\SL_3, \mathbb{A}^3).$$
This is related to the Example \ref{ex:SX3}.
  One can, in such cases, still apply the local and global conjectures
by reducing them to a case with an eigenform, as described earlier in 
\S \ref{ssseigencharinocuous} ; but it would be
good, in a further elaboration of this work, to do this in a more intrinsic way.  See e.g.\ Remark \ref{noeigenmeasure2}.
 \end{remark}

\begin{remark}

\begin{itemize}
\item[(a)]
The putative duality requires some kind of anomaly vanishing condition, unless one modifies the nature of the duality;
in the classical theory of automorphic forms the 
appearance of metaplectic groups necessitates a modification
even for the duality on reductive groups, cf. \cite{Weissman}. It is an 
important problem to study anomalous examples  --there are many in the automorphic literature, that is to say period integrals involving covering groups; let us point only to Shimura's integral \cite{ShimuraIntegral} -- and extend the proposals
here to that context.  

\item[(b)] We do not know of examples of dual $(M, \check{M})$ where {\em neither}
side admits a twisted polarization, but there seems to be no reason for such examples not to exist.  It would 
interesting to exhibit one.
\end{itemize}
\end{remark}

 \subsection{Hyperspherical dual pairs over arithmetic fields.}
  \label{GMdesiderata} \label{dhpFq}

  We have already seen in our previous discussions  of \S \ref{hdprings} and  \S \ref{CheckMRat}
  that it is important to consider 
   models for $(G, M)$
over other fields.  
We propose
that  in the non-anomalous
case, at least, there is a best one:
 \begin{quote} {\em  There exists a distinguished ``split'' form of  each hyperspherical $(G, M)$, defined over $\Z$.}
 \end{quote}

 In this subsection, after discussing this 
proposal, we will use it to give a working
definition of a class of dual pairs
$(G, M)_{/\FF} \leftrightarrow (\check{G}, \check{M})_{/\kk}$
where 
 $\mathbb{F}$ is either a finite field or $\C$, and $\kk$ is either $\CC$ or the closure of an $\ell$-adic field.  
  
Our primary motivation to understand the theory over general $\FF$ arises
from automorphic forms. 
We will discuss this motivation later, in 
 \S~\ref{pairs-motivation}. Perhaps the main takeaway is that, 
 for $\FF$ not algebraically closed, the hyperspherical
 datum over $\FF$
  best adapted to automorphic phenomena is not always the most obvious one.

 At this point, we should remind the reader what was done in  \S~\ref{hdprings} and \ref{CheckMRat}: In the former, 
 we used the structure theorem to define a notion of ``hyperspherical datum/scheme'' over (more) general rings, and to describe, in some cases, a ``distinguished split form'' of those. In the latter, we used the structure theorem to describe a split form of the dual Hamiltonian space of a spherical variety, depending on some mild choices. Here, we will combine these discussions into a wishlist for the ``distinguished split form of $M$,'' that will also make some forward references to our local conjecture.
We fix a pinning in order to rigidify $G$  (thus, the notion of a distinguished split form of $M$ is really to be understood with reference to a pinned group $G$). 

\begin{expectation}\label{GMdesiderata expectation}
Each nonanomalous hyperspherical pair $(G, M)$ over $\C$
admits a distinguished $\Z$-form $(G, M)_{\Z}$
with the following properties.

 \begin{itemize}
 \item[(a)]   
 Write $\Z'$ for the ring obtained from $\Z$
by inverting $N_G$, as in \eqref{NGdef}. Then
 $(G,M)_{\Z'}$ corresponds to a split $\Z'$-form of the linear algebra datum $\mathcal{D}(G, M)$
 (see Definition \ref{ringdatum}).

    \item[(b)]  Suppose that $(G, M), (\check{G},\check{M})$ form 
 a hyperspherical dual pair over $\C$, with $\check{M}=T^*X$ polarized. 
 Then, for any field $k$ not of characteristic $2$, the  base change $(G, M)_k$ to $k$ of the distinguished $\Z$-form
 belongs to the distinguished class of forms constructed  in
  \S \ref{CheckMRat},    switching $(G, M)$ and $(\check{G}, \check{M})$ in that discussion. 
  (cf. also Proposition \ref{invariantregular}).

    \item[(c)]   Suppose that $M, \check{M}$ form 
 a hyperspherical dual pair, with $\check{M}=T^*X$ polarized;
 then $\Z[\check{M}]$ is the local Plancherel algebra,
 see Remark \ref{localconjecturesplitform}.   

\end{itemize}
  \end{expectation}

 \begin{remark} \label{Whydesid}
 
 \begin{itemize}
 \item
 Points (a) and (b) are suggested   by the study of examples known to us. 
 The word ``split'' in (a)  is suggested by the work \cite{GaitsgoryNadler} wherein a maximal torus of $H$  -- the  reductive subgroup
  of $G$ appearing in the structural data for $M$, as in \eqref{basic0}  -- 
  appears from degenerating to a horospherical variety.  
  \item We are perhaps being overly cautious 
  in (a) and (b); perhaps it is unnecessary to invert {\em all} primes dividing $N_G$. 
 An example illustrating the difficulties at $p=2$ is the case $G=\mathrm{Sp}_4$
 when $M$ is defined by the datum $(H=\SL_2, S=\mathbb{A}^2)$
and the auxiliary $\SL_2$ is the centralizer of $H$ in $G$.

  \item
  
The strongest reason to believe in this expectation 
is the one of point (c), namely,  the local conjecture that we are about to formulate in
\S \ref{section-unramified-local}.  In favorable circumstances,
it  gives rise to an explicit ring with Poisson bracket
 which should be the coordinate ring of the distinguished  split form of $M$ over $\Z$. 
 
 Indeed, recall that one way of constructing the split form of a group is provided by the geometric Satake correspondence:
the split form of $G$ over $\mathbb{Q}$ and even over $\mathbb{Z}$, can be reconstructed from the category of sheaves
of the affine Grassmannian of $\check{G}$. Point (c) is 
an analogous proposal -- 
we will see in \S \ref{Poisson from loop} that $M$ can be  conjecturally reconstructed as a $G \times \GGm$-space  
given access to a polarized dual space $\check{M} =T^*(\check{X}, \check{\Psi})$
by considering a suitable category of constructible sheaves. 
This category of constructible sheaves
 can be defined with $\Z$ coefficients when $\Psi$ is trivial (simply by considering sheaves of abelian groups), in particular giving 
  a $\Z$-structure to $M$.   \footnote{In fact, a similar argument 
  applies when $\Psi$ is nontrivial. However,   it does not give rise to models over $\Z$, but
  rather over completions of various cyclotomic integer rings.  }
   \end{itemize}
  \end{remark}

\begin{remark}
For the purposes of applications to automorphic forms,  
the very best situation would not be 
to try to cherry-pick our favorite form, 
but rather find an enhanced 
 version of the duality
\begin{equation}\label{hyperdual}(G, M) \leftrightarrow (\check{G}, \check{M})\end{equation}
that took into account rational structures and Galois actions on both sides. 
In the case of $G \leftrightarrow \check{G}$, a form of 
 $G$ over a non-separably closed field $\FF$ provides an action of the absolute Galois group of $\FF$
 on the dual group $\check{G}$; this datum is used to define the $L$-group.  
   An analogous Galois action in the case of the duality \eqref{hyperdual} 
 is provided, in part, by the analysis of \S \ref{simple-M}, see Definition \ref{simpleaction}.
Such an action of the Galois group should naturally arise from the construction of the duality \eqref{hyperdual} in the general case.
 \end{remark}

\subsubsection{Some consequences of Expectation \ref{GMdesiderata expectation}} \label{GMdesidcons}
 Expectation \ref{GMdesiderata expectation}
has various explicit consequences which, again, we do not know how to prove in general,
but are valid in those cases we have examined. 

\begin{itemize}
\item[(i)] The existence of a $\Z$-form  proposed in Expectation \ref{GMdesiderata expectation}  implies that the isomorphism class of $(G, M)$ is stable by any field automorphism of $\C$,
and in particular $(G, M)$ can be unambiguously transferred to any algebraically closed field of characteristic zero
by the Lefschetz principle. 

\item[(ii)] Point (c) also has very strong implications (perhaps too strong?) The local conjecture relates the ring of regular functions on $M$ to certain constructible cohomology groups;
thus, for $\Z[\check{M}]$ to be flat over $\Z$,  these particular  cohomology groups should have no torsion. See
Remark \ref{localconjecturesplitform} for further discussion. 
\end{itemize}
 
\subsubsection{Characteristic $2$ subtleties, quadratic refinements}
\label{char2moment} 

  Let us consider the case when $M=V$ is a vector space
and look at what we expect about the integral model of $M$;  in this case, we expect, of course,
this to be simply a $\Z$-lattice $V_{\Z}$
equipped with its standard symplectic form, 
equipped with an action of the Chevalley form $G_{\Z}$.

However, we expect (hope?) this to have a further property, namely:

\begin{quote} (*)
the associated   representation of $\mathbb{F}_2$-algebraic groups 
$ G_{\mathbb{F}_2} \rightarrow \mathrm{Sp}(V_{\mathbb{F}_2})$ should preserve a quadratic refinement
of the symplectic form,
\end{quote}
i.e., $G_{\mathbb{F}_2}$ preserves a quadratic
function $Q$ on $V_{\mathbb{F}_2}$ such that $Q(x+y)-Q(x)-Q(y)$ gives the symplectic form. 

There are two reasons to expect this.
Firstly, this allows one to construct a moment
action for the action of $G_{\Z}$ on $V_{\Z}$.
The second is topological, arising
from the existence of extra mod $2$ operations
on the the cohomology of $E_3$-spaces (cf. \S \ref{spectral-factorization}).

Not all hyperspherical $G$-representations
will admit a $\Z$-lattice with this property. 
For example, the standard representation of $G=\SL_2$
on $M=\mathbb{A}^2$ does not have this property.
However, it is plausible that all anomaly-free
examples have property (*), and in some cases
it helps distinguish the appropriate form. 

For example, take $G=\SL_2$
acting on $M = \mathbb{A}^2 \times \mathbb{A}^2$,
there are two possible $\Z$ forms.
Namely  take
$$ M_{\Z} = (\Z^2, \omega) \otimes (\Z^2, B)$$
with $\omega$ the standard symplectic structure, and
for some symmetric bilinear unimodular form $B$.
 The only possibilities for $B$ are, up to equivalence,
 represented by the matrices $\left[{\small \begin{array}{cc} 0 & 1 \\ 1 & 0 \end{array}}\right]$ and the identity matrix; the first possibility satisfies (*) and the second does not.

 \subsubsection{Working definition of hyperspherical dual pairs over $(\FF, \kk)$}  \label{usf}  \index{distinguished split form of dual pair}

Although it is not clear that the different desiderata in Expectation \ref{GMdesiderata expectation}
are compatible with one another, they are quite restrictive when they apply.
We will focus here on (a) of that expectation. The ``distinguished split form'' introduced in Definition \ref{dhpFqdef}, whenever it applies, is necessarily the reduction of the conjectural distinguished split form to an appropriate finite field, and this permits us to transpose the notion of hyperspherical dual pairs, as in \S \ref{hdpC},
to pairs over arithmetic fields, suitable for our applications to the Langlands program.

 Suppose that $( G \times \GGm, M)$ and $(\check{G} \times \GGm, \check{M})$ is a hyperspherical dual pair over $\C$, in the sense of \S \ref{hdpC}.  
   Let $\mathbb{F}$ be either  $\mathbb{F}_q, \overline{\mathbb{F}_q}$ or $\C$, and let $\kk$ be either $\CC$ or the closure of an $\ell$-adic field. 
Then, by a 
  {\em distinguished split form} of the above pair over $(\mathbb{F}, \kk)$, we shall mean
 $$
(G \times \GGm, M)_{/\mathbb{F}} \mbox{ and }  (\check{G} \times \GGm,  \check{M})_{/\kk}$$  
where,   if $\mathbb{F} \neq \C$, the left hand side is a split form as defined in (a) above,
and the right hand side   is obtained from  $(\check{G} \times \GGm, \check{M})$  by means of an isomorphism $\kk \simeq \C$.\footnote{ It would follow
from \S \ref{GMdesidcons} (i) that this is independent of choice of isomorphism.}  
In the case that $M$ or $\check{M}$ come with distinguished polarizations, we can similarly define a split form of the pair equipped with their polarizations.

 \subsubsection{Some automorphic examples and motivation}\label{pairs-motivation}

The following discussion presupposes
 some familiarity with the conjectures in the remainder of the paper;
it is motivational and can be skipped.  The main takeaway from the section is that
the ``best'' form of $(G, M)$ may not be the most obvious one. 
A related discussion, but on the spectral side, is given in \S \ref{CheckMRat}. 
 
In the Langlands program
we are concerned with reductive groups $G$ over arithmetic fields $F$ of several types --
for example,
$F$ could be  the function field of a projective curve $\Sigma$ over a finite field, 
or a number field, or  a local field. Correspondingly, we want to be able
to work with $(G, M)$ over the same types of fields. 
Now, 
 the automorphic data corresponding to $M$ is sensitive to the form of $M$
 over $F$, and not only to its isomorphism class over $\bar{F}$. 
 Here is an example:

  \begin{example}
 Let $F$ be a local field. The space $X_d$ of $2\times 2$-matrices of determinant $d \in F^\times$, under the action of $G=\SL_2\times\SL_2$, is isomorphic to the space $X=X_1=\SL_2$ over the algebraic closure, but not necessarily over $F$. The set of $X$-distinguished irreducible representations $\pi\hookrightarrow C^\infty(X(F))$ consists of those of the form $\pi =\tilde\tau\otimes \tau$, while for $X_d$ they are of the form $\tilde\tau^d \otimes \tau$, where $\tilde\tau^d$ is the twist of $\tilde\tau$ (= the contragredient of $\tau$) by the automorphism given by conjugation by $\text{diag}(d,1)$. For the purposes of the Langlands parametrization, this automorphism does not change the $L$-parameter, but can act nontrivially on the elements of an $L$-packet. 
 
 In the language just introduced, both $X, X_d$ are defined by polarized hyperspherical data, with trivial $\SL_2$;
 what differs is the $F$-form of the embedding $H \rightarrow G$. 
\end{example}

A point of crucial interest to us is that  the form of $M$ or $X$ that interacts
in the cleanest way with the local Langlands conjecture may not be the obvious one. 
The following example can be considered a generalization of the $d=-1$ case of the previous one:

\begin{example} \label{group period example}
Consider the group $G=H\times H$, where $H$ is a quasisplit reductive group, and both copies of $H$ are assumed to carry the same pinning $f_H \in \mathfrak h^*$. Let $X$ be the $G$-space $H$, i.e.,
 $X$ is the quotient $\Delta H \backslash H \times H$
 and let $X'$ be the following form of $X$: 
 \[ X' = \Delta' H \backslash (H\times H),\]
 where 
 \[ \Delta'H= \ \text{graph of the inner automorphism} (-1)^{\check{\rho}_H}: H \rightarrow H.\]
 (The half sum of positive coroots $\check\rho_H$ is considered as a cocharacter into the group of inner automorphisms of $H$ via the pinning fixed.)
 Note that this inner automorphism (let us denote it by $\iota$) is the one connecting the Chevalley and duality involutions (see \S~\ref{dualityinvolution}): 
$h^d = (h^c)^\iota.$

As we now explain, it is functions on $X'$, rather than functions on $X$,
which look more natural on the spectral side of the Langlands correspondence:  

Consider a tempered local $L$-parameter $\phi$ for $H$; according to the local Langlands conjectures as refined by Vogan \cite{Vogan-LL},  the choice of a Whittaker model for $H$ (which is afforded by a combination of the pinning with a character of the additive group of the field) is supposed to be fixing a basepoint in the Vogan $L$-packet associated to $\phi$, and a parametrization of the elements $\tau$ of this $L$-packet by the irreducible representations $\eta$ of the component group of the stabilizer $\check H_\phi$ of $\phi$ in $\check H$. The pair $(\phi,\eta)$ can be called the Langlands--Vogan parameter of $\tau$. The (pinned) Chevalley involution $c$ of the dual group acts on the set of Langlands--Vogan parameters, and a naive hypothesis would be that $(\phi^c, \eta^c)$ is the Langlands--Vogan parameter of the contragredient representation $\tilde \tau$. This is wrong, however, and in \cite[Conjecture 2]{Prasad-contragredient} Prasad corrected this naive hypothesis by postulating that this is the Langlands--Vogan parameter of $\tilde\tau^\iota$, the twist of $\tilde\tau$ by the involution $\iota$. But the set of tempered representations $\tau\otimes \tilde\tau^\iota$ is precisely the set of representations appearing in the Plancherel formula for the space $L^2(X')$, where $T^*X'$ is the pinned hyperspherical space as above!

In terms of the language just introduced, both $X$ and $X'$ are defined by polarized hyperspherical data, with trivial $\SL_2$ and $\rho$,  and with $H=G$; what
 differs is the embedding $\iota: H_{/F} \hookrightarrow G_{/F}$.  As we just saw, $X'$ is in some respects more natural, although the associated
 embedding $\iota'$ involves a ``strange'' conjugation by an involution. 
 
To conclude discussion of this example, let us observe that  $X$ and $X'$ are distinguished from one another in the following abstract way:
$(f_H, -f_H)$ lies in the $F$-points of the image of the moment map for $T^*X$, but the point $f=(f_H, f_H)$ may not be in the image; on the other hand $M'=T^*X'$ contains $f$ in the ($F$-point) image of the moment map. 
This notion will be formalized later in \S \ref{CheckMRat}: we will say that $M'$ is a {\em pinned} hyperspherical space. 
\index{pinned hyperspherical space}

\end{example}

%% file: shearingPart2.tex
  In Part 2 of this work we formulate and study the local form of our conjecture in the unramified setting. 
  For an overview of this part, see
  page \pageref{part2intropage}.

  \section{Shearing and geometric Satake.}\label{shearingsec0}
  \index{shearing}
 In this section we discuss the operation of {\it shearing} with respect to $\Gm$-actions, in which the weights of the $\Gm$-action are paired with cohomological shifts. 
This concept will arise throughout the paper; it already arises in the geometric Satake isomorphism, 
 as we will recall in \S  \ref{sssderivedSatake}. 
Shearing is implicitly present throughout the Koszul duality literature. It is studied explicitly in the work of Arinkin and Gaitsgory -- see  ~\cite[Section A.2]{ArinkinGaitsgory}, but note that~\cite{ArinkinGaitsgory} shear by 2, while we shear by 1, hence the appearance of super-signs.

The contents of the section are as follows:
\begin{itemize}
\item  \S \ref{commutativity} discusses shearing of vector spaces. 
\item \S \ref{shearalgebra} discusses shearing of algebras.
\item \S \ref{shearcategory} discusses shearing of categories. 
\item \S \ref{shearing geometry} discusses various examples of shearing
on categories of geometric or representation-theoretic origin.

\item \S \ref{shearSatake} discusses abelian geometric Satake by way of example. Although
shearing does not appear in the usual formulations of this, we take
the opportunity, by way of illustration, to explain  how some of the subtle features there can be expressed in terms of shearing.
We will separately discuss an ``analytic'' and an ``arithmetic'' form, cf. \S \ref{analyticarithmetic}. 
\item \S  \ref{sssderivedSatake} discusses derived geometric Satake,
again, presented as an example of the shearing language. 
\item  \S \ref{Satake finite field} discusses  both abelian and derived Satake when
the base field is replaced by a finite field.
\item \S \ref{Mshear1} discusses the example of shearing the coordinate ring
of a hyperspherical variety. This could be subsumed in the previous section, but we isolate
it because of its later use in the paper. 
 \end{itemize}

\subsection{Shearing of vector spaces}\label{commutativity}

There is an autoequivalence of the category $\Rep(\Gm)$ 
of $\Gm$-equivariant complexes of $k$-vector spaces (i.e., complexes of graded vector spaces), 
  called \emph{shearing}
defined by
\[
M = 
\bigoplus_i M_i \mapsto M^{\shear} := \bigoplus_{i \in \Z} M_i[i]
\]
where $M_i$ is the $i$-isotypical space, 
upon which $\lambda \in \Gm$ acts by $\lambda^i$; 
this has  inverse (unshearing)
\[
N = \bigoplus N_i \mapsto N^\unshear := \bigoplus_{i \in \Z} N_i[-i]
\]

 Note that $\unshear$ has the property that it takes an ordinary graded vector space (i.e., a graded dg-vector space concentrated in cohomological degree $0$) to a dg-vector space for which the weight on the cohomology agrees with the cohomological degree.

The equivalence $M \mapsto M^{\shear}$ is a {\em monoidal} autoequivalence, i.e., there is a natural identification
$$ (M \otimes N)^{\shear} \simeq M^{\shear} \otimes N^{\shear}$$
using the corresponding property of translation (\S \ref{shift}). We moreover have 
$$ \Hom(M, N)^{\shear} \simeq \Hom(M^{\shear}, N^{\shear}),$$
 if we impose e.g.\ $N$ finite-dimensional (the $\Hom$s here are not $\Gm$-equivariant, hence carry a natural $\Gm$-action, besides the cohomological grading). 
 
The even iterates $M\mapsto M^{2n\shear}$ of the shearing have the natural structure of symmetric monoidal autoequivalences. 
The shearing operation itself $M\mapsto M^{\shear}$ is not symmetric monoidal (nor are its odd iterates), because of the Koszul rule of signs.
Thus we ``correct'' the shearing functor as in \S \ref{shearing notation} by replacing the shift $[i]$ with the parity-corrected shift $\Pi^i[i]$, giving an endofunctor of $\Rep^\super(\Gm)$, 
\[
M = 
\bigoplus_i M_i \mapsto M^{\shear} := \bigoplus_{i \in \Z} M_i\la i\ra 
\]
which does admit a natural symmetric monoidal structure. Restricting to even complexes $\Rep(\Gm)\subset \Rep^\super(\Gm)$ (applying ``Galois descent'', see Remark~\ref{super descent}) we find a symmetric monoidal equivalence 
$$(-)^\shear: \Rep(\Gm)\longrightarrow \Rep_\epsilon^\super(\Gm)\subset \Rep^\super(\Gm)$$
with the full subcategory of those super graded chain complexes in which the parity is given by the action of $\epsilon=-1\in \Gm$ (i.e.,   representations of odd $\Gm$-weight are odd vector spaces, and those of even $\Gm$-weight are even). 

\index{$\Rep_\epsilon^\super(\Gm)$}

With this convention, {\em shearing preserves traces}:  given an automorphism
$\alpha: M \rightarrow M$ of a bounded complex $M$, the
super-trace $\mathrm{tr}(M)$ satisfies $\mathrm{tr}(M) = \mathrm{tr}(M^{\shear})$,
since, in shearing the $i$th graded piece by $i$, we get a factor of $(-1)^i$ by the cohomological shift,
and another factor of $(-1)^i$ from switching between even and odd vector spaces.
As usual, the trace of an automorphism on a super-vector space
is understood to be the even trace minus the odd trace.

\begin{remark}\label{shear with Frobenius}
In \S \ref{shearing notation} we introduced a convention that, in a Frobenius-equivariant context, the shift operation $\la 1 \ra=\Pi [1](1/2)$ also includes
a Tate twist, and hence so does shearing. In the current context we would formalize this as follows: 
there is a monoidal autoequivalence of the category of representations
of $\GGm \times \langle \mathrm{Frobenius} \rangle$
where we additionally twist the action of Frobenius on the $n$-th graded piece by $q^{-\frac{n}{2}}$.   
One can readily transpose the discussion in this chapter to that setting. 
 \end{remark}

\subsection{Some motivation for shearing} \label{shearalgebra}
By way of motivation for what follows, let us describe a situation where shearing of algebras naturally arises  
(which indeed reflects the way it occurs in the main text), and also 
why it is natural to shear categories too.  The experienced reader can skip this section without loss. 

We will often encounter situations where there is an equivalence of triangulated categories
\begin{equation} \label{XexampleX} \mbox{a suitable category of constructible sheaves on a variety $X$} = D(\mathsf{A}),\end{equation}
where we have, on the right, a derived category of modules for
some   differential graded $k$-algebra
$\mathsf{A}$ -- e.g.,  $\mathsf{A}$ may arise as endomorphisms of a suitable object. 
Now, in a situation where $X$ is over a finite field, one often 
obtains a weight decomposition on the $\Hom$-spaces on the left, and in favorable cases,
this arises from a grading on  $\mathsf{A}$ itself, i.e., 
a decomposition
 $  \mathsf{A} =  \bigoplus \mathsf{A}^i_w$
 with differential increasing $i$ but preserving $w$.  
 In this situation one often has ``purity,'' that is to say, 
 \begin{quote} $H^i(\mathsf{A})$ has weight grading entirely in degree $w=i$\end{quote}

In this case, the shear  
$ \mathsf{B} := \mathsf{A}^{\shear}$ by the weight grading is {\em entirely in degree zero},
and therefore can be considered as a usual graded ring (``usual'' means that there is no differential to worry about). 
 We may then seek to describe the category of sheaves on $X$  in terms of the usual  ring $\mathsf{B}$.  

At the level of {\em graded} derived categories the answer is quite simple:
there is an equivalence
of the graded derived category of $\mathsf{A}$ and $\mathsf{B}$:
\begin{equation} \label{modelgradedequivalence} D^{\mathrm{gr}}(\mathsf{A}) \simeq D^{\mathrm{gr}}(\mathsf{B}), M \mapsto M^{\shear}. \end{equation}
To be explicit, this takes
a graded $\mathsf{A}$-module  $M = \bigoplus M^i_w$,
where $w$ is the grading variable so that again the differential increases $i$ and preserves $w$,
and then associating it to the graded $\mathsf{B}$-module
given by $M^{\shear} = \bigoplus M^i_w[w]$, i.e., 
  regrading $M^i_w$ to be in degree $i-w$.

 In the geometric situation described above,  $D^{gr}(\mathsf{A})$ 
 will typically describe a category of mixed sheaves on $X$. Correspondingly,
 this mixed category can be described as the category of graded modules for a (usual, underived) ring.

  Note  $\Hom_{D(\mathsf{A})}(M, N)$ and $
\Hom_{D(\mathsf{A}^{\shear})}(M^{\shear}, N^{\shear})$ do not coincide -- the latter 
is the shear of the former.
Indeed the {\em ungraded} categories $D(\mathsf{A})$ and $ D (\mathsf{B})$ need not be equivalent.
 They can be recoved from each other by an abstract categorical process of ``shearing,''
 which will be described in the next section, and will be written like this:
 $$ D(\mathsf{B}) = D(\mathsf{A})^{\shear}.$$
  In other words, there is an operation $\mathcal{C} \mapsto \mathcal{C}^{\shear}$  on dg categories (whose effect we're describing on underlying triangulated categories) which reflects the shearing operation on rings,
  informally obtained by
 \begin{quote}
 passing to graded objects, shearing all $\Hom$-spaces, and then passing
 back to all objects,
 \end{quote}
  and thus the category of constructible sheaves in \eqref{XexampleX} becomes, under the purity assumption above,  the shear of the derived category of a usual ring. 
        We now turn to formally constructing this 
        shearing operation on categories;       please fasten your seatbelts.

\subsection{Shearing of categories} \label{shearcategory}
 
 In this section we discuss different variants of the shearing operation on the level of categories.
 {\bf All categories in this section will be ``large'' (=presentable) dg categories \S \ref{HigherCatAppendix}.}
 This is essential for the frequent use of de-equivariantization and 1-affineness starting from \S \ref{shearing with Gm} (though in practice the categories we encounter are compactly generated so one can pass back to small categories).

\subsubsection{Shearing graded categories.}
We can use the shearing autoequivalence of $\Rep(\Gm)$ to shear {\em graded categories}, i.e., module categories over the rigid symmetric monoidal tensor category $\Rep(\Gm)$. Model example of such categories are:
\begin{itemize}
\item[(i)]  $\Gm$-equivariant sheaves $QC(X/\Gm)$ on a variety $X$ 
with $\Gm$ action;
\item[(ii)] Graded modules $\mathsf{A}-\textsf{mod}^{gr}$ for a graded ring $\mathsf{A}$. 
\end{itemize}

Such categories
are automatically enriched in graded vector spaces, and the shearing
operation does not change the underlying category but shears the graded $\Hom$ spaces
in the sense of \S \ref{commutativity}.

\begin{definition} 
We define an autoequivalence\footnote{As a model example for the definition that follows (and to normalize signs) that for $X$ a complex of graded vector spaces, we have
$X^{\shear} \otimes Y = (X \otimes Y^{\unshear})^{\shear}$;
that is to say, the action of $Y \in \Rep(\Gm)$ on $X^{\shear}$ corresponds
to the action of $Y^{\unshear}$ on $X$. } $(-)^\shear\actson \Rep(\Gm)\module$ of the category of graded categories by twisting the $\Rep(\Gm)$-action on a given category by the monoidal autoequivalence $(-)^\unshear$ of $\Rep(\Gm)$. Explicitly, 
the new action of $M \in \Rep(\Gm)$ corresponds to the old action of $M^{\unshear}$.

 In other words,
$$D^\shear = D\ot_{\Rep(\Gm)} \Rep(\Gm)^\shear$$ where $\Rep(\Gm)^\shear$ 
denotes $\Rep(\Gm)$ considered as a $\Rep(\Gm)$-bimodule, with left 
action of $Y$ given by $Y^{\shear} \otimes$ and right action given by $ \otimes Y$. 
 \end{definition}

For example, in our model example we have
\begin{equation}\label{shearing graded modules}
\left(\mathsf{A}-\textsf{mod}^{gr}\right)^\shear\simeq {\mathsf A}^\shear-\textsf{mod}^{gr},
\end{equation}
i.e., shearing of graded categories extends the notion of shearing of graded algebras.  

In particular, for a graded category $D$ as above,
there is a tautological equivalence of underlying categories $D \rightarrow D^{\shear}$,
which we will denote by $X \mapsto X^{\shear}$,
which has the property that 
\begin{equation} \label{enriched shear}  \Hom(X^{\shear}, Y^{\shear}) = \Hom(X, Y)^{\shear}\end{equation}
where the $\Hom$ is the enriched hom into graded vector spaces. 
We write out the argument for \eqref{enriched shear} to verify signs:
 The enriched $\Hom(X, Y)$
has graded degree $n$ component given by $\Hom_{\mathcal{C}}(X(n), Y)$, with $(n)$ denoting twisting
by the $\Gm$-representation of weight $n$.
By definition,  $X^{\shear}(n) = ( X(n))^{\shear}[-n]$, 
  and the corresponding graded degree $n$ component
  of $\Hom(X^{\shear}, Y^{\shear})$ is therefore given by
  $\Hom_{\mathcal{C}}(X(n), Y)[n]$.  

\subsubsection{Monoidal structure of shearing.}\label{monoidal shearing}
  Recall that modules over a commutative ring have a functorial symmetric monoidal structure. In our $\infty$-categorical setting~\cite{HA} this provides the category $\Rep(\Gm)\module$ of graded categories with a natural symmetric monoidal structure.
Since we will have occasion to consider the interaction of tensor structures with shearing, we note the following:

\begin{prop} 

\begin{itemize}
\item The autoequivalences $(-)^{2n\shear}\actson \Rep(\Gm)\module$ given by even shears naturally lift to symmetric monoidal autoequivalences of $\Rep(\Gm)\module$. 
\item The autoequivalence $(-)^\shear$ lifts to a symmetric monoidal autoequivalence of 
$\Rep^\super(\Gm)\module$.
\item The autoequivalence $(-)^\shear$ and its odd iterates lift to symmetric monoidal equivalences 
$$\Rep(\Gm)\module\longrightarrow \Rep_\epsilon^\super(\Gm)\module$$ 
(as do its odd iterates); see \S \ref{commutativity} for the notation $\Rep_\epsilon^\super(\Gm)$. \end{itemize}
\end{prop}

\index{$\Rep_\epsilon^\super(\Gm)$}

\subsubsection{Shearing categories with $\Gm$-action.}\label{shearing with Gm}
Next, we want to carry out a corresponding shearing for categories with a $\Gm$ action.

 We briefly recall the notion of categorical representation of an algebraic group $G$ (sometimes called a ``weak action'' of $G$), 
see~\cite{FrenkelGaitsgoryLLC} and~\cite{1affine} for a thorough study. 
Morally, an action of $G$ on a dg category $\cC$ is a family of autoequivalences of $\cC$ labelled by elements $g\in G$, satisying coherent associative composition laws and varying algebraically with $g$. Formally, this is captured by the notion of module category for the category $(QC(G),\ast)$ of quasicoherent sheaves on $G$ equipped with the monoidal structure given by convolution\footnote{Thanks to the self-duality of $QC(G)$, this is equivalent to the notion of {\em co}module category for $QC(G)$ with the convolution coalgebra structure, which is more immediately parallel to the notion of algebraic representation of an affine group scheme.}, or (thanks to the general formalism of descent) to the notion of quasicoherent sheaf of categories on the stack $BG$. We denote the category of $G$-categories by $Cat^G$.
\index{$Cat^G$} The main result in the subject is Gaitsgory's 1-affineness theorem ~\cite{1affine}, which (when applied to the stack $BG$) asserts that the notion of $G$-category (i.e., $(QC(G),\ast)$-module) for $G$ affine is equivalent to that of module category for the {\em symmetric} monoidal category $(Rep(G),\otimes)$. 

Let us focus on $\Gm$-categories, i.e.,
module categories for the convolution monoidal category $(QC(\Gm),\ast)$. 
Model examples of such categories are:
\begin{itemize}
\item[(i)]  Sheaves $QC(X)$ on  variety $X$ 
with $\Gm$ action;
\item[(ii)] (All) modules $\mathsf{A}\module$ for a graded ring $\mathsf{A}$. 
\end{itemize}

The relationship between this and our previous examples:
$$\mathsf{A}\module \leftrightarrow \mathsf{A}\module^{gr}, QC(X)\leftrightarrow QC(X/\Gm).$$
is a general one
relating 
\begin{equation} \label{Gmgraded} \mbox{categories with $\Gm$ actions} \leftrightarrow \mbox{graded categories}.\end{equation}
We detail this relationship in general; the discussion
that follows 
is a categorical version of the passage
from a $\Gm$-space $X$ to $\left( X/\Gm \to B\Gm\right)$ and,
in the other direction the passage from a space $Z$ over $B\Gm$ to 
the $\Gm$-space $Z \times_{B\Gm} pt$.

Given a category $\cC$ with $\Gm$ action we can apply equivariantization, i.e., pass to the category of equivariant objects
$$\cC\mapsto \cC^{\Gm}=\Hom_{QC(\Gm)}(\Vect, \cC)$$
using  the ``augmentation module'', i.e., the pushforward functor $QC(\Gm)\to \Vect$, which upgrades to a symmetric monoidal functor. The result has the natural structure of graded category, expressing the familiar fact that (nonequivariant) Hom spaces between equivariant sheaves carry representations of the group. Thus equivariantization defines a functor 
$$Cat^{\Gm}\longrightarrow \Rep(\Gm)\module.$$
Applying this construction to $\Vect$ itself we find 
$$\Vect^\Gm=\End_{\QC(\Gm)}(\Vect)\simeq \Rep(\Gm),$$
so that in particular we may consider $\Vect$ as a $(\QC(\Gm),\Rep(\Gm))$-bimodule category. 
Equivariantization has a left adjoint construction of de-equivariantization
$$\Rep(\Gm)\module\longrightarrow Cat^{\Gm}, \hskip.3in \D\mapsto \underline{\D}= \D\otimes_{\Rep(\Gm)} \Vect.$$

Gaitsgory's 1-affineness theorem ~\cite{1affine} (applied to the stack $B\Gm$) asserts that the bimodule $\Vect$
produces a Morita equivalence between the monoidal categories $(QC(\Gm),\ast)$ and $(\Rep(\Gm),\otimes)$, i.e., equivariantization and deequivariantization define inverse equivalences between $Cat^\Gm$ and $\Rep(\Gm)\module$. (Indeed this entire discussion holds with $\Gm$ replaced by any affine algebraic group $G$.)

We use this equivalence to transport the shearing operation from graded categories to $\Gm$-categories:

\begin{definition} 
The shear of categories with $\Gm$-action is the autoequivalence
$$(-)^\shear\actson Cat^\Gm,\hskip.3in \cC\mapsto \cC^\shear= (\cC^{\Gm})^\shear \ot_{\Rep(\Gm)} \Vect.$$
\end{definition}

It's important to note that shearing $\Gm$-categories {\em does} change the underlying category. However, it does not change the category of $\Gm$-equivariant objects: it induces an equivalence of categories \begin{equation} \label{dataeq} \shear:\cC^\Gm\to (\cC^{\shear})^\Gm.\end{equation}
 It is useful to think of the equivalence $\shear$ above -- which arises out of the construction -- as part of the data of a sheared category.

\begin{example}[Model example] \label{model example}
If $\mathsf{A}$ is a graded dg-ring, then
shearing carries $\mathsf{A}$-modules to $\mathsf{A}^{\shear}$-modules,
$$
\left(\mathsf{A}-\textsf{mod}\right)^\shear\simeq \mathsf{A}^\shear-\textsf{mod}.
$$  The associated equivalence \eqref{dataeq} is the equivalence~\ref{shearing graded modules} constructed 
 after \eqref{modelgradedequivalence} and is given simply by naive shearing. 
\end{example}

\begin{example}[Trivial action]\label{trivial action}
Any  linear category $\cC$ can be regarded as carrying a trivial $\Gm$-action. 
Shearing of graded objects (defined just as in \S \ref{commutativity})
defines an equivalence 
\begin{equation} \label{trivactionshear}  \unshear: (\cC^{\Gm})^{\shear} \rightarrow \cC^{\Gm},\end{equation}
which we will often denote later by  the ``unshear'' symbol $\unshear$ as it identifies the sheared with the usual category.

For $X$ an object of $\cC^{\Gm}$, this
equivalence sends the object $X^{\shear}$ on the left to the object regrettably\footnote{It could be worse.} also denoted $X^{\shear}$ on the right:
On the left, we use the notation described after \eqref{enriched shear}, i.e.
$X^{\shear}$ is the image of $X$ by the equivalence of underlying categories
which arises whenever we shear a graded category;
on the right, $X^{\shear}$ means the shear of the graded object $X$ defined in a way parallel to
\S \ref{commutativity}, i.e., we cohomologically
regrade $X$ 
using its $\Gm$-action.

In any case,  \eqref{trivactionshear} induces
$$\cC^{\shear} \simeq \cC.$$
 \end{example}

  \begin{remark}[Monoidal structure of shearing]
  As in \S \ref{supersloppy}, the sheared category comes also in a super-version,  
  i.e., as a $\Vect^{\super}$-linear category, which we will allow ourselves to denote by the same notation:
\begin{equation} \label{super shearing}   \cC\mapsto \cC^\shear= (\cC^{\Gm})^\shear \ot_{\Rep_{\epsilon}(\Gm)} \Vect^{\super}.\end{equation}
 It is better to use the super- version for considerations involving symmetric monoidal structure, and in particular
 considerations involving traces. 

 Indeed, just as representations of a group have a tensor product lifting the tensor of underlying vector spaces, $Cat^{\Gm}$ has a symmetric monoidal structure lifting the (Lurie) tensor product of categories, and the equivariantization equivalence $Cat^{\Gm}\simeq \Rep(\Gm)\module$ is naturally symmetric monoidal.
Then  \eqref{super shearing} describes a symmetric monoidal equivalence of 
 $\Vect^{\super}$-linear $\Gm$-categories.

 \end{remark}

\begin{example}[Functoriality of shearing] \label{functoriality of shearing}
A formal feature of shearing that will be useful in our applications is that shearing objects commutes with equivariant morphisms. 
Suppose $\pi_*:\cC\to \D$ is a morphism of $\Gm$-categories $$\pi_*\in Hom_{\Gm-cat}(\cC,\D)\simeq Hom_{\Rep(\Gm)-mod}(\cC^{\Gm},\D^{\Gm}).$$
Then we get a commuting square 

$$\xymatrix{\cC^\Gm\ar[d]^-{\pi_*}\ar[r]^-{\shear} & (\cC^\shear)^\Gm\ar[d]^-{\pi_*^\shear}\\
\D^{\Gm} \ar[r]^-{\shear} & (\D^\shear)^\Gm}$$

A case which will arise later is the situation when $\Gm$ acts trivially on $\D$. In this case
we can augment the square by means of the equivalence  of Example \ref{trivial action}. 
\begin{equation} \label{augsburg} \xymatrix{\cC^\Gm\ar[d]^-{\pi_*}\ar[r]^-{\shear} & (\cC^\shear)^\Gm\ar[d]^-{\pi_*^\shear}\\
\D^{\Gm} \ar[r]^-{\shear} & (\D^\shear)^\Gm  \ar[r]^{\sim}
& (\D)^{\Gm} 
}
\end{equation}
where $\sim$ comes from \eqref{trivactionshear}. 
The functor from $\cC^{\Gm}$ to the far right copy of $\D^{\Gm}$
is, therefore, given by $X \mapsto (\pi_* X)^{\shear}$: usual pushforward,
followed by shearing of a graded object through the $\Gm$ action.

\end{example}

\begin{example}[Automorphic shearing] \label{automorphic shearing} 
We present another point of view on the foregoing discussion
which can be seen as ``Cartier dual,'' and which appears as an automorphic counterpart to shearing. 
Here we replace the role of $\Gm$ by that of its Cartier dual $\Z$. 

We will apply the equivalences of symmetric monoidal categories 
$$( \Rep(\Z),\otimes)\simeq (\QC(\Gm), \ast) \mbox{ and } (\Vect(\Z),\ast)\simeq (\Rep(\Gm),\otimes).$$
 where $\Vect(\Z)$ means simply the category of $\Z$-graded vector spaces; e.g., the second equivalence
arises from the fact that a representation of $\Z$
 over the field $k$ is equivalent to a $k[t, t^{-1}]$-module. 
 It follows from \eqref{Gmgraded} that $(\Rep(\Z),\otimes)$ and $(\Vect(\Z),\ast)$ are Morita equivalent. 
 Modules $\cDD$ over $(\Vect(\Z),\ast)$, i.e., categories with $\Z$-action, are identified with local systems of categories over $S^1=B\Z$. 
 On the other hand, the global sections $\cC=\Hom_{Cat^\Z}(\Vect, \cDD)$ of such a local system of categories is a module for $End_{Cat^\Z}(Vect)=(Rep(\Z)=\Loc(B\Z),\otimes)$. (Equivalently, $\cC$ carries a $\Gm$-action where $\lambda\in \Gm$ acts by 
tensoring by a rank one local system on $S^1$ with monodromy $\lambda$.) The Morita equivalence $\cC\leftrightarrow \cDD$ then establishes the 1-affineness of $B\Z$, i.e., that we can recover a local system of categories $\cDD$ over the circle from its global sections $\cC$ as $\cDD=\cC\otimes_{\Rep(\Z)} \Vect$.
  
As was the case for (de)equivariantization of $\Gm$-actions (see \eqref{Gmgraded})
this Morita equivalence is a categorical shadow of the equivalence
 $$ \mbox{spaces $X$ over $S^1$} \leftrightarrow \mbox{spaces $Y$ with $\Z$-action}$$
 (where a space $X \to S^1$ determines a homotopy fiber  $Y=X\times_{B\Z} pt$ with  $\Z$-action, while
in the reverse direction, we pass from $Y$ to $X$ by taking homotopy quotient by $\Z$).\footnote{Note that these operations are opposite to those defining the Morita equivalence of $\Rep(\Gm)$ and $\QC(\Gm)$, in that the role of $\otimes$ and $\Hom$ or equivariantization and de-equivariantization are exchanged. This is made possible by the fact that $\Vect$ is canonically {\em self-dual}, allowing us to exchange tensors for Homs.}

Finally we ask, how does shearing look from this viewpoint?
If we consider $\Gm$-categories as (global sections of) local systems of categories over $B\Z$,
\begin{quote}
shearing = composing the monodromy automorphism with the shift $\la 1 \ra$.
\end{quote} 
 Explicitly, for a sheaf of categories over $B\Z$ with fiber $\cDD$ with monodromy automorphism $M:\cDD\to \cDD$, we obtain a sheared automorphism $M^\shear=M\circ\la 1\ra$ of $\cDD$. This defines a sheared local system of categories $\cDD^\shear$ over $B\Z$, whose global sections $\cC^\shear$ are the shear of the $(\Rep(\Z),\otimes)\simeq (\QC(\Gm),\ast)$-category $\cC=\cDD^\cZ$.

\end{example}

\begin{example}[Automorphic shearing of categorical representations]\label{shearing categorical reps}  We now explain a general pattern of shearing categorical representations as a {\em categorical analog of twisting representations by characters:}
 Given a group $G$ and a homomorphism $\iota: G \rightarrow \Z$, there is a way to shear a categorical representation $\mathcal{C}$ of $G$,
 which amounts to 
 \begin{quote} $\mathcal{C}^{\iota \shear} =\mathcal{C}$ with the action of  $g \in G$ sheared by $\langle \deg(g) \rangle$\end{quote}
We 
 will apply it to loop groups (and Hecke categories) in defining the normalized action of $G_F$ on $\SHV(X_F/G_O)$ in \S \ref{normalized-local}. 
Because this will come up for us often it will be useful to have a couple of different ways to think about it:

 Schematically, a categorical representation of $G$ is a sheaf of categories over $BG$ in a suitable sheaf theory. The category $SHVCAT(BG)$ of $G$-categories is linear over $\SHV(BG)$. Tensoring by an invertible sheaf of categories defines an autoequivalence of the category of categorical representations. A natural source of such invertible categorical representations is homomorphisms $ \iota: G\to \Z$: we pull back to $BG$ the sheaf of categories on $B\Z$ given by the categorical representation $\Vect\la 1 \ra$ of $\Z$ --
 this is the representation which, 
under the identification $(\Vect(\Z),\ast)\simeq (\Rep(\Gm),\otimes)$.  corresponds to the sheared fiber functor on $\Rep(\Gm)$. This defines the operation
denoted  $\cC\to \cC^{\shear \iota}$ on categorical representations of $G$. 

Another way to describe this shearing operation (along the lines of~\cite{BZG}) is as follows. Let's describe categorial representations of $G$ as modules for a monoidal category $(\SHV(G),\ast)$ of sheaves on $G$ under convolution in our fixed sheaf theory. Then the linearity of categorical representations over $\SHV(BG)$ amounts to a (braided monoidal) central functor $(\SHV(BG),\otimes) \to \cZ(\SHV(G)),$ which lifts the trivial monoidal functor $(\SHV(BG),\otimes)\to (Vect,\otimes) \to (\SHV(G),\ast),$ i.e., the action factors through the augmentation. This expresses the presentation $G\simeq pt\times_{BG} pt$ of $G$ as a groupoid over $BG$. Now given $\iota:G\to\Z$, we obtain a tensor functor $$(\Rep(\Z),\otimes)\simeq (\QC(\Gm),\ast)\to (\SHV(BG),\otimes)\longrightarrow \cZ(\SHV(G)),$$ which again is trivial (factors through the augmentation) as a plain monoidal functor. In other words, $G$-categories acquire functorial $\Gm$-actions, which are trivial as actions on the underlying categories. The operation $\cC\to \cC^{\shear \iota}$ above is given by shearing by this $\Gm$-action.
\end{example}

\subsection{Shearing in geometry} \label{shearing geometry}

We give several examples of shearing of categories of quasicoherent sheaves.

\begin{remark}[Shearing and coaffine stacks]\label{coaffine} 
For any stack $X$ with $\Gm$-action we obtain a $\Gm$-action $\QC(X)\in Cat^{\Gm}$ and consider the sheared category, which we denote  $\QCshear(X)$. 
We begin with a general warning about interpreting this category geometrically.

In general for $X=\Spec(A)$ with $A$ a non-negatively graded (discrete) algebra, we can think of the category $\QCshear(X)=A^\shear\module$ as a {\em variant} of the category of sheaves on the  ``coaffine stack''
$X^\shear=\rm{Spec}(A^\shear)$, the object of derived algebraic geometry represented by the {\em coconnective} commutative dga $A^\shear$.
The category of sheaves on the latter can be defined,  following the general formalism of derived algebraic geometry,  as a limit over maps from affines into $X^\shear$, i.e., in terms of {\em connective} cdgas. 

To illustrate how the two differ, take  $\mathbb{A}^1$ with squaring $\Gm$ action and $A$
the ring of functions on $\mathbb{A}^1$.  We denote the coaffine stack represented by $A^{\shear}$ as $\AA^1[2]$. 
The category of quasicoherent sheaves on this stack, which is equivalently described as $B^2 \Ga$, does not coincide 
with the category $$\QCshear(\AA^1)=A^\shear\module =\cO((\mathbb{A}^1)^\shear)\module.$$
 In terms of the Koszul dual exterior algebra of functions on the derived affine scheme $\mathbb A^1[-1]$, we  have $ \QCshear(\AA^1)\simeq QC^!(\AA^1[-1])$ while, 
 thinking of $\AA^1[2]$ as a coaffine stack, we have instead $QC(\AA^1[2])\simeq QC(\AA^1[-1])\simeq \cO(\AA^1[-1])\module$. 
\end{remark}

 \begin{example}[Shearing of $G$-representations]  \label{outershearexample} 
A cocharacter $\underline{\varpi}:\Gm\to Aut(G)$ gives an action of $\Gm$ on $BG$, hence on $QC(BG)=\Rep(G)$;
it makes  $k[G]$ into a graded Hopf algebra, $\fg$ into a graded Lie algebra and $G$ into a ``graded algebraic group'' in the sense of \S \ref{symplecticshear}.
We obtain a sheared category of representations $\mathsf{D} = \Rep(G)^{\underline{\varpi}\shear}$ which can be described as comodules for the sheared Hopf algebra $k(G)^{\underline{\varpi}\shear}$, \begin{equation} \label{Ccomodule} C \rightarrow C \otimes k[G]^{\underline{\varpi}\shear}.
\end{equation}

In particular, such a comodule inherits an action of the sheared Lie algebra $\mathfrak{g}^{\shear}$, i.e.
a map $$\mathfrak{g}^{\shear} \otimes C \rightarrow C.$$
 If we 
 grade the Lie algebra $\mathfrak{g} = \bigoplus \mathfrak{g}_j$
 via $\underline{\varpi}$, then $X \in \mathfrak{g}_j$ {\em decreases} degree by $j$
 (for $X \in \mathfrak{g}_j$ corresponds to an element of $\mathfrak{g}^{\shear}$ in cohomological degree $-j$). 
  That is to say,
  we can think of $\mathrm{Rep}(G)^{\shear}$ as being
  \begin{quote}
  ``complexes with a non-degree-preserving action of $G$,''
  \end{quote}
  (For a model example of such a complex, see  Example \ref{repGinnershear} below.)

The graded category obtained by equivariantizing $\Rep(G)$ is the category $\Rep(G\rtimes \Gm)$, with its natural $\Rep(\Gm)$ action: tensoring by representations inflated from $\Gm$. 
 \end{example}
 
 \begin{example}[Inner shearing of $G$-representations.] \label{repGinnershear}
 Continue with the prior example but now assume that $\underline{\varpi}$ in fact lifts to $\varpi: \Gm \rightarrow G$. Then we have in fact 
 \begin{equation} \label{eqqqq}\Rep(G) \simeq   \Rep(G)^{\shear \varpi}.\end{equation}
 Indeed  given $C \in \Rep(G)$
 the shear $C^{\shear}$ through $\varpi: \Gm \rightarrow G$ defines an object
 of the category $\mathsf{D}$ above.  Note that in the presentation \eqref{Ccomodule} of
 $\Rep(G)^{\shear}$ as comodules, 
 this category  comes with a fiber functor; the resulting (pulled back) fiber functor
 on $\Rep(G)$ via \eqref{eqqqq} is not  the usual
fiber functor on $\Rep(G)$ under \eqref{eqqqq}, but rather a shear of it. \footnote{
To say a different way:   for inner actions the resulting action of $\Gm$ on $BG$ can be trivialized;
but the   trivialization of the action on $BG$ doesn't preserve basepoints.}
 
 Alternately we can see \eqref{eqqqq} as arising from 
 an equivalence of the associated graded categories, which in turn arises from the group isomorphism
\begin{equation} \label{ab} (g, t) \in G \rtimes \Gm \mapsto (g \varpi(t), t) \in G \times \Gm.\end{equation}
 compatible with projection to $\Gm$.

\end{example}
 
\begin{example} \label{Vannoying}  
We generalize the foregoing example from the case of a point to a general $G$-space:
Suppose that $G$ acts on a space $X$
 and $\lambda: \Gm \rightarrow G$ is a one-parameter subgroup,
 which we regard as acting on $G$ via the adjoint action,
 and on $X$ through $G$. Then 
\begin{equation} \label{eqqq2} QC(X/G)^{\shear} \simeq QC(X/G).\end{equation}
This recovers
  \eqref{eqqqq} when $X$ is a point. 
  Again, it can be deduced from the equivalence of stacks
  $X/(G \rtimes \Gm) \simeq X / (G \times \Gm)$
  where $\Gm$ on the right is acting trivially;
  this equivalence arises from \eqref{ab} on the acting
  groups and the identity on $X$. 
  
We make the equivalence more explicit
  when $X = \mathrm{Spec}(R)$ is affine.

The category $QC(X/G)$ is  the category of $(G, R)$-modules, that is to say (complexes of) $R$-modules
with compatible $G$-action. $\Gm$ acts 
 on the pair $(G, R)$ via
$$x: g \mapsto \lambda(x) g \lambda(x)^{-1},
 r \mapsto \lambda(x) .r,$$ 
 and correspondingly acts on the category of $(G, R)$-modules.
 The category $QC(X/G)^{\shear}$ is now the category of $(G^{\shear}, R^{\shear})$-modules,
 i.e., complexes with an $R^{\shear}$-action and an action of $G$
 that does not preserve degree (just as in \eqref{outershearexample}). 
 Just as discussed after \eqref{eqqqq}, regrading 
 a $(G, R)$-complex by means of $\lambda$
 exhibits the equivalence with $(G^{\shear}, R^{\shear})$-modules.

\end{example}

\subsection{Abelian geometric Satake} \label{shearSatake}

We will now discuss the geometric Satake equivalence, 
both for reference in the rest of the paper, and
because it is a convenient example of shearing.

 \label{sssHecke}
 \index{Satake category}
 \index{affine Grassmannian}
 \index{Hecke category}
 \index{$\Gr_G$}
  Let $G$ be a reductive group over $\FF$. For now, we assume that $\FF$ is algebraically closed. 
  We will discuss the case of $\FF$ a finite field in \S \ref{Satake finite field}. 
 Attached to $G$ is the {\em affine Grassmannian} $\Gr_G$, an ind-variety over $\FF$.
 We call the {\em Satake category} of $G$ the abelian category
$$ \text{Sat}_G := \Perv_\kk(G_O\backslash \Gr_G)$$
\index{$\text{Sat}_G$ Satake category}
of $G_{O}$-equivariant perverse sheaves  on the Grassmannian with coefficients in $\kk$. 
This is a full subcategory of the abelian category $\Perv_\kk(\Gr_G)$ of all perverse sheaves on the Grassmannian, characterized by constructibility with respect to the stratification by $G_O$-orbits.  (The entire derived category of $G_O$-equivariant sheaves on $\Gr_G$ will be called the \emph{Hecke category} in this paper, in order not to have to distinguish between ``perverse'' and ``derived'' Satake categories.)

  This category has a monoidal structure defined by convolution, and a non-obvious
 commutativity structure that can be defined by fusion.
 These definitions can be found in \cite{MV}.
 As usual, there are several choices for sheaf theory. For our current purposes,
we will follow \cite{MV} and use {\'e}tale sheaf theory  
(see \cite[\S 14]{MV}).

 \subsubsection{Abelian Satake} \label{absat}
 
Roughly speaking, the abelian geometric Satake isomorphism identifies $\text{Sat}_G$, viewed as a Tannakian category with the fiber functor to $Vect_\kk$ given by total (non-equivariant) cohomology over $\Gr_G$, with the category of representations of the dual group $\check G$. However, this is slightly imprecise due to a fundamental parity issue. Namely, the cohomology functor from $\text{Sat}_G$ to vector spaces is a monoidal functor, but the 
commutativity constraint on $\text{Sat}_G$ is not compatible with that on vector spaces.
Rather, considering cohomology as a functor to {\em graded} vector spaces, 
the commutativity constraint on $\text{Sat}_G$ is intertwined
with the ``Koszul'' commutativity constraint on graded vector spaces
by \cite[Lemma 6.1]{MV}.

This parity issue is handled in Proposition 6.3 of {\em op.\ cit.}
by modifying the geometric commutativity constraint on the Satake category. One can instead modify the representation category in the following way, parallel to \S \ref{commutativity}: let us
replace $\Rep(\check{G})$ with the symmetric monoidal category $\Rep^{\super}(\Gv)$
of representations on super-vector spaces.
Concretly, objects are pairs $(V_+, V_-)$
of representations of $\check{G}$, and when 
  one swaps the tensor product of two odd vector spaces, one incurs a $-$ sign.
Then the abelian Satake correspondence gives a symmetric monoidal equivalence\index{$\Rep^{\super}(\Gv)$} 
\begin{equation} \label{abSatakeequivalence} \text{Sat}_G^{\super}\simeq \Rep^{\super}(\Gv)\end{equation}
 between the super-version of the Satake category of perverse sheaves on the Grassmannian  
 and super-representations of $\Gv$. Under this equivalence (applying descent, Remark~\ref{super descent}) even sheaves on the left correspond to representations whose parity
is given by the central element $(-1)^{2 \rho} \in \check{G}$. We thus get a symmetric monoidal equivalence for the usual Satake category in the form
\begin{equation} \label{Satake-unsheared} \text{Sat}_G \simeq \Rep^{super}_{2\rho}(\Gv):=\mbox{super-representations of $\check{G}$ whose parity is given by $(-1)^{2 \rho}$}.\end{equation}
\index{$\Rep^{super}_{2\rho}(\Gv)$} (Compare with the appearance of $\Rep^{super}_\epsilon(\Gm)$ in \S \ref{commutativity} and \S \ref{monoidal shearing}.)
 While this might seem overly elaborate, we have already discussed in  \S \ref{Sqrtqsuper} and \S \ref{supersloppy}
 why it is essentially inevitable to consider
 super-vector spaces if one wants to study numerical questions.  
 
In the following sections, use the following notation: \index{$\mathcal{T}_V$}
 For every representation $V$ of $\check G$, define $\mathcal{T}_V$ to be 
\begin{equation} \label{TVcaldef} \mathcal T_V = \mbox{ the correspondent to $V$ under \eqref{abSatakeequivalence}}\end{equation} 
 Hence, $\mathcal T_V$
 is ``analytically normalized,'' that is to say, it is Verdier self-dual,
 where the notion of Verdier duality on the affine Grassmannian
 is normalized to preserve the unit object.\footnote{We make a couple of normalization comments to avoid sign confusion later.   $\mathcal{T}_V$ is the IC sheaf of a certain stratum in $\mathrm{Gr}_G$. If that stratum is $k$-dimensional, then $\mathcal{T}_V$, restricted to the smooth locus, will be 
the constant local system shifted by $\langle k \rangle$, i.e., in cohomological degree $-k$ and with weight $-k/2$.
Moreover, as is relevant when computing trace functions, $\mathcal{T}_V$ is regarded as a super-sheaf
whose parity is given by the action of $e^{2\rho}(-1)$ on $V$.}

 \label{unsheared Satake example}
  \begin{example} 
  \begin{itemize}
  \item[(a)]
  Take $\lambda \in X_*(T)$, a cocharacter of the maximal torus of $G$. Under this equivalence, the IC sheaf of the closure $\mathcal{S}_{\lambda}$ of the $G_O$-orbit represented by $t^\lambda$ (considered as a super-vector space with parity $\left<2\rho,\lambda\right>$), is mapped to the representation with highest weight $\lambda$.

  \item[(b)]
  \label{oddTate} 
  Specializing further the example of (a), consider the case $G=\PGL_2$, and take 
  the first nontrivial stratum $\mathcal{S} \subset \mathrm{Gr}_G$ (the closed orbit on the ``odd'' or non-neutral component of the Grassmannian) is a copy of $\mathbb{P}^1$,
  corresponding to an elementary modification of a rank $2$ vector bundle. 
 
 The constant sheaf $\kk_{\mathcal{S}}\langle  1 \rangle$, placed in degree $-1$,
  defines an object of the Satake category, and under the above equivalence
  $$ \kk_{\mathcal{S}}\langle  1 \rangle \mapsto \mbox{standard representation of $\SL_2$ in even parity},$$
  that is to say, $\mathcal{T}_{\mathrm{std}}$ is simply $\kk_{\mathcal{S}}\langle  1 \rangle$. 
The associated trace function equals $q^{-1/2} 1_{\mathcal{S}}$.  Note there
is {\em no} minus sign despite the fact that the sheaf is in odd degree, for
the notation $\langle \dots \rangle$ includes a shift of super parity, according to our conventions. 
\end{itemize}
  \end{example}

 \subsubsection{Arithmetic shearing}  \label{sss:Satake-shearing} 
 There is another way of formulating the geometric Satake equivalence,
using shearing, rather than modifying the commutativity constraint
 or using super-vector spaces.

To explicate, we introduce a(n even) cohomological grading on $k[\check G]$ that is equal to the grading by the cocoharacter $\varpi=2\rho$
 , i.e., by the following action of $\Gm$:  
\begin{equation} \label{2rhoequation}  (t\cdot f)(g) = \text{Inn}(t^{2\rho})(f)(g) := f\left (\text{Inn}(t^{-2\rho}) (g) \right),\end{equation}
where $\text{Inn}$ is the \emph{left} action of $\check G$ on itself by conjugation.
The resulting graded Hopf algebra, considered as a dg-algebra with trivial differentials, will be denoted by $k[\check G]^\shear$ (as in \S \ref{outershearexample}); we will call this the ``arithmetic shear'' of $k[\check G]$ (or of $\check G$, by abuse of language).
Explicitly, if $t \cdot f = t^{a} f$, then $f$ defines an element of $k[\check{G}]^{\shear}$ in cohomological degree $-a$. 
 
Recall that the sheared category $\Rep(\check{G})^{\shear}$
of \S~\ref{shearcategory} from Example \ref{outershearexample}, can be described as complexes of vector spaces
with a comodule structure for $k[\check G]^{\shear}$. The element $e \in \mathfrak{g}$
has weight $-2$ for the action $\text{Inn}(t^{-2\rho})$ and therefore
{\em raises degree by $2$.}

 Now, inside the dg-category $\Rep(\check{G})^{\shear}$
 we can consider the full abelian subcategory $\Rep_{2\rho}(\check{G})^{\shear} $ 
 which comes from transporting the 
complexes supported in degree zero
under the equivalence $\shear: \Rep(\check{G})\to \Rep(\check{G})^{\shear}$ of  ~\eqref{eqqqq}.
  Equivalently,  $\Rep_{2\rho}(\check{G})^{\shear}$ coincides with the category of
 $k[\check G]^\shear$-comodules on finite-dimensional graded vector spaces, 
 where grading and the action of $2\rho$ coincide.  
 Then, $\Rep_{2\rho}(\check{G})^{\shear}$ is an abelian tensor category, 
 with monoidal structure
 and commutativity inherited from $\Rep(\check{G})^{\shear}$, 
 and there is an equivalence of abelian tensor categories
\begin{equation}\label{Satake-sheared}
 \text{Sat}_G \simeq \text{Rep}_{2\rho}(\check G)^{\shear}
\end{equation}
 compatible with tensor functors to graded vector spaces (cohomology over $Gr_G$ on the left, forgetting the $\check G^\shear$-action on the right).
 \index{$\mathrm{Rep}_{2\rho}(\check G)^{\shear}$}
 
  In other words, the action of $\check{G}$ 
 on  $V=$ the cohomology of a perverse sheaf on $\Gr_G$ is upgraded to a coaction
 $k[V] \to  k[\check G]^\shear \otimes k[V]$
respecting the cohomological grading. We will sometimes denote this coaction, by abuse of notation, as $\check G^\shear \times V \to V$.   
 \index{$k[\check G]^\shear$}
 
 \begin{example}
 
 We continue with Example \ref{unsheared Satake example} part (b). 
 In the same notation, the functor of \eqref{Satake-sheared}
 sends $k_{\mathcal{S}}$ to its total cohomology, i.e.
 the complex $k \oplus k[-2]$ with zero differential, 
 or equivalently the object  of $\Rep_{2\rho}(\SL_2)$ 
 obtained by 
 starting with $\mathrm{std}[-1] \in \Rep(\SL_2)$, the standard representation in cohomological
 degree one, and then shearing through $2 \rho$. \end{example}
 
 \subsection{Derived Geometric Satake} \label{sssderivedSatake}
 
 \index{Hecke category} \index{$\mathcal{H}_G$} \index{$\overline{\mathcal{H}}_G$}
 \index{$\overline{\mathcal{H}}$} \index{derived Hecke category} 

We now turn to the derived setting and discuss the full dg spherical Hecke category 
$$\mathcal{H}_G:=\Shv(G_O\backslash \Gr_G), \overline{\mathcal{H}_G} = \SHV(G_O \backslash \Gr_G) $$ \index{$\mathcal{H}_G$} \index{$\overline{\mathcal{H}_G}$}
of sheaves with $\kk$-coefficients on the affine Grassmannian defined over an algebraically closed field $\FF$.

As discussed in \S \ref{sheaf theory}, this comes in two versions, a small version $\mathcal{H}_G$  
of constructible sheaves, 
and a large version which we will denote by $\overline{\mathcal{H}_G}$. 
We note only here that while the small version is almost certainly what you think it is (if you think about such matters); 
the large version involves a choice about the order in which one takes certain limiting operations, and we use the ind-finite (i.e., renormalized) version, see
\S \ref{renormalization section}.

For $\FF=\kk=\CC$ the different sheaf theories from \S \ref{sheaf theory} -- constructible sheaves, $D$-modules, Betti sheaves or all sheaves -- give rise to equivalent small sheaf categories which we denote simply by $\Shv$ (i.e., $G_O$-equivariant coherent $D$-modules are forced to be regular holonomic, $G_O$-equivariant Betti sheaves are forced to be locally constant on the $G_O$-orbits and the compact objects constructible, and the two are identified by the Riemann-Hilbert correspondence compatibly with the various functors we consider). 
We will also be interested in the case $\FF=\overline{\FF_q}$, with $k = \overline{\Q_l}$ where $l$ is different from the characteristic of $\FF$, in which case $\Shv$ refers to $l$-adic {\'e}tale constructible complexes. 

We now state a mildly strengthened version of the derived geometric Satake theorem of  
Bezrukavnikov and Finkelberg~\cite[Theorem 5]{BezFink}. First, the original theorem is stated in the setting of triangulated rather than dg categories. However the technique of the proof extends to prove the stronger dg statement, and it is stated as such in~\cite[Theorem 11.3.3, Proposition 11.4.2]{ArinkinGaitsgory}. 
Next, the original theorem \cite[Theorem 5]{BezFink} was stated for $\FF=\CC$, $k= \CC$, but it is observed in \cite[Proposition 5]{BezFink} that it extends to the case $\FF=\overline{\FF_q}$, with $k = \overline{\Q_l}$, where $l$ is different from the characteristic of $\FF$.  We will outline below why the arguments of \cite{BezFink} work to establish it also with $k=\Q_l$, together with the action of Frobenius that follows (when $G$ is defined over $\FF_q$).

\begin{theorem}[Derived Geometric Satake] \label{Satake thm}  \label{large Satake}
For the shearing in the following statements, we regard $\fgxv$
as a $\Gm$-space via the squaring action,
and $\Gm$ is acting trivially on $\check{G}$. 

\begin{itemize}
\item Small version: 
 There is an equivalence of monoidal dg categories 
 $$ (\Sph_G, \ast) \simeq \mathrm{Perf}( k[\fgxv]^{\shear}/\Gv),$$
 where the right hand side denotes  dg derived category
 of $\Gv$-equivariant perfect\footnote{At the level of homotopy categories, this corresponds to 
 the smallest triangulated subcategory containing $\mathcal{O}^{\shear}$.} dg-modules over $k[\fgxv]^\shear$ and 
 with monoidal structure given by tensor product. 
 
 \item Large: There is an equivalence of monoidal large categories 
 $$ (\overline{\mathcal{H}_G}, *)\simeq \QC^{\shear}(\fgxv/\Gv)$$
  where the right hand side denotes 
the shearing of the category $\QC(\fgxv/\Gv)$ of
$\Gv$-equivariant coherent sheaves on $\fgxv$.
 \end{itemize}
 
  \end{theorem}

\begin{remark}

\begin{itemize}
\item 
Note that the theorem asserts an equivalence only of {\em monoidal} categories, even though the spectral side has an evident symmetric monoidal structure. 
This accounts for the absence of any ``parity'' correction. The derived substitute for the symmetric monoidal structure on the abelian spherical category is an $E_3$- or factorization structure on $\Sph$, see \S \ref{E3 spherical category}.

\item   (As discussed in \S \ref{supersloppy}): For issues involving symmetric monoidal structure or traces,  it is  preferable
to take the point of view  -- as in \S \ref{absat} -- that this is an equivalence
of categories of super-vector spaces,
wherein even sheaves in $\mathcal{H}_G$ 
are carried to sheaves of parity   determined by 
the central element $(-1)^{2\rho} \in \Gv$.

And just as in \S \ref{absat}, one can eliminate this issue
by use of a $2\rho$-shearing on $\check{G}$ itself; this has
the  advantage of being better suited to questions of rationality over a finite field;
we will describe this ``arithmetic version'', in the terminology
of \S \ref{analyticarithmetic},  in \S \ref{arithmetic dgs}.

\item 
The large version  is stated in \cite[Corollary 11.4.5]{ArinkinGaitsgory}.\footnote{Note that~\cite{ArinkinGaitsgory} also describes an unrenormalized form of the spherical category, which corresponds spectrally to the imposition of nilpotent singular support.}. We can pass to large categories simply by applying the functor ``Ind'' to the two categories in the small version.
On the automorphic side this produces the ind-finite category of sheaves (\S \ref{renormalization section}) on $G_O\backslash \Gr_G$ (again independent of sheaf-theoretic setting).
 \end{itemize}
\end{remark}

 \begin{remark}[Compatibility with Cartier duality of (co)characters]\label{Satake for cocharacters} 
  Consider the data of a character $\eta:G\to \Gm$, dual to that of a central cocharacter $\eta^\vee:\Gm\to Z(\Gv)$. 
  
Working over a local field, one classically knows that twisting an unramified representation 
by $\eta^* \lambda^{\val}$, for $\lambda \in \C^{\times}$, has
the effect of twisting  its Satake parameter through $\eta^{\vee}(\lambda) \in \check{G}$ (up to sign). 
We now describe a corresponding compatibility of the geometric Satake correspondence with Cartier duality for the center of $\Gv$ 
(an ``opposite'' compatibility, for the center of $G$, is mentioned in Remark~\ref{fluxes}).
 This will take the form of corresponding actions of $(\QC(\Gm), \ast)$ on the two sides of the Satake isomorphism. 
 
 On the automorphic side   the action of $(\QC(\Gm),\ast)\simeq (\Loc(B\Z),\otimes)$ comes from the map $$\xymatrix{
G_O\backslash G_F/G_O\ar[r] & BG_F\ar[r]^-{B\eta} & B(\Gm)_F \ar[r] & B\Z}$$
 which provides the spherical Hecke category the structure of sheaf of monoidal categories over the circle.
 This results in a central action of $(\Loc(B\Z),\otimes)\simeq (\QC(\Gm),\ast)$ on $\overline{\Hecke_G}$. Compare
Examples~\ref{automorphic shearing} and
 Example~\ref{shearing categorical reps}.
 
 On the spectral side, 
the cocharacter $\eta^\vee$ defines a braided monoidal functor from $(\QC(\Gm),\ast)$ to the Drinfeld center $\QC(\Gv/\Gv)$ of $\Rep(\Gv)$ (coming from pushforward along the map $\eta^\vee: \Gm\to \Gv/\Gv$), or better to the Drinfeld center $QC(\cLL(\fgxv/\Gv))^\shear$ of the Hecke category. 

 The fact these two actions correspond follows from the results of~\cite{BZG}:
 the automorphic action is part 
 of a larger central action, the Ng\^o action of~\cite{BZG}, by all of $(\Loc(BG_F),\otimes)$. The results of~\cite{BZG} identifying this categorical action with Ng\^o's spectral action of the group scheme of regular centralizers imply in particular the compatibility above.
\end{remark}

\subsection{Geometric Satake over a finite field} \label{Satake finite field}
 
Let us now consider the analogous story  when the group $G$ is defined over a finite field. 

We start by recalling the ``standard'' action of Frobenius on the dual group $\check G$, which we will term ``analytic''
in terms of the general division of \S \ref{analyticarithmetic},
as well as an action that will be termed ``arithmetic'':

\subsubsection*{Analytic Frobenius action on $\check G$}  \label{sss:Satake-shearing-arithmetic} 

When $G$ is defined over $\FF_q$, we have a Frobenius automorphism on the character group of its (universal) Cartan $A$, which gives rise to an action of Frobenius on the dual Cartan $\check A\subset\check G$. By construction ($A$ is defined as the torus quotient of any Borel subgroup), this action preserves the set of positive coroots, which define the standard Borel $\check B\subset \check G$.  The classical definition of the dual
group $\check G$ (over $k$)\footnote{Note that, if $k$ is not algebraically
closed, $\check G$ will always be taken to be split over $k$, as is necessary for the
existence of the pinning. } requires it to be pinned, allowing for a unique extension of the Frobenius action to $\check G$, by pinned
automorphisms.

\begin{definition}
The {\em analytic} action of Frobenius on $\check{G}$  is the pinned action dual to the Frobenius action on the root datum of $G$.
 \end{definition}

Now, we continue to denote by $\Gr_G$ the affine Grassmannian over the algebraic closure $\FF = \overline{\FF_q}$, and by $F, O$ the rings of Laurent and Taylor series over $\FF$, and we will use $\ff$, $\fo$ to denote $\FF_q((t))$, $\FF_q[[t]]$. Then, the category $\text{Sat}_G$ comes with an extra structure, which is the action of Frobenius by pullback of sheaves, $\mathcal F \mapsto \Fr^* \mathcal F$.
We then have the following:

\begin{prop}~\label{sheared Satake w Frobenius}
The equivalences  \eqref{Satake-unsheared} 
 and \eqref{Satake-sheared} of abelian tensor categories
 are compatible with Frobenius actions,
 using respectively analytic and arithmetic Frobenius actions on $\check{G}$:
 
 \index{$\mathrm{Sat}_G$} \index{$k[\Gv]^\shear$} \index{$\mathrm{Rep}_{2\rho}(\check G)^{\shear}$}
 \begin{itemize}
 \item[(a)] In the equivalence \eqref{Satake-unsheared}
we take the Frobenius action on $\check{G}$ to be the analytic action noted above. 
 
 \item[(b)] In the equivalence \eqref{Satake-sheared}
  $\mathrm{Sat}_G \simeq \mathrm{Rep}_{2\rho}(\check G)^{\shear}=k[\Gv]^\shear\comod_{fd}$
  we endow $k[\Gv]$ with the shear of the analytic Frobenius action as in \S \ref{shearing Frobenius notation} and Remark~\ref{shear with Frobenius}.
\end{itemize}
\end{prop}

\begin{remark}
\begin{itemize}
\item[(a)]
Explicitly, the ``arithmetic'' Frobenius action of (b)  
 on $k[\check{G}]^{\shear}$ is the analytic action, multiplied by the   the (left) inner automorphism obtained from the cyclotomic character $\Fr \mapsto q^{-1}$ composed with the cocharacter $(-\rho)$
  into $\check A_{\text{ad}}\subset \check G_{\text{ad}}$.
Thus we  have  $$k[\Gv]^\shear=\bigoplus k[\Gv]_i\la i \ra,$$
where the grading $i$ is via the left inner action of $2\rho$ as in \eqref{2rhoequation}, 
so the Frobenius action  on  the $i$-th (cohomologically) graded piece of $k[\check G]^\shear$ will be the analytic action twisted by $k(-\frac{i}{2})$. Note that, since the grading is even, this does not require choosing a square root of the Tate twist. 

The resulting sheared Lie algebra $\check{\mathfrak g}^\shear$ 
 possesses an element ``$e$'',
 which is in cohomological degree $2$ and is sent to $qe$ by Frobenius.
 Tensoring with any generator of $k(1)$, 
 $e$ therefore defines a Frobenius-invariant element of $\check{\mathfrak{g}}^{\shear}(1)$,

\item[(b)]  In the current context,  one motivation (see also \S \ref{analyticarithmetic}) for  introducing
an arithmetic action is that
the analytic action is {\em not} the action on $\check G$ that arises from the canonical isomorphism 
 of fiber functors \begin{equation}\label{fiberisom}H^\bullet(\Gr_G, \mathcal F) = H^\bullet(\Gr_G, \Fr^* \mathcal F).                                                                                                                                                                                                                                            \end{equation}
 
 To see this, recall that, in  the geometric Satake isomorphism, the pinning comes from cup product with the Chern class $c$ of the determinant line bundle;
that is to say, $c$ acts on the cohomology of perverse sheaves in a way 
that corresponds to the action of a principal nilpotent element  $e\in \check{\mathfrak g}$ in the associated $\check{G}$-representation.
However, this Chern class lives in 
 $H^2(\mathrm{Gr}_G, k(1))^{\Fr}$; thus, 
 without trivializing the Tate twist $k(1)$, this Chern class gives rise to
 an element of  $\check{\mathfrak g}(1)$. 
 This is the element noted in (a) above.

The natural action of Frobenius on $\check G$, that comes from the natural action on $\text{Sat}_G$ and the canonical isomorphism \eqref{fiberisom} via the Tannakian formalism and preserves this ``twisted pinning,'' has been described in \cite{xinwenramified,xinwenSatake}, where it is called ``geometric,''
and corresopnds to our arithmetic action of Frobenius.

\end{itemize}

\end{remark}

\subsubsection{Statement of derived geometric Satake} 
\label{arithmetic dgs}
We now consider  the action of Frobenius in the derived Satake correspondence.
We will leave the somewhat more straightforward ``analytic'' statement to the reader and describe the arithmetic version. 
 In the arithmetic version the shearings are as follows:

\begin{itemize}
\item As before, $\Gv$ is sheared by the left adjoint action of $e^{-2 \rho}$.
\item On $\fgxv$ we use the twist of the previous $\Gm$ action (through squaring) through the 
product of squaring and the left adjoint action of $e^{-2\rho}$. \end{itemize}

Just as in the statement of Proposition \ref{sheared Satake w Frobenius} we consider the shears 
$$ k[\fgxv]^\shear
 \mbox{ and }k[\Gv]^\shear$$ as coming with the ``arithmetic'' Frobenius action, which is to say, the shear of the standard (``analytic'') action on $\fgxv$ and $\Gv$. 
 We then have the following strengthened version of Theorem \ref{Satake thm} (we will
 restate the equivalences just to make shearing clear):

\begin{theorem}\label{Satake thm-arithmetic} 
(Arithmetic normalization of derived Satake, with Frobenius structures). 
\begin{itemize}
\item There is an  equivalence of monoidal dg categories 
$$(\Sph_G, \ast)\simeq (\Perf(k[\fgxv]^\shear)^{\Gv^{\shear}},\otimes)$$
between the spherical Hecke category and $\Gv^\shear$-equivariant perfect dg-modules over $k[\fgxv]^\shear$.

 \item There is an equivalence of monoidal dg categories
 $$(\overline{\mathcal{H}_G},\ast)\simeq (\mathrm{QCoh}(\fgxv/\Gv)^\shear,\otimes)$$ 
 between the large spherical Hecke category and the shear of quasicoherent sheaves on the coadjoint representation.
 
\item For $G$ defined over a finite field, this equivalence is Galois-equivariant, i.e., identifies the Frobenius action on the spherical category with 
the arithmetic Frobenius action on the shears.
\end{itemize}
 \end{theorem}
 
Explicitly, Frobenius acts on $\Gv^\shear$-equivariant dg-modules over $k[\fgxv]^\shear$ by pre-composing the $k[\fgxv]^\shear$-action by the inverse arithmetic Frobenius action on $k[\fgxv]^\shear$, and post-composing the $k[\Gv]^\shear$-coaction by the arithmetic Frobenius action on $k[\Gv]^\shear$.\footnote{Here ``arithmetic'' is referring to the distinction of \S \ref{analyticarithmetic} and not
the distinction between arithmetic and geometric Frobenius.}

The first two assertions are equivalent to those of Theorems  \ref{Satake thm} by the generalities of \S~\ref{shearcategory}. Namely, for every $\Gv^\shear$-equivariant dg-module $N$ over $k[\fgxv]^\shear$, we use the cocharacter $(-2\rho)$ into $\Gv^\shear$, in order to consider $N$ as a \emph{graded} $k[\fgxv]^\shear$-module. The equivalence \eqref{modelgradedequivalence}, then, gives rise to the unshear $N^\unshear$, which is naturally a $\Gv$-equivariant $k[\fgxv]$-module. Explicitly, we recall that if $N = \bigoplus_w N_w$ is the decomposition of $N$ into weight spaces for $(-2\rho)$, then $N^\unshear = \bigoplus_w N_w\la -w\ra$.

As stated, however, the above formulations of the Satake equivalence are not sufficient to pin it down uniquely. The theorem comes from a construction, which has the following additional properties:
\begin{enumerate}
 \item For every dominant coweight $\mu$ of $G$, the IC sheaf of the closed $\mu$-stratum on $\Gr_G$ (with its natural $G_O$-equivariant structure) corresponds to $V_\mu \otimes k[\fgxv]^\shear$, where $V_\mu$ is the irreducible $\check G$-module with highest weight $\mu$ (considered, in the formulation of Theorem \ref{Satake thm-arithmetic}, as a $k[\check G]^\shear$-comodule in degrees determined by the cocharacter $2\rho$.
 
 \item There is a choice of Kostant section $\mathfrak c = \check {\mathfrak c}^* = \check{\mathfrak g}^*\sslash \check G \to \fgxv$, such that the functor of $G_O$-equivariant cohomology $\mathcal F \mapsto H^\bullet_{G_O}(\Gr_G, \mathcal F)$, which is valued in modules\footnote{Of course, we could also consider it as valued in modules for $H^\bullet_{G_O}(\Gr_G,k)$, but we defer to \cite{BezFink} for a description of this structure.} for $H^\bullet(BG_O, k) = k[\mathfrak c]^\shear$, is given by the \emph{Kostant--Whittaker reduction}. (The shear on $k[\mathfrak c]$ will be described below.) 
 \index{Kostant section} \index{$\mathfrak c$}
 
 More precisely, $\fgxv$ contains a regular nilpotent element $e^*$ which is fixed by the $(2-2\rho)$-action of $\Gm$ (i.e., scaling by squares composed with the coadjoint action via $(-2\rho)$). Therefore, the resulting Kostant section, which, under an isomorphism $\fgxv\simeq \fgv$ sending $e^*$ to a nilpotent $e$ of an $\sl_2$-triple $(2\rho, e, f)$, identifies $\mathfrak c$ with $e+ \mathfrak g_f$, is preserved by this $\Gm$-action. Hence, at the level of sheared algebras it defines a homomorphism $k[\fgxv]^\shear \to k[\mathfrak c]^\shear$, where the shearing on $k[\mathfrak c]$ is the one corresponding to this $\Gm$-action.  
 
 In the setting of Theorem \ref{Satake thm-arithmetic}, now, the restriction of a given $k[\fgxv]^\shear$-module to this Kostant section is identified with the functor $H^\bullet_{G_O}$. Moreover, once such an identification of functors is fixed, the equivalence of the theorem is unique, by \cite[Theorems 2 and 5]{BezFink}.
\end{enumerate}

\begin{remark} Regarding the action of Frobenius, we again have canonical isomorphisms 
\begin{equation}\label{fiberisom-derived} H^\bullet_{G_O}(\Gr_G, \mathcal F) \simeq H^\bullet_{G_O}(\Gr_G, \Fr^*\mathcal F),
\end{equation}
 which are compatible with the action of Frobenius on $H^\bullet(BG_O, k)$. The latter is the analog of the arithmetic Frobenius action on $k[\mathfrak c]^\shear$, which simply combines the classical (``analytic'') action on $\mathfrak c$ (arising from its action on the root datum of $G$) with the twist by $(-\frac{i}{2})$ on the $i$-th cohomologically graded piece. This is compatible with the action of Frobenius on the Kostant section described above. 

As in the abelian case, we could also describe an action of Frobenius in the setting of  Theorem \ref{Satake thm} using the classical action on $\check{\mathfrak g}^*$, but we would again need to twist the canonical isomorphism \eqref{fiberisom-derived} by an operator that involves a choice of $q^\frac{1}{2}$, in general.
\end{remark}

{\small
We outline how the statement of Theorem \ref{Satake thm-arithmetic} (with the compatibilities indicated afterwords) follows from the arguments of \cite{BezFink}: The first step in the construction of op.\ cit.\ is the abelian Satake isomorphism \eqref{Satake-sheared}, which holds for arbitrary $k$. (``Arbitrary,'' in the setting of $\FF = \overline{\FF_{p^r}}$, should be interpreted as some ring suitable for $l$-adic cohomology, $l\ne p$.) Theorem 2 of op.\ cit.\ (specialized to $\hbar = 0$) extends this to a full embedding of the category of $\check G^\shear$-equivariant $k[\check{\mathfrak g}^*]^\shear$-modules of the form $V \otimes k[\check{\mathfrak g}^*]^\shear$, $V\in  \Rep_{2\rho}(\check G)^\shear$, into $\mathcal H$; this also immediately extends to arbitrary coefficients. (Note that here we are slightly reformulating by shearing, in order to keep track of cohomological degrees.) When $G$ is defined over a finite field, as \eqref{Satake-sheared} is equivariant with respect to the canonical identification of (non-equivariant) cohomology of a sheaf and its Frobenius pullback, the resulting functor of op.\ cit.\ is Frobenius-equivariant with respect to the identification of equivariant cohomology with the Kostant--Whittaker reduction. Indeed, for a $G_O$-equivariant perverse sheaf $\mathcal F$ on $\Gr_G$, the equivariant cohomology can be recovered from the non-equivariant one, i.e., we have 
\begin{equation}
 H^\bullet_{G_O}(\Gr_G, \mathcal F) = H^\bullet(\Gr_G, \mathcal F) \otimes H^\bullet(BG)
\end{equation}
as Frobenius-equivariant $H^\bullet(BG)$-modules. The extension of these isomorphisms to the entire derived category
 of $\Gv^\shear$-equivariant perfect dg-modules over the dg-ring $k[\fgxv]^\shear$, performed in Section 6 of op.\ cit., applies verbatim to arbitrary coefficient rings, and is unique, hence compatible with the action of Frobenius.

}

\subsection{The sheared coordinate ring of a  hyperspherical varieties} \label{Mshear1}

In our local  conjecture, there will be a particularly important role
played by shearing a hyperspherical variety by its associated neutral $\GGm$ action,
as well as the associated Frobenius structure. 
We explicate this   as a  a convenient reference for later parts of the paper.

Let  $(\check{G}, \check{M})$ be a hyperspherical variety over an algebraically closed field $\kk$ of characteristic zero.
(We use the notation $(\check{G}, \check{M})$ rather than $(G, M)$ simply to be suggestive: the discussion
that follows will be applied on the spectral side.)  We also note 
that the discussion that follows in the case of $\check{G}\times \Gv$ acting on $T^* \check{G}$,
where the $\Gm$ action is by squaring along cotangent fibers, 
is closely related to the discussion of the past sections \S \ref{sss:Satake-shearing} -- \S \ref{Satake finite field}.

\subsubsection{Analytic story}  \label{analytic-M} 
We may form
 the {\em shear of $k[\check{M}]$ by the neutral $\GGm$-action}, explicitly, 
  
  \begin{quote} $k[\check M]^\shear=$\emph{the algebra $k[\check M]$ considered as a super-dg-algebra with trivial differentials, in degrees and super-parity determined by the inverse of the action of $\GGm$. }
  \end{quote}
  \index{$k[\check M]^\shear$}

To avoid any sign confusion,  we emphasize that the action of $\GGm$ on $\mathcal{O}(\check{M})$ is defined
via the rule
 $$ \lambda \cdot f(m) = f(\lambda^{-1} m) \textrm{ or } f(m \lambda^{-1})$$
 depending on whether that action is written on the left or right on $\check{M}$
 (this looks wrong but 
  recall our left/right conventions from \S \ref{leftrightconventions}). 

Now let us moreover suppose that
$\check{M}$ is endowed with a finite order Frobenius action
commuting with $ \Gm$. 
For example, this could be taken trivial -- or, 
in
  \S~\ref{simple-M} we described an action of Frobenius\footnote{or, for that matter, of the absolute Galois group of the field of definition of $M$} for the dual $\check M$ of an {\em untwisted} polarized hyperspherical space, i.e., the dual of a  spherical variety.
In that case, we define the Frobenius action on $k[\check{M}]^{\shear}$
according to our general conventions on shearing  -- that is to say, 
$\langle n \rangle$ incorporates a Tate twist $(n)$. 
 
Therefore, if
 $f(\lambda  m) = \lambda^j f(m)$, 
 then
 $f$ lies in degree $-j$ for the $\Gm$ action; 
then it defines a class $f^{\shear}$ in $\mathcal{O}(\check{M})^{\shear}$
that occurs in cohomological degree $j$;
and if $f \in k[\check{M}]$ is fixed
by the finite order Frobenius action, 
 then the geometric Frobenius acts on $f^{\shear}$ by the scalar $q^{j/2}$.

\begin{example} \label{analytic-M-example}
If $\check M$ is a symplectic vector space, in which case the neutral $\GGm$-action is the scaling by the tautological cocharacter, then:
\begin{itemize}
\item  $\GGm$ acts on $k[\check M] = S^\bullet (\check M^*)$ by negative powers;
\item   $k[\check M]^\shear$ lives in positive cohomological degrees, and
\item The Frobenius action on $k[\check{M}]^{\shear}$ is pure with {\em non-negative} weights i.e.
eigenvalues of absolute value $q^{k/2}$ for various $k \geq 0$. 
\end{itemize}

\index{shifted symplectic geometry}\index{$P_3$ algebra}

 Note, however, that  that the degree $1$
 subspace of $k[\check{M}]^{\shear}$ is regarded as a  super-vector space with odd parity,  and as such $k[\check M]^\shear$  is indeed identified with symmetric powers of its dual. Here, and also more generally,
the symplectic structure on $\check{M}$ manifests itself in a {\em shifted symplectic structure}, in the sense of~\cite{PTVV} -- in particular, 
of a Poisson bracket $k[\check{M}]^{\shear} \times k[\check{M}]^{\shear} \rightarrow k[\check{M}]^{\shear}$ of degree $-2$ (a $P_3$ algebra).

\end{example}

 Moreover, all this structure 
  is compatible with the action of (non-sheared) $\check G$ on $\check M$, i.e.
 $\check{G}$ acts on  $ k[\check M]^\shear$ by transporting its action on the coordinate ring of $\check{M}$. 
Assuming that $\check{G}$ comes with a Frobenius action compatible
with that on $k[\check{M}]$, the same is true for this action;  in other words, we have a Frobenius-equivariant coaction of coalgebras
\begin{equation}
 k[\check M]^\shear \to k[\check M]^\shear \otimes k[\check G].
\end{equation}

 \subsubsection{Arithmetic shearing} \label{Marithshear}
 The above shearing by the neutral action is the ``analytic'' shear of $\check{M}$ in the parlance of \S \ref{analyticarithmetic}. 
The arithmetic version uses a different action of $\Gm$ to grade, which is 
  defined in the situation where $\check{M}$ 
  arises as the dual of a hyperspherical variety $(G, M=T^*(X, \Psi))$. 
  In this situation, after choosing an eigenmeasure on $X$,  we get a cocharacter $\eta: \Gm \rightarrow \check{G}$
  by dualizing  the eigencharacter. In  Definition \ref{arithmetic Gm action}
  we defined the {\em arithmetic action} to be the twist of the neutral action by $\Gm$,
  that is to say, the arithmetic action of $\lambda$ in $\Gm$ on $m \in M$
  is given by the composite of its neutral action and the action of $\lambda^{\eta} \in \check{G}$.
  
We correspondingly define the {\em arithmetic shear} of the coordinate ring:
\[k[\check M]^\shear = \mbox{ the shear of }k[\check M]\mbox{ by the twist of the neutral $\Gm$ action by $\eta+2\rho$.}.\]
  The arithmetic shearing has a very nice parity property. It follows
from Proposition \ref{zXparity}, or from the reinterpretation given directly thereafter: 
 the degrees of
 of  $k[\check M]^\shear$ are all even.
 Just as before, a Frobenius action on $\check{M}$ translates to one on
 $k[\check{M}]^{\shear}$; but here
 we observe that no choice of $\sqrt{q}$ is required because
 of this even-ness.

 Now $\check{G}$ does not act on $k[\check{M}]^{\shear}$
``but $\check{G}^{\shear}$ does,''
where the shearing on $\check{G}$ is through the right adjoint action of the cocharacter $e^{2 \rho}: \GGm \rightarrow \check{G}$
(equivalently, the left adjoint action of $e^{-2\rho}$, i.e., same as in  \S \ref{shearSatake}). 
 Formally, this means that $k[\check{M}]$ 
lies in the sheared category $\Rep(G)^{2 \rho \shear}$, or, equivalently,
it admits a coaction of the associated Hopf algebra:
\begin{equation}
 k[\check M]^\shear \to k[\check M]^\shear \otimes k[\check G]^\shear.
\end{equation}
\index{$k[\check G]^\shear$}

 \begin{example}\label{A1dualarithmeticshear}
 Let us take the case $\check{G}=\Gm, \check{M} =T^* \mathbb{A}^1$,
 which we think of as arising from the dual of $G =\Gm, X=\mathbb{A}^1$. 
 
 The eigenmeasure character is the tautological character of $G$, and corresponds thereby to the
 tautological cocharacter of $\check{G}$. 
 Consequently, if $\mathbb A^1\subset \check M$ is the eigenspace for the tautological character of $\check G$, the  
  arithmetic action of $\GGm$ on $\check{M}$  
 acts by squaring on  $\mathbb{A}^1$ and acts  trivially on the cotangent fiber. 
 
 Thus, writing $x$ for the coordinate on $\mathbb{A}^1$ and $\xi$ for the coordinate on the cotangent fiber,
 the $(\check{G} \times \GGm)$ degrees of $x$ and $\xi$ are given, 
 respectively, by $(-1, -2)$ and $(1, 0)$. For later use in checking signs, 
 the component of $\check{G}$ degree $1$ is spanned by
 $\xi, \xi^2 x^3, \dots$ and lies in $\GGm$ degrees $0, -2, \dots$.

\end{example}

%% file: local-geometric.tex
\newcommand{\GX}{\mathsf{GX}}
\newcommand{\std}{\mathrm{std}}
 
\section{Unramified local  duality} 
\label{section-unramified-local}

\subsubsection{Setup}  \label{uld setup}

Let $\FF$ be either the complex numbers $\C$ or the algebraic closure $\overline{\mathbb{F}_q}$ of a finite field. In this section we will work with a dual pair as in \S \ref{goodhypersphericalpairs}:
$(G, M)$ will be a split hyperspherical pair defined over $\FF$,
for which $M$ admits a distinguished polarization  $M = T^*(X,\Psi)$,
and $(\check{G}, \check{M})$ its dual defined over $\kk$;
this coefficient field $\kk$ will be almost always be taken to be $\C$ or the algebraic closure of an $\ell$-adic field
according to whether $\FF$ is $\C$ or of finite characteristic. 

Occasionally, we will implicitly assume that everything is defined over $\mathbb F_q$, and will introduce Weil structures that the reader can ignore over the algebraic closure. In particular, the appearance of $k\langle 1\rangle:= k[1](\frac{1}{2})$ suggests that we have chosen a square root of the cyclotomic twist, and the twist by $(\frac{1}{2})$ can be ignored by readers interested in statements over $\C$ or $\overline{\mathbb F_q}$.

As in \S \ref{localfieldnotn}, 
let $F = \FF((t))$, with integer ring $O = \FF[[t]]$. As before, when $\FF = \overline{\FF_q}$, we will be using $\ff, \fo$ to denote $\FF_q((t))$, $\FF_q[[t]]$. We will be using $X_O$ for the formal arc space, representing the functor $R\mapsto X(R[[t]])$, and $X_F$ for the formal loop space, representing $R\mapsto X(R((t)))$.
For $X$ affine, these are schemes and ind-schemes respectively.

\subsubsection{} 
In this part we introduce and study the local conjecture, Conjecture~\ref{local conjecture}, a counterpart of the geometric Satake correspondence for spherical varieties. The conjecture asserts an equivalence
$$\boxed{
\SHV(X_F/G_O)\simeq \QCshear(\Mv/\Gv)
}$$
between the categories of $G_O$-equivariant sheaves on $X_F$ and $\Gv$-equivariant quasicoherent sheaves on $\Mv$,
but with the latter category {\em sheared} (i.e., cohomologically regraded) using the $\GGm$-action on $\Mv$. We also state a version for ``small'' categories, matching constructible sheaves and perfect complexes. The equivalence is required to respect various structures on the two sides:
\begin{itemize}
\item[$\bullet$] the basic object (structure sheaf of the arc space $X_O$) is taken to the structure sheaf of $\Mv$ (\S\ref{pointings});
\item[$\bullet$] the Hecke action on $X_F/G_O$ is matched (under derived Satake) with the moment map $\Mv/\Gv\to \fgxv/\Gv$ (\S\ref{Hecke local conj});
\item[$\bullet$] the action of Frobenius on $X_F$ is matched with the action of the grading group $\GGm$ on $\Mv$ (\S \ref{Galois local conj}); and
\item[$\bullet$] the loop rotation structure on $X_F$ is matched with the Poisson structure on $\Mv$ (\S\ref{Poisson local conj},  \S \ref{Poisson from loop}).
\end{itemize}

Some significant aspects of the conjecture are discussed elsewhere. We defer to \S \ref{PlancherelCoulomb} the discussion of the various aspects of the local conjecture which are most concretely understood in terms of the Plancherel algebra (or relative Coulomb branch), a ring object in the spherical Hecke category. 
 In \S \ref{automorphic-factorization} and \S \ref{spectral-factorization} we discuss a crucial additional structure on the local conjecture, namely compatibility with {\em factorization structures} on the two sides (and the action of changes of coordinates) as well as with the global conjecture. 
  
 \begin{itemize}
 \item In \S \ref{spectral local category} we describe the local category on the spectral side.
 \item In \S \ref{arcspaces} we discuss loop spaces and their singularities and describe the local category on the automorphic side (using infinite type sheaf theory as in Appendix~\ref{infinite type}).
 \item In \S \ref{Sarcspace} we study the concrete description of the category of constructible sheaves on $X_F/G_O$ in the ``placid'' case. 

 \item In  \S \ref{normalized-local}  we discuss the geometric analogue of ``unitary normalization'' of the $G$-action on $X$,
 which simply involves including appropriate Tate and cohomological shifts. 
 
 \item  \S \ref{Slocalconjecture} formulates the local conjecture and some of its consequences, and
\item  \S \ref{basic examples}  discusses some examples in which the conjecture is known,  including pointers to the recent literature on the subject.
 
     \end{itemize}

\begin{remark} 
On the side of $X$, it is interesting to relax the requirement that $X$ be smooth, and work within the broader class of affine spherical varieties. In this setting, calculations of ``IC functions'' in \cite{BFGM, BNS, SaWang}  suggest a generalization of the conjecture that follows, by dropping the coisotropic condition on $\check M$ on the spectral side. However, we do not know how to formulate the categorical conjecture at this point -- or whether singular spherical varieties are even the correct objects to consider. (In the singular toric case, for example, it seems that one could replace singular toric varieties by smooth toric stacks.) 
\end{remark}

 \subsection{The spectral local category}\label{spectral local category}
 In this section we highlight some features of the local category $\QCshear(\Mv/\Gv)$ on the spectral side. Recall (\S \ref{uld setup})
 that we are considering a smooth affine variety $\Mv$ equipped with a graded Hamiltonian $\Gv$-action, with $\GGm$-equivariant moment map $\mu:\Mv\to \fgxv$. We will default to using the \emph{neutral} grading on $\check{M}$ (see \S \ref{Sneutral}), and comment later (\S \ref{unnormalized local})  on the effect of changing the grading.   Correspondingly, as is appropriate for the neutral grading, the $\GGm$ action on $\check{G}$ itself is trivial. 
 
 We will make free use of the notions of shearing from \S \ref{shearingsec0};
 see in particular \S \ref{Mshear1}, which we briefly recall: The $\GGm$-graded ring
 $\cO_{\Mv}$ can be sheared to give $\cO_{\Mv}^{\shear}$, 
 a differential graded ring with zero differentials.
 When it is relevant (for example, in taking traces, cf. \S \ref{supersloppy}), we always consider
$\mathcal{O}_{\Mv}$ as a {\em super-} DGA whose $\Z/2$ grading is the reduction
of the integral $\Gm$ grading.
 With this convention, $\mathcal{O}_{\Mv}^{\shear}$ is commutative {\em as a super-DGA};
 it may not be 
 be commutative as a DGA if we naively ignored the super-structure.
 
The local spectral category $\QCshear(\Mv)$ can be described as the associated category of differential graded modules, which
  by the discussion and notations of \S \ref{shearingsec0}, 
  can also be described as the sheared category of modules for $\cO_{\Mv}$ itself.
 Again, there is a shearing functor $E \mapsto E^{\shear}$: a $\GGm$-graded complex $E$ of $\cO_{\Mv}$-modules
 gives a complex 
 $E^{\shear}$ of $\cO_{\Mv}^{\shear}$-modules, by shifting the
cohomological grading by the $\GGm$-grading; and the same is true adding $\check{G}$-equivariance everywhere. 
 
 \begin{remark} \label{Vxlocalconjecture}  
 Recall  (\S \ref{Slodowy}, \S \ref{NTmov}) that
  $\check{M}$ has the structure of vector bundle
  over $\check{G}/\check{G}_X$
  with fiber $V_X$, a certain graded $\check{G}_X$-representation
  that has been discussed in the cited sections.
  It follows that the local category may be described in terms of $(V_X, \check{G}_X)$ rather than $(\check{M}, \check{G})$, 
  \begin{equation} \label{Vxlc2}
\QCshear(\Mv/\Gv)\simeq \QCshear(\check{V}_X/\GvX). 
  \end{equation}
where the grading on $\check{V}_X$ has been defined e.g.\ at the start of \S \ref{NTmov}.
   \end{remark}

\subsubsection{Hecke action}
The equivariant moment map $\mu:\Mv/\Gv\to \fgxv/\Gv$ endows the category of equivariant sheaves on $\Mv$ with a tensor action of sheaves on the coadjoint representation. The $\GGm$-equivariance of the moment map allows us to shear this structure, so we find that $\QCshear(\Mv/\Gv)$ is a module category for $(\QCshear(\fgxv/\Gv),\otimes)$ -- i.e., for the spherical Hecke category as it appears on the spectral side in Theorem~\ref{Satake thm}.

\subsubsection{Affineness} \label{Spectralaffineness}
The affineness of $\Mv$ ensures that $\Mv/\Gv$ is affine over $\fgxv/\Gv$. It follows that $\QC(\Mv/\Gv)$ is identified with the category of modules for the algebra object $\mu_*\cO\in \QC(\fgxv/\Gv)$, and similarly for the sheared version: 
$$\QCshear(\Mv/\Gv)\simeq \mu_*\cO^{\shear}\module_{\QCshear(\fgxv/\Gv)}.$$
Another perspective on this is that $\QCshear(\Mv/\Gv)$ is generated by the structure sheaf as a Hecke-module category. 
More concretely, 
a representation $V$ of $\check{G}$ defines a $\check G \times \GGm$-equivariant
sheaf $\underline{V}$ on $\check M$, i.e., we tensor $V$ (with trivial $\GGm$-action)
with the structure sheaf of $\check M$. The affineness then guarantees that the resulting objects generate $\QCshear(\Mv/\Gv)$. Moreover the spaces of morphisms are explicit and readily computable:
\begin{eqnarray*}
\Hom(\underline{V},\underline{W})&\simeq&  
  \Hom^{\Gv}(\Hom(V,W),\cO_{\Mv}^\shear)\\
&\simeq& \Hom^{\GvX}(\Hom(V, W), \mathcal{O}_{\check{V}_X}^{\shear}),\end{eqnarray*}
cf. \eqref{Vxlc2}.

\subsection{The automorphic local category}  \label{arcspaces}\label{automorphic local}

In this section we will assign to  a (untwisted) polarized hyperspherical variety $M=T^*X$, a category\index{$\mathcal{H}^X$} \index{$\mathcal{H}^X_G$}
$$  \mathcal{H}_G^X := \mbox{``$G_O$-equivariant sheaves on $X_F$''}$$
that we will call the ``$X$-Hecke'' or ``$X$-spherical'' category;
  sometimes we will write this simply as $\mathcal{H}^X$, although formally it depends both on $X$ and $G$
  as well as the specifics of the sheaf theory. 
  As usual it will have a small version $\mathcal{H}^X$ and a large version $\overline{\mathcal{H}^X}$.  
  This category will be equipped with a basic object and a Hecke action (\S \ref{Hecke action HGX}).
 \begin{quote}
   {\bf Caveat: some of the features of these sheaf categories, in particular safety/renormalization, are not adequately covered in the literature and should be taken as being provisional.}\footnote{At present, our
  reference for this material is the paper \cite{raskininfinite}. However, it is not yet published, and it  does not discuss either {\'e}tale/Betti contexts
   or safety and renormalization. We have formulated the Caveat simply to encourage further explication and study of the foundations of the theory.}
\end{quote}

The Caveat above is not ``serious,'' that is to say, we anticipate
that the relevant definitions and formalism can be filled in
with existing technology. 
However, it is adjacent to a ``serious'' issue: 
 we have given definitions that, in general, are, at least {\em a priori} completely impractical to compute with.
 The main issue is that $X_F$ is defined as a union of limits of schemes that are {\em singular}; 
 this problem arises in much work on arc spaces, when $X$ is singular, but it also appears for smooth $X$ when one considers strata in the loop space of $X$.  
 Although, as we will see below, there is a reasonable formal definition of the category of sheaves on $X_F$,
 this issue  prevents us from having a  concrete description of it, 
 as well as access to Verdier duality and the function-sheaf dictionary (see \S \ref{Sheafdesiderata} for some specific expectations which it would be good to prove).  
 We do not know how to resolve these issues in complete generality, 
 and regard this as a fundamental question for further study. 

However, there are two situations where many of these issues go away. Firstly, as described in the subsequent \S \ref{placid case},
$X_F$ is ``placid'' in many situations of interest, which permits
a much more explicit analysis of the sheaf categories. 
And, secondly, as we shall explain in the following section \S \ref{PlancherelCoulomb}, a substantial part of the
general conjecture (namely, everything related to the basic object) can always be reduced to a placid situation: arc spaces of smooth affine varieties.
 
 The contents of the current subsection are as follows:

\begin{itemize}
\item \S \ref{arc loop space} discusses arc and loop spaces and in particular illustrates
by example the singularities of $X_F$. 

\item \S \ref{XFsheaves}
 we outline the formal properties of categories of sheaves on loop spaces $X_F$ for $X$ an affine $G$-variety.
 These sheaf categories are defined
  formally, see Appendix~\ref{infinite type}.
  
  \item   \S \ref{XFGOsheaves}
describes the category of equivariant sheaves of interest, i.e., sheaves on $X_F/G_O$. 
  
  \item Definition \ref{basic object} and  \S \ref{Hecke action HGX} construct, respectively,  the basic object and the  Hecke action on sheaves on $X_F/G_O$. 
  
  These are  consequences of the general functoriality of sheaf categories on infinite type schemes and  
stacks~\cite{raskininfinite}, which we review in Appendix~\ref{infinite type}.
   \item  In \S \ref{subsec generalizations} we will briefly discuss the case of twisted polarizations and the unpolarized case.

\end{itemize}

\subsubsection{Arc and loop spaces} \label{arc loop space} For the moment let $X$ be an arbitrary affine $G$-variety defined
over $\FF$, although our case of primary interest is when $X$ arises from polarizing a hyperspherical $G$-variety.

 The
 arc space
 $X_O$ is the scheme representing the functor that
 sends a test $\FF$-ring $R$ to $X(R[[t]])$,
 see \cite[2.2.1]{KapranovVasserot}. 
  This is represented as an inverse limit of the schemes
$$X_O = \varprojlim X_n,$$
where 
$X_n$ represents the functor $X(R[t]/t^n)$.
For $X$ affine, $X_O$ and the $X_n$ are all affine.
For $X$ smooth, $X_O$ is pro-smooth: the morphisms $X_{n+1} \rightarrow X_n$
are in fact vector bundles of dimension $(\dim X)$.

 The loop space $X_F$ is the ind-scheme representing the functor that sends a test $\FF$-ring $R$ to $X(R((t)))$, see \cite[\S 2.5]{KapranovVasserot}. 
    We can write it as a colimit of schemes $X^l$, $l\in \mathbb Z$, by fixing a $G$-equivariant embedding $X \hookrightarrow V$ into the space of a $G$-representation, 
  and taking $X^l$ to be the points of $X_F$ that lie inside $t^{-l} V[[z]]$
  (where we write $V[[z]]$ for the $R[[z]]$ points of $V$).

  Unfortunately the behavior of the $X^l$ is in general much less nice than
  $X_O$, even if $X$ is itself smooth. 
  In particular, the natural presentation of 
  $X^l$ as a projective limit is not placid:
  if we write  
 $$ X^l = \varprojlim X^l_n$$
 where $X^l_n$ is analogously defined
 i.e., as elements $t^{-l} v$
 where $v \in V[[z]]/t^{n+1} V[[z]]$
satisfies the equations defining $X^l$, rewritten in terms of $v$, modulo terms of order $n+1$. 
Typically,  $X^l$ are not smooth and the transition
 maps are not smooth, no matter how large $n$ is. 
 To illustrate both the complexity and interest of the situation we discuss a simple example. 
 
 \begin{example} \label{quadric}
Consider the case of $X$ a
 sphere inside a vector space $V$, i.e., the level set $Q=1$ of a quadratic form. 
 Taking a basis $e_0, \dots, e_n$
 we take $Q$ to be given by
 $Q(\sum x_i e_i) = x_0^2   +q(x_1, \dots, x_n)$, 
 for a nondegenerate quadratic form $q$ in $n$ variables. 
 Then $X^{-1}_n$ is identified
 with  the locus of $(v_0, \dots, v_n) \in V$ satisfying
\begin{equation} \label{Qprez} Q(v_0 + v_1 t + \dots v_n t^n) = t^2 + O(t^{n+1}).\end{equation}
 For $n=1$ we get simply pairs
 $(v_0, v_1)$ where $Q(v_0) =0 $ and $v_0 \perp v_1$,
 in particular, a smooth variety, 
  but the situation becomes more complicated at the next stage:
  The fiber of $X^{-1}_2 \rightarrow X^{-1}_1$ is $(\dim X)$-dimensional
over each $(v_0, v_1)$ with $v_0 \neq 0$. Over the fiber with $v_0 = 0$,
however, the map is not even surjective at the level of points:
its image is the locus $Q(v_1)=1$, and each fiber above this codimension 
$1$ set is now $(\dim X+1)$-dimensional.  
The same general pattern for transition maps holds at each higher stage
$X^{-1}_{n+1} \rightarrow X^{-1}_n$ -- they fail to be surjective, and
the fibre dimension jumps along the image.

 Nonetheless, there are nice features to the situation.
 As proved in \cite{GK, DrFormalArc} 
 the singularities are ``locally finite dimensional''. 
For example, look at the singularity of $X^{-1}$
 at the constant arc $\gamma(t)=te_0$, 
which projects in \eqref{Qprez} to  $v_0=0, v_1 = \mbox{first coordinate vector}, v_2=v_3=\dots=0$. 
It can be checked that 
\begin{equation} \label{FNQ} \mbox{formal neighbourhood of $\gamma$ in $X^{-1}$ } \simeq \mathbb{A}^{\infty}
 \times \{q=0\},\end{equation}
 that is the local singularity structure is a quadric cone in one fewer dimension.
  \label{finitetypemodels}
  It is not a coincidence that this singularity is the
  same as the singularity of  the variety we get if we set $v_2=v_3=\dots= 0$, i.e., if
  we look just at  $\{v_0, v_1: Q(v_0 + v_1 t) = t^2\}$. This is a simple example
  of a global model (in this case, maps from $\mathbb{A}^1$ to $X$
  of degree $\leq 1$) modelling the singularities of a local situation (in this case,
  maps from the formal disc to $X$).
  It is possible that such techniques could provide concrete methods of accessing
  the sheaf theory on $X_F/G_O$ in the situations of interest to us.   
 \end{example}

 \subsubsection{Sheaves on $X_F$} \label{XFsheaves} 
As elsewhere in the paper, there are different choices for 
sheaf theory. Of relevance
for us in the current discussion are de Rham and {\'e}tale, although
in practice which one we choose makes very little difference in the local setting,
and for this reason we will not formally include superscripts ``dR, {\'e}t'' to 
indicate which one we are working with.\footnote{As mentioned in \S \ref{infinite type} the Betti sheaf theory presents challenges in infinite type, though in practice for stacks such as $X_F/G_O$ the Betti theory is expected to be well defined and produce the same categories as the de Rham theory.}
 
  As explained in  Appendix~\ref{infinite type},
working with either de Rham or {\'e}tale notions of sheaf theory\footnote{As explained in {\it loc. cit.}, to have small sheaf categories and ind-proper functoriality in the de Rham setting we have to restrict to a constructible setting, e.g., to ind-holonomic sheaves, though in our intended applications the discreteness of $G_O$-orbits discussed below makes the distinction between ind-holonomic and all $D$-modules disappear.},  we have -- by formally extending from the case of finite type schemes -- 
 \begin{itemize}
 \item Assignments  
 $$Z \leadsto \SHV(Z)=\SHV^!(Z), \SHV^!_s(Z), \Shv^!(Z)$$ for any scheme, stack or prestack of infinite type. Here we use $!$-sheaves by default and omit the $!$-notation when convenient. Recall that the categories of $*$-sheaves are canonically dual to those of $!$-sheaves whenever they are dualizable.  
\item  This is a contravariant functor under $!$-pullbacks, and covariant under ind-proper morphisms. We denote the pushforward functor as $f_*$ and note it is identified as a left adjoint of $f^!$.

 \item  These functors satisfy base change, in the strong sense that they define a functor out of the correspondence category of prestacks (with one leg ind-proper).
 \item For schemes $X$ and affine groups $G$ the category is compactly generated by finite objects (i.e., equivariant constructible sheaves) $\SHV^!(Z)=\Ind(\Shv^!(Z))$
 where $Z=X/G$. 
\end{itemize}

In particular this formalism can be applied to $Z=X_O$ or $Z=X_F$.
The former setting of $Z=X_O$ is particularly nice.
Since $X$ is smooth, $X_O$ is pro-smooth, in particular {\em placid} (see Appendix \ref{appendix placid setting}); it follows that $!$- and $*$-sheaves on $X_O$ are identified, in such a way that the dualizing sheaf $\omega_{X_O}\in \SHV^!(X_O)$ is sent to the constant sheaf $\kk_{X_O}\in \SHV^*(X_O)$. 

\subsubsection{Sheaves on $X_F/G_O$} \label{XFGOsheaves} 
    
 The loop group ind-scheme $G_F$ acts on $X_F$. We will primarily be interested in the quotient $X_F/G_O$ by the arc-group, and specifically when $X$ is an affine spherical variety.  
  Let us first discuss what kind of object $X_F/G_O$ actually is.
  
 There are [at least] two distinct objects that deserve to go by the name $X_F/G_O$: the quotient {\em prestack} and the quotient {\em stack}. Recall that a prestack is merely a functor on derived commutative rings valued in simplicial sets (or synonymously, at the $\infty$-category level at which we work, topological spaces or higher groupoids). There is a natural notion of quotient prestack of a functor such as $X_F$ by a group functor such as $G_O$, which is described as 
 pointwise taking the simplicial set determined by the $G_O(R)$ action on $X_F(R)$.

 On the other hand, we have the quotient {\em stack}, which is the (fppf) sheafification of this functor. Both objects are attached categories of $!$-sheaves by the general mechanism of right Kan extension. However it is reassuring to point out that categories of sheaves which satisfy descent for a given topology (say fppf) are not affected by replacing a prestack by its fppf sheafification. Thus in fact it follows from the results of~\cite{raskininfinite} - specifically the $h$-descent theorem for $!$-sheaves (~\cite[Proposition 3.8.1]{raskininfinite}) - that the two notions agree.
 
 From our point of view, the fundamental object of interest is the action itself and the resulting category of equivariant sheaves $\SHV(X_F)^{G_O}$ (as discussed in particular in~\cite[Section 3.9]{raskininfinite}). One can verify (extending~\cite[Proposition 3.9.2]{raskininfinite}) that this category agrees with that attached to the prestack or stack quotients $X_F/G_O$. Note that the morphism $X_F\to X_F/G_O$ is not itself finitely presented, but by factoring it as a quotient by a prounipotent group (for which equivariant $D$-modules form a full subcategory) and a reductive group (which is fppf) one can establish descent along this morphism directly.  
 
 In summary,
 sheaves on $X_F/G_O$ are understood to be ($!$-)sheaves on the quotient prestack $X_F/G_O$,
 and this agrees with all  other reasonable ways of defining the same concept:  
 $$ \overline{\mathcal{H}^X} = \SHV(X_F/G_O),
 \mathcal{H}^X =\Shv(X_F/G_O).$$
 
 As we discuss in Remark~\ref{! vs * conjecture}, we expect the theories of $*$- and $!$-sheaves on $X_F/G_O$ to be equivalent, and indeed such an equivalence is implied by Conjecture~\ref{local conjecture}. However we do not know how to see this directly.

 \begin{remark}
 Despite the fact that objects appearing seem (indeed, are) enormous, 
the problem of describing the entire category $\SHV(X_F/G_O)$ is actually rather concrete  and has very little of infinite nature in it. Namely, each of the countable strata of $X_F/G_O$ has the form $B \mathcal{G}$
 for a certain pro-group $\mathcal{G}_i$ which can be replaced (for the purpose of sheaf theory)
 with a finite dimensional reductive quotient $\overline{\mathcal{G}}_i$.  Thus
 the sheaf theory on each stratum is extremely simple; the only question
 is how these are to be glued. 
 \end{remark}
 
 \begin{remark}
   From the point of view of classical harmonic analysis, we might loosely think of
 $\Shv(X_F/G_O)$ as categorifying $G(\mathfrak o)$-invariant Schwartz functions,
 whereas $\SHV(X_F/G_O)$  categorifies distributions or generalized functions.
 See Example \ref{A1example1.5}. 
 \end{remark}

  A crucial feature of this setting is the {\em discreteness of $G_O$-orbits},   (see e.g.~\cite{GaitsgoryNadler}):   \begin{quote} \label{spherical is discrete} 
 {\em Discreteness of $G_O$-orbits:} For an affine {\em spherical} $G$-variety $X$, there are only countably many $G_O$-orbits on $X_F$. 
\end{quote}

  Thanks to this, we expect  all $D$-modules on $X_F/G_O$ to be ind-holonomic and all Betti sheaves to be ind-constructible, 
  and we expect   the de Rham and constructible sheaf theories to give equivalent categories (all three settings $\Shv, \SHV, \SHV_s$).
  Moreover, although we did not define Betti categories in 
 Appendix~\ref{infinite type}, it is reasonable to use the constructible categories
 here as the {\em definition} of the Betti category. 
Finally, we have access to the  full functoriality of constructible sheaves, which makes the situation considerably more tame than for example sheaf theory on $X_F$ itself.

 \subsubsection{Basic object and Hecke actions} \label{Hecke action HGX}
  We now describe some basic properties of the category $\mathcal{H}^X$,
  namely the basic object and the Hecke action. 
  
   \index{basic object}
 \begin{definition} \label{basic object}
 The {\em basic object} $\delta_X\in \Shv(X_F/G_O)$ is the pushforward to the loop space of the dualizing sheaf\footnote{
 Equivalently, with reference to the equivalence of $!$- and $*$-sheaves on the placid scheme $X_O$, the {\em constant} $*$-sheaf} on the arc space, i.e., for $i:X_O/G_O\hookrightarrow X_F/G_O$ we have $$\delta_X=i_*\underline{\omega_{X_O}}.$$  
 \end{definition}

Let us recall some general formalism about Hecke actions. 
The Hecke stack  $$Hecke:=BG_O \times_{BG_F} BG_O= G_O\backslash G_F/G_O$$ can be considered as an {\em ind-proper groupoid} over $BG_O$ (see~\cite[Sections II.2.5.1, III.3.6.3, V.3.4]{GR}).  Recall that a groupoid object over $S$ is a simplicial object $\cG_\bullet$ satisfying a Segal condition resulting in an identification of the simplices with iterated fiber products:
 \begin{equation*}
 \xymatrix{\cdots \ar[r]<1ex> \ar[r]<.5ex> \ar[r] \ar[r]<-.5ex> \ar[r]<-1ex> &
 \cG\times_S \cG\times_S \cG \ar[r]<.75ex> \ar[r]<.25ex> \ar[r]<-.25ex> \ar[r]<-.75ex> &
 \cG\times_S \cG \ar[r]<.5ex> \ar[r] \ar[r]<-.5ex> &
 \cG \ar[r]<.25ex> \ar[r]<-.25ex>&
 S}
 \end{equation*} 
 The Hecke groupoid is an ind-proper groupoid, meaning that all of the structure maps involved are ind-proper. For any prestack $Z$ with $G_F$-action, the quotient $Z/G_O\to BG_O$ carries a tautological action of the Hecke groupoid (the descent data describing $Z/G_O\to Z/G_F$). 
Our three flavors of sheaf theory ($\Shv,\SHV_s,$ and $\SHV$) evaluate on $Hecke$ to produce three flavors of Hecke categores. 
 
 The $G$-action on $X$ induces an action of the group ind-scheme $G_F$ on $X_F$, and thus the quotient stack $X_F/G_O\to BG_O$ carries an action of the Hecke groupoid. Again the functoriality of $\SHV$ automatically endows $\SHV(X_F/G_O)$ with the structure of $\overline{\Hecke_G}$-module category.

\begin{remark} 
[Self-adjointness]
It is useful to note that the action of the Hecke category on any module category, in particular on $\SHV(X_F/G_O)$, is self-adjoint, in the sense of isomorphisms
 $$ \Hom(T_V \star \mathcal{F}, \mathcal{G}) \simeq \Hom(\mathcal{F}, T_{V*}\star \mathcal{G})$$
 for $V \in \Rep(G)$, this isomorphism being natural in $\mathcal{F}$ and $\mathcal{G}$. This is a formal consequence of the dualizability of the objects $T_V$ (with duals $T_{V*}$) and more generally is part of the powerful duality package available for modules over {\em rigid} tensor categories such as $\Hecke_G$ (see e.g.~\cite[Section 1.9]{GR}).
\end{remark}

\begin{remark}[Signs]
The Hecke action corresponds, under sheaf function correspondence
when applicable, to the action arising from the action of $G_{\mathfrak{f}}$
on functions on $X_{\mathfrak{f}}$ wherein $g \in G_{\mathfrak{f}}$
sends $f$ to the function $x \mapsto f(xg)$.
\end{remark}

 \subsubsection{Generalizations} \label{subsec generalizations}
 
 It is our expectation that one can make  satisfactory definitions analogous to $\mathcal{H}^X$ in the twisted case 
 or in the unpolarized case using existing technology. We briefly discuss
 these in turn.

 \begin{remark} \label{twisted case}
(The case of twisted polarizations $M=T^*(X, \Psi)$).  
In the case of the Whittaker model, Gaitsgory \cite[\S 2]{Gaitsgory-Whittaker} has defined a local derived category of Whittaker sheaves on the affine Grassmannian (and, more generally, at quotients of arbitrary level of the loop space of $G$). The same definitions make sense in the general case of Whittaker induction. Namely given a $\Ga$-bundle $\Psi \rightarrow X$
one begins by constructing the category $\SHV(\Psi_F/G_O)$ of equivariant sheaves on $\Psi$; and
then, following \cite{Gaitsgory-Whittaker}, takes ``twisted coinvariants'' for the $\Ga_F$-action on this category. \end{remark}

\begin{remark} \label{unpolarized local}
 The non-polarizable case: it is evidently important to  extend the construction $X\leadsto \SHV(X_F/G_O)$ to non-polarizable 
 anomaly-free hyperspherical varieties $M$. 
Note that it is our expectation that the {\em normalized}  Hecke action, 
a twisting of \S \ref{Hecke action HGX} to be described in  \S \ref{normalized-local}, 
will extend to this situation, rather than the action of \S \ref{Hecke action HGX} itself.

In the case when $M$ is a symplectic vector space,
this amounts to studying 
 ``the unramified part of the geometric Weil representation.'' The possibility of such a construction is well-known to experts, although we do not know of an entirely satisfactory reference for our purposes. We describe two relevant works:
 
 For $G=$ the dual pair $\Sp_{2n}\times \SO_{2n}$ (with $M=$ the tensor product of their standard representations), the category analogous to $\Shv(X_F/G_O)$ was constructed by Lafforgue and Lysenko  \cite{LafforgueLysenko} (with an emphasis on the perverse objects), by geometrizing the Schr{\"o}dinger model and the canonical intertwining operators. On the other hand, Raskin \cite[\S 10]{Raskin-homological} has developed a general framework for describing $G_F$-actions on categories as (geometric analogs of) Harish-Chandra modules; this was applied in \cite[\S 2.3, 4]{RaskinNisyros} to modules for the Weyl algebra $\mathcal W$ associated to the symplectic space $M_F$, in order to define an action of $\mathcal H_G$ on $\mathcal W\mbox{-mod}^{G_O}$; here, $G$ can be any subgroup of $\Sp(M)$.
 
 \end{remark}

\subsection{Sheaves on loop spaces: the placid case} \label{Sarcspace} \label{placid case}

\subsubsection{Setup} \label{placid loop setup}
In many cases,
the geometry of $X_F/G_O$ is particularly nice,
and it admits the following presentation, enabling   
 one to avoid the subtle issues of singularities discussed in \S \ref{arcspaces}. There is:
 \begin{itemize}
 \item  an exhaustion $X_F = \varinjlim X^l$ of the
ind-scheme $X$, and a $G_O$-stable presentation of each $X^l$ as $\varprojlim_{n} X^l_n$, where
\item the $X^l_n$ are of finite type, and the transition maps $X^l_{n+1} \rightarrow X^l_n$ are
torsors for a unipotent group scheme, and 
\item   the $G_O$-action on $X^l_n$ 
factors through an action of some quotient $G_{N}$ on each $X^l_n$. 
\end{itemize}

Let us collectively refer to these properties as {\em the placid case}, although {\em a priori} it is slightly stronger
than $X_F$ being placid in the sense of the Appendix.

 For example, in the case of $X=$ a vector space, the $X^l_n$ were constructed in 
\S \ref{arc loop space}.\footnote{One should note that the situation
is quite rich even in the case when $X$ is a vector space,
because  
 $G_O$-orbit closures can be singular even though $X_F$ is very simple, cf. Example \ref{GJexample}. }When $X$ is a 
reductive group $H$ (and $G$ is a subgroup of $H \times H$ acting by left and right multiplication)
we can take the $X^l$ to come from strata of the affine Grassmannian,
and the $X^l_n$ to arise from the cover arising from the kernel of $G_O \rightarrow G_{n+1}$. 
Although these conditions are {\em a priori} quite restrictive,  a surprising number of  hyperspherical cases satisfy them:

 \begin{itemize}
  \item  The  Godement-Jacquet and Iwasawa--Tate periods, $X=M_n$ as a $G = \GL_n \times \GL_n$-space;
\item More generally, any vectorial case, e.g.\ the Hecke period $X=\mathbb{A}^2$ as a $G=\SL_2$-space; 
 \item  The group period, $X=H$ as an $H \times H$-space (possibly twisted by an automorphism in one factor).
 \item The Rankin--Selberg period, $X=\GL_n$ as a $\GL_n \times \GL_{n-1}$ space;
 \item the case of $X=\GL_n \times \mathbb{A}^n$ as a $\GL_n \times \GL_n$-space;
 \item the Gross--Prasad period, $X=\mathrm{SO}_n$ as a $\SO_n \times \SO_{n-1}$-space.
 \end{itemize}

Indeed, it seems possible that the placidity  of $X_F$ for all smooth $X$ can be proved by the techniques of  Drinfeld \cite{DrinfeldInfinite}, 
cf. \cite[\S 8.1]{ChenNadler}. The details of this have not appeared in print at the time of typing this sentence.

 The restriction to this placid setting allows us to describe the entire sheaf theory just in terms of 
 the notion of  sheaves on Artin stacks.
 The discussion of this section can be carried out in any of our sheaf-theoretic contexts, though we will be primarily concerned with the {\'e}tale sheaf theory of
$X_F/G_O$, since this version can be readily compared with numerical
statements in the finite characteristic case. Accordingly, we restrict the discussion to that setting to simplify notations. Thus the notation $\Shv$ denotes, by default, ``{\'e}tale sheaves,''
with constructibility or boundedness conditions to be specified.    However, we will freely use at later
points that the same discussion transposes to the de Rham theory.

\begin{quote}
{\bf Caveat: We expect
the proofs of the statements that follow in this section -- i.e.
various properties of the sheaf theory developed in \S 
\ref{automorphic local}, but restricted to the placid case -- 
 to be straighforward. However, we don't state them as theorems,   because of 
the fact, already noted, that sheaf theory in these contexts
is not clearly documented, and we prefer
to leave the credit for such statements 
to a paper which also develops the relevant foundations in detail,
across the various types of sheaf theory that we would like to access.}
\end{quote}

\subsubsection{Duality and integration}\label{Sheafdesiderata} 

The following properties should be straightforward to establish 
in the placid case. In fact, we expect that they hold in the general context
of spherical varieties,  but we don't know how to   establishing them
because of the issues of singularities mentioned.

\begin{itemize}
 
\item (Duality): The sheaf theory on $X_F/G_O$ is equipped with a Verdier duality satisfying $D \delta_X =\delta_X$
that is moreover compatible with that on the spherical category, i.e.
$D(T \star \mathcal{F}) \simeq (DT) \star (D \mathcal{F})$
(where the duality on the affine Grassmannian is normalized to 
preserve the identity of convolution).

\item (Integration): In the case where $\mathbb{F}$ is the algebraic closure of a finite field,
if $\mathcal{F}, D\mathcal{G}$ are Weil sheaves with trace functions $f,g: X_{\mathfrak{f}}/G_{\mathfrak{o}} \rightarrow \kk$ we have
\begin{equation} \label{int formula desid} [ \Hom(\mathcal{F}, \mathcal{G})^{\vee} ]= \int_{X_{\mathfrak{f}}/G_{\mathfrak{o}}} f(x) g(x) |\omega|,\end{equation} 
for motivation see Lemma \ref{Homlemma}; 
we follow the convention that  $\omega$ is a $G$-eigenmeasure (\S \ref{eigencharacter}), and normalized according to  
\begin{equation} \label{measXOGO}
\mu(X(\fo)/G(\fo)) = q^{-\dim(X)+\dim(G)} \frac{|X(\mathbb{F}_q)|}{|G(\mathbb{F}_q)|},
\end{equation} 
\end{itemize}

 We will sketch below an explicit model for objects in the sheaf category
 and explain, with reference to this explicit model, why duality and integration statements hold.

\subsubsection{Explication of $\mathcal{H}^X_G$.}
In the placid case, 
both objects and morphisms 
of $\mathcal{H}^X_G$ can be described in terms of finite-dimensional computations,
as we now recall. 

 With notation as in \S \ref{placid loop setup}
 write  $\bar{X}^l_n$ for the Artin stack
$$\bar{X}^l_n = X^l_n /G_{N}$$
 where $N$ is chosen so large that the $G_O$-action
factors through it -- the constructions that follow
will formally speaking depend on this choice, but
will be independent of it up to equivalence. 
Write  $\Shv^l_n$ for the category of constructible $l$-adic sheaves on this Artin stack.
While the notion of sheaves on an Artin stack can be of course
derived from the general formalism described in the Appendix,
it is also accessible in more direct ways as sheaves for the lisse-{\'e}tale topology,
see \cite{LO}. 
  
 Placidity
 implies, as in Appendix \ref{appendix placid setting}, that
 the category $\Shv(X_F/G_O)$ constructed previously is equivalent
 to the corresponding category defined via $*$-pullbacks. 
 Therefore, in what follows, we will freely work with the $*$-version.

  The formal construction of the category of $*$-sheaves $\Shv(X_F/G_O)$ implies the existence of functors
\begin{equation} \label{placidinclusion} \Shv^l_n \longrightarrow \Shv(X_F/G_O).\end{equation}
Explicitly, these functors are given by applying the functoriality of $*$-sheaves (these are dual to the already noted
functoriality for $!$-sheaves, namely, $*$-pullback over arbitrary morphisms and $!$-pushforward over ind-proper morphisms)
  to the diagram
$$\xymatrix{X^l_n/G_N & \ar[l] X^l_n/G_O & \ar[l] X^l/G_O \ar[r] & X_F/G_O}$$
 This permits us to think of objects of $\Shv^l_n$ as objects
 in $\Shv(X_F/G_O)$. 
  In fact,
 the morphisms can also be computed in $\Shv^l_n$, as we will explain.
 There are functors 
\begin{equation} \label{ShiftLN} \Shv^l_{n} \rightarrow  \Shv^l_{n+1}, \,\,\, 
 \Shv^l_n \hookrightarrow \Shv^{l+1}_n.\end{equation} 
 which arise from pullback, and from pushforward-pullback, along the
  diagrams of (space, group acting) 
$$ (X^l_{n+1},  G_{N'}) \rightarrow (X^{l}_n, G_{N}), (X^l_n, G_{N}) \leftarrow (X^{l}_n, G_{N''}) \to (X^{l+1}_n, G_{N''}),$$
where $N', N''$ are chosen so that the relevant $G_O$-action factors through these quotients.
Because the kernel of maps $G_{n+1} \rightarrow G_n$
are abelian unipotent group schemes,
the functors of \eqref{ShiftLN} are in fact {\em fully faithful},
so long as  $N, N', N''$ are chosen sufficiently large. 
The space of maps, computed in $\Shv(X_F/G_O)$,
between (the image of)
 a sheaf $\mathcal{F}$ on $\bar{X}^l_n$ and (the image of) a sheaf
 $\mathcal{G}$ on $\bar{X}^{l'}_{n'}$ is obtained 
 by pulling back both to some common $\bar{X}^K_N$ and then considering
 $$   \varinjlim_K \varinjlim_N \Hom_{\bar{X}^K_N}(\mathcal{F}, \mathcal{G}).$$
 Here ``$\varinjlim$''  refers strictly to the homotopy colimit of representing complexes;
 however in the inner limit all the maps are eventually quasi-isomorphisms
 and similarly for the outer colimit; therefore, this can be computed simply
 by a termwise direct limit of representing complexes.

 In terms of this discussion, the
 desiderata above 
 correspond to the following constructions:
 
 \begin{itemize}
\item
The Verdier duality functor on $\mathcal{H}^X_G$
is compatible with the Verdier duality on each $\Shv^l_n$
(where $n$ is sufficiently large relative to $n$), up to shift. 
  The shifts
 are determined up to an overall dimension constant fixed by requiring $$ D \delta_X = \delta_X.$$

\item 
 Also the integration formula \eqref{int formula desid}
follows, when $\mathcal{F}, \mathcal{G}$ arise from $\Shv^l_n$,  by applying  Lemma \ref{Homlemma}
applied to suitable $\bar{X}^K_N$. 
\end{itemize}

 \subsubsection{Examples}

 \begin{example} \label{A1example1}
  Take $X= \mathbb{A}^1$ and $G =\mathbb{G}_m$ (the ``Iwasawa--Tate case''). 
 For a test ring $R$, we have
 $X^l_n(R) = t^{-l} R[t]/t^{n}$
and the $G_O=R[[t]]^{\times}$-action on it factors through
 $(R[t]/t^{n+l})^{\times}$,
 with the transition maps in $n$ being the obvious projections.
 
  Lett $\delta_l$ denote the constant sheaf (with its trivial $G_O$-equivariant structure) on $X^l$, twisted by $\left<l\right>$, i.e., 
 \begin{equation} \label{deltadef} \delta_l = \underline{k}_{X^l}\left<l\right>,\end{equation}
 By this, we mean the $*$-sheaf arising from
 the constant sheaf on any one of the truncations $X^l_n$, by means of \eqref{placidinclusion}. 
  The twist by $\left< l\right>$ is, of course, the standard weight-zero perverse normalization of the sheaves, if we declare the zeroth stratum to have ``dimension $0$.''

 We then compute:
\begin{equation} \label{Deltadeltahom} \Hom({\delta_{l_1}}, {\delta_{l_2}}) \simeq  H^*(B\mathbb{G}_m)  \langle - d \rangle,\end{equation}
   where $d:=|\ell_1-\ell_2|$ is the absolute value of the difference. 
   Indeed, without the normalizing twist by $\langle l \rangle$, the $\Hom$-space
   would be $H^*(B\mathbb{G}_m)$ for $\ell_1 > \ell_2$ and $H^*(B\mathbb{G}_m) \langle -2d \rangle$ otherwise.

 \end{example}

\begin{example} \label{A1example1.5} ($ \overline{\mathcal{H}^X}$ as distributions).

Continue with the notation of Example \ref{A1example1}. 
There are morphisms
 $$ \underline{k}_{X^l} \rightarrow \underline{k}_{X^{l-1}}, \,\, \underline k_{X_l} \rightarrow \underline k_{X^{l+1}}\langle 2 \rangle,$$
the second arising from the inclusion of dualizing sheaves, which are trivial here (up to shifts).
There are corresponding maps of $G_O$-equivariant sheaves.
Let $\delta_0$ be the object of $\overline{\mathcal{H}^X}$ corresponding to 
taking the colimit along the first system, and $\mu$ the object that corresponds to colimit along the second system. 
Then $\delta_0$ and $\mu$ (with their evident Frobenius structures) are analogues of the Dirac delta function at $0$, and the Lebesgue measure, respectively.  

 Indeed, for a sheaf $\mathcal{F}$ on $X^l/G_O$
 considered in $\overline{\mathcal{H}^X}$ by means of \eqref{placidinclusion}, with corresponding function $f$ on $X^l(\mathfrak o)$, we 
 get by an application of \eqref{int formula desid}
 $$[ \Hom(\mathcal{F}, \delta_0)^{\vee} ]= f(0) \mu(X(\fo)/G(\fo)) \mbox{ and } [ \Hom(\mathcal{F}, \mu)^{\vee} ] = \int_{X(\ff)/G(\fo)} f(x) dx,$$
 where the latter measure is normalized according to \eqref{measXOGO}. 
(Of course, the price of this enlargement is that the large category $\SHV$ also contains many other ``large'' objects
that are much less familiar from harmonic analysis.)
\end{example}

  \begin{example}\label{GJexample} 
   We consider the case of   
 $$X=M_n=\mbox{$n \times n$ matrices}, G= \GL_n \times \GL_n,$$
 with action map $x\cdot (g_1, g_2) = g_1^{-1} x g_2$.
  This is, in the theory of automorphic forms, the ``Godement--Jacquet'' period.
  
Fix $0 \leq j \leq n$ and let $\mathcal T_j$ be the Hecke sheaf for $\GL_n$ 
  that corresponds to the representation $\wedge^j \mathrm{std}^*$ on the dual group, i.e.
  it is the intersection complex of the stratum of the affine Grassmannian for $\GL_n$
  that indexes subspaces of $t^{-1} O^n/O^n$ of dimension $j$. Of course, ``intersection complex,'' here, just means the constant sheaf twisted by $\langle j (n-j) \rangle $. 
  
  We think of it as a sheaf on the second copy of $\GL_n$ in $G$.  Then we readily compute  \begin{equation} \label{Ij1}  \mathcal T_j \star \delta_X = \mathrm{I}_j   \langle j(n-j) \rangle\end{equation} 
where $\mathrm{I}_j$ is a $\GL_n \times \GL_n$-equivariant sheaf on $M_{n,F}$, supported on $M_{n,O}$ and pulled back from $M_{n,O}/t M_{n,O} = M_n$,  
whose fiber at a matrix $S \in M_n$ is given by the cohomology of
the space of $j$-dimensional subspaces $E \subset \kk^n$
such that $S|_E = \{0\}$. 
 
  \end{example}

 \subsection{Normalized action of $G_F$} \label{normalized-local} 
 
 In order that our main conjecture 
 match appropriately with the  \emph{neutral} grading on $\check{M}$ (see \S \ref{Sneutral}), we
 need to introduce a cohomological shift to the action of $G_F$ on $\Shv(X_F)$ which, at the level of functions, corresponds to a normalized, unitary action on $L^2$-spaces. We will introduce this normalized version first, before explaining how to ``twist away'' the normalization by modifying the $\GGm$-action.

 Recall from the structure theorem (in the notation of \S \ref{dp}) that $X$ is an equivariant vector bundle over a homogeneous space of the form $HU\backslash G$:
  $$X = S^+\times^{HU} G.$$
  We will assume, as in \S \ref{eigencharacter}, 
  that $X$ admits a nowhere vanishing eigen-volume form $\omega$ with eigencharacter $$\eta: G \to \Gm.$$
  When $G, X$ are defined over a finite field $\FF_q$, taking complex absolute values induces a positive eigenmeasure $|\omega|$ on the points of $X$ over the local field $\ff = \FF_q((t))$. This gives rise to a Hilbert space completion $L^2(X(\ff), |\omega|)$ of the Schwartz space $\mathcal S(X(\ff))$, which furnishes a unitary representation of $G(\ff)$ under the \emph{normalized} action 
\begin{equation}\label{action-G-unitary}
 g\cdot f(x) = \sqrt{|\eta(g)|} f(xg).
\end{equation}

    We lift this normalized action to the level of loop spaces, as follows: the character $\eta$ gives a ``degree'' or ``valuation'' map 
   \begin{equation} \label{deg definition} \deg=\deg_\eta: G_F \to \Gm_F \xrightarrow{\val} \mathbb Z\end{equation}
   To be clear here about the convention: $t\in \Gm_F$ has valuation $1$.
At a categorical level, the normalized action can be described as the automorphic shearing $\cC\mapsto \cC^{\deg \shear}$ of categorical $G_F$-representations $\cC$ defined from the homomorphism $\deg$ as in Example~\ref{shearing categorical reps}.

   More explicitly,  this normalized action of $G_F$ on $\Shv(X_F)$ is defined by twisting the translation action of  $g\in G_F$ by $\left<\deg(g)\right>$, i.e., 
 \begin{equation}\label{eq:normalized-local}
  g\cdot \mathcal F:= g^*\mathcal F \left< \deg(g)\right>,
 \end{equation}
where $g^*$ denotes pullback by the (right) translation action,
and $\langle \dots \rangle$ is as in our standing conventions (cf. e.g.\ \S \ref{Sqrtqsuper}). 
 This induces a similar twisted action of the Hecke category $\mathcal H_G$ on $\mathcal H^X_G$. \footnote{
 The twisting of Hecke actions can be directly  described as a shearing operation as in Example~\ref{shearing categorical reps} and Remark~\ref{Satake for cocharacters}:
 Namely $\deg$ defines a map $G_O\backslash G_F/G_O\to B\Z$, hence endows a central action $$(\Loc(B\Z),\otimes)\simeq (\QC(\Gm),\ast)\longrightarrow \cZ(\Hecke_G)$$ of $\Gm$ on $\Hecke_G$, for which the underlying monoidal functor $\QC(\G_m)\to \Hecke_G$ is trivial. It follows that we obtain a shearing operation on $\Hecke_G$-modules, which doesn't change the underlying category. Applying this operation to the $\Hecke_G$-category $\Hecke^X$ recovers the normalized action.  }

 This twisted action is compatible with the numerical version.  If
 we are working over $\mathbb{F}_q$, and $\mathcal{F}$ is a Weil sheaf on $X_F/G_O$
 with associated trace function $f$, then
 $$ \mbox{trace function for $\mathcal T_V \star \mathcal{F}$} = T_V f$$
 where $\mathcal{T}_V$ is as in \eqref{TVcaldef}, 
 $T_V$ is the corresponding Hecke operator obtained by sheaf-function correspondence; on the right
 $T_V f$ is obtained by integrating the \emph{unitary} action \eqref{action-G-unitary} of $G(\ff)$ on $\mathcal S(X(\ff))$. 
  
 \begin{example} \label{A1example2} 

 In the Iwasawa--Tate case,  let $V_n$ be the representation $z\mapsto z^n$ of $\check G$;
 then
 $$ \mathcal T_{V_n} \star \delta_X = \delta_n$$
the latter object being described in Example \ref{A1example1}. 
 
  The associated
  function (under sheaf-function correspondence) is
  $q^{-n/2} 1_{t^{-n} \mathfrak o}$; it has $L^2$ norm independent of $n$. 
  Note in order to have this be valid without a sign $(-1)^n$
  the parity shifts embedded in $\langle n \rangle$ are important: $\delta_n$
is a  super-sheaf in parity $(-1)^n$.
  
By contrast, under the unnormalized action we would find instead
 that $\mathcal T_{V_n} \star \delta_X$ is  the sheaf $\underline k_{t^{-n} O}$.

\end{example}

\begin{example} \label{GJexample1.5}
In the Godement-Jacquet case of Example \ref{GJexample}
 the character $\eta$ is given by $\eta(g_1, g_2)= (\det g_1)^{-n} (\det g_2)^n$. 
Then we have
$\deg \mathcal{T}_j=- jn$ and
 for the normalized action
\begin{equation} \label{Inormdef} 
\mathcal T_j \star \delta_X=  I_j \langle  -j^2 \rangle
\end{equation}
compare with \eqref{Ij1}.

\end{example}

\begin{remark}\label{remark-choice-eta}
The foregoing definitions depend  on the choice of eigenmeasure $\omega$.
As discussed more generally in
\S \ref{ssseigencharinocuous}, the validity of the conjectures that we will formulate does not depend on this choice of extension. 
Indeed, any two eigenmeasures $\omega, \omega'$ will differ by $\nu: X \to \Gm$,
and this defines an equivalence of categories between sheaves on $X_F$ with the $\eta$-normalized $G_F$-action and sheaves on $X_F$ with the $\eta'$-normalized $G_F$-action, given by applying a  twist by $\left< \val\, \nu \right>$. 
\end{remark}

 \subsection{The local unramified conjecture} \label{Msheaves}
 \label{Slocalconjecture}

In this section we introduce the local conjecture. We are in the setting of \S \ref{uld setup}, with $M=T^*X$ polarized and the dual $\Mv$ is endowed with the neutral $\GGm$-action. 
 
       \begin{conjecture} \label{local conjecture}
There is an equivalence of categories\footnote{As per our general conventions, ``small categories'' are small $\kk$-linear idempotent complete stable categories with exact functors, and ``large categories'' are $\kk$-linear presentable stable categories with colimit preserving functors.}:
\begin{itemize}
\item (small version):
$$\LL_X: \Shv(X_F/G_O)\longrightarrow  \mbox{perfect $\Gv$-equivariant modules for $\mathcal{O}_{\Mv}^{\shear}$}$$
\item (large version):
$$\LL_X: \SHV(X_F/G_O) \longrightarrow\QCshear(\Mv/\Gv).$$
\end{itemize}
The equivalence is required to be compatible with {\em pointings} (\S\ref{pointings}), {\em Hecke actions} (\S\ref{Hecke local conj}), {\em Galois actions} (\S \ref{Galois local conj}) and {\em Poisson structures} (\S\ref{Poisson local conj},  \S \ref{Poisson from loop}). 
\end{conjecture}

We now enumerate and briefly discuss a sequence of structures we require to match on the two sides and some of the immediate consequences.
Note that the left hand side of the conjecture depends on a choice of de Rham or {\'e}tale sheaf theory, but
the conjecture  says that this depends on that choice only through the coefficient field $\kk$. 

\subsubsection{Pointings.}\label{pointings}

The first requirement we make of the conjecture is a matching of basic objects on both sides, i.e., a ``pointing'' of the equivalence:

\begin{itemize}
\item[$\bullet$][Pointing.] There is an identification $\LL_X(\delta_X)\simeq \cO_{\Mv/\Gv}^\shear$ between the image of the basic object and the sheared structure sheaf of $\Mv/\Gv$. 
\end{itemize}

 For instance if we compute $\Hom(\delta_X, \delta_X)$ inside the category of sheaves on $X_O/G_O$, we should have
\begin{equation} \label{invquot}  
\mbox{sheared algebra of functions on $\Mv\sslash \Gv$} \stackrel{?}{=}  \End(\delta_X)   = H^*_{G}(X),
\end{equation}
 the $G$-equivariant cohomology of $X$; the grading on the left hand side is by the $\GGm$-action on $\Mv$.
This can be readily verified by hand in many examples. \footnote{
 Note that this includes a claim on the \emph{formality} of the endomorphism dg-algebra, i.e., the fact that it is quasi-isomorphic to its cohomology. In the homogeneous case, $X=H\backslash G$, this is the same as the cochain complex of $BH$, which is well-known to be formal, and canonically isomorphic to the ring of $\check H$-invariant polynomial functions on $\check{\mathfrak h}^*$, the dual Lie algebra of the dual group. The general case reduces to the homogeneous case, using the structure of a vector bundle $X\to H\backslash G$, i.e., $\End_{\text{Shv}(X/G)} \underline{k}_X = \End_{\text{Shv}(BH)} \underline{k}_{BH} = H^*(BH) = H^*_G(X)$.}
See also Example~\ref{multiplicity-freeness} where an elaboration of this calculation, combined with the compatibility of the local conjecture with Poisson structure, forces the multiplicity-freeness of the Hamiltonian space $\Mv$.

\subsubsection{Hecke actions}\label{Hecke local conj}

Perhaps the most important requirement of the local conjecture is that it matches actions of the spherical Hecke category on both sides.
We first state a ``normalized'' form:
          
\begin{itemize}
\item[$\bullet$][Hecke actions]          
$\LL_X$ has the structure of equivalence of module categories for the spherical Hecke category $\Hecke_G$, where $\Hecke_G$ acts on the automorphic side by the normalized action as in \S \ref{normalized-local}\footnote{In particular, we are assuming a choice of square root of $q$ in $k$, whenever half-Tate twists appear. We will later present an unnormalized variant of our conjecture.} and on the spectral side via the
derived Satake isomorphism $\LL_G$, and the moment map
 $$ \check{M}/\check{G} \longrightarrow \fgv/\check{G}. $$
 \end{itemize}
 
 The Hecke-equivariance implies a very explicit knowledge of the automorphic category. Given a representation
$V$ of $\check{G}$ we have an equivariant vector bundle $\underline{V}\in \QCshear(\Mv/\Gv)$, and as noted in \S \ref{Spectralaffineness}, these objects generate and have explicitly computable Hom spaces. The matching objects on the automorphic sides are the 
sheaves $\mathcal T_V\star \delta_X$ 
on $X_F/G_O$ obtained by applying the corresponding Hecke operator to $\delta_X$ (under the unitarily normalized action). 
Under the sheaf-function correspondence, $\mathcal T_V \star \delta_X$
matches with the function $f_V$ obtained by applying the corresponding Hecke operator $T_{V}$
to $1_{X(\mathfrak o)}$. Thus the conjecture implies isomorphisms of complexes 
\begin{equation} \label{basic} \mathrm{Hom}(\mathcal T_V \star \delta_X, \mathcal T_W \star \delta_X) \simeq  \Hom^{\check G}(V, W \otimes \cO_{\Mv}^\shear),\end{equation} where the degree grading on the right hand side corresponds to the weight or cohomological grading on the left. 
Moreover these explicit Hom spaces (and their compositions) capture the entire automorphic category. In particular {\em all} sheaves on $X_F/G_O$ can be generated from the basic object through the action of Hecke functors. 

Note that the arc space $X_O$ is placid, and hence sheaf theory there is well behaved (for example satisfies Verdier duality). Thus the local conjecture guarantees that all sheaf theory on $X_F/G_O$ can be reduced using the Hecke action to inherently finite-dimensional calculations on $X_O$, and the potential subtleties of sheaf theory in the non-placid setting disappear.
 This is most succinctly expressed using the inner endomorphisms of the basic object -- {\em Plancherel algebra}, see \S\ref{PlancherelCoulomb}.
 
 \begin{remark}[$!$- vs $*$-sheaves]\label{! vs * conjecture} 
  Conjecture~\ref{local conjecture} provides an equivalence between the categories $\SHV^*(X_F/G_O)$ and $\SHV^!(X_F/G_O)$ of $*$- and $!$-sheaves on $X_F/G_O$, since the two sheaf theories are canonically dual and the category $\QCshear(\Mv/\Gv)$ is canonically self-dual. In fact this equivalence is naturally Hecke-equivariant, since the self-duality of $\QCshear(\Mv/\Gv)$ is linear over $\fgxv/\Gv$. It would be very interesting to have a direct proof of this equivalence, i.e., a construction of a $G_O$-equivariant renormalized dualizing sheaf on the (potentially) non-placid ind-scheme $X_F$. As we will discuss in \S \ref{automorphic-factorization} the $*$-sheaf theory is much more natural for the construction of the geometric form of $\Theta$-series and local-global compatibility for the ($!$-)period sheaf.   \end{remark}

\subsubsection{Changing $\GGm$-actions and the unnormalized conjecture.}\label{unnormalized local}
The local conjecture posits that the local automorphic category with the normalized Hecke action $\SHV(X_F/G_O)^{\deg \shear}$ --- sheared by the degree character of $G_F$ as in~\ref{shearing categorical reps} --- matches the local spectral category $\QC(\Mv/\Gv)^{\shear}$  where we shear using the {\em neutral} $\GGm$-action on $\Mv$.

This conjecture is ``analytic'' or ``normalized'' in the parlance of \S \ref{analyticarithmetic}. 
We can pass to the ``arithmetic'' or ``unnormalized'' version using the compatibility under geometric Satake between characters of $G$ and central cocharacters of $\Gv$ (see Remark~\ref{Satake for cocharacters}). 
In particular, in this way,  we deduce forms of the conjecture where the action of $G_F$  on $\Shv(X_F)$ is given by usual translations (i.e., the unnormalized case).

For explicitness in what follows, let us write the neutral action explicitly as $\textrm{neut}$, so that the neutral shear
will be denoted $\QC(\Mv/\Gv)^{\textrm{neut}\shear}$.  
Shearing by $\eta$ to undo the degree shift
in the conjecture, we find an identification of
 $ \SHV(X_F/G_O)$
 with $\QC(\Mv/\Gv)^{(\textrm{neut}+\check\eta)\shear}$,
 where on the right
 where we have modified the $\GGm$-action by the central cocharacter $\check \eta$, i.e., we define a non-neutral action by letting $z\in \GGm$ act instead by the neutral action of $(\check \eta(z), z)\in \check G\times \GGm$.
\footnote{Here is a sign check. In the Iwasawa--Tate example of Example \ref{A1example1},  let $V$ be the standard representation, 
then $T_V$ is the Hecke operator corresponding to a uniformizer, and the unnormalized
$T_V \star \delta_0 =\delta_1 \langle -1 \rangle$ with our previous notation;
thus $\Hom(T_V \star \delta_X, \delta_X) =  \Hom(\delta_1 \langle -1 \rangle, \delta_0)
\simeq H^*(B\Gm)$. On the other hand, we can compute
$\Hom(V, \mathcal{O}(\Mv)^{\shear})$ with the arithmetic shearing as in Example
\ref{A1dualarithmeticshear}: it is again in $\GGm$ degrees $0,2,\dots$.}

 Now, on the right hand side, we can further shear by the action of $2\rho$ on {\em both $\Mv$ and $\Gv$}
 by \eqref{eqqq2} to get
$$ \SHV(X_F/G_O)  \simeq \QC(\Mv/\Gv)^{(\textrm{neut}+\check\eta)+2\rho \shear}.$$
 Note that, above, the $\GGm$ action on $\Mv$ is through $\textrm{neut}+\check{\eta}+2\rho$ --
 that is to say, the arithmetic shear described in \S \ref{Marithshear} -- 
and the $\GGm$ action on $\Gv$ is through $2\rho$. So we can abridge the above statement to:
 $$ \SHV(X_F/G_O) \simeq \QC(\Mv/\Gv)^{\textrm{arith.}  \shear},$$
which  is the unnormalized form of the local conjecture, written in a form to be compatible with the arithmetic
form of geometric Satake.
Again, we emphasize that on the right $\GGm$ is acting on {\em both} $\Mv$ and $\Gv$. 
We remark that the shearing operations above do not change the underlying categories -- just the Hecke actions on them.

  \subsubsection{Galois-linearity}\label{Galois local conj}
   If $\FF$ is the algebraic closure of a finite field
of size $q$, then we require that $\LL_X$ is Galois-linear:

\begin{itemize}
\item[$\bullet$][Galois linearity]          
$\LL_X$ has the structure of Frobenius-equivariance, i.e., it intertwines the autoequivalence of $\mathcal H^X$ induced by Frobenius with  the ``analytic'' action of Frobenius (\S \ref{analytic-M}) on the spectral side.
 \end{itemize}

Note that, since we are supposing that $(G, M)$ is the ``distinguished split form'' hypothesized in \S \ref{goodhypersphericalpairs}, 
the expectation is that this analytic action on $\mathcal{O}_{\Mv}^{\shear}$ is simply the action which scales by $q^{i/2}$
in cohomological degree $i$. In the general case -- that is to say,
even if $(G, M)$ is not the distinguished split form -- this analytic action
 is a a finite order twist of this action, 
as constructed in \S \ref{CheckMRat}. 
We proceed with our discussion in the split case, leaving the modifications to the reader.

  The Galois-linearity  
  implies that $\LL_X$ takes Weil sheaves to $q^{1/2} \in \GGm$-equivariant sheaves, and for Weil sheaves $\cF,\cF'$ the Frobenius action
on the left hand side of
$$ \Hom(\mathcal{F},\mathcal{F}') \simeq \Hom(\LL_X(\cF), \LL_X(\cF'))$$
corresponds to the action of $q^{-1/2} \in \GGm$ on the right hand side.   
\footnote{
As a sign check consider the situation of \S \ref{pointings} when $\cF=\cF'=\delta_X$, 
  the basic
object with its trivial Frobenius-equivariant structure; the 
assertion is then that
the (geometric) Frobenius action on $H^*_G(X)$
coincides with the action of $1/\sqrt{q}$ on $\Gv$-invariant functions on $\Mv$.
The inverse seems strange but is in line with our conventions on shearing:
degree $k$ functions on $\Mv$ are sheared via a cohomological shift $[k]$
and a Tate twist of $(k/2)$. }

In fact,
the spectral side of the conjecture admits a natural graded lift $\QC(\Mv/\Gv\times \GGm)$, in which we impose in addition equivariance for $\GGm$. On the automorphic side, it is natural to expect this to match the category of constructible graded sheaves on $X_F/G_O$, as has been defined in the finite-dimensional setting in ~\cite{HoLiMixed}. The latter is roughly speaking obtained from the category of Weil sheaves by formally imposing that Frobenius acts by powers of $q$. This enhanced version of the Galois-linearity of the conjecture would generalize the ``Tate'' form of geometric Satake (as in~\cite{xinwenSatake}), and imply that the geometry of spherical varieties enjoys some of the very special mixed geometry (or ``Tate-ness'') of flag varieties.

\subsubsection{Parity and spin structures.}
Following the general discussion of \S \ref{supersloppy}
the conjecture has a version where, on both sides, 
  we take sheaves of super-vector spaces.  The parity constraint, where
we use the parity element $z_X$ from \eqref{zXdef}, should be the following:

\begin{itemize}
\item[$\bullet$][Parity]: 
  even sheaves on the left side are carried to sheaves of parity   $z_X$ (i.e.,
super-sheaves whose parity coincides with the action of $z_X$). 
\end{itemize}

That is to say: we consider  $\check{G}$-equivariant sheaf on $\check{M}$
whose super-grading coincides with the grading defined by the central involution $z_X \in \check{G}$.

Here are some examples:

\begin{itemize}
 \item The basic sheaf $\delta_X$. It is itself even,
 and it is sent to $\mathcal{O}_{\Mv}$; the parity 
 condition of \S \ref{parity} implies that $(z_X, -1)$ acts trivially on it;
 that says precisely that it has parity coinciding with the action of $z_X$.

 \item  The parity of $T_V \star \delta_X$
 is determined by the action of $z_X$ on $V$, e.g.
 if $z_X$ acts as $1$ on $V$ then $T_V \star \delta_X$
 is even. 
 
  To see this, recall that
 the parity of $T_V$ itself is, by our conventions, described
 by the action of $(-1)^{2 \rho}$ on $V$, and the convolution
 introduces a twist $\langle \deg T_V \rangle$ which alters the parity
 through the action of $(-1)^{\eta}$ on $V$.

 On the spectral side, the action of $z_X$
  on $V \otimes \mathcal{O}_{\Mv}^{\shear}$ 
  is the diagonal action. If $z_X$ acts by $1$ on $V$,
  its diagonal action here then coincides with the action of $-1 \in \GGm$, as desired.

  \end{itemize}

\begin{remark}[Dependence on spin structures.]\label{local conjecture depends on spin}
We have formulated the local conjecture in the presence of a local coordinate, i.e., an identification $O=\FF[[t]]$. It is natural to ask for a version that is independent of the choice of coordinate, i.e., equivariant for actions of the group of automorphisms of the disc $D=Spec(O)$. Perhaps surprisingly, we do {\em not} expect the most straightforward coordinate-independent formulation of the local conjecture to hold in general. This is due to a need to twist by spin structures on the curve which we will encounter globally in the definitions of both period- (\S \ref{bunX} and \ref{spindep}) and L-sheaves (\S \ref{LsheafXnormalized} and \ref{Lspinindep}). We defer to \S \ref{local-global} for a precise formulation of this spin-twisted and coordinate-independent local conjecture.  
\end{remark}

\subsubsection{Arthur Parameters.} \label{spectral local Arthur}\
 A striking feature of the $L^2$ theory of
spherical varieties is that they are expected (see \cite{SV} and also \cite{ClozelBurger})
 to contain
only unitary representations of a {\em single} Arthur type,  that is,
they are all associated to a single conjugacy class of maps $\SL_2 \rightarrow \check{G}$. 
This $\SL_2$ is the one from which the dual $\check{M}$ is Whittaker-induced, 
see \eqref{ArthurSL2}. 

Let us examine some related phenomena
im the local conjecture. 
Denote by $(\varpi,f)$ the corresponding $\mathfrak  sl_2$-pair (\S \ref{sl2pair}).  
\begin{itemize}
\item[(a)]
The  moment map for $\Mv$ misses $\cNN_{<f}$, the union of nilpotent orbits smaller than $f$. 
Therefore, 
the local conjecture
implies that 
\begin{quote} $\SHV(X_F/G_O)$ is {\em $f$-antitempered as a $\Hecke_G$-module category},
\end{quote} 
by which we mean that 
 it is annihilated by the sheared structure sheaf $\cO^\shear_{\cNN_{<f}}\in \QCshear(\fgxv/\Gv)$.
  
\item[(b)] The second  point is a more direct parallel to the fact, in the classical setting,
that the archimedean size of Hecke eigenvalues is essentially controlled by $\varpi$.   The $\Rep(\Gv)$-action on $\SHV(X_F/G_O)$, by means of the Hecke action,
 factors through the  sheared forgetful functor $(-)^{\varpi\shear}: \Rep(\Gv)\to \Rep(Z_\Gv(\varpi))$,
 obtained by restricting along the inclusion $Z_{\Gv}(\varphi) \hookrightarrow \Gv$
 and then shearing by $\varpi$. 
\end{itemize}

Point (b) comes  from the construction of the $\GGm$-action on the Whittaker induction.   It is for example
familiar from the usual Satake isomorphism when $X$ is a point. In that case, the assertion amounts to the fact
that the cohomology functor on $\mathrm{Gr}_G$ corresponds 
to the forgetful functor, taking a representation of $\Gv$ to its underlying vector space,
but sheared through the action of $2 \rho$.

\subsubsection{Singular support} 

We can also formulate a ``safe'' or ``nilpotent singular support'' version of the local conjecture. This version predicts an equivalence of categories
$$\SHV_s(X_F/G_O) \longrightarrow\QCshear_\cNN(\Mv/\Gv)$$ where we consider ind-safe sheaves on the automorphic side and sheaves supported over the nilpotent cone in $\fgxv/\Gv$ on the spectral side. This equivalence is required to be linear for the ``safe / nilpotent singular support'' form of the derived geometric Satake correspondence described in~\cite{ArinkinGaitsgory}.
 Indeed this version follows from the Hecke compatibility of the local conjecture as in Appendix \S \ref{safe Langlands}. Namely, we restrict both sides of the conjecture to the full subcategories which are torsion with respect to the action of the algebra $Z$ of endomorphisms of the unit in the Hecke category.

\subsubsection{Poisson structure.}\label{Poisson local conj}
We will describe two closely related structures on the automorphic side which are expected to reflect the induced Poisson bracket on functions on $\Mv$: loop rotation and factorization. We formulate the simpler statement about loop rotation here (namely that the deformation given by loop-rotation equivariance recovers the Poisson structure on $\Mv$), spell it out more concretely in \S \ref{Poisson from loop}, 
and defer to \S \ref{automorphic-factorization} and \S \ref{spectral-factorization} for a discussion of factorization.

\begin{remark}  
The local conjecture implies that the  the entire automorphic category $\SHV(X_F/G_O)$ carries a variant of a braided tensor structure:
  the category generated by the basic object carries a {\em locally constant} factorization structure, which in the Betti setting is equivalent to a braided monoidal ($E_2$) structure, or {\em balanced} braided structure ({\em framed} $E_2$-algebra) when combined with loop rotation. 
\end{remark}

The loop spaces $X_F$ and $G_F$ carry compatible actions of $\Gm$, 
 which arise from the action of $\mathbb{G}_m$ on $\FF((t))$ via rescaling of the parameter $t$. These actions
 fix the arc spaces $X_O$ and $G_O$, and induce compatible actions of $\Gm$ on the monoidal category $\Hecke_G$ and its module category $\SHV(X_F/G_O)$.
 
 Thus we may consider the $\Gm$-equivariant category $$\SHV_u(X_F/G_O):=\SHV(X_F/G_O\rtimes \Gm)=\SHV(X_F/G_O)^{\Gm},$$ i.e., the strongly $\Gm$-equivariant objects in the $X$-spherical category. This is a $H^*_{\Gm}(\mbox{pt}, k) \simeq \kk[u]$-linear category, where $u$ has cohomological degree two, which specializes at $u=0$ to a full subcategory (the {\em equivariantizable objects}) of the original category,
 $$\SHV_u(X_F/G_O)\otimes_{\kk[u]}k\hookrightarrow \SHV(X_F/G_O).$$ (The fullness is a consequence of the Koszul duality of~\cite{GKM}, see e.g.~\cite[Proposition A.6]{BZCHN2}.) In other words, we have a canonical deformation of a full subcategory of the $X$-spherical category. 
 
 \begin{lemma} Conjecture~\ref{local conjecture} implies that all of $\SHV(X_F/G_O)$ is equivariantizable, $$\SHV_u(X_F/G_O)\otimes_{\kk[u]}k\simeq \SHV(X_F/G_O).$$ 
 \end{lemma}
 
 \begin{proof} For the spherical Hecke category $\Hecke_G$ itself the full subcategory of equivariantizable objects is the entire category. Since the basic object is naturally $\Gm$-equivariant, it follows that the full subcategory it generates under the Hecke action is equivariantizable as well.  
 The Lemma follows from the observation (\S \ref{Hecke local conj}) that in the presence of the local conjecture, the entire automorphic category is generated in this way.
\end{proof} 
 
We now require that under the local conjecture this $u$-deformation defines a deformation quantization of $\Mv/\Gv$ with its (sheared) Poisson bracket. We formulate here a very weak form of this requirement:

\begin{itemize}
\item[$\bullet$][Poisson structure]          
The Hochschild cohomology class of $\SHV(X_F/G_O)$ defined by the first-order data of the $u$-deformation $\SHV_u(X_F/G_O)$ is matched with the Hochschild cohomology class of $\QCshear(\Mv/\Gv)$ defined by the (sheared, 2-shifted) symplectic form.
 \end{itemize}

Note that the Hochschild classes in question are in fact in degree 0. Automorphically, the $u$-deformation is a deformation over $k[u]$ where $u$ has cohomological degree 2. Spectrally, since the symplectic form on $\Mv$ has $\GGm$-weight 2, the sheared bracket on $\cO^{\shear}(\Mv)$ has degree $-2$ (i.e., defines a $P_3$-structure) whence again a degree 0 Hochschild class.

The structure theory of hyperspherical varieties can be used to prescribe a canonical global $k[u]$-deformation quantization of $\Mv/\Gv$.
A stronger form of this compatibility requires this entire deformation to match the $u$-deformation of the automorphic category. This is formulated (using the language of Plancherel algebras) as Conjecture~\ref{noncommutative poisson conjecture} in the next section.

 \subsection{Examples} \label{basic examples}
Let us discuss several examples where at least some parts of  Conjecture~\ref{local conjecture} are known.

\begin{example} \label{A1example3}
 We return to our  ``easy'' running example of $X=\mathbb{A}^1$ (see Example \ref{A1example1}, \ref{A1example1.5}, \ref{A1example2}). 
In this case one can verify the equivalence ``by hand,''
as we will essentially do below; we will also see what happens to the ``Dirac delta'' and ``Lebesgue'' sheafs 
from Example \ref{A1example2}. 

Here $\check{M}=T^* \mathbb{A}^1$.  
 In the normalized action, the $G^{\vee}=\Gm$ action
 on $\check{M}$ scales the $\mathbb{A}^1$-direction and antiscales the
 dual $T^* \mathbb{A}^1$-direction.
 Let $x$ be a coordinate on $\mathbb{A}^1$ and $\xi$ a dual coordinate on $T^* \mathbb{A}^1$, so that
 the $G^{\vee}$ action is given by $\lambda \cdot x=  \lambda^{-1} x, \lambda \cdot \xi = \lambda \xi$
 and the $\GGm$ action antiscales both; both $x,\xi$ are in cohomological degree $1$.
 
 We will write
 $[n]$ for the representation $z \mapsto z^n$ of $\check{G}$; the sheaf   $\delta_n$ introduced  in Exanple \ref{A1example1}
 corresponds to $[n] \otimes \mathcal{O}_{\check{M}}^{\shear}$.
 We have computed $\Hom(\delta_n, \delta_m)$ in \eqref{Deltadeltahom};
on the spectral side the same computation yields
 $$ \Hom([n] \otimes \mathcal{O}, [m] \otimes \mathcal{O}) = \Hom_{G^{\vee}}( [n-m] , \mathcal{O}).$$  
 The invariants in question are spanned by monomials $x^a \xi^b$ where $a,b \geq 0$ and $b-a = n-m$.
This lies in cohomological degree $a+b$.  That is to say, the $\Hom$-space
 is one-dimensional  in degrees $|n-m|, |n-m|+2, \dots$,  which matches with \eqref{Deltadeltahom}.
 
Let us also discuss what happens in the limit, to see where
the objects $\delta_0, \mu$ of Example \ref{A1example1.5} go. 
The {\em unnormalized} sheaf $k_{X^l}$ corresponds, on the spectral side, 
to $[l] \otimes \mathcal{O} \langle - l \rangle$,
which we can formally think of as $x^{-l} \mathcal{O}$. 
    The skyscraper at $0$ corresponds 
to taking a limit as $l \rightarrow -\infty$, which gives 
a $\mathcal{O}$-module which can be identified
with $\kk[x, x^{-1},\xi]$.  We can also take the limit of sheaves $k_{X^l}\langle 2l \rangle$
as $l \rightarrow \infty$, corresponding to $[l] \otimes \mathcal{O} \langle l \rangle$
or $\xi^{l} \mathcal{O}$, and get the object $\kk[x, \xi, \xi^{-1}]$.  It might be quite interesting
to study the images of other natural distributions on $X$, under the local conjecture,
in more general instances. 
 \end{example}

\subsubsection{The group case} \label{basic example: group}
 
 Here the space $X$ is a reductive group $H$ with a left and right action of $G=H \times H$; 
and the dual symplectic variety is  now $\Mv=T^\ast \Hv$ as a Hamiltonian $\Hv\times \Hv$-space
where the action of one factor is twisted by the duality involution. In particular, 
we have $\check{G}_X = \Hv$ and the space $V_X$ is 
the dual Lie algebra $(\mathfrak{h}^{\vee})^*$, placed in $\Gm$ weight $2$; 
thus, the conjecture 
 is precisely the derived Satake theorem of \cite{BezFink},
in the case of $\FF=\C$, 
which was recalled as Theorem \ref{Satake thm}.

Let us talk through the various desiderata: 
the requirement on pointings (\S\ref{pointings}) follows from the construction.
The statement about the equivalence
and Hecke actions (\S\ref{Hecke local conj})
{\em almost} follows from the monoidal nature of derived Satake,
but there is a subtlety:
one needs to know that the action of the
inversion map $H \to H$  is induced on the spectral side of $\LL_H$  by the {\em duality involution} 
$\check{H} \rightarrow \check{H}$ defined in \S \ref{dualityinvolution},
for which we don't know a reference. (This  assertion has content
 even in abelian Satake if we work over coefficient field $\kk=\Q$.)

The statement about Poisson structure 
(\S\ref{Poisson local conj} but even the stronger statement discussed in \S \ref{Poisson from loop}) was also proven by Bezrukavnikov and Finkelberg \cite{BezFink}. 
Finally, the 
statement about Galois actions does not apparently
appear in the literature, but
we have sketched a proof after \S \ref{Satake thm-arithmetic};
see also ~\cite{xinwenSatake}, 
~\cite{richarzscholbachmotivic}, and ~\cite{romaICM}.

\subsubsection{The trivial period}
We now take $X$ to be the point.
Thus $\Sph_X$ is simply the category of constructible sheaves on $BG$. That is to say (see \cite{BL}), 
  taking into account the formality of the cohomology of $BG$, 
$\Sph_X$ is simply the category of perfect modules under $H^*(BG)$.

 In this case, $\Mv$ is the Whittaker model, that is to say,
 the Hamiltonian $\Gv$-space obtained by Whittaker indcution
 of the point under the principal $\SL_2 \rightarrow \Gv$;
 and 
 $\Mv/\Gv\simeq \fc_G$ is the Kostant slice.
 The equivalence as abstract categories follows from the identification of the equivariant cohomology ring $H^*(BG)$ with the shear of the ring of invariant polynomials $\cO(\fgxv)^{\Gv}$. The equivalence as module categories for the Hecke category is a theorem of Bezrukavnikov-Finkelberg~\cite{BezFink}  (the same paper also proves the compatibility of the Poisson bracket with loop rotation equivariance).

\subsubsection{The Whittaker period} \label{sssWhittaker}
Dually, we consider the Whittaker model on the automorphic side. In other words, we take $M=T^\ast G\GIT_\psi U$. In this case the Hecke category $\Sph_M$ is the category $\Shv(G_F/G_O)^{U(F),\psi}$ of Whittaker sheaves on the Grassmannian. The dual period is $\Mv=pt$, so that $\QCshear(\Mv/\Gv)\simeq Rep(\Gv)^{\shear}.$ The resulting equivalence of categories is due to Frenkel-Gaitsgory-Vilonen~\cite{FGV} (see also~\cite{ABBGM}).

\subsubsection{Other examples} \label{OtherExamples}
Several other instances of the local conjecture have been proved recently --  verifying at least
the ``Hecke'' desideratum of the statement, although it is likely that the proofs give more. 
  
Hilburn and Raskin \cite{HilburnRaskin} prove a fully ramified, de Rham version of the local conjecture for the Iwasawa--Tate period.

 Braverman, Dhillon, Finkelberg, Raskin and Travkin \cite{RaskinNisyros}, and Teleman \cite{telemancoulomb} discuss geometric counterparts to the Weil representation and the construction of Coulomb branches (hence formulation of the local conjecture) in unpolarized settings.

Chen, Macerato, Nadler, and O'Brien \cite{CMNO} prove the local conjecture in the case $X = GL_{2n}/ Sp_{2n}$ and $\Mv=T^*(GL_{2n}/GL_n \ltimes U, \psi)\simeq GL_{2n}\times^{GL_n} M_n$
(in automorphic terminology, the ``symplectic period'' studied by Jacquet and Rallis \cite{JRsymplectic}). 

 Chen and Nadler \cite{ChenNadler} construct an equivalence, in the case of $X$ a symmetric variety, between the category  $\Hecke^X$ and a  category of sheaves
 arising in global {\em real} Langlands: sheaves on the space of principal bundles for a certain real form of $G$ over the real projective line. They show, among other results,
 the formality of the category $\mathcal{H}^X$ in this case, see \cite[Theorem 13.4]{ChenNadler}. 

The results of Macerato and Taylor \cite{MT} (when combined with \cite{ChenNadler}) imply that the local category for the symmetric variety $X=PSL_2/PO_2$ (which is outside of our hyperspherical setting due to the presence of Type N roots) contains a block equivalent to representations of quantum $SL_2$ at a fourth root of unity. This suggests that in this generality one should expect interesting braided ($E_2$) tensor categories to appear instead of the symmetric ones in our conjecture.

Devalapurkar \cite{Devalapurkar}  proves  the local conjecture for the case of homogeneous affine spherical varieties of rank $1$,   modulo
certain hypotheses relating to placidity and evenness.
Moreover, in \cite{sanathtriple}, under the same hypotheses, he proves the local conjecture in the case of the triple product $L$-function 
(the diagonal inclusion of $\PGL_2$ into $\PGL_2^{\times 3}$).

The papers~\cite{BravermanMirabolic} and ~\cite{BravermanOrthosymplectic} of Braverman, Finkelberg, Ginzburg and Travkin prove the local conjecture (in a Koszul dual formulation)
in cases corresponding to:
\begin{itemize}
\item[-] $X = \GL_n$ as $G = \GL_n \times \GL_{n-1}$-space,  
\item[-] $X=\GL_n\times \AA^n$ as $G = \GL_n \times \GL_n$-space
 \item[-] $X=\SO_n$ as $G =\SO_n \times \SO_{n-1}$-space.
\end{itemize}
 In automorphic terminology, these correspond to the Rankin--Selberg
periods for $\GL_n \times \GL_{n-1}$ and $\GL_n\times \GL_n$ ~\cite{JPSRS} and the Gross--Prasad period for $\SO_n \times \SO_{n-1}$.  Travkin and Yang \cite{TravkinYangI} prove
the local conjecture  in the case of the 
the Rankin-Selberg period for $\GL_m \times \GL_n$ when $m < n-1$; the dual is then
  $\Xv=\mathrm{Mat}_{m\times n}$ acted on by $\Gv=\GL_m\times \GL_n$.

\begin{remark}[The Gaiotto Conjecture] \index{quantum geometric Langlands}
The geometric Langlands correspondence, at least in de Rham and Betti forms, has a so-called {\em quantum} deformation introduced by Feigin, Frenkel and Stoyanovsky, in which the automorphic categories are replaced by categories of twisted sheaves on moduli of $G$-bundles, see~\cite{gaitsgoryquantum, FrenkelGaiotto,BettiLanglands} and references therein. The spectral side has a description, in the Betti version, in terms of categories of representations of quantum groups. It is a natural question which forms of the relative Langlands conjectures admit quantum deformations -- it is easy to see that many periods do not admit such deformations, and even in cases that do (such as the group case) the deformed local category becomes much smaller or trivial. The Gaiotto conjecture provides a quantum deformation of the local conjecture for a family of periods where (like the Whittaker case) the deformed category is of the same size. Remarkably, these cases correspond to quantum supergroups, with $\Gv$ as the even part and $\Mv$ as the odd part, and include all the Rankin-Selberg and Gross-Prasad type integrals in the previous paragraph. 

The Gaiotto conjecture is formulated in~\cite[\S 2]{BravermanMirabolic}, and the quantum analogs of the above cases have now been proven for generic level in~\cite{BravermanOrthosymplecticII},~\cite{TravkinYangI} and ~\cite{TravkinYangII} (the Whittaker case corresponds to the {\em Fundamental Local Equivalence}, proven for irrational level in~\cite{GaitsgoryFLE} and for rational level in~\cite{rationalFLE}). Crucially, these conjectural descriptions are equivalences of {\em monoidal} categories, with the monoidal structure coming from factorization corresponding to the braided tensor product coming from quantum groups. 
 It would be very interesting to compare these spectral models with ours, and to apply the techniques developed in these papers to construct monoidal (or even braided monoidal) structures on the categories $\Hecke^X$ in general, see \S ~\ref{factorizable Hecke} and Problem ~\ref{localconstancyHeckeX}.
 \end{remark}

%% file: local-geometric-Plancherel.tex
\newcommand{\ol}{\overline}
\newcommand{\one}{2}
\newcommand{\two}{1}
\section{The Plancherel algebra and the Coulomb branch} \label{PlancherelCoulomb}
 
In \S \ref{section-unramified-local} (whose notations we follow), we studied the category of constructible sheaves on $X_F/G_O$ as a module category for the spherical Hecke category, and formulated the local conjecture identifying this structure in terms of coherent sheaf theory on $\Mv$ and the moment map. In this section we will study a more concrete and less categorical manifestation of local geometric duality, connecting  the local conjecture to the study of {\em Coulomb branches}, as pioneered in the mathematical literature by Braverman, Finkelberg and Nakajima~\cite{BFN,BFNring}.
This concerns, roughly speaking, the {\em affine} aspects of the story, or 
``the part of the local conjecture that has to do with the basic object.'' \footnote{From the point of view of boundary conditions in topological field theory, we are passing from categories of line operators to vector spaces of local operators}

The key object here is what we call the
 {\em Plancherel algebra} (or {\em relative Coulomb branch}):   $$ \Planch_X =
\mbox{``endomorphisms of the basic object inside the Hecke category'' }$$
The Plancherel algebra is an algebra object in the Hecke category -- that is to say, 
it is equipped with an associative multiplication $\Pl_X \star \Pl_X \rightarrow \Pl_X$, where $\star$
is convolution of sheaves on the affine Grassmannian.
The definition will be explicated and spelled out in \S \ref{Coulomb1} -- it can be viewed as a gadget encoding not just endomorphisms of the basic object, but all Hom spaces between Hecke functors applied to the basic object. 
 The local conjecture predicts that, under the geometric Satake equivalence,
$\Pl_X$ corresponds to the ring of functions $\cO_{\Mv}^{\shear}$, as a representation of $\Gv$,
and compatibly with the moment map (i.e., as an algebra object in sheared $\Gv$-equivariant sheaves on the coadjoint representation $\fgxv$). See 
  \S \ref{OMloc} for more discussion of this point of view. 

The Plancherel algebra $\Planch_X$ is a more concrete object than the local category $\SHV(X_F/G_O)$, and also enjoys an important technical advantage.  As we saw,  ``constructible sheaf theory'' 
in the local conjecture involves serious subtleties related to the 
fact that loop spaces are singular; this makes it extremely difficult to actually compute anything.  However, $\Planch_X$ can
  be computed solely in terms of computations on smooth arc spaces
  and (therefore) its definition reduces entirely to computations on finite-dimensional varieties, as described in
Proposition
\ref{Planch mult def}.  As a result, it is possible to compute explicitly in examples, see e.g.\ \S \ref{Plancherel alg examples}.
  
  We derive this essentially finite-dimensional description of the Plancherel algebra from the general sheaf theory formalism.  
  That sheaf theory formalism is subject to the caveat noted in \S \ref{automorphic local}, namely, it is not well-documented in the literature at the current time.
  However, one could take the point of view that the
    finite-dimensional explication
       is the {\em definition} of the Plancherel algebra and bypass the categorical prerequisites. Moreover,  thanks to the affineness of $\Mv$, the local conjecture implies that one can recover the entire automorphic category from the Plancherel algebra (see Section~\ref{Hecke local conj}).
    
 The name ``Plancherel algebra'' emphasizes the fact, to be further explained in   \S \ref{section:numerical Plancherel},
 that it categorifies an object with a long mathematical history: the {\em Plancherel measure};
 and the 
 known numerical evaluations of the Plancherel measure provide
 supporting evidence for our categorical conjecture -- this will be explained in  \S \ref{section:numerical Plancherel}.

In the case when $X$ is a vector space, $\Planch_X$ is precisely the {\em relative Coulomb branch algebra}, the ring object
in the Hecke category introduced by Braverman, Finkelberg and Nakajima in~\cite{BFNring}
in the course of their mathematical study of Coulomb branches of 3d supersymmetric gauge theories. This ring object (or rather, its spectrum $\mu:\Mv\to \fgxv$) is a {\em relative} version of the Coulomb branch, corresponding physically not to a 3d field theory but to a boundary condition for 4d $\mathcal{N}=4$ supersymmetric gauge theory.
From the point of view of~\cite{BFNring}, perhaps the most novel part of the present chapter  
is to formulate a precise conjecture description of the relative Coulomb branch (and its noncommutative deformation) associated
to a spherical variety.

 We have included a liberal sprinkling of examples.
 Some of these examples correspond to already-proven cases of the local conjecture;
 our purpose in putting them here is, rather, ``for fun'' -- that is,  to illustrate   the pretty geometry
 and invariant theory  underlying the story.

 \begin{itemize}  
 \item \S \ref{OMloc}  motivates and defines the Plancherel algebra, and then formulates the ``Plancherel algebra conjecture,'' Conjecture
\ref{Plancherel algebra conjecture},
  identifying it with the ring of functions on $\Mv$.
 As explained above this is a ``reduced'' version of the full local conjecture which avoids
 the intricacies of sheaf theory on $X_F$.  
 
 \item \S \ref{GXgrassdef} introduces the ``relative Grassmannian''
 which will play a useful role;
 in the case of $X=H\backslash G$ it is essentially the affine Grassmannian for $H$. 
 
  \item \S \ref{Coulomb1} presents an explicit description of the Plancherel algebra, which in particular reduces to finite-dimensional calculations.
 \item \S \ref{Plancherel alg examples}  discusses several examples of Plancherel algebras. 
 
   \item \S \ref{Poisson to loop} studies loop rotation and formulates a conjecture on the corresponding noncommutative deformation of the Plancherel algebra. This explains, in particular, how the symplectic structure on $\Mv$ is relevant to the local conjecture.

   \end{itemize}

We follow general notations as in the last section
\S \ref{section-unramified-local}. In particular, $(G, X)$ will be defined over $\FF$,
either $\C$ or the algebraic closure of a finite field,
and the dual $(\check{G}, \check{M})$ is defined
over a coefficient field $\kk$. In line
with our conventions in
\S \ref{coefficients}, $\kk$ will generally be taken to be algebraically closed of characteristic zero,
although we will comment on more general settings at points.

  \subsection{The Plancherel algebra}\label{OMloc}
 
\subsubsection{Motivation: the ring of functions on $\Mv$ } 

As motivation, let us examine how to recover $\Mv$ from the coherent side of the local 
Conjecture
\ref{local conjecture},
that is to say, from the category of $\Gv$-equivariant sheaves on $\Mv$ (for simplicity,
let us ignore the shearing).

   Since 
 $\Mv$ is affine it is determined by its ring  $\mathcal{O}_{\Mv}$ of globally regular functions.
If we were not working with $\Gv$-equivariant sheaves, this can be recovered by taking endomorphisms of the (coherent) structure sheaf
of $\Mv$.
But, working instead in the category of sheaves on $\Mv/\Gv$, 
the endomorphisms of the structure sheaf recovers not $\mathcal{O}_{\Mv}$
but rather only the $\Gv$-invariant elements.
This can be remedied by computing
endomorphisms {\em as a $\Rep(\Gv)$-enriched category}.
In other words,
the category of $\Gv$-equivariant sheaves is a module category for the rigid tensor category $\Rep(\Gv)$, 
which permits us to compute ``internal endomorphisms in $\Rep(\Gv)$.'' 
(The phrase ``internal endomorphisms'' is most commonly used
in the case of a monoidal category acting {\em on itself} though it applies more generally.)
This is analogous to the discussion in Section~\ref{shearcategory} of $Rep(\Gm)$-categories as categories enriched in graded vector spaces (representations of $\Gm$).
\index{rigid tensor category}
 
Work, for a moment, with ordinary (underived) categories, and take an object $\mathcal{F}$ in the category of coherent sheaves on $\Mv/\Gv$.
Its internal endomorphisms, relative to $\Rep(\Gv)$, is the $\Gv$-representation $\underline{\End}(\mathcal{F})$ characterized
as such by the property
\begin{equation} \label{int end def}\Hom(V,  W \otimes \underline{\End}(\mathcal{F}))=  \Hom(V \otimes \mathcal{F},  W \otimes \mathcal{F}),\end{equation}
 for $V, W$ finite-dimensional $\Gv$-representations;
  the evident composition law on the right hand side, if we work with three representations $V, W, S$, gives
 on $\underline{\End}(\mathcal{F})$ the structure of $\Gv$-algebra. 
 
\begin{remark}[Categorical definition of inner endomorphisms]\label{inner endos}
 A more formal and categorical point of view, well adapted to the $\infty$-categorical setting, is based on the theory of rigid tensor categories briefly reviewed in Appendix ~\ref{rigid tensor categories} (see also~\cite[Section 1.9]{GR}). Namely, given a module category $\cM$ for a rigid tensor category $\cC$ and a compact object $F$, let $act_F:\cC\to \cM$ denote the $\cC$-module functor given by action on $F$. Then rigidity implies that $act_F$ admits a $\cC$-linear right adjoint, which is by definition the {\em inner hom} from $F$, $act_F^R=\underline{\Hom}(F,-)$. The algebra object $\underline{\End}(F)=\underline{\Hom}(F,F)\in Alg(\cC)$ represents the monad $act_F^R\circ act_F$ on $\cC$. Applying this to $\cM=\QC(\Mv/\Gv)$ and $\cC=\Rep(\Gv)$ recovers the construction above.
In particular Equation~\ref{int end def} applies in this generality (as a combination of the $\cC$-linearity of and adjunction between $act_F, act_F^R$), i.e., the internal endomorphism object of $F$ collects the data of all Hom spaces between the objects $V\ot F$ and $W\ot F$ for compact (hence dualizable) objects $V,W\in \cC$. 

Yet another formulation is that we can recover the category $\QC(\Mv)$ with its $\Gv$-action by de-equivariantization of $\QC(\Mv/\Gv)$ with its $\Rep(\Gv)$ action thanks to 1-affineness (as discussed in Section~\ref{shearcategory} for $\Gm$): i.e., we are computing endomorphisms of the $\Gv$-equivariant structure sheaf in the de-equivariantized category $\QC(\Mv)\simeq \QC(\Mv/\Gv)\otimes_{\Rep(\Gv)} Vect$. 
\end{remark} 

\subsubsection{The Plancherel algebra} 
 We have just explained how to reconstruct $\Mv$ from the coherent  sheaf category 
of the local conjecture  \ref{local conjecture}-- namely, by taking internal endomorphisms
of the structure sheaf of $\Mv$ relative to $\Rep(\check{G})$. Since that conjecture predicts that this is the
same as the ``automorphic'' sheaf category $\SHV(X_F/G_O)$, 
it also predicts that we can reconstruct $\Mv$  
by taking internal endomorphisms of the basic object in $\SHV(X_F/G_O)$ relative to the spherical Hecke category.

To spell out: we use that the $X$-spherical category is a module category under the (large) Hecke category $\HECKE_G \simeq \QCshear(\fgxv/\Gv)$, and just as above 
we can compute internal endomorphisms of an object and produce an algebra object in the Hecke category. Spectrally, this corresponds (after shearing) to thinking of $\mu_*\cO_{\Mv}$ as an algebra object in $\QC(\fgxv/\Gv)$, i.e., remembering both the moment map and $\Gv$-action. 

This motivates that we consider the following. 
 \begin{definition}\label{Pl def}   (Plancherel algebra):
The Plancherel algebra of the
$G$-variety $X$ is the algebra object in the Hecke category $\HECKE_G$ defined as the inner endomorphism algebra of the basic object
$\delta_X\in \HECKE^X$ in the $X$-Hecke module category\footnote{Recall that the superscripts over $\mathcal{H}$ refer to the large versions of the categories.}: 
\begin{equation} \label{PlanchXv1} 
\Planch_X :=\underline{End}_{\HECKE_G}(\delta_X)\in Alg(\HECKE_G)
\end{equation}

Via the  derived Satake equivalence,
Theorem \ref{large Satake}, this
$\Planch_X$ defines  (uniquely, up to homotopy)
a differential graded algebra over $\mathcal{O}(\fgxv)^{\shear}$ equipped with a $\check{G}$-action, 
and we will equally well use $\Planch_X$ to denote this algebra.  \footnote{In the cases of interest,
the local conjecture predicts that this differential graded algebra is formal, i.e., equivalent
to its cohomology. Thus, at least in principle, no information would be lost by passing to cohomology here.}

\end{definition}

Here ``internal endomorphisms'' are defined as explained after 
\eqref{int end def}, now working relative to the Hecke category rather than $\Rep(\check{G})$;  in particular, $\Planch_X$ is a unital algebra object of the Hecke category for $G$.   
As we have noted, the Plancherel algebra defined above 
 is {\em entirely reducible to finite-dimensional computations}, see 
 \S \ref{Coulomb1}.

 \begin{remark}
Just as with the previous section,
we should strictly speaking fix a choice of sheaf theory (de Rham or {\'e}tale) and 
include this in the notation, but, again since it matters very little, we don't do it that way.
A formal expression of this irrelevance is Proposition \ref{Planch mult def},
which computes $\Planch_X$ in terms of cohomology groups,  which can be
compared in a standard way between {\'e}tale and de Rham versions.
\end{remark}

 \begin{remark}
 We have only defined the sheaf category in the case when there is no
 $\mathbb{A}^1$-bundle   $\Psi \rightarrow X$. 
 However, once the foundational definitions
 have been suitably extended (cf. \S \ref{subsec generalizations})
the same definition should be adopted in both twisted and unpolarized cases
 (indeed, the paper \cite{RaskinNisyros} already 
 amounts to a definition of $\Planch_X$ for $X$ a symplectic vector space). 
\end{remark}

\subsubsection{The Plancherel algebra conjecture}      \label{Pl conj sec}
   
    Let notation be as in 
\ref{section-unramified-local}; in particular,
we take $(G, M=T^*(X, \Psi))$
and $(\check{G}, \check{M})$ a split hyperspherical pair
as in \ref{goodhypersphericalpairs}
defined over $\FF$ and $\kk$ respectively. We will now spell out the consequence of the local conjecture for the Plancherel algebra $\Planch_X$:
  
  \begin{conj} \label{Plancherel algebra conjecture}
 (A consequence of 
Conjecture~\ref{local conjecture}):

The Plancherel algebra $\Planch_X$ is isomorphic as a $\Gv$-equivariant algebra over $\mathcal{O}(\fgxv)^{\shear}$ with the ring of functions on $\mathcal{O}(\check{M})^{\shear}$:
\begin{equation} \label{plxom} \Planch_X \simeq \mathcal{O}(\check{M})^{\shear}\end{equation} 
(where the shearing on $\mathcal{O}(\check{M})$ is through the $\GGm$-action, and the  sheared moment map $\mathcal{O}(\fgxv)^{\shear} \rightarrow \mathcal{O}(\check{M})^{\shear}$ corresponds to
the unital structure on $\Planch_X$.)
 Moreover,  if $\FF$ is the algebraic closure of a  finite field, this isomorphism is Frobenius-equivariant,
where the Frobenius action on the left hand side
arises from its natural action on $\mathcal{H}_G$ and the trivial Frobenius structure on $\delta_X$; and 
Frobenius is acting on the right hand side according
to the shear of the analytic action (\S \ref{Marithshear}). 
\end{conj}
  
   In particular, $\Planch_X$, {\em a priori} a differential graded algebra, is in fact formal 
   and commutative. On the level of cohomology, the commutativity   follows 
   (again,  bearing in mind the general caveats of \S \ref{automorphic local}) from the compatibility of the product with factorization, as in Section~\ref{automorphic-factorization}, and convolution.
 We also remark that the statement \eqref{plxom} is  using
  analytic normalization of $\Planch_X$ and the \emph{analytic} shear on $\mathcal O_{\check M}$ (\S \ref{analytic-M}).
The corresponding assertion with the arithmetic normalization of $\Planch_X$, as in Remark \ref{remark-Plancherel-AA1},
is that the normalized version of $\Planch_X$ is isomorphic to the {\em arithmetic} shear   on $\mathcal O(\check M)$, \S \ref{Marithshear}.

\begin{remark} \label{Frob sign normalization}  
As a sanity check on signs let us compute the 
 case when $G=\Gm, X =\mathbb{A}^1$, so that $\check{M} =T^* \mathbb{A}^1$
 with scaling action of $\GGm$. 
 Frobenius structures were discussed in
Example \ref{analytic-M-example}.
 
 In this case, we have (cf. \eqref{Deltadeltahom})
 that for a character $\lambda \mapsto \lambda^j$
 the corresponding isotypical component $\Planch_X^{(j)}$, 
 the component by which $\check{G}$ acts by $\lambda \mapsto \lambda^j$, 
 is identified with $H^*(B\Gm)\langle -|j| \rangle$.\footnote{
 we emphasize that we have used the analytic normalization in this computation, 
 i.e., the Hecke operator $T_j$ corresponding
 to $\lambda \mapsto \lambda_j$ is supported in valuation $j$ and $\deg(T_j) = j$, cf. Example \ref{A1example2}.}
 This is one-dimensional exactly in degrees $j, j+2, \dots$,
 and the geometric Frobenius eigenvalue in degree $s$ is given
 by $q^{s/2}$. 
 
On the other hand, writing $x, \xi$ for the coordinates
on $\check{M}$, the isotypical space for $\lambda \mapsto \lambda^j$
on $\mathcal{O}(\check{M})$ is exactly
$x^j (x \xi)^m, m \geq 0$ or $\xi^{-j} (x\xi)^m, m \geq 0$
depending on whether $j$ is positive or negative. 
Again, this is one-dimensional exactly in $\GGm$ degrees
$-j, -j-2, \dots$,  corresponding to cohomological degrees $j, j+2, \dots$
(cf. Example \ref{analytic-M-example}), thus matching the geometric computation.

\end{remark}

   \begin{remark}   \label{localconjecturesplitform}
(Rationality issues:) 
  We would also conjecture that Conjecture
  \ref{Plancherel algebra conjecture}
  remains valid
  if $\kk$ is not algebraically closed,
   and even if $\kk$ is not a field, 
   where $\check{M}$ is taken to be the distinguished split form
   that was postulated in \S \ref{GMdesiderata}. 
   In the case where $\Psi$ is trivial we can even take $\kk$ to equal $\Z$.
 Indeed, in \S \ref{GMdesiderata}, we noted that compatibility
with this Plancherel algebra conjecture should be considered a defining property
of this split form.

  Note that even in the group case that the predicted $\Q$-form 
  of $\check{M}$ arising in the conjecture is not the ``pinned'' version;
  that is to say, $\check{M}$ is the cotangent bundle of $\check{G}$ as $\check{G} \times \check{G}$-space,
  but the action is twisted, on one factor, not by the Chevalley involution but rather by the duality
  involution (\S \ref{dualityinvolution}) which does not preserve a pinning.  
  
  It is extremely interesting to study these issues when $k$ has small characteristic. Work over $\FF=\C$   and consider, for example, the case when $X=G/H$ is homogeneous; we will see that, in this case, the invariants
  $\Planch_X^{\check{G}}$ are given simply by the cohomology $H^*(BH, k)$, the cohomology
  of the classifying space of $H$ with coefficients in $k$. Now, if the characteristic of $k$ is a bad prime
  for $H$, this cohomology can behave quite irregularly. On the other hand, the invariants $\Planch_X^{\check{G}}$
  can {\em also} behave badly: in finite characteristic, taking $\check{G}$-invariants is not exact, and one should use instead derived invariants.
  It is plausible that these pathologies match up with one another. Again, this is not even
  obvious in the group case, or in the case when $X=\mathrm{pt}$. 
  
   \end{remark}

\subsection{The relative Grassmannian} \label{GXgrassdef}
The following object (known as the ``variety of triples'' in the Coulomb branch literature~\cite{BFN}) will play a key role for us.
We use the term  ``relative Grassmannian'' and
 we have emphasized the intuition arising from the homogeneous case $X=H\backslash G$.
 where it is essentially the affine Grassmannian of the stabilizer $H$.  
 
Let $X$ be a quasi-affine $G$-variety. 
Informally, the relative Grassmannian $\mathrm{Gr}^X$ is the closed subscheme of $X_O \times \Gr_G  $ classifying  \begin{quote}
 pairs $x\in X_O, g \in G_F/G_O$ for which $xg \in X_O$.
\end{quote} 
where we identify $\Gr_G$ with $G_F/G_O$. 
In particular, it will come with a morphism 
 \begin{equation} \label{GrXdef} i^{\Gr}:\Gr^X \rightarrow  X_O \times  \Gr_G \end{equation}

Formally, we can define $\mathrm{Gr}^X$ as the  pullback 
\begin{equation} \label{GrXdef2} \xymatrix{\Gr^X\ar[r]^{b \qquad } \ar[d]^{a}& X_O \times_{G_O} G_F \ar[d]\\
X_O\ar[r]^i& X_F}\end{equation}  of the action of $G_F$ over the inclusion of the arc space. 
This pullback is taken in the category of schemes,\footnote{
where $G_F \times^{G_O} X_O$
is the scheme that represents the functor assigning to a ring $R$
the $G$-bundles on the formal disc $R[[t]]$ which are trivialized on $R((t))$
and moreover equipped with a section of the associated $X$-bundle.}
see \S \ref{explication grX} for further explication. 

It should be noted that the top horizontal arrow $b$ is not the same morphism
(indeed, does not even have the same codomain) as $i^{\Gr}$ from \eqref{GrXdef}.
Rather,  the  embedding $i^{\Gr}$ of \eqref{GrXdef}
is given by $\iota = (a, b_2^{-1})$, where $b_2: \Gr^X \rightarrow G_F/G_O$
is the second coordinate of $b$. 
Said differently, the horizontal arrow $b$ is represented,
in terms of pairs $(x,g)$ with $xg \in X_O$, 
 by ``$(x,g) \mapsto (xg, g^{-1})$.''

Next,   $\Gr^X$ is equipped with an action of $G_O$. 
If we consider $\Gr^X$ via \eqref{GrXdef} as classifying pairs $(x,g)$, the action
of $h \in G_O$ sends this to $(xh, h^{-1} g)$; if we think in 
terms of the pullback diagram   \eqref{GrXdef2} $G_O$ acts   both on $X_O$ at the bottom left,
and acts by $X_O \times_{G_O }G_F$ by right multiplication on $G_F$.

The quasi-affineness of $X$ guarantees (see below) that $\Gr^X$ defines a locally closed subscheme of $ X_O \times \Gr_G$. 
In the notation before \eqref{GrXdef}, the maps $(x,g) \mapsto g, (x,g)  \mapsto x,  (x,g) \mapsto xg, $   \index{$\Gamma$} \index{$p$} 
 define 
\begin{equation} \label{Gammadef} \Gamma: \Gr^X \rightarrow \mathrm{Gr}, p_1:   \Gr^X \rightarrow X_O, p_2: \Gr^X/G_O \rightarrow X_O/G_O\end{equation}
(note that $p_2$ is only defined after quotienting by $G_O$). 
\index{$\Gr^X$}

\subsubsection{Explication as a (representable) functor} \label{explication grX}

We will describe explicitly the functor represented by $\Gr^X$, and why it is represented
by a locally closed subscheme of $X_O \times \Gr_G$.

One can describe $\Gr^X$ as classifying a $G$-bundle $\cG$ on $\PP^1$ with a trivialization outside of $0$ and an arc in $X$ satisfying the condition that the arc, restricted to the punctured disc at $0$, extends to a section of the associated $X$-bundle of $\cG$ on the disc.

To say this symbolically,   $X_O \times \Gr_G$ is the functor assigning to a test ring $R$ the data of
\begin{itemize}
\item[(i)]
a $G$-bundle $\mathcal{G}$ on $\mathbb{P}^1_R$
trivialized on $\mathbb{A}^1_R$;
\item[(ii)] an element $x \in X(R[[t]])$.
\end{itemize}
Then, $\Gr^X$ is the subfunctor defined as follows: 
Let $\mathcal{X}$ be the $X$-bundle  on $\mathbb{P}^1_R$ associated to $\mathcal{G}$.
The trivialization of $\mathcal{G}$ on $\mathbb{A}^1$ gives
a corresponding trivialization of $\mathcal{X}$.
Therefore, $x$ determines a section $x^{\mathcal{G}}$ of $\mathcal{X}$ over $R((t))$, 
and $\Gr^X$ is determined by the condition that 
$x^{\mathcal{G}}$ extend to an $R[[t]]$-point of $\mathcal{X}$. 
Such an extension is unique for $X$ affine since
$R[[t]] \rightarrow R((t))$ induces an injection on points,
and moreover  the resulting functor is, in the case when $X$ is affine, in fact representable by a (ind-)closed subscheme
of $X_O \times \Gr_G$ (in what follows, we will omit the ``ind-''). 

To verify closedness, we reduce this to the case where $G=\GL_n$ and $X = \mathbb{A}^n$
by choosing an equivariant embedding of $X$ into a vector space;
and in that case
the $R$-points of $X_O \times \Gr_G$
are pairs:
\begin{itemize}
\item[-]
a locally free $R[[t]]$-submodule $M \subset R((t))^n$
such that $M[t^{-1}]$ is all of $R((t))^n$; 
\item[-] and   $x \in R[[t]]^n$. 
\end{itemize}
The condition defining $\Gr^X$ is then that $x \in M$, and this is readily verified to define a ind-closed subscheme. 

 If $X$ is {\em quasi-affine}, 
 the functor noted above  remains representable by a locally closed subscheme.
 Indeed, as before, the above  construction for  the affine closure $X^{\aff}$ of $X$
 yields a closed subscheme of $  X^{\mathrm{aff}}_O \times \Gr_G$;
 we then impose the additional {\em open} condition that $x^{\mathcal{G}}$
 reduce at the special point $R[[t]] \rightarrow R$
 to a point of $X^{\aff}(R)$ that in fact lies inside $X(R)$.

\subsubsection{The equivariant relative Grassmannian}

In the homogeneous case $X=H\backslash G$, we have an equivalence
\begin{equation} \label{Stabilizer Gr 1} \Gr^X/G_O \simeq H_O \backslash H_F / H_O,\end{equation}
  i.e., considered equivariantly relative to $G_O$, $\Gr^X$ 
  recovers  the $H_O$-equivariant affine Grassmannian of $H$. 
\footnote{This example perhaps motivates our name ``relative Grassmannian.''}
  
 More generally the equivariant version $\Gr^X/G_O$ can be described as the moduli stack of:
 \begin{itemize}
 \item pairs of $G$-bundles on the disc, together with
 \item  sections on the associated $X$-bundles on the discs,  and
 \item an identification of the bundles and sections on the punctured disc.
 \end{itemize}
Or, to put another way,  it is the groupoid object 
over $X_O/G_O$ obtained by restricting the Hecke action 
 of the equivariant Grasmannian  on $X_F/G_O$. This explains formally why the relative Grassmannian appears when studying the inner endomorphisms of the basic object (the constant sheaf on $X_O/G_O$) with respect to the spherical Hecke action.

 \subsubsection{The twisted case}
 If $X$ is equipped with an $\mathbb{A}^1$-torsor $\Psi \rightarrow X$ then 
 we get a morphism 
 \begin{equation} \label{GXpsi}A: \Gr^X \rightarrow \mathbb{A}^1, \end{equation}

Given a point of $\Gr^X$ represented by $(x \in X_O, g  \in G_F)$, 
 we lift $x$ to $\tilde{x} \in \Psi_O$,  and then consider $\tilde{x} g \in \Psi_F$;
 since it lies above $xg \in X_O$  it defines a canonical point
 in $F/O$ by comparison to an arbitrary lift $\widetilde{xg} \in \Psi_O$.
 Using the ``residue'' morphism $F/O \rightarrow \mathbb{A}^1$
 we get the desired map \eqref{GXpsi}. 
 For later use we observe the cocycle property
 of the map $A$ from \eqref{GXpsi}. 
\begin{equation} \label{Acocycle}
A(x, g_1 g_2) = A(x, g_1) + A(x g_1, g_2)
\end{equation}

\subsubsection{Essential finite-dimensionality} \label{fin dim model of Gr}
 
The relative Grassmannian is essentially a finite-dimensional construction, in the following sense.
Suppose we 
replace the role of the affine Grassmannian above by some
sufficiently large stratum $\mathrm{Gr}_{\leq n}$.
We obtain a corresponding locally closed subscheme $\Gr^{X}_{\leq n}$.
 The resulting
  $\Gr^X_{\leq n} \rightarrow X_O \times \mathrm{Gr}_{\leq n} $ is pulled back from a
locally closed immersion of  (finite dimensional!) varieties
\begin{equation} \label{GXfinite} i^{\Gr,N}:\Gr^X_{\leq n; N} \rightarrow X_{N} \times  \mathrm{Gr}_{\leq n} \end{equation}
   (for some $N$; here $X_N = X_{O/t^N}$ is the jet scheme of order $N$). 
    All this says is, if we fix the ``denominators''
   of $g$, the question of whether $xg \in X_O$ is entirely
   determined by the reduction of $x$ modulo a fixed power of $t$.   In all computations
   where we use $\Gr^X$ its role can be replaced by one of these finite-dimensional truncations.

 \subsection{Explication of the Plancherel algebra} \label{Coulomb1} 
  
 We are ready to give our explicit and inherently ``finite-dimensional''  description of the Plancherel algebra. 
  
   We describe the Plancherel algebra as 
 a direct sum
  $$\Planch_X = \bigoplus V \otimes \Planch_X^{(V)}$$ of multiplicity spaces for irreducible representations $V$ of $\Gv$, which we will describe individually.

  Let $V$ be a representation of $\check{G}$, let $\mathcal{T}_V$
be as in \eqref{TVcaldef}
  the associated perverse sheaf on $\mathrm{Gr}_G$. 
Fix a sufficiently large closed $G_O$-invariant subset $\Gr_{\leq n}$ on the affine
Grassmannian containing the support of $\mathcal{T}_V$.
 We present the formula for the {\em analytically normalized} Hecke action (see Remark \ref{remark-Plancherel-AA1} below for modifications in the other case).  Thus $\deg T_V$ is defined as in \S \ref{normalized-local} \index{$\deg T_V$} \index{degree of Hecke operator}
in terms of the eigenmeasure on $X$, or equivalently 
is the weight by which $\check{\eta}: \Gm \rightarrow \check{G}$
acts on $V$.     Finally, recall from \eqref{GXfinite} that 
$\iota: \Gr^X \rightarrow  X_O \times \Gr_G$ is  pulled back from a finite type 
$\iota: \Gr^X_{\leq n; N} \hookrightarrow X_{N} \times  \mathrm{Gr}_{\leq n},$ 
 where $X_N$ is the mod $t^N$ truncation of $X_O$ and $\iota$ a locally closed immersion.
Let $p: \Gr^X_{\leq n; N} \rightarrow X_N$ be the second coordinate of this truncated map. 

\begin{prop}  \label{Planch mult def} (Multiplicity spaces of the Plancherel algebra):
With notation as above, we have
\begin{equation}\label{shiftystuff0} 
\Planch_X^{(V)}\langle -\deg T_V \rangle = \Hom_G ( \mathcal{T}^V_X, p^{!} \kk_{X_N}) 
\simeq H^*_G(\mathbb{D} \mathcal{T}^V_X) \langle -2 \dim X_N \rangle
 \end{equation}
 where we compute the $\Hom$- on the right $G$-equivariantly\footnote{More naturally, we would say $G_O$-equivariantly;
 but since the map $G \rightarrow G_O$ induces a homotopy equivalence it amounts to the same thing, and 
 this way the computation visibly only involves finite-dimensional data.} on  the finite type scheme
 $\Gr^X_{\leq n; N}$;    $\mathcal{T}^X_V$ is the $*$-pullback of the Hecke operator
$\mathcal{T}_V$ to $\Gr^X_{\leq n;N}$, and $\mathbb{D}$
is computed for ordinary (not equivariant) sheaves on $\Gr^X_{\leq n;N}$.\footnote{More properly, of course,  $H^*_G$
should be replaced by a chain-level model, since $\Planch_X^V$ is formally a chain complex,
well-defined up to homotopy.}
  \end{prop}

   We will give the proof in \S \ref{PMDproof} and examples in the
  next section \S \ref{Plancherel alg examples}. At the moment we make several remarks.

  \begin{remark} [Comparison with~\cite{BFNring}]
The above statement was written in a way that was evidently finite-dimensional, that is to say,
it involves computing  cohomology on a finite-dimensional algebraic variety, equivariantly for the
action of a finite-dimensional algebraic group.  
 To facilitate comparison with  \cite{BFNring} we can rewrite   the right hand side of \eqref{shiftystuff0}, up to a shift that we will not explicate,  in the following way:
 \begin{equation} \label{PlanchXv2} 
\Planch_X^{(V)}\langle -\deg T_V \rangle = \Hom_{\Gr^X }( \mathcal{T}_V^X, \omega^{\textrm{ren}}_{\Gr^X}) \simeq\Hom_{\Gr}( \mathcal{T}_V, \Gamma_*\omega^{\textrm{ren}}_{\Gr^X}) 
\end{equation}
Here $\Gamma$ is the projection from $\Gr^X$ to $\Gr_G$. Strictly speaking, we first choose a stratum $\Gr_{\leq n}$ containing the support and compute everywhere on the corresponding restricted spaces, so
  the $\Hom$s are taken in $G_O$-equivariant $*$-sheaves on $\Gr^X_{\leq n}$ and  $\Gr_{\leq n}$, respectively; and $\omega^{ren}$ is the renormalized dualizing sheaf,
which, in the case at hand, is the system of sheaves $\omega \langle - 2 \dim \rangle$ 
on the various truncations $\Gr^X_{\leq n;N}$. \footnote{Recall that   sheaf theory on infinite type schemes comes in two variants, $*$ and $!$, according
to whether one takes, informally speaking, a system of   sheaves on finite type truncations
that are compatible with respect to $*$-pullbacks or with respect to $!$-pullbacks. In the prior Section
we used by default $!$-sheaves; in the current situation they are equivalent
by the placidity of $\Gr^X$, whch follows from \eqref{GXfinite}. }

  In~\cite{BFNring}, it is shown that (a suitable shift of) $\mathfrak{A}:=\Gamma_* \omega_{\Gr^X}$ has the structure of ring object in the spherical Hecke category, with $V$-isotypic components therefore given by \eqref{PlanchXv2}.  In fact \cite{BFNring} discuss
only the case of $X$ a vector space;  their construction
applies, however, without change for $X$ a smooth affine variety
(although not to the quasi-affine case).

It is natural to expect, and would be interesting to prove,  that this ring structure agrees with the algebra structure that arises 
via the definition of the Plancherel algebra as inner endomorphisms of the basic object (and such an identification should be completely formal, up to potential issues arising from
  subtleties of sheaf theory in infinite type).
 \end{remark}

  \begin{remark} \label{Whitdef} 

If $X$ is equipped with a $\mathbb{A}^1$-bundle $\Psi \rightarrow X$
and we work in either de Rham or {\'e}tale sheaf theory,
the formulas \eqref{PlanchXv2} or \eqref{shiftystuff0} still make sense by twisting  $p^{!} k_{X_n}$ or $\omega^{\textrm{ren}}_{\Gr^X}$ by $\mathcal{L}_{\psi}$,
the pull-back of an Artin-Schreier sheaf via \eqref{GXpsi}.

This will, presumably, coincide with Definition \ref{Pl def} after it is extended
to cover the twisted setting.  In the absence of this extension,  let us
 take \eqref{PlanchXv2} as the definition of $\Planch_X^{(V)}$ in the twisted case. 
 
We note that the direct construction of the product from \cite{BFNring}
  does not apply in the case of twisted polarizations
 because a certain map  
 fails to be proper in the quasi-affine case.
One expects that the twisting provided by 
 the Artin-Schreier sheaf makes it behave, from a cohomological standpoint,  as if it were proper,    but this remains to be worked out. 
\end{remark}

\begin{remark} \label{remark-Plancherel-AA1} (Normalized versus unnormalized version:) 
Note that \eqref{PlanchXv2} is an analytic normalization, in the sense
of \S \ref{analyticarithmetic}. For the ``unnormalized'' or ``arithmetic'' Plancherel algebra
we use the same definition, but omit the shear $\langle \deg T_V \rangle$. By default, statements
that follow will use the analytic normalization except where otherwise noted.
\end{remark}

 \begin{remark}[Trace of Frobenius]
 Restricting to the case of $\FF$ the algebraic closure of a finite field $\mathbb{F}_q$, 
 let us compute the trace of geometric Frobenius
 on the dual of $\Planch_X^{(V)}$ using Lemma \ref{Homlemma}.  
 The result of this computation will be used in the study of Plancherel measure. 
As before, we work with the analytic normalization; for the arithmetic version one ignores the $\sqrt{\eta}$. 
 
 The trace-function associated to $\mathcal{T}^X_V$
 is the pullback of the usual Hecke function $T_V$ associated to $\mathcal{T}_V$
 from the affine Grassmannian. The trace function
 associated to the $\mathbb{D} p^{!} \kk_{X_N}
 = p^* \omega_{X_N}$ is the pullback of the function $q^{-\dim X_N}$
 from $X_N$.  
 
Applying Lemma \ref{Homlemma},  the trace of Frobenius on the dual
 of $\Planch_X^V$ equals 
 the sum, over $(x,g) \in  \Gr^X_{\leq n; N}(\mathbb{F}_q)$,
 of $q^{\dim(G)} \frac{ T_V(g) \sqrt{\eta(g)}}{q^{\dim X_N} \# G(\mathbb{F}_q)},$
 where the order of $G(\mathbb{F}_q)$ in the denominator arises from
 the corresponding term in Lemma \ref{Homlemma},
 and $q^{\dim(G)}$ arises from the difference between $G$-equivariant
 duality on $\Gr^X_{\leq n; N}$, which is what appears in Lemma \ref{Homlemma}, 
 and ordinary non-equivariant duality which is what appears in the  formula
 \eqref{shiftystuff0}. 
  
 Recall our notations $\mathfrak{o}, \mathfrak{f}$
 for $\mathbb{F}_q[[t]]$ and $\mathbb{F}_q((t))$ respectively.
 Then $\Gr^X_{\leq n; N}$ consists of pairs $(x,g)$, 
 where $x \in X(\mathfrak{o}/t^N)$, 
 $g \in (G_{\mathfrak{f}}/G_{\mathfrak{o}})^{\leq n}$,
 and $xg$ is integral. 
 Rewriting in terms of $\mathfrak{o}$-points, we get 
\begin{multline}\label{Frob Inverse Trace} 
 \mbox{trace of  Frobenius
on dual of $\Planch_X^{(V)}$} 
  = \frac{q^{\dim(G)}}{\# G(\mathbb{F}_q)} \cdot   \int_{X_{\mathfrak{o}} \times  G_{\mathfrak{f}}/G_{\mathfrak{o}}}
  \sqrt{\eta(g)} 1_{X(\mathfrak{o})}(xg) T_V(g) dg, \
 \end{multline} 
and the measure normalizations assign mass $1$ to $G(\fo) \subset G(\ff)$, 
 and assign mass  $\#X(\mathbb{F}_q)/q^{\dim(X)}$ to $X(\fo)$. 
 
 In the case when $X$ is equipped with an affine bundle $\Psi$, 
 the same formula holds when we modify the definition of $\Planch_X$
 according to Remark \ref{Whitdef}, and including in the
integrand the twisting factor
 $\psi(A(x,g))$ where $A(x,g): X_{\mathfrak{o}} \times G_{\mathfrak{f}}/G_{\mathfrak{o}} \rightarrow \mathbb{F}$
 is obtained from the $\mathbb{F}$-points of \eqref{GXpsi}, and
 $\psi$ is an additive character of $\mathbb{F}$. 

\end{remark}

\subsubsection{Proof of  Proposition~\ref{Planch mult def}} \label{PMDproof} 
At a formal level 
we want to rewrite $\Hom(T \star \delta_X, \delta_X)$
by regarding  $\delta_X$
as the $*$-pushforward along $i: X_O \rightarrow X_F$,
and use $(i^*, i_*)$ adjunction  and base change in the diagram \eqref{GrXdef}.
However, we have to proceed with caution because
we can only use the adjunctions available in our infinite-dimensional $!$-sheaf theory. 
  Let us introduce some notation. 
 Along with $\Gr^X$ it will be convenient to define: 
 $$ \Gr^X_{F,O}=  X_O \times^{G_O} G_F  \mbox{ and } \Gr^X_F =  X_F \times^{G_O} G_F .$$
 We consider the Hecke action groupoid on $X_F/G_O$ and its restrictions (one leg at a time) to $X_O/G_O$. 
 
 $$\xymatrix{X_O/G_O \ar[dd]^-{i}& \ar[l]_-{p_{\one}} \Gr^{X}/G_O \ar[d]^-{i^{\Gr}} \ar[r]^-{p_{\two}} & X_O/G_O\ar[dd]^{i} \\
 & \ar[dl]_-{\wt{p}_{\one}}    X_O \times^{G_O} G_F/G_O  \ar[d]^-{\wt{i}} \ar[ur]^-{\wt{p}_{\two}}  & \\
 X_F/G_O & \ar[l]_{p_{\one}^F}  X_F \times^{G_O} G_F/G_O \ar[r]^-{p_{\two}^F} & X_F/G_O}
$$
 
 We use the notation indicated in the diagram and the additional notations $\Gamma,\wt{\Gamma},$ and $\Gamma_F$ for the projections to the Hecke stack $G_O\backslash G_F/G_O$ from $\Gr^X$, $\Gr^X_{F,O}$, and $\Gr^X_F$, respectively. 
  Recall that we work in the formalism of $!$-sheaves in infinite type as in Appendix~\ref{infinite type}, and that, correspondingly, the Hecke action is given by  $$\mathcal T\ast\delta_X:= (p_{\one}^F)_*((\Gamma^F)^! \mathcal T \ot^{!} (p_{\two}^F)^! \delta_X)$$ where the basic object $\delta_X$ is the pushforward of the dualizing ($!$-)sheaf on $X_O$.
We then calculate (using only base change, ind-proper adjunction, the projection formula and the unit for the $!$-tensor product) as follows--
where all $\Hom$s are computed $G_O$-equivariantly, but we drop this from the explicit notation after the first line:  \begin{eqnarray*}
 \Hom_{X_F/G_O}(\mathcal T\ast \delta_X, \delta_X) &\simeq & \Hom_{X_F/G_O}((p_{\one}^F)_*(  \Gamma_F^! \mathcal T \ot^{!} (p_{\two}^F)^! \delta_X), \delta_X)\\
 &\simeq &  \Hom_{\Gr_F^X}( \Gamma_F^! \mathcal T \ot^{!} (p_{\two}^F)^! i_*\omega_{X_O}, (p_{\one}^F)^! \delta_X)\\
  &\simeq & \Hom_{\Gr_F^X}( \Gamma_F^! \mathcal T \ot^{!} \wt{i}_* \omega_{\Gr_{F,O}^X}, (p_{\one}^F)^! \delta_X)\\
 &\simeq & \Hom_{\Gr_F^X}( \wt{i}_*( \wt{\Gamma}^! \mathcal T \ot^{!}  \omega_{\Gr_{F,O}^X}), (p_{\one}^F)^! \delta_X) \\
  &\simeq & \Hom_{\Gr_F^X}( \wt{i}_* \wt{\Gamma}^! \mathcal T, (p_{\one}^F)^! \delta_X) \\
 &\simeq &  \Hom_{\Gr_{F,O}^X}(  \wt{\Gamma}^! \mathcal T, \wt{i}^!(p_{\one}^F)^! \delta_X) \\
 &\simeq & \Hom_{\Gr_{F,O}^X}(  \wt{\Gamma}^! \mathcal T, (\wt{p}_{\one})^! \delta_X) \\
    &\simeq &  \Hom_{\Gr_{F,O}^X}(  \wt{\Gamma}^! \mathcal T, i^{\Gr}_*\omega_{\Gr^X})
 \end{eqnarray*}

 We will now descend to a finite-type computation. 
 Choosing a sufficiently large stratum $\Gr_{\leq n}$
 of the affine Grassmannian which supports $\mathcal{T}$,
 one can replace the role of $\Gr^X$ and $\Gr_{F,O}^X$ by the corresponding
cut-off versions $\Gr^X_{\leq n}$ and $\Gr_{F,O,\leq n}^X$. 
The pullback functor from sheaves on the truncated version
\begin{equation} \label{trunco} X_{N} \times^{G_O} \Gr_{\leq n}\end{equation}  to 
 $\Gr_{F,O,\leq n}^X$ is fully faithful. Therefore, we may
 compute the $\Hom$ on \eqref{trunco}.    $\tilde{\Gamma}^{!} \mathcal{T}$
 corresponds to the ordinary sheaf $\tilde{\Gamma}_N^{!} \mathcal{T}$
on \eqref{trunco}, 
 where $\tilde{\Gamma}_N$ is the analogous projection to $\Gr_{\leq n}$.
 Similarly (as explained in ~\cite[Remark 3.9.1]{raskininfinite}), thanks to base change the $!$-sheaf pushforward $i_*^{\Gr} \omega_{\Gr^X_{\leq n}}$ is pulled back
  from the ordinary sheaf $i_{N*}^{\Gr} \omega_{\Gr^X_{\leq n;N}}$
on \eqref{trunco}, where 
$i^{\Gr}_N: \Gr^X_{\leq n;N} \rightarrow \eqref{trunco}$
  is the level $N$ truncation of the top horizontal arrow of \eqref{GrXdef}.
  \footnote{Note that this is not the same as the truncation defined earlier,
  in \S \ref{fin dim model of Gr}.
  Namely, this truncation identifies elements $(x,g) \in \Gr^X$
  when the first $N$ ``digits'' of $xg$ agree, 
  whereas our previous truncation identifies elements $(x,g) \in \Gr^X$
  when the first $N$ digits of $x$ agree.   However, it will
  make no difference, because we can pass to a common smooth cover of both of the truncations --
  indeed, either truncation after sufficiently increasing $N$ provides such a smooth cover.}
 Therefore,  the $\Hom$-group in question becomes 
 $$ \Hom(\tilde{\Gamma}_N^! \mathcal{T}, i^{\Gr}_{N*} \omega_{\Gr^X_{\leq n,N}}),$$
 computed in $G_O$-equivariant sheaves on \eqref{trunco}.
 Upon using adjunction on the finite type schemes involved
 and the isomorphism $\tilde{\Gamma}^!_{\leq n} \simeq \tilde{\Gamma}^*_{\leq n} \langle 2 \dim X_n \rangle$,
this becomes
 $$ \Hom_{\Gr^X_{\leq n,N}}( i^{\Gr*}_N \tilde{\Gamma}_N^* \mathcal{T}, \omega_{\Gr^X_{\leq n}}). 
 \langle - 2 \dim X_n \rangle$$
 which agrees (recalling that we are working $G_O$-equivariantly everywhere) with \eqref{shiftystuff0}. 
  As noted earlier, this truncation $\Gr^X_{\leq n, N}$ does not
 coincide with that considered before but  the results coincide,
 because one can pass to a common smooth cover.

 \subsection{Some examples}\label{Plancherel alg examples}

The following examples are not given in the spirit of evidence,
but rather to illustrate that the objects appearing in the Conjecture
 are interesting and computable. {\em Everywhere here we will, for the purpose of computing,  work
 at the level of cohomology, i.e., $\Planch_X^V$ refers to the cohomology
 of the underlying chain complex.}
 
 \begin{example}\label{homogeneouscase}
  In the homogeneous case the discussion is particularly simple:
 we have
 for $X=H\backslash G$, with $H$ reductive, 
 the eigenmeasure character is trivial, and after 
 \eqref{Stabilizer Gr 1} we have the equality
\begin{equation} \label{GXhom} \Planch_X^{(V)} \simeq H^*_{H_O}(\Gr_H, \mathbb{D} \mathcal{T}_V),\end{equation}
 i.e., we take the Hecke operator on the $G$-affine Grassmannian,
 $*$-restrict to $\mathrm{Gr}_H$, dualize, and take equivariant cohomology. 
  ``Duality'' is normalized to preserve the basic object of $\mathrm{Gr}_H$,
 and then no further shifts are required; thus, for example, when $V$ is trivial,
 $\Planch_X^{(V)}$ gives the $H_O$-equivariant cohomology of $\Gr_H$ itself.

Here is a simple example where one can compute everything. Take $G=\PGL_2^3$ and $X = \PGL_2^3/\Delta \PGL_2$. Take
$V = V_a \otimes V_b \otimes V_c$
 an arbitrary irreducible of the dual $\SL_2^3$. 
 Then
 $ \Planch_X^{(V)}$ is  the cohomology of the $\PGL_2$-affine Grassmannian
 with coefficients in the tensor product of the  perverse sheaves associated by geometric Satake to $V_a, V_b, V_c$.
 This tensor product is zero unless $a,b,c$ all have the same parity;
 and in that case it is the constant sheaf $k \langle a + b + c \rangle$
 on the $m$-dimensional stratum of $\Gr^{\PGL_2}$ indexed by $m=\min(a,b,c)$.
Geometric Satake here tells us that the cohomology of $k$
on this stratum is identified with the cohomology of $\mathbb{P}^{\min(a,b,c)}$. 
Dualizing and computing cohomology we find 
a noncanonical isomorphism
of graded vector spaces (with Frobenius action when $\FF$ is of finite characteristic)
 \begin{equation} \label{chugchug}
\Planch^{(V)}_X \simeq 
 H^*(\mathbb{P}^{\min(a,b,c)}) \otimes H^*(B\PGL_2) \langle 2 \min(a,b,c)-a-b-c\rangle,\end{equation}
-- the noncanonicity arises from the fact that we used the degeneration of the pertinent spectral sequence. 
  Now let us   check the $V_a \otimes V_b \otimes V_c$-isotypical component of
 Conjecture \ref{plxom} including the statement about Frobenius structures.  
 Here we have $\check{G}=\SL_2^3$
and $\check{M} $ the standard representation, i.e., the tensor product
of the three standard representations, and so 
$$  \eqref{chugchug} \stackrel{?}{\simeq}
\Hom(V_a \otimes V_b \otimes V_c, 
\kk[\mathbb{A}^2 \otimes \mathbb{A}^2 \otimes \mathbb{A}^2]),$$
  On the right hand side, considering $\kk[\dots]$
 as polynomial on an $8$-dimensional vector space, the Frobenius action
 on degree $d$ polynomials is through a twist $\langle -d \rangle$, i.e.
 by $q^{d/2}$ (note that $d$ is, by our conventions, negative to the $\GGm$ grading). 
      To compute this  we note that  the $\SL_2$-action
 on $\mathbb{A}^8$ has a unique (up to scaling) quartic  invariant.
 This is well-known, particularly to those people who know it well. 
 The zero-locus of this semi-invariant has the form
$  U^0 \backslash \SL_2^3$
 where $U^0$ is the subgroup of upper triangular unipotent
 matrices $u_{x,y,z}$ satisfying $x+y+z=0$;
 the action of scaling corresponds here to the left
 action of the diagonal element $(1/t, t)$. 
  From this one can deduce an isomorphism 
 of graded vector spaces with $\SL_2$-action
 $$\kk[\mathbb{A}^8]  \simeq k[Q_4] 
 \otimes \kk[U^0 \backslash \SL_2^3],$$
 where $\SL_2$ acts trivially on $Q_4$ in degree $4$.  The factor $k[Q_4]$ with its grading
 matches with $H^*(B\PGL_2)$ above.
 The $V_a \otimes V_b \otimes V_c$-isotypical component of 
  $U^0 \backslash \SL_2^3$   
  is just the $U^0$-invariants on
 $ V_a \otimes V_b \otimes V_c$, 
 which are readily computed
 to be a $\min(a,b,c)+1$-dimensional
 space,  spanned by polynomials of degree
$$a+b+c -2j, 0 \leq j \leq \min(a,b,c),$$
which matches \eqref{chugchug}.

\end{example}

\begin{example}
We take the Whittaker case $X=U\backslash G$ with its affine bundle $\Psi$. 
In this case -- taking \eqref{PlanchXv2} as the definition
as in Remark \ref{Whitdef},   
$$ \Planch_X = \kk,$$ 
i.e., $\Planch_X^{(V)}$ is one-dimensional for $V$ trivial and zero otherwise. 
The computation that $\Planch_X^{(V)}$ is trivial for any nontrivial $V$ 
representation is precisely the geometric version of the Casselman--Shalika formula \cite{FGV}.

\end{example}

\begin{example} \label{moment map}
For $V =\mathbf{1}$, the trivial
representation, we have a distinguished morphism $k \rightarrow \Planch_X^{(\mathbf{1})}$.
This corresponds to the unit of the algebra $\Planch_X$. 
With reference to \eqref{Planch mult def},
it comes from the morphism $\mathcal{T}_{\mathbf{1}}^{X} \rightarrow p^{!} \delta_X$
that comes from the identification 
 $p_! \mathcal{T}_{\mathbf{1}}^X \simeq  \delta_X$.

There is always a map
from $\Hom_{\Gr_G/G_O}(\mathcal{T}_V, \mathcal{T}_{\mathbf{1}})$
to $\Planch_X^{(V)}$ (unnormalized version; for the normalized version add the shifts everywhere);
this encodes the $\mathcal{O}(\fgxv)^{\shear}$-module structure on $\Planch_X$. 
Again, from the point of view 
of \eqref{Planch mult def},
this arises as follows: 
 a homomorphism $f: \mathcal{T}_V \rightarrow \mathcal{T}_{\mathbf{1}}$
on $\Gr_G$ gives
a map $\mathcal{T}_V^X \rightarrow \mathcal{T}_{\mathbf{1}}^X \stackrel{unit}{\rightarrow} p^! \delta_X$,
i.e., a class in $\Planch_X^{(V)}$. 
\end{example}
  
\begin{example} (Product, in the homogeneous case):
Let us spell out what the product looks like in the homogeneous case $X=G/H$
(again, explicating only at the level of cohomology). 

  By \eqref{GXhom} we get  
 $\Planch_X^{(V)} = H^*_{H_O}(\Gr_H,  \mathbb{D} \mathcal{T}_V),$
where the dualization $\mathbb{D}$ takes place on $\mathrm{Gr}_H$
and is normalized to fix the basic object;
equivalently, $H_O$-equivariant homomorphisms from $\mathcal{T}_V$
to the  dualizing sheaf on $\Gr_H$. 
 Here there is a morphism 
\begin{equation} \label{VWmap}  \mathcal{T}_V \star_G \mathcal{T}_W|_{\Gr_H} \rightarrow (\mathcal{T}_V|_{\Gr_H} \star_H \mathcal{T}_W|_{\Gr_H}),\end{equation}
arising from a restriction map in cohomology;
on the left, $\star_G$ denotes convolution on the affine Grassmannian for $G$,
and on the right, $\star_H$ denotes
the same on the affine Grassmannian for $H$. 

Now we dualize on $\Gr_H$ (thus reversing the arrow) and take $H_O$-equivariant cohomology on both sides.
On the left hand side we get $\Planch_X^{(V \otimes W)}$.
On the right hand side  we   first of all get the $H_O$-equivariant cohomology
of 
$(\mathbb{D}    \mathcal{T}_V)|_{\Gr_H} \star_H (\mathbb{D} \mathcal{T}_W)|_{\Gr_H}$
(because the convolution of sheaves involves only
descent along a smooth map and a proper pushforward and thus is compatible with $\mathbb{D}$);
and in turn  we have a convolution product
in $H_O$-equivariant   cohomology
 $H^*(A) \otimes H^* (B) \rightarrow H^* (A \star_H B)$
 which shows that this admits a map from $\Planch_X^{(V)} \otimes \Planch_X^{(W)}$.

 For example, consider the case of $X=\PGL_2/T$, and take $V_n$
 the $n+1$-dimensional irreducible representation of the dual group $\check{G} = \SL_2$; then
 $\mathcal{T}_{1} := \mathcal{T}_{V_1}$ is the constant sheaf $\kk \langle 1 \rangle$ on $\mathbb{P}^1 \subset \mathrm{Gr}_{\PGL_2}$ cf. \S 
 \ref{unsheared Satake example}.
  Label the points of $\mathrm{Gr}_T$ as $x_n$,
 where $x_n$ corresponds to the cocharacter $t \mapsto \mathrm{diag}(t^n, 1)$ and
 write e.g.\  $\kk_n$ for the skyscraper on $\mathrm{Gr}_T$ with value $\kk$ at $x_n$.

First of all,  the restriction of 
$\mathcal{T}_{V_n}$ to $\Gr^T$ is the sum
$(\kk_n \oplus \kk_{n-2} \oplus \dots \oplus \kk_{-n})\langle n \rangle$,
from which we deduce that $\Planch_X^{(V_n)}$
is $n+1$-dimensional, with basis indexed by the points 
$x_i$ with $|i| \leq n$ and $i$ having the same parity as $n$. 
In particular, we get bases $d_0$ for  $ \Planch_X^{(V_0)}$, $e_{\pm 1}$ for 
 $\Planch_{X}^{(V_1)}$,   and $f_{-2}, f_0, f_2$ for $\Planch_X^{(V_2)}$,
 where the index is the same as the indexing of points of $\mathrm{Gr}_T$, 
We claim that the product $\Planch_X^{(V_1)} \otimes \Planch_X^{(V_1)} \rightarrow \Planch_X^{(V_2)}$  is given by
\begin{equation} \label{Tprodv1} e_1^2 = f_2, e_{-1}^2 = f_{-2}, e_1 e_{-1} = f_0 + \xi d_0.\end{equation}
 where  $\xi$ is a generator in $H^2(BT_O)$. 
 To see this, we compute
\begin{eqnarray*}
\mathcal{T}_1|_{T} &=& (\kk_1 \oplus \kk_{-1})\langle 1\rangle, \\ 
 \mathcal{T}_1|_T \star \mathcal{T}_1|_T &=& (\kk_2 \oplus \kk_0 \oplus \kk_0 \oplus \kk_{-2}) \langle 2 \rangle, \\
 \mathcal{T}_0 \oplus \mathcal{T}_2|_T &=&  (k_0) \oplus (\kk_2 \oplus \kk_0 \oplus \kk_{-2})\langle 2 \rangle. \end{eqnarray*}
 The morphism \eqref{VWmap} here, for $V=W=V_1$, gives
 $(\mathcal{T}_0 \oplus \mathcal{T}_2)|_T \rightarrow \mathcal{T}_1|_T \star \mathcal{T}_1|_T$. 
 This gives  a morphism of $\mathcal{T}_O$-equivariant sheaves on the points $\{x_i\}$,
and is an identity over $x_{\pm 2}$.   Over $x_0$,  the map comes from
 the map of $\C^{\times}$-equivariant sheaves on a point: 
\begin{equation} \label{sheafmapGR} \kk_0 \oplus \kk_0 \langle 2\rangle \rightarrow  (\kk_0\oplus \kk_0) \langle 2 \rangle,\end{equation}
 where all sheaves have trivial equivariant structure.
 Explicitly, this is the the restriction map from $\Gm$-equivariant cohomology of $\mathbb{P}^1$
 to $\Gm$-equivariant cohomology of the poles $\pm \infty$,
 where we consider both cohomology groups (via pushforward) as $\Gm$-equivariant sheaves on a point. 
 Thereore, in \eqref{sheafmapGR}, the induced map on $\kk_0 \langle 2 \rangle$ is the diagonal,
 and the induced map on $\kk_0$ is the map given by
 $(\xi, -\xi)$, where $\xi$ is a generator of $H^2(B \Gm)(1)$,
 considered as a map $\kk_0  \rightarrow \kk_0 \langle 2 \rangle$. 
 We recover \eqref{Tprodv1} from this, after taking Verdier dual and passing to cohomology.

Now, let us sketch how \eqref{Tprodv1} matches with the other side of Conjecture
 \ref{Plancherel algebra conjecture}
In this case,  $\check{G}=\SL_2, \check{M}=M_2$,
 the space of $2 \times 2$-matrices with its natural $\check{G}$-action
 by left multiplication, and with $\GGm$ acting by scaling, and the
 conjecture says
 $$ \Planch_X^{(V)} = \Hom(V, \kk[M_2]),$$
 where $M_2$ is the space of $2 \times 2$ matrices,
 with coordinates $(a,b; c,d)$ say.  
  Taking $V=V_1$, the standard representation,  the basis
  $e_{\pm 1}$ for $\Planch_X^{(V_1)}$ corresponds to the
  two embeddings $E_{\pm 1}$
sending the coordinate vectors $u,v$ of $V_1$ to 
$(a,b)$ and to $(c,d)$  respectively. 

Now, we may write $V_1 \otimes V_1 = V_0 \oplus V_2$, where the factors
are spanned by $u \otimes v - v \otimes u$
and $u \otimes u, v \otimes v, u \otimes v + v \otimes u$
respectively. 
We deduce that the $V_0$ component of $E_{\pm 1}^2$ 
equals zero, and the $V_0$ component of $E_1 E_{-1}$
corresponds to the embedding of $V_0$ to $k[M_2]$ given by the determinant;
the determinant on $M_2$ corresponds to  $\xi$ in the notation of \eqref{Tprodv1}. 
On the other hand, $E_1^2, E_1 E_{-1} , E_{-1}^2$
give maps $V_2 \rightarrow k[M_2]$ that 
correspond to $f_2, f_0, f_{-2}$ in the notation of 
  \eqref{Tprodv1}.

 \end{example}

\begin{example} \label{JHJG} (The regular locus and the scheme of regular centralizers): 
  In our prior computations we can replace the role of $\mathcal{T}_V$
  with the constant sheaf on $\Gr_G$; this leads to the original definition of Coulomb branch in \cite{BFN}. 
  The corresponding object in the Satake category, by
  \cite{BezFink}, is the pushforward of the structure sheaf from the Kostant slice. 
   
   Specialize now to the case $X=G/H$. 
The computation of (this analogue of, replacing $\mathcal{T}_V$ by the constant sheaf)  $\Planch_X^{(V)}$
  is well-known: it is the equivariant Borel--Moore homology
  of the affine Grassmannian for $H$, which, as computated
  by Bezrukavikov-Finkelberg, is precisely the ring of functions on 
  {\em group scheme of regular centralizers $J_H$}. 
  Working out the details 
  leads to the following consequence of the conjecture,  which we will study
  in later work:   $$ \mbox{fiber of $\check{M}$ above the Kostant slice for $\fgxv$} \simeq J_H.$$
 
   \end{example}

 \subsection{Noncommutative deformations; the symplectic structure on $\check{M}$} \label{Poisson to loop}
   \label{Poisson structure via loop equivariance} \label{Poisson from loop} 
 In Section~\ref{Poisson local conj} we discussed the canonical deformation of the automorphic category that arises from the structure of loop rotation, and which is required to relate to the Poisson structure on $\Mv$. We now formulate a more precise form of this compatibility in terms of a noncommutative deformation of the Plancherel algebra.
 
 Recall that we are considering the $H^*_{\Gm}(\mbox{pt}, k) \simeq \kk[u]$-linear category 
 $\SHV_u(X_F/G_O)=\SHV(X_F/G_O)^{\Gm}$ of $\Gm$-equivariant sheaves for the loop rotation action of $\Gm$ on $\FF((t))$.  
  The object  is naturally $\Gm$-equivariant
for this action.   We may thereby carry out the definition of $\Planch_X$ in $\Gm$-equivariant
cohomology. 
This produces, in place of a $\kk$-algebra, a $k[u]$-algebra  $\Planch_{X}^{\hbar}$ (with $u$ in cohomological degree $2$), specializing to $\Planch_X$ when $u=0$. 

The rigid structure of $\Mv$ provides a natural candidate for this algebra. Namely, $\Mv$ is built out of (twisted) cotangent bundles and symplectic vector spaces, which have canonical deformation quantizations. These deformation quantizations are {\em filtered}, i.e., their Rees algebras form $\Gm$-equivariant sheaves of algebras over $\AA^1$, and hence may be sheared to define $\kk[u]$-algebras, which we propose to identify with the deformed Plancherel algebra: 

\begin{conjecture} \label{noncommutative poisson conjecture}
Write, as per the construction (cf. \S \ref{hyp data checkM}), 
$\check{M}$ as a Whittaker induction of the symplectic
$\check{G}_X$-representation $S_X$
along a morphism $\check{G}_X \times \SL_2 \rightarrow \check{G}$,
and let $U$ be the unipotent radical of the parabolic associated to the $\Gm \subset \SL_2$. 

 Then $\Planch_X^{\hbar}$ is isomorphic to the sheared Rees algebra of the  
quantum Hamitlonian reduction of
$$ A :=  \mbox{differential operators on $\check{G}$} \otimes \mbox{Weyl algebra of $S_X \oplus \mathfrak{u}/\mathfrak{u}_+$}$$
by the twisted action of $\check{G}_X U$; that is to say, 
we first quotient by the ideal generated by the natural embedding 
$\check{\mathfrak{g}}_X \oplus \mathfrak{u} \rightarrow A$,
and then take invariants by $\check{G}_X U$. 
\end{conjecture}

As usual, a noncommutative deformation gives rise to a Poisson bracket, 
i.e., a bilinear mapping satisfying the Leibniz rule,  as follows: 
 If $A$ is a $\FF$-algebra, and $\tilde{A}$ is a flat deformation of $A$ over $\FF[u]$, 
then, for any $a,b \in A = \tilde{A}/u$, the commutator of any lifts $\tilde{a}, \tilde{b}$
to $\tilde{A}/u^2$  is of the form
$u x$, where the reduction $\bar{x} \in A$ is independent of choices.
In particular, we formulate the following expectation, which
is a corollary to the above conjecture, and plays
for us a more important role:

\begin{conjecture}\label{poisson conjecture}
(A consequence of Conjecture~\ref{noncommutative poisson conjecture}.) 

$\Planch_{X}^{\hbar}$ is  flat over $\FF[u]$; 
the associated Poisson bracket, explained above, 
  is identified with the reference to \eqref{plxom} with the 
  Poisson bracket on $\mathcal{O}(\Mv)^{\shear}$ induced by the given symplectic structure on $\Mv$. 
\end{conjecture}
We will revisit this from the point of view of factorization in \S \ref{spectral-factorization}, 
see in particular Remark \ref{Poissonloop2}. 

\begin{example} (Noncommutative deformation in the Iwasawa--Tate case).  
 Let us compute in the caes of $G=\Gm,X = \mathbb{A}^1$,
already discussed in Example \ref{A1example1};
we will sketch why the ``Weyl algebra''
of Conjecture \ref{poisson conjecture} appears.
\footnote{This example presents a glorious number of confusing possibilities
for sign normalization, and   
the reader should be sceptical of every single sign that follows.}
First of all, 
 $$\Planch^{(triv), \hbar}
= H^*(  B (G_O \rtimes \Gm^{loop}) ) = k[\hbar, \tau].$$
where we understand the first coordinate $\tau \in H^2(BG)$
and the second coordiante $\hbar \in H^2(B \Gm^{loop})$. 
Above, and in what follows, we will be able to harmlessly replace $G_O \rtimes \Gm^{loop}$
by its commutative subgroup $G \times \Gm^{loop}$
because the inclusion of one into the other is a homotopy equivalence. 

More generally, for any character $\lambda \mapsto \lambda^n$
of the dual group, the associated isotypical component
$\Planch^{(n),\hbar}$ is free over this algebra, i.e., 
$\Planch^{(n),\hbar} = k[\hbar, \tau] x_n,$
where $x_n$ is in degree $n$
and (viewed  as internal endomorphisms of the basic object) 
represents a generator for the degree $n$ map
of sheaves
$$x_n:  \kk_O \rightarrow \kk_{t^n O} \langle -n\rangle.$$
Then we have $x_n=x_1^n, x_{-n} = x_{-1}^n$, and
\begin{equation} \label{gilgamesh3} x_1 x_{-1} = \tau, x_{-1} x_1 = \tau-\hbar.\end{equation}
From this we see we can identify the whole algebra with  $k[x,y,\hbar]$ with $xy-yx=\hbar$
 and $xy=\tau$; $x_n$ corresponds to $x^n$ and $x_{-n}$ to $y^n$.

 The key computation \eqref{gilgamesh3} follows from the following analysis in the finite dimensional model:
\begin{itemize}
\item[(i)]   Given an inclusion $Y \rightarrow X$
of a divisor into a variety $X$,
we have a restrction map $\kk_X \rightarrow \kk_Y$
and also a map $\kk_Y \rightarrow \kk_X[2]$
coming from the natural map $\iota_! \iota^! \kk_X \rightarrow \kk_X$. 
The composites $\kk_X \rightarrow \kk_X[2], \kk_Y \rightarrow \kk_Y[2]$,
which represent classes in $H^2(X)$ and $H^2(Y)$ respectively, 
are given by (respectively)  by the fundamental class of $Y$
inside $H^2(X)$, and its restriction to $Y$.
\item[(ii)] The equivariant fundamental class of $\{0\} \hookrightarrow \mathbb{A}^1$
inside $H^2_{\Gm}(\mathbb{A}^1) \simeq H^2(B \Gm) = \kk[[t]]$
is given by $n t$ if $\Gm$ acts by $t \mapsto t^n$ on $\mathbb{A}^1$.
In particular, the fundamental class of $t^{n+1} O$ inside $t^n O$,
computed in $G_O \times \Gm^{loop}$-equivariant cohomology,
is given by $\tau + n \hbar$. 

\end{itemize}
 
\end{example}

   \begin{example}\label{multiplicity-freeness}  (Multiplicity freeness):
We can readily compute the $\check{G}$-invariants on $\Planch_X^{\hbar}$: \index{multiplicity free}
they arise from the loop-equivariant
  cohomology of $X_O/G_O$.
But the map $X_O/G_O \rightarrow X/G$ induces a cohomology isomorphism
where $X/G$ has trivial loop action.
That is to say, 
$\Planch_X^{\hbar}$ is the trivial deformation of $H^*(X/G)$,  
and (therefore, under Conjecture \ref{poisson conjecture}) 
  {\em all 
$\check G$-invariant functions Poisson-commute in $\check M$},
i.e., $\check{M}$ is a ``multiplicity free'' Hamiltonian $G$-space.\footnote{
This Poisson commutativity can also be seen from the point of view of factorization. Namely for a closed embedding $i:Z\to X$ the endomorphisms of the constructible complex $i_*k$ have a natural commutative ($E_\infty$) algebra structure (the ring structure of cochains on $Z$). In our case it follows that the endomorphisms of the basic object form a {\em commutative} ($E_\infty$) algebra compatibly with factorization, i.e., a commutative factorization algebra, whence the induced $P_3$ structure on cohomology is trivial.  }

This is precisely in line with our conjecture: it is condition
\ref{condcoisotropic} in our definition of hyperspherical (\S \ref{ssconditions}). 
 In the theory of automorphic forms, the central role 
of multiplicity one is well-known on the automorphic side,
and on the spectral side has also been observed experimentally.
\end{example}

%% file: local-numerical.tex
\section{The Plancherel algebra and the Plancherel formula for spherical functions} \label{section:numerical Plancherel}
\index{Plancherel algebra}

 In this section we describe our primary evidence for Conjecture \ref{local conjecture} -- via its derivative Conjecture \ref{Plancherel algebra conjecture} --  by explaining the relation to the numerical version. 
  Namely, working in the same setting as the previous section, 
 we will show that the Frobenius trace on Conjecture \ref{Plancherel algebra conjecture} recovers the {\em Plancherel formula} for spherical functions on $X(\ff)$, where $\ff$ is a nonarchimedean local field;
 indeed (somewhat informally)
\begin{center}
\emph{the Plancherel algebra categorifies the Plancherel measure.}
\end{center}
 More precisely, the Plancherel algebra is a $\check{G}$-representation, and we will
 see that we obtain the Plancherel measure essentially by taking its character (in the sense of distributions,
 and weighting by an element of $\GGm$). The Plancherel measure encodes the inner products of all Hecke operators applied to the basic vector, while the Plancherel algebra encodes the Hom spaces between all Hecke functors applied to the basic object. 
 The Plancherel measure for spherical varieties has been explicitly computed in many cases by the second-named author, partly in collaboration with Jonathan Wang \cite{SaSph, SaWang},
and we will see that this computation matches with the Frobenius trace of the local conjecture.

 We will now outline what we prove, although deferring precise
issues of normalization to the later subsections.  

 In general, the Plancherel formula -- for $\ff$ a nonarchimedean local field, $\fo$ its ring of integers,
  and $(G, X)$ defined over $\fo$ with $G$, for simplicity, split reductive -- gives  the decomposition of $L^2(X(\ff))$ as a $G(\ff)$-representation, i.e., 
 $L^2(X(\ff)) \stackrel{\sim}{\rightarrow} \int_{\lambda} \pi_{\lambda} $.
 Restricting to $G(\fo)$-invariants gives a corresponding direct integral representation over
 spherical representations, i.e.\ representations with a $G(\fo)$-fixed vector.  
 We fix the ``basic vector'',  \index{basic vector} \index{$\mathbf{e}$}
 $$ \mathbf{e} = \mbox{characteristic function of $X(\fo)$},$$
 and consider its Plancherel density $\mu$. This is a measure on the set of irreducible unitary, unramified representations of $G$ (that we identify with a subset of $\check A/W$ via the Satake isomorphism),\footnote{In this section, we take the dual group to be defined over $k=\CC$, and use $\check G, \check A$, etc., to denote the corresponding groups of $\CC$-points.} characterized by the property -- for Hecke operators $T_V, T_W$ indexed by representations $V, W \in \mathrm{Rep}(\check{G})$ -- 
\begin{equation} \label{Pldef}  \langle T_V \mathbf{e}, T_W \mathbf{e} \rangle_{L^2(X(\mathfrak{f}))} = \int_{t \in \check{A}/W} \chi_V(t) \overline{  \chi_W(t) } \mu(t).\end{equation}
 where $\chi_V, \chi_W$ are the characters of $V$ and $W$ respectively.

We have already seen the left hand side of \eqref{Pldef} arise geometrically, in \eqref{Frob Inverse Trace} -- 
at least for $W$ trivial, which is all that is needed to characterize $\mu$. Namely, an application of Lemma~\ref{Homlemma} identifies it (see Proposition \ref{Plcat} below) with
 \begin{equation} \label{Pldef2} [ \Hom(T_V \star \delta_X,   \delta_X)^{\vee} ]  =  [ \Hom(V,   \Planch_X )^\vee  ],\end{equation}
where we have used the ``trace of Frobenius'' notation as in the proof of Lemma \ref{Homlemma}, the left side denotes homomorphisms in $\Shv(X_F/G_O)$, and the right hand side denotes homomorphisms of $\check G$-modules;
here $\Planch_X$ is the Plancherel algebra introduced 
 in Definition \ref{Pl def}. 

 Since $\Planch_X$ corresponds, under Conjecture \ref{Plancherel algebra conjecture}, to $\mathcal{O}_{\check{M}}^\shear$
under geometric Satake, a simple computation will extract $\mu$ 
from the local conjecture.

The contents of this section, then, are as follows:
\begin{itemize}
\item \S \ref{numerical Plancherel setting} sets up notation. 
\item \S \ref{Cat To Num} we will prove Proposition \ref{Plcat}, which
computes the Plancherel density for $X(\ff)/G(\fo)$ conditional on Conjecture \ref{Plancherel algebra conjecture}. 
\item   \S \ref{ssPlancherel} we  will review the Plancherel formula for spherical functions on spherical varieties;
in Proposition \ref{VXcomp} it is shown to agree, under conditions on $X$, with the prediction
of Proposition \ref{Plcat}. 
\item  \S \ref{Hecke module structure} goes in a slightly different direction and discusses the {\em algebraic}
(rather than unitary) aspects of the Hecke module structure of functions on $X(\ff)/G(\fo)$. 
\end{itemize}

\subsection{Setup: $X, S_X$ and $V_X$} \label{numerical Plancherel setting} 
We will work now in the following setup. 

Let $\FF_q$ be the finite field with $q=p^?$ elements, $\kk$ be the algebraic closure of $\Q_l$ ($l\ne p$), 
and let $\ff=\FF_q((t))$ be the associated nonarchimedean local field,
with ring of integers $\fo= \FF_q[[t]]$. Much of the numerical part of the discussion
would apply to a general nonarchimedean local field, but this
setting is where we can compare with the local geometric conjecture. For the purposes of comparison, we also need to fix an isomorphism $k \simeq \CC$; in particular, this fixes a choice of $q^\frac{1}{2} \in k$, corresponding to the positive square root in $\CC$.

In this section, we take $G$ to be a split connected reductive group defined over $\FF_q$, 
with dual group $\check{G}$ defined over $\kk$. 
As in the prior sections, we will work with a dual pair
$$ M = T^*(X, \Psi) \mbox{ and } \check{M},$$
where we now suppose that $T^*X$ is a hyperspherical variety over $\FF_q$ (see \S~\ref{hdprings}), and $X$ is ``as split as possible,''
which we now formulate more precisely.

We already had a discussion of ``distinguished split forms'' for the dual of a spherical variety in \S~\ref{CheckMRat}\footnote{cf. also
Expectation \ref{GMdesiderata expectation} for a discussion in the broader setting of hyperspherical varieties} the conditions that we will impose on $X$ here are similar, but not necessarily identical: What we require is that, for every simple spherical root of even sphere type (\S~\ref{subsubsec:colors}), the associated rank-one subquotient $X^\circ P/\mathfrak R(P)$ is isomorphic to $\SO_{2n}\backslash \SO_{2n+1}$ \emph{with $\SO_{2n}$ split}. This, in particular, implies that \emph{all colors of $X$ are defined over $\FF_q$}. (Indeed, if colors of type $T$ -- corresponding to the subquotient $\SO_2\backslash\SO_3$ -- are defined over $\FF_q$, then it is easy to see that all colors are defined over $\FF_q$, given our assumption that $G$ is split.)
 We make these assumptions, in order for the ``simple'' action of the Galois group of $\FF_q$ on the symplectic representation $S_X$ of Definition \ref{simpleaction} to be trivial; however, this is just a simplifying assumption, and the more general Galois actions of that definition, together with the modified ``analytic'' and ``algebraic'' actions on $\check M$, described in \S~\ref{Mshear1}, should produce the correct answer in every case.

 We recall, again, in the notation of \S \ref{dp}, that $$X \simeq S^+\times^{HU} G,$$ where $S^+$ is a representation of the reductive subgroup $H$. Assuming, as in \S~\ref{eigencharacter} and \S~\ref{normalized-local}, that the modular character by which $H(\ff)$ acts on the Haar measure on $S^+(\ff)$ extends to a positive character $\eta: G(\ff) \to \mathbb R_+^\times$, which we fix (hence, here, $\eta$ is the absolute value of what was denoted by the same letter in \S \ref{normalized-local}), the space $X(\ff)$ carries a unique up to scalar $(G(\ff),\eta)$-eigenmeasure with factorization $\int_X f(x) dx= \int_{H\backslash G} \int_{S^+} f(sg) ds dg$, for $ds=$ a Haar measure on $S^+$ and $dg$ an eigenmeasure valued in the appropriate line bundle over $H\backslash G$. We normalize this measure in such a way that the measure of $X(\fo)$ equals
\begin{equation}\label{eq:measureX}\textrm{vol} \ X(\fo) =  \frac{  |X(\mathbb{F}_q)|}{|G(\mathbb{F}_q)|} q^{\dim(G)-\dim(X)},
\end{equation}
and normalize the action of $G(\ff)$ on $L^2(X(\ff))$ to be unitary, as in \eqref{action-G-unitary}.

We recall that $\check{M}$ was constructed in 
 \S \ref{omspherical}
as the Whittaker induction of a triple
\begin{equation}  (\check{G}_X \subset \check{G},  
\mathfrak{sl}_2 \to \check{\mathfrak{g}}, \mbox{$S_X$ a self-dual representation  of $\check{G}_X$})\end{equation}
where the $\sl_2$ and $\check{G}_X$ commute, and defined over $\kk$.  
What is more important for our considerations in this section is the space $V_X$ defined in  \S \ref{NTmov}
as
\[V_X = S_X \oplus [\mathfrak{g}_X^{\perp} \cap \check{\mathfrak{g}}_e]\]

Recall that this space carries an action of the group $\check G_X\times \Gm$
for which the $\Gm$ degrees on $V_X$ are all  positive.
See the discussion after \S \ref{NTmov}, where the $\Gm$ was denoted 
by $\GGm'$ for reasons mentioned there.

\subsection{Categorical to numerical} \label{Cat To Num}

  In this section we will prove  the following. 
  
 \index{$U_X$} \index{$\check{A}_X^{(1)}$}
 \begin{proposition} \label{Plcat} 
 Assume the Plancherel Algebra Conjecture \ref{Plancherel algebra conjecture}. Let $U_X\subset \check G_X$ be a maximal compact subgroup.

 Then, with measures normalized as in \S \ref{numerical Plancherel setting},
  and using the notation recalled above,  the unramified Plancherel measure of \eqref{Pldef} can be written
  as follows, for $f$ a continuous function on $\check{A} \sslash W$: 
\begin{equation} \label{muform2} \mu(f) = 
|W_X|^{-1} \int_{\check{A}_X^{(1)}} f(q^{-\rho_{L(X)}} t)
  \frac{\det(I-t|{\check{\mathfrak g}_X/\check{\mathfrak a}_X})}{\det(I- (t,q^{-\frac{1}{2}})|V_X)} dt, \end{equation}
where $\check{A}_X^{(1)}$ is the maximal compact subgroup of the maximal torus $\check{A}_X \subset \check{G}_X$, 
$W_X$ the Weyl group for $\check{G}_X$, 
  $dt$ is the probability Haar measure and $\rho_{L(X)}$ is as in e.g.\ \eqref{rhorho1}.

  \end{proposition}
   In particular,  $\mu$ is supported on the image of the coset $U_X \cdot q^{-\rho_{L(X)}} \subset \check{G}$ under $\check G\to \check A\sslash W$. (Recall again that the cocharacter $2\rho_{L(X)}$ commutes with $\check G_X$.)
   It admits the following alternate description:
Let $\mu_0$ be  the character  of $U_X \times q^{-1/2} \subset \check G_X \times \Gm$ acting on the symmetric algebra $\Sym V_X$;
that is to say, the function on $U_X \times q^{-1/2}$ whose value at $(u, q^{-1/2})$ is the  sum $\sum \chi_n(u) q^{-n/2}$, with $\chi_n$ the trace
of $u$ acting on the degree $n$ component of $\Sym V_X$;
this defines a smooth measure on the compact manifold $U_X$. Then we have
\begin{equation} \label{measmudef}  \mu = \mbox{pushforward of $\mu_0$ by $\check{G}_X \times \Gm \stackrel{\Id, 2 \rho_{L(X)}}{\longrightarrow} \check{G}$}, \end{equation}

\proof  
$\mu$ is uniquely characterized by \eqref{Pldef},
and is indeed characterized by only the cases when $W$ is trivial. Taking into account the normalized action of $G$, we have 
\[\langle T_V \mathbf{e},   \mathbf{e} \rangle = \int_{X(\fo)} \int_{G(\ff)} \sqrt{\eta(g)} T_V(g) 1_{X(\fo)}(xg) dg dx.\] 
This
was computed in \eqref{Frob Inverse Trace}
to equal
the trace of geometric 
Frobenius on the linear dual of $\Planch_X^{(V)}$. (Note that we have modified the measure on $X(\fo)$ here, to absorb an extra factor from \eqref{Frob Inverse Trace}.)
Assuming Conjecture \ref{Plancherel algebra conjecture}
this is the trace of Frobenius on $\Hom(V, \mathcal{O}_{\Mv}^\shear)^\vee$,
where the relevant shear is the 
analytic sheaf of \S~\ref{analytic-M}, that is to say:
$$ \langle T_V \mathbf{e}, \mathbf{e} \rangle = [\Hom(V, \mathcal{O}_{\Mv}^\shear)^\vee].$$

Now, we compute more explicitly using the isomorphism
  $ \check{M} =  V_X \times^{\check{G}_X} \check{G}$;
 see \S \ref{NTmov}. The representation $\mathcal{O}_{\Mv}^\vee$ -- by
 which we mean the space of $\check{G}$-finite linear functionals on $\mathcal{O}_{\Mv}$,
 i.e., the ``algebraic'' dual -- 
 is thereby induced as an algebraic $\check{G}$-representation from the $\check{G}_X$-representation $\mathcal O_{V_X}^\vee = $ the symmetric algebra on $V_X$.
Ignoring for a moment the $\GGm$-action (and the shear), we have isomorphisms of vector spaces
\[\Hom(V, \mathcal{O}_{\Mv})^\vee = (V\otimes \mathcal{O}_{\Mv}^\vee)^{\check G} = (V\otimes \Sym V_X)^{\check G_X}.\]
 One can easily check that the $\GGm=\Gm$-action on the right hand side corresponds to the combination of the action on $V$ via the cocharacter $2\rho_{L(X)}$, and the action on $V_X$ described in \S~\ref{NTmov} (where $\Gm$ is denoted by $\GGm'$).
 
 A vector in degree $n$ for this $\Gm$ action (i.e., $\lambda \cdot v = \lambda^n v$) will
 correspond to a vector in the sheared space with Frobenius eigenvalue $q^{-n/2}$. 
  Using the fixed isomorphism $\kk \simeq \C$, and taking $U_X = $ a maximal compact subgroup of $\check G_X$, we can employ the Weyl character formula to compute invariants. This gives
\[ \int \chi_V d\mu = \langle T_V \mathbf{e},   \mathbf{e} \rangle = \int_{U_X} \chi_V(q^{-\rho_{L(X)}} g) \tr((g, q^{-\frac{1}{2}}) | \Sym V_X) dg\]
(with $q^{-\frac{1}{2}}$ acting on $\Sym V_X$ via $\GGm'$), which is
\eqref{muform2}.

\qed

 \subsection{Known computations of the Plancherel density} \label{ssPlancherel}
In this section we explain the role of the $\check{G}_X$-representation $V_X$ in
the unramified Plancherel formula for $X$ over a nonarchimedean local field.

 The density $\mu$ was computed in many cases by the second-named author \cite{SaSph} and extended
 in joint work with J. Wang \cite{SaWang}.  These papers are not restricted to the smooth case (indeed, the main objective of \cite{SaWang} was to study the ``IC functions'' of possibly singular spherical varieties), but in the smooth case we can summarize the findings as follows (still, with $G$ split):
 
\begin{quote} \emph{ Suppose that $X$ is a spherical affine $G$-variety such that $T^*X$ is hypespherical.  If $X$ is homogeneous or $\check G_X=\check G$, and under some combinatorial assumptions that are true in every example that we know (see Proposition \ref{VXcomp}), 
the unramified Plancherel measure for $X$ is given by \eqref{muform2}; that is to say,
the ``Frobenius trace of the local conjecture'' is correct. }
\end{quote}

To see this, we must recall the setup of the quoted papers and translate it a form
  where it can be readily compared with \eqref{muform2}.  As we shall
recall here, the papers \cite{SaSph, SaWang}
 prove (under assumptions on the spherical variety) that the Plancherel 
 measure is given as the pushforward, under $\check A_X^{(1)}\ni \chi \mapsto \chi q^{-\rho_{L(X)}} \in \check A \to \check A\sslash W$, 
of a measure of the form
\begin{equation}\label{Plancherel-spherical} 
 d\mu(t) =  |W_X|^{-1} 
  \frac{\det(I-t|{\check{\mathfrak g}_X/\check{\mathfrak a}_X})}{\det(I- (t,q^{-\frac{1}{2}})|V_X')} dt, 
\end{equation}
where $V_X'$ is a graded representation of $\check A_X$, equivalently, a representation of $\check A_X\times \Gm$. The multiset of weights of the representation $V_X'$ is $W_X$-invariant and will be described below. 
It is not manifest
in the general situation of\cite{SaSph, SaWang}
that this is the restriction of $\check G_X$-representation.
Hence, to compare with \eqref{muform2}, we need to recall the conditions under which \eqref{Plancherel-spherical} is proven, and to compare the graded $\check A_X$-representations $V_X|_{\check A_X}$ and $V_X'$.

Let us first mention the conditions on $X$ -- in fact, also on twisted cases $(X,\Psi)$ addressed in the aforementioned papers, in which case the assumptions below apply to the Whittaker induction datum: $X$ is a spherical affine $G$-variety, or ``Whittaker-induced,'' in the sense of \eqref{Whittakerinduction}, from one satisfying the following assumptions: 
\begin{itemize}
 \item $M=T^*X$ is hyperspherical. This, in particular, implies that $X$ is smooth and affine and has no ``roots of type $N$'' (see Proposition \ref{hypersphericalspherical}) -- the latest being a necessary assumption for the validity of \eqref{Plancherel-spherical} in the aforementioned papers.
 \item $X$ is either homogeneous, with the support of every spherical root having the type of a classical group, or $\check G_X=\check G$. Moreover, in the first case, there is a combinatorial condition \cite[Statement 7.1.5]{SaSph} on certain data that will be recalled in \S~\ref{VX-homogeneous} below; under this condition, \eqref{Plancherel-spherical} is \cite[Theorem 9.0.1]{SaSph}. In the second case, it is \cite[(1.11)]{SaWang}. We expect the condition on the support of spherical roots to be removed, once a few more cases of simple spherical varieties are checked along the lines of \cite[\S~6]{SaSph}, and the combinatorial condition \cite[Conjecture 7.1.5]{SaSph} to always be true.
\end{itemize}

\emph{For the remainder of this subsection, we assume the conditions above, and compare the 
$\check A_X$-representations $V_X|_{\check A_X}$ and $V_X'$.} The reader who is not familiar with the aforementioned papers is advised to read the statement of Proposition \ref{VXcomp}, and skip its proof. 

\begin{remark}
 The conditions of \cite[Theorem 9.0.1]{SaSph} are a little bit more permissive than our current assumptions: they allow for spherical roots of even sphere type, where for the associated subquotient $\SO_{2n}\backslash\SO_{2n+1}$ the subgroup $\SO_{2n}$ is not split. The results from that paper motivated the painful definition of the ``simple'' Galois action that we provided for those cases in Definition \ref{simpleaction}. However, we will leave it to the reader to compare with the results of \cite{SaSph} for those cases.
\end{remark}

 \subsubsection{The space $V_X'$} \label{VXprime}

For the proposition that follows, and for a given multiset $A$ in a set $B$ (i.e., function $B\to \mathbb N$), where $B$ carries an action of a group $W$, we understand ``the multiset of $W$-translates of $A$'' to be the smallest multiset that is $W$-stable and contains $A$. In other words, its elements are the $W$-translates of elements of $A$, and the multiplicity of each element is the maximum multiplicity of a $W$-translate in $A$. Many natural questions about multiplicities, in the discussion that follows, are easily resolved by applying Lemma \ref{freecolors} to eliminate multiplicities, and we will not make further comments on those.

\begin{proposition}\label{VXcomp}
 Assume the conditions above for $X$, so that the Plancherel density is given by \eqref{Plancherel-spherical}, for a graded $\check A_X$-representation $V_X'$ described in  \cite[Theorem 9.0.1]{SaSph} or \cite[(1.11)]{SaWang}. The space $V_X'$ admits a decomposition $V_X' = S_X' \oplus V_X''$, which compares to the decomposition $V_X = S_X \oplus [\mathfrak{g}_X^{\perp} \cap \check{\mathfrak{g}}_e]$ of \eqref{vx} as follows:
 \begin{itemize}
  \item When $X$ is affine homogeneous (or Whittaker-induced from such), the multiset of weights of the space $S_X'$ is the multiset of $W_X$-translates of valuations associated to colors of even sphere type, with grading 1; in particular, $S_X'\subset S_X$. If those weights are minuscule, we have $S_X'=S_X$.
  \item When $X$ is the affine closure of its open $G$-orbit, and $\check G_X=\check G$,\footnote{The first two cases intersect for affine homogeneous spherical varieties, with $\check G_X=\check G$. In those cases, the corrected version of \cite[Corollary 7.3.4]{SaWang} implies that the weights of $S_X$ and the $W_X$-translates of colors coincide (without counting multiplicities, except for the highest weights). Thus, our two statements about $S_X'$ agree in those cases.} the multiset of weights of the space $S_X'$ contains, with the same multiplicities and grading 1, the multiset of $W_X$-translates of valuations associated to colors of type $T$, and, ignoring multiplicities, coincides with the set of weights of $S_X$. Again, if those weights are minuscule, we have $S_X'=S_X$.
  \item In the general case (with $\check G_X=\check G$), both $S_X$ and $S_X'$ are obtained from the corresponding spaces $S_Y$, $S_Y'$ associated to the affine closures of their open $G$-orbits by the recipe of Definition \ref{def:SX2}.In particular, if $S_Y'=S_Y$ then $S_X'=S_X$.
  \item Finally, the graded $\check A_X$-representation $V_X''$ can be identified with the smallest $W_X$-invariant subspace of $\mathfrak{g}_X^{\perp} \cap \check{\mathfrak{g}}_e$ which contains its zero-weight subspace, as well as its intersection with the span of the support of every spherical root. (Here we identify $\check{\mathfrak g}$ with $\check{\mathfrak g}^*$.)
\end{itemize}

\end{proposition}

We clarify the meaning of the last condition: Every (simple) spherical root $\gamma$ can be written as a sum of simple roots of $G$; those appearing nontrivially in the sum form its \emph{support}. The support of each spherical root defines a standard Levi subgroup of $\check G$, and since the dual group $\check G_X$ is defined, as a subgroup of $\check G$, uniquely at least up to conjugation by $\check A$, its intersection with that Levi is well-defined up to conjugation by $\check A$. Choosing an invariant bilinear form to identify $\check{\mathfrak g}$ with its dual, the image of that Levi subalgebra in $\check{\mathfrak g}$ does not depend on the choice of form, and its intersection with $\mathfrak{g}_X^{\perp} \cap \check{\mathfrak{g}}_e$ gives rise to a sub-$\check A_X$-representation of $\mathfrak{g}_X^{\perp} \cap \check{\mathfrak{g}}_e$. The smallest subspace invariant under the action of the normalizer, in $\check G_X$, of  $\check A_X$ is what is meant by ``smallest $W_X$-invariant subspace'' in the last item.

\begin{remark}
In all examples of smooth affine spherical varieties that we know, the weights associated to colors are minuscule (hence, $S_X'=S_X$), and $V_X''=\mathfrak{g}_X^{\perp} \cap \check{\mathfrak{g}}_e$.  Of course, we expect that $V_X=V_X'$, always. 
\end{remark}

\subsubsection{Recollection of \cite{SaSph} and proof of Proposition \ref{VXcomp} in the homogeneous case.}\label{VX-homogeneous}

We start with the unramified Plancherel formula of \cite[Theorem 9.0.1]{SaSph}, for the cases of homogeneous affine spherical varieties satisfying the conditions recalled in the beginning of \S~\ref{ssPlancherel}. 
To write the formula given in \cite{SaSph} in the form \eqref{Plancherel-spherical}, one needs to take a number of steps, which we detail here.

 First of all, one needs to modify the measures; since we are talking about the case when $X=H\backslash G$ with $H$ reductive (or a Whittaker-induction of such), we will be working with $G$-invariant measures. The normalization \eqref{eq:measureX} means that the measure on $X$ is Tamagawa measure divided by the factor $q^{-\dim G} |G(\mathbb{F}_q)|$, where by ``Tamagawa measure'' we mean the measure obtained by a $G$-invariant, residually nonvanishing, integral volume form on $X$ (and the probability Haar measure on $\fo$).

 \begin{remark} \label{remark-point}
For usage later in the proof, we will reformulate this volume when $X$ is a point. In that case, $\vol(X(\fo)) = \frac{  q^{\dim(G)}}{|G(\mathbb{F}_q)|}$, which is the reciprocal Tamagawa measure of $G(\fo$). The formula of Steinberg \cite{SteinbergEnd, Gross-motive}  states that this is equal to $\det(I - \Fr |V_X)^{-1}$, where $\Fr$ is the (geometric) Frobenius element, and $V_X$ is the Galois representation 
 \[ V_X = \bigoplus_i (T_0 \cc)_i (i),\]
 where $\cc = \check{\mathfrak t}^*\sslash W$, $T_0$ is the tangent space at the image of $0\in \check{\mathfrak t}^*$, the index $i$ denotes its $i$-th graded piece by of the natural $\Gm$-action descending from the action on $\mathfrak t$, and $(d)$ denotes the $d$-th cyclotomic twist (multiplying the action of $\Fr$ by $q^{-d}$). 
 
 In our setting, $\cc$ is interpreted as the quotient $\check M/\check G$, where $\check M$ is the Whittaker cotangent space for the dual group, $\check M = (f+\check{\mathfrak u}^\perp) \times^{\check N} \check G$. By the Kostant section, $\cc$ can also be identified with $\check{\mathfrak g}_e$, the centralizer of a principal nilpotent $e$, and the double of the above grading is obtained by the action of the element $h$ of a corresponding principal $\mathfrak{sl}_2$-triple \emph{plus 2}:
\begin{equation} \label{trivcase}  \vol(X(\fo)) = \frac{  q^{\dim(G)}}{|G(\mathbb{F}_q)|} =  \det(1- \Fr| \check{\mathfrak{g}}_e)^{-1}.\end{equation}
\end{remark}
 
 We now return to our main concern. 
 On the other hand, the formula of \cite[Theorem 9.0.1]{SaSph} is using on $X$ $(1-q^{-1})^{-\dim A_X}$ times the Tamagawa measure (see the remark at the end of Section 9 there). In particular, writing $\vol'$ for the volumes of that paper and $\vol$ for our current normalization, and taking $\Phi = \mathbf e$, equation (9.2) of \cite{SaSph} reads:
 \[ \vol'(X(\fo)) = \frac{1}{Q |W_X|} \int_{\check A_X^1} \vol'(X(\fo))^2 L_X(\chi) d\chi,\]
 hence, in our current normalization,
 \[ \Vert \mathbf e\Vert^2 = \vol(X(\fo)) = \frac{1}{Q |W_X|} \int_{\check A_X^1} \vol'(X(\fo)) \vol(X(\fo)) L_X(\chi) d\chi = \]
 \[ = q^{\dim G} \frac{(1-q^{-1})^{\dim A_X}}{|G(\mathbb{F}_q)| }  \frac{1}{Q |W_X|} \int_{\check A_X^1} \vol'(X(\fo))^2 L_X(\chi) d\chi,\]
 Thus the Plancherel density of $\mathbf e$ has the form 
 \[ q^{\dim G} \frac{(1-q^{-1})^{\dim A_X}}{|G(\mathbb{F}_q)| }  \frac{1}{Q |W_X|}  \vol'(X(\fo))^2 L_X(\chi),\]
 (where we omit the implied probability Haar measure $d\chi$ on $\check A_X^1$).
 Now, we combine the definition of $L_X$ given in \emph{op.\ cit}.\ 7.2.3 with Theorem 9.0.3 which computes $\vol'(X(\fo))$. The former has the form
 \[ L_X(\chi) = c^2 \frac{\det(1-t|_{\check{\mathfrak g}_X/\check{\mathfrak a}_X})}{\prod_{\theta\in \Theta} (1-q^{-r_\theta} e^\theta)}\]
 (we will comment on the set $\Theta$ below, but we mention that the simplifying assumption that all colors are defined over $F$ means that we can ignore the signs $\sigma_\theta$ of \emph{loc.cit.}), and the latter has the form
 \[ \vol'(X(\fo)) = Q c^{-1}.\]
 This simplifies the Plancherel density of $\mathbf e$ to  
 \begin{equation}\label{Planchereldensity-intermediate}  \frac{q^{\dim G} Q\cdot (1-q^{-1})^{\dim A_X}}{ |G(\mathbb{F}_q) | } \cdot \frac{\det(1-t|_{\check{\mathfrak g}_X/\check{\mathfrak a}_X})}{|W_X| \prod_{\theta\in \Theta} (1-q^{-r_\theta} e^\theta)}.
 \end{equation}
 The constant $Q$ depends on the parabolic $P(X)$, and is equal to 
 \[Q=\frac{\vol(G(\fo))}{\vol(P(X)^-(\fo) P(X)(\fo))},\] where $P(X)^-$ is opposite to $P(X)$ (and the  volume  appearing above explicated
 e.g.\ in the statement of \emph{op.\ cit.}  9.0.3.) 
  One then computes the first factor (the first fraction) of \eqref{Planchereldensity-intermediate} is equal to 
 \[\frac{(1-q^{-1})^{\dim A_X}}{\vol L(X)(\fo)} = \frac{1}{\vol \ker(L(X)\to A_X)(\fo)},\]
 where the volumes are taken with respect to Plancherel measure. 
 
 Hence, we can write the Plancherel density of $\mathbf e$ as
 \[ \frac{1}{|W_X|}  \cdot \frac{\det(1-t|_{\check{\mathfrak g}_X/\check{\mathfrak a}_X})}{\vol \ker(L(X)\to A_X)(\fo)\prod_{\theta\in \Theta} (1-q^{-r_\theta} e^\theta)}.\]

 Note that the kernel of the map $L(X)\to A_X$ is connected (by condition \ref{condtorsion} on hyperspherical varieties, \S~\ref{condALL}
 and examining the proof of Proposition
\ref{hypersphericalspherical} ) and its dual is the cokernel of the map $\check A_X \to \check L(X)$. By Remark \ref{remark-point},
and in particular \eqref{trivcase},  $\vol \ker(L(X)\to A_X)$ is the (alternating) trace of $q^{-\frac{1}{2}} \in \Gm$ acting on the exterior algebra of
   \[ \check{\mathfrak l}(X)^e / \check{\mathfrak a}_X,\]
  where $\check{\mathfrak l}(X)^e$ denotes the part of the Levi algebra $\check{\mathfrak l}(X)$ dual to $P(X)$ annihilated by the element $e$ of the $\mathfrak{sl}_2$-triple; equivalently, we could have written this as $\check{\mathfrak l}(X)_e$, the centralizer of $e$ in $\check{\mathfrak l}(X)$.
   This part is graded by the action of the element $h$ of the $\mathfrak{sl}_2$-triple \emph{plus 2}. In particular, its $2$-graded piece is $\mathfrak z(\check{\mathfrak l}(X)) / \check{\mathfrak a}_X$, and its $>2$-graded piece is $[\check{\mathfrak l}(X),\check{\mathfrak l}(X)]^e$.

  Hence, we arrive at the following Plancherel density for $\mathbf e$.
  \begin{equation}\label{Planchereldensity}
  \frac{1}{|W_X|}  \cdot \frac{\det(I-t|_{\check{\mathfrak g}_X/\check{\mathfrak a}_X})}{\det(I- (q^{-\frac{1}{2}}) | \check{\mathfrak l}(X)^e / \check{\mathfrak a}_X) \cdot \prod_{\theta\in \Theta} (1-q^{-r_\theta} e^\theta)}.
 \end{equation}

 The denominator is what is denoted by $\det(I- (t,q^{-\frac{1}{2}})|V_X')$ in \eqref{Plancherel-spherical}. The product over $\Theta$ expresses the nonzero $\check A_X$-weights of $V_X'$, while the factor $\det(I- (q^{-\frac{1}{2}}) | \check{\mathfrak l}(X)^e / \check{\mathfrak a}_X)$ contains its zero weights. To arrive at the claim of Proposition \ref{VXcomp}, we need to recall the definition of the pairs $(r_\theta,\theta)$. The reader is warned that there is nothing pleasant about their description. 
 
 The set $\Theta$ of \cite{SaSph} consists of pairs $(r_\theta, \theta)$ (although, for notational simplicity, we write $\theta\in\Theta$), comprised of a half-integer $r_\theta$ and an element $\theta$ of the cocharacter lattice of $A_X$. For our purposes, $2r_\theta$ should be interpreted as the $\GGm$-weight. This set is $W_X$-stable, and stable under multiplication of the $\theta$'s by $\pm 1$, hence corresponds to a graded representation of the dual torus $\check A_X$, where each graded piece is self-dual (as an ungraded representation), and isomorphic to its $W_X$-conjugates. Its multiset of weights is the smallest $W_X$- and $(\pm 1)$-invariant multiset of pairs $(r_\theta, \theta)$ containing the ``virtual colors'' introduced in \cite[\S~7.1]{SaSph}; we will repeat the definition, for convenience of the reader. Note that in \cite{SaSph} there was a third piece of data for virtual colors -- a sign -- that here we may ignore, because we are for simplicity working with the ``split form'' of the space. 
 
The virtual colors are the smallest multiset $\mathcal D_v$ of pairs $(r_\theta, \theta)$  that 
\begin{itemize}
 \item contains the pairs defined by colors, as follows: if $D$ is a color, then $\theta=\check v_D \in X_*(A_X)$, the valuation defined by $D$, and \footnote{The formula for the grading $2r_D$ of a color needs to be modified in the non-homogeneous case; see Remark \ref{remark-colors-grading}.}
  \begin{equation}\label{grading-color}
  2r_D = \left<\check v_D,  2\rho - 2\rho_{L(X)}\right>;
 \end{equation} 
 in the cases of Whittaker induction, in the notation of \eqref{Whittakerinduction}, \emph{this applies only to the colors induced from $X_L$} -- the rest of the $B$-stable divisors on $X$ are ignored;\footnote{This point is not stated in \cite{SaSph}, making the definition of relevant colors imprecise in the cases of Whittaker induction; however, it readily follows from the arguments that this is the correct definition. }
 \item 
 if $\gamma$ is a root of $G$ which is also a simple coroot of $\check G_X$ (i.e., $\gamma$ is a \emph{spherical root} which also happens to be a root of $G$), and $(r_\theta, \theta) \in \mathcal D_v$ with $\left<\theta,\gamma\right>>0$, then there is a distinct $(r_{\theta'},\theta') \in \mathcal D_v$ with $\theta' = -w_\gamma \theta$, where $w_\gamma$ is the simple reflection associated to $\gamma$, and 
  \begin{equation} \label{grading-color-virtual} 
 2r_{\theta'} = \left<2\check\rho,\gamma\right> - \left<\theta, 2\rho - 2\rho_{L(X)}\right>.
 \end{equation}
 This pair $(r_{\theta'},\theta')$ may or may not come from another color; if not, it is added to the multiset artificially, hence ``virtual'' colors. We observe that the grading of all virtual colors is positive, as follows (at least) by the case-by-case analysis of rank-one and rank-two spherical varieties in \cite{SaSph}. 
 
 \end{itemize}

 We should specifically discuss the case of colors $D$ of type $T$ (see \S~\ref{subsubsec:colors}): If $D, D'$ are two colors of type $T$ and $\alpha$ is a simple root such that $D, D'\subset X^\circ P_\alpha$ (recall that $X^\circ$ denotes the open Borel orbit), then $(r_\theta, \theta)= (\frac{1}{2}, \check v_D)$ and $(r_{\theta'},\theta') = (\frac{1}{2},\check v_{D'})$. This follows from \eqref{coloreq} and \eqref{colorparity}. (Note that in the affine homogeneous case, $\eta$ can be taken to be trivial.)

  \begin{example} We revisit Example \ref{ex:5sphere}. Here
there is only one color $D$ with $\check v_D = \frac{\check\gamma}{2}$, so $r_D = \frac{3}{2}$, but we also have a virtual color $D'$ with $\check v_{D'} = \check v_D $ and $r_{D'} = \frac{1}{2}$. The representation $V_X$ is two copies of the standard representation of $\SL_2$, one with grading $1$ and the other with grading $3$.
\end{example}

 We can now divide the (multi)set $\Theta$ into a disjoint union $\Theta_1 \sqcup\Theta_2$, where $\Theta_1$ consists of those pairs $(r_\theta, \theta)$ with $r_\theta = \frac{1}{2}$. We define $S_X'$ to be the graded representation of $\check A_X$ with weights in $\Theta_1$, and $V_X''$ to be the direct sum of $\check{\mathfrak l}(X)^e / \check{\mathfrak a}_X$ with the one with weights in $\Theta_2$. 
 
 \begin{lemma}
  $\Theta_1$ consists precisely of the $W_X$-translates of
  \begin{itemize}
   \item colors of type $T$; 
   \item virtual, but not actual, colors associated to spherical roots of even sphere type $\SO_{2n}\backslash\SO_{2n+1}$ with $n\ge 2$. 
  \end{itemize}
 \end{lemma}

Note that, here, we don't need to additionally say ``$(\pm 1)$-translates,'' because the set of $W_X$-translates as in the lemma is automatically closed under $(\pm 1)$, by \eqref{coloreq}, \eqref{colorSeq}.

\begin{proof}
 This is by inspection of the cases of rank-one spherical varieties listed in \cite[Theorem 6.11.1]{SaSph}, and the subsequent calculations in that paper.
\end{proof}

Hence, the representation $S_X'$ satisfies the statement of Proposition \ref{VXcomp}. To prove the statement on $V_X''$, we check the zero- and nonzero-weight spaces for $\check A_X$ separately. The zero weight space clearly coincides with the zero weight space of $\mathfrak{g}_X^{\perp} \cap \check{\mathfrak{g}}_e$, since $\check{\mathfrak l}(X)$ is the centralizer of $\check{\mathfrak a}_X$ in $\check{\mathfrak g}$.
 
 For the nonzero-weight spaces, we have unfortunately been unable to identify them with those of $\mathfrak{g}_X^{\perp} \cap \check{\mathfrak{g}}_e$, but one can check the statement of Proposition \ref{VXcomp} about the support of  spherical roots, as follows: Let $\gamma$ be a (simple) spherical root, and $P_\gamma$ the parabolic defined by the support of its spherical roots. One can then consider the variety $X^\circ P_\gamma$, and its quotient $X_\gamma$ by the unipotent radical of $P_\gamma$, which is a homogeneous spherical variety for the Levi $L_\gamma$. Colors of $X$ contained in $X^\circ P_\gamma$ are in bijection with colors of $X_\gamma$, and one can easily check by hand, for each of the cases of such simple roots (appearing in Section 6 of \cite{SaSph}) that $V_{X_\gamma} = V_{X_\gamma}'$. In particular, $V_{X_\gamma}'' = \mathfrak l_\gamma\cap  \mathfrak{g}_X^{\perp} \cap \check{\mathfrak{g}}_e$, and therefore $V_X''$ can be identified with the smallest $W_X$-invariant $\check A_X$-subrepresentation of $\mathfrak{g}_X^{\perp} \cap \check{\mathfrak{g}}_e$ which contains its intersections with $\mathfrak l_\gamma$, for every simple spherical root $\gamma$.

\subsubsection{Recollection of \cite{SaWang} and proof of Proposition \ref{VXcomp} in the non-homogeneous case.}

In the cases considered in \cite{SaWang}, we have $\check G=\check G_X$, and therefore (since this does not include twisted -- i.e., Whittaker-induced -- cases) every simple root is of ``type $T$,'' and $\check M = S_X$. In that case, the formula \cite[(1.11)]{SaWang} on the Plancherel density of the basic function is of the form \eqref{Plancherel-spherical}, with $V_X'=S_X'$ described by a multiset $\mathfrak B_X$ of weights -- a crystal -- with the following properties (see \cite[Theorem 7.1.9]{SaWang}): 
\begin{enumerate}
 \item It is obtained from the corresponding crystal of the affine closure of the open $G$-orbit, according to the recipe of Definition \ref{def:SX2}. We are employing, here, Proposition \ref{Gvaluations}, to identify the set $\mathcal D^G_{\text{sat}}(X)$ of \emph{op.cit}.\ with the set $\mathcal D^G(X)$ of $G$-invariant divisors.
 \item In the case of $X = \overline{X^\bullet}^\aff$, without counting multiplicities, the weights in $\mathfrak B_X$ coincide with the weights of $S_X$; moreover, $W_X$-translates of colors appear in $\mathfrak B_X$ with the same multiplicity as the corresponding colors. (This is unambiguous, because Lemma \ref{freecolors} reduces us to the case where all color multiplicities are $1$.)
\end{enumerate}

The $\GGm$-grading of the space $S_X'$ is equal to $1$.

 \begin{remark}\label{remark-colors-grading}
To compromise the recipes for the grading, here and in the previous case, we remark that the formulas \eqref{grading-color} and \eqref{grading-color-virtual}, in the non-homogeneous case where there is an eigenmeasure with eigencharacter $\eta$, need to be shifted by the character $\eta$ of \eqref{eigencharacter}; that is, \eqref{grading-color} should become  
 \begin{equation}\label{grading-color-eta}
  2r_D = \left<\check v_D,  \eta + 2\rho - 2\rho_{L(X)}\right>,
 \end{equation}
and \eqref{grading-color-virtual} should become 
 \begin{equation} \label{grading-color-virtual-eta} 
 2r_{D'} = \left<2\check\rho,\gamma\right> - \left<\check v_D, \eta + 2\rho - 2\rho_{L(X)}\right>.
 \end{equation}

Note that, by \eqref{rhorho1}, $\eta + 2\rho - 2\rho_{L(X)}$ is a character of $A_X$, hence these definitions make sense.  
The fact that \eqref{grading-color-eta} gives $2r_D=1$ when $\check G_X=\check G$ follows, as in the homogeneous case, from \eqref{coloreq} and \eqref{colorparity}.

\end{remark}

\subsection{Questions about the Hecke module structure of spherical functions} \label{Hecke module structure}
 
  The Plancherel formula describes the structure of $L^2(X_{\ff})^{G_{\fo}}$. One may also be interested in more algebraic versions of the same question:
 in particular, the module structure of the space of functions $X_{\ff} \rightarrow \kk$ (with $\kk$ a coefficient ring)
 that are invariant by $G_{\fo}$; this question becomes particularly interesting for $k$ of finite characteristic.
We now sketch  what our conjecture suggests about this, assuming throughout that 
   the residue characteristic should be invertible in $\kk$.

One hopes that a suitable ``trace of Frobenius'' on the category of sheavees,
 recovers the vector space of $G_{\fo}$-invariant functions on $X_{\ff}$.  This is studied in the group case in \cite{xinwentrace};
 we do not know it is true for $X_{\ff}/G_{\fo}$. Nonetheless, it is reasonable to believe it is so.
 
 On the spectral side,  the trace of an automorphism   $F$
 on a category of coherent sheaves is the space of functions on the fixed locus of $F$;
 or, more abstractly,  the (derived) intersection of the diagonal  and the graph of $F$.
In the case of an automorphism of $\check{M}$ commuting with $\check{G}$,  
 the fixed locus of $F$ on $\check{M}/\check{G}$ is then
  the twisted inertia stack:
$    \{ (g \in \check{G}, m\in \check{M}): gm=Fm\}/\check{G}$.  
This reasoning suggests, then, that 
\begin{equation} \label{CXk}   C^{\infty}_c(X_{\ff}; \kk)^{G_{\fo}} \stackrel{?}{=} 
\left( \kk \left[ g \in \check{G}, m \in \check{M}: (g, q^{-1/2}) \cdot m = m  \right]^{\check{G}}  \right)^{\shear}
\end{equation}

On both sides we take invariants in the {\em derived} sense -- that is to say,
each $G_{\fo}$-orbit $G_{\fo} x \subset X_{\ff}$ contributes not a copy of $\kk$,
but of the cohomology $H^*(G_{\fo, x}, \kk)$ of the stabilizer.  
The right hand side should  also be interpreted in a derived sense. 

If $\kk = \C$, these derived phenomena are irrelevant, and
we get simply the usual function space of $\kk$-valued $G_{\fo}$-invariant functions on $X_{\ff}$. 
Then \eqref{CXk}  follows from the work \cite{SaSph} of the second-named author
in the situations to which that work is applicable. 
Indeed, we can rewrite the right hand side as
$ \kk \left[ g \in \check{G}_X, v \in V_X: g \cdot v = q^{1/2} v \right]^{\check{G}_X}$.
Now  for any $g \in \check{G}_X$ the space of solutions to $g v = q^{1/2} v$
 corresponds to the $1/\sqrt{q}$-eigenspace for $g$, and a $g$-invariant function on this space 
 is clearly constant. Therefore, the ring on the right is simply the ring of class functions of $\GvX$, and one can apply the results of \cite{SaSph}   to deduce
 the desired isomorphism of Hecke modules (for $\kk=\C$)
 $$C^{\infty}_c(X_F/G_O)\simeq \C[\GvX]^{\GvX}$$  In particular, for $\kk=\C$, the module structure is insensitive to $V_X$, but this is  presumably false for general $\kk$. 
  This is an interesting
 point to study further.

To summarize parts of our previous discussion:
while, as we just noted,  $V_X$ is not reflected in the ``mere'' module structure of $C^{\infty}_c(X_F/G_O)$,
it becomes visible if we consider also the basic function and inner product.  
Indeed, our discussion implies
that, when
Proposition \ref{Plcat}  applies,
 there is an identification of Hecke modules equipped with distinguished vector and inner product:
 $$(C^{\infty}_c(X_{\mathfrak{f}}/G_{\mathfrak{o}}) \ni \mathbf{e}, \langle -, - \rangle) \simeq 
 (\C[\GvX]^{\GvX} \ni  1, \langle - , - \rangle_{L^2(\mu)} )$$
 where $\mu$ is the Haar measure on a translate of the compact form of $\GvX$, multiplied by    a $q$-deformation of the character of the representation $V_X$.

%% file: global-geometric.tex
  In this Part we study the global story -- that is, starting from a 
 dual pair $(G, M)$ and $(\check{G}, \check{M})$
  of hyperspherical varieties, we examine
  the matching of  associated automorphic and spectral quantizations
in the setting of global Langlands -- more specifically, for 
a curve over either a finite field or the complex numbers.

We refer to page \pageref{part3intropage} for a summary of the contents of the various subsections.
To briefly reprise: we begin in \S \ref{section:global-geometric} by examining automorphic
quantization in global geometric Langlands,
giving ``period sheaves'' and ``period functions.''
   The spectral side of the story will be entirely parallel, 
  giving what we call the $L$-sheaf; it will 
  be treated in \S \ref{Lsheaf};
and in \S \ref{GGC} we will put the two together and  formulate the geometric form of the conjecture, which asserts
  that period and $L$-sheaves match under geometric Langlands. 
  After carrying out some sanity checks for the case of $\mathbb{P}^1$
  in  \S \ref{P1}, we then turn in \S \ref{L2numerical}
  to the arithmetic manifestation of the same phenomena -- that is, the corresponding
  statements concerning equalities of numerical periods with $L$-values.

\section{Period functions and period sheaves}
\label{section:global-geometric}

According to the general picture explained in \S \ref{AFT2}, 
periods should give 
{\em objects} in the categories on the two sides (automorphic
and spectral) of the global geometric Langlands conjecture,
just as they give {\em vectors} in the vector space of automorphic functions
in global Langlands and {\em objects} of the category of local representations
in local Langlands. 

In the current section we shall describe the global automorphic objects and the global automorphic vectors arising from a   
polarized hyperspherical $G$-variety $M$ - we will be more specific about the setting in a moment.  These will be called ``period sheaves'' and ``period functions.''
For example, when $M=T^*X$ (that is to say, there is no twist in the polarization),  we get a morphism
$$\pi: X/G \longrightarrow BG, $$
By taking maps from an algebraic curve $\Sigma$, and then pushing forward
the constant sheaf   along the maps induced by $\pi$, we obtain 
  the ``period sheaf''  on the space of $G$-bundles.
  Frobenius trace extracts the usual period function.

The contents of the  present section are as follows:
\begin{itemize}
\item \S \ref{notn} sets up the considerable amount of notation needed.

\item \S \ref{bunX} introduces the space of bundles with an $X$-section, together
with the necessary twists needed for our conjecture. 

\item \S \ref{periodX} defines the basic period sheaf and period function. 

\item \S \ref{normalizedperiod} introduces a crucial ``unitary'' (analytic) normalization
of the period sheaf.  The normalization process is a half-twisting that depends on the $\GGm$-action.

\item \S \ref{Affcase} describes the modifications needed when $M$ admits a twisted polarization (\S \ref{dp}). 

\item 
\S \ref{BunXex} gives a number of examples, emphasizing
the roles of the various twists.

\item \S \ref{spindep} discusses how the statements
can be reformulated to be manifestly independent of spin structure. 

\item \S \ref{Pinduction} explains the compatibility of period sheaves with Whittaker induction of Hamiltonian spaces, reducing their study to the vectorial case. 

\item \S \ref{Pindep}
explains why the construction is independent of the choice of polarization --
an incarnation of the Fourier transform.   (Therefore, 
the period sheaf can indeed be viewed as attached to $M$, rather than $X$; we have however preferred
to denote it as $\mathcal{P}_X$ rather than $\mathcal{P}_M$).

 \end{itemize}

\subsection{Notation: $\Sigma$, $G$ and $M$} \label{notn}

\subsubsection{Coefficient fields}

We work over an algebraically closed base field $\FF$, which is either
the algebraic closure of $\FF_q$, or $\CC$. We fix a smooth projective curve
$\Sigma$ over $\FF$.  \index{$\FF=$ base field of curve.} For any statement that entails a Frobenius morphism (e.g., any mention of Weil structures on sheaves), or adeles, it is understood that $\Sigma$ is defined over $\FF_q$; in those cases ``adeles'' and ``function field'' of the curve will always mean over $\FF_q$.

We also have a ring of coefficients $\kk$,
which we take to be $\overline{\Q_{\ell}}$ or $\C$ according
to whether  $\FF$ is of finite characteristic or complex. \index{$\kk=$ field of coefficients.}
In the finite case, we fix a square root $\sqrt{q} \in \kk$
of the cardinality of $\FF_q$. As in the local case, we will present ``normalized'' or ``analytic'' versions of our conjectures, which make use of this choice, and arithmetic versions, which do not.

\subsubsection{Spin structure on the automorphic side} \label{spinstructure-automorphic} 
\index{$\cK^{1/2}$ spin structure} \index{$\mathfrak{d}$ different} \index{$\mathfrak{d}^{1/2}$ half-different}

It will be convenient to choose a spin structure on $\Sigma$, i.e., a square root $\cK^{1/2}$ of its canonical bundle  (though we will keep track of the dependence of our constructions on the choice and indicate how to formulate statements independent of it, see in particular \S \ref{spindep}.).  For convenience, we will fix a rational section of this square root as well. 
Squaring this gives rise to  \index{spin structure} 
a meromorphic differential form $\omega$ \index{$\diff$=different} \index{$\partial=$ id\`ele associated to different}
whose zero divisor (the ``different'')   \begin{equation} \label{diffdef} \mathfrak{d} := \sum n_v v, \ \ \sum n_v =  (2g-2).\end{equation} 
has all even multiplicities; we put $\mathfrak{d}^{1/2} = \sum \frac{n_v}{2} v$ and
write $\mathcal{K}^{1/2} =\mathcal{O}(\mathfrak{d}^{1/2})$ for the associated line bundle. 
We also write $\partial$ for an id\`ele associated to $\mathfrak{d}$
so that $\partial = (\partial_v)_v$ with $\partial_v = \pi_v^{n_v}$
and $\pi_v$ a local uniformizer.  

\index{$F=$ function field of curve.}
   Let $F$ be the function field of $\Sigma$,
 $\mathbb{A}_F$ its ring of adeles, and $\widehat{\mathfrak{o}} \subset \mathbb{A}_F$
   its integral subring. 
In the finite field setting, $\omega$ gives rise to a homomorphism \index{$\psi=$ adelic additive character}
\begin{equation} \label{psidef}  \psi: \mathbb{A}_F/F \rightarrow \FF_q,\end{equation}
whose restriction to the copy of $F_v$ inside $\mathbb{A}_F$
comes from the pairing $(f, \omega) \mapsto \mathrm{Res}_v(f \omega)$. 
Fixing  once
  and for all an additive character $\psi$ of $\FF_q$
  valued in $\kk^{\times}$, 
   we get a character of $\mathbb{A}_F/F$, also to be denoted by $\psi$. 
   Notice that
   varying the choice of  rational section 
   while fixing the spin structure $\mathcal{K}^{1/2}$  varies $\psi$
by a {\em square} in $F^{\times}$.    
    From our definitions, the character $x \mapsto \psi( x)$ of $\mathbb{A}_F$ 
has the property that on each completion $F_v$ it is trivial on $\partial_v^{-1} \mathfrak{o}_v$
but not on $\varpi_v^{-1} \partial_v^{-1} \mathfrak{o}_v$;
thus $x_v \mapsto \psi(\partial_v^{-1} x_v)$ is an  unramified character of $F_v$. 

\subsubsection{$G, M$, and polarizations}  \label{periodsheafgeneralsetup}

$(G, M = T^*(X,\Psi))$ will be a distinguished split
form over $\FF$,  in the sense of
Definition \ref{dhpFqdef}\footnote{The word ``split'' is relevant
here only in the context where $\FF$ has finite characteristic,  and there because we did not 
define a notion of hyperspherical in finite characteristic.  For the purposes 
of this chapter,  the reader can ignore the word ``split'' entirely, and instead take the data of $(X, \Psi)_{/\FF}$
as a starting point, as in Remark \ref{XnotM}.},  of a hyperspherical variety admitting a distinguished polarization; where we recall $G$ is to act on the right. 
See \S \ref{Pindep} for discussion of independence
of polarization as well as removing the assumption that $M$ is polarized.

We will frequently allow ourselves to assume that $X$ admits an eigenmeasure, see
\S \ref{eigencharacter}. 
This assumption is ``harmless'' for global applications, for reasons outlined in \S \ref{ssseigencharinocuous},
 and 
it should also be possible to formulate the discussion to avoid it entirely, see
preliminary discussion along these lines in Remark \ref{noeigenmeasure2}.  However,
we find it extremely helpful in thinking about how to normalize to maintain this assumption.

We will make use of 
the quantity
\begin{equation} \label{betaXdef} \beta_X = (g-1) (\dim G + \gamma_X - \dim X ). \end{equation} \index{$\beta_X=$normalizing shift}
where $\gamma_X$ is, as in \eqref{gammaXdef}, the character through which $\GGm$ scales the eigenmeasure.  This will appear as
a normalizing shift in our period sheaf;
a corresponding shift will also intervene on the spectral side.

 \begin{remark} \label{XnotM}
Note that, while we make the hyperspherical assumption for global coherence of the paper,
{\em all} the considerations of this chapter can be applied to an arbitrary
$G \times \GGm$-space $X$
with an eigenmeasure,  or for that matter
the same situation allowing a $\mathbb{A}^1$-torsor $\Psi \rightarrow X$, 
and there are certainly
examples where duality theory seems to work well that land outside our hyperspherical framework. A particularly important
example is the ``Eisenstein case'' $X=U\backslash G$ considered as a $T \times G$ space, which we will
at times consider by way of contrast.  
 \end{remark}

\subsubsection{Context for the Langlands program and sheaf theory} \label{BunLocintro}
 \label{aut sheaf notation}

We briefly recall the main outline of the geometric Langlands program and the underlying sheaf theory, see Appendices~\ref{sheaf theory} and~\ref{geometric Langlands} for a more thorough overview.

 Attached to   $\Sigma$ there are two basic spaces of interest for the geometric Langlands story:
\begin{itemize}
\item  The space   $\Bun_G$ of coherent $G$-bundles
on $\Sigma$.  \index{$\Bun_G$} 
\item The space $\Loc_{\check G}$ of  $\check{G}$-bundles on $\Sigma$.  \index{$\Loc_G$}
 \end{itemize}
 
 We have used the word ``space'' loosely; more precisely, $\Bun_G$ is an algebraic stack over $\FF$
 and $\Loc_{\check{G}}$ is a derived algebraic stack over $\kk$. 
   The dimensions of these spaces will often come up, and we abridge them
   in the following way:  \index{$b_G$} \index{$b_H$}
\begin{equation} \label{5bdef} b_G = \dim \Bun_G= (g-1) \dim G, \; b_H = \dim \Bun_H = (g-1) \dim H, \dots\end{equation}
 and so on, where $g$ is the genus of the curve $\Sigma$;
 the dimensions of the corresponding $\Loc$ spaces are obtained by doubling these.

 The geometric Langlands correspondence posits an equivalence
 between a category of ``constructible''-type sheaves on the former space, 
  and a category of ``coherent'' sheaves on the latter. 
  This general vision has
  been formulated in at least three different contexts,
  with varying specifics. 
 $\Bun_G$ is the space of algebraic $G$-bundles in all cases, but the
 category of sheaves on it varies, and the definition of $\Loc_{\check{G}}$ also varies. 
 We give a r{\'e}sum{\'e} of these constructions in Appendix \ref{geometric Langlands};
for now we summarize sheaf theory on the automorphic side.

  In all cases, we will write simply \index{$\Autshv$} 
\begin{equation} \label{automorphic sheaf category definition} 
 \autshv(\Bun_G) \supset \autshvspec(\Bun_G),
  \bigautshv(\Bun_G) \supset \bigautshvspec(\Bun_G),
 \end{equation}
  for the ``small'' and ``big'' categories $\autshv$ or $\bigautshv$ of sheaves on $\Bun_G$, and, inside it, the ``spectrally decomposable'' subcategories $\autshvspec$ or $\bigautshvspec$ i.e.
  the largest category on which it is reasonable to think about Hecke actions.   
  We describe $\autshv(\Bun_G)$ more explicitly in each case,  but refer again to Appendix \ref{geometric Langlands} for details.

  \index{$\autshv(\Bun_G)$}
  
 \index{$\FF$} 
 \begin{itemize}
 \item  Finite context: $\FF=\overline{\FF_q}$ and $\kk = \overline{\mathbb{Q}_{\ell}}$.   
$\autshv(\Bun_G)$ consists of   {\'e}tale constructible  sheaves on $\Bun_G$ with coefficients in $\kk$.
This has a Frobenius action (when $\Sigma$ is defined over $\FF_q$); Frobenius-equivariant objects are
then Weil sheaves, and here we can talk about ``Tate twists.''

 \item Betti context: $\FF = \C$, and $\kk=$  an algebraically closed  field of characteristic zero.\footnote{Betti sheaves, actually, make sense over any coefficient ring}
 $\autshv(\Bun_G)$ consists of (certain) sheaves on $\Bun_G$ with coefficients in $\kk$ and Lagrangian singular support.   In many ways this is technically the simplest setup.

  \item De Rham context: $\FF=\kk=\C$.  Objects of $\autshv(\Bun_G)$
are coherent $D$-modules on $\Bun_G$.

  \end{itemize}

\begin{remark} 
Still another context, 
closely related to the ``finite'' setting,  is to take $\FF=\C$ but to use  constructible sheaves for the classical, instead of the \'etale, topology (and arbitrary coefficients), which also affects the definition of $\Loc_{\check G}$ on the spectral side. We will, somewhat abusively, use the ``\'etale setting'' to describe 
either this setting or the finite setting, since many of the same conclusions will apply in both cases. 
\end{remark}

{\em An important technical point is that we will work in all cases with the {\em ind-finite} or {\em renormalized} category of sheaves; \index{renormalized category}
  see  \S \ref{renormalization section}}
  and Appendix \ref{geometric Langlands} for discussion;
  while this is important in trying to get various categories to line up, the
  reader unfamiliar with this notion will not lose a lot by ignoring these words
  at a first reading (and this choice can be adjusted at the cost of adding ``nilpotent singular support'' requirements to the sheaf theory on the spectral side).

The constructions presented  will have slightly different interpretations according to context.  
We will explain most constructions in the Betti context, where they are particularly simple to describe. 
 When the translations to other contexts are not straightforward we will add a description of them.

\index{angle bracket twists $\langle d \rangle$}
  In each context we will use the twist notation $\mathcal{F} \langle d \rangle$
  following \S \ref{anglebracketnotation},
  and in particular \eqref{ultimate shearing}. Recall 
  that this means three simultaneous twists (which may or may not apply, according to context):
  a cohomological twist $[d]$, a Tate twist $(d/2)$, 
using the fixed choice $\sqrt{q}$ inside our field of coefficients $\kk$,
and a parity twist.

\subsection{The space $\Bun_G^X$ of bundles with an $X$-section} \label{bunX}
\index{$\Bun_G^X$} \index{Bundles with $X$-section}
 We discuss first the case of untwisted cotangent bundles, $M= T^*X$, with the twisted case to be discussed in \S \ref{Affcase}.  The space of primary interest 
 for defining the period sheaf is informally
  $$ \textrm{$\Bun_G^X$} := \mbox{``$G$-bundles  with a section of the associated $X \otimes \cK^{1/2}$-bundle.''}$$
For example, when $G=\GL_n$ and $X$ its standard representation with scaling $\GGm$ action,  the fiber of $\Bun_G^X$ over a vector bundle $V$ is the space of sections of $V\otimes \cK^{1/2}$; when
$X=H\backslash G$ with trivial $\GGm$ action
then $\Bun^X_G \rightarrow \Bun_G$ is identified with $\Bun_H \rightarrow \Bun_G$. 

 Formally, 
  $\Bun_G^X$ -- or $\Bun^X$ when the group $G$ is clear -- is the algebraic stack defined as the pullback of mapping stacks 
\begin{equation} \label{diagdog} \xymatrix{\Bun_G^X\ar[r]\ar[d] &\Map(\Sigma, \frac{X}{G\times \GGm})\ar[d]\\
\Bun_G\ar[r]^-{\mathrm{id} \boxtimes \cK^{1/2}}&\Bun_{G \times \GGm}}.\end{equation}

 This is in fact an Artin stack,  cf. \cite[Theorem 1.1]{Olsson} or see (a) below for a sketch.
 Note that, for the considerations that follow, it does not matter whether we consider
 these as classical or  derived stacks.
 The reason is as follows: Although one can meaningfully enrich $\Bun^X_G$ to a derived stack, this will have no effect on
 constructible sheaves -- in particular the period sheaf -- which are sensitive only to topology. 
 
\begin{remark} \label{alg-geom-tech1}

 \begin{itemize}
\item[(a)]  $\Bun_G^X$ and the  geometry of $\Bun_G^X \rightarrow \Bun_G$ is quite tame, as we now explain:

$\Bun_G$ is a union of open substacks, each of which 
is a global quotient of a scheme; we exhibit the global quotient structure by ``adding level structure,'' i.e.
fixing a point on the curve and trivializing the bundle up to some order at that point, and then descending; the group involved
is thus a pro-unipotent extension of $G$. Upon pullback to each such open substack, the morphism $\Bun_G^X \rightarrow \Bun_G$
is  a global quotient of a morphism of schemes, as we can see by choosing a $G$-equivariant embedding of $X$ into a vector space
and thus reducing to the case when $G=\GL_n, X=\mathbb{A}^n$.

 \item[(b)]  For $X$ homogeneous, $\Bun_G^X$ is smooth (although  its morphism to $\Bun_G$ need  not  be smooth). This is true more generally
     for $X$ smooth on the locus of maps
     which generically land in the open $G$-orbit on $X$. For an example of nonsmoothness see \S \ref{TateBunX}.
     
 \end{itemize}
 \end{remark}

 \subsection{The period sheaf and the period function.} \label{periodX} 
 The compactly supported (i.e., $!$-) pushforward of the constant sheaf along $\Bun_G^X \rightarrow \Bun_G$
will be called the  {\em unnormalized period sheaf} $\mathcal{P}_X$.  
 $$ \mathcal P_X = \mbox{compactly supported pushforward of constants along $\Bun_G^X \rightarrow \Bun_G$.}$$
 We emphasize that this is a $!$ pushforward; a
dual $*$ pushforward
 will appear on some rare and interesting occasions e.g.\ Remark \ref{starperiod} and \S \ref{starperiods}; and in that case we will use $*$ explicitly in the notation.

 In the finite context $\mathcal P_X$ is considered as a Weil sheaf, i.e.,
 with a Frobenius equivariant structure.  We take the trivial Frobenius action on the constant sheaf here, to be ``corrected'' later, when we introduce the normalized period sheaf.

Let us compute the function associated to this Weil sheaf.  

Recall that we can uniformize $\Bun_G(\FF_q)$ as the quotient $G(F) \backslash G(\mathbb{A})/
G(\widehat{\mathfrak{o}})$; the sections of the bundle parameterized
by $g \in G(\mathbb{A})$ are identified with the elements of $x \in G(F)$ with the property that $xg \in G(\widehat{\mathfrak{o}})$.  Thus, for example in the case $G=\Gm$, 
the element of $\Gm(\mathbb{A})$ that is given by the uniformizer $\varpi_x$ at a single
point $x$ of the curve $\Sigma$ parameterizes the line bundle $\mathcal{O}(x)$ (we spell this out to avoid possible sign confusion). 
With reference to these adelic uniformizations the chosen spin structure of \S \ref{diffdef} can be identified with the class of $\partial^{1/2} =\prod_v\varpi_v^{n_v/2}$ in $F^{\times} \backslash \mathbb A_F^\times/\prod_v \mathfrak o_v^\times=
\Bun_{\Gm}(\FF_q)$, 
and the set of $\FF_q$-points in the fiber of $\Bun_G^X$ over 
the $G$-bundle represented by $g \in G(\mathbb{A}_F)$
is identified with 
\begin{equation} \label{bunxGm} X(F) \cap \prod_v X(\mathfrak{o}_v) \cdot( g^{-1}, \partial^{-1/2}).
\end{equation}
The Frobenius trace on the period sheaf recovers the {\em period function} (or theta series)
associated to $X$, to be denoted by regular font,  \index{period function} \index{$P_X$}
\[ P_X : \Bun_G(\FF_q) \rightarrow \kk,\]
which sends
a $G$-bundle to the number of sections of the associated $X \otimes \cK^{1/2}$-bundle, equivalently:
\begin{equation} \label{PXform1} P_X(x)  : g \in G(\mathbb{A}) \mapsto \sum_{x \in X(F)} (g, \partial^{1/2}) \cdot \Phi(  x).\end{equation}
 Here, $\Phi$
is the characteristic function of  $X(\widehat{\mathfrak{o}})$
 inside the adelic points of $X$, and the action of $G \times \GGm$
 on such functions
 is understood as $(g, \lambda) \Phi (x) = \Phi(x (g, \lambda))$ (i.e., it is an \emph{unnormalized} action, as opposed the normalized action that we will introduce in \eqref{PXstardefinition} below).

As is often the case (see \S \ref{analyticarithmetic}), 
 it will be convenient to also have a version of $P_X$
 that incorporates half-twists, because it will relate
 more clearly to unitary structures in the theory of automorphic forms.
 We turn to this next, in \S \ref{normalizedperiod}. 
 
\index{star period}
 \begin{remark} (The star period) \index{$\mathcal{P}_X^*=$ star period} \label{starperiod}
  There is a second variant of the period sheaf, which   is of great interest, although
  it will not play such a role in this paper. This is the $*$-variant $\mathcal P_X^*$ -- 
  understood as the *-push-forward of the {\em dualizing sheaf} from $\Bun_G^X$.  Equivalently, this is the ``naive Verdier dual'' of $\mathcal{P}_X$; that is to say,
 a sheaf on $\Bun_G$ being a compatible system of sheaves on a system
 of open truncations, we apply usual Verdier duality at each level; note however
that \cite{DrinfeldGaitsgorycompact} (cf. \S \ref{tensor products of sheaves}, \S \ref{tensor and miraculous}) this naive Verdier duality is not an equivalence of categories. 
 The classical meaning of this sheaf is not so easy to understand, 
and its existence is an interesting puzzle in the classical theory.    See \S \ref{starperiods}.
\end{remark}

 \begin{remark} 
 Let us discuss technical issues in the de Rham context. As we have seen  (\S \ref{alg-geom-tech1}), 
the map $\Bun_G^X \rightarrow \Bun_G$ has a very simple nature.
  In particular, there is no difficulty in defining either $!$ or $*$ pushforward of the constant sheaf
   along these maps in any of our sheaf-theoretic contexts.
   
Indeed (by definition, cf. \S \ref{renormalization section}) a sheaf on $\Bun_G$  is a compatible system of sheaves
on open quasicompact substacks of $\Bun_G$; and in turn, on each open quasicompact substack,  
sheaves are defined as the ind-completion of
the category obtained by taking the limit of small sheaf categories over 
maps of an affine $Y$ into that open quasicompact substack.
When pulled back to such $Y$, the map $\Bun^X_G \rightarrow \Bun_G$
becomes a morphism of schemes.  

Therefore, $\cP_X$ and $\cP_X^*$ are, locally on $\Bun_G$, even compact objects
of the automorphic category, i.e., objects of the small category of sheaves, though they are not compact themselves since they don't have quasicompact support. 
By means of the  functor \eqref{unrendef} from the ind-finite to the usual category, they can also
be considered elements of the latter.

\end{remark}

 \subsection{Normalized periods and normalized period sheaves} \label{normalizedperiod}
 $\mathcal P_X$ has an important variant.
This is the {\em normalized period sheaf}, which corresponds   
 to formula \eqref{PXform1} using the {\em unitarily normalized} action of $G$. 
  To define a normalized period sheaf  we 
assume the existence of an eigen-volume form as in \eqref{eigencharacter}
with eigencharacter $\eta: G \rightarrow \Gm$.

  \index{$P_X^{\norm}$}
The normalized version of \eqref{PXform1}  
is 
 \begin{equation} \label{PXform2} P_X^{\norm} : g \in G(\mathbb{A}) \mapsto q^{-\beta_X/2} \sum_{x \in X(F)} g \star  ( \partial^{1/2} \cdot \Phi(  x)),\end{equation}
 \begin{equation} \label{PXform3} =  q^{\frac{g-1}{2}(\dim X-\dim G)} \sum_{x \in X(F)}  (g, \partial^{1/2}) \star  \Phi(  x).\end{equation}
 where  $\star$ denotes normalized actions, defined thus:
\begin{equation} \label{PXstardefinition} g  \star\Psi =|\eta(g)|^{1/2}    (g\cdot \Psi ), (g, \lambda) \star \Psi = |\eta(g)|^{1/2} |\lambda|^{\gamma_X/2} (g,\lambda) \cdot \Psi\end{equation}
 are the unitary
  $\beta_X$ is as in \eqref{betaXdef},
 and $\gamma_X$ is as in \eqref{gammaXdef}.   That is to say,  \eqref{PXform2} uses only unitary normalization of $G$, whereas \eqref{PXform3}
 uses unitary normalization of both $G$ and $\GGm$; to deduce \eqref{PXform3} from \eqref{PXform2} we use
 $|\mathfrak{d}^{1/2}| = q^{1-g}$  and the definition \eqref{betaXdef} 
of  $\beta_X$.

To make the sheaf version 
note that $\eta: G\rightarrow \Gm$
induces  $\Bun_G \rightarrow \Bun_{\Gm}$
and in particular a degree function $\deg: \Bun_G \rightarrow \Z$.  
We put 
\begin{equation} \label{PXnormdef} \mathcal P_X^{\norm} := \mathcal P_X  \langle \deg + \beta_X\rangle. \end{equation}
Note that $\mathcal P_X^{\norm}$ is not, in general, a global twist of $\mathcal P_X$
 because of the non-constant function $\deg$. The star period sheaf of Remark \ref{starperiod} also has a normalized variant, obtained by taking naive Verdier duality, or equivalently $\mathcal{P}_X^{* \norm} = \mathcal{P}_X^{*} \langle - \deg - \beta_X \rangle$.

 Let us try to describe the origin of the two twists in \eqref{PXnormdef}, by degree and by $\beta_X$.  \footnote{\label{footnote-deg}{A better, but notationally cumbersome, way to separate the twists is to split off the term $(g-1)\gamma_X$ from $\beta_X$; this term plays exactly the same role for $\GGm$ as the degree plays for $G$ -- remember from \eqref{diagdog} that $\Bun_G^X$ is defined as a fiber over $\cK^\frac{1}{2}\in \Bun_{\GGm}$. In the discussion that follows, references to $\beta_X$ really are meant for the remaining terms $(g-1)(\dim G-\dim X)$.}} 
 \begin{itemize}
 \item[-] The twist by  $\langle \beta_X \rangle$ appears already in the homogeneous case $X=H\backslash G$, in which case $\Bun_G^X=\Bun_H$, and splits the difference between the constant and the dualizing sheaf of $\Bun_H$; in this case, 
 the period sheaf  can be thought of as the push-forward of the intersection complex of $\Bun_H$, but we emphasize that our definition does {\em not} agree with the intersection complex of $\Bun^X_G$ in general. 
The numerical explanation  of the $\beta_X$ twist is that the factor  $q^{-\beta_X/2}$ is an attempt to render $P_X^{\norm}$ approximately $L^2$-normalized
\begin{equation} \label{norm goal} \| P_X^{\norm}\|_{L^2} \approx 1.\end{equation}
\item[-]
 On the other hand, the degree twist $\langle \deg \rangle$ has to do with the linear fiber $S^+$, and is a standard twist in the geometrization of the Weil representation \cite{Lysenko}. If $S^+$ is nontrivial, this twist cannot be interpreted through the intersection complex of $\Bun_G^X$ (which is, in general, singular), but can be thought of as the sheaf-theoretic analog of half-densities in the 
Schr{\"o}dinger model of the Weil representation.
\end{itemize}

 \begin{remark}[The degree sheaf]\label{degree sheaf} \index{degree sheaf} \index{$\underline{\deg}$}
It will be convenient to represent the degree shift $\cF\mapsto \cF\langle \deg \rangle$ as tensor product with a sheaf. Namely, this is achieved 
(in any of our contexts) by the ($\eta$-pullback of the) locally constant sheaf $$\ul{\deg}\in \Autshv(\Bun_\Gm),$$ \index{degree sheaf} \index{$\ul{\deg}$} which is given by $\kk \langle d \rangle$ on $\mathrm{Pic}^d$;
in the language of \S \ref{shearingsec0}, $\underline{\deg}$ is the shear of the constant sheaf by the $\Gm$-action on $\Autshv(\Bun_\Gm)$ corresponding to the $\Z$-grading by degree. 
Under the sheaf-function correspondence  (see \S \ref{anglebracketnotation}) 
 this corresponds to the function $\mathcal{L} \mapsto q^{-\mathrm{deg} \mathcal{L}/2}$.
 If we uniformize line bundles ad\`elically
 via \eqref{adelicuniformization}  this matches with $x \in \mathbb{A}^{\times} \mapsto |x|^{1/2}_{\mathbb{A}}$. In other words, $\ul{\deg}$ is an automorphic avatar of the square root of the cyclotomic character $\varpi^{1/2}$ under class field theory.
\end{remark}

 \begin{remark} \label{Heisenberg comment}
It is worth noting that there are three shifting processes appearing implicitly in the above discussion.

\begin{itemize}
\item[(a)]   the shift by $\cK^{1/2}$ embedded in the definition of $\Bun^X_G$, which reflects the $\GGm$ action;
\item[(b)]  the twist by $\langle \deg \rangle$ in \eqref{PXnormdef}; this reflects the failure of $X$ to be unimodular; 
\item[(c)] the twist by $\langle \beta_X\rangle$.
\end{itemize}
Roughly speaking, the first twist is related to a translation on the automorphic side,
whereas the second twist is related to a translation on the spectral side -- see Remark \ref{Heisenberg comment2}. 
As such, these twists do not commute with one another: as is usual in Fourier analysis, 
spectral and automorphic translations do not commute. 
 In physics automorphic and spectral translations correspond to shifts of magnetic and electric flux, respectively, and the lack of commutativity is an aspect of the ``uncertainty of fluxes'' studied in~\cite{FreedMooreSegal1,FreedMooreSegal2}; see also Remark~\ref{fluxes}. 
\index{fluxes}
  
  \end{remark}

\subsubsection{Changing the grading}\label{changing grading automorphic}
The following remarks are not essential to understanding the main conjecture but
will be used later in analysis of parity issues. 
Given a central cocharacter $\lambda:\Gm\to Z(G)$ we denote by $X[\lambda]$ the $G\times \GGm$-space $X$ with $\GGm$-action twisted by $\lambda$.
The effect of passing from $X$ to $X[\lambda]$ is to translate the unnormalized period function/sheaf by the translation action of the $Z(G)$-torsor $\lambda(\cK^{1/2})$ (or numerically by the central element $\lambda(\mathfrak{d}^{1/2})$ of $G$):
$$\mathcal{P}_{X[\lambda]} = \lambda(\cK^{1/2})\ast \mathcal{P}_X.$$
Here, 
to normalize signs, ``translation by $\lambda(\cK^{1/2})$'' means that the delta function at a point would be sent to the delta function at its translate by $\lambda(\cK^{-1/2})$.   We will prove that the normalized period sheaf is transformed similarly:

\begin{lemma}  \label{PXGmshift}  The operation $X\mapsto X[\lambda]$ affects the normalized period sheaf as follows:
$$ \mathcal{P}_{X[\lambda]}^{\norm} \simeq  
\lambda(\cK^{1/2})\ast
 \mathcal{P}_X^{\norm}.$$
\end{lemma}
 
 \proof
Observe first of all that, see \eqref{betaXdef}, 
\begin{equation} \label{betachange} \beta_{X[\lambda]} =  \beta_X +  (g-1) \langle \lambda, \eta \rangle.\end{equation}
If we translate the sheaf $\langle \deg \rangle$ on $\Bun_{\Gm}$ by $\cK^{1/2}$
we get $\langle \deg + (g-1) \rangle$. Correspondingly,  since we understand $\langle \deg \rangle$ on $\Bun_G$
as the pullback of the $\Gm$ sheaf via $\eta: G \rightarrow \Gm$, 
\begin{equation} \label{Tshift} \mathrm{T}( \mathcal{F} \langle \deg \rangle) = (\mathrm{T} \mathcal{F}) \langle \deg + (g-1) \langle \eta, \lambda \rangle \rangle.\end{equation}
where $\mathrm{T}$ refers to a translation by $\lambda(\cK^{1/2})$. Thus, from the definition \eqref{PXnormdef}, 
 $$ \mathcal{P}_{X[\lambda]}^{\norm} = ( \mathrm{T}   \mathcal{P}_X) \langle \deg +\beta_{X[\lambda]} \rangle, \mathrm{T} 
 \mathcal{P}_X^{\norm} =   \mathrm{T} \mathcal{P}_X \langle \deg + \beta_X + (g-1) \langle \eta, \lambda \rangle \rangle.$$
 and comparing with \eqref{betachange} we deduce the lemma. 
\qed

\subsection{Modification for the twisted case} \label{Affcase} 
\index{$\mathcal{P}_X$ (twisted case)}
We are now going to explain how to define the period sheaf
and its variants in the case of twisted polarizations $(X,\Psi)$ (see \S 
\ref{twisted cotangent bundle section}). 
At a high level, the modification is simply
\begin{quote}
{\em twist by a rank one ``Artin--Schreier'' local system on $\Bun^X_G$}.
\end{quote}

Here, the Artin--Schreier sheaf is pulled back from a sheaf on $\Ga$
in the finite and de Rham context, and a slight modification in the Betti context.
We start by defining it on $\Ga$, and then explain how it is to be pulled back to $\Bun^X_G$. 

\begin{definition}\label{def:ArtinSchreier} \index{Artin--Schreier sheaf}
We understand the ``Artin--Schreier sheaf'' on thus:
\begin{itemize}
\item[(a)] In the finite context, it is the {\'e}tale sheaf of rank one $\kk$-modules on  $\Ga$ whose trace function is the fixed additive character $\FF_q\to \kk^\times$, obtained, in the usual way, from the covering $x^q-x$, see e.g.\ \cite[Sommes. trig.]{SGA4.5}.
\item[(b)] In the de Rham context, it is the exponential $D$-module on $\Ga$, i.e.,
the sheaf associated to the differential equation $f'=2\pi i f$. 
\item[(c)]
  In the Betti context we understand the ``Artin--Schreier sheaf''
  to be a locally constant sheaf in the analytic topology on  $\Ga/\Gm$, where $\Gm$ 
  acts by squaring on $\Ga$, defined as
    $$\left( j_! \kk^- \oplus j_* \kk \right) [ -1]$$ where  
  $j: \Gm \hookrightarrow \Ga$,
  and $\kk, \kk^{-}$ are respectively the trivial
  and nontrivial one-dimensional local systems on $\Gm/\Gm \simeq B\mu_2$.

\end{itemize}
\end{definition}

The point in (c) is that  although we cannot
 make sense of a Artin--Schreier sheaf on $\Ga$ itself, its pushforward to $\Ga/\Gm$
 makes sense (cf.  ~\cite[Section 2.5.2]{NadlerYun3pts}), 
 and this will be sufficient for our purposes.  To compare with \cite{NadlerYun3pts}, note that in our case
 the action of $\Gm$ on $\Ga$ is the squaring action; computing
 the pushforward of the sheaf used in op.cit.\ leads to the formula above.

 Now, in our current situation, $M=T^*(X, \Psi)$ with $\Psi$ a $\Ga$-bundle over $X$, 
and the $\Gm$ action scales $\Ga$ via squaring. 
Recall from \S \ref{bunX} that we define $\Bun_G^X$ as the fiber of $\Map(\Sigma, \frac{X}{G\times \Gm})$ over $\cK^{\frac{1}{2}}$.
An affine bundle $\Psi \rightarrow X$ as in \S \ref{dp} defines a map
 $X \rightarrow B\Ga$
equivariant for the action of $G \times \GGm$ where the $\GGm$-action scales $\Ga$
by the square character, and $G$ acts trivially. In particular, we have a map 
$$ \Map(\Sigma, \frac{X}{G\times \Gm})\to \Map(\Sigma,  \frac{\mbox{pt}}{\Ga \rtimes \Gm}),$$
wherein the action of $\Gm$ on $\Ga$ is by squaring.  

Correspondingly, 
$\Bun^X_G$ maps to the stack of $\Ga \rtimes \Gm$ torsors reducing to $\cK^{1/2}$; 
 said differently,
this is the stack of torsors for $\cK=(\cK^{1/2})^{\otimes 2}$ considered as a vector bundle; 
this stack is identified with
the quotient of $H^1(\Sigma, \Omega^1) \simeq \mathbb{G}_a$ by the trivial action of $H^0(\Sigma, \Omega^1)$, and in particular maps to the affine line $\mathbb{G}_a$. This gives us a (non-schematic)
 morphism: 
\begin{equation} \label{GaGmT} \Bun^X_G \longrightarrow \mathbb{G}_a\end{equation}
which measures the obstruction of lifting the $X$-section to a $\Psi$-section.
 This already suffices to define the Artin--Schreier sheaf on $\Bun^X_G$ in the de Rham and finite contexts as the pullback of the corresponding sheaf on $\mathbb{G}_a$,\footnote{Or, to only use a schematic morphism,  one first notes that the Artin--Schreier sheaf on $\mathbb{G}_a$
descends to its quotient by the trivial action of $H^0(\Sigma, \Omega^1)$, and then pulls back to $\Bun^X_G$.}
and we define $\mathcal{P}_X$   as the compactly supported pushforward of this sheaf:
 $$ \mathcal P_X = \mbox{$!$-pushforward of the Artin--Schreier sheaf along $\Bun_G^X \rightarrow \Bun_G$.}$$
We also define the $*$-period sheaf $\mathcal{P}_X^*$ as the ordinary pushforward of the Verdier dual of the Artin--Schreier sheaf.

In the Betti context we note that we get also
$$\Bun^X_G/\Gm \rightarrow \mathbb{G}_a/ \mathbb{G}_m$$
where the $\Gm$ action on the left arises from that on $X$; and on the right it is squaring. 
We can correspondingly define the Artin--Schreier sheaf on $\Bun^X_G/\Gm$ by pullback.  Since 
 the morphism $\Bun^X_G \rightarrow \Bun_G$ factors through the quotient
$\Bun^X_G/\Gm$ we then define $\mathcal{P}_X$  by pushing
forward this Artin--Schreier sheaf via $\Bun^X_G/\Gm \rightarrow \Bun_G$. 

We introduce normalized versions of these sheaves according to precisely
the same shift $\langle \deg + \beta_X \rangle$ that occurred previously in \eqref{PXnormdef}.

\begin{remark}\label{dependenceAS}
 Note that, in the twisted case, the definition of the period sheaf depends on the choice of an Artin--Schreier sheaf,
 which we fixed in
Definition \ref{def:ArtinSchreier}; recall that in the finite case, this depends  on the choice of an additive character of $\mathbb{F}_q$. 
 We consider this choice fixed throughout the paper.
\end{remark}

 \subsubsection{Period function}
 We will now describe explicitly the associated period function, i.e., trace of Frobenius, which
recovers well-known ``Fourier--Whittaker periods'' in the theory of automorphic forms. 

Let  $\psi'$ be the character of $\adele$ given by $\psi'(x) = \psi(\partial^{-1} x)$.  This need not be trivial on $F$;
as we have seen however it is ``unramified'' on each $F_v$. 
Recall that $\Psi$ is the total space of a $\mathbb{A}^1$-bundle over $X$. 
 Consider then the induced space of functions
\begin{equation} \label{Whittaker function} \Phi: \Psi(\mathbb{A}_F) \rightarrow k, \ \
\Psi(\tilde{x} \cdot t) =  \Psi(\tilde{x}) \psi'(t) \ \ (\tilde{x} \in \Psi(\mathbb{A}_F), t \in \mathbb{G}_a(\mathbb{A}_F)),\end{equation}
and inside it consider the ``basic function'' which is characterized
by satisfying \eqref{Whittaker function}, being supported on the preimage of $X(\widehat{\mathfrak{o}})$,
and being identically $1$ on $\Psi(\widehat{\mathfrak{o}})$. 

For such a  function $\Phi$ the translate $ \partial^{1/2} \cdot \Phi: \Psi(\mathbb{A}_F) \rightarrow k$
satisfies 
$$ (\partial^{1/2} \cdot \Phi) ( \tilde{x} t) =  \psi'((\partial^{1/2})^2 t) \Phi(\tilde{x} \partial) = \psi(t) (\partial^{1/2} \cdot \Phi)(\tilde{x})
\ \ (\tilde{x} \in \Psi(\adele), t \in \mathbb{G}_a(\adele)),$$
and therefore equalling $\partial^{1/2 }\cdot \Phi(x)$ if $t \in F$. 
In particular, for $x \in X(F)$, the value of $\partial^{1/2} \cdot \Phi$ on any lift $\tilde{x} \in \Psi(F)$
is independent of choice of $\tilde{x}$; let us call this simply $\partial^{1/2} \cdot \Phi(x)$.  
With these conventions \eqref{PXform1} still holds.

 \begin{example} \label{Whit explicit computation example}
 For later usage, we will write out a  formula in the Whittaker case, i.e., $X=U\backslash G$
 with the $\mathbb{G}_a$-torsor defined by the sum $U \rightarrow \mathbb{G}_a$ of identifications of the simple root spaces with $\Ga$.  There is nothing
novel here, but it will be useful to have explicit formulas for the various shifts and constants.
 Let $f$ be a function on $\Bun_G(\FF_q)$. We will compute
\begin{equation} \label{FTXcollapse} \sum_{\Bun_G(\FF_q)} P_{X}(g) f(g) = \int_{G(F) \backslash G(\mathbb{A})} P_X(g) f(g), \end{equation}
 where the sum is taken over $G$-bundles weighted by inverse size of their automorphism group (i.e., the sum over $\Bun_G(\FF_q)$ considered as a groupoid), 
 and the latter integral is taken with respect to  the measure with $\vol (G(\hat{\mathfrak o})) = 1$.  
 
 Unfolding via \eqref{PXform1}, taking into account the prior discussion, the above equals
 $ \int_{U(F) \backslash G(\mathbb{A})} f(g)  (g, \partial^{1/2}) \cdot \Phi( g).$
 Now, $\partial^{1/2}$ is acting on $X=U \backslash G$ by means of (see \S \ref{ggmm}) left translation through the element
 $a_0^{-1}$, where  \index{$a_0$}
\begin{equation} \label{a0def} a_0 := e^{-2 \check{\rho}}(\partial^{1/2}) \in T(\mathbb{A}_F)\end{equation}
 and correspondingly $\partial^{1/2} \cdot \Phi$ is 
  is supported on $ U(\mathbb{A})   a_0 \cdot G(\widehat{\mathfrak{o}})$. 
  Write $du$ for the measure on $U(\mathbb{A})$ 
  where $U(\mathfrak{o})$ has measure $1$;
  it assigns to $U(\mathbb{A})/U(F)$ the measure $q^{(g-1) \dim U}$. 
  The measure  $d( a_0^{-1} u a_0)$
  assigns to $a_0 U(\mathfrak{o}) a_0^{-1}$ the mass $1$, 
  and equals $du$ multiplied by $|e^{2 \rho}(a_0^{-1})|=  q^{-(g-1)\langle 2\rho, 2 \check{\rho}  \rangle}$. 
  Thus \eqref{FTXcollapse}   equals
  $$ q^{-(g-1) \langle 2\rho, 2 \check{\rho}   \rangle} \int_{U_F \backslash U_{\mathbb{A}}}  \psi(u) f(u a_0 ) du
  =q^{(g-1)(\dim U - \langle 2 \rho, 2 \check{\rho} \rangle)} \int_{U_F \backslash U_{\mathbb{A}}} \psi(u) f(u a_0 ) d'u.$$
  where in the latter  integral we use the Haar probability measure.  Observe
  that the exponent of $q$ is simply $\beta_X = (g-1) (\dim U - \langle 2 \rho, 2\check{\rho} \rangle)$
 (compute via   \eqref{betaXdef}  and \eqref{gammadef}), 
   so we get 
  \begin{equation} \label{FTXcollapse2} \sum_{\Bun_G(\FF_q)} P_{X}(g) f(g) = q^{\beta_X}
  \int_{U_F \backslash U_{\mathbb{A}}} \psi(u) f(ua_0) d'u,  
  \end{equation}
  and the analog for the normalized period where we replace $\beta_X$ by $\beta_X/2$, 
  see \eqref{PXform2}.  In the last formula, $d'u$ is again the invariant measure with total mass $1$.

  \end{example}

\subsection{$\Bun_G^X$ and the period sheaf : examples}  \label{BunXex}

The following collection of examples are intended to give some indication
of the geometry involved with both period sheaves and the spaces $\Bun_G^X$. 
Note that we will sometimes consider examples of $X$ with non-neutral $\GGm$-action,
i.e., with a $\GGm$-action that differs from that specified in \S \ref{dp}; 
the point to note here is that the choice of $\GGm$-action really affects $\Bun_G^X$! 

 \subsubsection{ Homogeneous spaces:}  \label{homspace} 
 We consider the case of $X$ with trivial vectorial and $\SL_2$ component, i.e., 
  $X=H\backslash G$ for $H$  reductive. 
  Here we will use the neutral $\GGm$ action, as specified in \S \ref{dp};
  what makes this case particularly easy think about is that
  \begin{quote} the neutral $\GGm$ action is {\em trivial}.
  \end{quote}
 Also, in this case, there exists a $G$-invariant volume form on $X$ with $\eta$ and $\gamma$ both trivial.

The map
 $\Bun_G^X \rightarrow \Bun_G$ is simply identified with the natural map $\Bun_H \rightarrow \Bun_G$;
 the period sheaf $\mathcal P_X$ assigns to a $G$-bundle the compactly supported cohomology 
 of the space of reductions to an $H$-torsor, and the period function $P_X$ counts
 the number of such reductions (i.e., the size of the fibers of $\Bun_H \rightarrow \Bun_G$). 
 The normalized period sheaf twists by $\langle b_H \rangle$, and 
  the normalized period function takes the value
 \begin{equation} \label{PXexp} P_X^{\norm}(x) = P_X(x) \cdot q^{-b_H/2}.\end{equation}

 \subsubsection{The Iwasawa--Tate period.} \label{TateBunX}
 
 Take $X=\mathbb{A}^1$ as a $G=\Gm$-space, and let us start
 by considering the {\em trivial} $\GGm$-action.   The unnormalized period function is equal to
      $$P_X: \mathcal{L} \rightarrow q^{h^0(\mathcal{L})} \ (\textrm{trivial $\GGm$-action})$$
      (for $\mathcal{L}$ a line bundle on $\Sigma$). 
      The star period $P_X^*$ is  more interesting and is described in \S \ref{starTate}.

 Let us describe the geometry of $\Bun_G^X$ and the period {\em sheaf} in this case.
  There is a tautological map
         $$ \pi_r: \mathrm{Sym}^r \Sigma \rightarrow \mathrm{\Bun}_{\Gm}, \ \ (Q_1, \dots, Q_r) \mapsto \mathcal{O}(\sum Q_i).$$
      The fiber of $\pi_r$
      above a line bundle $\mathcal{L}$ is the space of such effective divisors $\{Q_i\}$ together with an isomorphism $\mathcal{O}(\sum Q_i)\simeq \mathcal{L}$, which is the same as the
      space of nonzero global sections of $\mathcal{L}$; that is to say, the fibers of $\pi_r$ are punctured affine spaces. 
  Note here that we are really regarding $\Bun_{\Gm}$ as a stack; if, in our discussion,
  we were to replace its  role by the Picard  {\em scheme}, the analogous fibers would be projective spaces.

With this in hand, we can describe
 $$ \{\Bun_G^X \rightarrow \Bun_G \} = \mbox{partial compactification of $\pi_r$} $$
where we  allow the zero sections of line bundles, i.e., omitting the phrase ``nonzero'' in the above description.
      Note that, unlike \S \ref{homspace}, $\Bun_G^X$ is not smooth.

 Now let us switch to the {\em scaling} $\GGm$-action, i.e., the ``neutral'' action. 
 Now
 $\Bun_G^X\to \Bun_G$ parametrizes line bundles $\Ll$ with a section of $\Ll\ot \cK^{1/2}$ (and its geometric description
is exactly parallel to that given in the previous paragraph). The normalized period function is
\begin{equation} \label{PXTatenorm} P_X^{\norm} \mbox{ with scaling $\GGm$}: \mathcal{L} \mapsto q^{h^0(\mathcal{L}\otimes \cK^{1/2}) - \frac{1}{2} \deg(\mathcal{L} \otimes \cK^{1/2})}.\end{equation}
Observe this now has a pleasing symmetry $\mathcal{L} \leftrightarrow \mathcal{L}^{-1}$. 

 Let us spell out the twists.
 Here, $\beta_X$ as in \eqref{betaXdef} is given by $(g-1)$ 
 and the factor $q^{-\beta_X/2}$ of   \eqref{PXform2},  coincides with $q^{-\frac{1}{2} \deg(\cK^{1/2})}$; also, 
 $|\eta(g)|^{1/2}$ contributes $q^{-1/2 \deg(\mathcal{L})}$. Geometrically,
 the normalized period sheaf twists the period sheaf
 by $\langle d + g-1 \rangle$ 
 on the component where $\deg(\mathcal L)=d$, categorifying the factor  $q^{-\frac{1}{2}(d+g-1)}$. We can think of $g-1$ as the dimension of $\Bun_G$ and
 $d$ as the Euler characteristic of $\mathcal L \otimes \cK^{1/2}$, i.e.
the expected dimension of fibers of $\Bun^X_G \rightarrow \Bun_G$, so that all in all
this is the analog of the twist by $b_H$ appearing in the case \S \ref{homspace}.

  \subsubsection{The Eisenstein case $U \backslash G$} \label{A2weirdness}

  For $G$ arbitrary take 
    $$X=U\backslash G$$ 
    as a $G \times T$-space    i.e., $(g,t): Ux \mapsto U t^{-1} x g$.
  This case does not fall in our general setup,  for $X$ is not affine, but nonetheless our definitions 
  of period sheaf and period function make sense, and it will be valuable for us to examine them.

    There are two $\GGm$-actions we shall consider in this paper; one is trivial,
    and the other, which we shall examine here, is 
    where $\GGm$ acts via the restriction of  the $G \times T$-action 
    via $(1,e^{-2\rho})$. (A discussion of the relation between these two actions is given, in a more general context, 
    in \S \ref{mitch}). 
  Explicitly  in the latter action $\lambda \in \GGm$
   acts through left multiplication on $U \backslash G$ by $\lambda^{2\rho} \in G$. 
  
   For example, in the case  $G=\SL_2$, this $X$ is the
  punctured affine plane via $g \in \SL_2 \mapsto (0,1) g$;
  if we take the nontrivial action of $\GGm$
  it amounts to inverse scaling on this punctured plane.
Then the fiber $\Bun_G^X \rightarrow \Bun_G$ over a rank $2$ unimodular vector bundle $V$ is 
the space of
everywhere injective maps $\cK^{1/2} \rightarrow V$, equivalently,
the space of extensions  $$ \cK^{1/2} \rightarrow V \rightarrow \cK^{-1/2}.$$
 The total space of $\Bun_G^X$ can thereby
  identified with the affine space $H^1(\Sigma, \cK) \simeq \mathbb{A}^1$
modulo automorphisms $H^0(\Sigma, \cK)$.    
 
 More generally, if the base curve has genus $\geq 2$, and we take $X=U \backslash G$
 with the nontrivial $\GGm$ action just described,  we may identify 
   \begin{equation} \label{ArBun} \mathrm{\Bun}^X_G= \mathbb{A}^r/\mathcal{U}\end{equation} as the quotient of
   the affine space $\mathbb{A}^r$ (where $r=$ the semisimple rank, and the coordinates are indexed by simple roots for $G$)
by a {\em trivial action} of a certain unipotent group scheme $\mathcal{U}$.
The map  $\Bun_G^X \rightarrow \mathrm{\Bun}_G$
is not a closed immersion, but rather  factors 
through the quotient of $\Bun_G^X$ by 
the action of a torus in $G$; as a result, at the level of points over an algebraically closed field, 
the image of $\mathbb{G}_m^r \subset (\mathbb{A}^1)^r$ is a single
point of $\Bun_G$.

\subsubsection{The Whittaker case. }  
The Whittaker case was already discussed in Example \ref{Whit explicit computation example}
and we just make a couple of minor additional remarks. 
With reference
to the presentation \eqref{ArBun}, the period function is given by pushing forward $\psi(\sum_{i=1}^r x_i)$ on $\mathbb{A}^r$
to $\Bun_G$,  and the period sheaf
the geometrization of that construction via Artin--Schreier sheaves.
By definition and \eqref{sigmaXdef}, 
the normalized period function is related to the unnormalized one by \begin{equation}\label{Whittakernormalized} P_X^\norm = 
q^{-\beta_X/2} P_X, 
\beta_X = (g-1) [-\langle 2 \rho, 2 \rho^{\vee} \rangle +\dim U].
\end{equation}
 
 \begin{remark}\label{remarkbetaX}\index{$\beta_X$}
The quantity $\beta_X$ is closely related to dimension of $\Bun^X_G$;
this is easy to see in the homogeneous case $H \backslash G$,
but also remains true in the twisted case. 
For example the Whittaker case just described, the dimension of $\Bun_G^X$, i.e., the space of \eqref{ArBun}, is  the sum of  negated Euler characteristics
of bundles $\langle \alpha, \rho \rangle \cdot \cK$
over positive roots $\alpha$, with $\cK$ the canonical divisor;  
{\small \begin{equation}  \label{mudef}  \dim \Bun_G^X =   (g-1) \sum_{\alpha} \left[ 1- 2\langle \alpha, \rho \rangle    \right] =   (g-1) [-\langle 2 \rho, 2 \rho^{\vee} \rangle +\dim U]
= \beta_X.
\end{equation}}
 \end{remark}

\subsection{Dependence on spin structure} \label{spindep} \label{period sheaf and spin}
  
As mentioned, we have felt free to choose a spin structure. 
 However, it is sometimes desirable to have a formulation
which is manifestly independent of spin structure. 
We will discuss such a formulation now, which will use 
  the extended dual group (\S \ref{subsection-extended-group}, \S \ref{extended-group appendix}).
Note, however, we will make little use of this formulation and include it for completeness. 

For the discussion that follows, we assume that $M$ satisfies the parity condition discussed
in \eqref{P2}. Namely, we assume that we specify a \begin{quote} \label{parityglobalgeometric}   central involution $z \in G$ 
whose action on $M$ coincides with $-1 \in \GGm$.  
\end{quote}
 As we observed in the discussion surrounding \eqref{P2},
if $M$ admits a dual hyperspherical pair $(\check{G}, \check{M})$
with $\check{M}$ polarized, then $z$ is the product of $e^{2 \rho}$
with the dual of the character by which $\check{G}$
acts on an eigenmeasure on $\check{X}$.

\index{$\zBun_G$}
 Assuming \eqref{parityglobalgeometric}, the action of $G \times \GGm$ then factors through its quotient by $(z, -1)$, which
is precisely the extended group  ${}^CG_z$,   (\S \ref{subsection-extended-group},\ref{extended-group appendix}). 
 Now define 
 $\zBun_G$ as
the fiber of $\Bun_{{}^CG_z} \rightarrow \Bun_{\Gm}$ above the
canonical bundle. Equivalently we can write $$\zBun_G\simeq \Bun_G \times^{\Bun_{\Z/2}} \mathrm{Spin}_\Sigma,$$
where $\mathrm{Spin}_\Sigma$  \index{$\mathrm{Spin}_\Sigma$ stack of spin structures on $\Sigma$} denotes the stack of spin structures (square roots of $\mathcal K$) on $\Sigma$, and the group stack $\Bun_{\Z/2}$ of $\Z/2$-torsors on $\Sigma$ acts both on $\Bun_G$ via the central embedding $z:\Z/2\to Z(G)$, and simply transitively on the 
$\mathrm{Spin}_{\Sigma}$. 
The choice of a spin structure $\cK^{1/2}$ gives rise to an identification
\begin{equation} \label{zbunbun} \zBun_G \simeq \Bun_G\end{equation}
and changing this choice $\cK^{1/2}\mapsto \cK^{1/2}\otimes \mathcal{L}$ by a 2-torsion line bundle changes the identification by the translation action of 
$\mathcal{L}$ on $\Bun_G$. 

Next we define a period sheaf on $\zBun_G$ by replacing pushforward along $\Bun_G^X\to \Bun_G$ by pushforward along its twisted version
that renders the following square Cartesian:  \begin{equation*} \xymatrix{\zBun_G^X\ar[r]\ar[d] &\Map(\Sigma, \frac{X}{{}^C G_z})\ar[d]\\
\zBun_G\ar[r]& \Bun_{{}^C G_z}}.\end{equation*}
 Our above definition of $\Bun_G^X$ is obtained from this one by  transporting   via  \eqref{zbunbun}
 after fixing a spin structure. 
  An equivalent way of formulating this is that the
  period sheaf is (independently of the choice of spin structure) defined as an object in a twisted version of sheaves on 
  $\Bun_G$:
\begin{equation} \label{pxtwisted} \cP_X\in \Hom_{\Bun_{\Z/2}}(\mathrm{Spin}_\Sigma, \Autshv(\Bun_G))\end{equation} where
$\Bun_{\Z/2}$ acts on $\Bun_G$ via the central homomorphism
$\{\pm 1\} \rightarrow Z(G)$ sending the nontrivial element to $z$.

\subsection{Reduction to the vectorial case}  \label{Pinduction}  

Hyperspherical varieties are built (as in Theorem \ref{thm:structure})
by a process of Whittaker-induction from the special case of $M$ a vector space.
This gives rise to a corresponding structure for period sheaves:
\begin{quote}
the period sheaf for general $M$ is a Whittaker induction for the period sheaf in the case of $M$ a vector space.
\end{quote}
 
 We will make this explicit. 
This explication is completely straightforward, and can be referred to only as necessary; 
we note it mainly for reference and as a comparison point for a similar  (but less straightforward) discussion in the spectral case. 
  For simplicity we restrict ourselves to unnormalized period sheaves in this section. 

Recall that,  in the case $\FF=\C$,  a polarized hyperspherical variety $M$
has the structure (see Theorem \ref{thm:structure}) of a 
  Whittaker-induction
along $H \times \SL_2 \rightarrow G$ of a  polarized symplectic $H$-representation
$T^* S$, with associated twisted polarization $X = S \times^{HU} G$ with $\Ga$-bundle $\Psi \rightarrow X$. 
 In what follows, the fact that $M$ is hyperspherical will not matter;
all that matters is the homomorphism $H \times \SL_2 \rightarrow G$
and the  $H$-space $S$.  Let us review the notation in more detail. 
 
We fix a homomorphism $H\times SL_2\to G$ with underlying cocharacter $\varpi:\Gm\to G$.
We will also restrict to the situation (automatic in the polarized case
by Definition  \ref{distinguishedpolarizationdef})
where all the    $\varpi$-weights on the Lie algebra are even.\footnote{There is 
 an analogue of the construction that follows without this requirement, but now involving a geometric
  version of a Jacobi $\theta$ function. To simplify our life, we simply exclude this situation.}
Let $U=U_+\subset G$ be the unipotent subgroup defined by the positive part of the grading.
To this data we can associate
the $\mathbb{A}^1$-bundle $(\Psi \longrightarrow U \backslash G)$
where $\Psi=U_0 \backslash G$, 
$U_0$ being the kernel of $U \rightarrow \mathbb{G}_a$.
 
 These are $G \times H$-spaces: $G \times H$ acts on compatibly on  $\Psi$ and $U \backslash G$ by the rule $(g,h): U_0 x \mapsto U_0 h^{-1} xg$. 
Correspondingly,  the twisted cotangent   bundle $T^*_{\Psi}(U\backslash G)$ 
 (see  \S \ref{twisted cotangent bundle section}  for definition)  is a graded  Hamiltonian $G \times H$-space, cf.
 Example \ref{Slodowy slice example}; it can be regarded as the Whittaker induction along $H \times \SL_2 \rightarrow G$
 of $T^* H$.

Associated to the $G\times H$-space $U\backslash G$, endowed with the affine bundle $\Psi$, we obtain
 the Whittaker period sheaf $\cP_{U\backslash G,\Psi}\in \Autshv(\Bun_G\times \Bun_H)$ following \S \ref{Affcase}, explicitly,
 the pushforward under 
 $\Bun_{HU}\to \Bun_G\times \Bun_H$ of the corresponding pulled-back Artin--Schreier sheaf;
 this can be regarded as the quantization of $T^*_{\Psi}(U\backslash G)$. 
 We use this sheaf as an integral transform to define automorphic Whittaker functoriality.
 To avoid difficulties with *-pullback and !-pushforward in the de Rham setting, we 
restrict to {\'e}tale or Betti settings, although it would be interesting to give a uniform treatment along the lines
of \S \ref{Pshkernels}.
 
\index{$\mathsf{WI}$ Whittaker induction}
\begin{definition}\label{automorphic Whittaker induction def}  
 The automorphic Whittaker induction functor $$\mathsf{WI}: \bigautshv(\Bun_H) \longrightarrow \bigautshv(\Bun_G)$$ 
is the integral transform given by the Whittaker period sheaf $\cP_{U\backslash G,\Psi}$:
$$\mathsf{WI}(\cF)=\pi_{1!}(\pi_2^*\cF\otimes \cP_{U\backslash G,\Psi}).$$ 
\end{definition}

By the projection formula, Whittaker induction is equivalently described as the integral transform given by the Artin--Schreier sheaf on 
$\Bun_{HU}$ (as a correspondence between $\Bun_G$ and $\Bun_H$).
It follows that the Whittaker induction of the constant sheaf on $\Bun_H$ recovers the period sheaf associated to the $G$-space $(X=H U\backslash G,\Psi)$.\
More generally, using base change and the projection formula one checks that Whittaker induction commutes with the formation of period sheaves:

\begin{lemma}  \label{automorphic 1022}
Given a homomorphism $H\times SL_2\to G$, where $\mathbb{G}_m \subset \SL_2$
has only even weights in its action on the Lie algebra of $G$, and $S=T^*Y$ a polarized Hamiltonian $H$-space, the period sheaf of the Whittaker induction $(X=Y\times^{H U} G,\Psi)$ of $S$ is naturally identified with the Whittaker induction of the period sheaf of $S$:
$$\mathsf{WI}(\cP_{Y})\simeq \cP_{X,\Psi}.$$ 
\end{lemma}

This Lemma allows us to reduce certain questions
about period sheaves  
to the case of symplectic representations, 
see for example the next section \S \ref{Pindep}.

\begin{remark}[Whittaker reduction]
We can also use the Whittaker period sheaf $\cP_{U\backslash G,\Psi}$ as an integral transform in the opposite direction
to define a Whittaker restriction (or ``Whittaker-Jacquet'') functor $$\mathsf{WJ}: \bigautshv(\Bun_G) \longrightarrow \bigautshv(\Bun_H).$$ 
An analogous argument shows that $\mathsf{WJ}$ performs Whittaker reduction on period sheaves, i.e., takes the period sheaf for a polarized $G$-space $M$ to that of the (twisted-polarized) Hamiltonian $H$-space given by its reduction $M\GIT_\psi U$.
It would be interesting to verify if $\mathsf{WJ}$ is identified with the left adjoint of $\mathsf{WI}$ -- this seems technically nontrivial
because of the absence of smoothness or properness above. 
\end{remark}

\subsection{Independence of polarization.}\label{Pindep}  
We have described the construction of period sheaves for polarized hyperspherical varieties $M$. This leaves two natural questions: show that the period sheaf is independent of the polarization, i.e., depends only on $(G \times \GGm, M)$ --  and extend the definition to $M$ which don't admit a polarization. We address the sheaf-theoretic question in the finite setting; probably a version works in the other contexts too, but we did not check.

\begin{proposition} \label{prop:Pindep}
 Let $M = T^*(X_1,\Psi_1) = T^*(X_2,\Psi_2)$ be two distinguished polarizations (\S~\ref{dp}, \S~\ref{hdprings}) of an $\FF$-hyperspherical variety $M$,
 defined by a  completely reducible datum $\mathcal{D}_\FF$ over $\FF$ (see discussion below). 
 
  There is a  isomorphism $\mathcal P^\norm_{X_1}\simeq \mathcal P^\norm_{X_2}$ between their normalized period sheaves, in the finite setting.
  When the hyperspherical data specifying $M$ and the polarizations are defined over a finite field $\FF_q\subset \FF$, the corresponding normalized period functions $P^\norm_{X_1}$ and $P^\norm_{X_2}$ are equal. 
\end{proposition}

Note that this independence is asserted \emph{having fixed an Artin--Schreier sheaf}, which affects the definitions in the twisted case (see Remark \ref{dependenceAS}). Recall that hyperspherical varieties (and polarizations) over $\FF$ have been defined in Definition \ref{ringdatum} using the notion of a hyperspherical datum $\mathcal D_\FF$ which includes a symplectic or usual representation of a reductive subgroup $H$ of $G$ (over $\FF$). We say that the datum is completely reducible if this is the case for that representation of $H$.\footnote{ This complete reducibility is, as usual, automatic in  characteristic
	that is large relative to the weights of $\rho$ or $\rho^+$; 
	this can deduced from \cite[Part II, Chapter 6]{JantzenAlgebraicGroups}.}

\begin{proof}

The only difference between two polarizations of $M$ arises from
the possibility of two different $\rho^+$ polarizing the same representation $\rho$.
Using 
 Lemma \ref{automorphic 1022}, the sheaf-theoretic statement is reduced to the vectorial case, i.e., when $M=S$ is a symplectic representation of $G$,
 polarized in two different ways:
  $$ S = X_1 \oplus X_1^* = X_2 \oplus X_2^*.$$
The statement about period functions similarly reduces to this case, as well. 
By virtue of our assumption of conclude reducibility, 
the conclusion of Remark \ref{Xdatumuniqueness}
applies, 
and  this permits us to further reduce
 to the case where $X_1=X_2^*$ and $X_2=X_1^*$.

 Hence, assume that $M=S$, and that the symplectic 
 induces a perfect pairing on $X_1 \times X_2$. 
 Fourier transform (obtained by integrating with reference to the self-dual measure the kernel
$\psi(\langle x, x^* \rangle)$) gives rise to an isomorphism of Schwartz spaces
 \[ \mathcal S(X_1(\mathbb A)) \rightarrow \mathcal S(X_2(\mathbb A)),\]
 The isomorphisms are equivariant with respect to the normalized $G(\mathbb A)$-action denoted by $\star$ in \S~\ref{normalizedperiod}.
 The functions  $\partial^{1/2} \cdot \Phi(x)$ 
in the notation of \eqref{PXform2}, for $X_1$ and $X_2$, are mapped to each other under Fourier transform, and the Poisson summation formula  implies the equality of normalized period functions.

  Let us now sketch independence of polarization for the period sheaf.    It follows from a geometric version of the previous argument,
as has been given in the work of Braverman and Gaitsgory \cite[Lemma 7.3.6]{BravermanGaitsgory},\footnote{We thank Tony Feng and Jonathan Wang for explaining this reference and argument. }
where it is related to the (sheaf-theoretic) functional equation for Eisenstein series. We note that this argument
is given in the finite case, where it uses properties of Artin--Schreier sheaves. Probably a version works in the other contexts too, but we did not check.
Specifically, the quoted Lemma is
to be 
applied to a $2$-term complex on $\Bun_G$ computing the cohomology of the vector bundle associated to $X$;
we can find such a complex at least on any quasicompact open substack, and the resulting isomorphisms 
can be glued by Lemma 7.3.7 (b) of the same reference.

\end{proof}

\subsubsection{Unpolarized periods} \label{remark:unpol}
Of course, we would like to define the period sheaf  and period functions without recourse to a polarization; in particular,
for $M$ that do not admit a polarization.   Now we discuss this unpolarized setting,  where, in short, 
 the ingredients all exist but a more detailed study of certain technical issues is required 
 to formulate them in the level of generality considered in this paper. 
 See \S \ref{Lindep} for the spectral counterpart of this discussion.

 In general terms, this should  follow from the theory of 
   $\theta$ series or its geometric analogue, Lysenko's geometrization of $\theta$-series \cite{Lysenko}.
   (A related topic is the recent construction of Coulomb branches for general symplectic representations~\cite{RaskinNisyros}.)
However, to carry this out in a way that is sufficiently detailed
for our needs, one needs
an analysis of the issue of splittings, which
should be closely related to the issue of the anomaly. 
 
Let us restrict, for what follows, to
 the case when $\FF$ has finite characteristic; we suppose 
(as we expect, see Expectation \ref{GMdesiderata expectation})
  that $(G, M)_{/\FF}$ arises
from the base change  to $\FF$ of some
  split hyperspherical datum, as in Definition \ref{ringdatum} , defined over some subring 
$R \subset \C$.

 In that case -- using the notations of Part \ref{part:structure} -- we begin
by constructing a $\theta$-function on $HU$ associated to the symplectic space $S \oplus \mathfrak{u}/\mathfrak{u}_+$;
  this depends on the choice of spin structure, through the choice of the additive character $\psi$. At the level of functions, 
Weil's theory constructs a ``Jacobi'' $\theta$-function 
in the space of automorphic functions for the semidirect product of the metaplectic group $\widetilde{\Sp}$
and the Heisenberg group on which it acts. Fixing a splitting of $\widetilde{\Sp} \rightarrow \Sp$ over $H$
permits us
to pull back this $\theta$-function to $HU$; we then use   the $\Theta$-series (summation over rational points)
associated to $HU \backslash G$ to construct an automorphic function on $G$. 
At the level of sheaves, the geometrization of the $\theta$-series 
(at least on $\widetilde{\Sp}$, rather
than the larger semidirect product by the Heisenberg group)
 was carried out by Lysenko~\cite{Lysenko}.

To carry this out one must have a splitting of 
the extension  $\widetilde{\Sp} \rightarrow \Sp$ 
over the ad{\`e}lic points of $H$. 
It is our hope that 
the vanishing of anomaly -- understood   in the sense of Definition \ref{anomalyfree}
applied to $(G, M)_{\C}$ --  should provide, in fact, a  distinguished splitting
(this is a reasonable hope at least  for   the distinguished split form postulated in  \S \ref{dhpFq}; otherwise
one may need to modify the anomaly vanishing condition to take into account issues of rationality).  
 The most favorable case is where $H$ is simply connected; in that case,
 the vanishing of anomaly for $(G, M)_{\C}$ as in Definition \ref{anomalyfree} implies that
the metaplectic cover of $\Sp(\mathbb{A}_F)$ splits uniquely over $H(\mathbb{A}_F)$,
in a fashion that is compatible with the splitting of the metaplectic cover on $F$-points:

\begin{itemize}
\item[(a)] The splitting of the cover can be deduced from \S \ref{anomalousautomorphic}.
Here one uses the fact that  $H$ is simply connected to pass
statements from an algebraically closed field to $\FF_q$
as in \cite[1.10]{Deligne}. 

\item[(b)] The resulting splitting is unique, because $H(\mathbb{A}_F)$
has trivial abelianization -- again, this uses that $H$ is simply connected.   
 \end{itemize}

In general -- that is to say, when $H$ is not simply connected -- both points
become less clear.  Lemma \ref{anomaly rationality} is a partial result 
in the direction of (a). More interesting, however, is (b): 
the question of {\em choice of splitting}.   In the classical theory of $\theta$ correspondence
Kudla  \cite{KudlaSplitting} has introduced a certain set of functional splittings, whose 
significance on the dual side is understood; what is needed is to abstract these examples. 

Our local conjecture suggests the following proposal for how to split
metaplectic covers over local fields:  For 
$V$ a representation of the dual group $\hat{H}$ with associated Hecke operator 
$T_V$, we  should have 
\begin{equation} \label{TVX} \langle T_V \delta_X, \delta_X \rangle \geq 0\end{equation}
when $\delta_X$ is a spherical
vector in the metaplectic representation. 
If such a splitting exists it is unique.

\begin{remark}
Implicitly, the condition \eqref{TVX} depends on a choice of $q^{1/2}$,
 which enters through the definition of the metaplectic representation;
 by default we take the positive choice.
 If we used its negative, 
  the splitting is modified  through $H(F) \stackrel{\theta}{\rightarrow} F^{\times} \stackrel{\textrm{val}}{\rightarrow} \Z/2\Z$, 
where $\theta$ is a character of $H$ as in Lemma \ref{calculateanomaly}. 
\end{remark}

 \section{$L$-functions and $L$-sheaves} \label{Lsheaf}
This section is the spectral analogue of \S \ref{section:global-geometric}: starting
with a hyperspherical $(\check{G}, \check{M})$, with $\check{M}$ polarized, 
 we will define an ``$L$-sheaf'' on the spectral side
 of the Langlands program, and 
 explain the sense in which it geometrizes an $L$-function. {\bf Note   that in this section we write the $\check{G} \times \GGm$-action on $\check{M}$ on the left. We
recall that the convention for passing from left to right actions is to invert the action of $\check{G}$
but not of $\GGm$, cf. \S~\ref{leftrightconventions}.}

\begin{itemize}
\item \S \ref{Lsetup} sets up some general notation supplementing that of \S \ref{notn}.

\item\S\ref{epsilon section} sets up notation on $\epsilon$-factors.
  
\item  \S \ref{locX}: we define
the space $\Loc_{\check{G}}^{\check{X}}$ of $\check{G}$-local systems with $\check{X}$-section, which is an analogue
of $\Bun_G^X$ defined previously. 

\item \S \ref{LsheafX}: we define the $L$-sheaf, which again comes in normalized and unnormalized forms;
the normalized form will be discussed in \S \ref{LsheafXnormalized}.

\item\S \ref{case4} defines the $L$-sheaf in the case of a {\em twisted} polarization, which is quite subtle
and requires the idea of shearing (\S \ref{shearingsec}). Roughly, the role of twisting
is to shift cohomological degrees in the $L$-sheaf.  In particular, the role of the Artin--Schreier sheaf on the automorphic is played on the spectral side by the spectral exponential sheaf introduced in~\S\ref{spectral exponential}.

\item \S \ref{Lspinindep} discusses the role of spin structures.

\item  \S \ref{locXex}: we compute fibers of the $L$-sheaf and show that these give
(geometrizations of) $L$-functions, thus the name ``$L$-sheaf.'' 

\item \S \ref{Linduction} explains the process of {\em spectral Whittaker} or {\em Arthur} induction, which can be used to reduce the study of L-sheaves to the vectorial case (parallel to \S \ref{Pinduction}). 

\item   \S \ref{Lindep}  explains independence of polarization, a categorified form of the functional equation for $L$-functions  (parallel to \S \ref{Pindep} but the computations
are less familiar).  
 \end{itemize}

  \subsection{Setup} \label{Lsetup}
 
 We will follow the general notation set up in \S \ref{notn}, but will fix some extra notation related to the spectral side.
 
  \subsubsection{Derived stacks}
Since the $L$-sheaf
 involves  algebraic rather than topological constructions, it is sensitive to
 the derived structures on the spaces involved.  
 Therefore,  although these ideas are similar to those of the last chapter, the level 
 of technicality involved in implementing them is greater.  
  The foundational
 theory is quite involved, and we will have to use it as a black box, most notably   
 the theory of quasi-coherent sheaves (and their variants, ind-coherent sheaves) on such spaces.
A standard reference for this material is the book~\cite{GR} of Gaitsgory and Rozenblyum.
To avoid being buried in a mountain of technicality, and to help preserve the
sanity of the authors, we will often take the liberty of either sketching certain constructions, 
or proving them under specific assumptions, with the understanding that
we expect their extension to other cases to be routine for the experts, and that
we will clearly flag any issues that do not seem to be straightforward.

 A {\em prestack} over $\kk$ means 
a functor of $\infty$-categories
from ``derived commutative rings'' to anima (the homotopy theory of simplicial sets, topological spaces or $\infty$-groupoids).
 There are different models for ``derived commutative rings''; since we work in characteristic zero it is convenient to take differential graded commutative rings which are {\em connective} (in degree $\leq 0$
 with degree-increasing differential). A {\em derived stack} is a prestack satisfying a sheaf condition.
 Most of our derived stacks
 will be in fact quotients of derived schemes by an affine algebraic group.

We briefly recall  informally some features of sheaf theory, and refer
to  \S \ref{coherent sheaf theories} in the Appendix for more details.  A ``quasicoherent sheaf'' on a derived stack means, informally, a compatible system of quasi-coherent
 sheaves on affines $\mathrm{Spec}(A) \rightarrow X$ mapping to $X$ -- i.e., the $\infty$-category of quasi-coherent sheaves is defined as the limit of (the $\infty$-derived category of) $A$-modules over all affines over $X$.
  On singular schemes and stacks $X$ (such as the stacks of local systems arising as Langlands parameters) it's crucial to enlarge the category $QC(X)$ of quasi-coherent sheaves to that of {\em ind-coherent} sheaves $QC^!(X)$, which account more fully for the singularities (by replacing the role of perfect complexes by that of arbitrary coherent complexes). See~\cite[Part II]{GR} for a detailed study. The definition of $QC^!$ is more subtle than that of $QC$, and in particular requires that $X$ satisfy finite type assumptions. See \S \ref{IndCoh section} for a summary. 

The theory of ind-coherent sheaves is the natural home for Serre duality on singular spaces, and in particular the dualizing sheaf $\omega_X\in QC^!(X)$ is naturally ind-coherent, as are the $L$-sheaves we introduce in this section as pushforwards of dualizing sheaves. These sheaves may lose crucial information if we try to project them to $QC$ (in fact this projection {\em vanishes} in the presence of a nontrivial Arthur $SL_2$.) 
We will use the $!$-tensor product structure $\otimes^!$ on $QC^!$ (for which $\omega$ is the unit) but also the tensor product action of $QC$ on $QC^!$, which we denote by a plain $\otimes$. 
  
 \subsubsection{$\Loc_{\check{G}}$ in the different contexts} \label{BunLocintro2}

We continue the discussion of \S \ref{BunLocintro}, now on the spectral side,
and again pointing to Appendix \S \ref{geometric Langlands} for details.
In all cases
 $\Loc_{\check{G}}$ will be a derived  stack over $\kk$.

 \begin{itemize}
 \item  Finite context: $\FF$ is the algebraic closure of a finite field and $\kk=\overline{\mathbb{Q}_{\ell}}$.     The space $\Loc_{\check{G}}$ is taken to be the space of ``restricted local systems'' (or local systems with ``restricted variation''), 
 defined as a pre-stack in \cite[\S 1.3]{AGKRRV1}. It classifies a certain class of 
  representations of the {\em geometric} fundamental group of $\Sigma$
 and in particular comes with a Frobenius action, when the curve is defined over $\FF_q$. (Again, any mention of Frobenius in the text will, naturally, assume this.)   In {\em op. cit.} Theorem 1.3.2
 various geometric properties are given; it is in particular ``locally of finite type'' and one can talk
 of ind-coherent sheaves as in   \cite{GR}.   
   
   {\bf Warning:}  This situation comes with technical details not encountered in the situations below. We have not examined
   these issues in detail but will flag
   them at relevant points in the text, e.g.
   (ii) of \S \ref{locX}.
 
  \item De Rham context: $\FF=\kk=\C$.  We take $\Loc_{\check{G}}$ to be the space
   of de Rham local systems (i.e., $\Gv$-bundles with flat connections); this is defined as a mapping stack (from the de Rham space $\Sigma_{dR}$ to $B\Gv$) in \cite[10.1.1]{ArinkinGaitsgory}
   and also studied in \cite[\S 2]{BD}.
   As observed in the later reference \cite[2.11.2]{BD}, if the genus is $\geq 2$ and $\check{G}$ is semisimple, this is a classical
stack, i.e., it is representable by the quotient of a (non-affine!) underived scheme by $\check{G}$. 

 \item Betti context: $\FF = \C, \kk =$ an algebraically closed field of characteristic zero.   Here, $\Loc_{\check{G}}$
  is  the space of Betti $\check{G}$-local systems on $\Sigma$, which can be again defined as a mapping stack, now from the homotopy type $\Sigma_{Betti}$ of $\Sigma$ to $B\Gv$.
In the case of genus $\geq 1$ (see \S \ref{P1} for the other case) we fix a basepoint and can consider this as a 
  space parameterizing $\pi_1(\Sigma)$-representations, where we retain stack and derived structure;
  it can therefore be presented for a curve of genus $g$ as the conjugacy quotient of the space
  of representations of $\pi_1$ into $\check{G}$ 
\begin{equation} \label{locGprez} \Rep_{\check{G}} :=  \{ (x_1, y_1, \dots, x_g, y_g) \in \check{G}^{2g}: [x_1, y_1] \dots [x_g, y_g]=e\}.\end{equation}
  Observe that its algebraic structure does not depend on the algebraic structure of
  $\Sigma$.  As in the de Rham case, if $g \geq 2$ and $\check{G}$ is semisimple, 
  $\Rep_{\check{G}}$ is in fact a locally complete intersection affine ring, and $\Loc_{\check{G}}$ is a classical(=underived) Artin stack.\footnote{This 
  can be deduced
  from the corresponding
  assertion in the de Rham case. }  \end{itemize}

  \begin{remark}
The space $\Loc_{\check{G}}$ of local systems with restricted variation that
we use in the finite context in fact makes sense for any field, and for $\FF=\C$ sits inside both the Betti and de Rham spaces of local systems (see \S \ref{spectral side section}). This space only sees formal neighborhoods of irreducible representations. In general, {\em restricted variation} means the semisimplification of a local system is fixed in any family, see~\cite[0.5.3]{AGKRRV1}. Still, this is sufficient to 
compare with numerical predictions, and, more to the point, there is no known way to go beyond formal neighborhoods of semisimple parts in the finite setting.   \end{remark}

 \subsubsection{Cohomology of $\Sigma$} \label{sigmacoh}
 
 We understand $H^*(\Sigma, -)$ to mean singular cohomology of a Betti local system,
 de Rham cohomology of a de Rham local system, or (geometric) {\'e}tale cohomology
 of an {\'e}tale sheaf, according to context. (Recall ``geometric'' means that, even if our curve is defined over $\FF_q$, we base change to $\FF=\overline{\FF_q}$.)

 \subsubsection{The Frobenius action on $\Loc$} \label{Frobaction}
  In the finite context (with $\Sigma$ defined over $\FF_q \subset \FF=\overline{\FF_q}$),
there is an action of Frobenius on $\Loc_{\check G}$:
\begin{equation} \label{LocFr} \Fr: \Loc_{\check G} \rightarrow \Loc_{\check G}.\end{equation}
We will write it out  as part of our running battle with   signs. 
We will understand this to be defined 
by means of (equivalently):
\begin{itemize}
\item pullback of {\'e}tale sheaves by the 
morphism 
\begin{equation} \label{arithFrobenius} \mathrm{id} \otimes (\lambda \mapsto \lambda^q)^*\end{equation} 
on 
$\Sigma_{\overline{\FF}} = \Sigma \times_{\Spec \  \FF} \Spec \  \overline{\FF}$, or  
\item The {\em inverse} to the pullback  by geometric Frobenius acting on {\'e}tale sheaves.
Recall that geometric Frobenius is 
  the morphism which raises coordinates to the $q$th power
with respect to a fixed $\FF$-projective embedding.
\end{itemize}

A  $\kk$-point of $\Loc_{\check{G}}$ fixed by this action
amounts to  giving a $\check{G}$-local system $\rho$ with $\kk$
 coefficients equipped with an isomorphism $\Fr \ \rho \simeq \rho$.
Suppose in addition that $\mathcal{F}$ is a Frobenius-equivariant coherent sheaf on $\Loc_{\check{G}}$. 
  Then there is an induced ``Frobenius''
 on the $\rho$-fiber  $\mathcal{F}_{\rho}$. 
By definition, we understand this to mean 
the composite
\begin{equation} \label{compo} \mathcal{F}_{\rho} = \mathcal{F}_{\Fr \rho} =  \Fr^* \mathcal{F}_{\rho} \stackrel{\sim}{\rightarrow} \mathcal{F}_{\rho}, \end{equation} where the first map uses the Frobenius-equivariant structure on $\rho$, i.e.
the structure that renders it fixed by Frobenius, and the final map uses
the Frobebius-equivariant structure on the sheaf $\mathcal{F}$; moreover, in \eqref{compo}, $\Fr$
is as in \eqref{LocFr}, and should not be confused with the geometric Frobenius on $\Sigma$ itself. 

 For example, if $\check{G}=\GL_n$, the 
  the cohomology (or rather cochains) of the $n$-dimensional local system 
  associated to each $\rho \in \Loc_{\GL_n}$ 
  can be regarded as the fibers of a certain coherent sheaf $\mathcal{F}$ on $\Loc_{\check{G}}$. 
  With our conventions, 
  the action of Frobenius on $\mathcal{F}_{\rho}$
  is naturally identified with the pullback action of geometric Frobenius
  on $H^*(\Sigma_{\bar{\FF}}, \rho)$.

\subsubsection{Tate twists} \label{Tatetwistrecall} 
\index{angle bracket twists $\langle d \rangle$}
As explained in \eqref{ultimate shearing}, $\langle 1 \rangle$ denotes 
the simultaneous application of the following three shifts:
\begin{itemize}
\item a cohomological shift $[1]$;
\item a Tate twist by $1/2$, where applicable (e.g., if we are dealing with Frobenius equivariant objects);
\item a change of parity, where applicable  (i.e., if we apply it to a sheaf, we regard that sheaf
as a super-sheaf and change its parity). 
\end{itemize}

\subsubsection{Conventions for arithmetic class field theory} \label{CFTGLconventions}  
 \label{CFTnormalization}
 We will use class field theory, i.e., the Langlands correspondence for $\GL_1$,
 and  again we will {\em try} to get signs right,  for which
 reason we briefly recall it here.  Suppose we are in the finite context. 
We normalize local and global class field theory so that geometric Frobenius elements  in the Galois group are carried to uniformizers in
the local fields or adeles. 
In this version, the cyclotomic character of a nonarchimedean local field with residue
characteristic $q$, which sends geometric Frobenius to $q^{-1}$,
is matched with the normalized valuation character
 $x \mapsto |x|$ (sending a uniformizer to $q^{-1}$).  Globally the  cyclotomic
 character is also matched with $x \mapsto |x|_{\mathbb A^\times}$.
  This convention coincides
 with Tate's in \cite{Tate}. As discussed in \S \ref{periodX}, the adelic uniformization
\begin{equation} \label{adelicuniformization} \mathbb{A}^{\times} \rightarrow \mathrm{Pic}(\Sigma)(\FF)\end{equation}
carries the adele corresponding to a uniformizer $\varpi_x$
at a closed point $x \in \Sigma$, to the line bundle $\mathcal{O}(x)$.

 \subsubsection{Conventions for geometric class field theory} 
 \label{GLGMnormalization} \index{geometric class field theory} \index{GCFT}
 Continue in the setting of \S \ref{CFTnormalization}, i.e., we are in the finite context.  \index{$\Pic$ vs. $\Bun_{\Gm}$}
To each $\Gm$-local system $\rho$ on $\Sigma$ is associated a Hecke eigensheaf $\chi_{\rho}$ on $\Bun_{\Gm}$,  which 
  \footnote{We use $\mathrm{Bun}_{\Gm}$ 
   as opposed to $\mathrm{Pic}_{\Gm}$
to emphasize that we are interested also in the stacky aspect of its structure - the notation $\mathrm{Pic}$
is often used to refer to a scheme.} 
for $r \gg 1$  descends from  $\rho^{\boxtimes r}$ via $\Sigma^r \rightarrow \mathrm{Bun}_{\mathbb{G}_m}$,
 $(P_1, \dots, P_r) \rightarrow \mathcal{O}(\sum P_i)$. 
 In particular, the function $\Bun_{\Gm}(\mathbb{F}_q) \rightarrow \kk^{\times}$,  attached to $\chi_{\rho}$ (the trace of geometric Frobenius), 
after pullback to the id\`eles $\mathbb{A}^{\times}$, gives the id\`ele class character
associated by class field theory to $\rho$.

 The construction just described gives rise, more generally, to an equivalence of categories:  
 \begin{equation} \label{GCFT} \mathsf{GCFT}: \bigautshvspec(\Bun_\Gm) \simeq QC^{!}(\Loc_\Gm)
\end{equation}
which holds in all contexts, but with appropriate conditions for sheaf theory on both sides. In the de Rham setting, the equivalence $\mathsf{GCFT}$ is a mild extension of the Fourier-Mukai transform of Laumon~\cite{LaumonFourier} and Rothstein 
\cite{Rothstein}
 identifying $D$-modules on the Jacobian of $\Sigma$ with quasicoherent sheaves on the moduli scheme of flat line bundles on $\Sigma$; in the Betti case~\cite[4.3]{BettiLanglands} it is a simple consequence of the identification of the first homology of $\Sigma$ and the fundamental group of its Jacobian. The {\'e}tale version is likewise elementary but currently missing an explicit reference. 
See the discussion of Conjecture \ref{GLC} for more details on the geometric Langlands conjecture of which this is the $GL_1$ case.\index{$\mathsf{GCFT}=$geometric class field theory equivalence}

We take the opportunity to fix some signs. The Hecke operator $T_x$ at 
  $x \in \Sigma$ is ``translation by $x$'' 
arising the map $D \mapsto D+x$ on divisors. It acts by pullback on sheaves, thus  
sending a sheaf supported on $\mathrm{Pic}^0$ to a sheaf supported on $\mathrm{Pic}^{-1}$.
The action of $T_x$ on $\chi_{\rho}$
corresponds, on the right hand side of \eqref{GCFT}, 
to  tensoring with the fiber $\rho_x$ of $\rho$ at the point $x$.

 \begin{remark}  
  As a (rather minor) warning,  while this normalization of geometric class field theory  
 is (up to sign) the standard one, it will not coincide with the normalization of the Langlands
 correspondence for $G=\Gm$ posited in the global
 geometric duality conjecture Conjecture \ref{GlobalGeometricConjecture}.
Including these twists here would be needlessly heavy for the minor way in which we use it. 
 \end{remark}

\subsubsection{$\check{G}$ and $\check{M}$}  \label{Lsheafgeneralsetup}

$\check{G}$ and $\check{M}$ will be   hyperspherical over $\kk$;
since  $\kk$ is algebraically closed
of characteristic zero, there are no issues of rationality to consider -- this is defined 
as in Part \ref{part:structure}.

 We will moreover restrict in the current section 
to the case that $(\check{G}, \check{M})$ is polarized, possibly with a twisting (see the end of \S \ref{Lindep}
and further discussion for discussion of the general case).  That is to say,

$$ \check{M}=T^*\check{X} \mbox { or } T^*(\check{X}, \check{\Psi}),$$
where $\check{X}$ is a $\Gv \times \GGm$-space, and, if applicable,
$\Psi$ is a $\Gv \times \GGm$-equivariant $\Ga$-torsor over $\check{X}$,
where $\Gv$ acts trivially on $\Ga$ and $\GGm$ acts by squaring on $\Ga$. 
\index{left versus right actions}

Finally,  we will allow ourselves
to assume that $\check{X}$ has an eigenform. 
As previously discussed (see \S \ref{ssseigencharinocuous}, and 
 \S \ref{periodsheafgeneralsetup} for the corresponding computation for period sheaves) this
should be considered a matter of convenience.

\subsection{Epsilon factors}\label{epsilon section} 
In this section we set up basic notions regarding $\epsilon$-factors that will be used to define normalized $L$-sheaves.
 
 \subsubsection{Recollections on $L$ and $\epsilon$-factors}  \label{Leps}\index{functional equation for $L$-function}
  Now, restrict to the case of $\Sigma$ defined over $\FF_q$,
  and let $T$ be an {\'e}tale local system of $\kk$-vector spaces on $\Sigma$.  Recall that we denote by $\Gamma$ the Weil group of the function field of the curve.
   In this setting we have an $L$-function and
  an $\epsilon$-factor 
  $$ L(s, T)\mbox{ and }\epsilon(s,\psi,T) \in \kk(q^s),$$
   defined using  
  $\psi$ as in \eqref{psidef}, which, we recall, depends in particular
  on a fixed spin structure on $\Sigma$.  We follow the conventions of number theory
  in writing this as a function of $q^s$, although it would be more reasonable in our current
  situation to treat $q^s$ as a formal variable. 
  The spin structure being fixed,  we abridge $\epsilon(s, \psi, T)$  to $\epsilon(s, T)$.

  \index{$\epsilon$-factor}
  We understand the $\epsilon$-factor to
  be as defined by Tate \cite[\S 3.6]{Tate} taking the measure $dx$ therein to be self-dual Haar measure. 
 We then have \index{$\epsilon(s, T^{\shear})$} 
\begin{equation} \label{epsform} \epsilon(s, T) = \det(T) (\partial) \cdot q^{(1/2-s) (2g-2) \dim T},\end{equation}
where, in writing $\det(T)(\partial)$, we identify the determinant of $T$ with an id{\`e}le class character via class field theory,
and the element $\partial$ was defined after \eqref{diffdef}.
  We will also use notation such as $\epsilon(0, T^{\shear})$ in a way similar
to the use \eqref{Lsheardef} for $L$-functions. Thus,  if $T$ is a representation of $\Gamma \times \GGm$
 then we set $$\epsilon(s, T^{\shear}) := \prod_k \epsilon(s+k/2, T_k)$$ with $T_k$
 the $k$th graded component.

Now, our choice of a square root of the different distinguishes a square root $\sqrt{\epsilon}(s,T)$
  associated to the spin structure, namely  
\begin{equation} \label{epsilondef} \sqrt{\epsilon}(s, T)  = (\det T)(\partial^{1/2}) \cdot q^{(1/2-s) (g-1) \dim T}.\end{equation}
   We have 
 $ L(s,T) = \epsilon(s,T) L(1-s, T^\vee) \mbox{ and } \epsilon(s, T) \epsilon(1-s, T^\vee) = 1$. In  particular the ``normalized $L$-function''  \index{$L^{\norm}$, normalized $L$-function}
\begin{equation} \label{Lnormdef}  L^{\norm}(s,T) := L\epsilon^{-1/2}(s,T) := \epsilon(s,T)^{-1/2} L(s,T)\end{equation}
  is actually invariant under $(T,s) \leftrightarrow (T^\vee,1-s)$, 
\begin{equation} \label{LnormFE}  L^{\norm}(s,T) = L^\norm(1-s, T^\vee).
\end{equation}

  \begin{remark}
  Let us note that the Grothendieck--Lefschetz trace formula interprets $L(s, T)$ as the alternating product
  $\prod_{i} \det(1-q^{-s} \Fr | H^i)^{(-1)^{i+1}}$, with $H^i$ the cohomology of the geometric curve
  $\Sigma_{\overline{\mathbb{F}_q}}$ with coefficients in $T$,  and $\mathrm{Fr}$ the  pullback action of geometric Frobenius, 
  which gives us the interpretation
\begin{equation} \label{firstepsilondet} \epsilon(0,T) = \prod_{i} \mathrm{det}( \Fr | H^i)^{(-1)^{i+1}},\end{equation}
that is to say, the inverse of the trace of Frobenius on the determinant of cohomology. 
  \end{remark}

\subsubsection{Geometric class field theory for line bundles}
We will next  introduce the analogue of this $\sqrt{\epsilon}(s,T)$ in the geometric setting. It will
 be a line bundle on $\Loc_{\Gm}$ denoted by $\varepsilon_{1/2}$, depending on a choice of spin structure $\cK^{1/2}$ (see \S \ref{Lspinindep} for a
 formulation independent of this choice).
 Although our numerical discussion of $\epsilon$ was restricted to the finite context, $\varepsilon_{1/2}$
 will be defined in all contexts.
 
 The construction is based on the following basic feature of geometric class field theory, which we describe from several perspectives:

\begin{prop}  There is a unique homomorphism of groups\footnote{In this statement, $\mathrm{Pic}$ is to
be understood as a ``mere'' abelian group, i.e., isomorphism classes of line bundles, without additional algebraic structure: this is all we need.} (to be denoted $\mathcal{L} \mapsto \mathcal{[L]}$)
\begin{equation} \label{Lprimedef} \mathrm{Pic}( \Sigma) \longrightarrow \mathrm{Pic}( \Loc_{\Gm}), \ \ \mathcal{L} \mapsto [\mathcal{L}].\end{equation} 
for $\mathcal{L}$ a line bundle on $\Sigma$, with the property that  $\mathcal{O}(x)$ is sent to
 the line bundle whose fiber
 over $\rho \in \Loc_{\Gm}$ is the fiber of $\rho$ at $x$.
 
 \index{$\mathcal{[L]}=$class field theory for line bundles}
 
 \end{prop}

 Indeed, the definition uniquely specifies what $[\mathcal{O}(D)]$ is,
 and one checks independence of $D$ using the fact that, for fixed $\mathcal{L}_0$,
 the space of $D$ for which $\mathcal{O}(D) \simeq \mathcal{L}_0$ is a projective
 space and in particular simply connected if nonempty.  
 The map is easy to describe concretely in each context. 
  \begin{itemize}
  \item[--]  In the de Rham setting this homomorphism comes from the  
  pullback of the Poincar\'e line bundle (expressing the self-duality of $\Pic$) under the projection  $\Loc_{\Gm} \rightarrow \Bun_{\Gm}$ from rank one flat connections to degree zero line bundles.

 \item[--] In the Betti context the bundle $[\mathcal{L}]$
 is obtained by taking the bundle associated to the representation
 $z \mapsto z^{\deg\mathcal{L}}$ on $B\Gm$
 and pulling it back via the map $\Loc_{\Gm} \rightarrow B\Gm$ given by taking fiber at a fixed point of $\Sigma$

 \item[--]  In the finite context, the bundle $[\mathcal{L}]$
 is again  pulled back from $B\Gm$ just as in the Betti case.\footnote{This does not depend on the 
 point at which we take fiber.}
 \end{itemize}

 Note that in the Betti and finite cases the class of  $[\mathcal O(x)]$ is actually independent of $x$.

In the finite case, if we suppose that $\mathcal{L}$ to be defined over $\mathbb{F}_q$,
 then  there is a natural way to equip $[\mathcal{L}]$ with the structure
of Frobenius-equivariant line bundle on $\Loc_{\Gm}$; 
if, for example, $\mathcal{L}=\mathcal{O}(x)$
 for $x \in \Sigma(\mathbb{F}_q)$, then
 this arises from the tautological identification of
 the fibers of $\rho$ and $\mathrm{Fr}(\rho)$ at $x$. 
 
For $\rho$ a  $\Gm$-local system 
defined over $\FF_q$, in particular defining a Frobenius fixed point on $\Loc_{\Gm}$, 
we  have for this equivariant structure
 $$ \mbox{ 
  trace of geom. Frobenius on $[\mathcal{L}]_{\rho}$} =\chi_{\rho}(\mathcal{L})
 $$
 where $\chi_{\rho}: \Pic(\Sigma) \rightarrow \kk^{\times}$
 is the character associated to $\rho$ by class field theory, normalized as in
 \S \ref{CFTGLconventions}.  For example, if $\mathcal{L}=\mathcal{O}(x)$
 for $x \in \Sigma(\mathbb{F}_q)$, then 
  trace of geometric Frobenius on $[\mathcal{L}]_{\rho}$
is then the trace of geometric Frobenius on the fiber of $\rho$ at $x$, i.e.
by the normalization of class field theory (\S \ref{CFTGLconventions}), $\chi_{\rho}$ evaluated
at the uniformizer $\pi_x$, which by \eqref{adelicuniformization}
uniformizes $\mathcal{L}$.

 \begin{remark}[Construction via geometric class field theory] \label{GCFTpicconstruction}
 A more structured approach to the homomorphism $[-]$ is given by 
  categorical geometric class field theory
-- that is to say, the  abelian case of the geometric Langlands correspondence.

 We discuss first the de Rham situation. 
 In that case,  with our normalizations (cf. \S \ref{GLGMnormalization})
 to each $\Ll\in \Pic(\Sigma)$  we may consider the  skyscraper sheaf $i_{\mathcal{L}^{-1},*}\kk$
 at $\mathcal{L}^{-1}$ (the inverse is an artifact of our normalizations).
 We 
  take its image under 
 geometric Langlands correspondence
 \eqref{GCFT}, 
 one version of which is   an equivalence  
 of symmetric monoidal categories\footnote{Here, in contrast to \eqref{GCFT}, it is convenient to take the \unrenormalized category
 on the left, and quasi-coherent rather than ind-coherent sheaves on the right.}
  $$(\D(\Bun_{\Gm}),\ast)\simeq (QC(\Loc_\Gm),\otimes).$$ 
where  the monoidal structure on the source
 is given by convolution; 
 and the skyscraper is an invertible
 object of the source with respect to this monoidal structure.  Therefore, the resulting
 sheaf on $\Loc_{\Gm}$ is an invertible object of the category of quasi-coherent sheaves, i.e.,
 a line bundle; this is just our $[\mathcal{L}]$. 

 In the Betti and \'etale settings, we must first apply the spectral projection before applying geometric Langlands. 
 We will discuss the notion of spectral projection at more length in
 \S \ref{nilpotent projection is good}. For example in the Betti setting, we simply replace $\Pic$ by its homotopy type, 
and replace $D$-modules on $\Pic$ by local systems. In particular we replace the skyscraper $i_{\mathcal{L}^{-1},*}\kk$ by the corresponding ``universal cover local system'', where the inclusion
of $\mathcal{L}^{-1}$ is replaced by the path fibration at $\{\mathcal{L}^{-1}\}$.
 \end{remark}

 \subsubsection{The spectral bundle $\varepsilon_{1/2}$} \label{spinstructure-spectral} 
 \index{half-epsilon}
 \index{$\varepsilon_{1/2}$}
 \begin{definition}
The {\em half-epsilon line bundle} associated to a choice of spin structure $\cK^{1/2}$  is defined as
 $$\varepsilon_{1/2} = \mbox{ line bundle $[\cK^{1/2}]$  on $\Loc_{\Gm}$ associated to $\cK^{1/2}$ via \eqref{Lprimedef}.}
  $$

  \end{definition}
  
    Having fixed a half-different $\mathfrak{d}^{1/2}=\sum_{v\in \Sigma} \frac{n_v}{2} v$, we have an explicit description 
\begin{equation} \label{epsilondescription} 
 \textrm{fiber of $\varepsilon_{1/2}$ at $L$} \simeq  \bigotimes_{v \in \Sigma} L_v^{n_v/2} \end{equation} 
  We will sometimes write, for $T$ a local system of vector spaces,  \index{$\varepsilon_{1/2}(T)$}
 $$\varepsilon_{1/2}(T) := \mbox{fiber of $\varepsilon_{1/2}$ at $\det(T)$.} $$

The bundle $\varepsilon_{1/2}$ categorifies
 the square root of a {\em central} $\varepsilon$ factor, i.e., a $\epsilon$-factor as in \S \ref{Leps}
 evaluated at its center of symmetry $s=1/2$. 
On the other hand, the square root of an  $\epsilon$-factor at $s=0$
is categorified by
  \begin{equation} \label{epsilon0description}  \varepsilon_{1/2}(T) \langle (1-g) \dim(T)\rangle.\end{equation}
  From our comments above, in the finite context, the trace of geometric Frobenius on the line $\varepsilon_{1/2}(T)$ resp. $ \varepsilon_{1/2}(T) \langle (1-g) \dim(T)\rangle$ is given by
 $\sqrt{\epsilon}(\frac{1}{2}, T)$ resp. $\sqrt{\epsilon}(0, T)$, where the choice of square root is determined as in 
 \eqref{epsilondef} via the spin structure.

 \begin{remark}   \label{det of coh remark}  
These constructions are closely related to determinant of cohomology,
cf. \eqref{firstepsilondet}: $\varepsilon_{1/2}$ is, up to a global twist by a line,
identified with the inverse square root of the determinant of cohomology on $\Loc_{\Gm}$.

Let us explicate this as directly as possible in the de Rham case.
Take a  vector bundle $E$
equipped with flat connection. Computing
its cohomology by the de Rham complex, we
find that the determinant of cohomology $\mathcal{D}$
for the associated local system is then the product
of the determinant of cohomology for $E$, and the inverse of the determinant
of cohomology for $E \otimes \Omega^1$:
$$ \mathcal{D} \simeq \det H^*(E) \otimes  \left( \det H^*(E \otimes \Omega^1)  \right)^{\vee},$$
where on the right we have coherent cohomology.   Fixing a
rational section $s$ of $\Omega^1$
with divisor $\sum_{i} n_i P_i$
gives rise to an isomorphism  of line bundles
$\mathcal{O}(\sum n_i P_i) \stackrel{\times s}{\longrightarrow} \Omega^1$. 
We compute $\det(E \otimes \Omega^1)$ by repeatedly using the exact sequence
$E \rightarrow E(p_i) \rightarrow E_{p_i} \otimes \mathcal{K}_{p_i}$, where $E_{p_i}$ is the stalk
of $E$ at $p_i$, and $\mathcal{K}_{p_i}$ the stalk of the canonical bundle.
This leads to an identification  $\mathcal{D}  \simeq \bigotimes (\det E_{p_i})^{ -n_i} \otimes \ell$,
where the line $\ell$ depends on the tangent spaces at the various $p_i$s,
but not on $E$  -- in fact, it is  identified with
the determinant of cohomology for the trivial local system.
Therefore, if all $n_i$s are even,
the line $\otimes_{i} (\det E)^{-n_i/2}$ gives a square root of $\mathcal{D}$, 
at least up to a line that depends only on the curve and choice of section $s$,
but not on the local system.
\end{remark}

\index{spin gerbe}

\subsection{$\Loc_{\Gv}^{\check{X}}$ and the $L$-sheaf}  \label{locX}
Let $\check{X}$ be a $\check{G} \times \GGm$-space as in \S \ref{periodsheafgeneralsetup}.
 The discussion that follows will largely apply unchanged to a general $\check{G} \times \GGm$-space
 $\check{X}$ but, with very rare and clearly noted exceptions, we will only use it in
 the narrow situation just quoted, that is to say, spaces derived from hyperspherical varieties. 

We define the spectral analogue of ``$G$-bundles with $X$-section,'' namely, 
we define  
 $\Loc_{\check G}^{\check X}$ to be the moduli space of $\check G$-local systems
together with a (flat or locally constant) section of the associated $\check X$-bundle.
 In all cases,  if $\check{X}=\check G_X\backslash \check{G}$, then this space $\Loc_{\check{G}}^{\check{X}}$ will be
      $\Loc_{\check G_X}$ together with its natural map to $\Loc_{\check{G}}$.

Before giving a precise definition we note that in the Betti and de Rham cases the geometry is quite tame:
 $\Loc_{\check{G}}^{\check{X}} \rightarrow \Loc_{\check{G}}$ is  a global quotient by $\check{G}$ of a  morphism of derived schemes. 
    In the finite context the morphism remains schematic, but the geometry of $\Loc_{\check{G}}$
  itself is complicated, as we have already noted.

\begin{itemize}
\item[(i)]  
We will discuss the Betti context in the most detail, as it is most explicit.  

One can construct $\Loc_{\check G}^{\check X}$ as a mapping stack in derived algebraic geometry (as in \cite[Appendix B]{ArinkinGaitsgory})  from the Betti space associated to $\Sigma$ (its homotopy type) to $\Xv/\Gv$:
$$\Loc_{\check G}^{\check X}=\Map(\Sigma_{\textrm{Betti}}, \Xv/\Gv)\longrightarrow \Loc_\Gv=\Map(\Sigma_{\textrm{Betti}},\textrm{pt}/\Gv).$$

This $\Loc_{\check{G}}^{\check{X}}$ is in fact representable as the quotient of an affine derived scheme by an algebraic group, as we shall explicitly sketch now
in the case when the genus of $\Sigma$ is $\geq 2$ and $\check{G}$ is semisimple (for the case when $\Sigma$ has genus zero, see Remark \ref{P1locprez}). 

Let $\Rep_{\check{G}}$ be the space of representations of $\pi_1(\Sigma)$
into $\check{G}$. Specifically, after fixing
a basepoint $\star \in \Sigma$,
homomorphisms $\pi_1(\Sigma, \star) \rightarrow \check{G}(S)$
for $S$ a $\kk$-algebra correspond
to $S$-points on the fiber
of the ``product of commutators'' mapping
$ \check{G}^{2g} \longrightarrow \check{G},$
considered as a map  of $\kk$-varieties.
We take $\Rep_{\check{G}}$ to be this fiber; {\em a priori}
one this could be taken as a derived scheme, but since we suppose that $\check{G}$ is semisimple
and the genus is $\geq 2$, 
it coincides with a usual scheme. Let $R$ be the ring of functions on the affine $\kk$-variety $\Rep_{\check{G}}$;
then we get a universal representation $\pi_1(\Sigma) \rightarrow \check{G}(R)$.

 The product  $\check{X} \times \Rep_{\check{G}}$
 correspondingly carries an action of $\pi_1(\Sigma)$, using 
 this universal representation, and we can consider the derived fixed points, 
i.e., we take derived $\pi_1(\Sigma)$-coinvariants on the ring of functions
$\mathcal{O}[\check{X} \times \Rep_{\check{G}}]$.\footnote{
The model example is as follows: 
$\pi_1$ acts on a vector space $\check{X}$, then the derived invariants
of $\pi_1$ on the ring of functions $k[\check{X}]$ will be the symmetric
algebra on the complex computing homology of $\pi_1$ with coefficients in $\check{X}$.
An explicit model for derived invariants in general can be given as follows: 
Use the equivalence between commutative connective differential graded $k$-algebras
and simplicial commutative algebras over $k$, then pass to a ``free resolution''
  replacing $\mathcal{O}$ with 
  a cofibrant simplicial $k$-algebra,  and then finally passing to $\mathcal{O}^{\otimes E \pi_1}$,
where $E \pi_1$ is the usual contractible simplicial set with $\pi_1$-action, and
taking levelwise $\pi_1$-coinvariants. }

Call the resulting scheme $(\check{X} \times \Rep_{\check{G}})^{\pi_1}$;
 and explicitly we may present  $$\Loc_{\Gv}^{\check{X}} := (\check{X} \times \Rep_{\check{G}})^{\pi_1}/\check{G}$$
 as the quotient of the resulting derived scheme  by $\check{G}$. 

 An even more explicit presentation 
when $\check{X}$ is a vector space will be given in \S \ref{Lindep}.

The fiber of $\Loc^{\check X}_{\check G}$ above
 a classical local system $\rho: \pi_1(\Sigma) \longrightarrow \check G(\kk)$ 
 is then the ``space of flat sections'' of $\rho$. When the genus of $\Sigma$ is $\geq 1$, 
 so that $\Sigma$ is topologically a $K(\pi_1, 1)$,  this can be equivalently
 described as  the  {\em derived} locus of fixed points
 of the $\pi_1$-action of $\rho$ on $\check X$.  These ``derived fixed points'' represent  
 the functor sending a   
 differential graded $k$-algebra $A$
 to the  homotopy fixed points of $\pi_1$ acting on the simplicial set $\check{X}(A)$
 by means of $\rho$.

\item[(ii)]   In the {\'e}tale context 
the definition of $\Loc$ is more complicated. We will define $\Loc_{\Gv}^{\check{X}}$
only in the homogeneous case (i.e., trivial vectorial part) where, if we have $\check{X}= \check{G}/\check{G}_X$,
we take $\Loc_{\Gv}^{\check{X}} = \Loc_{\check{G}_X}$ with its natural map to $\Loc_{\check{G}}$.

We do not expect any essential difficulty in transposing definitions similar to (i) or (iii) to this setting, which has been studied by \cite{AGKRRV1}, 
in greater generality,   but  we have not checked the  details.

 \item[(iii)]
 In  the de Rham context, the space 
of flat sections is similarly  constructed as a mapping space replacing the role of $\Sigma_{\textrm{Betti}}$ above by the
de Rham functor $\Sigma_{\textrm{dR}}$
 so that the fiber over a flat $\Gv$-bundle is the space of flat sections of the associated flat $\Xv$-bundle.

 Again, this stack is representable as the quotient of a derived scheme by an algebraic group.
The construction of this scheme can be carried out using the theory of jet schemes (thought of as commutative $\D_\Sigma$-algebras or equivalently commutative chiral algebras) which we explain pointwise: 
  associated to $\check{X}$ and a de Rham $\check{G}$ local system $\rho$ 
 is a commutative $\D_\Sigma$-algebra, the sections of the jet scheme for the associated flat $\check{X}$-bundle,
 and the construction of spaces of horizontal sections (the commutative case of chiral homology, see  \cite[\S 2.4.1 and \S 4.6]{BD}) provides a dg ring representing
 the space of flat sections of $\rho$.   We expect this discussion to extend
 to families without difficulty.
  
\end{itemize}

 The definition as mapping stack in the Betti and de Rham context formally implies a computation of the tangent complex, as in \cite[Appendix B]{ArinkinGaitsgory}. 
 We will use the following case: for a $\kk$-point of $\Loc^{\check{X}}_{\check{G}}$, the pullback of the tangent complex of $\Loc^{\check{X}}_{\check{G}}$
 to $\Spec(\kk)$ is given by the cochain complex of (Betti or de Rham) $\Sigma$ with coefficients in the pullback of the tangent bundle of $\check{X}/\check{G}$. 
For example, in the Betti case, take a representation $\rho: \pi_1 \rightarrow \check{G}(k)$
and a $k$-point $x \in \check{X}(k)$ fixed via $\rho$;
the pair $(x, \rho)$ gives rise to a $\kk$-point of $\Loc^{\check{X}}_{\check{G}}$
lifting $\rho$. 
The cohomology of the tangent complex, pulled back to $\Spec(\kk)$,  computes
\begin{equation} \label{TC} H^*(\Sigma, [\check{\mathfrak{g}} \rightarrow T_{x}]_{\rho}),\end{equation}
where the subscript $\rho$ indicates that  both the Lie algebra of $\check{G}$ and the tangent space $T_x$
are considered as $\pi_1(\Sigma)$-representations by means of the representation $\rho$,
and $[\check{\mathfrak g} \rightarrow T_x]$ is referring to a $2$-term complex
with $T_x$ in degree zero. 
 In the finite context, a similar result can  be deduced in the homogeneous case using the computation of \cite{AGKRRV1} of the
 tangent complex.

 \subsection{The $L$-sheaf} \label{LsheafX}

\index{$\mathcal{L}_{\Xv}$}
\index{$L$-sheaf}

The spectral analog of the period sheaf is the push-forward of the dualizing sheaf along
$\pi: \Loc_{\Gv}^{\check{X}} \rightarrow \Loc_{\check{G}}$;  we define
the {\em unnormalized  $L$-sheaf} as the sheared $*$-pushforward of the dualizing sheaf along $\pi$:
\begin{equation} \label{Lsheafdef}  \mathcal L_\Xv := (\pi_* \omega_{\Loc^\Xv})^\shear \in QC^!(\Loc_\Gv). \end{equation}
Here the shear $\shear$ is obtained by shifting cohomological degrees thus:\footnote{Again, signs: we regard the $\GGm$ action as a {\em left} action, and the shearing
 on the dualizing sheaf arises from  regarding sections as equipped with the {\em left} $\GGm$ action.
 } the $\GGm$-action on $\check{X}$ induces  
 a $\GGm$-action on $\Loc^{\check{X}}_{\check{G}}$ covering the trivial
 action on $\Loc_{\check{G}}$. Accordingly,  the $*$-pushforward of the dualizing sheaf
 obtains a $\GGm$-action that can be used to shear it.
 That is to say, we may write $\pi_* \omega$
 as a direct sum $\bigoplus_{n} (\pi_* \omega)_n$
 of weight spaces, where $\GGm$ acts by the character $x \mapsto x^n$
 on $(\dots)_n$, and the shear is then defined as
 $\bigoplus_n (\pi_* \omega)_n \langle n \rangle$.  
 For the meaning of the angle bracket notation see \S \ref{Tatetwistrecall};
 for much more than you want about shearing, see
\S \ref{shearingsec0}.

 \begin{remark} \label{shearing0}
 
 Let us try to compare this definition with the automorphic period sheaf in \S \ref{periodX}.
 The automorphic period sheaf involves the stack $\Bun^X_G$ whose definition in \eqref{diagdog} involves a twist by $\cK^{1/2}$,
reflecting the $\Gm$-action on $X$.
 A parallel definition might be to define $\Loc_{\Gv}^{\check{X}}$ not as we have in \S \ref{locX}, 
 but instead a ``sheared'' version where its ring of functions is cohomologically sheared \index{nonconnective}
 by means of the $\Gm$-action. Such a construction will land us in general outside of ordinary derived algebraic geometry, which is built on affines given by
  {\em connective} (nonpositively graded) commutative dg algebras, and much care needs to be taken
  because  naive generalizations of usual constructions often fail. However, one can reasonably define the relevant category of sheaves 
 on such a space, using the general formalism of shearing categories. So, 
 one might think of the construction above as a substitute for actually constructing the shear of $\Loc_{\check G}^{\check X}$. 
 
 In fact, in this spirit, we can reinterpret the above definition in the following way.
 Let $p_{\Loc^\Xv}: \Loc_{\check G}^{\Xv} \rightarrow \mathrm{pt}$ be the morphism to a point; then
 the $L$-sheaf is the image of the dualizing
 sheaf of a point under the following:
    \begin{equation}   \label{quercus2}
\xymatrix{ 
 &  \IndCoh(\Loc_{\check G}^\Xv)^{\shear}\ar[rd]^-{\pi_*^\shear} \\
\IndCoh(\mathrm{pt})^{\shear} \ar[ru]^-{(p_{\Loc^\Xv}^!)^\shear} & & \IndCoh(\Loc_\Gv)^{\shear} \ar[r]^-{\unshear}& \IndCoh(\Loc_\Gv)}\end{equation}
The last map $\unshear$ refers to the identification of sheared and ordinary categories arising from the trivialization of the $\Gm$-action on $\Loc_{\Gv}$,
as in Example \ref{trivial action}.
To see this actually gives the same result as \eqref{Lsheafdef}, we use functoriality of shearing, as in Example \ref{functoriality of shearing},
and apply the discussion of \eqref{augsburg} to the diagram:  
    \begin{equation}   
\xymatrix{ 
 \IndCoh(\mathrm{pt})^{\Gm} \ar[r]^{\shear}  \ar[d]^{\pi_* p_{\Loc^\Xv}^!}  & \IndCoh(\mathrm{pt})^{\Gm\shear} \ar[d]^-{(\pi_* p_{\Loc^\Xv}^!)^\shear} \\
   \IndCoh(\Loc_{\Gv})^{\Gm} \ar[r] &
   \IndCoh(\Loc_\Gv)^{\Gm\shear} \ar[r]^-{\unshear}& \IndCoh(\Loc_\Gv)} \\
\end{equation}

 In our setting, the ``explicit'' construction of \eqref{Lsheafdef} seems to be the correct one, but it may be in more general
 situations (for example, the study of more general boundary conditions,
 in the language of \S \ref{AFT2}),  that the appropriate way to proceed is through some sheared nonconnective geometry of this kind.
This may also be related to puzzling unexplained shifts that we find when we pass outside the case of $X$ affine, cf.   \S \ref{eis period numerical comp}.
 \end{remark}

 \subsection{Normalized $L$-sheaf} \label{LsheafXnormalized}
 
 As in \S \ref{normalizedperiod}, the $L$-sheaf has a normalized version. We give a definition that depends on a choice of spin structure; see \S \ref{Lspinindep} for an invariant definition.
 
Again we suppose that 
  $\Xv$ admits a $\Gv$-eigenform with scaling character defined by \eqref{gammaXdef}, which we will denote here by $\check\eta$ to avoid confusion with the automorphic side
 \index{$\check\eta$}  \begin{equation}\label{etaspecdef} \check\eta: \Gv \rightarrow \Gm,
  \end{equation}
 This gives a map $\underline{\check\eta}: \Loc_{\Gv} \rightarrow \Loc_{\Gm}$
and we will denote also by $\varepsilon_{1/2}$ the pull-back line bundle 
$(\underline{\check\eta})^* \varepsilon_{1/2}$ on $\Loc_{\Gv}$ (see  \S \ref{spinstructure-spectral}). 
 We put
\begin{equation} \label{LXnormdef} \mathcal L_{\Xv}^{\norm} = (\pi_*  \omega)^{\shear} \otimes \varepsilon_{1/2}^{\vee} \langle -\beta_{\check{X}} \rangle.\end{equation} 
 This can be compared with   \eqref{PXnormdef} -- see Remark \ref{Heisenberg comment2} below.
 As we will see in \S \ref{Lsheafdef}, the $L$-sheaf itself can be regarded as a geometric version of the $L$-function;
 and, with respect to this,   $\mathcal L_X^{\norm}$ imitates the ``normalized $L$-function'' $L^{\norm} = L\epsilon^{-1/2}$
discussed after \eqref{epsilondef}.

\begin{remark} \label{remarksignwarning} {\bf Sign warning!} Our conventions about left and right become confusing at this point. We define 
 $\check\eta$ via \eqref{gammaXdef}, but  we are now using a {\em left} action of $\check{G}$. 
  It may be helpful to note, in tracking signs,    
that $\underline{\check\eta}$ coincides with the determinant of the $\Gv$-action on the tangent space to any fixed point,  with $\Gv$ acting by pushforward of tangent vectors.
This conclusion is very similar
to the discussion  before \eqref{gammadef}; the difference is that various left- and right- actions have been switched. 
 
   \end{remark}

\begin{remark} \label{Heisenberg comment2}
(Compare also to Remark \ref{Heisenberg comment}). The normalized $L$-sheaf
involves three twists, which are:
\begin{itemize} \item  the shearing twist $\shear$,
which depends on the $\GGm$ action on $\Xv$;
\item the $\varepsilon_{1/2}$-twist,
which reflects the failure of $X^{\vee}$ to be unimodular;
\item the twist by $\langle -\beta_{\check{X}} \rangle$.
\end{itemize}
These are roughly analogous to the three twists of Remark \ref{Heisenberg comment} but with the order of the first two twists   reversed. The negative sign of the third one
contrasts with the corresponding sign in Remark \ref{Heisenberg comment}; this is due to the fact that, before, we were shifting the constant sheaf, while here we are shifting the dualizing sheaf.

\end{remark}

 \subsubsection{Change of grading}\label{changing grading spectral}
 
 We now discuss the effect of changing the grading ($\GGm$-action) on $\Xv$ through a central modification $\lambda: \Gm \rightarrow \check G$, in parallel with \S \ref{changing grading automorphic}. 
Again, these remarks are not essential to understanding the main conjecture. Recall that, starting from a fixed $\GGm$-action, we denote by $\check X[\lambda]$ the space $\check X$ with the $\GGm$-action twisted by $\lambda$, which now denotes a central cocharacter into $\check G$.
We will prove in analogy with  Lemma~\ref{PXGmshift} that 
\begin{equation} \label{LXGmshift} \mathcal{L}_{\check{X}[\lambda]}^{\norm}    
\simeq \check{T} \mathcal{L}_{\check{X}}^{\norm} \end{equation}
with $\check{T}$ a suitable spectral translation. To express this translation we first 
digress back to geometric class field theory, cf. \S \ref{GLGMnormalization}.  (Note that, while geometric class field theory
is very helpful in interpreting the result, the proof will
not use it in any essential way; the reader can simply read the formal
proof starting from \eqref{checkTdef2}.)

Recall the degree sheaf $\ul{\deg}$ on $\Bun_\Gm$ (Remark~\ref{degree sheaf});
it is an avatar of the  square root of the cyclotomic character.
 Since $\ul{\deg}$ is locally constant, we may directly apply the equivalence of geometric class field theory 
 from \eqref{GCFT}. 

We will now identify  $\mathsf{GCFT}(\ul{\deg})$.
Let $\iota: \mbox{pt} \rightarrow \Loc_{\Gm}$
be the inclusion of the trivial local system. The skyscraper $\delta_0=\iota_*\kk$ at the trivial local system carries an inertial action of $\Gm$ by automorphisms -- it corresponds to the regular representation of $\Gm$ under the identification of the trivial bundle locus with the closed immersion $B\Gm \subset \Loc_{\Gm}$. This corresponds under class field theory to the decomposition of the constant sheaf on $\Bun_\Gm$ by components --- except with component of degree $n$ corresponding to the weight  $\lambda \mapsto \lambda^{-n}$.\footnote{
To see why the sign arises with our normalizations, note that translation by the Hecke operator $T_x$ for $x \in \Sigma$ carries sheaves on $\Bun_{\Gm}^{(n)}$ 
to sheaves on $\Bun_{\Gm}^{(n-1)}$, cf. \S \ref{GLGMnormalization}, and
on the automorphic side tensors by the sheaf that sends a local system $\rho$ to its fiber $\rho_x$, which is
in $\Gm$ degree $1$. Correspondingly, increasing the degree in $\Bun_{\Gm}$
corresponds to reducing the $\Gm$-weight on $\Loc$.}  

Therefore,  we can shear  $\iota_* \kk$  by this {\em inverted} $\Gm$-action obtaining an object $\delta_0^{\shear}$ corresponding to the $\Gm$-representation $W := \bigoplus_{n} \kk_{-n} \langle n \rangle$,
wherein $\kk_{-n}$ is in $\Gm$-weight $-n$. 
This matches the description of the degree sheaf as a shear of the constant sheaf, so that we
 obtain 
 \begin{equation} \label{deltasheafconvo} \mathsf{GCFT}(\degsheaf)\simeq \delta_0^{\shear}.\end{equation}

Thanks to the symmetric monoidal property of $\mathsf{GCFT}$, the spectral counterpart to tensoring with $\degsheaf$ is given by convolution $\delta_0^{\shear} \star -$ with the sheared skyscraper.
For example, since $\varepsilon_{1/2}$ is in $\Gm$-degree
$(g-1)$ and $\varepsilon_{1/2}^{\vee}$ in $\Gm$-degree $1-g$ we get:  
\begin{equation} \label{Delta0shear} \delta_0^{\shear} \star \varepsilon_{1/2}^{\vee} = \varepsilon_{1/2}^{\vee} \langle (g-1) \rangle,\end{equation}
The shift $\langle g-1 \rangle $ at the level of functions corresponds to the value of the square root of cyclotomic
on $\cK^{1/2}$, which equals $q^{(1-g)/2}$.

Let us now return to the nonabelian situation with a given central cocharacter $\lambda: \Gm \rightarrow \check G$.
 Note that $\lambda$ induces an action of the group object $\Loc_{\Gm}$ on $\Loc_{\check G}$,
with respect to which we shall consider the convolution action $\check{T}$ of the inverse of $\delta_0^{\shear}$, that is to say
\begin{equation} \label{checkTdef2} \check{T} = \check T_\lambda: \mathcal{F} \mapsto \mbox{pushforward of $(\delta_0^{\shear})^{-1} \boxtimes \mathcal{F}$
via $\Loc_{\Gm} \otimes \Loc_{\check G} \rightarrow \Loc_{\check G}$},\end{equation}

Equivalently, we observe that $B\Gm \subset \Loc_{\Gm}$
also acts on $\Loc_{\check G}$ through the embedding $\lambda$,  and we may write
\begin{equation} \label{checkTdef3} \check{T} := \mbox{convolution action
of $(\delta_0^{\shear})^{-1}$, as a sheaf on $B\Gm$, on $\Loc_{\check G}$}\end{equation}
Then $\check{T}$ has the effect of regrading sheaves on $\Loc_{\check G}$
according to the action of $\lambda: \Gm \rightarrow \check G$. 
In particular, this describes the action of the substitution $X \mapsto X[\lambda]$ on  $(\pi_* \omega)^{\shear}$:
 $$ \mathcal{L}_{X[\lambda]} = \check{T} \mathcal{L}_X.$$
 
Now by \eqref{Delta0shear}, taking account that we now have inverted $\delta_0^{\shear}$, 
and that $\varepsilon_{1/2}^{\vee}$ is pulled back via $\underline{\check\eta}$ (see before \eqref{LXnormdef}), 
we get for an arbitrary sheaf $\mathcal{F}$ on $\Loc_{\check{G}}$ the equality
$$\check{T} (\mathcal F \otimes \varepsilon_{1/2}^{\vee}) = \check{T} (\mathcal F) \otimes \varepsilon_{1/2}^{\vee} \langle (1-g) \langle \check\eta, \lambda \rangle \rangle.$$
For the normalized $L$-sheaf we find
$$  \mathcal{L}_{X[\lambda]}^{\norm} =  \mathcal{L}_{X[\lambda]} \otimes \varepsilon_{1/2}^{\vee} \langle - \beta_{\check{X}[\lambda]} \rangle =   (\check{T} \mathcal{L}_X) \otimes \varepsilon_{1/2}^{\vee} \langle - \beta_{\check{X}[\lambda]} \rangle$$
  $$ =\check{T} (\mathcal{L}_X \otimes \varepsilon_{1/2}^{\vee}) \langle -\beta_{\check{X}[\lambda]} + (g-1) \langle \check\eta, \lambda \rangle \rangle
= \check{T} \mathcal{L}_X^{\norm}$$
for just as in \eqref{betachange} we have $\beta_{\check{X}[\lambda]} = \beta_{\check{X}} + (g-1) \langle \check\eta, \lambda \rangle$. This 
confirms \eqref{LXGmshift}.

  \subsection{$L$-sheaves for twisted polarizations.} \label{case4} 
We now discuss the construction of $L$-sheaves for twisted cotangents and Whittaker inductions, a spectral counterpart to the twisted period sheaves in \S \ref{Affcase}. Let us recall that the construction of these twisted period sheaves
amounted to ``twist by a rank one Artin--Schreier local system on $\Bun^X_G$,''
and moreover that Artin--Schreier sheaf was pulled back from 
the space of torsors for the canonical bundle $\mathcal{K}$. 
 We will do something similar, now, on the spectral side. 
 
 \begin{remark}
 We note that Whittaker inductions on the spectral side correspond to what one might call  ``Arthur functoriality'' or ``Arthur lifting'' in the classical theory of automorphic forms, as discussed in \S \ref{L2numerical} (specifically \S \ref{nontempered} and \S \ref{Arthur functoriality}). 
\end{remark}

 Let $\Gv\times \Gm\actson (\Xv,\Psi)$ be as in \S \ref{twisted cotangent bundle section}. 
The affine bundle $\Psi\to \Xv$ defines a $\Gv\times \Gm$-equivariant
map from $\Xv$ to the classifying space $B\Ga$ (with the squaring action of $\GGm$ and trivial action of $\Gv$), whence a map 
$$ \Loc_{\Gv}^{\Xv}    \rightarrow \Loc_{\Ga}.$$
Our goal is to define the $L$-sheaf
by twisting the previous construction by the pullback to $\Loc_{\Gv}^\Xv$ of a ``spectral exponential'' (or spectral Artin--Schreier) sheaf on $\Loc_{\Ga}$.
The discussion will be parallel to the corresponding discussion of period sheaves, 
but much more confusing, because
\begin{itemize}
\item[(i)] of the presence of shearing in the definition of the $L$-sheaf, and
\item[(ii)] it is not clear what the spectral exponential sheaf on $\Loc_{\Ga}$ should be.
\end{itemize}
The key point is that difficulty (i) and (ii) cancel each other out:
our construction of the spectral exponential sheaf (spelled out in \S \ref{spectral exponential}) exists only after shearing.

\begin{remark} 
In the case $\Xv=\check U\backslash \Gv$ with
the standard $\Psi$, the only work that we are aware of is due to V. Lafforgue \cite{LafforgueP1} --
 it  can be checked that our definition matches with that suggested therein -- 
 and work in preparation of Hilburn-Yoo \cite{HilburnYoo}.
In fact, Lafforgue's computation makes clear the following striking point:
the spectral Whittaker sheaf is an object
of  ind-coherent category $QC^{!}$ {\em whose projection to $QC$ is not bounded below.}
\end{remark}

\subsubsection{Construction of the spectral exponential sheaf} 
We will describe the construction in the Betti case, leaving the modifications to the reader. A choice of orientation on $\Sigma$ determines a map
\begin{equation} \label{loca1morph} \Loc_{\Ga} \rightarrow \AA^1[-1]\end{equation}
which will play a role similar to \eqref{GaGmT}.
 Here $\AA^1[-1]$ is the derived scheme 
 whose ring of functions is $\kk[x_{-1}]$, see \S \ref{shiftingofschemes}.

Here is how the map \eqref{loca1morph} arises. Given a   complex $V^{-1} \rightarrow V^0 \rightarrow \dots$ of vector spaces over $\kk$, 
we can functorially associate a derived stack over $\kk$, 
which can be understood as the quotient of 
the derived scheme whose ring of functions is the dual symmetric algebra on $[V^0 \rightarrow V^1 \rightarrow \dots]$ 
 by the action of $V^{-1}$ considered as a vector group. The isomorphism
 class of the resulting derived stack depends only on the quasi-isomorphism class of the original complex.   
 (For the functor of points in more generality see  \cite[Section 3.3]{ToenEMS}.)
Now the space of local systems $\Loc_{\Ga}$ for $\Ga$
is the vectorial derived stack associated, in this fashion,  to the cochain complex $C^*(\Sigma)[1]$;
mapping this complex to  its truncation $H^2(\Sigma)[1]$ in degrees $\geq 1$ gives
a morphism from $\Loc_{\Ga}$ to $H^2(\Sigma, \kk)[-1] \simeq \AA^1[-1].$

For example, in the Betti case,  the morphism \eqref{loca1morph}
follows
 immediately from the presentation \eqref{locGprez}, taking
 into account that the commutator is trivial and that the conjugation action is trivial;
 namely, the generator $x_{-1}$ is sent to the degree $-1$
 element of the ring of functions associated to the relation from \eqref{locGprez}.

On $\AA^1[-1]$ there is an object that plays the role of the exponential sheaf,
but exists only after some fiddling: 
 $$\exp\in QC^!(\AA^1[-1])^\shear$$ 
 It is defined in \S \ref{spectral exponential};
 the shearing is for the squaring action of $\Gm$ on $\AA^1$. 
Roughly, it is ``Koszul dual to a skyscraper at $1$'' and the various adornments $!, \shear$
are formal adjustments to the category that allow this to make sense:
 
 \begin{itemize}
 \item By Koszul duality, the category $QC^{!}$ for $\AA^1[-1]$ is identified
 with the ``category of sheaves on $\AA^1[2]$,''
 which, formally speaking, is defined as a sheared version of the category of sheaves
 on $\AA^1$. 
 
 \item The shearing has the effect of replacing $\AA^1[2]$ by $\AA^1$.
 
 \item Finally, $\exp$ corresponds to the skyscraper at $1 \in \AA^1$.
 \end{itemize}

\subsubsection{Construction of the $L$-sheaf}
Before reading the following, the reader may want to glance at
the reformulation of the untwisted definition  given in \eqref{quercus2}. 
It is this reformulation that the twisted definition will parallel. 

Let $\overline{\Psi}_C$ denote the composite $\Gm$-equivariant morphism $ \Loc^\Xv \rightarrow  \Loc_{\mathbb{A}^1} \rightarrow \AA^1[-1]$
(squaring action on $\AA^1$). 
The resulting diagram
\begin{equation} \label{first} \xymatrix{
& \ar[ld]_-{\overline{\Psi_C}}  \Loc^\Xv  \ar[rd]^-{q}
\\
\AA^1[-1] && \Loc_\Gv}\end{equation} is also $\Gm$-equivariant, with respect to the trivial action on $\Loc_\Gv$. It induces $\Gm$-equivariant functors on categories of ind-coherent sheaves
 and then (by the functoriality of the shearing process) also on sheared categories (see \S \ref{shearingsec} for background)  \begin{equation} \label
{shear crux} 
\xymatrix{ 
 &  \IndCoh(\Loc_\Xv)^{\shear}\ar[rd]^-{q_*^\shear} \\
\IndCoh(\AA^1[-1])^{\shear} \ar[ru]^-{(\overline{\Psi_C}^!)^\shear} & & \IndCoh(\Loc_\Gv)^{\shear} \ar[r]^-{\unshear}& \IndCoh(\Loc_\Gv)}\end{equation}
The last map $\unshear$ refers to the identification of sheared and ordinary categories arising from the trivialization of the $\Gm$-action on $\Loc_{\Gv}$,
as in Example \ref{trivial action}.

We are finally ready for the definition.

\begin{definition}\label{spectral Whittaker def}
Let $\Gv\actson (\Xv,\Psi)$. We define the associated $L$-sheaf by pushing forward the Whittaker-exponential sheaf  from $\Loc_\Gv^\Xv$
along the sequence  \eqref{shear crux}. In symbols, 
\begin{equation} \label{Whitsymbols} \mathcal L_{\Xv,\Psi}= \unshear \circ q_*^\shear\circ (\overline{\Psi}_C^!)^\shear (\mathrm{exp}) .\end{equation}

The normalized $L$-sheaf is defined by twisting 
the definition of $\exp$ by $\varepsilon_{1/2}^{\vee}$ and 
shifting the end result by $-\beta_{\check{X}}$, 
cf. \eqref{LXnormdef}.
\end{definition}
 
\index{spectral Whittaker sheaf}
 In the case $\Xv=\Gv/\check U$ with $\Psi$ a nondegenerate character of $U$ this object can be reasonably termed the {\em spectral Whittaker sheaf},
 i.e., a spectral analogue of the automorphic Whittaker sheaf that, in turn, geometrizes the Whittaker period in the theory of automorphic forms. Unfortunately,  it is quite difficult to compute with.  Our primary evidence that is the right definition comes -- besides the parallel with the automorphic definition --  from our computations in the $\mathbb{P}^1$ case (see \S \ref{P1})
 since we do not do any direct numerical computations.  See \S \ref{Arthur functoriality} and \S \ref{WhitArth} for further discussion of the spectral Whittaker construction in the context of Arthur parameters.

 \subsection{Dependence on spin structures} \label{Lspinindep}
 
 In this short subsection,
 parallel to \S \ref{spindep}, we collect a few observations on
 the dependence of constructions on spin structures.  See also \S \ref{extended-group appendix}. 
 Again, we emphasize we will make little use of this formulation and include it for completeness. 
  \index{spectral spin structure}
 
 Let us first say what we are {\em not} doing. 
One can argue
 that the spectral analogue of a spin structure, i.e., the spectral analogue of a
 square root of the canonical bundle, is given 
 -- at least in the finite context -- by a  square root of the cyclotomic character. 
 To make our discussion parallel to that of periods, then,
 we might like to more systematically 
 choose a square root of the cyclotomic character,
 which roughly speaking will index a choice of $\sqrt{q}$
 that ``varies over $\Sigma$,''
 and with such a choice  we could treat
spectral spin structures and automorphic spin structures on precisely the same footing.
However, we have chosen {\em not} to do so -- in effect, the choice of $\sqrt{q} \in \kk$
leads to a specific
  choice of a square root of the cyclotomic character, which in the 
  notation of \S 
\ref{CFTGLconventions}  corresponds to the square root $|x|^{1/2}$
fixed by the given choice of $\sqrt{q}$.   To our knowledge other choices are  never considered in the automorphic literature.

 The remarks below are, therefore, not about this question, but about spin structures in 
precisely the same sense as \S \ref{spindep}, that is to say, square roots of the coherent
dualizing bundle.   We will use
the same notation as in \S \ref{spindep}.

 First, let us discuss the dependence of the half-epsilon line bundle $\varepsilon_{1/2}$ on spin structures. Changing $\cK^{1/2}\mapsto \cK^{1/2}\otimes \mathcal{L}$ for $\mathcal{L}\in \Bun_{\Z/2}$ replaces $\varepsilon_{1/2}=[\cK^{1/2}]$ by
$$[\cK^{1/2}\otimes \mathcal{L}]\simeq \varepsilon_{1/2}\otimes [\mathcal{L}]$$ 
where $[\mathcal{L}]$ is a 2-torsion line bundle on $\Loc_\Gm$. In other words, if we consider $\Bun_{\Z/2}$ as acting on $QC(\Loc_{\Gm})$ via the homomorphism $[-]$
of \eqref{Lprimedef},  then we have a canonically defined object in a twisted version of sheaves, $$\varepsilon_{1/2}\in QC(\Loc_\Gm)^{\mathrm{Spin}}:=\Hom_{\Bun_{\Z/2}}(\mathrm{Spin}_\Sigma, QC(\Loc_\Gm))$$ (note that this fomulation requires defining $[-]$ in a structured fashion rather than just the level of isomorphism classes --
see Remark \ref{Abelian double covers} below). 

 Next, the map $\underline{\check\eta}: \Loc_{\Gv} \rightarrow \Loc_{\Gm}$ allows us to pull back the $\mathrm{Spin}_\Sigma$-twisted object $\varepsilon_{1/2}\in QC(\Loc_\Gm)^{\mathrm{Spin}}$
just defined, to give a twisted sheaf on $\Loc_\Gv$, which we can tensor with the unnormalized $L$-sheaf to account for dependence on spin structures: 
$$ \mathcal{L}_{\Xv}^{\norm}  =\mathcal{L}_\Xv \langle -\beta_{\check{X}} \rangle  \otimes \varepsilon_{1/2}^{\vee} \in QC^!(\Loc_\Gv)^{\mathrm{Spin}}:=\Hom_{\Bun_{\Z/2}}(\mathrm{Spin}_\Sigma^\vee, QC^!(\Loc_\Gv))$$
where $\Bun_{\Z/2}$ acts on $QC^!(\Loc_\Gv)$ by tensor product with 2-torsion line bundles pulled back from $\Loc_\Gm$.  The duality exponent on $\Spin$ reflects
the fact that $\varepsilon_{1/2}$ is dualized in the definition of the normalized $L$-sheaf. 
Compare with \eqref{pxtwisted}.

\begin{remark} \label{Abelian double covers} For later use, we note that the action of $\Bun_{\Z/2}$ on $QC^!(\Loc_\Gv)$
that has just appeared can be described directly in terms of double covers of $\Gv$
 (see also Remark~\ref{fluxes} in the Appendix).\index{flux}
 For this,  we describe in the Betti or de Rham case an explicit version of the map  $[-]$ for $\Z/2$-torsors.
The map $\Gm\to B(\Z/2)$ classifying the double cover of $\Gm$ determines a cover of $\Loc_\Gm$ with Galois group $H^1(\Sigma, \Z/2)$,
whose fiber at a local system consists of its lifts along $\Gm\stackrel{x\mapsto x^2}{\longrightarrow} \Gm$.  Cup product with the class of a $\Z/2$-bundle (i.e., the Weil pairing self-duality of $H^1(\Sigma,\Z/2)$)
  then gives us a map
 from 
 \begin{equation} \label{abelian double cover eqn} \Bun_{\Z/2} \longrightarrow \mbox{$\Z/2$-local systems on $\Loc_\Gm$ }\end{equation}
Pullback via $\underline{\check\eta}$ gives  the
desired morphism $\Bun_{\Z/2}\to \Pic(\Loc_\Gv)$. \end{remark}

\subsection{$L$-sheaves and $L$-functions} \label{locXex}
  
  In order  to better understand the meaning of the definitions, 
  we are going to compute some fibers of the $L$-sheaf, and see
  that the resulting vector spaces categorify $L$-funtions.  The main results
  are summarized in Table 
 \ref{LsheafLfunctioncomputation}. The results here are not used in any formal way
 in the study of the geometric conjecture, but they are part of the motivation for the numerical conjecture
 enunciated later. Finally, in \S \ref{derived volume forms} we sketch another point of view on $L$-functions coming from the theory of categorical traces in derived algebraic geometry, as expressing the derived volume of Galois fixed points on $\Xv$.
 
   Let $\rho$ be a $\check{G}$-local system with coefficients in $\kk$. In order to discuss simultaneously the geometric and the arithmetic contexts, in the finite case we will be thinking of $\rho$ as a representation
    $$\rho: \pi_1(\Sigma) \longrightarrow \check{G}(\kk),$$
    where $\pi_1(\Sigma)$ is the \'etale fundamental group (with respect to a fixed geometric base point) of the curve $\Sigma$ \emph{over $\FF_q$}, and 
  we will write 
  $\rho^{\geom}$
when we need to emphasize that we restrict
  $\rho$  to geometric $\pi_1$. From $\rho$ we get 
  a $\kk$-point of $\Loc_{\check{G}}$, 
  $$ \iota_{\rho}: \mbox{pt} \rightarrow \Loc_{\check{G}},$$
which depends only on $\rho^{\geom}$ in the finite context. 
    
  In what follows we will describe two types of results. Firstly, we will compute 
  various stalks of the $L$-sheaf; these results are valid in all contexts.
  We will then compute Frobenius traces on these stalks; these
  results should be understood as applying only to the finite context.
  When we talk about Frobenius trace, here, we are implicitly
  using an extra structure, namely, 
  the fact that the $L$-sheaf has a natural equivariant structure for the action 
  of Frobenius on $\Loc$ described in
   \S \ref{Frobaction}.

  The main takeaways are, firstly, that these traces are $L$-functions, and, secondly, that
  the $\GGm$-action shifts the point of evaluation of those $L$-functions. 
  The reader accustomed to arithmetic settings will note with some bemusement the signs $(-1)^d$
  that occur below due to cohomological shifts by $d$;
  somewhat surprisingly,  these reflect numerical phenomena that are quite complicated, see \S \ref{starperiods}!

  \subsubsection{Notation and the conditions (a), (b), (c)} \label{SSSrestrictions}
  
  $\rho^\geom$ defines  a $\kk$-point of $\Loc_{\check{G}}$, 
 and 
 let $$\delta_{\rho}=\iota_{\rho*}(\kk)$$ be the associated skyscraper sheaf.  
 $\rho^\geom$ may have automorphisms;  we write $Z(\rho)$ for this algebraic group,  
 i.e., the centralizer of $\rho^\geom$ in $\check{G}$, and put
$$ d = \mbox{dimension of  $Z(\rho)$.} $$

 Consider the following three restrictions (a), (b), (c) 
 that can be placed on the situation.
We will always (in 
the current section \S \ref{locXex})
 impose (a) and (b) below, and sometimes also (c) which will
force $d=0$:
\begin{itemize}
 \item[(a)]  
 We assume that $\check{M} =T^*\check{X}$, i.e., the affine bundle $\Psi$ is trivial. 
 However, our analysis
 will apply to arbitrary smooth $G \times \GGm$-spaces $\Xv$ (e.g.\ not necessarily affine, or
 satisfying any other conditions). \footnote{  In \S \ref{P1} we will carry out an analysis
 of the general type performed here in cases that allow twisted polarization, but restricted to the case of the base curve $\mathbb{P}^1$.}
 
  \item[(b)]  The classical fixed point locus of
  $\rho^{\geom}$ on $\Xv$  is a singleton $x_0$ -- that is to say, 
  $$ \Xv^{\rho_{\geom}} = \{x_0\},$$
  where we consider $x_0$ as a reduced scheme.
  This implies also -- in the finite context --  that $x_0$ will be fixed by $\rho$. 
  We write $T$ for the tangent space at $x_0$; then our assumption entails that
 \begin{equation} \label{tanvanish} H^0(\Sigma, T) = H^2(\Sigma, T) = 0,
 \end{equation}
where $T$ is regarded a local system of vector spaces on $\Sigma$
via $\rho^{\geom}$, and (see \S \ref{sigmacoh}) cohomology in the finite
case is always geometric {\'e}tale cohomology. 
See also
 Remark \ref{nondiscrete}.
 
   \item[(c)]   $\rho$ is a smooth point of the moduli space $\Loc_{\check{G}}$,
  equivalently, $$H^0(\Sigma, \ad \  \rho^{\geom}) = H^2(\Sigma, \ad \ \rho^{\geom}) = 0,$$
  with $\mathrm{ad}$ the adjoint action of $\check{G}$ on its Lie algebra, i.e., $\ad \ \rho^{\geom}$
  is a local system of vector spaces of dimension $\dim(G)$. 
    \end{itemize}

 Now, the $\GGm$-action on $\check{X}$ fixes $x_0 \in \check{X}(k)$ by the assumed uniqueness in (b). 
Therefore, $\Gamma_F$ fixes $x_0$ not only
through its action via $\rho$, but through its action via
the corresponding extended Langlands parameter:
$$  \pi_1(\Sigma) \stackrel{(\rho, \varpi)}{\longrightarrow} \check{G} \times \Gm \rightarrow \mathrm{GL}(T), \ \ T=T_{x_0} X.$$
 We may therefore define the $L$-function $L(s, T^{\shear})$ according to the prescription of \eqref{Ldef};
 thus explicitly if 
     $T = \bigoplus_{k} T_k$ we have  $L(0,T^{\shear}) := \prod_{k} L(k/2, T_k)$ and if $T$ has trivial grading we get just $L(0,T)$.
     We have similarly defined the normalized $L$-function, see \eqref{Lnormdef}, and
     will e.g.\ talk of $L^{\norm}(0, T^{\shear})$ defined similarly. 

To keep the typography simple, we define
$$H^* T := \mbox{ cohomology of $\Sigma$ with coefficients in the local
system defined by $T$,}$$
(again, {\em geometric} {\'e}tale cohomology). Moreover, we define
$H^* T^{\shear}$ to be the sheared version of this cohomology, which coincides with the cohomology
of $T^{\shear} = \bigoplus_k T_k\langle k \rangle$. 
 Thus, for example, if $H^* T$ is concentrated in degree $1$, so $H^* T = H^1 T[-1]$, 
we have $(H^* T)^{\shear} = H^*(T^{\shear}) =\bigoplus H^{1} T_k[k-1] (k/2)$.

  \begin{remark} \label{caveat} Recall that,  
    in the finite context, we have not defined $\Loc^{\check{X}}_{\check{G}}$
  in all cases; for the purpose of interpreting the results, one may either
  restrict to the case in which $\check{X}$ is defined, or assume
  that there exists a definition  
  for which the analysis below goes through (the only key points is that it should
  satisfy base change and has the correct tangent complex).  We do not  see any essential difficulty in establishing these points,
  but we have not done so.
   \end{remark}

    \begin{remark} \label{nondiscrete}
In the finite context we are supposing that the fixed locus
 of $\rho$ restricted to {\em geometric} $\pi_1$
 is a singleton. However, one may hope that a
 similar result remains valid under the restriction only to {\em arithmetic} $\pi_1$,
 because of ``localization.''    
  \end{remark}

   \begin{remark} \label{completion caveat}
  The results that follow, in the case when $\rho^\geom$ has positive dimensional centralizer $Z(\rho)$,
  will refer to certain ``completions'' where we replace
  a sum over $Z(\rho)$-isotypical components by a product.
  More precisely,  for $\Zrep$ a representation of $Z(\rho)$, decomposed into isotypical components \index{$\Zrep^{\wedge}$}
as $\Zrep =\bigoplus T_{\alpha}$, we write $\Zrep^{\wedge}$ for the ``completion'' 
\begin{equation} \label{completionX} \Zrep^{\wedge} =  \Hom(k[Z(\rho)], \Zrep) = \prod \Zrep_{\alpha}.\end{equation}
  \index{completion} \index{$\Zrep^{\wedge}$}
 Let us describe where this infinitude arises from in a simple example, just
 so as to emphasize that it does not reflect any actual pathology of the situation.
Take $\check{G} = \Gm$ and $\check X=\check{G}$. 
In this case, the $L$-sheaf is the push-forward $\iota_{\rho*} \kk$ from
$\rho=\mbox{trivial}$ and corresponds, automorphically, to the constant
sheaf.   
Then 
$$  \Hom(\delta_{\rho}, \mathcal{L}_{\check{X}}),$$
for $\rho$ trivial,   corresponds to the  automorphic computation of
$\Hom(\kk, \kk)$ on $\Bun_{\Gm}$, i.e., we are computing
$H^*(\Bun_{\Gm})$, which is a product over the infinitely many connected components.
For this reason, the above $\Hom$-space is an infinite direct product.
This infinitude reflects the non-compactness of $\Bun_{\Gm}$,
and can be avoided by restricting the computation above to 
   a finite collection of connected components of $\Bun_{\Gm}$.
  \end{remark}

 {\small 
  \begin{table}[!htbp]
  \caption{The $L$-sheaf and which $L$-function it categorifies. In the case $d>0$
  the results must be ``completed''; see Remark \ref{completion caveat}. 
  Also $\det H^1(\ad \rho)$ is placed in degree $\dim(\ad \rho)$.}
  
  \label{LsheafLfunctioncomputation} \centering 
  \begin{tabular}{|c|c|c|c|}
 \hline
assumption			&  		computation							&  			 result							&			trace of Frobenius on {\em dual}						 	\\ 
\hline
(ab)			&  $\Hom( \delta_{\rho}, \mathcal{L}_{\Xv})	$					&  $\left(\Sym^* H^*(T^{\shear})\right) [-d] $										& $(-1)^d L(1,T^{\shear\vee})$			  \\
(ab)			&	 $ \Hom( \delta_{\rho}, \mathcal{L}_{\Xv}^{\norm})$				& $ \left(\textrm{same}\right) \otimes \varepsilon^{\vee}_{1/2}(T)\langle -\beta_X  \rangle $										&  $(-1)^d q^{-b_{\check G}/2} L^{\norm}(1, T^{\shear\vee})$.			  \\
 \hline
 (abc)			&  $\Hom( \mathcal{L}_{\Xv},\delta_{\rho})	$			&   $(\Sym^* H^*T^{\shear})^{\vee}  \otimes   \det  H^1(\ad \rho) $										&  
  $ q^{-b_{G}/2}  L(0,T^{\shear})$				  \\
  (abc)			&   $  \Hom(\mathcal{L}_{\Xv}^{\norm}, \delta_{\rho})$		& $  \left(\textrm{same} \right) \varepsilon_{1/2}(T) \langle \beta_X \rangle $										& 
 $ q^{-b_{G}/2}  L^{\norm}(0,T^{\shear}).$		  \\
 \hline
 \end{tabular}
 \end{table}
}

 \subsubsection{Results} \label{LsheafLfunctionresults}
  
  We summarize the results that follow in Table \ref{LsheafLfunctioncomputation}. 
  It lists (under the assumptions listed in the far left) what 
  $\Hom$s from $L$-sheaves to and from skyscrapers are. 
  The final column is relevant only in the
  finite context, and computes the trace of Frobenius
  on the {\em dual} of this $\Hom$-space, which
  will help guide us in our later numerical discussions. 
 We now spell out what the table says a little more carefully.
 
 First two lines of the table:
Assuming  conditions (a), (b) from \S \ref{SSSrestrictions},  and that the centralizer of $\rho$ is
  $d$-dimensional we will compute below that 
    \begin{multline} \label{LX1} \Hom( \delta_{\rho}, \mathcal{L}_{\Xv}) =   
\mbox{completion of } \left(\Sym \ H^*(T^{\shear}) \right) [-d]  \mbox{ and } \\  
 [\mbox{uncompleted RHS above}^{\vee}] =  (-1)^d L(1,T^{\shear\vee}). 
 \end{multline}
 The first isomorphism is an isomorphism of Frobenius modules in the finite case.
 The completion is discussed above in Remark \ref{completion caveat} but is necessary only in the case $d>0$.   The meaning of the notation $[\dots]$ is to take the trace of geometric Frobenius, see \S \ref{anglebracketnotation}. 
 The reason that we are computing with the dual space on the second line of \eqref{LX1} comes from Lemma \ref{Homlemma}, which expresses the inner product of two functions in terms of the dual of a Hom-space of sheaves;
 we anticipate, therefore, that it is not $\Hom(\delta_{\rho}, \mathcal{L}_{\check{X}})$
 but its dual  whose trace corresponds to a meaningful numerical computation on the automorphic side. 
 
This and the definition \eqref{LXnormdef} of the normalized $L$-sheaf
will imply  corresponding results in the normalized case:
\begin{multline}    \Hom( \delta_{\rho}, \mathcal{L}_{\Xv}^{\norm}) =   
\mbox{completion of } \left(\Sym^* H^*(T)\right) [-d] \otimes \varepsilon^{\vee}_{1/2}(T)\langle -\beta_X  \rangle \mbox{ and }
 \\ \label{LX1b}  [\mbox{uncompleted RHS above}^{\vee}] = 
 (-1)^d q^{-b_{\check G}/2} L^{\norm}(1, T^{\shear\vee}).
 \end{multline}
 Note that the symmetric algebra here simply amounts to the exterior algebra on $H^1$
 since assumption (a) entails that $H^0$ and $H^2$ is zero.
 
 Second two lines of the table:
{\em Additionally supposing} (c) from \S \ref{SSSrestrictions},  we will also compute:
    \begin{multline} \label{LX0} \Hom(\mathcal{L}_{\check X}, \delta_{\rho}) =
(\Sym \  H^*T^{\shear})^{\vee}  \otimes   \det  H^1(\ad \rho) \implies  
[ \Hom(\mathcal{L}_{\check X}, \delta_{\rho}) ^{\vee}] =  q^{-b_G} L(0,T^{\shear}), \end{multline}
    \begin{multline} \label{LX0b} \Hom(\mathcal{L}_{\check X}^{\norm}, \delta_{\rho}) =
\left(\textrm{same} \right) \varepsilon_{1/2}(T) \langle \beta_X \rangle \implies  \\
 [ \Hom(\mathcal{L}_{\check X}^{\norm}, \delta_{\rho})^{\vee}] \sim
 q^{-b_{G}/2}  L^{\norm}(0,T^{\shear}).\end{multline}

\begin{remark} \label{interesting}
It is an important and interesting question to relax
restriction (c) of \S \ref{SSSrestrictions}, i.e., to prove versions of \eqref{LX0} and \eqref{LX0b}
when $\rho$ is not a smooth point. The reason is that
it is understood that ``$!$ period paired with an automorphic form gives $L$-function''
is valid not merely for cusp forms: it remains valid in the Eisenstein case, away
from the polar locus of the $L$-function. What happens at the polar locus of the $L$-function
is a subtle question, even at the purely numerical level -- relaxing (c) would help understand this. 
\end{remark}

  \subsubsection{Proofs of the statements about fibers} \label{numnum}
  
 We now give the proof of the above statements about fibers of the $L$-sheaf that were stated in \S \ref{LsheafLfunctionresults}; the numerical
 statements will be proved separately below in \S \ref{numnum}. 
 
 We will compute this in the Betti model, the analysis of the other cases being similar (see Remark \ref{caveat}). 
  We shall consider the pullback diagram
 $$ 
\xymatrix{
\Loc^{\Xv}_{\rho} \ar[r]^{\iota_{\Xv}} \ar[d]^{\pi_{\rho}} & \Loc_{\Gv}^{\Xv} \ar[d]^{\pi} \\
\pt \ar[r]^{i_{\rho}} & \Loc_{\Gv}
}
$$

Now we observe that   the fiber of $\Loc_{\Gv}^{\check{X}}$ above $\rho$ is the derived scheme $H^1(\Sigma,  T)[-1]$.  
Indeed, by our assumption (b), this fiber is an affine derived scheme with just one classical point, and its tangent
complex is $H^1(T)$ concentrated in degree $1$.

 Now,  in $QC^{!}$ one has base change for $*$-pushforward and $!$-pullback, without restrictions.     
 On the other hand, the $*$-pullback of ind-coherent sheaves, satisfying base-change with $*$-pushforward, is only defined for morphisms which are eventually co-connective \cite[4.3.2]{GR};  this is certainly the case for the inclusion $\iota_{\rho}$ of a smooth point of $\Loc$. 

 Suppose first that $\rho$ has finite centralizer, so that $d=0$, $i_\rho$ is proper and $(i_{\rho*},i_\rho^!)$ form an adjoint pair; we then compute by base change  that 
\begin{equation} \label{irho!} \Hom(\delta_{\rho}, \mathcal L_{\Xv}) = \Hom(i_{\rho*} \kk, \mathcal L_{\Xv}) =  \iota_{\rho}^{!} \mathcal L_{\Xv}
= \pi_{\rho*} \omega_{\Loc^\Xv_{\rho}}^{\shear} \end{equation}
so we get sheared sections of the dualizing sheaf of $\Loc_X^{\rho}$. 
 As we saw $\Loc_X^{\rho}$ is simply 
  the vector space $\mathsf{V} = H^1(T)$
shifted via $[-1]$ to a derived scheme; as a derived
scheme its ring of functions is $\mathrm{Sym} \  (\mathsf{V}^\vee[1])$ and 
forms are the dual of functions:
\begin{equation} \label{omega sec}\mbox{sections of $\omega$} =  \mathrm{Sym} \  (\mathsf{V}[-1]) = \Sym\ (H^1 T[-1]) =  \Sym \  (H^* T).\end{equation} 
{\em  The symmetric algebra is taken in the graded sense}, i.e.
   at the level of underlying vector spaces this is the exterior algebra on $\mathsf{V}$. 
   In combination with \eqref{irho!} 
 this establishes \eqref{LX1} in the case when $\rho$ has finite centralizer, i.e., $d=0$; we will discuss the numerical version below.

Now we examine the case when $\rho$ has $d$-dimensional centralizer $Z$. 
Here we will use the completion from \eqref{completionX}. 
Now we can factor $\iota_{\rho}$ as 
 $\mathrm{pt} \stackrel{a}\rightarrow \mathrm{pt}/Z \stackrel{b}\rightarrow \Loc_{\Gv}$
 with $b$ proper.
 Observe \begin{eqnarray} \label{BZrev} \Hom(b_* a_* \kk, \mathcal L_X) &=& \Hom(a_* \kk, b^{!} \mathcal L_\Xv)\\ 
 & \stackrel{\cdot}{=} &
 ( a^* b^{!} \mathcal L_{\Xv}  )^{\wedge}\\
 &=&  (a^{!} b^{!} \mathcal L_{\Xv} [-\dim Z])^{\wedge} \\
 &=& (i_{\rho}^{!} \mathcal L_{\Xv}[-\dim Z])^{\wedge},\end{eqnarray} 
where we identified sheaves on $\mathrm{pt}/Z$ with $Z$-representations. 
  The dotted equality
arises as follows:
the sheaf $a_* \kk$ is identified with the regular representation $\kk[Z]$
of $Z$, and then we use  $\Hom(\kk[Z], W) \simeq W^{\wedge}$ for any $Z$-representation $W$,
with $W^{\wedge}$ as in \eqref{completionX}. 
We also used $a^{!} = a^*[\dim Z]$ for the smooth morphism $a$. 
Equation  \eqref{irho!}  therefore holds with a further shift of $-\dim Z$, giving again \eqref{LX1} in this case. 
The version with the normalized $L$-sheaf follows at once from the definitions.

Next, let us impose assumption (c), so that  $\Loc$ is smooth at the point $\rho$.
 In this case the map $i_{\rho}$ is ``eventually co-connective'' and 
 also LCI; so 
 $i_X^{!}$ and $i_X^*$ differ by a shift. 
By adjunction we have
$$\Hom(\mathcal{L}_X, \delta_{\rho}) = \Hom_{\mathrm{pt}}(\iota_{\rho}^* \mathcal{L}_X, \delta)$$
Then 
 $ i_{\rho}^* \pi_* \omega = \pi_{\rho*} i_{\rho}^* \omega$
 (here and below $\omega$ is dualizing on $\Loc^X$)
; canonically  
 $ i_{\rho}^{!} = i_{\rho}^* \otimes   (\det T_{\rho})$
 where $T_{\rho}$ is the determinant of the tangent space at $\rho$,
 {\em which we understand to be placed in dimension $\dim T_{\rho}$}; 
 and -- since $i_{\rho}^!$ carries dualizing to dualizing --  $i_{\rho}^* \omega = \omega_{X,\rho} (\det T_{\rho})^{-1}$ where  $\omega_{X,\rho}$ is dualizing of $\Loc^X_{\rho}$. 
  We deduce that $$\Hom(\mathcal{L}_X, \delta_{\rho}) =  \left[
  \mbox{sections of $\omega_{X, \rho}^{\shear}$}
  \right]^{\vee} \otimes \det  T_{\rho}.$$
  Note that $T_{\rho}= H^1(\mathrm{ad} \rho)$. This in combination with \eqref{omega sec} proves \eqref{LX0}
  and again \eqref{LX0b} follows simply by taking into accounts the shift in the definition of the normalized $L$-sheaf.   
  
  \subsubsection{Proofs of statements about Frobenius traces}
  
 Now we examine the trace of Frobenius statement from \eqref{LsheafLfunctionresults}, restricting, of course, to the finite context.
 Our conventions about Frobenius morphisms have been
 given in \S \ref{Frobaction}.  
 In the case $d>0$ we are going to be computing Frobenius traces only on the duals of the uncompleted spaces, see
Remark \ref{completion caveat}; we anticipate these to be the relevant quantities that will match with
the computation of automorphic inner products.

 For \eqref{LX1} we note that the external shift by $[-d]$ has no effect
 except multiplication by a factor $(-1)^d$. We will therefore compute without it. In what follows, the cohomology space $H^1$ is {\em considered as a vector space in cohomological degree $1$} (in particular, its symmetric powers, in the graded sense, are exterior, in the ungraded sense).  
 We assume,
 at first, that the $\GGm$ action on $T$ is trivial.
For \eqref{LX1}, we have:
 \begin{multline*} \mathrm{tr}(\Fr |  \Sym  \ H^1(T)^\vee) = \det(1-\Fr| H^1(T)^{\vee}) \\ =
    \det (1-q^{-1} \Fr| H^1(T^{\vee})) = L(1, T^{\vee}),\end{multline*}
 since by assumption $T$ has no $H^0$ or $H^2$
 and we have a perfect pairing $H^1(T) \times H^1(T^{\vee}) \rightarrow \kk(-1)$,
 where $\Fr$ is acting by $q^{-1}$ on the target. 
 For \eqref{LX0}, we use 
   $$ \mathrm{tr}(\Fr |  \Sym \ H^1(T)) = \det(1-\Fr| H^1(T)) = L(0,T),$$
  using again vanishing of $H^0$ and $H^2$.   
  We additionally 
use the fact that the determinant of Frobenius on $H^1(\ad \rho)$ equals $q^{b_G}$, by \eqref{epsform} and \eqref{firstepsilondet};
indeed, the relevant orthogonal $\epsilon$-factor is trivial since $\ad \rho$  has trivial determinant;
moreover, the dimension of $H^1(\ad \rho)$ is even, so there is no sign shift. 

The modifications for nontrivial action of $\GGm$ are immediate, as the shears ``come along for the ride.''
  
  For \eqref{LX1b} and \eqref{LX0b} we use  the fact that $\varepsilon_{1/2}$
  geometrizes the square root of the central $\epsilon$-factor, see discussion before \eqref{epsilon0description}. 
  We will now pay more attention to the shears when computing. 
Decompose $T = \bigoplus T_k$ according to the $\GGm$-grading. We get
\begin{equation} \label{epsilon0k} 
\epsilon(0, T^{\shear})= \prod \epsilon(k/2, T_k) 
  \stackrel{\ref{epsform}}{=}
  \epsilon(1/2, T) q^{(  (g-1) \alpha}, \end{equation}

 \begin{equation} \epsilon(1, T^{\vee\shear})  \stackrel{\eqref{epsform}}{=}  \epsilon(1/2, T^{\vee}) q^{-(g-1)\alpha}\end{equation}
   where $\alpha =  \sum   (1-k) \dim(T_k) = \dim(X) - \gamma$,
   because $\GGm$ is acting on $\det(T)$ by the character $\sum k \dim(T_k)$,
   and this must equal $\gamma$, cf. the discussion preceding \eqref{gammadef}.    For \eqref{LX1b} we reason thus using \eqref{sigmaXdef} (and recalling that
  we are dualizing at the start!) 
  {\small 
\begin{eqnarray*}  [ \left( \left(\Sym^* H^*(T)\right) [-d] \otimes \varepsilon^{\vee}_{1/2}(T)\langle -\beta_X  \rangle  \right)^{\vee}] = 
 (-1)^d q^{-\beta_{\check{X}}/2} [\epsilon_{1/2}(T)] L(1, T^{\shear\vee})   \\
 = (-1)^d q^{-\beta_{\check{X}}/2} [\epsilon_{1/2}(T^{\vee})]^{-1} L(1, T^{\shear\vee})  
  \\ =
 (-1)^d   q^{-\beta_{\check{X}}/2} q^{-(g-1) (\dim X-\gamma)/2} \sqrt{\epsilon}(1, T^{\shear \vee})^{-1} L(1, T^{\shear \vee}) 
 \\ \stackrel{\eqref{Lnormdef}}{=} (-1)^d q^{-b_G/2} L^{\norm}(1, T^{\shear \vee}),\end{eqnarray*}}
 where we used $\beta_{\check{X}} + (g-1) (\dim X-\gamma) = (g-1) \dim G$,
 see \eqref{betaXdef}, and the square root of $\varepsilon$ is chosen based on the fixed spin structure, as in \eqref{epsilondef}.
 For \eqref{LX0b} we reason thus:

\begin{multline} 
 [ (\mbox{ghastly mess})^{\vee}] =q^{-b_G} q^{\beta_X/2} \epsilon(1/2, T)^{-1} L(0, T^{\shear}) 
 \\ = q^{-b_G} q^{\beta_X/2} q^{(g-1) (\dim X-\gamma)/2} L^{\norm}(0, T^{\shear}) =q^{-b_G/2} L^{\norm}(0, T^{\shear}). 
\end{multline}

  \qed

\subsubsection{$L$-functions, algebraic distributions and volume forms}\label{derived volume forms}\index{algebraic distributions}
We briefly sketch a geometric point of view on $L$-functions suggested by $L$-sheaves and the theory of categorical traces, as an interesting direction for future research.

 In summary, we formally expect applying the trace of Frobenius to the $L$-sheaf $\Ll_\Xv$ defines an algebraic distribution -- the ``$L$-distribution'' on the fixed points of Frobenius on $\Loc_\Gv$, i.e., the stack of arithmetic local systems. This trace is only defined (at best) after restricting to an open locus where we require in particular that fixed points are isolated, i.e., away from poles of the $L$-function. The $L$-distribution can be thought of as the relative form of the volume of the derived Galois fixed points on $\Xv$. (See Remark \ref{AGKRRV Arthur} for a parallel discussion of Arthur parameters.)

First let us recall some of the functoriality of traces and fixed points in derived algebraic geometry, as explained in~\cite{nonlinear,HSS,GKRV}. It is a consequence of the functoriality of ind-coherent sheaves that a morphism $f:Z\to Z$ defines a  vector space (in the derived sense, so, strictly a chain complex): the {\em categorical trace} of $f_*$
as an endomorphism of $\QC^!(Z)$. This  vector space is identified with algebraic distributions -- derived section of the dualizing sheaf -- on 
the derived fixed points of $f$,
$$\mathrm{Tr}(f_*)\simeq \Gamma(Z^f, \omega).$$ 
Furthermore a coherent sheaf $\cF$ (i.e., a compact object of $\QC^!(Z)$) equipped with an $f$-equivariant structure defines a trace distribution  -- an object
of the vector space $\mathrm{Tr}(f_*)$.

This setup is realized in the Betti setting as follows. We consider a topological surface $\Sigma$ with a diffeomorphism $F:\Sigma\to \Sigma$, with mapping torus a 3-manifold $\Sigma_F$ (which plays the role of a global field in the classical Langlands correspondence).
We will apply the foregoing discussion to 
$$Z=\Loc^B_\Gv(\Sigma)$$
  the Betti moduli of local systems on $\Sigma$,
  and we take the map $f$ to be that induced by $F$.
  The derived fixed point locus   
  is the stack of local systems $\Loc^B_\Gv(\Sigma_F)$ on the 3-manifold $\Sigma_F$.  
  The $L$-sheaf attached to a $\check{G}$-space $\check{X}$ gives an $F$-equivariant sheaf on $Z$.
  It
  is not in general coherent, but it will be
  so if we restrict to any subset where 
$Z$ is smooth {\em and} the  map
$\Loc^{\check{X}}(\Sigma) \rightarrow \Loc_{\Gv}(\Sigma)$
is finite.  So we get  an algebraic distribution at
least on the corresponding subset of
$\Loc^B_{\Gv}(\Sigma_F)$.
  
    This algebraic distribution is   a geometric  avatar of the $L$-function.
  To make a clearer connection
  to $L$-functions in the arithmetic sense,  we would like to apply this formalism in the setting of the \'etale geometric Langlands conjecture~\cite{AGKRRV1}, as reviewed in \S \ref{geometric to arithmetic Langlands}. Namely, we would like to take $(Z,f)$ to be the formal algebraic stack $\Loc_\Gv^{\et}$ of restricted local systems with its Frobenius action, so that $Z^f=\Loc_\Gv^{\textrm{arith}}$ is the stack of arithmetic local systems. The formalism of traces is applied in this setting in~\cite[24.7]{AGKRRV1} and the categorical trace of Frobenius is identified as expected with algebraic distributions $\Gamma(\Loc_\Gv^{\textrm{arith}},\omega)$.  Thus the formalism
  above   produces an $L$-distribution in this space, 
but, again, only after restricting  to a suitable subset on which the $L$-sheaf is coherent (rather than merely ind-coherent).
  We leave a clearer understanding of this subset, and what happens away from it,
  as a question for later study.

The discussion above is  related to the well known analogy between $L$-functions and Reidemeister torsion and its generalizations (see, for example, ~\cite{AV}, where this analogy plays a key role). 
A  particularly relevant perspective on the  theory of Reidemeister torsion volume forms has been developed in~\cite{naefsafronov}: it provides various mapping spaces in derived algebraic geometry with canonical volume forms -- i.e., sections of the determinant line bundle of the cotangent complex.

\subsection{Reduction to the vectorial case}
\label{Linduction} 
  \label{Whittaker functoriality} 
  
We now consider the spectral counterpart of \S\ref{Pinduction}:
$L$-sheaves in general can be expressed as  ``spectral Whittaker inductions'' (or better, Arthur inductions) of $L$-sheaves
in the vectorial case. 
  We set up this functor as in the automorphic case, with the spectral exponential sheaf
 replacing the Artin--Schreier sheaf. 
 
 Our notation will be parallel to that of \S \ref{Pinduction}. Thus we fix a homomorphism $\Hv\times SL_2\to \Gv$ with underlying cocharacter $\varpi:\Gm\to \Gv$, which we assume to have even weights on $\fgv$. 
Let $U=U_+\subset \Gv$ be the unipotent subgroup defined by the positive part of the grading.
Just as in \S \ref{Pinduction} we then consider the graded Hamiltonian $\Gv\times \Hv$-space $T^*_\Psi \Gv/\check U$ associated to 
\begin{equation} \label{XvLindef}(\Xv =U \backslash \Gv,\Psi:\ U\longrightarrow \AA^1)\end{equation} and its $L$-sheaf $\mathcal{L}_{\Gv/\check U,\Psi}\in QC^!(\Loc_\Gv\times \Loc_\Hv)$ following Definition~\ref{spectral Whittaker def}.

\index{Arthur induction}
\begin{definition}\label{Whittaker induction def} 
(Compare with Definition
 \ref{automorphic Whittaker induction def}): \index{$\mathsf{AI}$ Arthur induction}
Let notation be as above. The Arthur induction functor $$\mathsf{AI}: QC^!(\Loc_\Hv) \longrightarrow QC^!(\Loc_\Gv)$$ is the $(!-)$integral transform given by the spectral Whittaker $L$-sheaf:
$$\mathsf{AI}(\cF)=\pi_{1*}(\pi_2^! \cF\otimes^{!} \mathcal{L}_{\Gv/\check U,\Psi}).$$
\end{definition}

Let us spell this out. First note we are applying throughout the general functoriality of $QC^!$ from~\cite{GR} (see \S \ref{IndCoh section} for a summary), in particular the $!$-tensor product for which $!$-pullbacks are symmetric monoidal and pushforwards satisfy the projection formula~\cite[III]{GR}.

The stack $\Loc^\Xv$ in the setting of \eqref{XvLindef}  is identified with $\Loc_{\Hv U}$ with its natural map to $\Loc_\Gv\times \Loc_\Hv$. 
 Then using the morphism $$\xymatrix{\overline{\Psi}_C:\Loc_{\Hv U}\ar[r]^{\Psi_C}&  \Loc_{\AA^1}\ar[r] & \AA^1[-1]},$$ we can identify the Arthur induction as 
$$\mathsf{AI}(\cF)=\left(q_*^\shear( (\overline{\Psi}_C^!)^\shear (\exp)\ot^! p_\Hv^!\cF)\right)^\unshear.$$
with $q$ the projection to $\Loc_{\Gv}$, and $p_{\Hv}$ the projection from $\Loc_{\Hv U}$ to $\Loc_{\Hv}$.

 The Arthur induction of $\omega_{\Loc_\Hv}$ recovers the $L$-sheaf associated to $(\Xv=\Gv/\Hv \check U,\Psi)$, and more generally Arthur induction of (polarized) L-sheaves realizes Whittaker induction on Hamiltonian spaces:

\begin{lemma}  \label{1022}
Given a homomorphism $\Hv\times SL_2\to \Gv$ (with even $SL_2$) and $S=T^*\Yv$ a polarized Hamiltonian $\Hv$-space, the $L$-sheaf of the Whittaker induction $(\Xv=\Yv\times^{\Hv U} \Gv,\Psi)$ of $S$ is naturally identified with the Arthur induction of the $L$-sheaf of $S$:
$$\mathsf{AI}(L_{\Yv})\simeq L_{\Xv,\Psi}.$$ 
\end{lemma}

\begin{proof}
The identification comes from considering locally constant maps from $\Sigma$ into the following commutative diagram with Cartesian diamond:
 $$\xymatrix{  & pt/\Ga & \ar[l]\Xv/\Gv \ar[dl]\ar[dr]&\\
  & pt/\Hv U \ar[u]\ar[dr]\ar[dl]&&\Yv/\Hv\ar[dl]\\
  pt/\Gv&&pt/\Hv&}$$
\end{proof}

\begin{remark}[Arthur restriction]  \label{ArthurResRemark}
As in \S \ref{Pinduction} we can also use the spectral Whittaker sheaf as an integral transform in the opposite direction
to define an Arthur restriction (or ``Arthur-Jacquet'') functor $$\mathsf{AJ}: QC^!(\Loc_\Gv) \longrightarrow QC^!(\Loc_\Hv)$$ which (by the same argument) takes the $L$-sheaf for a polarized $\Gv$-space $\Mv$ to that of the (twisted-polarized) Hamiltonian $\Hv$-space given by its reduction $\Mv\GIT_\psi U$.
It would be interesting, in particular in light of the heuristic discussion of \S \ref{Arthur functoriality}, to verify if $\mathsf{AJ}$ is identified with the left adjoint of $\mathsf{AI}$. 
\end{remark}

  \subsection{Independence of polarization.}\label{Lindep}
  We now discuss the spectral counterpart of \S\ref{Pindep}. An ultimate goal is to reformulate our entire study in terms of $(\check{G}, \check{M})$ and as in \S \ref{Pindep} we can think of this in two parts: independence of polarization (a form of the functional equation for $L$-functions), and constructing the $L$-sheaf even if a distinguished polarization does not exist.
\index{functional equation for $L$-function}

 We will now examine independence of polarization in the Betti case
and we shall prove  independence of polarization after projection via $QC^{!} \rightarrow QC$
  under condition \eqref{condB}. These restrictions are to make the discussion as simple and explicit as possible: it is likely
  that a suitably framed version of the same argument works in the other contexts, applies to $QC^{!}$,
  and does not require \eqref{condB}.

   The reader can compare these points with the corresponding constructions in \S \ref{Pindep}. From the point
of view of the classical theory of periods, the parallel between these situations is quite surprising.  
As in the automorphic case, we will only briefly discuss here the unpolarized situation, 
\S \ref{specWeil-start}.

We now put ourselves in the polarized vectorial setting, thanks to the discussion of \S\ref{Linduction}.
We are also going to assume that  
\begin{equation} \label{condB} \mbox{ $\check{G}$ is semisimple and the genus is $\geq 2$. } 
\end{equation}
There should be no difficulty in removing this assumption (in fact, the corresponding assertion in genus $0$
follows from our computations in \S \ref{P1}). 

Under this assumption and working in the Betti context we are going to construct an isomorphism (as mentioned, inside $QC$):
\begin{equation} \label{CandC}\mathcal L_{\check{X}} \simeq \mathcal L_{\check{X}^{\vee}} \otimes \det H^*(\check{X})^{-1} \end{equation}
\begin{equation} \label{CandCnorm} \mathcal L_{\check{X}}^{\norm}
 \simeq \mathcal L_{\check{X}^{\vee}}^{\norm}\end{equation}
 where $\det H^*(\check{X})$ denotes the determinant of cohomology, considered as a line bundle
 on $\Loc_{\check{G}}$.  
 These isomorphism could be seen, in light of \S \ref{locXex}, as categorifying the functional equation of the $L$-function,
 but the situation is a bit more subtle: the point is that the $L$-sheaf is only a reasonable
geometrization of the $L$-function over a certain nice locus in the moduli space, 
but this ``categorical functional equation'' for the $L$-sheaf nonetheless holds everywhere.

 We will prove \eqref{CandC} in a very straightforward way: we write down a presentation for the left hand side that is manifestly invariant under duality
(up to the twist).  Before we do that, let us see how to go from \eqref{CandC} to \eqref{CandCnorm}.
Note that  $\beta$ for $\check{X}$ and $\check{X}^{\vee}$ coincide by definition, see \eqref{betaXdef}. On the other hand,
by definition (\S \ref{LsheafXnormalized}), 
$\varepsilon_{1/2}$ for $\check{X}$ and for $\check{X}^{\vee}$ are pulled  
back from the same bundle $\varepsilon_{1/2}^{\vee}$ via (respectively) $\check\eta$
and $\check\eta^{-1}$.
Taking into account (see Remark  \ref{det of coh remark} for the argument in the harder de Rham case)  that $\varepsilon_{1/2}^{\vee}$
on $\Loc_{\Gm}$ is a square root of the determinant of cohomology line bundle, we get 
$$ \varepsilon_{1/2, \check{X}}^{\vee} \otimes \varepsilon_{1/2, \check{X}^{\vee}}
\simeq \mbox{determinant of cohomology for $\check{X}$},$$
and so from \eqref{CandC} we will get the desired
invariance under duality for the normalized $L$-sheaf \eqref{CandCnorm}.

\begin{remark}

     The identity  that we are about to see is analogous to the following simple numerical fact. Given a complex of finite dimensional
 vector spaces $V$ with endomorphism $F$,  the formal $L$-value
$L(V) := \prod_k  \det(1- F | H^k V)^{(-1)^{k+1}}$
can be categorified in two different ways,  namely both by
$\Sym V \mbox{ and } \Sym V^\vee \otimes (\det V^\vee)$ -- the determinant
is regarded as placed in cohomological degree equal to the dimension of $V$. 
This amounts to the identity
 $$ \sum_{k \geq 0} x^k= \frac{1}{1-x} =  -\frac{x^{-1}}{1-x^{-1}} = - \sum_{k=-\infty}^{-1} x^{k}$$   
 Nonetheless,
the reader may like to convince themselves before reading further that \eqref{CandC} is not a {\em formality}, i.e,
the definitions of both sides do not obviously align.  \end{remark}

\subsubsection{The vectorial $L$-sheaf}

Let $R$ be the ring of functions on homomorphisms from $\pi_1(\Sigma)$ to $\check{G}$,
as in \S \ref{locGprez}.  
It is an affine complete intersection ring; we have
$$ \Loc_{\check{G}}= \Spec(R) / \check{G},$$
and sheaves on $\Loc_{\check{G}}$ are then identified with $\check{G}$-equivariant
sheaves on the spectrum of $R$. It is in this model that
we will compute the $L$-sheaf.

As we are supposing the genus of $\Sigma$ is $2$ or greater,
this $R$ is a {\em usual} (underived) ring, which is even a complete intersection;
some   readers may follow the example of some  authors and  sigh with relief.

We have a representation
$$ \pi_1 \rightarrow \check{G}(R) \rightarrow \GL_{\check{X}}(R),$$
i.e., a local system of free $R$-modules on $\Sigma$. 
Proceeding, a cell decomposition of the Riemann surface associated to $\Sigma$
gives a complex of free $\check{R}$-modules that computes the cohomology of $\Sigma$
with coefficients in this local system; explicitly,
we lift this cell decomposition to a universal cover $\tilde{\Sigma}$
 and take the complex  of $\pi_1$-invariant cochains. 
 The result is a complex of finite rank free $R$-modules 
 in degrees $0,1,2$:
 \begin{equation} \label{Cdef} \mathcal{C}:  Q  \rightarrow Y \rightarrow P\end{equation}
with bases indexed by $0$, $1$- and $2$-cells; differentials are the usual differentials of singular cohomology, but twisted
by the universal local system.

We can think of $Q, Y, P$ as sections of vector bundles over $\Spec(R)$, 
to be denoted by the same letters.  
The pullback of $\Loc^{\check{X}}_{\Gv}$ to $\Spec(R)$
is  ``the kernel of the map
 from $Q$ to $\mathrm{ker}(Y \rightarrow P)$,'' except we interpret this in the derived sense.
 To make this more precise,
 let us call  $\widetilde{\Loc}_{\Gv}^{\check{X}}$  this pullback of $\Loc_{\Gv}^{\check{X}}$ to the spectrum of $R$, 
 so that we have a diagram, where the horizontal arrows take quotient by $\Gv$: 
$$
\xymatrix{
\widetilde{\Loc}_{\Gv}^{\check{X}}   \ar[d]^{\pi} \ar[r] & \Loc_{\Gv}^{\check{X}}  \ar[d]^{\pi} \\
\Spec(R) \ar[r] & \Loc_{\check{G}}
}
$$

Then $\widetilde{\Loc}_{\Gv}^{\check{X}}   $
is the spectrum of the dual symmetric algebra to the complex $\mathcal{C}$ of \eqref{Cdef}. 
Above, $Q, Y, P$ are in degrees $0,1,2$ respectively; the dual symmetric algebra is thereby in degrees $ \leq 0$;
ignoring differentials, it is identified with
$$ \Sym \mathcal{C}^{\vee} = \Sym \  Q^* \otimes \Sym \  Y^*[1] \otimes \Sym \ P^*[2],$$
where the symmetric algebras are taken in the graded sense so that, for example, $Y$ is in fact contributing an abstract exterior algebra.
 Also the $\GGm$-action scales $Q,Y,P$ by the tautological character and acts dually on the ring of functions. 
 
 To compute the $L$-sheaf, we will factorize the map $\pi$ as 
 $$ \widetilde{\Loc}_{\Gv}^{\check{X}} \stackrel{\pi_0}{\rightarrow} Q \rightarrow \mathrm{Spec}(R)$$
 where $Q$ is now considered as a vector bundle over the spectrum of $R$. 
   Let 
  $$ \omega = \mbox{dualizing sheaf of $\widetilde{\Loc}_{\Gv}^{\check{X}}$}.$$ 
We have \begin{equation} \label{LsheafOmega} \mathcal{L} = \Omega_{\check{X}}[-\dim G]^{\shear}, \end{equation}
where
$$ \Omega_{\check{X}} := \mbox{ sections of $\pi_* \omega$ on $\Spec R$},
$$  as a $\check{G}$-equivariant
 $R$-module.  This $\Omega_{\check{X}}$ is equivalently
 described as sections of $\omega$ on $\widetilde{\Loc}_{\Gv}^{\check{X}}$,
 i.e., as homomorphisms $\mathcal{O} \rightarrow \omega$ computed there, 
and this can bee computed by  adjunction along the proper
map $\pi_0: \widetilde{\Loc}^{\tilde{X}} \rightarrow Q$:
  $$\Omega_{\check{X}}= \Hom_{\widetilde{\Loc}_{\Gv}^{\check{X}}}(\mathcal{O}, \pi_{0}^! \omega_Q) = \Hom_{Q}(\pi_{0*} \mathcal{O}_{\widetilde{\Loc}_{\Gv}^{\check{X}}}, \omega_Q).$$
  Here, as above, we are by an abuse of notation regarding $Q$
as a vector bundle over $\Spec(R)$. 
Now, the dualizing sheaf of $\omega_Q$ is identified with 
$$ \omega_Q \simeq  \omega_{\mathrm{Spec} \ R} \otimes \Sym(Q^{\vee}) \otimes (\det Q)^{\vee}[\dim Q]$$
(the fact that $R$ is LCI implies that $\omega_{R} := \omega_{\Spec  R}$ is in fact locally free of rank $1$)
and {\em ignoring differentials at the moment}
we have an identification of complexes

$$ \omega_R^{-1} \Omega_{\check{X}} \otimes (\det Q)[-\dim Q] \simeq
\Hom_{\Sym \ Q^{\vee}}(\Sym \ \mathcal{C}^{\vee}, \Sym \ Q^{\vee} )$$
which, after shearing everything in sight, becomes
$$ \omega_R^{-1} \Omega_{\check{X}}^{\shear} \otimes (\det Q)
\simeq \Hom_{\Sym Q^{\vee}[-1]}( \Sym \mathcal{C}^{\vee}[-1], \Sym \ Q^{\vee}[-1]).$$
Recall here that the shearing includes a parity shift on $\mathcal{C}^{\vee}$,
so that the shear of the symmetric algebra coincides with the symmetric algebra of the shear,  
that is to say, there is no funny business switching between exterior and symmetric algebras upon shearing. 
Now $\mathcal{C}^{\vee}[-1]$ is the complex  $P^{\vee} \rightarrow Y^{\vee} \rightarrow Q^{\vee}$ in degrees $-1,0,1$, and
 {\em forgetting differentials} we get 
 \begin{equation} \label{X01}  \omega_R^{-1} \Omega_{\check{X}}^{\shear} \otimes (\det Q) \simeq 
\wedge^*Y \otimes \Sym (P[-1] \oplus Q^{\vee}[-1]) .\end{equation}
Here, to avoid unwarranted suffering,  we write simply $\wedge^* Y$ where $Y$
is now regarded as a usual module in even parity and degree zero,  in place of than $\Sym Y$ where  $Y$ is regarded as a fermionic module in degree zero.

By Poincar{\'e} duality, fixing an orientation, the corresponding construction, replacing $\check{X}$ by $\check{X}^{\vee}$,
arises from the complex $\mathcal{C}^{\vee}[-2]$, that is to say, $P^{\vee} \rightarrow Y^{\vee} \rightarrow Q^{\vee}$; and then we would similarly obtain
 after tensoring by the determinant of $Y$,  {\em again forgetting differentials}
\begin{equation} \label{X12}  \omega_R^{-1} \Omega_{\check{X}^{\vee}}^{\shear} \otimes (\det P^{\vee}) 
    \otimes (\det Y) \simeq  \ \wedge^*Y \otimes \Sym^* (P[-1] \oplus Q^{\vee}[-1]),\end{equation}
  
  Since the determinant of cohomology of $H^*(\check{X})$ is identified with $(\det Q) \otimes (\det P) \otimes (\det Y^*)$, 
 comparing \eqref{X01} and \eqref{X12} and checking the differentials match will shows the desired identification 
 $$\Omega_{\check{X}}^{\shear} \otimes (\det H^*(\check{X})) \simeq \Omega_{\check{X}^{\vee}}^{\shear}$$ of \eqref{CandC}.
 To verify that the differentials match, 
 we will write the differentials in \eqref{X01} in a way that is manifestly symmetric. To do 
 this we will write out things in a basis, which, in this case, shows a certain beauty of the construction.

 \begin{lemma} \label{symmetric_presentation}
  As above let $\mathcal{C}: Q \rightarrow Y \rightarrow P$  be a complex of free $R$-modules in degrees $[0,2]$
  and write $$d: Y \rightarrow P, d^*: Y^{\vee} \rightarrow Q^{\vee}$$
for the differential and dual differential from \eqref{Cdef}. 
 Here $Y^{\vee}=\Hom(Y, R)$ and similarly for $Q^{\vee}$.   Fix a basis $x_1, \dots, x_k, \dots$ for $Y$ as an $R$-module with dual basis $x_i^*$ for $Y^{\vee}$.

 Then, with reference to the identification just made
\begin{equation} \label{explforml}  \Hom_{\Sym Q^{\vee}[-1]}( \Sym \mathcal{C}^{\vee}[-1], \Sym \ Q^{\vee}) 
 \simeq \wedge^* Y \otimes \Sym(P[-1] \oplus Q^{\vee}[-1])  \end{equation}
 the differential on the right-hand complex  may be characterized thus:
  
  It is   $\Sym(P[-1] \oplus Q^\vee[-1])$-linear and given
on a term $x_1 \wedge x_2  \wedge \widehat{x_3} \wedge \dots$ 
(where, as usual, a hat on top of $x_k$ means that it is omitted) by the following rule:
  \begin{itemize}
 \item replacing any term $x_j$ by $\pm \widehat{x_j} (dx_j)$; here $dx_j \in P$; 
 \item replacing any term $\widehat{x_j}$ by $\pm x_j (d^* x_j^*)$; here $d^* x_j^* \in Q^*$,
  \end{itemize}
 where in both cases the sign $\pm$ is given by $(-1)^{j-1}$.
 
 \end{lemma}
 
  The presentation of the differential is manifestly symmetric under duality, concluding the proof of \eqref{CandC}.

 \proof
 We just need to write out \eqref{explforml} and compute. 
 \qed

 \subsubsection{The spectral Weil representation}
  \label{specWeil-start}
  The isomorphism \eqref{CandC} can be seen as a kind of ``spectral Fourier transform,''
and we may expect it to be related to the Fourier transform that appears in the Weil representation;
indeed, that
 \begin{quote} {\em the $L$-sheaf in the vectorial case
is characterized, up to a projective ambiguity, as the unique  irreducible representation of a certain algebra}.\end{quote}
i.e., a kind of Stone-Von Neumann theorem. 
This is indeed the case, but since it will take us a little way from the central concerns
of this paper,   we will report on the construction in 
separate work.
  This points us a way towards both a more intrinsic way of understanding the self-duality, and more broadly a way
 of understanding $L$-sheaf.  
 This is closely related to the
 ``algebra of $L$-observables'' that is discussed in \S \ref{case5} and at length in \S \ref{spectral geometric quantization}.

 Let us at least see where the algebra comes from. 
  In the situation above, let  
\begin{equation} \label{VVdef} \mathsf{V} = \mathcal{C}[1] \oplus \mathcal{C}^\vee[-1], \end{equation}
where $\mathcal{C}$ is as in \eqref{Cdef} a complex of sheaves on $\Loc_{\check{G}}$ in degrees $[0,1,2]$
computing the cohomology with coefficients in $\check{X}$;
the complex  $\mathsf{V}$
 is represented by a complex in degrees $[-1,1]$. 
$\mathsf{V}$ has moreover a self-duality structure for it computes the cohomology of  the $\check{M}=T^* \check{X}$-local system associated to a
  $\check{G}$-local system, just as $\mathcal{C}$ computed the cohomology of the $\check{X}$-local system,
  and the self-duality structure arises from Poincar{\'e} duality pairing on cohomology.
 In particular, the resulting pairing  $\mathsf{V} \rightarrow \mathsf{V}^{\vee}$ gives $\kk \oplus \mathsf{V}$ the structure of a
 a differential graded Lie superalgebra.    Moreover, this Lie  superalgebra acts on the $L$-sheaf constructed above
 and can be used to characterize it.

%% file: global-geometric-2.tex
\newcommand{\nilPX}{\overline{P_X^{\norm}}}
\section{The global geometric conjecture} \label{GGC}

 We will formulate the global geometric conjecture in multiple versions.
 The general nature of the conjecture is that
 \begin{quote} the normalized period sheaf
 attached to $(G, M)$ matches under geometric Langlands with the normalized
 $L$-sheaf attached to $(\check{G}, \check{M})$ for a hyperspherical dual pair.
 \end{quote}

  Before we go further, it may be helpful
 to say specifically the `data points' on which are conjectures are based. The most important are:
 \begin{itemize}
 \item[(i)] agreement with known numerical statements, to be discussed in  \S \ref{L2numerical},  \item[(ii)] the case of $\mathbb{P}^1$, to be discussed in \S \ref{P1}.
 \end{itemize}
 What we present below should be regarded as a first attempt to formulate global conjectures
 that are consistent with these data points and the known formalism of geometric Langlands.
We expect that, especially in regard to technical details concerning the appropriate categories (see e.g.\ \S \ref{nilpotent projection is good}), they may need further modification.

\begin{itemize}

\item  In \S \ref{case1} we formulate the conjecture in the case where $M, \check{M}$ are both polarized (after discussing in what categories to place the period and $L$-sheaves).

\item In  \S \ref{6examples} we examine a few simple examples of \S \ref{case1}. 
 The examples we present here are the much more restricted class where one can
compute explicitly with the sheaves, and should be regarded more as illustrations than as evidence. 
\item In \S \ref{groupcase} we give a somewhat extended discussion of the
``group case,'' and reformulate some cases of the conjecture in terms of functorial transfers
between different groups. 

\item  In \S \ref{nilpotent projection is good} we discuss the role of spectral projection in the Betti and \'etale forms of the global geometric conjecture.

\item \S \ref{parity returns}  -- \S \ref{Parity Examples}
are an investigation of the role of parity  (see \S \ref{analyticarithmetic} and also \S \ref{parity}).
In \S \ref{parity returns} we explain why the parity condition 
 \S \ref{parity} implies that the conjecture is independent of the choice of spin structure. 
 \S \ref{mitch} is probably the most detailed test of the parity condition. We  first explain
 how to shift between normalized and unnormalized version of the conjecture. This in particular leads
 in 
 \S \ref{Parity Examples}
 to another implication of the conjecture related to parity, \eqref{parity prediction}, which we verify by hand in many examples

 \item 
In \S \ref{case5} we study the situation with polarized $M$ and unpolarized $\check{M}$. 
 We achieve this generality at the cost of giving a less precise conjecture, describing
 not the $L$-function but its ``square'' (this is familiar in the automorphic literature). 
 The version given in this section
 works only over a certain subset of the space of local systems;  we discuss how to extend this form of the conjecture  to the whole space in the concluding section \S
\ref{spectral geometric quantization}
 of this paper. 

 \end{itemize}

\begin{remark}  
  \label{metastuff}
 In a certain sense, the conjecture that follows is a geometrization of the 
  ``$L$-value equals period'' philosophy well known in number theory.  Let us describe
  how one would arrive at the conjecture, starting from that philosophy;
note that  the purely numerical version of the conjecture is studied 
  in \S \ref{L2numerical},
  and the relationship between geometric and numerical is discussed in \S \ref{mumerical}.

  \begin{itemize}
  
  \item[(a)] Firstly, one must decide {\em which} $L$-function
  is attached to a given period. The guiding philosophy of our paper
  is that this fits into the duality framework of  hyperspherical varieties.

  \item[(b)] Next, one must understand the correct normalization of the conjecture; i.e., 
describe it in a way where the ``fudge factors'' that appear
in a typical ``$L$-value equals period'' formula come from.    This shows up, for us, in the precise twists
and shifts in the definition of normalized $L$- and period sheaves. 

\item[(c)] Thirdly, even in 
very simple settings such as the Hecke or Iwasawa--Tate period,
it is not in fact true that the relation ``$L$-value equals period''
holds when paired against all automorphic forms:
this is related to subtleties of Eisenstein periods where the $L$-function has a pole (see Remark \ref{interesting}). This is a feature that is not even understood well in the numerical world. 
The correct conjecture must 
reproduce this behavior. This  point was emphasized to us by Tony Feng and Jonathan Wang, who have
studied the Hecke case in unpublished work. We do not discuss this point in isolation
in this paper, although it comes up implicitly a few times; however,
to the extent that we have tested it on this front,  the conjecture passes. 

\item[(d)] Finally, in geometrizing
numerical statements, we must decide where to use the dualizing sheaf,
where to use the constant sheaf, and where to use something in the middle.
  For example, these choices (in the group case) must recover
the fact (``miraculous duality'') of Drinfeld--Gaitsgory that
duality on automorphic and spectral sides are related in a somewhat subtle way; see \S \ref{diagdiag}.

 \end{itemize}
\end{remark}
 
\subsection{Normalized period conjecture: statement}\label{case1}

We follow the general notation set up in \S \ref{notn} and \S \ref{Lsetup};
in particular, $\Sigma$ is a projective smooth curve over the field $\FF$,
and we consider sheaf theory with coefficients in $\kk$.

Before stating the conjecture, let us discuss the setting: what categories do the $L$- and period sheaves live in?
 As has been studied in the work \cite{ArinkinGaitsgory} of Arinkin and Gaitsgory in the de Rham setting (and more recently in~\cite{BettiLanglands} in the Betti setting and~\cite{AGKRRV1} in the \'etale setting), the precise formulation of a geometric Langlands conjecture is quite a subtle matter,
 because one must choose carefully which categories to work with on both sides. We summarize the different forms of the geometric Langlands correspondence available in different contexts in Appendix~\ref{geometric Langlands},
cf. \S \ref{BunLocintro}, \S \ref{BunLocintro2}. We recall here that the de Rham and Betti settings are formulated only over $\FF=$ the complex numbers;
for us, the \'etale setting will be applied with $\FF$ the algebraic closure of a finite field. (In fact   
   the \'etale setting applies in greater generality, providing a common language for the finite setting of $\ell$-adic sheaves and for the ``common part'' of the de Rham and Betti contexts over $\FF=\C$.)

 One of the curious features of our $L$- and period- sheaves is that, although 
 both definitions are very natural and parallel to one another, 
  {\em they do not always live in the standard categories in terms of which geometric Langlands has been formulated.}

\begin{itemize}
\item[-]  The geometric Langlands correspondence, as reviewed in \S \ref{geometric Langlands}, only takes as input {\em automorphic} sheaves on $\Bun_G$. In the de Rham setting, all sheaves are automorphic, but in the Betti and \'etale settings one restricts to {\em nilpotent} sheaves on $\Bun_G$, a class that includes Hecke eigensheaves but not in general period sheaves. This comes from the fact (Section~\ref{spectral action}) that the Hecke action is not automatically locally constant over the curve, so only a subcategory of all sheaves can spectrally decompose over $\Loc_\Gv$. As a result, to apply the correspondence we need to consider period sheaves as {\em functionals} on automorphic sheaves. Equivalently, to apply the correspondence to period sheaves we must first take them out of their native habitat and apply spectral projection $\mathcal P\mapsto \mathcal P^{spec}$, \index{spectral projection}
which is a functor 
\begin{equation} \label{spec proj} \bigautshv(\Bun_G) \rightarrow \bigautshvspec(\Bun_G)\end{equation}
from all sheaves to automorphic sheaves, i.e., the spectrally decomposable subcategory. The nature of the Betti and \'etale categories is substantially different (for example because $\Loc^B$ is close to an affine variety while $\Loc^{et}$ is close to an affine formal scheme) -- as we illustrate in some detail in Section~\ref{nilpotent projection is good}, in the former the spectral projector gives the {\em left} nilpotent projection (left adjoint to the inclusion), while in the latter it is the {\em right} nilpotent projection. Thus in the Betti setting the period functional is given by Hom {\em from} the period sheaf, while in the \'etale setting it is given by Hom {\em to} the period sheaf. 
\item[-] In all settings, the spectral side of the Langlands correspondence is usually formulated (as in~\cite{ArinkinGaitsgory,AGKRRV1}) for ind-coherent sheaves with nilpotent singular support on the corresponding stacks $\Loc_\Gv$ of Langlands parameters, a class that does not include most $L$-sheaves. To accommodate them we choose to work in all settings  with the larger (but less familiar)  \renormalized\footnote{which we could also call ``renormalized'' or ``cohomological''}
  form of the correspondence, whose spectral side consists of arbitrary ind-coherent sheaves. 
  
  As explained in Appendix~\ref{geometric Langlands} the distinction between the \renormalized and more standard ``safe'' forms of the Langlands correspondence, on the automorphic side, boils down to whether we treat sheaves on $BG$ ``homologically'', as modules for $H_*(G)$ (the ind-safe category, as is the norm in the setting of $D$-modules and thus for example in~\cite{ArinkinGaitsgory}),  or ``cohomologically'', as modules for $H^*(BG)$ (the ind-finite category). The cohomological (ind-finite) conjecture is strictly stronger than the homological (ind-safe) one but is not well documented or supported. In particular, the analogue of the projector \eqref{spec proj} is not documented in this setting, but we will proceed assuming it is valid;
  in case of any concern, all of our statements can be ``projected'' to the safe category, see \S \ref{safe Langlands}.

 \end{itemize}

 We now formulate the matching of period and $L$-sheaves in all three settings, using the spectral projectors that we will further discuss in \S \ref{nilpotent projection is good}. 
  Before stating the conjecture, we recall
  issues of rationality, only of relevance in the finite case:
  In \S~\ref{dhpFqdef}  we  specified
   a rather clumsy working definition of a ``distinguished split form of a hyperspherical dual pair'',
  where the two sides are defined over   fields $\FF,\kk$
  covering all our present options for coefficients.
   We will
  use this notion here, but it can be avoided  by setting up the situation up a little differently; 
see Remark \ref{rationality conjecture}.
  
\begin{conjecture}  \label{GlobalGeometricConjecture} 
Take $(G, M = T^*(X, \Psi))_{/\FF}$ and $(\check{G}, \check{M} = T^*(X, \Psi))_{/\kk}$
be a distinguished split form of a hyperspherical dual pair with distinguished polarizations, both sides admitting eigenmeasures\footnote{See Remark
\ref{noeigenmeasure2} for some first remarks on removing this.}.  Write
$$\mathcal P = \mathcal P_X^{\norm}, \mathcal L= \mathcal L_{\check{X}}^{\norm}$$
for the associated period and $L$-sheaves as defined in
\S \ref{section:global-geometric} and \S \ref{Lsheaf} respectively,  $(-)^{spec}$ for the spectral projection to automorphic sheaves, and $d$ for ``dualizing twist'', i.e., for the effect of applying the dualizing involution (\S \ref{dualityinvolution})
on $G$ or $\check{G}$. Then in all three settings we conjecture 
\begin{center} the spectral projection $\mathcal P^{spec}$ of the period sheaf and the dual $L$-sheaf $\mathcal L^d$ match under the geometric Langlands correspondence.
\end{center}
In the finite case, when $\Sigma$ is defined over $\FF_q$ and
assuming $(G, M=T^*(X, \Psi))$ to have a distinguished split form over $\FF_q$, this is moreover compatible 
with the  natural Frobenius-equivariant structures on both sides.
\footnote{On the left, this
arises from the definition \eqref{PXnormdef}, taking into account that the Tate twist modifies the equivariant structure; on the right,it similarly arising from the Frobenius action on $\Loc^{\check{X}}_{\check{G}}$
covering the action \eqref{Frobaction} on $\Loc_{\check{G}}$,  again taking account that the Tate twists implicit in \eqref{Lsheafdef} modify the equivariant structure.}
\end{conjecture}

Since the Beilinson spectral projector has very different features in the three settings, let us spell out separately what the conjecture amounts to:
  
\begin{description}
\item[de Rham] The period sheaf itself $\mathcal P=\mathcal P^{spec}$ and dual $L$-sheaf $\mathcal L^d$ match under the de Rham geometric Langlands correspondence of~\cite{ArinkinGaitsgory} -- or, rather, 
 its proposed ind-finite variant, see discussion above. 

\item[Betti] The {\em left} nilpotent projection of the period sheaf $\mathcal P^{spec}=\bar{\mathcal P}^B_l$ and dual $L$-sheaf $\mathcal L^d$ match under the Betti geometric Langlands correspondence of~\cite{BettiLanglands}.

\item[\'etale]
The {\em right} nilpotent projection of the period sheaf $\mathcal P^{spec}=\bar{\mathcal P}_r^{res}$ and dual $L$-sheaf $\mathcal L^d$ match under the \'etale geometric Langlands correspondence of~\cite{AGKRRV1}. 
   
\end{description}
``Left'' and ``right'' refer to left and right adjoints to the inclusion of
the subcategory of nilpotent sheaves, see \S \ref{igd2}.

   \begin{remark} 
   A few initial remarks on the statement:
   \begin{itemize}
   \item[(a)]
   The $P$-and $L$-sheaves are independent of polarization
 (see \S \ref{Pindep}, \S \ref{Lindep}), and correspondingly the conjecture depends only on $(M, \check{M})$.
 
 \item[(b)] At first sight it might appear unpleasant to have the dual $\mathcal L^d$ of the $L$-sheaf appear in the conjecture rather than the $L$-sheaf itself. 
Indeed, this is merely an issue of our (somewhat standard) choice of normalizations and can be reversed in one of several ways. One could change the normalization of the Hecke action in geometric Langlands, or switch the roles of arithmetic and geometric Frobenius in passing to the numerical conjecture, or absorb the twist in the duality between hyperspherical varieties, or switch left or right conventions... Of course,
each one of these choices would cause changes somewhere else.

 \item[(c)] Let us comment on the implied normalization of the geometric Langlands correspondence.
The above conjecture includes the matching between automorphic Whittaker and spectral dualizing sheaves up to an explicit twist,
which is, up to this explicit twist,  the standard Whittaker normalization of the correspondence, see \S \ref{Whittakertwist}. 
We have taken the point of view, however, that one should not privilege Whittaker normalization over the matching of other periods. 

\item[(d)] We expect that in the \'etale setting there is a natural pro-$L$-sheaf which matches the pro-automorphic sheaf given by the left nilpotent projection $\bar{\mathcal P}_l^{res}$ (a pro-version of the spectral projection of $\mathcal P$). This ``left'' version is important for matching the conjecture with numerical predictions, but we do not develop this version here.  
  \end{itemize}
    \end{remark}
\begin{remark}  \label{rationality conjecture}  (Alternate approach to rationality issues:)  The conjecture used the still tentative notion of distinguished split form in the case when $\FF$ has finite characteristic.
  This can be avoided  by taking {\em as starting point} a spherical variety over $\FF$, giving
  a conjecture that is in some ways more general:
    
  As usual, let $G$ be the split group over $\FF$, and take a spherical $G$-variety  $X$ over $\FF$, possibly endowed with a torsor $\Psi$, and without roots of type $N$.  
  Assuming that the characteristic of $\FF$ is sufficiently large,\footnote{It would be nice to have a more explicit statement here, for which
  we need to examine the various structural results that have been used in \S \ref{CheckMRat} and verify their validity away from an explicit set of bad charcteristics. This can at least be done straightforwardly in any specific case.}
the discussion of \S \ref{dualofX} can be applied directly in this situation 
 to obtain
  $(\check{G}, \check{M})_{/\kk}$, which we assume moreover to be hyperspherical,
  and to admit a distinguished polarization
  $\check{M} = T^*(\check{X}, \check{\Psi})$. In this setting, we may again conjecture  that 
 the period and $L$-sheaves of Conjecture
 \ref{GlobalGeometricConjecture}  match. Under further assumptions, one 
formulates a statement about Frobenius equivariance using the Frobenius action on $\check{M}$ that has been specified in \S \ref{CheckMRat}.

  \end{remark}

 \subsection{Some illustrative examples} \label{6examples}

We discuss a few illustrative examples. 
As we have already mentioned, the strongest evidence for the global conjecture of course comes from the
study of {\em numerical} examples, which we examine in later sections.  
Both for this reason and because of the technical difficulty of working out details in this situation,
we make no attempt at complete or rigorous treatment, simply sketching some computations
that seem to us  likely representative of the true situation.

To avoid repeatedly having to mutter ``spectral projection'' in each example, {\em let us restrict ourselves to the de Rham
context for this subsection \S \ref{6examples}.}

\subsubsection{The automorphic Whittaker model}
  \label{Whittakertwist}
Take $M = T^*_{\Psi}(U\backslash G)$ the Whittaker model
with dual $\check{M}=\mbox{pt}$.  The conjecture says that, under the Langlands correspondence,
the Whittaker sheaf  is exchanged with the dualizing sheaf of $\Loc_{\check{G}}$
(up to a shift that will be examined below). 
This is a known prediction (``Whittaker normalization''), due to Drinfeld. It often plays a distinguished role because it can be used
to normalize the Langlands correspondence; from our point of view it is only one of many matching pairs. 
The only point of relevance to examine is the normalizations,
which we now spell out. 
$\mathcal P_X^{\norm}$ differs from the ``standard'' Whittaker sheaf $\mathsf{W}$
by the shift of \eqref{PXnormdef}, which
is given here by \eqref{betaXdef} and \eqref{sigmaXdef}: 
 $$\mathcal P_X^{\norm} = \mathsf{W} \langle    (g-1) \left( \dim U - \left( \langle 2 \rho, 2 \check\rho \rangle  \right)  \right)   \rangle.$$
On the other hand, the $L$-sheaf is  
the dualizing sheaf now shifted by \eqref{LXnormdef}:
$$\mathcal L_{\check X}^{\norm} = \omega \langle - (g-1) \dim G \rangle$$
In other words, in our normalization, 
the usual Whittaker sheaf and the dualizing sheaf match after including a shift of $\langle Q \rangle$ on the spectral side with
$$ Q = (g-1) \left(  \langle 2 \rho, 2\check\rho  \rangle - \dim U - \dim G \right).$$
In other words, this is the shift by which the normalization implicit in 
Conjecture
\ref{GlobalGeometricConjecture}  differs from the standard normalization.

\subsubsection{The spectral Whittaker model} 
We now take $M =\mbox{pt}$ and $\check{M}= T^*_{\Psi}(\Gv/\check U)$ the Whittaker model. Recall that we have defined the corresponding $L$-sheaf in \S \ref{case4}. The conjecture says that, under the Langlands correspondence, the constant sheaf on $\Bun_G$ is exchanged with the ``spectral Whittaker sheaf on $\Loc_\Gv$'' again up to an explicit shift. 
The work of V.\ Lafforgue \cite{LafforgueP1} can be seen as proving a version of this statement
in the case of the projective line; see also \S \ref{P1}. 

\subsubsection{The Iwasawa--Tate case}  \label{Taterevisited}

  We consider the Iwasawa--Tate case: $M=T^*\mathbb{A}^1$, polarized by $X=\mathbb{A}^1$, as a $\Gm$-space,
with the ``neutral'' $\GGm$-action which acts by scaling on $M$
\footnote{Note, however, that $G=\Gm$ is not  acting by scaling on $M$: it scales $X$ and acts by inverse scaling in the cotangent direction.}
 We already described the geometry of $\pi: \Bun_\Gm^X \rightarrow \Bun_\Gm$ in \S \ref{TateBunX}.
            Let $z: \Bun_\Gm \rightarrow \Bun_\Gm^X$ be the zero section and $j: \Bun_\Gm^{X,0} \rightarrow \Bun_\Gm^X$ its complement;
            we have an exact triangle (denoting here by $k$ just the constant sheaf -- it should be obvious on which set) $ j_{!} k \rightarrow k \rightarrow z_* k$,             and correspondingly
           \begin{equation} \label{Verybasicsequence} \pi_{!}^{0} k \rightarrow \mathcal \Pcal_X  \rightarrow k.\end{equation}
           where $\pi^{0}$ is the restriction of $\pi$ to $\Bun_\Gm^{X,0}$. 
   Note also that the fibers of the map $\pi^{0}$  are $\Gm$-torsors
over projective spaces; and so $\pi^{0}_{!}  k \simeq \pi^{0}_* k [-1]$.
                       
       Taking account of the various twists to pass to the unnormalized sheaves, the conjecture amounts to the assertion that
 $ \Pcal_X \langle g-1 \rangle \mbox{ and } \Lcal_X^d$
 match under the Langlands equivalence (for a more general discussion of how to
 pass to unnormalized sheaves from their normalized versions, see Lemma \ref{shiftpred}). 
           Let us check this stalk-by-stalk on $\Loc_{\Gm}$, away from the trivial local system:

  Let $\rho$ be a {\em nontrivial} local system on $\Sigma$. 
   Associated to $\rho$, there is (cf. \S \ref{GLGMnormalization}) a  locally constant sheaf $\mathcal{G}_{\rho}$
      on the space $\Bun_{\Gm}(\Sigma)$, which has the property that 
      it pulls back to $\rho^{\boxtimes r}$ under $\Sigma^r \rightarrow \Bun_{\Gm}$
      sending $(P_1, \dots, P_r)$ to the line bundle $\mathcal{O}(\sum P_i)$. 
 The normalization of the Langlands correspondence relevant in the conjecture
 sends the perversely shifted version of $\mathcal{G}_{\rho}$, that is to say
 $\mathcal{F}_{\rho} := \mathcal{G}_{\rho} \langle g-1 \rangle$, 
 to the skyscraper sheaf $\delta_{\rho}$  on $\Loc_{\Gm}$ (i.e., the $*$ pushforward from 
 the point $\rho$).  We then compute:
\begin{multline} \label{ml1} 
\Hom(\mathcal{F}_{\rho}, \Pcal_X \langle g-1 \rangle) =  \Hom(\mathcal{G}_{\rho}, \Pcal_X) = \Hom(\mathcal{G}_{\rho}, \pi^{0}_*[-1] k )  \\ =
 \Hom((\pi^{0})^* \mathcal{G}_{\rho}, k)[-1]  = 
\left( \Sym H^*(\Sigma, \check\rho) \right) [-1]. 
 \end{multline}
 or strictly speaking a completion of this symmetric algebra,
 taking into account that the $\Hom$ is a product over components over $\Bun_{\Gm}$; 
 and by \eqref{LX1}
 \begin{equation}   \label{ml2} \Hom( \delta_{\rho}, \mathcal{L}_X) =   
 \mbox{ a completion of } \Sym^* H^*(\Sigma, \rho) [-1].
 \end{equation} 
 Taking into account the duality twist in the conjecture, which inverts $\rho$, we 
 see that \eqref{ml1} and \eqref{ml2} do indeed match.
 It is not difficult to carry this analysis out in families of $\rho$, except in the neighbourhood of trivial $\rho$,
 where the situation is more  interesting;
   Tony Feng, Jonathan Wang
 and one of the authors (A.V.) have verified that the conjecture
holds there too.

   \subsection{The group case and functoriality} \label{groupcase}
Continuing our study of examples, we now discuss the group case, and then explain how to utilize the group case to translate between the theory of periods and Langlands functoriality. We then briefly discuss a couple of instances of functoriality encoded in the global period conjecture.

\subsubsection{The group case and miraculous duality} \label{diagdiag}
Let us consider the group case: $M=T^*G$ as a $G \times G$-space.
Unwinding the definition in this case, we find that the un-normalized period sheaf $\cP_X\simeq \Delta_!\kk$ is the $!$-pushforward of the constant sheaf (equivalently: twisted dualizing sheaf) $\kk\simeq \omega\langle -2b_G\rangle$ of the smooth stack $\Bun_G$ under the diagonal embedding $\Delta:\Bun_G\to \Bun_{G\times G}$. 

To obtain the normalized period sheaf we further twist by $\beta_X=(g-1)\dim G$, i.e., $$\cP^{norm}\simeq \Delta_!\kk \langle b_G\rangle\simeq \Delta_!\omega \langle -b_G\rangle.$$

On the other hand, on the spectral side the dual Hamiltonian space is the  dualizing twist $\Mv=T^*\Gv$, i.e., the action of one copy of $\check G$ is twisted by the duality involution (\S \ref{dualityinvolution}).\footnote{Since we are taking $k$ to be algebraically closed here, the duality involution is the ``same'' (i.e., conjugate) as the pinned Chevalley involution.}  Thus, denoting by $\check\Delta$ the diagonal of $\Loc_{\check G}$, the un-normalized $L$-sheaf is
$\check\Delta_*\omega^{d}$ where the superscript $\mathrm{d}$ means that we twist the diagonal inclusion $\check\Delta$ (or, equivalently, twisting the sheaf $\omega$) by the dualizing involution in one factor.
The normalized $L$-sheaf is given by shifting this by $b_{\Gv}=b_G$:
$$\cLL\simeq \check\Delta_*\omega^{d}\langle -b_\Gv\rangle.$$
Thus, the global conjecture in the group case is the prediction that under the geometric Langlands correspondence for $G\times G$, the sheaf $\Delta_!\omega$ corresponds to the dualizing twist of $\check\Delta_*\omega$. 
 
 We now recall from \S \ref{duality section} the notion of duality for categories, and from \S \ref{tensor products of sheaves} and \S \ref{tensor and miraculous} its explicit realization for categories of ind-coherent sheaves. Namely, the sheaf $\Delta_*\omega\in QC^!(X\times X)\simeq QC^!(X)\otimes QC^!(X)$ for any QCA stack (quasicompact with affine diagonal, such as the stacks of local systems of either de Rham or Betti flavor) encodes Serre duality $QC^!(X)\simeq QC^!(X)^\vee$ as its unit (and an analogous duality applies to the stack of \'etale local systems). On the other hand, as explained in Section~\ref{tensor and miraculous} it is a very special feature of the stack $\Bun_G(\Sigma)$, Gaitsgory's miraculous duality Theorem~\ref{miraculous thm}, that the very sheaf $\Delta_!\omega$ encodes a self-duality of the category of automorphic sheaves (of either de Rham, Betti or \'etale flavor). We thus deduce the following:

\begin{quote}
 \label{group=miraculous}
The group period conjecture is equivalent to the assertion (formulated as~\cite[Conjecture 0.2.5]{gaitsgorystrange} in the de Rham setting) that the geometric Langlands correspondence (in de Rham, Betti or \'etale versions) intertwines the miraculous self-duality of automorphic sheaves and Chevalley-twisted Serre duality of spectral sheaves.
\end{quote}

\subsubsection{Period sheaves as kernels for functoriality} \label{Pshkernels}

We  can use functions or sheaves on a product space as kernels for ``integral transforms.''
   In particular, this can be done on either the automorphic or spectral
side of the Langlands correspondence. But  -- starting with matching automorphic and spectral
kernels --  it is not clear how the resulting integral transforms are related.

The group period conjecture, in the form stated at the end of \S \ref{diagdiag} permits
us to analyze this issue.  It asserts the compatibility of geometric Langlands with ``inner products'' (i.e., self-dualities),  and
thereby  allows us to pass more readily  between periods on product groups, and functoriality statements.

As noted in Remark~\ref{Langlands for product groups}, the geometric Langlands correspondence is expected to be tensorial under products of groups, i.e., to produce a commutative diagram of equivalences~\ref{product Langlands}.
Applying the matching self-dualities of automorphic and spectral categories, we find a commutative diagram
\begin{equation}
\label{BZfunct} \xymatrix{\bigautshv(\Bun_{G\times H})\ar[r]\ar[d]& QC^!(\Loc_{\Gv\times \Hv})\ar[d]\\
 \Hom(\bigautshv(\Bun_H), \bigautshv(\Bun_G)) \ar[r] & \Hom(QC^!(\Loc_\Hv), QC^!(\Loc_\Gv)).}
 \end{equation}
Here the horizontal arrows are induced by the geometric Langlands correspondences for $G\times H$, $G$ and $H$, while the vertical arrows are constructions of integral transforms from kernel sheaves 
(see Remark below for more discussion), provided by 
provided by the matching self-dualities of automorphic and spectral sheaves for $H$:
on the automorphic side, miraculous duality; on the spectral side,
Serre duality together with the dualizing involution.  
 Under this dictionary the group period is taken to the identity functor on Langlands parameters.

 \begin{remark}[Relation to integral transforms]\label{strange integral transforms}
 Explicitly, on the spectral side of \eqref{BZfunct} we have the usual construction of $!$-integral transforms  associated to a kernel sheaf $$\cK\in QC^!(X\times Y) \leadsto \left(\cF\mapsto p_{Y*}(p_X^!\cF^{d_H} \otimes^! \cK)\right),$$
 where we include a pre-composition with the dualizing involution on $\Loc_{\Hv}$;  while on the automorphic side we have a subtle ``miraculous'' modification (see ~\cite{gaitsgorykernels} for the $D$-module setting and~\cite{AGKRRV2} for the subtler ``enhanced" version in the \'etale setting).
 \end{remark}
 
\medskip 

In particular suppose 
$$ G\times H\actson M=T^*(X, \Psi), \Gv\times \Hv\actson \Mv=T^*(\check{X}, \check{\Psi})$$ are   dual hyperspherical varieties, for which the (normalized) period and $L$-sheaves $\mathcal P_X\in \autshv(\Bun_{G\times H}), \mathcal L_{\Xv}\in QC^!(\Loc_{\Gv\times \Hv})$ are defined. In other words, we are in the setting of  Conjecture~\ref{GlobalGeometricConjecture} for product groups. Then the conjecture predicts that $\mathcal P_X^{\spec}$ and $\mathcal L_{\Xv}^d$ match under geometric Langlands. Then we have the following corollary of Conjecture~\ref{GlobalGeometricConjecture}:
\begin{conj}\label{functoriality conjecture}
The geometric Langlands correspondence intertwines the functors given as in \eqref{BZfunct} by the integral kernel $\mathcal P_X^{\spec}$ and by $\mathcal L_{\Xv}^{d}$, i.e., 
we have a commutative diagram $$\xymatrix{ \bigautshv(\Bun_{H})\ar[r]\ar[d]^-{\mathcal P_X^{spec}}& QC^!(\Loc_{\Hv})\ar[d]^-{\mathcal L_{\Xv}^{d}}\\
 \bigautshv(\Bun_G) \ar[r] & QC^!(\Loc_\Gv)},$$
where $d$ denotes the dualizing involution.\end{conj}

 \begin{example}
  \label{Eisgeom} 
  An interesting example is the  ``Eisenstein case'',
by which we mean a  putative duality between
  $$X = U \backslash G, \hskip.3in \check{X}= \check{U} \backslash \check{G}$$
as $G \times T$ and $\check{G} \times \check T$ spaces, both with trivial $\GGm$-action.
  In this case $M,\Mv$ are {\em not affine} hence do not fit our definition of hyperspherical variety. In many ways this example fits well with the formalism of this paper nonetheless,
  but it does present peculiarities, which we will draw attention to. 
 The actions are given by 
\begin{equation} \label{numact} (g,t): U x \mapsto U t^{-1} x g, \ \ (\check{g}, \check{t}): \check{U} x \rightarrow \check{U} \check{t} x \check{g}\end{equation}
  (note the inverse; to see why  something of this nature, take $G$ a torus, where this  should reduce to the group case).
  
Conjecture~\ref{GlobalGeometricConjecture} in this case predicts that the normalized Eisenstein period sheaf $\mathcal P_{\Eis}$ and Eisenstein $L$-sheaf $\mathcal L_{\Eis}$ match
(after applying the duality involution $d$ to the latter). 
However, numerical computation  (\S \ref{Eisappendix})
suggests that the shifts (Tate and cohomological) of the global Conjecture
are not correct, and rather there should be an additional shift of $(g-1)\dim(U)$.
This discrepancy will also manifest itself in the study of parity  (\S \ref{Parity Examples}).

In any case, let us ignore this for the moment, 
 and describe the relation of this conjecture with the theory of geometric Eisenstein series. 
  After switching $\check{X}$ to a left action by our general conventions \S \ref{leftrightconventions},
  we can write the spaces as:
\begin{equation} \label{bminusid} X \simeq B \backslash (G \times T), \check{X}  \simeq (\check{G} \times \check{T})/\check{B},\end{equation}
   and the map $B \rightarrow T$ is the standard one, but the map
  $ \check{B} \rightarrow \check{T}$ is the inverse of the standard one. 
  The duality twist of $\check{X}$ is identified with
  $(\check{G} \times \check{T})/\check{B}^-$,
  again with inverted map $\check{B}^- \rightarrow \check{T}$.
 Now consider the diagrams:

$$
\xymatrix{ 
&\Bun_{B} \ar[dl]^-q\ar[dr]^-p& \\
\Bun_T && \Bun_G}.  
\hskip.3in 
\xymatrix{ 
&\Loc_{\Bv^-}\ar[dl]^-{\qv}\ar[dr]^-{\pv}& \\
\Loc_{\check{T}} && \Loc_{\Gv}
}$$

We have a spectral  Eisenstein series functor $\Eis_{\spec}=\pv_*\qv^!$, which is the functor defined by the $L$-sheaf $\mathcal L_{}^{d}$ using Serre duality,
taking into account the remarks before the diagram
and \eqref{bminusid}.   (In the above diagram, the map $\check{B}^- \rightarrow \check{T}$
is just the standard map: in passing from kernels to functors,
we have implicitly used the dualizing involution on the source, and this
has the effect of removing the previous sign.)

We have two versions of automorphic Eisenstein series
$\Eis_*=p_*q^!$ and $\Eis_!=p_!q^*$, which are (formally speaking) the $!$- and $*$-integral transforms defined by the period sheaf $\mathcal P_{\Eis}$;\footnote{It is a nontrivial result of~\cite{DrinfeldGeometricconstantterm2016} that the latter functor is well defined in the de Rham setting, where a priori only 
$!$-pullbacks and $*$-pushforwards exist in general.}
up to sign issues mentioned below, the Eisenstein desideratum formulated by Arinkin and Gaitsgory \cite[\S 13]{ArinkinGaitsgory} is that the geometric Langlands correspondence intertwines the functors $\Eis_!$ and $\Eis_{\spec}$.
The identification of $\Eis_!$ with the integral transform associated to $\mathcal P_{\Eis}$ under miraculous duality -- hence the
compatibility between the Eisenstein desideratum and the Eisenstein functoriality provided by Conjecture~\ref{functoriality conjecture}, is precisely the subject of~\cite[Theorem 4.1.2]{gaitsgorystrange} (specifically the first half of the proof).

As we see there are some issues that remain to be studied here. Most significant is the unexplained shift remarked above
and computed in \S \ref{Eisappendix}. 
There is also a matter of signs, i.e., on the right diagram $\Bv^-$ appears
rather than $\Bv$; while this does not naively  align with \cite{ArinkinGaitsgory}, it is probably a minor issue
related to the various possibilities for different sign normalizations 
implicit in our discussion.
  \end{example}
\subsubsection{Geometric Gan-Gross-Prasad Period (GGGPP) and the Geometric Theta Correspondence} 
An important example of using period sheaves as functorial kernels
 is given by the $\theta$-correspondence. This can be used not only transfer  automorphic forms, but to transfer
interesting periods. \footnote{We thank Wee Teck Gan for this suggestion.}

 We discuss briefly perhaps  the simplest example 
 of this; it  relates, in classical language, to the relationship between
 Fourier coefficients of $\theta$-series and representation numbers of quadratic forms.

 Consider  the hyperspherical dual pair
\begin{equation} \label{Besseltheta} (\SL_2 \times \mathrm{SO}_{2n}, \mathrm{Std} \otimes \mathrm{Std}) \mbox{ and } (\SO_3 \times \SO_{2n},  \mathrm{Bessel}),\end{equation}
 The Bessel period on the left is, more precisely, defined by the subgroup $\Delta (\SO_3) \cdot (V, \Psi)$ for a suitable unipotent subgroup.
 
 Using the respective period and $L$-sheaves to define
 transforms as in \eqref{BZfunct}, we obtain
\begin{equation} \label{thetafunct}  \mbox{sheaves on $\Bun_{\SL_2}$} \rightarrow \mbox{sheaves on $\Bun_{\SO_{2n}}$}.\end{equation}
and similarly on the spectral side. 
It is a familiar phenomenon in the theory of $\Theta$-correspondence that \eqref{thetafunct}
``carries the Whittaker period on $\SL_2$ to the period on $\SO_{2n}$ defined by $X=\SO_{2n}/\SO_{2n-1}$.''
The numerical
statement is that the $\Theta$-lift adjoint to \eqref{thetafunct} pulls back
the Whittaker coefficient to the $\SO_{2n-1}$-period;
we have not verified the analogous phenomenon geometrically.

On the dual side, this presumably manifests itself as follows, with reference to \eqref{Besseltheta}: the symplectic reduction of the Bessel
space by $\SO_3$ gives $T^*\check X$, with  $\check X$
the $\SO_{2n}$-space defined by $\SO_{2n}/\SO_3 (U, \Psi)$;
this is just the dual of the $\SO_{2n}$-space $X$.

It will be interesting to study this further. 

 \subsection{Spectral projections}\label{nilpotent projection is good}
In this section we return to the point that has been raised at the start of \S \ref{case1}, and explore the role of projection of automorphic sheaves to nilpotent singular support in the global conjecture; after an informal general discussion (\S \ref{igd} and \S \ref{igd2}) we explain in detail the simplest case -- the duality of the automorphic Whittaker model for $G=\Gm$  -- i.e., $G=X=\Gm$,   ``Dirichlet boundary conditions,'' and the trivial period $\Gv=\Gm$
acting on $\Xv=\pt$,  ``Neumann boundary condition''. \index{Dirichlet boundary condition}
\index{Neumann boundary condition}

We refer to Section~\ref{automorphic categories} for an overview of the main features of the different automorphic sheaf theories,
and in particular to \S \ref{spectral action} for a discussion of spectral projection. 

\subsubsection{General discussion: why do we want to project?} \label{igd}
One might want a direct comparison between period and $L$-sheaves under the geometric Langlands correspondence, and indeed we predict such a comparison in the de Rham setting. However,  we'll see below that -- even in very simple cases --  the endomorphisms of the period sheaf in finite or Betti contexts are much smaller than those of the $L$-sheaf, so no naive comparison is possible. More generally, in the constructible world (in any characteristic) there is no known version of geometric Langlands that takes as input the entire category of (ind-)constructible sheaves on $\Bun_G$, so we must project the period sheaf into the ``spectrally decomposable'' category in order to apply the correspondence. 

Equivalently, we can only ``test'' -- take homomorphisms to/from --  the period sheaf against a suitable subcategory of all automorphic sheaves including in particular Hecke eigensheaves (which correspond to skyscrapers on the spectral side). 
\begin{itemize}
\item[-] In the finite characteristic situation, and more generally in the setting of the \'etale geometric Langlands correspondence, the available automorphic test objects are ind-constructible sheaves with nilpotent singular support, which correspond spectrally to sheaves with finite (or equivalently proper) support.
\item[-] In the Betti situation we have access to arbitrary automorphic $\C$-sheaves with nilpotent singular support,  which correspond to arbitrary support on the Betti stack of local systems (i.e., the Hecke eigenvalues can vary in algebraic families). 
\end{itemize}

The subtleties between the behavior of period sheaves in the different formulations concern behavior ``at infinity'' in $\Loc$, for example since the de Rham and Betti functions on $\Loc$ \footnote{That is: the ring of regular functions on $\Loc$, but with respect to the differing algebraic structures
corresponding to the Betti and de Rham moduli spaces.} differ by their growth at infinity. Dually, on the automorphic side from a microlocal or symplectic perspective the [compact] nilpotent sheaves correspond to objects living over {\em finite} subsets of the base of the Hitchin system (in the semiclassical limit).
Period sheaves are typically very different -- for example, the Whittaker period sheaf corresponds to a {\em section} of the Hitchin fibration. Thus, spectral projection on period sheaves is a violent operation some of whose properties remain mysterious. 
 
\subsubsection{Left, right, and spectral projections}  \label{igd2}
We are interested in projecting automorphic sheaves
into the subcategory of sheaves with nilpotent singular support. 

In general, there are two ways in which we can attempt
to project a category into a subcategory: by taking the left or right adjoint
of the inclusion. We will refer to these as ``left projection'' and ``right projection,''
and will denote them by $\mathcal P \mapsto \bar{\mathcal P}_l$ and $\mathcal P \mapsto \bar{\mathcal P}_r$ when they exist.
  There are tautological maps $\mathcal P \rightarrow \bar{\mathcal P}_l$ and $\bar{\mathcal P}_r \rightarrow \mathcal P$. 
  In our specific case, there is another projection of interest, the {\em spectral projector} \index{spec} \index{$l$=subscript of left projection} \index{$r=$ subscript of right projection}
  $$\mathcal P \rightarrow \mathcal P^{spec},$$
  which, as we will see, corresponds to either a left or right projection according to context:

 The behavior of the two projections in our settings is as follows:
 
 \begin{itemize}
 \item In the Betti setting, both left and right nilpotent projection exist. 
   This left adjoint exists in the Betti setting, and in fact agrees with the spectral projector, resulting in the Betti period sheaf $$\mathcal P\mapsto \mathcal P^{spec}= \bar{\mathcal P}_l, \ \ \mbox{Betti setting}. $$
   
 \item 
 In the \'etale setting,  the Beilinson spectral projector $\mathcal P\mapsto \bar{\mathcal P}=\mathcal P^{spec}$ produces a nilpotent sheaf which conjecturally identifies with the right projection $\mathcal P \rightarrow \bar{\mathcal P}_r$:
 $$ \mathcal P \mapsto \mathcal P^{spec} \stackrel{?}{=} \bar{\mathcal P}_r, \ \ \ \mbox{\'etale setting}.$$

On the other hand,   the left adjoint does not exist except as a pro-functor\footnote{Note that the left adjoint to a colimit preserving functor must preserve compact objects, thus is much harder to represent in the \'etale setting where only constructible sheaves can be compact. This is perhaps easier to see on the spectral side where $QC(\Loc_\Gv^{Betti})$ has lots of perfect complexes (like the structure sheaf) while in the \'etale setting (where $\Loc$ is close to being a formal scheme) only sheaves with finite support (on the coarse moduli space) can be compact.
} in the \'etale setting.
We obtain thus a pro-object $\bar{\mathcal P}_l$ in the \'etale setting, 
which can  also be  identified with the natural pro-version of the Beilinson spectral projector coming from the ind-structure of $\Loc^{res}$. 
Thus, both (pro-)left and  right projections in the \'etale setting are given in terms of the spectral projection, so are ``not too far apart'' in a precise sense.

\end{itemize}

  Matching the geometric conjecture with the numerical conjecture, as discussed in \S \ref{mumerical}, requires
studying homomorphisms {\em from} $\mathcal P$ and $\mathcal L$ {\em to}  Hecke eigensheaves. 
For this purpose it is useful to understand the {\em left} nilpotent projection $\bar{\mathcal P}_l$, corepresenting the functor $\Hom(\mathcal P, -)$ 
 on the category of sheaves with nilpotent singular support.
On the other hand, right nilpotent projection $\bar{\mathcal P}_r$ will control, instead,
homomorphisms {\em to} $\mathcal P$ and $\mathcal L$ from Hecke eigensheaves; such homomorphisms are also classically interesting, and are captured numerically by the more exotic star periods, see \S \ref{starperiods}. 
   \index{star period}

 This plethora of projections gives rise to several possible candidates for the geometric conjecture. The candidate
 we have chosen (to use the spectral projection in each context) is based on some simple
 plausibility checks that rule out the other options. 
As we discuss below, in the Betti setting the right period sheaf $\bar{\mathcal P}_r^{B}$ fails to match the $L$-sheaf already in the first nontrivial example, while the left version passes that test;
in the finite context the situation is reversed. 
It remains of interest --
in order to reproduce standard numerical statements, as in Example \ref{mumerical} -- 
to understand
what corresponds spectrally to the pro-object $\bar{\mathcal P}_l$ in the finite context. For this one would need a pro-version of the $L$-sheaf; we do not attempt this in this paper.

All in all, then, the situation involves many complications and is still an evolving one!

\begin{remark} 
Also, from the above discussion, we see
that, although the Betti right-projected period sheaf $\bar{\mathcal P}_r^{B}$ does not 
match appropriately with the $L$-sheaf, 
this issue is rectified by projecting further to the \'etale category of sheaves. In effect,
this latter projection amounts to  working only over finite subsets of $\Loc$ (up to unipotents). Speaking informally, the problem
with the Betti  right-projected period sheaf $\bar{\mathcal P}_r^B$ is therefore that its behavior at infinity in $\Loc_{\check{G}}$
is bad.

 \end{remark}

\subsubsection{The abelian Whittaker period: de Rham and Betti} \label{WhitBetti}
We consider in some detail the case $G = \Gm$ and $X=\Gm$ the Whittaker space.

Note that $\Bun_{\Gm}$ is a product of the Picard scheme with $B\Gm$, while the derived stack of local systems, in either Betti or de Rham context,  is the product of the classical stack of local systems with the spectrum of an derived exterior algebra, and the geometric Langlands correspondence respects this decomposition (reducing to Koszul duality on the ``nonclassical'' factor).

For convenience, {\em we will ignore the factor of $B\Gm$ automorphically
and the derived exterior algebra spectrally.} 
As we will explain in Remark \ref{missing derived} below, both the period
sheaf and the $L$-sheaf will have the structure of an external product
along the decomposition
$$ \Bun_{\Gm} = \Pic \times B\Gm \mbox{ and }\Loc_{\Gm} = \Loc_{\Gm}^{\mathrm{classical}} \times \Spec k[x_{-1}],$$
and the second factors of these decompositions will match.
So we restrict our attention to what happens on the $\Pic$ and $\Loc_{\Gm}^{\mathrm{classical}}$. 
 
 The Whittaker period sheaf  on $\Pic$ is simply a skyscraper ($\delta$-function $D$-module)
at the identity, i.e., the pushforward of the constant sheaf under $\mathrm{pt} \rightarrow \Pic$.
Thus, the endomorphisms of the period sheaf are simply scalars in either setting. 
The singular support (or semiclassical limit) of the period sheaf is the cotangent fiber to $\Pic$ at the trivial bundle.

On the other hand the $L$-sheaf  on $\Loc_{\Gm}^{\mathrm{classical}}$
is, up to a shift, the structure sheaf. Since $\Loc_{\Gm}^{\mathrm{classical}}$  is
  the product of a smooth
   variety by $B\Gm$ these
 endomorphisms  are
$$ \mathcal{O} = \mbox{algebraic functions on $\Loc_\Gm$}.$$

These functions however differ drastically between the Betti and de Rham setting. In the Betti setting we find $ \mathcal{O} = \kk[H_1]$ the group algebra of the first homology of the curve (abelianization of $\pi_1$) with coefficients in our structure ring $\kk$; this is, equivalently, the ring of functions on the algebraic torus which is the representation variety of homomorphisms $\pi_1 \rightarrow \Gm$
(considered now without derived or stack structure). In the de Rham setting on the other hand $\Loc_\Gm$  is the universal vector extension of the Jacobian and admits no nonconstant algebraic functions, so the endomorphisms of the $L$-sheaf are just scalars. 

Thus, in the de Rham situation the endomorphisms of period and $L$-sheaves match, while in the Betti setting the former is much smaller than the latter.
 
 \begin{remark}[Derived version] \label{missing derived} 
Let us mention the (matching) behavior of period and L-sheaves over the ``derived'' factor we ignored above. 
The period sheaf in both Betti and de Rham settings is the external product of the skyscraper above with the ``regular'' sheaf on $B\Gm$
(the pushforward of $\kk$ under $\mathrm{pt} \rightarrow B\Gm$.) This regular sheaf corresponds to the free module over the homology $H_*(\Gm)\simeq 
\kk[x_{-1}]$.  

On the spectral side the derived ring of functions on $\Loc$ carries an extra tensor factor $\kk[x_{-1}]$.  This contributes a free $\kk[x_{-1}]$-module
to the dualizing sheaf,
matching the period sheaf.

\end{remark}

\subsubsection{Projection to local systems.}
Restricting now to the Betti context, we now consider the left and right projections of the skyscraper sheaf $\mathcal P$ to local systems on $\Pic$, i.e., to automorphic sheaves with nilpotent singular support.  

  {\em Left} nilpotent projection $\bar{\mathcal P}_l$ replaces the role of $\delta$ in the discussion above by the regular $\pi_1(\Bun_{\Gm})$-module $\mathcal{O}[2g]$
  where $\mathcal{O}=\bigoplus_{\pi_1} \kk$ is the group algebra
(considered as a local system on $\Bun_{\Gm}$ in the obvious way). In other words, to make $\delta$ into a local system we replace the inclusion of the trivial local system by the path fibration (in this case universal cover) of $\Bun_{\Gm}$ and thus the compactly supported pushforward of the constant sheaf becomes the ``regular'' local system.

{\em Right} nilpotent projection $\bar{\mathcal P}_r$ replaces $\delta$ by
 $M\simeq \prod_{\pi_1} \kk$, where 
 $M$ is the space of all functions $\pi_1 \rightarrow \kk$.  
 These statements simply record the fact that taking
 fiber at the identity exhibits an equivalence of locally constant sheaves on $\Pic$ with
 the  category of $\pi_1$-representations on $\kk$-vector spaces, and for such a representation $V$ we have
 $$V \simeq \Hom_{\pi_1}( \mathcal{O},V), \hskip.3in V^* \simeq \Hom_{\pi_1}(V, M).$$

   Note, computing endomorphisms as $\mathcal{O}$-module, 
\begin{equation} \label{yuk} \mathrm{End}(\mathcal{O}[2g]) = \mathcal{O},   \ \ \mathrm{End}(M) = \textrm{``next question, please!''}.\end{equation}
 Corresponding, the endomorphisms of $\bar{\mathcal P}_{l}$ exactly match functions on $\Loc_\Gm$
 and indeed   $\mathcal P^{spec}=\bar{\mathcal P}_l$ is exactly what corresponds to $\mathcal{L}$
 under the Langlands correspondence. However the endomorphisms of $\bar{\mathcal P}_r$ are far too large.

 \begin{remark} \label{rl at point}
Despite their difference in size, $\bar{\mathcal P}_r$ and $\bar{\mathcal P}_l$
in the Betti case are, in a certain sense, not far apart, as we now sketch.
 For $\rho$ a point of $\Loc$,  the restriction $\iota_{\rho}^* \mu$ of the composite morphism $\mu: \bar{\mathcal P}_{r}^B \rightarrow \bar{\mathcal P}_l^B$ arising from
$$\bar{\mathcal P}_r^B \rightarrow \mathcal P \rightarrow \bar{\mathcal P}_l^B$$
is an isomorphism. 
Let us illustrate the phenomenon in the case $\pi_1=\Z$
leaving the real-world
case  $\pi_1=\Z^{2g}$ to the reader. 
The short exact sequence of $k[x^{\pm 1}]$-modules
 $$ \kk[x^{\pm 1}] \rightarrow \kk((x)) \oplus \kk((x^{-1})) \rightarrow  \underbrace{\mbox{formal series $\sum_{n \in \Z} b_n x^n$}}_{R}$$
gives rise to a morphism of the right hand group $R$ to $k[x^{\pm 1}][1]$
in the derived category of $k[x^{\pm 1}]$-modules. 
Clearly, this morphism -- which is analogous to $\bar{\mathcal P}_r^B \rightarrow \bar{\mathcal P}_l^B$ --  is not an isomorphism, but it is an isomorphism
when pulled back to any $\kk$-point of $\Gm$,
because $\kk((x))$ and $\kk((x^{-1}))$ do not have support there. 
The  error is ``supported at $0$ and $\infty$.''
It would be interesting to study the analog of this example for some semisimple $G$.
\end{remark}

 \subsubsection{Projection to \'etale categories of sheaves}
 
  Let us carry through the same discussion ($G=\Gm$ and $X=\mbox{Whittaker}$) in the \'etale case.
 Recall here that 
   $\kk$ is either $\CC$ or the algebraic closure of $\Q_{\ell}$, and 
 the allowable sheaves on the automorphic side are ind-(locally constant constructible) sheaves.  On an irreducible $\FF$-variety $X$, such sheaves correspond to representations of $\pi_1(X)$
 on a $\kk$-vector space which are {\em locally finite}, i.e.
 each vector lies in a finite dimensional $\pi_1$-stable subspace.  We are primarily interested in this context
 when $\FF$ is the algebraic closure of a finite field, but we will also remark on what happens in the \'etale setting over $\C$.\footnote{This should be distinguished from our
 default usage of the word Betti, where we consider allow locally constant rather than only ind-(locally constant constructible) sheaves.}

 Each $\kk$-point of $\Loc_{\Gm}$, i.e., each rank one local system $\rho: \pi_1 \rightarrow \kk^{\times}$,
 admits a universal deformation to a map $\pi_1 \rightarrow R_{\rho}^{\times}$
 for a certain smooth complete $\kk$-algebra $R_{\rho}$ abstractly isomorphic to  $ \kk[[x_1, \dots, x_{2g}]]$ ; and the \'etale classical moduli space $\Loc_{\Gm}^{\mathrm{classical}}$ is a disjoint union of the spectra of 
 $R_{\rho}[x_{-1}]$, each quotiented again by the trivial action of $\Gm$. 
 (In fact, all of these $R_{\rho}$
 are isomorphic to one another by twisting.)
 In particular, we again have  
\begin{equation} \label{Lhom2} \mathrm{End}(\mathcal{L}) = \prod R_{\rho}\end{equation}
 or, as in Remark \ref{missing derived}, if we are to include
 the contribution of the derived structure, we would additionally add
 a generator in degree $-1$ in each factor. 
 
 In this case the left nilpotent projection of $\mathcal P$ is (up to a shift) the pro-object 
 corresponding to the complete $\kk$-algebra $\bigoplus R_{\rho}$, i.e., the structure pro-sheaf of the formal stack $\Loc_\Gv$. Of course $\Loc_\Gv$ (like any stack) carries a structure sheaf as an object in $QC$, 
which is given by calculating the limit of the pro-sheaf version of the structure sheaf ($QC$, being presentable -- see Appendix \ref{HigherCatAppendix} -- is closed under all small limits by ~\cite[Corollary 5.5.2.4]{HTT}), however this object doesn't corepresent the given left adjoint (as can be seen for example by its lack of compactness). Thus, it is clear in this case how to modify the $L$-sheaf to get a pro-version that matches $\bar{\mathcal P}_l$. Dually, this is the pro-local system associated to the complete (not locally finite) $\pi_1$-representation $\bigoplus R_{\rho}$. 
 
 On the other hand, the 
   right nilpotent projection of $\mathcal P$
corresponds, on the spectral side,
to the object whose $\rho$th component is given by 
 $E_{\rho}=$ the top local cohomology of the local ring $R_{\rho}$. 
 In coordinates $R_{\rho} =  \kk[[x_1, \dots, x_{2g}]]$   
 the module $E_{\rho}$ is given by  ``$\kk[x_i^{-1}]$''
 where, to define the module structure,  any expression involving $x_i$ to a positive exponent is regarded as zero.
Observe that for $R_{\rho}$-module $M$  of  finite length we have
 $\Hom_{R_{\rho}}(M, E_{\rho}) \simeq \Hom_k(M, k)$
 via the ``constant coefficient'' map $E_{\rho} \rightarrow k$
 (by exactness of the left-hand  functor we reduce to the case of $M$ the augmentation). 
 
 In particular $\mathrm{End}(E_{\rho}) = R_{\rho}$ because $R_{\rho}$ is complete 
 (see for example \cite[Theorem 3.10]{Huneke})
and therefore the endomorphisms of $\bar{\mathcal P}_r$
are precisely  $\prod R_{\rho}$ (or the same adjoining $x_{-1}$
if we compute on $\Bun_{\Gm}$ rather than $\Pic$) which, in contrast to  \eqref{yuk}, match the endomorphisms of the $L$-sheaf.

\begin{remark}
Over $\FF=\C$,  we can think of 
the right projection $\bar{\mathcal P}_r$
to the  category of ind-constructible sheaves 
 as ``fixing'' the behavior of the right Betti period sheaf $\bar{\mathcal P}^B_r$, 
 by (speaking in dual terms)   
 restricting to finite subsets of $\kk$-points $\rho$ of $\Loc_{\check{G}}$.
This makes sense: we  already saw in Remark \ref{rl at point} that restricted to points
 $\bar{\mathcal P}^B_r$ and $\bar{\mathcal P}^B_l$ are actually the same thing.
 \end{remark}

\subsection{Parity and independence of spin structures} \label{parity returns}
We now observe that  parity conditions on our hyperspherical spaces implies that the validity of the global geometric conjecture, Conjecture~\ref{GlobalGeometricConjecture}, is independent of the choice of spin structure, or equivalently allows us to formulate it independently of this choice. 
In the setting of Conjecture \ref{GlobalGeometricConjecture}, we denote by \begin{equation} \label{Autspecnotn}  \eta: G \rightarrow \Gm,\;  
\check\eta: \check{G} \rightarrow \Gm,\; 
 \eta : \Gm \rightarrow \check{G},  \;
  \check\eta: \Gm \rightarrow G\end{equation}
 the two eigencharacters, cf.  \eqref{gammaXdef}, considered also as central cocharacters into the dual groups.   We get corresponding maps
  $\Bun_G \rightarrow \Bun_{\Gm}$ and $\Loc_{\check{G}} \rightarrow \Loc_{\Gm}$.

In Remark~\ref{period sheaf and spin} we noted that the dependence of the period sheaf on spin structures 
can be 
summarized 
by defining the period sheaf
inside  $$\cP_X\in \Hom_{\Bun_{\Z/2}}(\mathrm{Spin}_\Sigma, \autshv(\Bun_G))$$ where
$\Bun_{\Z/2}$ acts on $\Bun_G$ through multiplication via the central homomorphism
$$z=\check\eta e^{2 \rho}: \Z/2 \rightarrow Z(G).$$ The same is true of the normalized period sheaf, which involves no additional spin twist. 
That is to say: modifying the choice of $\mathcal K^{1/2}$
by a $2$-torsion line bundle $\mathcal{L}$ on
$\Sigma$ modifies the period sheaf by ``translation through $\mathcal{L}$'' via $z$. 

Likewise in \S \ref{Lspinindep} we saw that the dependence of the normalized $L$-sheaf on spin structures is captured by
considering 
the normalized $L$-sheaf
as an element of
 $$
 \mathcal L_{\Xv}^{\norm}  \in  \Hom_{\Bun_{\Z/2}}(\Spin_\Sigma^\vee, QC^!(\Loc_\Gv)).$$
Here $\Bun_{\Z/2}$ acts on $QC^{!}(\Loc_\Gv)$ by tensor product with 2-torsion line bundles which can be described as pulled back from $\Loc_\Gm$ via $\check\eta: \check{G} \rightarrow \Gm$.
That is to say, modifying the choice of $K^{1/2}$ by a $2$-torsion line bundle $\mathcal{L}$
has the result of tensoring the $L$-sheaf with the pullback of $\mathcal{L}$ via $\check\eta$.

We now need to appeal to the compatiblity between the Langlands correspondence and abelian duality for the center of $G$ expressed in Remark~\ref{fluxes}.
Namely, the translation action of torsors for the center $\Bun_{Z(G)}$ on $\autshv(\Bun_G)$ matches under the geometric Langlands correspondence
with the tensor product action on $QC^!(\Loc_\Gv)$ given by the canonical dual homomorphism
$\Bun_{Z(G)}\to  \mathrm{Pic}(\Loc_\Gv)$
(see discussion before
\eqref{abelian double cover eqn}). 
 In particular for a central involution $z:\Z/2\to Z(G)$   the induced actions of $\Bun_{\Z/2}$ agree.

Now observe that  the twists appearing in the period and $L$-sheaf differ by a universal amount, i.e.
{\em independent of the period under consideration} --  namely, the shift by the canonical parity element $e^{2\rho}:\Z/2\to Z(G)$.
Therefore, the validity of the conjecture does not depend on the choice of $\mathcal K^{1/2}$: If we change the choice of $\mathcal K^{1/2}$,
it can be compensated by changing the normalization of the Langlands correspondence, twisting
by a line bundle on the spectral side. See also Section~\ref{extended-group appendix}.

\subsection{Parity and change of grading} \label{mitch} 
In this section we study the interaction of change of grading (as in Sections~\ref{changing grading automorphic} and ~\ref{changing grading spectral}) with Conjecture~\ref{GlobalGeometricConjecture}. This is, more or less, just book-keeping. 

Specifically, suppose that $(G, M)$ and $(\check{G}, \check{M})$ are a hyperspherical pair in the sense of \S \ref{goodhypersphericalpairs},
with both $M, \check{M}$ polarized. What we will explain here is how to twist the conjecture to avoid the normalizing factors on the $L$- and period sheaf. This will introduce an extra and rather unenlightening Tate twist
as well as altering the $\GGm$-action on $M$ and $\check{M}$. 
 The main feature of the resulting ``unnormalized''
conjecture is that it has a more transparent parity condition
in the sense of \S \ref{analyticarithmetic},  see 
\eqref{mtpcap}.
We use notation as in \eqref{Autspecnotn}; we will
use the same letters for the induced maps   $\Bun_G \rightarrow \Bun_{\Gm}$ and $\Loc_{\check{G}} \rightarrow \Loc_{\Gm}$

  \begin{definition} \label{unnormalized action def}
If $(G, M=T^*X)$  and $(\check{G}, \check{M} = T^* \check{X})$
are a dual hyperspherical pair, both sides equipped with a polarization and eigenmeasure, 
we call the $\GGm$-actions on $X$ and $\check{X}$
obtained by twisting the neutral actions by $\check{\eta}^{-1} \mbox{ and } \eta^{-1}$
respectively the {\em unnormalized $\GGm$-actions} on $X, \check{X}$.

We denote the twisted spaces by $$X^{\tw} := X[\check{\eta}^{-1}] \mbox{ and }\check{X}^{\tw} := \check{X}[\eta^{-1}],$$ i.e.
   $X[\dots]$ means that we twist the $\GGm$-action on $X$ by 
the stated character.  
\end{definition}

\begin{example}
It should be noted that the {\em unnormalized actions depend on the choice of polarizations.}

Here is an example. 
 Take $(G=\mathbb{G}_m, M=T^*\mathbb{A}^1)$, with the scaling $\GGm$-action.
Then the dual space is the same $(\check{G} = \mathbb{G}_m, \check{M} =T^* \mathbb{A}^1)$. 
 With the standard polarization we have $$\eta = \check\eta= \mbox{standard character of $\Gm$}$$
  and correspondingly the twisted spaces $X^{\tw}, \check{X}^{\tw}$ are given by $\mathbb{A}^1$
 with {\em trivial} action of $\GGm$. 
 
The twisting process, however, {\em depends on a choice of polarization}. Had we chosen, for example,
the ``other'' polarization of $T^*\mathbb{A}^1$,  thus inverting the $G$-action, then $\eta$ would be inverted, and correspondingly
the normalized action on the dual $\check{X} =\mathbb{A}^1$ would be given by $\lambda \mapsto \lambda^2$. 
In words this corresponds to various choices of ``functions'' or ``forms'' on either side.

\end{example}

\begin{remark} ``Un-normalized'' refers loosely to the fact that these actions are adapted
to consideration of sheaves and functions without incorporating $L^2$-twists,
that is to say, it is ``arithmetic'' in the general parlance of \S \ref{analyticarithmetic}. 

The definition is somewhat strange: the twist involved in defining
the normalized action on $X$ involves the eigencharacter on volume forms on $\check{X}$
and vice versa. 
Now    the eigenform for $X$ is not unique, but determined only up to characters of $X$,  
  i.e., $G$-eigenfunctions in $\mathcal{O}(X)$; 
  this dependence does not matter by
  \S \ref{ssseigencharinocuous}.

Regrettably, as far as global coherence of notation,
although these ``unnormalized'' actions fall in the ``arithmetic''
side of the arithmetic/analytic divide of  \S \ref{analyticarithmetic},
they do not quite coincide with the ``arithmetic actions'' defined in \S \ref{Marithshear}. 
\end{remark}

 \begin{proposition} \label{shiftpredlem}
 The following are equivalent:
 \begin{itemize}
 \item[(a)] 
$\mathcal P_X^{\norm}$ and the dualizing twist of $\mathcal L_{\check{X}}^{\norm}$ are dual
  to one another;
  \item[(b)]  The corresponding statement holds for $X, \check{X}$ with {\em unnormalized}
  $\GGm$-action and {\em unnormalized} period and $L$-sheaves, i.e.:
  
        $\mathcal P_{X^{\tw}} \langle  r \rangle$ and  the dualizing twist of $\mathcal L_{\check{X}^{\tw}}$
are dual to one another.
\end{itemize}
Here, with notation as in \eqref{betaXdef}, 
 \begin{equation} \label{shiftpred} r= \beta_{X^{\tw}} + \beta_{\check{X}^{\tw}} - (g-1) \tau.\end{equation}
and    $ \tau = \langle \eta, \check\eta \rangle \in \Z$.  
\end{proposition}

 \proof
 We start with the statement $\mathcal L_{\check{X}}^{d,\norm} \iff  \mathcal P_X^{\norm}$ of the normalized conjecture
 (the symbol $\iff$ here means ``matches under the Langlands correspondence''), 
 and use  \eqref{PXGmshift} and \eqref{LXGmshift}. 
 
 If we abridge $\check{X}', X'$ for the spaces with twisted $\GGm$-actions we have
 $T \mathcal{\mathcal P}_{X}^{\norm} =  \mathcal{\mathcal P}_{X'}^{\norm}$
 where $T$ is a translation by $\check\eta(\mathcal{K}^{1/2})$;   and similarly $\mathcal{L}_{X'}^{\norm} = \check{T} \mathcal{L}_{X}^{\norm}$
 with $\check{T}$ now convolution with $\delta_{\shear}^{-1}$
 through $\eta^{-1}$ (see \eqref{checkTdef3}). 
 
 So  the normalized conjecture is equivalent to
 \begin{equation}
 ( (\eta \delta_{\shear}^{-1}) \boxtimes \mathcal L^{\norm}_{\check{X}'})^{d} \iff  (T\mathcal P^{\norm}_{X'})
 \end{equation}  
 or, taking account of the effect of $d$ on $\delta_{\shear}$, 
  \begin{equation}
 (\eta \delta_{\shear}) \boxtimes (\mathcal L^{\norm}_{\check{X}'})^{d} \iff  (T\mathcal P^{\norm}_{X'})
 \end{equation}  
  Now spectral convolution with $\eta\delta_\shear$
 corresponds on the automorphic side by tensoring with
 $\langle \deg \rangle$, see after \eqref{deltasheafconvo},  so the above is equivalent to
 \begin{equation} \label{snb}
 (\mathcal L^{\norm}_{\check{X}'})^{d} \iff (T \mathcal P^{\norm}_{X'}) \langle -\deg \rangle 
\stackrel{\eqref{Tshift}}{=} T(\mathcal P^{\norm}_{X'} \langle -\deg \rangle) \langle  (g-1)\tau \rangle.
 \end{equation}
 where we applied \eqref{Tshift} with $\lambda= \check\eta$, noting
 the presence of another minus sign from the fact we have $-\deg$.

 Now, $T$   corresponds on the spectral side to tensoring by  the pullback of $[K^{1/2}]$ via $\check\eta: \Loc_{\check{G}} \rightarrow \Loc_{\Gm}$, 
 which in our notation is the line bundle $\varepsilon_{1/2}$ on $\Loc_{\check{G}}$,   see \eqref{epsilondescription} and \eqref{etaspecdef},  and
the inverse of $T$ corresponds to with tensoring with $\varepsilon_{1/2}^d$,  so \eqref{snb} becomes
  \begin{equation}
  (\varepsilon_{1/2} \mathcal L^{\norm}_{\check{X}'})^{d} \iff  \mathcal P^{\norm}_{X'} \langle -\deg  + (g-1)\tau \rangle.
 \end{equation}
  
Taking into account (\eqref{PXnormdef}) that  $\mathcal P^{\norm}_{X'} = \mathcal P_{X'} \langle \deg + \beta_{X'} \rangle$ and using \eqref{LXnormdef}, 
 $$
 \mathcal{L}_{\check{X'}}^d  \langle -\beta_{\check{X}'} \rangle \iff    \mathcal P_{X'} \langle \beta_{X'} + (g-1) \tau \rangle.$$
thus 
   \begin{equation}
 (\mathcal{L}_{\check{X}'})^{d} \iff  \mathcal P_{X'} \langle \beta_{X'} + \beta_{\check{X}'} + (g-1) \tau \rangle
 \end{equation}
 See also \eqref{betaXdef}.  By \eqref{betachange}
 $\beta_{X'} =\beta_X - (g-1) \tau$ and  $\beta_{\check{X}'} = \beta_{\check{X}} - (g-1) \tau$
 and thus 
  $$(\mathcal{L}_{\check{X}'})^{d} \iff  \mathcal P_{X'} \langle \beta_{X} + \beta_{\check{X}} - (g-1) \tau \rangle$$
 
 \qed

  \begin{remark} \label{noeigenmeasure2}   
By comparing to a suitable
cover as in \S~\ref{ssseigencharinocuous},  one can directly define the
unnormalized $\GGm$-action for the dual of a spherical variety $(X, \Psi)$
without reference to Definition \ref{unnormalized action def} -- in particular, whether
or not $(X, \Psi)$ has an eigenmeasure.  Then part (b) above gives
rise to the appropriate  ``arithmetically normalized'' formulation of the global conjecture, valid whether or not $(X, \Psi)$ has an eigenform. 
    
For instance, returning to the example  of $(X, \check{X})$ 
mentioned in Example \ref{noeigenmeasure},
statement (a) of   Proposition \ref{shiftpredlem} is undefined; but statement (b) 
makes sense and can be taken as the statement of the global conjecture,
where one takes the unnormalized $\GGm$-action on both sides to be trivial.

  This is related to the fact, already commented,
  that the parity condition is ``more transparent'' when one works with unnormalized $\GGm$-action. 
 Namely, 
Proposition
\ref{zXparity}, if applicable, proves that $(e^{2\rho}(-1),-1)$ acts trivially on $X^{\tw}$ and similarly for $\check{X}^{\tw}$.
That is to say:
\begin{equation}  \label{mtpcap}
\mbox{the action of $\check{G} \times \Gm$ on $X^{\tw}$ and $\check{X}^{\tw}$
factor through  ${}^CG$} 
\end{equation}
where ${}^CG$ is as in 
\S \ref{subsection-extended-group}.  
 \end{remark}
   
 \begin{example}
 
Take the ``Godement-Jacquet'' example, where  $(G, X) = (\GL_n \times \GL_n,  M_{n,n})$,
with action by $m \cdot (g_1, g_2) = g_1^{-1} m g_2$.
On the dual side, $\check{X} = \GL_n \times \mathbb{A}^n$
with action $(g_1, g_2) \cdot (h, v) = (g_1 h g_2^T, g_2 v)$.
\footnote{up to possibly permuting
some indices, which we didn't check, but makes no difference to the point of the example.}
 Here the eigenform character on $X$ is given by $\det(g_2 g_1^{-1})^n$, and correspondingly
 $\eta: z \mapsto (z^{-n}, z^{n}) \in \check{G}$.
   In particular, the normalized action of $ z \in \GGm$ on $\check{X}$ 
   is given by $z \cdot (h, v) = ( h, z^{1+n}v )$. 
 We see in this example that the normalized $\GGm$-action has no evident
  positivity property along the vectorial fibers of $\check{X}$.

  \end{example} 
 
 \begin{example}  Here are a few more examples comparing normalized and unnormalized $\GGm$-actions;
 the notation (u) refers to an unnormalized example. 
 
   {\small 
  \begin{table}[!ht]
\centering 
\label{tab:output2}

 \begin{tabular}{|c|c|c|c|c|}
 \hline
 name		& $(G,X)$										& $\GGm$								&			$(G,X)$									&	$\GGm$		\\
\hline
Iwasawa-Tate   	& $(\Gm, \mathbb{A}^1)$							& scaling									&  	 		$(\Gm, \mathbb{A}^1)$ 						&	scaling			\\ 
Iwasawa-Tate (u) 	& $(\Gm, \mathbb{A}^1)$							& trivial   									& 			$(\Gm, \mathbb{A}^1)$ 						&  	 trivial 				\\ 
\hline
 R.-S.   	& $(\GL_n \times \GL_{n+1},  \mathrm{Std}_{n(n+1)})$ 	& scaling									& 			$(\GL_n \times \GL_{n+1}, \GL_{n+1})$			&   trivial  \\
  R.-S. (u) 	& $(\GL_n \times \GL_{n+1},  \mathrm{Std}_{n(n+1)})$ 	& scaling									& 			$(\GL_n \times \GL_{n+1}, \GL_{n+1})$			&  $ (x^{n+1}, x^n)$ \\
  \hline 
Whit &			 $(G, G/(U, \psi))$							& $(1, e^{-2\rho})$										& $(G, G)$ & trivial \\
\hline
$\mathbb{A}^2$ norm.& $(\SL_2 \times \Gm, \mathbb{A}^2)$					& scaling 									& $(\PGL_2 \times \Gm, \PGL_2)$									& trivial 	\\
$\mathbb{A}^2$ (u) & $(\SL_2 \times \Gm, \mathbb{A}^2)$				& scaling 									& $(\PGL_2 \times \Gm, \PGL_2)$									& $e^{2\rho}$	\\
\hline
 \hline
 \end{tabular}
 \caption{Some examples comparing the normalized and unnormalized $\GGm$-actions.}

 \end{table}
}
\end{example}

 \subsection{Parity phenomena} \label{Parity Examples}
 
Proposition \ref{shiftpredlem} gives a useful mod $2$ check on our various shifts.

Namely, when doing computations involving  extracting Frobenius traces, it is best to 
(if we are using the analytic normalization of Satake)
 formulate the global Langlands equivalence as an equivalence of supersheaves, 
see Remark \ref{GLCsuper} and the discussion of \S \ref{supersloppy}.
Then   both period and $L$-sheaves should really be understood as {\em super}sheaves,
these parity twists being incurred in the normalization process on both sides
and reflecting half-integral powers of $q$ in the 
numerical theory, cf. \S \ref{Sqrtqsuper}). 

Both normalized period sheaf and normalized $L$-sheaf have in general nontrivial parity.
Let us examine the relationship between these parities. 

By construction,  
 the unnormalized period sheaf has even parity.
 On the other hand, Proposition \ref{zXparity} implies that 
the action of $-1 \in \GGm$ on $\check{X}$
coincides with $e^{2 \rho}(-1) \in \check{G}$. 
Therefore the unnormalized $L$-sheaf
has parity determined by the action of $e^{2\rho}(-1) \in Z(\hat{G})$,
this arising from the shearing operation of \eqref{Lsheafdef}. 
 The only reasonable way for
 this to be valid for all hyperspherical dual pairs simultaneously
 is the following:
 
  \begin{quote}
{\em Prediction:}   The parity of the integer $r$ of Proposition \ref{shiftpredlem}
  is independent of the choice of $(X, \check X)$. 
  \end{quote}
  which, in turn, implies  the following statement
\begin{equation} \label{parity prediction}
\dim(X) + \gamma_X  + \dim(\check X) + \gamma_{\check X}  +\dim(G/U) \stackrel{?}{\equiv} \langle \eta, \check\eta \rangle \mbox{ mod $2$.}
\end{equation}
which results by comparing a   general $(X, \check X)$ to the Whittaker case $(X=\mathrm{pt}, \check X=\Gv/\check U)$. 
 
We do not have a general proof of this but it certainly holds in all examples that we computed. See the table below. 
An interesting {\em nonexample} is the Eisenstein case -- formally outside the validity of our conjectures --
$X=G/U, \check X=\Gv/\check U$, considered as $G \times T$ and $\Gv \times \check T$-spaces.
This is another indication that
there remain some interesting issues to resolve in this case, cf.  Example \ref{Eisgeom}. 
 
 {\small 
  \begin{table}[!h] \label{tab:output3}
\centering 
 \begin{tabular}{|c|c|c|c|c|c|c|}
 \hline
 name	&  $\dim(X)$	&  	$\gamma_X$	& $\dim(\check X)$	& $\gamma_{\check X}$ & $\langle \eta, \check\eta \rangle$ & $\dim(G/U)$ \\
 \hline
  Iwasawa-Tate	 (l.2)&  $1$		&  	 $1$			& $1$			& $1$ 				& $1$ 	& $1$	\\
  Group  (l.3)	& $g$ 		& 	$0$			& $g$			& $0$ 				& $0$  	& $0$ \\
  Whittaker	(l.8)& $g-u$		&  	$0$			& $0$			&$0$					& $0$  	& $g-u$\\
Godement-Jacquet	(l.6)& $n$		&	$n$			&	$0$			&	$n$				&	$n$	&   $0$ \\
 $\GL_{2n+1}/\GL_n \times \GL_{n+1}$ (l.4)
 		& $0$ 		& $0$ 			& $n+1$ 			& $0$				& $0$ 	& $n+1$ \\
Jacquet-Shalika (l.9)	& $0$ 		& $0$ 				& $n$			& $0$			& $0$	& $n$ \\ 
Hecke.  (l.1)
	 	& $0$		& $0$			& $0$			& $0$				&	$0$	& $0$ \\
$G\times T, G/U$ 	& $g-u$ 		& $0$		& $g-u$			& $0$					&  $0$.	
																				& $g+t-u$ \\
 \hline
 \end{tabular}
 \caption{ Some examples examining \eqref{parity prediction};  dimensions taken mod $2$. $g=\dim(G)$ etc. $\gamma$ as in \eqref{gammadef}.  
 Line references, e.g.\ l.3, are to Table \ref{tab:examples}. }
 \end{table}
}

 \subsection{The  $L^2$ conjecture and the algebra of $L$-observables} \label{case5} 
 
 Our conjecture thus far has required access to a polarization, i.e., $\check{M} = T^*\Xv$ or a twisted version thereof, so as to construct an $L$-sheaf to match with the period sheaf. 
 
 However, in \S \ref{specWeil-start} --
we hinted that there is a natural algebra acting by endomorphisms on the $L$-sheaf.
This algebra is a deformation of the ``doubled'' $L$-sheaf obtained
by substituting $\check{M}$ for $\check{X}$. 

  This has a manifestation
at the level of the period sheaf.
Instead of attempting to describe $\mathcal P_X^{\norm}$ or its spectral projection,
  we instead describe its ``square,'' namely, the endomorphisms $\mathrm{End}(\mathcal P_X^{\norm})$;
  we anticipate this object can often be described solely in terms of $\check{M}$ (even when $\Mv$ is not polarizable). 
  This corresponds to the fact, familiar in the theory of periods,
 that access to $\check{M}$ alone still permits one to describe
  the square of the period. In the physics language, it is also a manifestation of the passage from geometric quantization (which requires a polarization)
 to deformation quantization (which doesn't), or from states to observables. 
 
  It is therefore reasonable to ask:
 
 \begin{quote}
 Can we give a spectral description of the endomorphisms $\mathrm{End}(\mathcal P^{spec})$
 of the period sheaf in terms of $\check{M}$?
 \end{quote}

 To discuss this in more detail, let us 
   assume (as in \S \ref{specWeil-start}) that {\em $\check{M}$ is a symplectic vector space}.
We will remark after on the formulation in the general case, and the discussion will be revisited from a more general and structured point of view in \S 
\ref{spectral geometric quantization}.

 Let us moreover restrict ourselves to the locus
 $\Loc_G^{\circ}$ of representations that have a unique classical fixed point
 on the symplectic vector space $\check{M}$. 
 This should be understand as the complement of the locus where the relevant $L$-function has a pole.
   This means that the  complex $\mathsf{V}$ encountered in \eqref{VVdef} is in fact cohomologically concentrated in a single degree and renders our discussion extremely concrete.
 
Over $\Loc_G^{\circ}$ we can form a vector bundle $\mathsf{H}$ whose fiber at a point $\rho$ is given by the cohomology $H^1(\check{M})$
 where $\check{M}$ is considered as a symplectic flat bundle over the curve by means of $\rho$.
 Over the locus $\Loc_G^{\circ}$, then, $\Loc^{\check{M}}$ is  
 the derived scheme obtained by $(-1)$-shifting of the total space of $\mathsf{H}$. 

Now, since the vector bundle $\mathsf{H}$ just mentioned carries an orthogonal structure $\mathsf{H} \simeq \mathsf{H}^{\vee}$ coming from Poincar{\'e} duality, we may form its Clifford algebra:
\begin{equation} \label{Adef} \bO_{\Mv,\Sigma}^\circ := \mbox{Clifford algebra of $\mathsf{H}$ with its natural quadratic form},\end{equation}
 which is now a sheaf of algebras on $\Loc_G^{\circ}$. \footnote{This is closely related to the 
 universal enveloping algebra of the dg Lie algebra mentioned in \S \ref{specWeil-start}.}
  It is very natural to suppose that this $\bO_{\Mv,\Sigma}^\circ$ is the ``spectrally decomposed algebra of endomorphisms
 of the period sheaf,'' a notion we now explain.  
\index{$\bO_{\Mv,\Sigma}^\circ$ $L$-observable Clifford algebra}

 The spectral action of quasi-coherent sheaves on $\Loc_\Gv$ on automorphic sheaves (the ``automorphic to Galois'' direction, see Section~\ref{spectral action} of Appendix~\ref{geometric Langlands}) implies that 
  the $\Hom$-space $\Hom(F, G)$ of spectrally decomposable automorphic sheaves $F,G\in \autshvspec(\Bun_G)$
  can be disintegrated (or spectrally decomposed) over $\Loc_{\check{G}}$, or, more formally,
  enriches to a quasi-coherent sheaf $\mathsf{Hom}(F,G)$ on  $\Loc_{\check{G}}$.  
Namely, this sheaf is defined as inner Hom in $QC(\Loc_\Gv)$, i.e., by the universal property that
for  $Q\in QC(\Loc_\Gv)$  we should have
$$ \Hom(Q, \mathsf{Hom}(F, G)) = \Hom(Q \star F, G).$$
(Of course given the full geometric Langlands conjecture we may simply take $Q$ to be the internal sheaf Hom between the Langlands transforms of $F$ and $G$ -- note that the
 $\Hom$ of ind-coherent sheaves can be enriched to take values in quasi-coherent sheaves.)
In particular the endomorphism algebra of $\mathcal P^{spec}$ enriches from a mere algebra to a quasi-coherent sheaf of algebras over $\Loc_\Gv$. We now give a conjectural description of the restriction of this sheaf of algebras to $\Loc^\circ$.
  \begin{conjecture}  (Algebra of $L$-observables):  \label{L2conj}
Suppose that we are in the setting of the global geometric 
Conjecture \ref{GlobalGeometricConjecture} and $\check{M}$ is a symplectic vector space,
which we do not require to be polarized.

Then there is an isomorphism of quasicoherent sheaves of algebras
  on $\Loc_\Gv^{\circ}$
 $$ 
 \bO_{\Mv,\Sigma}^\circ \stackrel{\sim}{\rightarrow} \mathsf{End}_{QC(\Loc_\Gv^{\circ})}(\mathcal P^{spec}),$$
 where $\bO_{\Mv,\Sigma}^\circ$ is the sheaf of algebras constructed from $\check{M}$ in \eqref{Adef}, 
 $\mathcal P^{spec}$ is the spectral projection of $\mathcal P_X^{\norm}$ (and thus agrees with $\mathcal P_X^{norm}$ in the de Rham setting),
 and $\mathsf{End}$ is the internal endomorphisms valued in quasi-coherent sheaves, as noted above. 
  \end{conjecture}

 \medskip
 
 In other words, over the locus $\Loc_\Gv^{\circ}$ we have 
 deformation quantized the fixed points of $\Mv$ to the Clifford algebra $\bO_{\Mv,\Sigma}^\circ$ which makes sense independently of polarization data.
 Given a polarization,  the spectrally decomposed period sheaf provides a compatible geometric quantization.  
 
  For example, if $\check{M} = T^* \check{X}$ for a vector space $\check{X}$,
 with scaling $\GGm$-action, 
 let us compute as in \S \ref{locXex} (and using the same notation)
  some fibers of the above statement.
 Let $\rho$ be, as there, a $\check{G}$-valued local system
with a unique fixed point $0$ on $\check{X}$.
Under the global conjecture the
fiber of the right hand side is the same as
as 
$ \Hom(\mathcal{L}^{\norm}_{\rho}, \mathcal{L}^{\norm}_{\rho})$
with $\mathcal L^{\norm}_{\rho}$ the fiber of the normalized $L$-sheaf at $\rho$. 
 As in \eqref{LX0} 
 this space is identified up to twist with endomorphisms
 of $\Sym^* H^*(T) = \wedge^* H^1(T)$,
 where $T \simeq \check{X}$ is just the tangent
 space to $\check{X}$ at $0$. \footnote{Note that the $H^1(T)$ appearing
 in the above statement really appears in degree zero, because of shearing.}
  The assertion of the conjecture \ref{L2conj} just amounts to the fact
 that the Clifford algebra of the orthogonal  space $H \oplus H^*$ -- with $H=H^1(T_x)$ the cohomology group appearing above --
 maps isomorphically to the endomorphisms of $\wedge^* H$;
 this is the usual realization of the spin representation in presence of a polarization.

\begin{remark} (The anomaly and spectral quantization) \label{anomalyspectral}
The conjecture implies that Clifford algebra $\bO_{\Mv,\Sigma}^\circ$ should in fact split (on $\Loc_{\check{G}}^{\circ}$) --
i.e., it is isomorphic as algebra to the endomorphisms of a vector bundle.
 This issue is precisely the spectral analogue of the automorphic discussion in  \S \ref{quantization}.
Let us spell out our expectations on this issue:

The obstruction to splitting the Clifford algebra is
an element of $H^2(\Loc_{\check{G}}^{\circ}, \Z/2)$, 
which is closely related to the second Stiefel-Whitney class
of the quadratic vector bundle (see e.g.\ \cite[V.3]{Lam}). 
It is reasonable to suppose that this obstruction class arises  from the second Chern class 
$c_2 \in H^4(BG)$ from the embedding $\check{G} \rightarrow \mathrm{Sp}(\check{M})$, i.e.
$c_2$ yields by pullback a cohomology class in the space of maps from $\Sigma$ to $B\check{G}$ which can be integrated over fibers
to give a degree $2$ class in $\Loc_{\check{G}}$.   
This is related to the study of the ``Maslov cocycle.'' \footnote{It should
be able to prove  a version of this statement, at least on any field-valued point of $\Loc_{\check{G}}^{\circ}$, 
using the results of \cite{ParimalaMaslov} and the methods of Meyer \cite{Meyer}
(who worked only over the real numbers).
 Cf. also \cite{AV}.}

In particular, this would indeed imply that the Clifford algebra is split
so long as $\check{M}$ is  anomaly free in the sense
of \S \ref{quantization}, that is to say, if $c_2 \in H^4(\check{M}/\check{G})$, considered modulo $2$,
is the square of an integral class in $H^2$.\footnote{ 
Then the degree $4$ class on 
$\Sigma \times \Loc_{\check{G}}$ 
is the square of  a degree $2$ class.
All that matters for us is the $(1,1)$-component
of this class -- call it 
$c_{11} \in H^1(\Sigma, \Z/2) \otimes H^1(\Loc_{\check{G}}, \Z/2) $.
Because the squaring map $H^1(\Sigma, \Z/2) \rightarrow H^2(\Sigma, \Z/2)$
is trivial we get $c_{11}^2=0$ cf. Lemma \ref{anomalousautomorphic},
and this implies that the integrated class in $H^2(\Loc_{\check{G}},\Z/2)$ vanishes, too.
All we used, in fact,  was that $c_2$ was the square of a mod $2$ class. }

\end{remark}
 
  \begin{remark} 
  What if $\check M$ is not a symplectic vector space? There is a similar discussion but with the cohomological degrees shifted, so it looks less familiar from a classical viewpoint. For more details see \S \ref{polarized spectral}. 
 
  A simple example is $\check M=T^*(\check G/\check H)$ where (possibly after restricting to an open subset of $\Loc_G$)
  we suppose that $i: Y := \Loc_{\check H} \rightarrow Z:= \Loc_{\check G}$ is a closed immersion.  In this case, 
  the endomorphisms $\End_{\mathcal{O}_Z}(i_* \mathcal{O}_Y)$
  can be considered as a deformation quantization of the {\em relative} cotangent bundle of $Y \rightarrow Z$.
  This relative cotangent bundle is in turn identified with $\Loc^{\check M}_{\check G}$; thus we may regard the endomorphisms as a deformation of its structure sheaf. We return to this in \S \ref{spectral geometric quantization}. 
        \end{remark}

%% file: global-geometric-P1.tex
\section{The case of the projective line} \label{P1}

We consider now the case of $\mathbb{P}^1$ 
and discuss the relationship between the local conjecture and the polarized global conjecture. 
It is our expectation that the global conjecture should actually be a consequence of the local conjecture,
but we do not attempt to push this through here. Our goal is a more modest one: 
 we verify in various cases that
\begin{equation} \label{P1desid} \Hom(\mathbf{e}', \textrm{period sheaf}) = \Hom(\mathbf{0}, \textrm{$L$-sheaf}),\end{equation}
where $\mathbf{e}'$ and $\mathbf{0}$ are corresponding ``basic objects'' in the automorphic
and spectral category.\footnote{The prime on the $\mathbf{e}$ is meant to remind of certain shifts in the normalizations.}  

These computations are sufficient to support 
 various subtler points of the conjecture (the choice of $\mathcal{O}$ versus $\omega$,
the precise Tate and cohomological twists, etc.)  In particular, in cases involving
twisted polarizations on the spectral sides, the computations here are the {\em only} evidence 
that we currently have that the proposed shift in the global conjecture is correct.

We will work in the {\'e}tale setting throughout this section. 
 Both sides of \eqref{P1desid} are {\em a priori} complexes   of $\kk$-vector spaces with a Frobenius action. However, as we will see,
 this Frobenius action arises from a natural $\GGm$ action (by 
 letting Frobenius act by $q^{-a/2}$ in $\GGm$ degree $a$)
 and therefore we will prove \eqref{P1desid} as an isomorphism
 of graded complexes of $\kk$-vector spaces.

\begin{remark}
 The reader with a background in the arithmetic Langlands program
 who is only interested in numerical consequences may ask:  Why should one spend any time or effort on the ``degenerate'' case of $\mathbb{P}^1$?
The reason is that the numerical computations in the classical Langlands  program provide
good evidence for our conjecture on the cuspidal locus.  However, from this point of view, the situation with, e.g., the constant automorphic form, or other forms constructed from residues of  Eisenstein series, remains
murky. The study of $\mathbb{P}^1$  provides a toy example
where these latter complications are still present,
and can be studied without the appearance of cusp forms.   In particular,
the issues we see in this Chapter are not only geometrical in nature -- corresponding complications
would also appear in the numerical study of automorphic forms on $\mathbb{P}^1$. 
 \end{remark}

 \begin{remark}[The UFO, a.k.a. the raviolo]\index{UFO}
The global geometric conjecture on $\mathbb{P}^1$ has a close variant (which we do not attempt to state formally) provided by the global geometric conjecture on the
``UFO'' (or ``raviolo''), the non-separated curve ${\mathcal R}=D\coprod_{D^*} D$ given by two formal discs glued away from the origin.  Indeed on the spectral side the stacks of local systems on ${\mathcal R}$ and $\mathbb{P}^1$ coincide.
 The stack $\Bun_G({\mathcal R})$ is simply the equivariant affine Grassmannian, and its category of sheaves is the Hecke category. The geometric Langlands conjecture on ${\mathcal R}$ becomes the derived geometric Satake correspondence, the period sheaf on ${\mathcal R}$ recovers the Plancherel algebra, and the global geometric conjecture recovers the Plancherel algebra conjecture -- except that in all cases the algebra structures (convolution and factorization) are encoded separately in special features of the UFO. 
 \end{remark}

The contents of the section are as follows:

 \begin{itemize}
 \item In \S \ref{P1example} we discuss the geometry of bundles on $\mathbb{P}^1$. 
 \item In \S \ref{P1GL} we summarize the explicit geometric Langlands correspondence in the case of $\mathbb{P}^1$.
 \item In \S \ref{Koszul dual via volume forms} we review Koszul duality in the form that we will use it.
 \item In \S \ref{LocXstack} we describe the geometry of $\Loc^{\check{X}}$. 
 \item In \S \ref{LunP1} we compute the unnormalized $L$-sheaf.  More precisely, we compute its Koszul dual,
 and will find a very pleasant phenomenon: this Koszul dual depends on $\check{X}$ only through its contangent bundle $T^* \check{X}$.
 \item In \S \ref{Whittaker P1} we complete the computation of the $L$-sheaf in the Whittaker case, highlighting the perspective of Atiyah bundles.
  \item In \S \ref{P1Lnorm} we use the foregoing to compare normalized period and $L$-sheaves in 
  the case when $X$ has no Whittaker twist. 
  
\end{itemize}

Throughout this section we work in the {\'e}tale framework, since  
the main subtlety we hope to examine involve tracking cohomological and Tate twists. 
Accordingly  throughout this section,  we use the notation \index{$(a, b]$ twist}
\begin{equation} \label{doublebracketnotn} \Q_{\ell}(a,b] := \Q_{\ell}(\frac{a}{2}) [b]\end{equation}
i.e., a Tate twist by $a$ and a cohomological shift by $b$.  
These Tate twists will actually be kept track of by a
$\GGm$-action.

We also  denote by lower case letters the dimension of the associated algebraic variety, e.g.:
$$g=\dim(G),  b=\dim(B), x =\dim(X).$$
In other contexts we have used $g$ to mean the genus of the curve, but hopefully this will not
cause confusion in the current section since we are working only with $\mathbb{P}^1$.

\subsection{Example: $\PGL_2$-bundles on $\mathbb{P}^1$} \label{P1example}
To orient the reader (and, more importantly, the author) who has not studied this situation before, we briefly describe the simplest instance of the geometry,
when $G=\PGL_2$:

The $\F_q$-points of $\Bun_{\PGL_2}(\mathbb{P}^1)$ are parameterized by non-negative integers:
 up to twisting, each bundle is of the form 
 $$[n] := \mathcal{O} \oplus \mathcal{O}(n)$$ for a {\em unique} $n \geq 0$. 

But although the picture at the level of points is straightforward, the algebraic geometry of this situation is already somewhat
 nontrivial; in fact the closure of $[0]$ contains $[2]$, the closure of $[2]$ contains $[4]$ and so on. \
To draw pictures here we can pass to a smooth cover,
and a convenient one is the map
 $$ \mathrm{Ext}^1(\mathcal{O}(n), \mathcal{O}) \rightarrow \Bun_{\PGL_2}.$$
 Then, for $n \geq 2$ the preimage of the closures of $[n], [n-2], \dots$ 
 gives an increasing stratification of this $n-1$-dimensional vector space $\mathrm{Ext}^1$ by varieties of dimension 
 $0, 2, 4, \dots$, which captures the ``stratified topology'' of $\Bun_G$. 
 
For example take $n=4$.  We can identify $\mathrm{Ext}^1(\mathcal{O}(4), \mathcal{O})$ with
  sections 
   $ P := x^{-2} y^{-2}( a \frac{y}{x} + b + c \frac{x}{y})$
   of $ \mathcal{O}(-4)$ 
  on $\mathbb{G}_m$: use ${\tiny \left( \begin{array}{cc} 1 & P \\ 0 & 1 \end{array}\right) }$
  to glue $\mathcal{O} \oplus \mathcal{O}(4)$ on the complements
  of $0$ and $\infty$. 
   In this identification, the  closure of the
   preimages of $[4], [2], [0]$ are:\footnote{ To see this, the vector bundle $V$ attached to a polynomial $P$ belongs to the closure of
 $[2]$ exactly when $V(-3)$ has a section, equivalently, when there are sections
 $f_{1}$ of $\mathcal{O}(1)$ and $f_{-3}$ of $\mathcal{O}(-3)$
 over $\mathbb{A}^1$, 
 such that both $f_{-3} + P f_1$ and $f_1$
 extend over $\infty$. This forces $(a,b)$ and $(b,c)$ to be linearly dependent, 
 giving the cone $b^2=ac$. 
}
 $$ \mbox{origin} \subset \{b^2=ac\} \subset \mathbb{A}^3.$$
 Thus, various computations  with
  on $\Bun_{\PGL_2}$ will reflect the topology of the singularity of this cone at the origin.
If we go deeper, we see more involved varieties defined by determinants.

\subsection{Geometric Langlands for $\mathbb{P}^1$}  \label{P1GL}

\subsubsection{The automorphic side}

As usual, the automorphic side will be a suitable version of the category of {\'e}tale sheaves on $\Bun_G$, see
\S   \ref{aut sheaf notation}.
The class of the trivial bundle defines
an open inclusion $j: BG \rightarrow \Bun_G$.
Set $\mathbf{e} = j_{!} \Q_{\ell}$;
it corresponds to the automorphic function whose value is $1$ on the trivial
bundle and zero off it.  
We also put 
   \begin{equation} \label{e'def} \mathbf{e}' := \mathbf{e}(b/2-g, -g]\end{equation}
   using the notation of \eqref{doublebracketnotn}, i.e., the first coordinate is a Tate twist
   and the second a cohomological twist.
   The shifts in $\mathbf{e}'$ are adapted to the normalization of the Langlands correspondence needed for the conjecture -- recall
that our conjecture asserts there is a choice of shifts which makes period and $L$-sheaves match
for any choice of period.

We will use $\infty \in \mathbb{P}^1(\mathbb{F})$ as a basepoint.

\subsubsection{The spectral side} \label{SSP1}   
The stack $\Loc_{\Gv}$ of local systems on $\PP^1$ is identified with 
the quotient of the derived scheme $\cg[-1]$  by the adjoint action of $\check{G}$: 
$$ \Loc_{\check{G}} = \check{\mathfrak{g}}[-1]/\check{G}\simeq T[-2] (\mathrm{pt}/\Gv)$$
where   $T[-2] (\dots)$ denotes the shifted tangent bundle of $\mathrm{pt}/\Gv$; recall that the tangent complex of $\pt/\Gv$ is given by the adjoint representation in cohomological degree $-1$, i.e., $\fgv[1]/\Gv$.
  The category of coherent sheaves on $\Loc_{\Gv}$ is then 
 the category of  $\check{G}$-equivariant differential graded modules for $\Sym \  \cg^*[1]$
 with finite dimensional cohomology; here, 
 $\Sym \ \cg^*[1]$ is understood to have trivial differential.
 
The derived geometric Satake correspondence then further identifies these categories with the derived ``small'' spherical Hecke category $\Sph=\Sph_G$ for $G$.
 This is Koszul dual to the formulation of \cite{BezFink} that
 was recalled in Theorem \ref{Satake thm}; see \S \ref{Koszul dual via volume forms} for discussion of Koszul duality in this context.
 
With reference to this, the augmentation $\mathbf{0}$ of  $\Sym( \cg^*[1])$, 
considered  as a $\Sym(\cg^*[1])$ -module with its trivial $\check{G}$-equivariant structure, 
corresponds to the unit of the spherical category. 
 
 \subsubsection{The statement of the Langlands correspondence}
    The action by Hecke operators at $\infty$ on the sheaf $\mathbf{e}'$ defines a functor
$$\Sph\to \autshv(\Bun_G),\hskip.3in T\mapsto T \mathbf{e}'$$
which is an equivalence of categories \cite{LafforgueP1, ArinkinGaitsgory} -- see~\cite[Theorem 3.1.9]{darioBM} for a detailed treatment;
that is to say, in this case, 
   we have a factorization of the Langlands correspondence, denoted by the dashed arrow below: 
\begin{equation} \label{P1LD} 
\xymatrix{
  \Sph \ar[d]^{act_{\infty}e} \ar[rr] & & \mbox{Coherent sheaves on $\cg[-1]/G$}  \ar[d] \\ 
 \mbox{$\autshv(\Bun_G)$}    \ar@{.>}[rr] && \Coh(\Loc_{\check{G}})
 }
 \end{equation}
where the top vertical map is that of geometric Satake. On the left side, 
we have the small  versions of the categories of constructible sheaves\footnote{i.e., constructible sheaves $!$-extended
from quasicompact substacks} and on the right side
we have categories of coherent sheaves; we can then pass to Ind-categories
everywhere to recover the ``large'' equivalence.

\subsubsection{Incorporating Frobenius}
The diagram \eqref{P1LD} is compatible with Frobenius structures, i.e., 
we have a corresponding diagram where 
we impose mixed
Weil sheaves on the left\footnote{In other words, we consider on the left only those complexes of sheaves on $\Bun_G$ or
the affine Grassmannian 
which are equipped with an isomorphism $\mathrm{Frob}^* V \simeq V$,
with the property that each cohomology sheaf has a filtration whose
associated graded sheaves are pointwise pure, that is to say,
with reference to an isomorphism $\kk \simeq \C$
all the Frobenius eigenvalues on stalks have absolute value $q^{w/2}$ for some {\em integer}
$w$. 
}, and $\GGm$-equivariant coherent sheaves at the right.  
Here $\GGm$ is acting by $x \mapsto x^{-2}$ on $\cg$ (this
will become more familiar after we pass to the Koszul dual picture,
where it will be acting by $x \mapsto x^2$ on $\fgxv$). And 
with respect to this equivalence, Frobenius on the left corresponds 
to $q^{-1/2} \in \GGm$ on the right; said differently, it acts by shearing.

\subsection{Koszul duality and volume forms}\label{Koszul dual via volume forms}
\index{Koszul duality} In our later computations, we will see that the $L$-sheaf
becomes much simpler viewed through the lens of Koszul duality: 
Let us write $\mathbf{0}$ for the augmentation of   $\Sym \  \cg^*[1]$.
Then
$\mathrm{End}(\mathbf{0}) \simeq \mathrm{Sym} \  \cg[-2]$.
More generally, 
the functor $$M \rightsquigarrow \Hom(\mathbf{0}, M),$$
 of Koszul duality
 carries  differential graded $\Sym( \cg^*[1])$ modules to differential graded $ \mathrm{Sym}(\cg[-2])$ modules.
 We refer to \S \ref{shearing affine line} for a discussion of the extent to which it is an equivalence
 (this depends on the precise finiteness conditions on both sides). 
 
 It induces 
 a similar functor on $\check{G}$-equivariant modules. 
It will be convenient to adopt the convention that {\em $\fgxv$ lies
in $\GGm$-degree $2$.} With this convention,
we can write 
$$ \Sym \cg[-2] = \mathcal{O}(\fgxv)^{\shear}.$$

It is useful to note that the Koszul dual symmetric algebra $\mathrm{Sym} \cg[-2]\simeq \cO(\fgxv)^{\shear}$ to the exterior algebra  
$\cO(\fgv[-1])$ is naturally realized as the convolution algebra of volume forms $\omega(\fgv[-2])$ on the group-stack $\Omega(\fgv[-1]) \simeq \fgv[-2]$
of loops in $\fgv[-1]$.
Indeed from the point of view of derived algebraic geometry, the Ext-algebra of a skyscraper always appears as the convolution algebra of volume forms on the based loop space.
 This realization comes from proper adjunction and base-change for the diagram
$$\xymatrix{
\fgv[-2]\ar[r]^-{p}\ar[d]^-{p} & 0 \ar[d]^-{i}\\
0\ar[r]^-{i} & \fgv[-1]
}
$$
as follows:
\begin{eqnarray*}
\Hom(i_* k, i_* k) &\simeq& \Hom(k, i^! i_* k)\\
&\simeq& \Hom(k,p_*p^! k)\\
&\simeq& \omega(\fgv[-2])
\end{eqnarray*}
 We will need the following basic Koszul duality calculation  
that generalizes the above:

\begin{lemma}\label{family Koszul duality}
Let $p:W\to Z$ be a vector bundle over a smooth stack, $p:V=W[-1]\to Z$ its shift (a derived vector bundle) 
and $Z\times_V Z\simeq W[-2]$ the self-intersection of the $0$ section $i:Z\to V$. 
Then we have an equivalence $$\omega(Z\times_V Z)\simeq \Gamma(Z, p_* \underline{\End}(i_* \cO_Z) \otimes \omega_Z).$$
\end{lemma}

In other words, algebraic distributions on the total space of $W[-2]$ are identified with an $\omega_Z$-twist of the Koszul dual $\underline{\End}(i_* \cO_Z)$ of the exterior algebra of functions on $W[-1]$, which is the graded symmetric algebra of functions on $W^*[2]$.

\begin{proof}
The lemma follows from the relative version of the Koszul duality calculation of \S \ref{indcoh for exterior algebra} over the smooth base $Z$ (we use the smoothness of $Z$ to identify $\QC(Z)$ and $\QC^!(Z)$ via the inverse equivalences $\Xi,\Psi$).
In other words, we describe the category $QC^!(W[-1])$ by Barr-Beck-Lurie as modules for the algebra object $S^\shear=\underline{\End}(i_* \cO_Z) \in Alg(QC(Z))$, where $S=Sym(W)$ is the graded symmetric algebra on $W$. This equivalence identifies the adjunction $(i_*,i^!)$ is with the tensor-hom adjunction, tensoring with $S^\shear$ and forgetting the $S^\shear$-action. 

Consider the pullback diagram
 $$\xymatrix{ Z\times_V Z \ar[r]^-{\Pi}\ar[d]^-{\Pi} & Z\ar[d]^-{i}\\
 Z \ar[r]^-{i} & V}.$$
The volume forms $\omega(Z\times_V Z)$ are calculated as global sections on $Z$ of  $$\Pi_*\omega_{Z\times_V Z}\simeq \Pi_*\Pi^!\omega_Z \simeq i^!i_*\omega_Z.$$
Thus we can rewrite
\begin{equation} \label{desco} \omega(Z\times_V Z)\simeq \Hom_{QC(Z)}(\cO_Z, i^!i_*\omega_Z).\end{equation}

Under the Koszul duality equivalence, $ i^!i_*\omega_Z$ is identified with the $S^\shear$-module $S^\shear \otimes \omega_Z$ with $S^\shear$-action forgotten, whence the result follows. \footnote{Alternately, we can proceed
from \eqref{desco}   by rewriting the right hand as  $\Hom_{V} (i_! \cO_Z, i_* \omega_Z)$,
and using adjunction to express this as   $\Hom_V(i_*\cO_Z, i_* \omega_Z)= \Hom_{\QC^!(V)}(i_* \cO_Z, i_* \cO_Z) \otimes p^* \omega_Z$.
At the last step, we use the fact that $\omega_Z$ is a line bundle.}
  \end{proof}

\subsection{The stack $\Loc^{\Xv}$}  \label{LocXstack} 
Let $\Xv$ be a smooth affine $\Gv \times \GGm$-variety and $\Psi\to \Xv$ an equivariant $\mathbb{A}^1$-bundle
where as usual $\GGm$ acts on $\mathbb{A}^1$ by squaring. As before, we may consider the derived stack  
$$\pi:\Loc^\Xv=\Map_{\mathrm{lc}}(\PP^1, \Xv/\Gv)\longrightarrow \Loc_\Gv $$ of locally constant maps from $\PP^1$ to the quotient stack $\Xv/\Gv$, or equivalently, of
$\Gv$-local systems equipped with a locally constant section of the associated $\Xv$-bundle.

A useful description of $\Loc^\Xv$, well adapted to the Koszul dual description of sheaves on $\Loc_\Gv$, comes by taking the fiber over the trivializable local system locus $i:\pt/\Gv\to \Loc_\Gv$:
we find a pullback diagram 
 \begin{equation}\label{shifted tangent diagram}
 \xymatrix{(T[-2]\Xv)/\Gv\ar[r]^-{\wt{i}}\ar[d]^{\tilde{\pi}}& \Loc^\Xv\ar[d]^-{\pi}\\
 \pt/\Gv \ar[r]^-{i}& \Loc_\Gv }
 \end{equation}
 where the shifted tangent complex $T[-2]\Xv$ is identified with the functor of locally constant maps from $\PP^1$ to $\Xv$. 

 \begin{remark} \label{P1locprez}
Although we won't use it, an alternate presentation  which makes clear the geometric nature is 
that  $\Loc^\Xv$ is the quotient by $\Gv$ of the fiber product  
$$\mathcal{Y} =  \Xv\times_{I(\Xv)} \Xv,$$
where  $I(\Xv)$ is the locus (in the derived sense) of pairs $(g,x)$ with $gx=x$;
and the map $\Xv \rightarrow I(\Xv)$ sends $x$ to $(\mathrm{id}_G, x)$. 
Informally, the two copies of $\Xv$ give the section above lower and upper hemispheres
and the class in $I(\Xv)$ gives the gluing datum.
\end{remark}
 \subsection{Computation of the unnormalized $L$-sheaf} \label{LunP1}

  We continue in the generality of the previous section \S \ref{LocXstack}.
Our goal in this section is to describe the unnormalized $L$-sheaf of $\Xv$, i.e., $(\pi_*\omega_{\Loc^{\check{X}}})^{\shear}$, and its Whittaker version in which $\omega_{\Loc^\Xv}$ is replaced by a $\Psi$-twisted version (the pullback of the spectral exponential sheaf). The answer is given explicitly in terms of the geometry of the associated Hamiltonian $\Gv$-spaces $T^*\Xv$ and $T^*_\Psi\Xv$.

Let us first of all observe that the moment map $T^* \Xv \rightarrow \fgxv$
is equivariant for $\GGm$ actions, where we  modify the natural $\GGm$ action by the action by squaring along fibers of $T^* \Xv$. 
Therefore, we obtain a morphism 
$$ \mathcal{O}_{T^*\Xv[2]} = \mathcal{O}(T^* \Xv)^{\shear} \leftarrow \cO(\fgxv)^{\shear},$$
where, on the left, the ring of functions on $T^*\Xv$
but sheared so that linear functions on the fiber are taken in cohomological dimension $2$.
 This will be implicitly used in the following statement:

 \begin{proposition} \label{Lsheafcomputation}  
 (The Koszul dual to the $L$-sheaf is volume forms on $T[-2] \Xv$): 
  The Koszul dual $\underline{\Hom}(\mathbf{0}, L_{\Xv,\Psi})$ of the $L$-sheaf is given as the $(0,-g]$-shift of the following $\Gv$-equivariant module for $\cO(\fgxv)^{\shear}$:
\begin{equation} \label{KoszulLresult} \Gamma(\Xv, \mathcal{O}_{T^*_\Psi \Xv} \otimes_{\mathcal{O}_\Xv} \omega_\Xv)^{\shear}\end{equation}
 that is to say,  global sections of the line bundle $p^*\omega_\Xv$ on $T^*[2]\Xv$
 or the twisted version thereof.  
 
 Moreover, the above isomorphism is Frobenius equivariant,
where we have endowed \eqref{KoszulLresult}
with a Frobenius action through the shearing of
the trivial Frobenius action  on  $\mathcal{O}_{T^*_{\Psi}} \Xv, \omega_{\Xv}, \mathcal{O}_{\Xv}$
through the $\GGm$-action (see  \S \ref{shearing Frobenius notation} for generalities).

 \end{proposition} 
 
 Note that, when we write $\underline{\Hom}(\mathbf{0}, \dots)$, we mean
that we take $\Hom$ ``relative to $\mathrm{pt}/\check{G}$,''
i.e., we take the sheaf-Hom
and regard it  by pushforward as a sheaf on $\mathrm{pt}/\check{G}$,
so that the result
is a $\check{G}$-representation.
 
Although the twisted case of the Proposition of course includes the untwisted one, we will prove them separately to try
to distinguish twistedness from other aspects. 

\subsubsection{Proof of Proposition~\ref{Lsheafcomputation}: the untwisted case}

  Consider the pullback diagram~\ref{shifted tangent diagram}.
 Since $i$ is proper  we may rewrite  the desired $\underline{\Hom}(\mathbf{0}, L_{\check{X}})$ as
\begin{equation} \label{first_computation} \underline{\Hom}(k, i^! \pi_*\omega_{\Loc^\Xv})\simeq  \wt{\pi}_*\wt{i}^! \omega_{\Loc^\Xv}\simeq  \wt{\pi}_* \omega_{T[-2]\Xv / \Gv}\end{equation}
 since inner Hom from the trivial representation $k_{\pt/\Gv}$ is the identity, and applying base change. 
 In other words, {\em  the Koszul dual of the L-sheaf is given by the $(0,-g]$-shift of volume forms on $T[-2]\Xv$, as a $\Gv$-representation. }
 The shift arises from the fact 
 that what appears in \eqref{first_computation} is the dualizing sheaf of $T[-2] \Xv/\Gv$,
 rather than $T[-2] \Xv$, and   
 the dualizing sheaf of $BG$ is not $\kk$ but rather $\kk(0,-g]$.

This proves the claim \eqref{KoszulLresult} as an isomorphism of $\check{G} \times \GGm$ modules
but we must also check the structure  as module for $\cO(\fgxv)^{\shear}$. For this,
we must compute the module structure under  $$\underline{\End}(i_*k)\simeq (\omega(\fgv[ -2]),\ast)\simeq (\cO(\fgxv[ +2],\cdot),$$
 where the isomorphisms are as discussed in \S \ref{Koszul dual via volume forms}. 
  In order to compute this we now perform Koszul duality on $\Xv$, repeating the argument of Lemma~\ref{family Koszul duality} with $V=(T[-1]\Xv)/\Gv\to Z=\Xv/\Gv$, but now working equivariantly for the action of additive groups over $\pt/\Gv$ as follows:

$$\xymatrix{\fgv[-2]/\Gv\ar[r]\ar[d] & \pt/\Gv \ar[d]&\actson&T[-2]\Xv/\Gv \ar[r]^-{p}\ar[d]^-{p} & \Xv/\Gv\ar[d]^-{i}\\
\pt/\Gv \ar[r] & \fgv[-1]/\Gv&\actson& \Xv/\Gv \ar[r]^-{i} & T[-1]\Xv/\Gv}.$$

Here the vector group $\fgv[n]$ (i.e., $\fgv[n]/\Gv\to \pt/\Gv$) acts on $T[n]\Xv$ (again connoting $T[n]\Xv/\Gv\to \pt/\Gv$) by the (shifted) derivative of the $\Gv$-action.
Each step of the identification of Lemma~\ref{family Koszul duality}
$$ \omega(T[-2]\Xv) \simeq p_* \underline{Hom}(i_*\cO_\Xv, i_*\cO_\Xv) \otimes \omega_\Xv$$ is now compatible with the action of $\omega(\fgv[-2])$ as endomorphisms of the functor $i_*$.   In particular $\omega(\fgv[-2])$ acts through its action on the first factor,
$$\cO(T^*[2]\Xv)= p_* \underline{Hom}(i_*\cO_\Xv, i_*\cO_\Xv),$$ 
where it is identified (through the $\Gv$-equivariant algebra isomorphism   
 $(\omega(\fgv[-2]),\ast)\simeq (\cO(\fgxv)^{\shear},\cdot)$ of Section~\ref{Koszul dual via volume forms}) with the moment map action of $(\cO(\fgxv)^{\shear},\cdot)$, as claimed.

\subsection{The Whittaker $L$-sheaf on ${\mathbb P}^1$.}\label{Whittaker P1}

In the twisted case, we  are going to apply the equivariant version of Lemma~\ref{family Koszul duality} as in the untwisted case, but now with the total space $\underline{\Psi}$ taking the role of $\Xv$ and the group $\Gv\times \Ga$ taking the role of $\Gv$. 
To guide us in this argument, we will recall some facts about Atiyah bundles: the Atiyah bundle $At_\cP=(T\cP)/H$ of a principal $H$-bundle $\cP\to \Xv$ is, by definition,  the quotient of the tangent bundle of $\cP$ by $H$.
It fits
 into a fiber sequence  $$\ad_\cP\to At_\cP \to T\check{X}.$$
 By rotating this triangle we can realize the Atiyah bundle as the fiber of the tangent map $T\check{X}\to ad_\cP[1]=\psi^* T(\pt/H)$ to the map $\psi:X\to \pt/H$ classifying $\cP$, i.e., as the following pullback:
\begin{equation} \label{pbd} \xymatrix{ At_\cP\ar[r]\ar[d] & T\Xv\ar[d]^-{\psi}\\
\pt/H \ar[r] & T(\pt/H)}.\end{equation}

\subsubsection{Proof of Proposition~\ref{Lsheafcomputation}: the twisted case}

 Recall the definition of L-sheaves in the Whittaker setting, Definition~\ref{spectral Whittaker def}, in which the role of $\omega_{\Loc^\Xv}$ is 
taken up by the pullback $(\overline{\Psi}^!)^\shear exp$ of the spectral exponential sheaf (combined with shearing). 

As in the untwisted case, we first realize this $L$-sheaf as the sections of the pulled back exponential sheaf (rather than dualizing sheaf) on $T[-2]\Xv$.
Referring again to the diagram
 $$ \xymatrix{(T[-2]\Xv)/\Gv\ar[r]^-{\wt{i}}\ar[d]^{\tilde{\pi}}& \Loc^\Xv\ar[d]^-{\pi} \ar[r]^{\bar{\Psi}}  &\mathbb{A}^1[-1] \\
 \pt/\Gv \ar[r]^-{i}& \Loc_\Gv  }
$$
we  compute, as in \eqref{first_computation}, 
\begin{equation}  \label{third computation}
 \underline{\Hom}(\mathbf{0}, L_{\check{X}})= 
 \unshear \ \circ ( i^! )^{\shear}  \pi_*^{\shear} (\overline{\Psi}^!)^\shear \exp  \simeq  \unshear \circ \wt{\pi}_*^{\shear}   (\overline{\Psi} \circ \tilde{i}^!)^{\shear} \exp.\end{equation}
 The $\unshear$ at the end is as in   Definition~\ref{spectral Whittaker def}: to remind us that the
 term $\left[  \wt{\pi}_*^{\shear}   (\overline{\Psi} \circ \tilde{i}^!)^{\shear} \exp \right]$
 formally speaking lives in $QC^{!}(\pt/G)^{\shear}$, but
 since the $\GGm$ action is trivial, this is identified with $QC^{!}(\pt/G)$ itself (cf. Example \ref{trivial action}); 
 we will however not explicitly mention this in the analysis that follows. 
 
 The right hand side of \eqref{third computation}, said in words, 
describes the sections of $(\bar{\Psi} \circ \tilde{i}^!)^{\shear}  \exp$ on $T[-2] \check{X}$;  we shall prove
\begin{equation} \label{baltimore} \mbox{ sections of } (\bar{\Psi} \circ \tilde{i}^!)^{\shear}  \exp  \mbox{ on $T[-2] \check{X}$} 
  = \Gamma(\Xv, \mathcal{O}_{T^*_\Psi \Xv} \otimes_{\mathcal{O}_\Xv} \omega_\Xv)^{\shear}.\end{equation}
We are going to compute the right hand side
in terms of volume forms on a suitable Atiyah bundle.  Notice  first that $\mathcal{O}_{T^*_{\Psi} \Xv}^{\shear}$  
 correspond to functions on  the shifted hamiltonian reduction  of $T^* \Psi$ by $\Ga$
 and  can be computed as follows:  
 \begin{equation} \label{exp via shifted skyscraper}
 \cO(T^*_{\Psi} \Xv)^{\shear}  \simeq \cO(T^*  \underline{\Psi})^{\shear, \Ga}\otimes_{\cO(\fg^*_a[2])^{\shear}} k_1.\end{equation}  
   where $k_1$ is a skyscraper sheaf at $1 \in \fg^*_a$, cf. Definition \ref{sesdefn}.

 Consider  now the 
 derived vector bundle on $\Xv=
  \underline{\Psi}/\Ga$ defined as
  $$At_\Psi[-2]:=(T[-2] \underline{\Psi})/\Ga.$$
  We may think of it as a
 shifted version of the Atiyah bundle  
associated to the principal $\Ga$-bundle $\Psi \rightarrow \check{X}$.
  Using Lemma \ref{family Koszul duality}, applied
  to  $At_{\Psi}[-1] \rightarrow \Xv$ and its zero section,    we find an identification of volume forms on $At_\Psi[-2]$ with sections of a line bundle on the shifted cotangent bundle
$$\omega(At_\Psi[-2])^{\shear} \simeq ( \mathcal{O}(T^*\underline{\Psi})^{\Ga } \otimes_{\cO_\Xv} \omega_\Xv)^{\shear}.$$
Moreover this equivalence respects not only the shifted Hamiltonian $\Gv$-action as before, but also the action of $(\omega({\fg}_a[-2]),\ast)\simeq (\cO({\fg}_a^\ast[2]),\cdot)$. 
 Therefore, the right hand side of \eqref{baltimore} equals 
\begin{equation} \label{LHSRHS}   \omega(At_\Psi[-2])^{\shear} \otimes_{\omega(\fg_a[-2])^{\shear}} k_1\end{equation} (compatibly with shifted hamiltonian $\Gv$-actions)
and we will now  prove that the LHS of \eqref{baltimore} is expressed
by the same formula. 

Via   (a shifted version of) the discussion of \eqref{pbd},   we have the pullback diagram  
 $$\xymatrix{ At_\Psi[-2]\ar[r]^-{\wt{q}}\ar[d]^-{\pi} & T[-2]\Xv\ar[d]^-{\psi}\\
 pt \ar[r]^-{q} & \AA^1[-1]/\Ga}.$$
 The $\Ga$-action on $\AA^1[-1]$ is trivial, hence $\AA^1[-1]/\Ga=\AA^1[-1] \times B\Ga$, and we will now discuss the exponential sheaf (an object in $\QC^!(\AA^1[-1])^{\shear}$) as living on $\AA^1[-1] \times B\Ga$, by tensoring with the trivial sheaf on $B\Ga$. The exponential sheaf  is defined as the twisted coinvariants $q_*^{\shear} k \otimes_{\omega(\fg_a[-2])^{\shear}} k_1$ of the sheared skyscraper $q_*^{\shear} k
$ (considered as an object of the sheared, ind-coherent category) by the action of its endomorphisms $\omega(\fg_a[-2])^{\shear}$. 
By base-change  we have
$$( \psi^{!} q_{*} )^{\shear} k = (\wt{q}_* \pi^{!})^{\shear} k = \wt{q}_*^{\shear} \omega_{At_{\Psi}[-2]}^{\shear}$$
and passing to twisted coinvariants we get:
 $$ \psi^{!} \exp = ( \wt{q}_*^{\shear} \omega_{At_{\Psi}[-2]}^{\shear})
 \otimes_{\omega(\fg_a[-2])^{\shear}} k_1.$$
 The left-hand side of \eqref{baltimore}
 also coincides with sections of $\psi^{!} \exp$, and therefore
 coincides with   \eqref{LHSRHS}, concluding the proof.

 \subsection{Normalized period and $L$-functions} \label{P1Lnorm}
 
 We are now ready to compare the normalized period and $L$-sheaves, or, more precisely,
   we want to  compare \index{$P$  (case of $\mathbb{P}^1$)}
   \index{$L$  (case of $\mathbb{P}^1$)}\index{$P'$  (case of $\mathbb{P}^1$)}
\begin{equation} \label{Pprimedef} \underbrace{ \Hom(\mathbf{e}',  \textrm{normalized period sheaf})}_{P'} \mbox{ and } \underbrace{ \Hom(\mathbf{0}, \textrm{normalized $L$-sheaf})}_L.\end{equation}
We will also use $P$ to denote the analogue of the left hand side defined with $\mathbf{e}$ instead of $\mathbf{e}'$. 
We also use a superscript or subscript ``norm'' for the normalized analogue.
 
We will consider the following setting: 
 \begin{itemize}
 \item $X$ a vector bundle over a homogeneous affine $G$-variety, i.e., 
 $X = G \times_{H} V$
 for an $H$-representation $V$; here
 we suppose that $G, H$ are split reductive over $\FF_q$.
 \item $(\check{X}, \check{\Psi})$  is as in
 \S \ref{LocXstack}, a smooth affine $\check{G}$-variety and $\AA^1$-bundle,
 with $\GGm$ actions acting by squaring along $\AA^1$; 
 \item  The consequence \eqref{invquot} of the local conjecture  holds:
the $G$-equivariant cohomology of $X$ is identified 
with the $\check{G}$-invariant functions on $T^*(\check{X}, \check{\Psi})^{\shear}$,
compatibly with Frobenius where the Frobenius action on invariant functions
comes from shearing the trivial action. 
   \end{itemize}
 
 We are, of course,  interested in the case when $(G, T^*X)$ and $(\check{G}, T^*(\check{X}, \check{\Psi}))$
 form a hyperspherical dual pair; but the above is all we will actually use.
 
  In terms of the larger scheme of this paper,  the most important
 case of the following computations   is the case of spectral Whittaker (i.e., $X=\pt$): this gives an important
 data point that the shift is correct, whereas for automorphic Whittaker the corresponding fact is well-attested by numerical computations e.g.\ \S \ref{Whitexample}.
  
 \subsubsection{The normalized period sheaf} \label{shift0}
   
The fiber of $\Bun_G^X$ over the trivial bundle is  ``the 
space of sections of $X \otimes K^{1/2}$,''
which, here,  amounts to a (necessarily constant) map $\mathbb{P}^1 \rightarrow G/H$, 
and then a section of $V \otimes K^{1/2}$ over $\mathbb{P}^1$,
where $V$ really means the pullback of the vector bundle $X\rightarrow G/H$ to $\mathbb{P}^1$,
and is therefore a trivial bundle. 
Since $K^{1/2}$ has degree $-1$
the bundle $V \otimes K^{1/2}$ has no nonzero sections,  
and so   the fiber of $\Bun_G^X$ is simply
$\bar{X} := G/H$. Consequently,  by the definition
\eqref{PXnormdef}, we have
\begin{equation} \label{Pnaive} P := \Hom(\mathbf{e}, P_X^{\norm} ) = H^*_{c,G}(\bar{X}) \langle -h \rangle \end{equation} where
we used $\beta_X = -h (=-\dim H)$ in the case at hand,  see \eqref{betaXdef}, and
we write $H^*_{c,G}(\bar{X}) $ for $G$-equivariant cohomology with compact support ``along $X$'', which is formally defined as
the cohomology of $BG$ with coefficients in  the compactly supported pushforward $\pi_{!} \kk$
along $\pi: \bar{X}/G \rightarrow BG$.
Here, and below, identifications such as that of \eqref{Pnaive} are understood to be Frobenius-equivariant. 
 
 We  will need to relate this $H^*_{c,G}$ to the usual equivariant cohomology,
 i.e., without compact support conditions:

 \begin{lemma}
Continue with the notation $\bar{X}=G/H, h=\dim(H), g=\dim(G)$. So long that
the residue characteristic of $\FF$ is sufficiently large, we have:
\begin{equation} \label{Ppnormcalc} H^*_{c,G}(\bar{X}) \simeq H^*_G(\bar{X})   (\frac{h-g}{2} +\frac{r_G-r_H}{2}, h-g].\end{equation}
Here $r_G, r_H$ are the ranks (dimension of maximal tori) of $G, H$ respectively.
\end{lemma}

\proof 
  $H^*_c(\bar{X})$ is a $H^*(\bar{X})$-module, and we will first of all prove it is free of rank one.
 By standard comparison arguments it suffices to verify this in the singular setting and with $\bar{X}$  considered over $\C$.  
 (This is where the assumption that the characteristic of $\FF$ is sufficiently large comes in.
 Presumably we could sharpen this by using the compactification  theory of $\bar{X}$.)

By Poincar{\'e} duality -- now working with singular cohomology of complex points -- 
$$H^*_c(\bar{X}) \simeq H_{2 \dim(\bar{X})-*}(\bar{X})$$
 so it is sufficient to verify that homology of $\bar{X}$  is free as a cohomology module.
Here we can replace $X$   by the homotopy equivalent  $G^{\textrm{cpct}}/H^{\textrm{cpct}}$,
a compact manifold (here a superscript cmpt marks the compact form of an algebraic group over $\C$),
where the claim is simply Poincar{\'e} duality. 

  Returning now to the {\'e}tale setting, let  $\alpha_X$ be a generator
  for $H^*_c(\bar{X})$ over $H^*(\bar{X})$. Then $\alpha_X$
  is in degree $g-h$;  for $H^{g-h}_c(\bar{X})$ is one-dimensional   and there is no compactly supported cohomology below degree $g-h$. 
  
  In particular $\alpha$ is an eigenvector for Frobenius.
  We will show that its eigenvalue is $q^w$ with 
  \begin{equation} \label{wcomp} w =  \frac{1}{2} (g-h + r_H-r_G).\end{equation}
In fact since cupping with $\alpha_X$ gives an isomorphism of compactly supported
and ordinary cohomology we deduce, after taking Frobenius traces, 
\begin{equation} \label{msrtp}
\mathrm{tr}(\Fr| H^*_c(X, \Q_{\ell})) =  (-1)^{g-h} q^w \mathrm{tr}(\Fr | H^*(X, \Q_{\ell})),\end{equation}
 i.e., if we write $X(q)$ for the number of points of $X$ over the finite field with $q$ elements, we 
 obtain by \eqref{msrtp} and Poincar{\'e} duality
 the numerical consequence
\begin{equation} \label{numcons} X(q)= (-1)^{g-h} q^w q^{\dim X} X(q^{-1}) \implies q^{w} = (-1)^{g-h} \frac{X(q)}{q^{\dim X} X(q^{-1})}.\end{equation}
To compute $w$ from 
\eqref{numcons},  we compare \eqref{numcons} 
with the  parallel computations for $G$ and $H$ separately.  
  We see that the similarly defined $w_G$ and $w_H$  (i.e.
  if we replace the role of $X$ above by $G$ or $H$)
  are the dimensions
  of the respective unipotent radicals,  and then we get $w=w_G-w_H$
  proving \eqref{wcomp}. 
 
Observe that $\alpha_X$ is represented by an equivariant class: it extends uniquely to $\tilde{\alpha}_X \in H^{g-h}_{G, c}(\bar{X})$; 
this follows from  the Leray spectral sequence associated to $X/G \rightarrow BG$,
because, for   degree reasons, there is no differential that can kill $\alpha_X$,
and no other $E^2$ term that can contribute. 
 Let $\pi : \bar{X}/G \rightarrow BG$; the action of compactly supported cohomology of fibers on cohomology of fibers corresponds 
at the sheaf level to a product
$$ \pi_{!} \kk \otimes \pi_* \kk \rightarrow \pi_{!} \kk$$
Therefore $\tilde{\alpha}_X \in H^{g-h}(\pi_{!} k)$ defines a map $\pi_* k(-w, h-g] \rightarrow \pi_{!} k$ of sheaves on $BG$, which is in fact an isomorphism. 
  Taking cohomology gives the desired result.
 \qed

 \subsubsection{The normalized $L$-sheaf}

 We must first take into account the effect of normalization on the $L$-sheaf,
 a truly depressing process because it is all about signs. 
 Recall the definition from \eqref{LXnormdef}: $ \mathcal L_{\Xv}^{\norm} =
 \mathcal L_{\Xv} \otimes   \varepsilon_{1/2}^{\vee} \langle -\beta_{\check{X}} \rangle$. 
 
 This  
 $\varepsilon_{1/2}^{\vee}$
 is pulled back,  via the  eigenmeasure character $\eta = \eta_{\spec}: \check{G} \rightarrow \Gm$, from a line bundle on 
  $\Loc_{\Gm} \simeq \mathfrak{g}_m[-1]/\Gm$, the quotient taken
  with trivial action. This line bundle    is seen  to be the trivial bundle on $\mathfrak{g}_m$ with the scaling action of $\Gm$ 
 along fibers (interpreting e.g.\  $\epsilon_{1/2}^{\vee}$  as the square root of the determinant of cohomology,
 see Remark \ref{det of coh remark}).   From \eqref{Lsheafcomputation} we  deduce that
  $$\underline{\Hom}(\mathbf{0}, \mathcal{L}_X^{\norm})  =
 \Gamma(\Xv, \mathcal{O}_{\check{M}} \otimes_{\mathcal{O}_\Xv} \omega_\Xv)^{\shear} \otimes \varepsilon_{1/2}^{\vee} \langle - \beta_{\Xv} \rangle  (0,-g],$$
 where as usual 
 $$\check{M} = T^*_{\Psi}\Xv.$$
 
The existence
of a global differential form means that $\omega_{\Xv}$ is trivial,
but it is not trivial equivariantly for the action of $G$.   Indeed, fixing  a global differential form $\omega$ on $\Xv$;
 we have $g^* \omega = \eta(g) \omega$; the {\em left} action of $\check{G}$ on forms
 is via $(g^*)^{-1} =\eta(g)^{-1}$, since $\check{G}$ is acting on the left on $\check{X}$.
On the other hand $\check{G}$ is acting on $\varepsilon_{1/2}^{\vee}$
through $\eta$. 
  Therefore, as far as the $G$-action is concerned, the twists coming from $\omega_{\Xv}$ and $\varepsilon_{1/2}^{\vee}$ cancel with one another.

  Now we   consider the $\GGm$ action.   
The left action of $\lambda \in \GGm$ on forms on $\check{X}$
is given by  pullback through $\lambda^{-1}$. 
Therefore,  cf. \eqref{gammaXdef}, 
 the  $\GGm$ action on $\omega_{\check{X}}$ is via $\lambda \mapsto \lambda^{-\check{\gamma}}$.
 
 The above discussion allows us to eliminate the role of $\omega_{\check{X}}$: 
 $$\uHom(\mathbf{0}, \mathcal{L}_X^{\norm})  =
\mathcal{O}_{\check{M}}^{\shear} (0, \check{x}-g] \langle-\gamma_{\check{X}} -\beta_{\check{X}} \rangle.$$
where the cohomological twist $[\check{x}]$, with $\check{x} = \dim\check{X}$,  arises from the cohomological shift in $\omega_{\Xv}$ itself.

Recall from \eqref{betaXdef} that $\beta_X = - (\dim G + \gamma_X - \dim(X))$
in the case of genus zero
 and so
$ - \beta_{\check{X}} - \gamma_{\check{X}} = g-\check{x} $.
 Thus
$$
\uHom(\mathbf{0}, \mathcal{L}_X^{\norm})  = \mathcal{O}_{\check{M}}^{\shear} (\frac{g-\check{x}}{2}, 0].$$
Now, pass to $\check{G}$ invariants and use 
the fact that, as stated in the setup at the start of \S \ref{P1Lnorm}, we are assuming \eqref{invquot} to be true,
in order to compute these invariants in terms of $X$. The result is 
\begin{equation} \label{Lpnormcalc}
 L := \Hom(\mathbf{0}, \mathcal{L}_X^{\norm})  =  
 H^*_G(X) (\frac{g-\check{x}}{2}, 0].\end{equation}

 \subsection{Comparison of automorphic and spectral sides} \label{P1finalcomparison}

Combining \eqref{Pnaive} and \eqref{Ppnormcalc}  we get
 $$ P := H^*_G(G/H)   (\frac{-g}{2} +\frac{r_G-r_H}{2}, -g]$$
Because of our normalization, we want to instead use the twist $P'$ 
of \eqref{Pprimedef}, i.e., $\Hom(\mathbf{e}', P_X^{\norm})$
 rather than $\Hom(\mathbf{e}, P_X^{\norm})$.  The two are related
 via \eqref{e'def}:
 $$P'  =  P (g-b/2,g] =H^*_G(G/H) (\frac{g-b + r_G-r_H}{2}, 0].$$
 The global conjecture asserts that this should coincide with $L$
 as computed in \eqref{Lpnormcalc}:
 \begin{equation}  
 L := \Hom(\mathbf{0}, \mathcal{L}_X^{\norm})  =  
 H^*_G(X) (\frac{g-\check{x}}{2}, 0].\end{equation}

 The cohomological shifts match
 and so all that remains is to check
 the Tate twist:
   \begin{equation} \label{crux} \check{x} \stackrel{?}{=} b-r_G+r_H.\end{equation}
 We do not have a direct proof of this, but it can be checked case-by-case in all examples.
 In fact, it would follow from the validity of our general duality proposal (Expectation \ref{anomaly expectation})
 i.e., the conjecture that applying the construction of \S \ref{dualofX} to $(\check{G}, \check{X}, \check{\Psi})$, when applicable, 
reconstructs $(G, X)$.   This would imply that the unipotent
 radical of $\check{G}$ acts freely on the open orbit of $\check{X}$, and  
 the quotient is identified with an $r_H$-dimensional torus, proving
 \eqref{crux}. 
 
  To be precise, then, we have verified that, starting from a pair $(G, X)_{/\FF_q}$ and 
$ (\check{G}, (\check{X}, \Psi))_{/\kk}$ as described after \eqref{Pprimedef}, 
and further assuming that:
\begin{itemize}
\item[-] the characteristic of $\FF$ is sufficiently large 
{\em and}
\item[-]
  the equality \eqref{crux} of numerical invariants is valid -- which, as we recall, 
 would follow from the fact that our construction of $\check{M}$
 from $X$ is in fact  symmetric, cf. \S \ref{anomexpcon} (ii),
 \end{itemize}
 {\em then}, the consequence \eqref{P1desid} is valid: 
 the period and $L$-sheaf match ``when tested against $\mathbf{e}$,''
 i.e., $$\Hom(\mathbf{e}', \textrm{period sheaf}_{X}) \simeq \Hom(\mathbf{0}, \textrm{$L$-sheaf}_{\check{X}})$$
 are Frobenius-equivariantly isomorphic. 
  
 While this seems very conditional, we regard it  rather as a significant check of the self-consistency
 of our overall pictures, in particular in relation to Tate and cohomological shifts. 
It is also plausible that by very similar methods one could prove
 that the local conjecture implies the global conjecture in the case of $\mathbb{P}^1$,
 by testing against arbitrary Hecke translates of $\mathbf{e}'$ and not just $\mathbf{e}'$ itself. (See~\cite[\S 3]{darioBM} where similar calculations are carried out in the spherical Hecke category itself.)

%% file: global-numerical.tex
\section{Numerical conjecture}
\label{L2numerical}

 Here we shall explicitly formulate, as conjectures in their own right,
 the numerical consequences suggested by the global conjectures,
 and     compare them to known  
  statements in the theory of automorphic forms.
  
   These conjectures avoid the
  various technicalities of derived geometry that we have encountered in the previous section.
  They are, to some extent, consequences of the global geometric conjecture but we prefer to regard them as free-standing
  statements with somewhat different ranges of applicability.
  For further discussion of this point, see
   \S \ref{mumerical}.
   
  We restrict to the case of {\em everywhere unramified automorphic forms over a function field,}
  i.e., eigenfunctions of all unramified Hecke operators.
 The reader familiar with automorphic forms will be disappointed with this restriction; but the
 general picture here is already sufficiently complicated to suggest that it would
 be foolhardy to go beyond this at an early stage, and in fact we believe that 
even here the story offers several interesting features that have not
 been properly explored (e.g.\ \S \ref{nontempered}, \S \ref{realperiods}, \S \ref{starperiods}) 
 in the classical theory.
 Of course, it is a fundamental question 
  to develop the conjectures in greater generality, which should go hand-in-hand with a deeper development
  of the ramified local conjecture. 
  
  The contents are as follows:
  \begin{itemize}
  \item \S \ref{Lconven} sets up general notation.
  \item \S \ref{conjtempered} gives the numerical conjecture in the tempered case
  (the most novel aspect of this for number theorists is the version about the $*$ period) and
 \item \S \ref{nontempered} discusses the conjecture in the nontempered case. 
  \item \S \ref{realperiods} discusses the questions related to whether periods are real-valued. 
  \item \S \ref{tempex} studies examples of the tempered conjecture (the main interest
  here is to make sure the constants are right) and \S \ref{NTexample} studies examples of the nontempered conjecture. 
  \item \S \ref{mumerical} discusses the relationship between the geometric and numerical conjecture. 
  \item \S \ref{starperiods} gives an introduction to the question of {\em star periods}. 
  The star period $P_X^*$  is an unfamiliar object
  in the classical theory of automorphic forms, but has been studied
  in the geometric Langlands context by Drinfeld, Gaitsgory, Schieder and Wang,
  and it is likely it plays an interesting role in the numerical theory also.
  \item \S \ref{Arthur functoriality} discusses the role of Arthur functorality and suggests a geometric interpretation of the nontempered conjecture. 
   \end{itemize}

Throughout this section we write:
$$b_G =(g-1) \dim G$$
for the dimension of the smooth stack $\Bun_G$.

  \subsection{Some conventions about $L$-functions} 
  \label{Lconven}

Our general notation will follow that of the previous sections, in the 
finite context.  In particular we have a finite field $\FF_q$ and
take for coefficient field 
  the algebraic closure $\kk$ of an $\ell$-adic field. 
We work with a  projective smooth curve $\Sigma$ over   $\FF_q$;
we denote by $F$ the function field of $\Sigma$, and use other notation as in \S \ref{section:global-geometric}
and \S \ref{Lsheaf}. 

To handle issues of rationality, 
we will work  with split forms of  hyperspherical dual pairs, a still somewhat tentative notion 
  postulated and discussed in  \S~\ref{GMdesiderata}; 
  a working definition of a class of such split forms is specified in  and after
Definition \ref{dhpFqdef}.
However, just as in \S \ref{GGC}, there is an alternate way to handle issues of rationality which avoids this notion: 
see  Remark   \ref{rationality conjecture}.

We suppose that the $G$-side admits a twisted polarization over $\FF_q$:
$$(G \times \GGm, M = T^*(X, \Psi))_{/\FF_q} \mbox{ and } (\check{G} \times \GGm, \check{M})_{/\kk},$$
and, as elsewhere in the global part of this paper, we assume that $X$ admits an eigenmeasure;  
when $\check{M}$ is polarized, we will also implicitly assume that this polarization, too, has an eigenmeasure 
(i.e., a top volume form which is preserved up to scalars by the group action). 
  We will be writing $\check G$ for $\check G(k)$.
  
  We will also fix an isomorphism $k \simeq \C$ and use it
  to transfer results to the complex numbers without comment;
  it also fixes a ``positive'' square root of $q$ in $k$, which gives us a root
  $\varpi^{1/2}$ of the cyclotomic character.   We  will therefore apply language from complex coefficients
 to automorphic forms valued in $\kk$, e.g., we will speak of an automorphic representation
 $\pi$ being ``tempered,'' which we will understand
 to mean that it is so when transported by any isomorphism $\kk \simeq \mathbb{C}$. 
\footnote{An eigenform is ``tempered'' when the associated
automorphic representation is tempered, that is to say, when its matrix coefficients
lie in $L^{2+\epsilon}$ modulo center.}

Recall (\S \ref{Ldef}) the notion 
of extended Langlands parameter, which we will here understand to be a  Frobenius-semisimple morphism  $$\phi_E = \phi_L \times \varpi^\frac{1}{2}:  \Gamma_F \longrightarrow  \check{G} \times \GGm$$ from the Weil group $\Gamma=\Gamma_F$
to an extended version of the dual group.
  For such $\phi_E$ we may consider the scheme of fixed points
$$\check{X}^{\phi} := \{x \in \check X: \phi_E(\gamma) x = x\}, \check M^{\phi}:= \{x \in \check M: \phi_E(\gamma) x =  x.\}.$$
considered as \emph{classical} schemes. (Hence, our convention here is different from the one of Section \ref{Lsheaf}, where we used similar notation to denote derived schemes of fixed points.)  Note that by default the superscript $\phi$ {\em always refers to fixed points
of the extended Langlands parameter.}

To each fixed point $x$ of $\phi_E$ on $\check M$ or $\check{X}$, let $T_x \check{X}$ or $T_x \check{M}$
be the tangent space to $\check M$  or $\check{X}$ at $x$; for our current discussion
we abridge both cases to $T_x$.  It carries a representation $(\phi_{x,E},T_x)$ of $\Gamma$
obtained by linearizing the $\phi_E$-action of $\Gamma$ at $x$, 
and we thereby obtain by differentiation  a homomorphism  
$$ \phi_{x,E}:  \Gamma \longrightarrow \GL(T_x)$$
which has an associated    $L$-function \index{$L$-functions with $\shear$}
which we shall denote by $L(s, T_x^{\shear})$
\begin{equation} \label{Lsshearlocaldef} L(s, T_x^{\shear}) := \mbox{$L$-function  for the $\phi_{x,E}$-action on $T_x$};
 \end{equation}
we similarly define the normalized version $L^{\norm}(s, T_x^{\shear})$, according to the general conventions introduced in
\S \ref{Leps}. Note that the parameter $\phi$ will be implicit in our notation for the $L$-value, except when necessary.

We use the same notation if the role of $T_x$ is replaced here by an arbitrary $\check{G} \times \GGm$ representation:
If $\rho: \check{G}  \times \GGm \rightarrow \mathrm{GL}(W)$ then 
\begin{equation} \label{shearedLfunctiondirect} L(s, W^{\shear}) := \mbox{$L$-function for $\Gamma_F$-action on $W$ via $\rho \circ \phi_{E}$.}\end{equation}

This usage of the $\shear$ notation is compatible with that introduced in \eqref{Lsheardef} and \eqref{Lnormdef}
and reminds us that we are dealing with an extended parameter, or, equivalently, the $L$-functions above  implicitly include shifts related to the $\GGm$-action.  
Indeed, suppose that $x$ is isolated, thus fixed by $\GGm$;
then $x$ is fixed by $\phi_L$ and  its tangent space is graded by the $\GGm$-action. Writing $\phi_{x,L}^{(i)}$
 for the induced representation on the $i$-th graded piece of the tangent space, we get
 \begin{equation}
 \label{Lxdecomp}
 L(s,T_x^{\shear})=  \prod_i L(\phi_{x,L}^{(i)}, s+\frac{i}{2}), \ \ \phi_{x,E} \simeq \bigoplus_{i}  \phi_{x,L}^{(i)} \otimes \varpi^{i/2} . \end{equation}

In the case of $\check M$, continuing to suppose that $x$ is isolated, we have a symplectic form on $T_x$ that 
 pairs the $i$- and $(2-i)$-eigenspaces, and in particular one gets (by the definition \eqref{Lnormdef} and the functional equation)
\begin{equation} \label{dropnorm} L^{\norm}(s, T_x M^{\shear}) = L^{\norm}(-s, T_x M^{\shear}), \,\,\,
 L^{\norm}(0, T_x M^{\shear}) = L(0, T_x M^{\shear}).\end{equation}
See Example \ref{exponents} \eqref{exponents-three} below for explication of the last equality.

 To help decipher  these hieroglyphs we include:
 \begin{example} \label{exponents}

 \begin{itemize}
 \item[(1)] Take $G=\GL_n$ and let $W$ be the standard representation
 with $\GGm$-action given by $\lambda \mapsto \lambda^n$. 
 Then for $\phi_L: \Gamma_F \rightarrow \check{G}$ we have
 $$ L(s, W^{\shear}) = L(s+\frac{n}{2}, \phi),$$
 $$ L(s, (W \oplus W^*\langle 2 \rangle)^{\shear}) = L(s+\frac{n}{2},\phi) L(s+1-\frac{n}{2}, \check{\phi}).$$
 where $W^* \langle 2 \rangle$ here just means that we modify the $\GGm$-action on $W^*$
by squaring, i.e., $W^* \langle 2 \rangle = W^*$ with $\GGm$-action given by $\lambda \mapsto \lambda^{2-n}$.   
 
 \item[(2)]
 The following example looks a bit arbitrary at the moment but will come up multiple times. 
Let \[ \phi_L = e^{ 2\rho} \circ \varpi^\frac{1}{2}, \]
 let $\check{G}$ act on its Lie algebra $\check{\mathfrak{g}}$ by the (left) adjoint action,
 and let $\GGm$ act on it by squared scaling. Then, writing $e$ for a regular  nilpotent
 compatible with $\rho$, its centralizer $\check{\mathfrak{g}}_e$ becomes a Galois module, and we have 
\begin{equation} \label{bhs} L(s, \check{\mathfrak{g}}_e^{\shear}) = \prod \zeta(s+d_i),\end{equation}
 where the $d_i$'s are the {\em exponents} of $\check{G}$, i.e., \index{exponents} \index{$d_i=$ exponents}
 the ring of invariant polynomials for the $\check{G}$-action on its Lie algebra
 is generated by homogeneous polynomials of degree $d_i$. 
 
 \item[(3)] \label{exponents-three} Unravelling notation in \ref{dropnorm}: $T_x M$ is graded by the $\GGm$-action
 as $T_x M = \oplus T_x M^{(i)}$, and
 then 
 \begin{multline*} L^{\norm}(s, T_x M^{\shear}) = \prod_{i} L^{\norm}(s+\frac{i}{2}, T_xM^{(i)}) = \\  \prod_{i} \sqrt{\epsilon}(s+\frac{i}{2}, T_x M^{(i)})\cdot L(s, T_x M^{\shear}).
 \end{multline*}
 When $s=0$, we can split the $\epsilon$-factors into pairs
 $\sqrt{\epsilon}(\frac{i}{2}, T_x M^{(i)}) \sqrt{\epsilon}(1-\frac{i}{2}, T_x M^{(2-i)})$
 and $\sqrt{\epsilon}(\frac{1}{2}, T_x M^{(1)})$,  both of which are identically equal to $1$ --
 the first by duality, between $T_x M^{(i)}$ and $T_x M^{(2-i)}$,  and the second because,
 being symplectic, $T_x M^{(1)}$ has trivial determinant, see \S~\ref{epsilon section}.
  \end{itemize}
 
\end{example}

 \subsection{The conjecture in the tempered case} \label{conjtempered}
   The conjecture comes in several variants, which apply in overlapping
   but slightly different situations. 
  We will formulate below all three forms for tempered Langlands parameters. The nontempered conjecture
 will be discussed in \S \ref{nontempered}.

 As elsewhere in this paper,
we write $b_G =(g-1) \dim G$ for the dimension of $\Bun_G$, and are only considering everywhere unramified automorphic forms and representations.
 We suppose that $M=T^*(X, \Psi)$ is (possibly, twisted-)polarized. Let
$P_X, P_X^*, P_X^{\norm}, P_X^{*\norm}$
be as previously defined (see \eqref{PXform1}, \eqref{PXform2}, Remark \ref{starperiod} and after \eqref{PXnormdef}) and put
for $f$ an automorphic function (i.e., a $\kk$-valued function on $\Bun_G$)  \begin{equation} \label{!def} P_X(f) = \int_{[G]} P_X(g) \cdot f(g) \, dg = \sum_{x \in \Bun_G(\FF_q)} \frac{1}{\# \mathrm{Aut}_x} P_X(x) f(x).
\end{equation} 

Thus, e.g., $P_X(f)$
in the case of $X=H\backslash G$ is the sum of $f$ over $H$-bundles weighted by inverse automorphisms as $H$-bundle.  

Note that, for the equality above to be true, we have fixed the measure \index{$\int_{[G]}$} on $[G]= G(F)\backslash G(\mathbb A)$ such that the preimage of a $G$-bundle under $[G] \to \Bun_G(\FF_q)$ has measure
equal to the inverse of the number of its automorphisms. We fix this measure on $[G]$ throughout, unless where stated otherwise.

The sum occurring in \eqref{!def} may not
be, in general, convergent, if $f$ is not compactly supported, e.g., a cuspidal automorphic form when $G$ is semisimple. In the divergent case, when $f$ is an automorphic form, there are often standard techniques to regularize those integrals, as integrals of asymptotically finite functions, see \cite[\S 5--6]{SaStacks}. The possibility to extend those integrals, as $G(\mathbb A)$-invariant functionals, from compactly supported functions to a space of asymptotically finite functions depends on the \emph{exponents} of such functions, i.e., of $P_X$ and of $f$; we will regard the integral as undefined otherwise.

The conjecture that follows relies on the Langlands parametrization of tempered automorphic representations.  Recall that tempered automorphic representations are irreducible summands of representations unitarily induced from cuspidal tempered representations, and therefore their Langlands parameters are provided by the work of V.\ Lafforgue \cite{LafforgueChtoucas, LafforgueICM}.

 \begin{conjecture} \label{numconj} (Global conjecture, tempered case). 
Suppose that $\pi$ is an everywhere unramified, tempered automorphic representation with   Langlands parameter $\phi$ (in particular, $\pi$ has unitary central character). 
 Then we may choose   a spherical vector $f=f_{\phi} \in \pi$ in such a way that  
$f^d=\overline{f}$ (where $d$ is the duality involution, as in \eqref{dualityinvolution};
and $\overline{f}$ refers to complex conjugation, transported to $\kk \simeq \C$),
 and  moreover the following properties  
hold with reference to any distinguished split form of a dual hyperspherical  pair $(M = T^*(X, \Psi), \check{M})$ as in \S \ref{usf}): 
 
 \begin{itemize}
 
 \item[(i)]
 Suppose that $\check{M} =  T^* \check{X}$ is polarized, without any twisting,
 and that the fixed points of the extended Langlands parameter $\phi_E$ on  $\check{M}$ is a finite (reduced) set of points; let  $\{x_1, \dots, x_r\} \in \check X$ be the 
 corresponding fixed set on $\check{X}$. 
 \footnote{Although the statement below refers only to fixed points on $\check{X}$, discreteness of
of fixed set on $\check M$, rather than merely $\check X$, is needed to avoid poles of the $L$-function. }
Then, for $f$ cuspidal, 
\begin{equation}   \label{PXshriekprediction}   P_X^{\norm}(f) \stackrel{?}{=}  q^{-b_G/2} \sum_{i} L^{\norm}(0, (T_{x_i} \check X)^{\shear})\end{equation}

Without cuspidality conditions on $f$,  the normalized star-period is given by 
\begin{equation} \label{PXstarprediction} P_X^{* \norm}(f) \stackrel{?}{=}  
(-1)^s  q^{-b_G/2} \sum_i L^{\norm}(1,\phi^d, (T_{x_i}  \check X)^{\vee \shear})\end{equation}
where $s$ is the dimension of the centralizer of the Langlands parameter of $f$,
the $x_i$ are now the fixed points for the {\em dualized} parameter $\phi^d$,
and the $T_{x_i}$ are considered as $\Gamma_F$-modules through $\phi^d$.
(See also \eqref{PXstarpredictionunshear} below.)

\item[(ii)]   Without assumptions on $\check{M}$, but assuming
that $f$ is cuspidal
  \begin{equation}\label{L^2 normalized equation}    P_X^{\norm}(f)   \stackrel{?}{=} q^{-b_G/2} \sum_{i} \sqrt{  L^{\norm}}(0, (T_{m_i} \check M)^{\shear}) \end{equation} 
  the sum now being taken over fixed points  $\{m_1, \dots, m_r\}$ of $\phi_E$ on $\check{M}$, again assumed finite. Here, 
  $\sqrt{  L^{\norm}}(0, (T_{m_i} \check M)^{\shear})$ refers to some square root of the quantity
 $L(0,T_m M^{\shear})$; note that these latter are real-valued and positive\footnote{This assertion is assuming -- as is expected -- the purity of the Langlands parameter; see \S \ref{realperiods}.}.
Moreover, this choice is   invariant by the action of the centralizer  $Z(\phi) \in \check{G}$ of the Langlands parameter,
which permutes the fixed points.

 \end{itemize}

\end{conjecture}

\begin{remark}
A few comments on the statements:
\begin{itemize}

\item[(a)]    As we will see, in each individual case where the $X$-period has been previously analyzed in number theory,   the statements
   (at least the ones about the usual period, rather than its $*$-counterpart) boil down to known theorems or conjectures. 
  However, even restricted to these cases, there is some value in formulating it as above: 
It illustrates how these known theorems or conjectures
 fit into the uniform duality formalism of this paper, in particular it is far from obvious
 that the various constants in the individual examples admit a uniform description. 

\item[(b)] To help process the $*$-period statement, it is helpful to restrict to the 
(very common) case that all fixed points of $\phi_{E}$ on $\check{X}$ are in fact $\check{G}$-fixed points.
In this case, the statement becomes (see Remark \ref{unshearremark}): 
\begin{equation} \label{PXstarpredictionunshear} 
P_X^{* \norm}(f) \stackrel{?}{=}  
(-1)^s  q^{-b_G/2} \sum_i L^{\norm}(1, T_{x_i}  \check X^\unshear)\end{equation}
where $\unshear$ denotes shearing for the {\em inverted} $\GGm$-action. 
This now looks more like \eqref{PXshriekprediction}, but observe that:
\begin{itemize}
\item[(i)]   the $L$-value has been
shifted to $1$;
\item[(ii)] there is an appearance of the interesting sign $(-1)^s$
\item[(iii)] The $\GGm$ shift has been inverted.
\end{itemize}

\item[(c)] 
Regarding the condition $f^d =\bar{f}$: At least if we suppose multiplicity one, in that the line through $f$ is uniquely specified
by the Hecke eigenvalues, we have 
  $f^d = t \overline{f}$ for some $t \in \mathbf{C}^{\times}$;
and since $t \bar{t} = 1$ we can  modify $f$ so that $t=1$. The resulting choice is determined up to a real scalar.
Moreover, in this situation, the line of $f$ is defined over the field generated by Hecke eigenvalues; being a CM field,
the choice of isomorphism $\kk \simeq \mathbb{C}$  does
not affect the validity of $f^d = \bar{f}$. 

\item[(d)] The normalization of $f_{\phi}$ can be compared with standard ones by using $X = \mbox{Whittaker}$,
or (up to sign) by using $X$ the group case, see Example \ref{L2Whitnorm}. 
We have preferred however not to single these cases out, regarding them as   special (if particularly useful) cases of the general duality phenomenon. 

\item[(e)]   It is a straightforward matter
 to pass between this conjecture and a corresponding one for unnormalized periods,
 using the same type of discussion as   \eqref{shiftpred}. The unnormalized period
 then involves various powers of $q$ depending on $\beta_X, \beta_{\check{X}}$. 
 
\item[(f)] Our discreteness assumptions in the theorem statement mean that the $L$-functions appearing on the right hand side 
are always defined, i.e., they are never evaluated at a pole point.   \footnote{Indeed,
the quantity
$ L^{\norm}(0, (T_{x_i} \check X)^{\shear})$ 
can be evaluted by means of \eqref{Lxdecomp};
it has a pole precisely  when 
one of the $\phi^{(i)}_{x,L} \otimes \varpi^{i/2}$
contains a copy of the trivial representation or the cyclotomic character $\varpi$;
assuming, as one expects, each $\phi_{x,L}^{(i)}$ to be pure of weight zero,
this can happen only for there is a $\phi_{x,L}$-fixed vector
in $T_{x_i} \check{X}$
lying in $\GGm$ degree $i=0$ or $i=2$.  In the former case,
there is a $\phi_E$-fixed tangent vector
at $x_i$, contradicting our supposed reducedness; in the latter case,
the same happens for the ``vertical'' tangent space for $T^* \check{X}$
above $x_i$. }
 
Correspondingly, 
we would expect that the left hand side can always be assigned an unambiguous regularization, cf. discussion after \eqref{!def}.

\item[(g)]  The phenomenon of obtaining a sum of $L$-functions is not common but has been observed. As we will see in the discussion of the group case below, it is sometimes ``hidden'' in sizes of centralizers of Langlands parameters that appear in conjectures about periods.

One example is the 
``Eisenstein case'' of $X=G/U$,  which
we discussed from the geometric viewpoint in \S \ref{Eisgeom}. It does not fit into the hyperspherical umbrella but has many features in common with it.  The relevant period computation is the constant term of Eisenstein series, which involves a sum of $L$-functions indexed by the Weyl group;  see \S \ref{eis period numerical comp} for an examination of how this fits with the conjecture.

A more interesting example
is the appearance of centralizer groups in period formulas, as in Example \ref{L2Whitnorm} below; and 
a yet more interesting example was given
in the case of (a form of) $X=\mathrm{GL}_2 \backslash \mathrm{SO}_5$ in the PhD thesis
of X.\ Wan \cite{WanThesis}.
 
 \item[(i)] One case in which the integral  diverges for   trivial reasons
 is where a central torus in $G$ acts trivially on $X$.
 This is usually handled, in the automorphic literature, by integrating
 on $G$ modulo this central torus.  In the current case, this ``trivial divergence''
 matches with a similar trivial divergence on the dual side, and can 
be handled by passing to an isogenous situation.
 \end{itemize}
\end{remark}

 \begin{example} \label{L2Whitnorm}
 The conjecture  entails the normalizations
 \begin{equation} \label{Whit} \mbox{normalized Whittaker period of $f$} = q^{-b_G/2} \end{equation}
 (unnormalized Whittaker period $=q^{\frac{(g-1) [-\langle 2 \rho, 2 \rho^{\vee} \rangle +\dim U - \dim G]}{2}} $),
 \begin{equation} \label{l2} \int_{[G]}|f(x)|^2 \, dx = (\# Z_{\phi}) L(1, \check{\mathfrak{g}})
\end{equation} 
with $Z_{\phi}$  the centralizer of the Langlands parameter $\phi$ inside $\hat{G}$.

 To see the first statement, we apply the conjecture to the case of $X$ the Whittaker period;
 then $\check{X}$ is a point, and the result is obvious. For the unnormalized Whittaker period, see \eqref{Whittakernormalized}.

For the latter,  we apply the conjecture to the case of the group period
 $X=G$ as a $G \times G$-space, where we twist
 the second factor of $G$ to act through the duality involution (we expect this to be the distinguished split form in general, 
 see Example \ref{group period example}). The normalized period equals $q^{-b_G/2} \sum |f(x)|^2$ because $f^d = \overline{f}$. 
 
 The dual space $\check X$
 is the standard $\check{G}$ as $\check{G} \times \check{G}$-space. 
 The fixed points are then precisely $z \in \check{G}$ which centralize
 $\phi_L$, i.e., the centralizer.   The $L$-function appearing is orthogonal;
 recalling that the group is now $G \times G$, the power of $q$ on the right hand side of
 \eqref{PXshriekprediction} is $q^{-b_G}$, and finally 
the term $\epsilon^{-1/2}$ that appears in the normalized $L$-function \eqref{Lnormdef},
 after using the functional equation to switch evaluation point to $s=1$,  equals
$q^{b_G/2}$; thus \eqref{l2}.  

As is clear from this example, the sizes of centralizer groups appear naturally from fixed point counts.
The appearance of these sizes has been noted in period conjectures before;
see, in particular, the work of Ichino and Ikeda \cite{II}. \footnote{In the literature, one usually normalizes the ratio of local and global $L^2$-norms to be equal to $1$, and this results in the sizes of these centralizers showing up in the formulas for other periods (such as the Gross--Prasad and Whittaker periods). In our formulation of Conjecture \ref{numconj}, the eigenform $f$ is normalized so that the sizes of centralizers play a role in its $L^2$ norm, and do not appear in Whittaker or Gross--Prasad periods -- and this fits in more naturally with relative Langlands duality.}
\end{example}

\begin{remark}[Nonlinear $L$-functions]  
 
The right hand side of \eqref{PXshriekprediction} can be considered  as a ``nonlinear'' $L$-function attached to $\check{X}$.
It is a function on $\check{G}$-valued Galois representations  -- satisfying a mild discreteness condition. In the case that $\check{X}$ is a vector space,  it recovers the usual $L$-function, evaluated at a point determined by the $\GGm$-action.
See \cite{ChenVenkatesh}  for a more explicit development of this viewpoint in some cases.
  \end{remark}

\begin{remark}[The case of number fields and ramified representations] \label{remark-global-ramified}
	Although in this paper we work over function fields and everywhere unramified representations, the numerical versions of the conjectures are well-suited for generalization to number fields and the ramified setting.
  At present, however, we cannot formulate them with the same level of precision as the conjectures of this paper. 
  
The  primary difficulty  is that we cannot, in general, pick a single theta series and a single vector $f$ in the space of an automorphic representation $\pi$ for which nice numerical formulas for various periods will hold, but we rather have to reformulate the equalities of Conjecture \ref{numconj} as equalities of two pairings between the Schwartz space $\mathcal S(X(\mathbb A))$ of the adelic points of the spherical variety, and the space of an automorphic representation $\pi$. The sums on the right hand side of \eqref{PXshriekprediction}, \eqref{PXstarprediction}, \eqref{L^2 normalized equation} will be over fixed points of the hypothetical (extended) global Langlands parameter of $\pi$, and one can replace the $L$-functions that appear by partial $L$-functions (away from a finite set of places that contain all ramified and archimedean places), multiplied by ``local zeta integrals.'' The latter would be local pairings of the same form as the global period pairings, and while the literature is abundant with examples of those, we do not, in general, know how to encode them into a general recipe, except in the case of \eqref{L^2 normalized equation}, where, under an additional ``multiplicity one'' assumption, such functionals were described in \cite[\S~17]{SV}, generalizing the local Ichino--Ikeda periods \cite{II}. (Without this multiplicity-free assumption, the paper \cite{FLO} suggests that the local multiplicity space should have a basis that is related to the fixed points of the local Langlands parameter on $\check M$; this basis would provide the local Euler factors in these conjectures.)
\end{remark}

\subsection{Nontempered representations} \label{nontempered}

We shall analyze how to modify part (iii) of our prior conjecture in nontempered cases (see also \S \ref{Arthur functoriality} for a more geometric perspective). 

 \index{Arthur parameter}

\subsubsection{Arthur parameters} \label{apsetup}
We first set up notation for later use.  Assuming the Arthur conjectures \cite{Arthur-unipotent-conjectures} on the parametrization of the automorphic discrete spectrum, let $f_{\phi}$\footnote{
 Although
elsewhere in this section we use $f$ for an automorphic form, we
will try to consistently use $f_{\phi}$ here to avoid any confusion with the $f$ of the 
$\mathfrak{sl}_2$-triple.
}
be an $L^2$ Hecke eigenform
attached to the Arthur parameter
  $$ \phi_A:  \Gamma_F \times \SL_2 \rightarrow \check G(k).$$   We will assume that $\phi_A|\Gamma_F$ is pure of weight zero, i.e.
  all eigenvalues of Frobenius elements, taken after any fixed embedding $\check{G} \rightarrow \GL_N$,
  have absolute value $1$ after transport to $\CC$.\footnote{ Note that, in \cite{Arthur-unipotent-conjectures}, Arthur assumes the existence of the hypothetical Langlands group $L_F$
 and works with $\C$ coefficients. Above we have allowed ourself, in line with our general setup, to work with $k$ coefficients.
Arthur assumes that the image of $L_F$ is bounded; with $k$ coefficients  the purity
condition is a reasonable substitute.}
 We write $(h, e, f) \in \check{\mathfrak{g}}$ for the $\mathfrak{sl}_2$-triple.
In particular, for $f_\phi$ to be $L^2$, the parameter is \emph{discrete}, that is, its centralizer $Z_{\phi}$ of $\phi_A$ is finite.

  Associated to $\phi_A$ are a Langlands parameter $\phi_L$, 
which
should reproduce the Hecke eigenvalues of $f_{\phi}$,  and an extended Langlands parameter $\phi_E$: 
\begin{equation} \label{phiLdef}  \phi_L = 
 \phi_A \circ  \left( \mathrm{id} \times \left[\begin{array}{cc} \varpi^{1/2} & 0 \\ 0 & \varpi^{-1/2} \end{array} \right]\right):   \Gamma_F \rightarrow \check G(k),\end{equation}
\begin{equation} \label{phiEdef}  \phi_E := \phi_L \times \varpi^{1/2} : \Gamma_F \rightarrow \check{G} \times \GGm(k),\end{equation}
where, as usual, $\varpi$ is cyclotomic
and $\varpi^{1/2}$ the square root defined using the chosen square root $\sqrt{q}$ (see \eqref{cyclotomicdef}). 
 Let us also write  $\iota$ for the restriction of $\phi_A$ to $\SL_2$
  and $a=\iota|_{\Gm}$ for the associated cocharacter $a: \Gm \rightarrow \check{G}$
  given by
  \begin{equation} \label{adef}
  \lambda \mapsto \iota \left[ \begin{array}{cc} \lambda & 0 \\ 0 & \lambda^{-1} \end{array} \right]. 
  \end{equation} 
In other words, $a$ is the composition of $\iota$ with the cocharacter denoted by $\varpi$ in \S~\ref{sl2pair}, but we refrain from using $\varpi$ for that purpose here, since it has been reserved for the cyclotomic character.

By definition, the restriction of $\phi_A$ to $\Gamma_F$ is pure of weight zero and this implies the following useful fact:

\begin{lemma}\label{apstar}  
Assume that $\phi_A|\Gamma_F$ is pure of weight zero. Then  
the Zariski closure of $\phi_E$ contains
the group\footnote{In \eqref{Gmprime} we defined $\GGm'$ by the inverse cocharacter into $\check G$; the inversion is due to the switch from right to left actions, see \S~\ref{leftrightconventions}. The notation is compatible if we consider $\GGm'$ as a copy of $\Gm$ with a map to the automorphism group of $\check M$, and it is only through its action on $\check M$ that this group plays a role here.} 
 \begin{equation}\label{GGmprime} \GGm':= \mbox{image of }(a, \mathrm{id}):\Gm \rightarrow \check{G} \times \GGm. 
 \end{equation}
\end{lemma}
 
 \begin{proof}
Consider 
\begin{equation} \label{morgan} \phi_A|_{\Gamma_F} \times \varpi^{1/2}:
\Gamma_F \rightarrow \check{G} \times \GGm\end{equation}
and let $\Lambda$ be its image. 
Suppose that the Zariski closure of $\Lambda$ did not contain $\GGm$.
Then this Zariski closure intersects $\GGm$ in a finite subgroup, say $\mu_N \subset \GGm$, and then the Zariski closure itself 
 must be contained in a subgroup of the form $\{(q, \lambda): \lambda^N = \chi(q)\}$
where $Q \leqslant \check{G}$ is algebraic and $\chi: Q \rightarrow \GGm$ a character.
This contradicts purity, which asserts that all
  Frobenius eigenvalues, for some fixed embedding $Q \subset \GL_N$
and after transport via $k \simeq \C$, all have absolute value $1$. 
Indeed, fixing a maximal torus $T_Q \subset Q$,
we can write $\chi|_{T_Q}$
as a linear combination of various characters occurring
in the embededing $T_Q \hookrightarrow \GL_N$,
and it follows from this that $\chi$ at any (semisimplified) Frobenius element
also has absolute value $1$ in $\C$, which is false -- $\chi$ on this Frobenius coincides with 
  the $N/2$th power of the cyclotomic character.
\end{proof}

\subsubsection{}

 The overall slogan is that to compute the period
 for an Arthur form $f_{\phi}$ one uses not $\check{M}$
 but a certain ``Kostant'' or ``Slodowy'' slice of $\check{M}$;
 that is to say, rather than considering $\Gamma_F$-fixed points on $\check{M}$
 we consider $\Gamma_F$-fixed points on
 the (smooth) slice of $\check{M}$ given by
\begin{equation} \label{checkMslice} \check M_{\slice} := \mbox{preimage of $f+\check{\mathfrak g}_e$ under the moment map $\mu: \check{M} \rightarrow \check{\mathfrak{g}}^*$},\end{equation}
 where $e,f$ is the $\mathfrak{sl}_2$-triple associated to the Arthur parameter;
 and the action is through $\phi_E$. 
 Here, as in \S~\ref{Slodowy}, we write $f+\check{\mathfrak g}_e$ for the affine subspace of $\check{\mathfrak g}^*$ that would more canonically be denoted by $f+\check{\mathfrak g}^{*,e}$.
 
 \begin{remark}  \label{msliceform}
 $\check{M}_{\slice}$ is not $\GGm$-stable,
 because $\GGm$ doesn't preserve $f+\check{\mathfrak g}_e$.
 However, it {\em is} stable under the action of $\Gamma_F$
 through the extended Langlands parameter $\phi_E$  associated to $\phi_A$ (see \eqref{phiEdef}). 
This follows from the fact that, with our previous notation, $(a, \mathrm{id}): \GGm \rightarrow \check{G} \times \GGm$
preserves $f+\check{\mathfrak g}_e$ under the adjoint action -- this is the action of the group $\GGm'$ from \eqref{Gmprime}, taking into account the switch from right to left actions (\S~\ref{leftrightconventions}).
We also observe that the $\Gamma_F$-action on $\check{M}_{\slice}$
scales the symplectic form through the cyclotomic character,
i.e., $\gamma^* \omega = \varpi(\gamma) \omega$. 
  \end{remark}

  In \S\ref{Arthur functoriality}, we will suggest a heuristic reason
 for the appearance of this $\check{M}_{\slice}$ (see also Remark \ref{whynontempered} (b) for another point of view). By the theory of the Slodowy slice,
 the morphism $\check{G} \times (f+\check{\mathfrak g}_e) \rightarrow \check{\mathfrak{g}}^*$
 is smooth (see e.g.\ \cite[\S 2.2]{Gan-Ginzburg}) and from this we readily deduce that $\check{M}_{\slice}$ is smooth. 
 Indeed,  $\check M_{\slice}$ is a  {\em twisted Hamiltonian reduction}
 of $\check M$ by the unipotent subgroup determined by the $\mathfrak{sl}_2$-triple,
 and this is the point of view from which it will appear in our heuristic discussion \S \ref{Arthur functoriality}.

 Observe that the fixed points of $\Gamma_F$
 on $\check M_{\slice}$ are, by Lemma \ref{apstar}, fixed by  
 the action of $\GGm'$,  and in particular {\em all such fixed points map to $f \in \check{\mathfrak{g}}^*$. }

 \begin{conjecture} \label{numerical nontempered} 
 (Nontempered periods are obtained by Slodowy-slicing $\check{M}$:)  
Take a distinguished split form of a dual hyperspherical  pair $(M = T^*(X, \Psi), \check{M})$ as in \S \ref{usf}), with $M$ polarized.

   Let $f_{\phi}$ be an  everywhere unramified automorphic form belonging to the discrete series
 with Arthur parameter $\phi_A$, as above, and associated Langlands and extended Langlands
 parameters $\phi_L$ and $\phi_E$.  
 Assume that the fixed points of $\phi_E$ on  $\check M_{\slice}$ form a finite set $m_1, \dots, m_r$.
 Then for a normalization of $f_{\phi}$ with $f_{\phi}^d= \overline{f_{\phi}}$, independent of $\check{M}$,
 we have the equality:
   \begin{equation} \label{l2conjecturenontempered}   P_X^{\norm}(f_{\phi}) =  q^{-b_G/2} \sum_{i} \sqrt{ L^{\norm}}( 0, T_{m_i} \check M_{\slice}^{\shear})
 \end{equation}
 where the $L$-function on the right is to be interpreted as in \eqref{Lsshearlocaldef}, and the square root should be interpreted as in \eqref{L^2 normalized equation}, i.e., as a choice of signs invariant under the centralizer of $\phi_A$.  \end{conjecture}
 
 The role of the shear in the notation above was explained in \eqref{Lsshearlocaldef}; to reformulate, the space $T_m \check M_{\slice}^{\shear}$ is considered as a graded Galois representation via the restriction of the Arthur parameter $\phi_A$ to $\Gamma_F$, with the grading coming from the action of $\GGm'$ (see Lemma \ref{apstar}); and this action, composed with the square root $\varpi^\frac{1}{2}$ of the cyclotomic character, introduces shifts to the point $0$ of evaluation.

 \begin{remark}
 \label{whynontempered} ~
 \begin{itemize}
  \item[(a)]
As a parallel to \eqref{l2} we will compute in \S~\ref{L2case} that  
 \eqref{l2conjecturenontempered} implies in particular for the $L^2$ norm that
  \begin{equation} \label{l2nontempered} \int_{[G]} |f_{\phi}(x)|^2 \, dx   = (\# Z_{\phi_A})    L(0, T_f (f+ \check{\mathfrak{g}}_e)^\shear) =  (\# Z_{\phi_A})    L(1,  \check{\mathfrak{g}}_e^\shear),\end{equation}
 where $\check{\mathfrak g}_e$ is the centralizer of $e$ and the shear on it in the last expression comes only from the action of $\Gm\hookrightarrow\SL_2$. On the other hand, the shear on the tangent space $T_f (f+ \check{\mathfrak{g}}_e)$ comes from the action of $\GGm'$, which is the combination of the action of $\Gm\subset\SL_2$ and the usual (square) action of $\GGm$ on $\check{\mathfrak{g}}_e\subset \check{\mathfrak g}^*$. The most vivid example is when $f_{\phi}$ is the trivial function; then 
the normalization relevant for the conjecture is just $f_{\phi}=1$, and 
\eqref{l2nontempered} is giving the usual ``Tamagawa number''
formula for the volume of $\Bun_G$ -- the product of the order of the center of $\check{G}$,
giving the number of components, and a product of $\zeta$-functions.  See \eqref{bhs} or \S \ref{BZGHstuff} for that last deduction.

\item[(b)] Conjecture \ref{numerical nontempered}
 can be considered a ``regularization'' of the prior Conjecture \ref{numconj}.
 That is to say if we  naively take the case (ii) of Conjecture \ref{numconj}, 
and  apply it to the nontempered case,  we find that  the right hand side of \eqref{L^2 normalized equation} 
  diverges. Nonetheless (as is familiar from examples, such as those examined
  in \cite{II}, 
 or from the consideration of Eisenstein series),
 if we take the {\em ratio} between \eqref{L^2 normalized equation} and \eqref{l2}, 
 and  formally cancel divergent factors, the result is compatible
 with Conjecture \ref{numerical nontempered}. The algebraic
 explanation lies in the following diagram:
 Let $O = [\check{\mathfrak{g}}, e]$, so that 
  $O \subset \check{\mathfrak{g}}$
is a complement to $\check{\mathfrak{g}}_f$ inside $\check{\mathfrak{g}}$.  Then  with $o_x: \mathfrak{g} \rightarrow T_x$ the orbit map,
one has a commuting diagram of isomorphisms of vector spaces:

$$ \xymatrix{
 O \ar[d]^{\sim} & \oplus &  T_x \check M_{\slice} \ar[rrr]^{(X, Y) \mapsto o_x(X) +Y }  \ar[d]_{\sim}^{d\mu} &&& T_x \check M \ar[d]^{\mu}   \\
  \check{\mathfrak{g}}/\check{\mathfrak{g}}_f &  \oplus  &  Z(e) \ar[rrr]^{(X,Y) \mapsto [X,f] + Y} && &  \check{\mathfrak{g}}^*,
}
$$
 so the ratio  of \eqref{l2conjecturenontempered} and \eqref{l2nontempered}  
 gives the same result as the ratio of \eqref{L^2 normalized equation} and  \eqref{l2}.
  \end{itemize}
 \end{remark}

\subsection{Real structures on normalized periods} \label{realperiods}  
 We examine the $L$-functions
 on the right hand side of \eqref{L^2 normalized equation}.
 
 The only case of interest here is where the $\SL_2$
 associated to $\check{M}$ is trivial; otherwise, one expects
 the period of all tempered forms to vanish identically. 
 Assuming this is so, $\check{M}$ has the form
 $\check{G} \times^{\check{G}_X} (S \oplus  \check{\mathfrak{g}}_X^{\perp})$ by \eqref{fund_id}.  
 We can, without loss of generality, suppose
that a $\Gamma_F$-fixed point $m_i \in \check{M}$
lies above the identity coset of $\check{G}/\check{G}_X$,
and therefore that $\Gamma_F \rightarrow \check{G}$
factors through $\check{G}_X$. 
 The tangent space $T_{m_i} M$ decomposes,
 as a $\Gamma_F \times \GGm$-representation,
 as 
 $$T_{m_i} \check{M} = (\check{\mathfrak{g}}/\check{\mathfrak{g}}_X)_0 \oplus S_1 \oplus (\check{\mathfrak{g}}/\check{\mathfrak{g}}_X)_2$$ 
 where the subscripts denote the $\GGm$-weight. Note that  $(\check{\mathfrak{g}}/\check{\mathfrak{g}}_X)$
 is an orthogonal representation of $\Gamma_F$, $S_1$ is a symplectic representation, and 
   temperedness of the automorphic form (for any isomorphism $\kk \simeq \C$)
 implies that the representation of $\Gamma_F$ on $T_{m_i} \check{M}$ is pure of weight zero. 
 
With this notation, invoking the functional equation \eqref{LnormFE}, the $L$-value  $L^\norm(0, T_{m_i} \check M^{\shear})$ appearing in \eqref{L^2 normalized equation} 
can be rewritten as
$$L^{\norm}(1,  (\check{\mathfrak{g}}/\check{\mathfrak{g}}_X))^2 L^{\norm}(\frac{1}{2}, S_1).$$ 
This is real  (because each Euler factor is real) and non-negative (by the Riemann hypothesis). 
In particular, \eqref{L^2 normalized equation} entails that:
\begin{quote} $P_X^{\norm}(f)$   is {\em real-valued},
when $f$ is chosen so that $f^d = \bar{f}$.
\end{quote}
(since we have taken $f$ to be $\kk$-valued, this statement
should be interpreted as holding after transporting via $\kk \simeq \mathbb{C}$.)

Curiously, this statement is not obvious, and not discussed explicitly in the literature, at least in any generality. We think that it represents an interesting phenomenon, and 
 pause to discuss it here. 
 In what follows, we will often refer to $e^{\rho}(-1)$; it is 
 an element of the adjoint group of $G$ and its adjoint action defines an involution of $G$.

With $f$ normalized so that $f^d =\bar{f}$, we have 
$$ \overline{  P_X^{\norm}(f) } =  \int_{[G]} \overline{ P_X^{\norm}} \cdot f^d    = \int_{[G]} \overline{(P_X^{\norm})^d} \cdot f.$$
 and so $P_X(f)$ will be real-valued for such $f$ if
  \begin{equation} \label{key reality assertion} \overline{P_X^{\norm}} \stackrel{?}{=}  (P_X^{\norm})^d.\end{equation}

 Note that the right hand side is independent of the choice of isomorphism $k\simeq \mathbb C$, while the left hand side depends on it; in other words, this stament entails the assertion that $P_X^{\norm}+ (P_X^{\norm})^d$ is totally real.
 
\subsubsection{Some examples and corollaries}   
  We have not verified \eqref{key reality assertion} in general, but there are several cases where we can confirm it. We also discuss some consequences of the expectation.  
  For the discussion that follows, note that,
  by transport of structure, 
$(P_X)^d = P_{X^d}$, and 
the same holds for the normalized version.   
 
  \begin{enumerate}[(a)] 
 
 \item In the Whittaker case,  $\overline{P_X^d} = P_X$
by the very way $d$ is constructed. This can be seen as one reason to prefer
the duality involution over the pinned Chevalley involution.

 \item  In the case of ``distinguished split forms of hyperspherical varieties,'' in the sense of Definition \ref{dhpFqdef}, \eqref{key reality assertion} follows from Lemma \ref{dualitystableX} . Namely, 
 $(P_X^{\norm})^d$ and $\overline{P_X^{\norm}}$ should be considered, respectively,  as quantizations of the Hamiltonian space $M^d$ (obtained from $M$ by composing the $G$-action and the moment map with the duality involution, with the same $\GGm$-action), and the space $\overline{M}$ (obtained from $M$ by negating the symplectic form and moment map). 
 The independence of the normalized period function from polarization (Proposition \ref{prop:Pindep}) 
 shows that Lemma \ref{dualitystableX} 
 then implies
 \eqref{key reality assertion}.  

\item    (Relationship between the normalized vectors for $f$
and that for the dual representation:)  If $f$ is normalized so that $f^d = \bar{f}$ and $W(f)=q^{-b_G/2}$
for the normalized Whittaker period $W$,  as in \eqref{L2Whitnorm}, 
then $\bar{f}$ belongs to the dual (=conjugate) representation, satisfies 
the same relationship with its $d$-twist, and $W(\bar{f}) =  \pm q^{-b_G/2}$,
where the sign is chosen according to the action of $\Ad(e^{\rho}(-1))$
upon $f$. 

   In other words, the ``good'' normalization of a vector
in the dual representation is $\pm \bar{f}$, the sign being taken according
to the action of $e^{\rho}(-1)$ upon $f$.

To see this, note that the distinguished split form of the group period is  (presumably, see Example \ref{group period example}) $G \backslash (G \times G)$,
where the embedding is not the diagonal, but rather the graph of $\Ad(e^{\rho}(-1))$;
and correspondingly the period of $(f, \pm \bar{f})$ 
is given by $\int_{[G]} |f|^2$, which is real and positive as desired.  The distinguished split form of the Chevalley-twisted group period is now $G \backslash (G \times G)$,
where the embedding is via the duality involution, and correspondingly the period of $(f, f)$
is given by $\int f f^d = \int |f|^2$, again real and positive.

\item
  (The role of $\sqrt{-1} \in \GGm$). 
Above we have considered reality of the ``spectral transform'' of $P_X^{\norm}$;
but indeed  $P_X^{\norm}$ is often real-valued itself. 
 For example, in the Whittaker case,
this is so if $-1$ is a square in $\FF_q$ or if $\check{\rho}$ belongs to the cocharacter lattice.   The role of a square root of $-1$ here is of interest
and it seems to occur in several related ways.  We do not understand this at a deeper level,
but we observe that  the action  of $\sqrt{-1} \in \GGm$
 gives an equivalence between $M$ and the same space with negated symplectic form and negated moment map. 
 Since $P_X^{\norm}$ and $\overline{P_X^{\norm}}$ should be considered as quantizations
 of $M$ and the negated space, respectively, this reality is not unexpected. 
 
    \end{enumerate}

 \subsection{Tempered examples: Whittaker, Gross-Prasad, Eisenstein, Tate} \label{tempex}
 We shall briefly examine (i) or (iii) of Conjecture 
\ref{numconj} in various cases. As observed at the start of the section, in each case
the $L$-function formulae are familiar in number theory; the main point is to check that the constants are right. 
 
\subsubsection{The Iwasawa-Tate period}
Take $X= \mathbb{A}^1$ and $G=\mathbb{G}_m$. 
We have computed the normalized period in \eqref{PXTatenorm}. 
An automorphic form is simply an idele class character $\chi$. The relevant normalization
is taking   $f  = q^{-(g-1)/2} \chi$, as can be seen from \eqref{Whit}.   It is clear that $f^d =\overline{f}$, and  we get
$$ \langle P_X^{ \norm}, f \rangle = q^{-(g-1)/2} \chi(K^{1/2})^{-1} L(\frac{1}{2}, \chi)  = q^{-b_G/2} L^{\norm}(1/2, \chi),$$
as required.

 \subsubsection{The Whittaker case: compatibility with Lapid--Mao} 
   \label{Whitexample}    

Let $X$ be the Whittaker space with its twisted polarization.   We will sketch that the ratio of \eqref{l2} and \eqref{Whit} is as predicted by the work of Lapid and Mao.
   The main point here is to verify that the power of $q$ is correct.
   To avoid some minor issues with center, suppose $G$ is semisimple. We follow notation as in the conjecture,
   so that $f$ is an automorphic form with Langlands parameter $\phi$.

 Let $W(f)$ denote the unnormalized Whittaker period of $f$;  with reference to adelic uniformization:
 \[W(f) = \int_{n \in U(F) \backslash U(\adele)} \psi(u) f(ua_0) du, \,\,\,
 a_0 = e^{2 \check\rho}(\partial^{-1/2}),\]
where the measure on $U(\adele)$ is normalized to give mass one to the quotient.
Now, by   Example \ref{Whit explicit computation example}
 $$P_X^{\norm}(f) = q^{-\beta_X/2}  q^{\beta_X} W(f), \beta_X = (g-1) (\dim(U) - \langle 2 \rho, 2 \check\rho \rangle). $$
 and combining \eqref{Whit} and \eqref{l2} reads
   $$ \frac{|W(f)|^2}{  \int_{[G]} |f|^2} =  \frac{q^{\beta_X} q^{-b_G}}{(\# Z_{\phi}) L(1, \check{\mathfrak g})}
   = 
    \frac{q^{-(g-1) u + (g-1) \langle 2 \rho,2  \check\rho\rangle }}
    { \# Z_{\phi} q^{b_G} L(1, \check{\mathfrak g})}  $$ 
and so our conjecture implies
\[ \frac{|W(f)|^2}{ \int_{[G]} |f|^2 }  =
  \frac{q^{-(g-1) u + (g-1) \langle 2 \rho,2 \check\rho  \rangle}}{ \# Z_{\phi}   q^{b_G} L(1, \check{\mathfrak g})} ,\]
  where the measure on $[G]$ is normalized to give the maximal compact
  of $G(\mathbb{A})$ the volume $1$.

Let us now explain why this conjecture is compatible with that of Lapid and Mao \cite{LM}
(they also prove many cases). Their assertion,
specialized to the unramified case, says that 
the left-hand side above equals 
\[ ([Z_{\phi}: Z(\hat{G})])^{-1} \vol([G])^{-1}
\int_{n\in N(\mathbb{A})} \langle n a_0 f, a_0 f \rangle \psi(n)\, dn,\]
where  the integral at the end is to be expressed as a product
and regularized in a standard way.
Almost all the local factors in this product coincide with 
the local factors of 
 $
\frac{\prod \zeta_v(d_i)}{L_v(1, \mathrm{Ad})},$
where $d_i$ are the exponents of $G$, cf. \S \ref{exponents}. 
The measure  is chosen so that  $dn = q^{-(g-1) u} \prod_{v} dn_v$
where each $dn_v$ assigns mass $1$ to the integral points of $N$.
Note that in fact $\psi(n) = \psi_0(a_0^{-1} n a_0)$
where $\psi_0$ is everywhere ``unramified,'' and $a_0$
is as in \eqref{a0def}. 
Making the substitution $n \leftarrow a_0 n a_0^{-1}$, 
 so that $d(a_0n a_0^{-1}) =|e^{2\rho(a_0)}| dn$, we compute the integral above to equal
 \[q^{-(g-1) u}  |e^{2\rho(a_0)}| \int \langle \pi(n) v,   v \rangle \psi_0(n)\, dn = q^{-(g-1) u + (2g-2) \langle 2 \rho, \check{\rho} \rangle} \frac{\prod \zeta(d_i)}{L(1, \mathrm{ad})},\]
 and now taking into account $\vol([G]) =  \# Z(\hat{G}) q^{b_G} \prod \zeta(d_i)$
 gives the claimed equivalence between our conjecture and that of Lapid and Mao.

\subsubsection{The reductive case and the Ichino-Ikeda conjecture}
Again, to avoid minor issues with the center, suppose that $G$ is semisimple.
  Assume now that $X=H\backslash G$ and 
  that $\check{M}$ is defined by 
  data $\check{G}_X \subset \check{G}$
  and a symplectic $\check{G}_X$-representation $S_X$,
  with {\em trivial} $\SL_2$. 
  
In this case the numerical conjecture asserts that the $H$-period $ \int_{[H]} f_{\phi}  $
vanishes if the parameter $\phi$ doesn't factor through $\check{G}_X$; if it so factors uniquely up to $\check G_X$-conjugacy, say as
$\phi_X: \Gamma_F \rightarrow \check{G}_X$, then  $\int_{[H]} f_{\phi}$ is real and
the conjecture implies 
\begin{equation}\label{IIhat} q^{b_G/2-b_H} \frac{ \left| \int_{[H]} f_{\phi} \right|^2}{ \int_{[G]}
   |f_{\phi}|^2 } =  
\frac{ q^{b_G/2-b_{G_X}} }{\# Z_{\phi}} \frac{L(\frac{1}{2}, S_X) L(1, \check{\mathfrak{g}}/\check{\mathfrak{g}}_X)}{L(1, \check{\mathfrak g}_X) },    \end{equation}
  where we write $b_{G_X} = (g-1) \dim \check G_X$, and the $L$-functions are those associated to the parameter $\phi_X$.  The ratio of $L$-functions arises
  as \[\frac{L^{\norm}(0, T_x \check{M})}{L^{\norm}(0, \check{\mathfrak g})},\] where $x$ is the unique fixed point determined by the factorization $\phi_X$.   It is often in essentially this form that the conjecture
  appears in the literature taking account that $S_X \oplus \check{\mathfrak{g}}/\check{\mathfrak{g}}_X \simeq V_X$;
  the paper of Ichino and Ikeda \cite{II} was particularly influential.
 One must check the constants, which we do not do here. 

Many interesting examples fall in this case.
  See also \S~\ref{SVcompare} for the more complicated analogous discussion in the nontempered case.

\subsubsection{Eisenstein periods} \label{Eisperiodmaintext}
As mentioned  \S \ref{Eisgeom} the Eisenstein case (which does not lie
inside our general framework of dual hyperspherical pairs, but has many formal similarities)
presents peculiarities. We will   compute in \S \ref{Eisappendix} 
the ratio of Eisenstein and Whittaker periods and show that it coincides (formally, because there
are some trivially regularizable infinities) with the ratio of the right-hand sides of \eqref{PXshriekprediction},
{\em multiplied by} $q^{-b_U/2}$. As remarked in \S \ref{Eisgeom} this is an interesting discrepancy which requires study.

\subsection{Nontempered examples: trivial, diagonal, polarized homogeneous}  \label{NTexample}
   
   We now analyze several examples of the nontempered conjecture \S \ref{nontempered}. 
   A particularly interesting class of cases is when the $\SL_2$-type of the form and the space $X$ coincide.
         
      The following Lemma will be useful at several points. 
  \begin{lemma} \label{Compremark} Take
  an orthogonal vector space $W$ with
  $\Gamma_F \times \SL_2$-action respecting the quadratic form. 
  Then, with notation as in \S \ref{Lconven}, 
\begin{equation} \label{Wfnl} q^{(g-1) \dim W}  L(1,  (\mbox{$e$-invariants on $W$})^\shear) =L(0, (\mbox{$e$-coinvariants on $W$})^{\shear}) \end{equation} 
where $e$ is the standard element of a $\mathfrak{sl}_2$-triple; 
moreover, the shearing in the $L$-function refers to the extended Langlands parameters
obtained via  the standard embedding $\Gm \stackrel{\eqref{adef}}{\rightarrow} \SL_2$.
\end{lemma}

\begin{proof}
We write
$W$ as a sum of spaces
$E \otimes V_{\ell}$
where $V_{\ell}$
is the $\SL_2$ representation of dimension $\ell$ with
weights $\{(\ell-1)/2, \dots, -(\ell-1)/2\}$;
and $E$ is self-dual. 
 The $L$-values appearing on the  left and right sides of \eqref{Wfnl} are 
 now  $L(\frac{\ell+1}{2}, E)$ and $L(-\frac{\ell-1}{2}, E)$
and the   functional equation relating these takes  the desired form
 $L(\frac{\ell+1}{2}, E) = q^{ \ell (g-1) \dim E} L(-\frac{\ell-1}{2}, E)$.
 Note that the root number that intervenes, being a global $\epsilon$-factor attached
 to a self-dual unramified representation of $\Gamma_F$, is trivial, because the different is a square. 
 \end{proof}

  \subsubsection{The diagonal case and $L^2$ norms.} \label{L2case} 
We explain why the predicted formula for $L^2$ norms   \eqref{l2nontempered} follows from the general nontempered conjecture \eqref{numerical nontempered}.
Here, suppose that $G$ is semisimple and let us work with the period   (cf.(c) of \S \ref{realperiods})
\[ X = \Delta^d G \backslash G^2,\]
where $\Delta^d$ is the diagonal twisted by the duality involution. Let us take an (unramified) Arthur parameter $\phi_A$ for $G$, and ``double'' it to obtain an Arthur parameter $\Phi_A$ for $G^2=G \times G$
  and associated doubled extended parameter $\Phi_E$. 
We will compute the period of $(f_{\phi}, f_{\phi})$ in the setting of Conjecture \ref{numerical nontempered};
since $f_{\phi}^d =\bar{f_\phi}$, this computes the  square of the $L^2$ norm of $f_\phi$.

 The dual period is $\check{X}=\check{G}$ as $\check{G } \times \check{G}$-space.
 Identify the tangent space $T_g \check{G}$ with $\check{\mathfrak{g}}$
by means of left-invariant vector fields; thus $Z \in \check{ \mathfrak{g}}$
is associated to the derivative of $g e^{tZ}$;
we make a corresponding identification $T^*\check{G} \simeq \check{G} \times \check{\mathfrak{g}}^*$. 
We take the
right action of $\check{G} \times \check{G}$ on $\check{G}$, viz. $(g,h) x = g^{-1}  xh$,
and this induces the action on $T^* \check{G}$ given by $Z \mapsto  \Ad(h^{-1}) Z$ in the second cordinate. 
The moment map 
is given up to sign in each factor by $(g, Z) \mapsto (\mathrm{Ad}(g^{-1}) Z, Z)$.
  
 Take a $\Phi_E $-fixed point
 \[x = (g, Z) \in (T^* \check{G})^{\slice},\]
 where the slice is the one associated to the $\sl_2$-triple of the Arthur parameter by \eqref{checkMslice}.

The moment image equals $(f,f)$, i.e.,  $Z=f$
and $g$ centralizes $f$.  Now (by Lemma \ref{apstar}) 
$g$ also commutes with $\Gm \subset \SL_2$;
so $g$ commutes with $\SL_2$ and lies in the centralizer of the Arthur parameter;
that is to say, the whole fixed space in the slice is
\[ Z(\phi_A) \times f \subset G \times \check{\mathfrak{g}}^* = T^* \check G.\]

We focus on the fixed point $x = (\mathrm{id}_G,f)$; 
 for all other elements of the centralizer will contribute in exactly the same way. 
The tangent space at $x$  to the (isotropic) fiber $\mathrm{id}_{\check G} \times \check{\mathfrak{g}}^* \subset T^*\check G$ intersects
  $T_x^{\slice}$ in the Slodowy slice $f + (\check{\mathfrak{g}}^*)^e$, which is Lagrangian by dimensional considerations:   
its dimension is that of $(\check{\mathfrak{g}}^*)^e$, which is half of the dimension of 
the preimage of $f + (\check{\mathfrak{g}}^*)^e$ inside $T^*\check G$. Thus, the sliced tangent space $T_x^{\slice}$ can be identified with the direct sum of the $e$-invariants\footnote{Note that in the group case we have been using the notation $\check{\mathfrak g}_e$ for the centralizer Lie algebra of $e$; this conflicts with the invariant/coinvariant notation used here. We will temporarily adopt this new notation for the present example, but it should not cause any confusion in the remainder of the paper.} on $\check{\mathfrak{g}}^*$ and its dual, the $e$-coinvariants on $\check{\mathfrak{g}}$.
The relevant $\GGm$-action on the former is through the combination of the action of $\Gm\hookrightarrow \SL_2$ and the squaring $\Gm$-action, while the action on the latter is only through $\Gm\hookrightarrow \SL_2$. For the calculation that follows, let us use the action of $\Gm\hookrightarrow \SL_2$ to shear both spaces,
absorbing the squaring twist into the point of evaluation of the $L$-function.  
Then, \eqref{l2conjecturenontempered} gives upon squaring
\begin{multline*}  q^{-b_G} \left[ \int_{[G]}|f_{\phi}(x)|^2 \right]^2 \\ =  \frac{ (\# Z_{\phi_A})^2}{q^{2b_G}} L(0, (\mbox{$e$-coinvariants on $\check{\mathfrak{g}}$})^\shear) L(1, 
(\mbox{$e$-invariants on $\check{\mathfrak{g}}^{*}$})^\shear
) \\ \stackrel{\eqref{Wfnl}}{=} q^{-b_G} (\# Z_{\phi_A})^2  L(1,  \check{\mathfrak{g}}_e^\shear)^2
\end{multline*}
implying the desired formula \eqref{l2nontempered}, taking into account that all the signs in \eqref{l2conjecturenontempered} are the same,
and that the period is positive.

\subsubsection{$f_{\phi}=1$ and $X=H\backslash G$ homogeneous: } \label{BZGHstuff}
We suppose now that $H$ is semisimple. 
Let $f_{\phi}$ now be the trivial (=constant) form with Arthur packet $\phi_A$.  
We know from the computation of \S \ref{L2case},
combined with the Tamagawa number formula
to compute the volume of $[G]$, 
  that the appropriate normalization is $f_{\phi} = \pm 1$.\footnote{In the nontempered case,  Conjecture \ref{numerical nontempered} pins down
the normalization of forms at best up to  sign, because of the square roots on the right hand. }

It is possible to verify that  the local conjecture implies
 an identification
\[ \check M_{\slice} \simeq \check{J}_H,\]
the {\em group scheme of  dual regular centralizers for $H$,} that is to say,
the group scheme over the Kostant slice whose fiber is the centralizer of a regular element in the dual group. 
This identification will be discussed further in a sequel to this paper (it is related to Proposition \ref{regnilpotent} as well as Example
\ref{JHJG}). 

More precisely, writing $\mathfrak{c}_H, \mathfrak{c}_G$ for the invariant-theoretic quotients of the  Lie algebras of $H, G$ respectively by conjugation action,
the inclusion of $H$ into $G$ induces a morphism $f: \mathfrak{c}_H \rightarrow \mathfrak{c}_G$.
Now, $\mathfrak{c}_H$ is identified with the invariant-theoretic quotient of the {\em dual}
Lie algebra for $\check H$, and via this identification we obtain
a group scheme of regular centralizers $\check{J}_H \rightarrow \check H$ over $\mathfrak{c}_H$,
and similarly for $\check{J}_G$.  There is a morphism of group schemes $f^* \check{J}_G \rightarrow \check{J}_H$ over $\mathfrak{c}_H$,
which in particular gives rise to an action of $\check{J}_G$ on $\check M_{\slice}$.
 
For a parameter of $\check M_{\slice} \simeq \check{J}_H$ to be fixed by the parameter of the trivial representation for $G$, 
it must have image $f$ under the moment map, so map to $0 \in \mathfrak{c}_G$,
and so must also have mapped to $0 \in \mathfrak{c}_H$; that is to say, the fixed points of $\Gamma_F$
acting via $\phi_A$
on $\check M_{\slice}$ are contained in the centralizer of a regular nilpotent on $\check{H}$,
 and in fact are precisely given by the center of $\check{H}$
 (this by consideration of $\mathbb{G}_m$-actions, see Lemma \ref{apstar}). 
 
At each such point $m$,
the tangent space $T_m \check{M}_{\slice}$
admits the centralizer $\mathfrak{v}$ of a regular nilpotent in $H$
as a Lagrangian space; $\phi_E$ 
 acts on $\mathfrak{v}$ through the cyclotomic character raised to the power $2d_i$,
the $d_i$ being the exponents of $H$ (see \S \ref{reductive group notation}). 

Now,
the normalized period of $f_{\phi}=1$ is here given by $q^{-b_H/2} \# \Bun_H$
and using  $\sum (2d_i-1) = \dim H$  we get   
\[ q^{-b_H} (\# \Bun_H)^2   =   (\# Z_{\check{H}})^2 q^{b_H} (\prod \zeta(d_i))^2 = (\# Z_{\check H})^2 \cdot  \prod \zeta(d_i) \zeta(1-d_i). \]
This verifies Conjecture  \ref{numerical nontempered} here.

 \subsubsection{The case when the $\SL_2$-types of $X$ and $f_{\phi}$ coincide; comparison with \cite{SV}} \label{SVcompare}
 We now examine the situation where   $\SL_2$-type of the space $X$ and the automorphic form $f$ are the same. 
 This situation has several simplifying features -- it was for example
 the case in which \cite{SV} proposed a general conjecture, with arbitrary ramification.
 Let us see how the Conjecture \ref{numerical nontempered} recovers the period conjecture
 in a form close to that of \cite{SV} in the unramified case.  We will again suppose that $G$ is semisimple.
   
  We will follow the setup as in \S \ref{omspherical}. 
Let the dual data for $X$
be $\check{G}_X \subset \check{G}, \iota: \SL_2 \rightarrow \check{G},$
and the symplectic $\check{G}_X$-representation $S_X$,
so that 
\begin{equation} \label{cMdef} \check{M} = \mbox{Whittaker induction of $S_X$ along $\check{G}_X \times \SL_2 \rightarrow \check{G}$.}
\end{equation}
As before we put
\begin{equation} \label{VXrecall}
V_X = S_X \oplus \check{\mathfrak{g}}_e/\check{\mathfrak{g}}_X.\end{equation}
Here $V_X$ is graded, i.e., $V_X = \bigoplus_{i} V_X^{(i)}$,
where $S_X$ lies in weight $1$, and $\check{\mathfrak{g}}_e$ is
graded via $\Gm \subset \SL_2$ {\em plus two}, see  e.g.\ discussion in \S \ref{Slodowy}.

To simplify our considerations somewhat we will assume that 
{\em $\check{G}_X$ is the centralizer of $\SL_2$.}  (Otherwise the
considerations below can be modified in a fairly straightforward way
involving sums over possibly more fixed points.)

We are going to consider a parameter 
\[\phi^0_L: \Gamma_F \rightarrow \check{G}_X,\]
giving a Langlands parameter $\phi_L: \Gamma_F \rightarrow \check{G}$ and 
  an Arthur parameter $\phi_A = \phi_L \times \iota:  \Gamma_F \rightarrow\check{G}$, and will
derive from the conjecture the following explicit formula, which is essentially the proposal of \cite{SV} (cf. \eqref{IIhat}): 
\begin{equation}\label{XperiodSV} \frac{\left| \mbox{normalized $X$-period of $f_{\phi}$}\right|^2}{\langle f_{\phi}, f_{\phi} \rangle}  =
\frac{ q^{-b_{G_X}} }{\# Z_{\phi_A}}
\frac{  \prod_{i} L(\phi_L^0, V_X^{(i)}; i/2)}{L(\phi_L^0,  \check{\mathfrak g}_X; 1 )},
\end{equation}

\begin{proof}[Proof that Conjecture \ref{numerical nontempered} implies \eqref{XperiodSV}] 
  $\check{M}$ is Whittaker-induced from $S_X$ 
      along $\check{G}_X \times \SL_2 \rightarrow \check{G}$. 
     In  particular there is a morphism  
 $\check{M}   \rightarrow \check{G}/\check{G}_X \check{U}$
 where  $\check{P} = \check{L} \check{U}$ is the Levi subgroup determined by $\mathfrak{sl}_2$.  We are going to check
 
 \begin{quote} 
 {\em Claim:}  Any $\Gamma_F$-fixed point
  in $\check{M}_{\slice}$   maps to the trivial coset
  in $\check{G}/\check{G}_X \check{U}$. 
 \begin{equation}\label{claimfixedpoint}
~   
  \end{equation}
  \end{quote}
  
  Assuming this for a moment we derive \eqref{XperiodSV}.
Since any fixed point for $\check{M}_{\slice}$ has moment
image equal to $f$,  we see (see e.g.\ \eqref{cfp3} and \eqref{cfp4})
that the fixed points of $\Gamma_F$ on $\check{M}_{\slice}$
correspond precisely to fixed points of $\Gamma_F$ 
on $S_X$ under $\phi_E^0$, the extended parameter corresponding
to $\phi_L^0$. The only such fixed point is the origin
by similar reasoning to the proof of Lemma \ref{apstar}. 
It remains to understand the tangent space $T$ to $\check{M}_{\slice}$ 
at this fixed point $m_0$ corresponding to the origin of $S_X$.

Referring to \eqref{cfp5} 
we see that $T_{\slice}$ has a composition series
 whose associated graded factors are\footnote{
  The latter factor arises, if we follow the notation of
  \eqref{cfp5}, by noting that $\{X \in \mathfrak{g}: [X, f] \in f +\mathfrak{g}_e \}$ 
  is precisely the centralizer of $f$,
  and taken modulo $\mathfrak{h}$ (in the same notation) is dual to  
  $\mathfrak{g}_e/\mathfrak{h}$. }
 \[   
  W := \check{\mathfrak{g}}_e/\check{\mathfrak{g}}_X, 
  S_X, W^{\vee}.\]
  
  Here:
  \begin{itemize}
  \item 
  The $\phi_E$-action on $S_X$ is the action of $\phi_L^0$ on $S_X$
  multiplied by   $\varpi^{1/2}$.
  \item 
  The $\phi_E$-action on $W$  is $\phi_L^0$  multiplied by the 
  action through $\Gamma_F \stackrel{\varpi^{1/2}}{\rightarrow} \GGm(k)$.
 Here $\GGm$ acts through $\Gm \subset \SL_2$ with a further shift by $2$, just as in
 \S \ref{Slodowy}.
  \item The $\phi_E$-action on $W^{\vee}$ is determined from that on $W$ by means of the duality,
  recalling  (Remark \ref{msliceform}) that the pairing $W \otimes W^{\vee}$ is valued in $\kk(1)$. 
  \end{itemize}
  
 Conjecture \ref{numerical nontempered} and
\eqref{l2nontempered} therefore show that 
\begin{multline}   \frac{ P_X^{\norm}(f_{\phi})^2}{\langle f_{\phi}, f_{\phi} \rangle }  = \frac{q^{-b_G}   }{\# Z_{\phi_A}}  
    \frac{ L (0, W^{\vee \shear}) L(1, W^{\shear})   L (1/2, S_X)}{L(1, \check{\mathfrak{g}}_e^\shear)}
 \\    \stackrel{\eqref{Wfnl}}{=}
      \frac{q^{-b_G}   }{\# Z_{\phi_A}} q^{(b_G-b_{G_X})}  \frac{ L(1, (\check{\mathfrak{g}}_e/\check{\mathfrak{g}}_X)^\shear)  L (1/2, S_X)}{L(1, \check{\mathfrak{g}}_X)} 
\end{multline}
where the shearing on $W$ is now considered via the $\Gm \subset \SL_2$
without further shift. 
That is precisely \eqref{XperiodSV}.

We now give the proof of 
Claim \ref{claimfixedpoint}:

  By Lemma \ref{apstar}, any fixed point $x \in \check{M}$  
 is also fixed by the action of
the group $\GGm'$ of \eqref{GGmprime}. Since that action is contracting the Slodowy slice $f+\check{\mathfrak g}_e$ to $f$, any fixed point must lie over $f$ under the moment map. Recall \cite{Gan-Ginzburg} that the product $\check M_f:= G\times (f+\check{\mathfrak g}_e)$  is the Hamiltonian reduction of $T^*\check G$ by $\check U$ over coadjoint orbit of $\check{\mathfrak u}$ of those elements which restrict to the functional $f$ on $\check{\mathfrak{u}}_+=$ the subalgebra 
		whose $\Gm$-weights are $\geq 2$. Since\footnote{If $\Ad(g) f = f+X$ for some nonzero $X \in  \check{\mathfrak{g}}_e$,
			then $\Ad(a_{\lambda} g a_{\lambda}^{-1}) f =  f + \lambda^2 \mathrm{Ad}(a_{\lambda}) X$,
			and by taking $\lambda \rightarrow 0$ we would contradict the fact that   $G$-orbits meet  $f + \check{\mathfrak{g}}_e$ transversally (see \cite[2.2]{Gan-Ginzburg}). } the $G$-orbit of $f$ intersects $f + \check{\mathfrak{g}}_e$
		precisely in $f$, the only points of $\check G/\check U$ over which the moment image of $\check M_f$ contains $f$ are those represented by the centralizer $Z(\check f)$. By construction of the map $\check M\to \check G/\check G_X\check U$, the moment image of a fiber of $\check M$ is contained in the moment image of the corresponding fiber of $\check M_f$ (e.g,  the fiber over the identity coset in $\check M$ has moment image contained in $f+\check{\mathfrak{u}}_+^\perp$), therefore any fixed point of $\check M_{\slice}$ has to live over $Z(f)\check G_X\check U/\check G_X\check U$. On the other hand, if $\check{P}_-$, the parabolic opposite to $\check{P}$
		with respect to the Levi subgroup $\check{L} = $ the centralizer of $a(\Gm)$, the subset $Z(f)\check G_X\check U/\check G_X\check U$ lives over\footnote{To see that an element $g\in Z(f)$ belongs to $\check{P}_-$, we can 
			use \cite[Corollary 20]{McNinch} or
			proceed as follows:
			$h':=\Ad(g) h, f$ form part of a $\mathfrak{sl}_2$-triple, and 
			by \cite[Theorem 3.6]{KostantTDS}, $h'$ must have the form $h+u^{-}$
			where $u^-$ has has negative $h$ weight. In particular,  
			$\mathrm{ad}(h')$ acting on $\mathfrak{p}^-$
			is triangular with respect to a basis of $h$-eigenspaces,
			and all its eigenvalues are therefore $\leq 0$;
			thus, $\mathfrak{p}^-$ is also the sum of negative weight spaces for $\mathrm{ad}(h')$.
			This implies that $g$ normalizes $\check{P}_-$, and therefore $g \in \check{P}_-$.}
		 the open $\check P_-$-orbit on the flag variety $\check G/\check P$, and the only $a$-fixed point in that Bruhat cell is the coset of $1 \check P$. Hence, any $\GGm'$-fixed point on $\check M_\slice$ has to live over a right coset of $\check G_X\check U$ represented by $\check L\cap Z(f)=$ the set of elements that centralize both  $a$ and $f$,
		hence centralize $\SL_2$ (cf. \cite[Cor 3.5]{KostantTDS}). By virtue of our assumption on centralizers, stated after \eqref{cMdef}, it lives over the identity coset of $\check{G}_X \check U$ (while without this assumption, we would have to sum over the set of $\phi_A|_{\Gamma_F}$-fixed points of the centralizer of $\SL_2$ mod $\check G_X$, analogously to the tempered case).

\comment{
With $a$ the cocharacter as in \eqref{GGmprime}, if 
 $\phi_E$ fixes a point on the slice over $g \check{G}_X \check{U}$  (some $g \in \check{G}$) it  thus follows that
\begin{equation} \label{aeq}  a g \check{G}_X \check{U}  = g \check{G}_X \check{U} a\end{equation} 
and we want to show $g \in \check{G}_X \check{U}$.

The moment image of any point of $\check{M}$ mapping to $g \check{G}_X \check{U}$   
lies in $\Ad(g) (f + \mathfrak{u}_{+}^{\perp})$,
with $\mathfrak{u}_+$ the subalgebra of $\mathfrak{u}$
whose $\Gm$-weights are $\geq 2$. 
Note that : 
Therefore, 
\[ \Ad(g^{-1}) f \in f + \check{\mathfrak{u}}_+^{\perp} \implies g \in  Z(\check{f}) \check{U}\]
where we use the fact that the orbit map
from $\check{U} \times (f +\check{\mathfrak{g}}_e)$ 
is an isomorphism onto its image  $f + \check{\mathfrak{u}}_{+}^{\perp}$
  \cite[Lemma 2.1]{Gan-Ginzburg}.
Modifying $g$ on the right by $\check{U}$ as we may harmlessly do,
we may suppose, then that $g$ centralizes $f$.  
In particular, $g \in \check{P}_-$, the parabolic opposite to $\check{P}$
with respect to the Levi subgroup $\check{L} = $ the centralizer of $a(\Gm)$. 

From \eqref{aeq} we deduce that
\[g^{-1} a g a^{-1} \in \check{G}_X \check{U} \subset \check{P}.\]
Given that $\check{g} \in \check{P}_- = \check{L} \check{U}_-$ this inclusion forces  $g \in \check{L}$,
 and so centralizes $a$.  Then $g$ centralizes both $a$ and $f$,
 and so centralizes $\SL_2$ (cf. \cite[Cor 3.5]{KostantTDS}) so actually belongs to $\check{G}_X \subset \check{P}$
 by virtue of our assumption on centralizers, stated after \eqref{cMdef}. 
}
 \end{proof}

   \subsection{How are the geometric and numerical conjectures related?} \label{mumerical}
   
   Let us now return to the question of how the geometric and numerical conjectures are related.
In short, the numerical conjecture should be a consequence of the geometric one,
 but the deduction involves technical issues that we have not studied (and, particularly in the nontempered case,
 may also involve some new structures of independent interest). 
 
    We restrict ourselves to the polarized case: $M=T^*(X, \Psi), \check{M} =T^* \check{X}$. 
 In the finite context, the geometric form of the conjecture, Conjecture \ref{GlobalGeometricConjecture}, 
  asserts
 that a certain period sheaf $P_X^{\spec}$  (the spectral projection of $P_X$)
 matches, under a suitable form of the Langlands equivalence,
 with the $L$-sheaf $L_{\check{X}}$. 
 
   Let $f=f_{\phi}$ be a Hecke cuspidal eigensheaf on $\Bun_G$,
  with Langlands parameter $\phi$, a $\kk$-point of $\Loc_{\check{G}}$.
  Let us suppose that:
  \begin{itemize}
  \item[(a)]
   $\phi$ ,
  restricted to geometric $\pi_1$, fixes a single point on $\check{X}$,
  and the same is true for the dualized parameter $\phi^d$. 
  \item[(b)]  $f_{\phi}$
    is a   pure self-dual perverse sheaf.  
    \end{itemize}

Let $f(x): \Bun_G(\FF_q) \rightarrow \kk$ be the associated function, and let $\widetilde{f}$
be the function arising from the Verdier dual, i.e., the conjugate of $f$
after fixing an isomorphism $\kk \simeq \mathbb{C}$.
Numerically (b) has the effect that the $L^2$-norm of $f$ is of size about $1$. 
As in \S \ref{anglebracketnotation}
   let us write $[\dots]$ for geometric Frobenius trace
   on a vector space.

Since spectral projection in our context
amounts (at least conjecturally, see \S \ref{spectral action} and references therein) to the
right adjoint to the inclusion of nilpotent sheaves, we find
\begin{equation} \label{joeb} \Hom(f, P_X^{\norm}) = \Hom(f, \textrm{spectral projection of $P_X^{\norm}$}) \simeq
\Hom(\delta_{\phi}, \mathcal{L}_{\check{X}}^{\norm}).\end{equation} 
Passing to Frobenius trace, and using
Lemma \ref{Homlemma} to evaluate the left-hand side
and \eqref{LX1b} to evaluate the right hand side, 
   \[ \sum f_{\phi}(x) P_X^{*,\norm}(x) = q^{-b_G/2} L^{\norm}(1, \phi^{d}, T^{\vee\shear}),\]
     where on the  left  the sum is over $\Bun_G(\FF_q)$ and we weight by inverse-automorphisms;
and  $T$  is the tangent space to $\check{X}$ at the unique fixed point
for the dualized parameter $\phi^d$. On the far right,  $T^{\vee \shear}$ is considered
as a $\Gamma_F$-module through $\phi^d$.
 
That is precisely \eqref{PXstarprediction} of the numerical conjecture Conjecture \ref{numconj}.
That is to say, the geometric conjecture, together
with the anticipated identification of the spectral projection, implies the {\em star period}
part of the numerical statement in the case that $\phi$
has finite centralizer.
 
 \begin{remark}
The same reasoning also suggests why the $(-1)^s$
in the statement of \eqref{PXstarprediction} needs to be there
when $f_{\phi}$ is no longer assumed cuspidal: 
it arises (eventually) from a cohomological shift by $s$. 
(In more detail, it should arise for the same reason as the $(-1)^d$ in the second line of \eqref{LX1b},
whose source can be seen at \eqref{BZrev}.)
\end{remark}

\begin{remark}   \label{unshearremark}
Suppose that $\check{G}$ fixes the unique fixed point above. Then we may rewrite
\begin{equation} \label{unshearequation}  L(1, \phi^d, T^{\vee \shear}) = L(1, \phi, T^{\unshear}),\end{equation}
where on the right $T^{\unshear}$ means that we shear by the negated $\Gm$-action on $T$. 

Indeed,  the representation $\phi^d_E$ of $\Gamma_E$ on  $T^{\vee}$ is obtained as the composite:
\[ \Gamma_F \stackrel{(\phi, \varpi^{1/2})}{\longrightarrow} \check{G} \times \GGm \rightarrow \GL(T^{\vee}),\]
where $\check{G}$ acts on $\GL(T^{\vee})$ through
its standard action {\em precomposed with the dualizing involution}, 
and $\GGm$ acts on $T^{\vee}$ through its action arising 
from the $\GGm$-action on $\check{X}$. 

Now precomposition with the dualizing involution {\em on both $\check{G}$ and $\Gm$}
switches the isomorphism class of $T$ and that of $T^{\vee}$.
Therefore, the action of $\check{G} \times \Gm$ on $T^{\vee}$
that occurs above is isomorphic to the action of $\check{G} \times \Gm$
on $T$, which is the standard action on the first factor,
but the {\em inverse} of the action on the second factor. This explains \eqref{unshearequation}. 
 \end{remark}

However, the star period part of the numerical statement is also the least well-attested part
by classical computations! 
One certainly wants to carry out the same deduction for the $!$-period, and for this
we would want to argue that
\eqref{joeb} holds in the opposite direction, i.e., $\Hom(P_X^{\norm}, f)$
is unchanged if we replace $P_X$ by its spectral projection.
Assuming that is valid, 
the argument above goes through to say

\begin{equation} \label{Xperiod}  
  \sum_{x} P_X^{\norm}(x) \widetilde{f_{\phi}}(x)=q^{-b_G/2} L^{\norm}(0, \phi^d,  T^{\shear}),
  \end{equation}
  where we now used \eqref{LX0b}.
    But on the other hand 
the Langlands parameter of $\tilde{f}$ should be the image of $\phi$
 under the dualizing involution,  and therefore, replacing $\phi$ by $\phi^d$, we get
\begin{equation} \label{op}  \sum_{x} P_X^{\norm}(x) f(x) = q^{-b_G/2} L(0, T^{\shear}),\end{equation}\
which is now \eqref{PXshriekprediction} of Conjecture \ref{numconj}. 
On the right, $T^{\shear}$ is considered as a $\Gamma_F$-module in the ``obvious'' way,
i.e., through $\phi$ and not through $\phi^d$.

 Consequently,  for $f$ normalized as above, and
 {\em if  we were to assume  that the left nilpotent projection exists and coincides with the right nilpotent projection}, the 
 geometric conjecture implies the numerical statement  \eqref{op} about the period of the form $f$.
 In general, this left nilpotent projection need not exist,
but we may hope that some suitable ``interpretation'' of it does (cf. \S \ref{igd2}). 
 At the moment, then, the numerical conjecture is not a consequence of the geometric conjecture,
 but rather a parallel statement.

\begin{remark}[Period and $L$ distributions]\label{AGKRRV numerical} 
\index{algebraic distributions}
Recall from \S \ref{derived volume forms} (following \S \ref{geometric to arithmetic Langlands}) the proposed formulation of L-functions as meromorphic algebraic distributions (i.e., meromorphic sections of the dualizing complex) on $\Loc_\Gv^{\textrm{arith}}$ in the setting of~\cite{AGKRRV1}. This suggests a parallel formulation of the numerical period conjecture, which is what one might hope to obtain by taking the categorical trace of Frobenius directly on the geometric period conjecture. Namely, we expect $L$-distributions to be identified under the unramified arithmetic Langlands correspondence with the corresponding dual period functionals, represented as elements in localizations (via the spectral action) of the space of automorphic functions $\kk_c[\Bun_G(\F_q)]$ over open subsets of the stack of arithmetic local systems. 
See also \S \ref{AGKRRV Arthur} for an analogous discussion of Arthur parameters.
\end{remark}

\subsection{Star periods and asymptotics} \label{starperiods}

We will now discuss the star period function $P_X^*$, formulating some conjectures.\footnote{Our focus here will be not so much on Conjecture \ref{numconj} itself,
but rather about even more basic formulas relating to  the $*$-period.
  Tony Feng, Jonathan Wang and the third-named author have
 carried out some computations supporting the $*$-period assertion in Conjecture \ref{numconj}, but we will not describe  such computations here, 
 except for Example \ref{starTate}.}

Recall from Remark \ref{starperiod} that the star version of the period sheaf is defined
as the $*$-pushforward of the dualizing sheaf from $\Bun^X_G$, and, by definition,
$P_X^*$ is the function obtained by taking trace of Frobenius.

{\em To simplify we assume that the $\GGm$-action on $X$ is trivial throughout this section.}
 The work of Schieder and Wang suggests that $P_X^*$ can generally be computed in terms 
of the theory of asymptotics on $X$, which
we now briefly recall.  

Assume, here, that $X$ is affine and homogeneous. Recall that to each subset $\Theta$ of the roots of $\check{G}_X$ we may attach
a boundary degeneration $X_{\Theta}$ of $X$ (terminology of \cite{SV} -- in particular, $X_\Theta$ is homogeneous). For every (nonarchimedean) place $F$, there is an ``asymptotics'' map
\[ \textrm{asymp}_{\Theta}^*: C^{\infty}(X(F)) \rightarrow C^{\infty}(X_{\Theta}(F)),\]
so that a function is equal to its image close enough to infinity ``in the $\Theta$ direction.'' 

Restricted to compactly supported functions $C_c^\infty(X(F))$, the asymptotics map is known to have image in a space of functions of moderate growth and bounded support -- ``bounded'' means that it has compact closure in an affine embedding of $X_\Theta$. Taken over all places together, we obtain a map 
\[ \textrm{asymp}_{\Theta}^*: C_c^{\infty}(X(\mathbb A)) \rightarrow C^{\infty}(X_{\Theta}(\mathbb A)),\]
whose image consists of functions of bounded support. 
Now let $P_{\Theta}: C_c^{\infty}(X_{\Theta}(\mathbb A)) \rightarrow C^\infty([G])$
be the theta series,  that is to say, sending a function $f$
to $\sum_{X_{\Theta}(F)} f(xg)$; 
it extends to smooth functions of bounded support.  
 \begin{conjecture}\label{conjecture-star}
 Suppose that $X$ is affine homogeneous with point stabilizer $H$ (and trivial $\GGm$-action).  Then the $*$-period $P_X^*$ is obtained by evaluating 
\begin{equation}  q^{-b_H} \label{altsum} \sum_{\Theta} (-1)^{|\Theta|} P_{\Theta}  \circ \mathrm{asymp}_{\Theta}: C^{\infty}(X) \rightarrow C^\infty([G])\end{equation}
 at the basic vector $\delta_X \in C^{\infty}_c(X)$, that is to say, at the characteristic function of integral points.
 \end{conjecture}
 
 The factor $q^{b_H}$ arises from the fact that the star period is obtained from the dualizing sheaf and not the constant sheaf. 
 It may be possible to use a suitable variant of \eqref{altsum}  as a {\em definition} of the $*$-period in the number field case. One needs to give a suitable definition to $\mathrm{asymp}_{\Theta}$ at archimedean places.

  In any case, Conjecture \ref{conjecture-star} allows one
  to explicitly compute the star period fairly readily,  and, where we have looked, it appears to be  compatible with our numerical conjecture
Conjecture \ref{numconj}.  The strange sign $(-1)^s$ in Conjecture \ref{conjecture-star} is related to the alternating sum above, but in a subtle and beautiful way, since several $\Theta$ will contribute to a given pairing. 

 We discuss here only the simplest example.  Assume that $X$ is wavefront (terminology of \cite{SV}); 
then the only term in \eqref{altsum} that can have a nontrivial pairing with cusp forms is the term for $\Theta=\emptyset$, and we get
\[ \langle P_X^{*}, f_{\phi} \rangle = q^{-b_H} \langle P_X, f_{\phi} \rangle\]
for $f_{\phi}$ a cusp form.  Our assumptions guarantee that the eigencharacter $\eta$ is trivial,
and therefore $P_X^{*\norm} = P_X^*[-\beta_X], P_X^{\norm} = P_X[\beta_X]$; since $\beta_X=b_H$ here,
we would also get:
\begin{equation} \label{Sdphenom} \langle P_X^{*\norm}, f_{\phi} \rangle =  \langle P_X^{\norm}, f_{\phi} \rangle\end{equation}
According to Conjecture \ref{numconj}, the left and right hand are given
(in the unique fixed point case) in the general form
$L^{\norm}(0, T^{\shear})$ and $L^{\norm}(1, T^{\unshear})$. 
Although we do not have a general analysis,
this equality arises in the examples we have looked at
in the following way:    $T$ is a sum
of components $T_0 \oplus T_1$ in $\GGm$ weights $0$ and $1$.
Then (again, in examples we studied) $T_0$ is self-dual as a $G$-representation, which
gives $L^{\norm}(0, T_0) = L^{\norm}(1, T_0)$;
and for $T_1$ we have \[L^{\norm}(0, T_1^{\shear}) = L^{\norm}(\frac{1}{2}, T_1) 
= L^{\norm}(1, T_1^{\unshear}).\]

 \begin{remark}
Conjecture \ref{conjecture-star} generalizes \cite[Theorem C.7.2]{wangbilinear}, and we hope that it can be proven along the lines of Wang's argument, by \emph{compactifying} the morphism $\Bun^X\to \Bun_G$, using an $X$-analog of the Vinberg monoid.

Namely, consider the \emph{affine degeneration} of $X$: this is an affine family $\mathcal X \to \overline{A_{X,\text{ad}}}$, where $A_{X,\text{ad}}$ is the quotient of $A_X$ whose character group is spanned by the spherical roots, and $\overline{A_{X,\text{ad}}}$ is its toric embedding corresponding to the dual of the cone of spherical roots. The family carries an action of $G\times A_X$, and contains both $X$ and an affine embedding of $X_\Theta$ (for all $\Theta\subset \Delta_X$) as special fibers. See \cite{PopovContractionsactionsreductive1987}, \cite[\S 5.1]{GaitsgoryNadler}, \cite[\S 2.5]{SV} for the construction; the precise base $A_{X,\text{ad}}$ is (probably) not very important for the argument we are outlining, and one can replace $A_{X,\text{ad}}$ by a torus that is isogenous to it.
Let $\mathcal X^\bullet$ be the open subset whose fiber over every point on the base is the open $G$-orbit in the corresponding fiber of $\mathcal X$. The torus $A_X$ acts freely on it, and the quotient $\mathcal X^\bullet/A_X$ is a compactification of $X/\mathcal Z(X)$ (sometimes called the \emph{wonderful compactification}, although this term is usually reserved for the cases when it is smooth). As in \cite[Lemma C.8.2]{wangbilinear}, we expect that the morphism
\[ \Map^\bullet(\Sigma, \mathcal X/G)/A_X \to \Bun_G,\]
where the bullet denotes maps which generically land in $\mathcal X^\bullet$,  compactifies the map $\Bun^X\to\Bun_G$, and can be used to address Conjecture \ref{conjecture-star}.
\end{remark}

\medskip

Conjecture \ref{conjecture-star} fails when $X$ is not affine homogeneous, as the argument that we outlined breaks down (e.g., the requirement for objects of $\Map^\bullet(\Sigma, \mathcal X/G)/A_X$ to lie generically in $\mathcal X^\bullet$ misses out a part of $\Bun^X$, when $X$ is not homogeneous).  
An interesting example is the Iwasawa-Tate case, with which we conclude:

\begin{example}
\label{star1} \label{starTate}
 The Iwasawa-Tate case is of interest precisely because the self-duality phenomenon
described after \eqref{Sdphenom} fails and the mechanism for compatibility between star and $!$ conjectures
is somewhat different.

In the Iwasawa-Tate case with neutral $\GGm$-action we have, for $L$ a point of $\Bun_{\GGm}$,
\begin{equation} \label{PX!Tate} P_X^{\norm} (L)=   q^{h^0(L \otimes K^{1/2}) -\frac{1}{2} \deg(L \otimes K^{1/2})},\end{equation}
\begin{equation} \label{PX*Tate} P_X^{*\norm} (L)=  q^{\mathrm{deg}(L)/2- (g-1)/2}  + q^{-\mathrm{deg}(L)/2-(g-1)/2} - P_X^{\norm}(L) \end{equation}
\end{example}
 
 We readily verify the predictions of Conjecture \ref{numconj}:
 if we pair $P_X^{\norm}$ with a character $\chi$ we get 
 $\chi(K^{1/2})^{-1} L(\frac{1}{2}, \chi)= L^{\norm}(1/2, \chi)$,
 where we regularized the pairing in the unique $G$-invariant way; 
 if we do the same for $P_X^{*\norm}$ we get $-L^{\norm}(1/2, \chi)$. 
 It is instructive to consider the asymptotic behavior. Write $r=\deg(L)$.  As $|r| \rightarrow \infty$ we have
 \[P_X^{\norm} = q^{\frac{1-g+|r|}{2}}, P_X^{* \norm} = q^{\frac{1-g-|r|}{2}}.\]
 Therefore $P_X^{\norm}$ blows up at $\infty$ whereas $P_X^{*\norm}$ decays.

We already proved \eqref{PX!Tate}, see \eqref{PXTatenorm}. 
We sketch the argument to check \eqref{PX*Tate}: we can split $\Bun_X$
here into the open $\Bun_X^{\circ}$ and a closed zero-section that is identified with $\Bun_{\Gm}$.
Each fiber of $\Bun_X \rightarrow \Bun_{\Gm}$ is a $\Gm$-torsor over a projective space,
and correspondingly the   $*$-pushforward of $\omega$ to $\Bun_X^{\circ}$ can be checked to be a shift of the $!$ pushforward
of the constant sheaf, which can be numerically computed fiber by fiber.
Observe that  the symmetric form of \eqref{PX*Tate} is not seen by this way of computing: the first term comes from the zero-section, 
and the remaining two terms come from $\Bun_X^{\circ}$.

\subsection{Arthur functoriality} \label{Arthur functoriality}

 In this speculative subsection we will make some suggestions\index{Arthur parameter}
of geometric interpretations of the formulas presented in 
the nontempered case (\S \ref{nontempered} and \S \ref{NTexample}) and the role of Arthur functoriality. We postpone a somewhat sharper discussion in the geometric setting  to \S \ref{WhitArth}, where we discuss also  Arthur functoriality as an operation on arithmetic field theories.

The basic situation for this section (as in other discussions of spectral Whittaker data such as \S \ref{Whittaker induction} and \S \ref{Whittaker functoriality}) is that we are given a homomorphism 
\[ \iota: \check{H} \times \SL_2 \rightarrow \check{G}\] (where we often further assume that $\check{H} \subset \check{G}$ is the centralizer of the $\SL_2$). We restrict ourselves to even $\SL_2$'s, i.e., we demand that the corresponding cocharacter $\varpi_\iota$ acts on $\fgv$ with only even weights (see Remark~\ref{odd SL2s} for a discussion of the odd case). 

Arthur's conjectures \cite{Arthur-unipotent-conjectures} suggest very broadly that the automorphic representation theory of $G$ is built out of the {\em tempered} automorphic representation theories of the group $H$ Langlands dual to $\Hv$. The process
by which nontempered representations are built out of tempered ones, which we might informally call ``Arthur lifting'' or ``Arthur functoriality,'' is analogous to the conjectural Langlands functoriality relating automorphic forms on groups $H$ and $G$ coming from an inclusion $\Hv\subset \Gv$ of their duals --
indeed, that is the special case of Arthur functoriality when $\iota$ is trivial. 

Continuing with the setup from \S \ref{apsetup}, an Arthur parameter $\phi_A$ associated to $\iota$ defines a $\Gamma_F$-representation into the centralizer $\Hv$ of the $\SL_2$ (i.e., a Langlands parameter into $\Hv$) which is pure of weight zero, and the associated Langlands parameter $\phi_L$ into $\Gv$ is a shift of it by $\iota|_{\Gm}: \Gm \rightarrow \SL_2 \stackrel{\iota}{\rightarrow} \Gv$ (composed with the square root of the cyclotomic character).
(cf. ~\eqref{phiLdef}).

This passage $\phi_A\mapsto \phi_L$ has a natural geometric version which we explain in \S \ref{sheared Loc}  --- namely, given an $\Hv$-local system, we can shear the induced $\Gv$-local system by its $\Gm$-symmetry coming from $\iota|_{\Gm}$. The result is not an ordinary $\Gv$-local system but rather a {\em derived} or sheared local system -- its ``associated vector bundles'' are cohomologically graded through $\iota|_{\Gm}$, a geometric counterpart of the nontempered nature of the Langlands parameter $\phi_L$. Moreover the resulting local system comes equipped with a Lefschetz operator (endomorphism of cohomological degree 2). The model example of such an object is the cohomology of an algebraic variety over the curve $\Sigma$ equipped with its Lefschetz 
operator. 
 
 In the geometric setting we also have access to a much sharper version of this pointwise construction of Langlands parameters into $\Gv$ from (tempered) Langlands parameters into $\Hv$.
Namely, in \S \ref{Linduction} we have constructed an Arthur (or spectral Whittaker) induction functor
 \[\mathsf{AI}: \mbox{spectral category for $\check{H}$}  \longrightarrow \mbox{spectral category for $\check{G}$}\]
 which by Lemma \ref{1022} interacts nicely with $L$-sheaves: 
 \[ \mbox{$L$-sheaf of $\check{Y}$} 
\stackrel{\mathsf{AI}}{\rightsquigarrow} \mbox{$L$-sheaf of Whittaker induction $\check{X}$ of $\check{Y}$}.\]
Moreover, in Corollary~\ref{sheared eigenvalues} we show, in the geometric setting, that $\mathsf{AI}$ interacts in the expected way with the underlying operation $\phi_L\leftrightarrow \phi_A$  on Langlands parameters.  Applying the geometric Langlands correspondence, the Arthur induction induces (conjecturally) a functor 
 \[\check{\mathsf{AI}}: \mbox{automorphic category for $H$}  \longrightarrow \mbox{automorphic category for $G$}\]
which provides a geometric analog of Arthur lifting --- it will send Hecke eigenobjects  in the ordinary sense to Hecke eigenobjects with $\iota$-sheared derived local systems as eigenvalues. 
As we discuss in \S \ref{WhitArth} this suggests a strong geometric form of Arthur's conjectures, in the form of a 
semi-orthogonal decomposition of the automorphic category for $G$, indexed by $\SL_2$-parameters,
and with ``associated graded'' pieces
generated by tempered automorphic sheaves for the  various $\SL_2$-centralizers.

\begin{remark}[Arthur lifting]\label{AGKRRV Arthur}

In the formulation of the unramified arithmetic Langlands correspondence of~\cite{AGKRRV1} reviewed in
\S \ref{geometric to arithmetic Langlands}, one can formally consider taking the Frobenius trace of the functor $\mathsf{AI}$ to obtain an arithmetic Arthur induction map
 $ \Gamma(\Loc^{arith}_{\Hv},\omega)\longrightarrow \Gamma(\Loc^{arith}_{\Gv},\omega)$ on distributions on the space of Langlands parameters. (Such a map will only exist meromorphically, i.e., after restriction to a suitable open substack. See \S \ref{derived volume forms} and \S \ref{AGKRRV numerical} for a parallel discussion of $L$-functions and periods as meromorphic algebraic distributions.) Dually, assuming the unramified Langlands correspondence, this would provide an Arthur lifting map on the full space of unramified automorphic forms, after a suitable localization with respect to the Hecke action;
 ignoring this latter subtlety, this is a map
 \[\check{AI}: \kk[\Bun_H(\F_q)]_c \longrightarrow \kk[\Bun_G(\F_q)].\]

In the arithmetic Langlands program the Arthur parameterization is not usually thought of in 
terms of such a map $\check{AI}$, because the Arthur parameterization
concerns automorphic representations, which at best pin down individual functions up to scaling. 
The point is, however, that 
the theory of periods suggests that there is a distinguished way to normalize
automorphic forms (although our conjecture often only normalizes them up to sign, 
it seems likely that the sign ambiguity can also be resolved).  The map $\check{AI}$
above then takes  normalized tempered forms on $H$ to normalized nontempered forms on $G$.

Notice that we have assumed above that the $\SL_2$ has even weights. This is in fact an important assumption here;
in the odd case we can encounter issues of {\em anomaly}, 
and it is plausible that there is no natural way to make a map $\check{{AI}}$
without it annihilating some tempered eigenforms (this is related to the existence of CAP forms). 

 \end{remark}

\subsubsection{Nontempered periods, revisited}
We now revisit the conjectural description of nontempered periods, Conjecture \ref{numerical nontempered}, whose notation we keep;
however, to simplify our notation, we will just for now index $L$-sheaves by Hamiltonian spaces $\check{M}$ rather
than by their polarizations $\check{X}$. 

Thus we would like to describe the $M$-period of an Arthur form $f_\phi$ for a polarized hyperspherical variety $M=T^*(X, \Psi)$.
The geometric counterpart to this is the Hom pairing 
\[\Hom(\cP_X, \check{\mathsf{AI}}(\cF_{\phi}))\]
 between the period sheaf associated to $M$, and the geometric Arthur lift of a Hecke eigensheaf $\cF_\phi$ on $\Bun_H$, i.e., the automorphic sheaf corresponding to the Arthur induction of a skyscraper sheaf $\mathsf{AI}(\cO_{\phi})$ on $\Loc_\Hv$.

Now, recall from  Remark \ref{ArthurResRemark} the operation of Arthur restriction (or Arthur-Jacquet) $\mathsf{AJ}$, which performs Whittaker reduction (adjoint to Whittaker induction) on the level of Hamiltonian spaces :
\begin{equation} \label{WJnum} \mathsf{AJ}: \mbox{$L$-sheaf of $\Mv$} \rightsquigarrow \mbox{$L$-sheaf of Whittaker reduction $\Mv\GIT_\psi \check{U}$}
 \end{equation} 
where $\check{U}$ and $\psi$ are the unipotent subgroup and character determined by $\iota$.

To calculate the periods of Arthur forms geometrically, we need a guess (based on similar phenomena in representation theory rather than explicit calculation) that $\mathsf{AJ}$ can be identified with the left adjoint of $\mathsf{AI}$.
It is by no means clear that such an identification is even morally correct, 
but we may hope that it is close enough to true for numerical purposes, and conversely take known evidence for the numerical statements about nontempered periods as suggestive of this statement, on the level of Frobenius traces.
Given this highly optimistic setup,  and ignoring issues of normalization, we
compute
\begin{eqnarray*}
\Hom_{\Bun_G}(\cP_X, \check{\mathsf{AI}}(\cF_{\phi}))&\simeq& \Hom_{\Loc_\Gv}(\Ll_{\Mv}, \mathsf{AI}(\cO_{\phi}))\\  
&\simeq& \Hom_{\Loc_\Hv}(\mathsf{AJ}(\Ll_{\Mv}), \cO_{\phi})\\
&\simeq& \Hom_{\Loc_\Hv}(\Ll_{\Mv\GIT_\psi \check{U}}, \cO_{\phi})
\end{eqnarray*}
In other words, the $M$-period of an Arthur sheaf on $G$ is calculated spectrally by the Whittaker reduction, or Slodowy slice, to the dual hyperspherical variety $\Mv$.
Passing from sheaves to functions as discussed in \S \ref{mumerical}, this discussion suggests that $X$-periods of Arthur lifts on $G$ are given by $L$-functions on $\Hv$ associated to the Slodowy slice to $\Mv$, which -- modulo finer issues such as normalization -- is precisely what Conjecture \ref{numerical nontempered} says.  

If we further assume that the Hamiltonian $\Hv$-variety $\Mv\GIT_\psi \check{U}$ has a Hamiltonian dual $H$-variety $M_H=T^*(X_H, \Psi_H)$, then we can further describe the interpret the final term above automorphically:
\[\Hom_{\Bun_G}(\cP_X, \check{\mathsf{AI}}(\cF_{\phi}))\stackrel{?}{\simeq} \Hom_{\Bun_H}(\cP_{X_H}, \cF_\phi).\]
where the $?$ reminds that this is not a theorem, but based on our optimistic
speculations about adjointness of $\mathsf{AI}, \mathsf{AJ}$. 

\begin{remark}[The $L^2$ picture]
We can apply the above analysis in the group case, as in \S \ref{L2case} and \S \ref{diagdiag}, to describe the $L^2$ norm of an Arthur form in the arithmetic setting or the endomorphisms of an Arthur sheaf
$\check{\mathsf{AI}}(\cF_{\phi})$ in the geometric setting.
 First let us 
 consider the formula \eqref{l2nontempered} for $L^2$ norms:  \[ \langle f_{\phi}, f_{\phi} \rangle= (\# Z_{\phi}) L(1,  \check{\mathfrak{g}}_e^{\shear}).\]
In the case $\# Z_{\phi}=1$ the right hand side
 can be interpreted as trace of Frobenius
on the self-Hom of a $\check{G}$-derived local system equipped with a Lefschetz operator.
This is, at least in an informal sense, compatible with the  geometric picture discussed earlier in this section.
\end{remark}

%% file: automorphic-factorization.tex
\newcommand{\basic}{\delta}

\label{part4longintropage}

Our main goal in Part 4 of this work is to explain how to relate the local theory from Part 2 with the global theory from Part 3. A one-sentence summary is 
 \begin{quote} {\em the local categories (automorphic and spectral) provide Hecke constraints satisfied by the period and $L$-sheaves, and these constraints are intertwined by the local and global conjectures.}
 \end{quote}
This can be compared with the picture developed in~\cite{SV} in which $X$-periods of automorphic forms are given by Euler products whose local factors describe the Plancherel measure on spherical functions on $X_F$.

We first describe the one-point form of local-global compatibility in \S \ref{local-global}, asserting that the local and global period conjectures, Conjectures~\ref{local conjecture} and~\ref{GlobalGeometricConjecture}, are intertwined by Hecke-linear functors from the local to the global categories: the unramified automorphic and spectral {\em $\Theta$-series}.

This form of local-global compatibility can be significantly sharpened by inserting the local category at many points $\{x_i\}\subset \Sigma$ and allowing these points to vary and collide. This is captured by the notion of {\em factorization}, which we apply in \S \ref{automorphic-factorization}-\ref{spectral geometric quantization}.

{\bf
The ideas of this Part are to a much greater extent than before obstructed by technical issues, in particular issues of sheaf theory in infinite type on the automorphic side and the de Rham spectral side, so that the most complete picture we present is on the Betti spectral side. }

In \S \ref{automorphic-factorization} we introduce a factorizable form of the local category and of the Plancherel algebra. The mechanism of factorization homology -- a geometric analog of Euler products -- produces a global counterpart to the Plancherel algebra, the {\em RTF algebra} $\RTF_{X,\Sigma}$, an algebra in the global Hecke category. We also sketch the idea of a factorizable form of $\Theta$-series, which realizes the RTF algebra as a ``locally defined'' source of maps between the period sheaf and its Hecke transforms, and indicate its relation to the relative trace formula.

In \S \ref{spectral-factorization} we discuss the spectral counterpart to the factorizable Plancherel algebra, the {\em $L$-algebra} $\bO_\Mv$, and its relation to the hyperspherical variety $\Mv$ via the mechanism of {\em spectral deformation quantization}. This gives rise to a factorizable form of the local conjecture, in the setting where both $M$ and $\Mv$ are polarized.

In \S \ref{spectral geometric quantization} we study the spectral counterpart of the RTF-algebra, the {\em $L$-observables} $\bO_{\Mv,\Sigma}$ and its action on the $L$-sheaf. We are guided by the analogy that the $L$-observables are to the $L$-sheaf as deformation quantization is to geometric quantization, an analogy which we flesh out in a couple of ways. 
The $L$-observables are a geometrization of the $L$-function of $\Mv$ (the square of the $L$-function of $\Xv$) and its 
factorization homology construction is a geometric counterpart of the corresponding Euler product.
We conclude with a foray into the geometric study of Arthur parameters, whose construction is closely related to the theory of $L$-sheaves for twisted polarizations.

 \section{Theta series and local-global compatibility}\label{local-global}
 \index{$\Theta$-series}

In \S \ref{section-unramified-local} we studied the local automorphic and spectral categories $\SHV(X_F/G_O)$ and $\QCshear(\Mv/\Gv)$ and formulated the local conjecture, Conjecture~\ref{local conjecture}, identifying the two. Among the consequences of this conjecture is that the local automorphic category is controlled by the Plancherel algebra $\Planch_X$ studied in \S \ref{PlancherelCoulomb}, which is identified with the sheared coordinate ring $\bO_\Mv:=\cO^{\shear}(\Mv/\Gv)$ (the ``$L$-algebra'') as algebra objects in the spherical Hecke category.

On the other hand we formulated   in \S \ref{GGC} a global conjecture on an algebraic curve $\Sigma$
when $\Mv=T^*\Xv$ is polarized (or more generally twisted polarized $\Mv=T^*_\Psi\Xv$).
This  statement, 
Conjecture~\ref{GlobalGeometricConjecture},  relates the period sheaf $\cP_X$ on $\Bun_G(\Sigma)$ studied in \S \ref{section:global-geometric} with the L-sheaf $\Ll_\Xv$ on $\Loc_\Gv(\Sigma)$ studied in \S \ref{Lsheaf}.\footnote{The conjecture depends on a choice of spin structure $\cK^{1/2}$ on $\Sigma$, which can be eliminated by the use of C-groups as in \S \ref{spindep},~\ref{Lspinindep},~\ref{extended-group appendix}.}

In order to relate the local and global conjectures we first remove the coordinate dependence in the local conjecture. We write the stack $X_F/G_O=\Bun_G^X(D,D^*)$ as the moduli of $G$-bundles on the disc with a section of the associated $X$-bundle on the punctured disc; when $\GGm\actson X$ is nontrivial, we rather work with a ``normalized'' form where we twist the section by $\cK^{1/2}$. Likewise the Koszul dual form $\QCshear(\Mv/\Gv)\simeq \QC^!(\Ll \Xv/\Gv)$ of the local category admits a parallel coordinate independent formulation via the identification $\Ll \Xv/\Gv\simeq \Loc_\Gv^{\Xv}(D,D^*)$. 
\footnote{For simplicity we are assuming here that we are in the polarized rather than twisted-polarized case, which we also handle. Moreover, in the case when the eigencharacter $\check \eta:\Gv\to \Gm$ associated to $\Xv$ is nontrivial, we must also twist the spectral category $\QC^!(\Ll \Xv/\Gv)$ by a ``half-epsilon gerbe'' $\varepsilon_{1/2,D}$, the local analog of the epsilon line bundle correcting the normalized $L$-sheaf in \S \ref{LsheafXnormalized}, see Remark~\ref{epsilon gerbe remark} for a brief discussion.} \index{$\check\eta$} 

The compatibility between the local and global conjectures captures the relation between the four categories of sheaves involved (local and global, automorphic and spectral) as modules for the unramified Hecke operators. This compatibility has both 
{\em one-point} and {\em factorizable} formulations; this chapter will study the one-point version. 

 Fixing a point $x\in \Sigma$, we show that the Plancherel algebra $\Planch_X$ (as an algebra in the Hecke category) acts on the period sheaf via the Hecke action on the global automorphic category. This action is described in \S \ref{Theta section} as a formal consequence of the construction of the (one-point, unramified) geometric $\Theta$-series functor, a Hecke-linear functor
 $$\Theta_{X,x}:\SHV(\Bun_G^X(D,D^*))\longrightarrow \AUT(\Bun_G(\Sigma))$$ sending the basic object to the period sheaf. 
 Likewise the $L$-algebra $\bO_\Mv$ acts on the $L$-sheaf via the spectral Hecke action. This action is described analogously in \S \ref{enhance} using a spectral $\Theta$-series construction
$$\Ll_{X,x}:\QC^!(\Loc_\Gv^{\Xv}(D,D^*))\longrightarrow QC^!(\Loc_\Gv(\Sigma)).$$ 
 This naturally leads to a compatibility between the local and global conjectures. This enhanced conjecture is formulated\footnote{Again we omit for simplicity the case of twisted polarizations of $\Mv$, and refer to Remark~\ref{epsilon gerbe remark} for the normalized version for nontrivial $\check\eta$ involving epsilon lines.}
 in Conjecture~\ref{local-global conjecture}
as the construction of a Hecke-linear commutative diagram intertwining the automorphic and spectral $\Theta$-series functors
\begin{equation} \label{1pt Hecke conjecture} 
 \xymatrix{
 \AUT(\Bun_G(\Sigma))  \ar[r]^{\mathsf{GL}} & \QC^!(\Loc_\Gv(\Sigma))\\
 \SHV(\Bun_G^X(D,D^*)) \ar[u]^{\Theta_{X,x}}  \ar[r]^{\mathbb{L}_X} & \QC^!(\Loc_\Gv^{\Xv}(D,D^*)) \ar[u]^{\Ll_{\Xv,x} } 
   }
 \end{equation}
 in which the horizontal arrows are the equivalences of the local conjecture and geometric Langlands.
 
 Perhaps more concretely, the essential content of the local-global compatibility is the assertion that the global conjecture respects the actions of the algebras $\Planch_X\simeq \bO_\Mv$, which we express symbolically as an action of the bottom row on the top row in the following diagram: 

$$
\xymatrix{
{ \text{Modules:} }&\cP_X \in \SHV(\Bun_G(\Sigma))\ar[rr]^-{\mathsf{GL}} && \QC^!(\Loc_\Gv) \ni \Ll_{\Xv} \\
{\text{Algebras:}}&\Planch_X \in \HECKE_G\ar[rr]^-{\tiny \text{Satake}} && \QCshear(\fgxv/\Gv) \ni \bO_\Mv}
$$
 
\subsection{Automorphic $\Theta$ series}\label{Theta section}
In this section we introduce (the one-point, unramified version of) the $\Theta$-series functor. It is in essence
a routine transposition of the notion of $\theta$- or Poincar{\'e} series to a sheaf-theoretic context.
The discussion applies in de Rham, \'etale or Betti sheaf theories. 

\subsubsection{Setup.} Fix a point $x \in \Sigma$ on a smooth projective curve, and let $O=O_x \subset F=F_x$
be the completed local ring at $x$ and its quotient field, respectively.  We denote by $D=D_x\subset D^*=D^*_x$ the spectra of these rings, i.e., the formal disc and formal punctured disc at $x$. We omit the subscript $x$ when possible, i.e., until we begin to vary the point $x\in \Sigma$. We restrict to the case of $X$ a smooth affine $G$-variety. The local unramified theory concerns the stack $X_F/G_O$ of
$G$-bundles on the disc with a meromorphic section of the associated $X$-bundle.

{\bf Caveat: It will be important for this discussion that we assume that the $X$-spherical category $\HECKE^X$ of $!$-sheaves on $X_F/G_O$ is identified with the category of $*$-sheaves $\SHV_*(X_F/G_O)$.}

The Caveat is here since we will need the theory of $*$-pullbacks, which is only available for $*$-sheaves. Recall from Remark~\ref{! vs * conjecture} that the theories of $!$- and $*$-sheaves are expected to be equivalent on $X_F/G_O$ and specifically such an equivalence is provided by the Local Conjecture. The $*$-sheaf category is pointed by 
 the basic ($*$-)sheaf $\basic_X=\ul{k}_{X_O}$, the constant sheaf on the substack $X_O/G_O$ where the $X$-section extends to the whole disc.
 
\subsubsection{Normalizing the local category.} \index{normalized local category}
Recall that the period sheaf $\cP_X$ was defined (\S \ref{periodX}) using the stack $\Bun^X=\Bun^X_G(\Sigma)$ of $G$-bundles with a section of the associated $X$-bundle {\em twisted by a chosen spin structure} $\cK^{1/2}$. To compare this construction with the local category $\SHV(X_F/G_O)$ we need to twist the latter by a spin structure as well. This correction is not necessary when the $\GGm$-action on $X$ is trivial. In general this normalization does not affect the local category $\SHV(X_F/G_O)$ up to equivalence (even as a Hecke module), but is necessary to formulate coordinate-independent or factorizable versions of the local conjecture as well as the $\Theta$-series and local-global compatibility.
\index{$\Bun_G^X(D,D^*)$ twisted version of $X_F/G_O$}

Recall that the stack $X_F/G_O$ (respectively, $X_O/G_O$) parametrizes $G$-bundles on the disc $D$ with a section of the associated $X$-bundle on the punctured disc $D^*$ (respectively, the disc $D$).
We fix a spin structure $\cK^{1/2}$ (on the curve $\Sigma$, or for the purposes of the local category alone, just on the disc).
The group-scheme $\GGm(O)$ acts on $X_F/G_O$ and we let $$\Bun_G^X(D,D^*):=X_F/G_O\times^{\GGm(O)} \cK^{1/2}(O)^\times$$ 
be the twist of $X_F/G_O$ by the $\GGm(O)$-torsor of trivializations of $\cK^{1/2}$ over the disc. We could, equally, describe this as a twist by the induced torsor of sections of $\cK^{1/2}$ on the punctured disc for the ind-group-scheme $\GGm(F)$.

 This definition follows the global paradigm from \S \ref{bunX}, i.e.,  $\Bun_G^X(D,D^*)$ parametrizes $G$-bundles on the disc with a section of the associated $X\ot \cK^{1/2}$-bundle on the punctured disc, i.e., informally the fiber product
\begin{equation*} \xymatrix{\Bun_G^X(D,D^*)\ar[r]\ar[d] &\Map(D^*, \frac{X}{G\times \GGm})\ar[d]\\
\Bun_G(D) \ar[r]^-{\mathrm{id} \boxtimes \cK^{1/2}}&\Bun_{G \times \GGm}(D)}.\end{equation*}

\index{$\Bun_G^X(D)$} \index{$\Bun_G^X(D,D^*)$} 
Likewise we let $$\Bun_G^X(D):=X_O/G_O\times^{\GGm(O)} \cK^{1/2}(O)^\times,$$ which is a closed substack $i:\Bun_G^X(D)\to \Bun_G^X(D,D^*)$.

We define the {\em normalized local automorphic category} $$\HECKE^{X,\cK^{1/2}}:=\SHV(\Bun_G^X(D,D^*)),$$ with the basic object $\delta_X=i_*\underline{\kk}$, and considered with the normalized Hecke action of \S \ref{normalized-local}. 
\index{$\HECKE^{X,\cK^{1/2}}$ normalized local automorphic category}

\subsubsection{$\Theta$-series.} We now define the $X$-$\Theta$-series at $x$, $\Theta_{X,x}$, as a functor $$\Theta_{X,x}:\HECKE^{X,\cK^{1/2}}=\SHV(\Bun_G^X(D,D^*))\longrightarrow \SHV(\Bun_G(\Sigma))$$ 
which categorifies the classical numerical theta series  $\Phi \mapsto \sum_{x \in X_F} \Phi(xg)$.
As with the classical version, this will respect the action of $\HECKE$ by Hecke modifications at $x$.

For this let us introduce the stack 
$\Bun^X_G(\Sigma,\Sigma\setminus x)$ of $G$-bundles together with an $X \otimes \cK^{1/2}$-section {\em away from $x$},
where we use the same language as in \S \ref{bunX}. 
 We may restrict such a section to the formal neighborhood of $x$;
 fixing a trivialization of the $G$-bundle there, and so also of the $X$-bundle, we get an $F$-point of $X \otimes \cK^{1/2}$ that is well-defined up to the action of $G(O)$.
 This discussion gives the horizontal mappings in the diagram
 $$\xymatrix{\Bun_G^X(D)\ar[d]^-{i} & \ar[l]^-{\pi} \Bun^X_G(\Sigma)\ar[d]^-{i_x}\ar[dr]^-{q}&\\
\Bun^X_G(D,D^*)&\ar[l]^-{\pi_x}  \Bun^X_G(\Sigma,\Sigma\setminus x)\ar[r]^-{q_x}& \Bun_G(\Sigma)}.$$

Note that the unramified morphism $q$ is locally of finite type, and $!$-pushforward along it was used to define the period sheaf. Its ramified variant $q_x$ is of {\em ind-}finite type -- i.e., if we bound the poles at $x$ of the twisted map to $X$ we obtain finite type morphisms. As a result, the functor 
$$q_{x!}:\SHV(\Bun^X_G(\Sigma,\Sigma\setminus x))\longrightarrow \SHV(\Bun_G(\Sigma))$$ is well defined, independently of the theory of sheaves in infinite type. (In fact it agrees with the colimit preserving extension of the $!$-pushforward along the finite type closed substacks with bounded poles.) 
We may then compose with the (Hecke-linear) spectral projection (\S \ref{spectral action}, \S \ref{nilpotent projection is good}) to land in the ``automorphic'' global category:

\begin{definition} The unramified $\Theta$-series at $x$ is the composition
$$\xymatrix{\Theta_{X,x}:\HECKE^{X,\cK^{1/2}}\ar[r]^-{\pi_x^*}&\SHV(\Bun^X_G(\Sigma,\Sigma\setminus x))\ar[r]^-{(q_{x!})^{spec}}&\AUT(\Bun_G(\Sigma))}$$ 
\end{definition}

By base change on the pullback square above there is a natural identification
 \begin{eqnarray}
\label{PXThetaX} \mathcal{P}_X &\simeq & q_!\underline{k}\\
\nonumber &\simeq& q_! \pi^*\underline{k}\\
\nonumber  &\simeq& q_{x!} i_{x!}\pi^*\underline{k}\\
\nonumber  &\simeq& q_{x!} \pi_{x}^*i_{!}\pi^*\underline{k}\\
\nonumber   &\simeq& \Theta_{X,x}(\basic_X) \end{eqnarray}
between the period sheaf and the $\Theta$-series of the basic object $\basic_X\in \HECKE^{X,\cK^{1/2}}$. 
Moreover, $\Theta_{X,x}$ admits a natural $\HECKE$-linear structure,
for it ``arises from $G(O)$ invariants on a $G(F)$-equivariant diagram.''

The functor $\Theta_{X,x}$ has a normalized version \begin{equation}\label{normalized automorphic theta}
\Theta_{X, x}^{\norm}=\Theta_{X,x}\langle \deg+\beta_X \rangle 
\end{equation} 
obtained by post-composing wuth the shift $\langle \beta_X\rangle$ as well as the $\deg$-shear as in \eqref{PXnormdef}. The constant shift by $\beta_X$ commutes with Hecke actions, and the $\deg$ shift has the effect of making $\Theta_{X,x}^{\norm}$ Hecke-linear for the standard action on the target and {\em normalized} ($\deg$-shifted) action on $\SHV(X_F/G_O)$ of \S \ref{normalized-local}. There is a similar identification of the image of the basic sheaf by $\Theta_{X,x}^{\norm}$ and the normalized spectrally-projected period sheaf $(\mathcal{P}_X^{\norm})^{spec}$.

The $\Theta$-series description of the period sheaf provides it with an action of the Plancherel algebra:\index{Plancherel algebra}

\begin{prop}\label{plancherel symmetry} The identification $(\mathcal{P}_X)^{spec}\simeq \Theta_{X,x}^{\norm}(\delta_X)$ endows the spectrally projected period sheaf $(\mathcal{P}_X)^{spec}$, as an object in the $\HECKE_G$-module category $\AUT(\Bun_G(\Sigma))$ with the structure of module for the Plancherel algebra $\Planch_X\in Alg(\overline{\HECKE_G})$, the inner automorphisms of the basic object.
\end{prop}

Note a similar assertion holds before spectral projection, but we focus on the spectrally projected version for comparison with the spectral side.

\begin{remark}(Twisting by spin structures.) \index{spin structure}
We can remove the dependence of the local category $\HECKE^{X,\cK^{1/2}}$ and the $\Theta$-series on the choice of spin structure following the model of \S \ref{spindep} (whose conventions we follow). Namely we consider the stacks $\zBun_G^X(D)\hookrightarrow \zBun_G^X(D,D^*)$ of ${}^CG$-bundles on the disc with associated $\Gm$-bundle $K$ and sections (on $D$ or $D^*$) of the associated $X$-bundle. The $\Theta$-series construction defines a functor $$\SHV(\zBun_G(D,D^*))\longrightarrow \SHV(\zBun_G(\Sigma))$$ sending the basic object to the period sheaf, which is identified with the version defined above for any choice of $\cK^{1/2}$, and any two identifications differ by the translation action of $\Bun_{\Z/2}$ on $\Bun_G(\Sigma)$. 
\index{$\zBun_G$}\end{remark}

 \subsection{Spectral $\Theta$-series and local-global compatibility}\label{enhance}
We now discuss the spectral counterpart of the one-point $\Theta$-series construction of Section~\ref{Theta section}. 

Fix a smooth affine $\check{G} \times \GGm$-variety $\check{X}$. Rather than working with the symplectic space $\Mv/\Gv$ with its shearing, i.e., morally on the shifted cotangent $T^*[2]\Xv/\Gv$, we work with a Koszul dual form. Namely consider the stack $$(L \Xv)/\Gv\simeq (T[-1] \Xv)/\Gv$$ which classifies $\Gv$-local systems on the disc with a flat section of the associated $\Xv$-bundle on the punctured disc (here $L\Xv=Map(S^1,\Xv)\simeq T[-1]\Xv$ is the derived loop space of $\Xv$ \index{loop space}). This is the natural spectral counterpart of the stack $X_F/G_O$ of $G$-bundles on the disc with a section of the associated $X$-bundle on the punctured disc (as in Section~\ref{Theta section}). 

The spectral category $\QCshear(\Mv/\Gv)$ has a Koszul dual description:
we identify [ind-]coherent sheaves on $T[-1]\Xv/\Gv$, sheared to take into account the $\GGm$-action on $\Xv$,
with sheared quasicoherent  sheaves on $T^*\Xv$ by the shear of Koszul duality   
\begin{equation} \label{KoszuldualLfunctor} {\mathbb K}_\Xv=\Hom(i_*\omega_X,-):\QC^{!}(T[-1]\Xv/\Gv)\longrightarrow \QC(T^*[2]\Xv/\Gv)\end{equation}
  $${\mathbb K}_\Xv^{\shear}: \QC^{!}(T[-1]\Xv/\Gv)^{\shear} \longrightarrow \QC(T^*\Xv/\Gv)^{\shear}.$$
  On the right of the first equation, $\QC(T^* [2] \Xv/\Gv)$ really
  means that we shear $\QC(T^* \Xv/\Gv)$  by the rescaling action on the cotangent bundle.
  The shearing in the second equation additionally involves the correction from the action on the base. Koszul duality has been normalized here to identify the pushforward  $i_* \omega_\Xv$ under the zero-section with the ring of functions $\mathcal{O}_{T^*[2]\Xv}$. 

The spectral category, in either Koszul dual form, carries an action of the (spectral) spherical Hecke category. Starting from Section~\ref{spectral local category} we have been using the description of this action as the action of $\QCshear(\fgxv/\Gv)$ on $\QCshear(\Mv/\Gv)$ given by pullback along the moment map. On the other hand the spectral Hecke category has a Koszul dual description (discussed in Section~\ref{P1GL}) as the category ind-coherent sheaves on the stack $\fgv[-1]/\Gv\simeq (0\times_{\fgv} 0)/\Gv$ of pairs of local systems on the disc equipped with an identification of the punctured disc. This category is monoidal under convolution, and acts on sheaves on $L\Xv/\Gv$ by modifications of local systems at the origin. 
 
 \index{$\Theta$-series, spectral}
 
\subsubsection{The spectral $\Theta$-series} We define a functor of $\HECKE$-modules, the (one-point unramified) spectral $\Theta$-series
\begin{equation} \label{LXLoc} \Ll_{\Xv,x}:\QCshear(\Mv/\Gv)\stackrel{{\mathbb K}_\Xv^{\shear}}{\rightarrow} QC^!(L \Xv/\Gv)^{\shear} \longrightarrow QC^!(\Loc_\Gv),\end{equation}
where in the middle we have the category of ind-coherent sheaves on $L \Xv$,
but {\em sheared through the $\GGm$ action} and ${\mathbb K}_\Xv$ is the Koszul duality functor 
of~\eqref{KoszuldualLfunctor}.

 The functor $\Ll_{\Xv, x}$ is defined using the spectral counterpart of 
 $\Bun_G^X(\Sigma,\Sigma\setminus x)$ from Section~\ref{Theta section}, namely the stack $\Loc^\Xv_\Gv(\Sigma,\Sigma\setminus x)$ of $\Gv$-local systems with a flat section of the associated $\Xv$ bundle {\em away from $x\in \Sigma$}, via the resulting correspondence
 $$\xymatrix{\Loc^\Xv_\Gv(D)\ar[d]^-{i} & \ar[l]^-{\pi} \Loc^\Xv_\Gv(\Sigma)\ar[d]^-{i_x}\ar[dr]^-{q}&\\
\Loc^\Xv_\Gv(D,D^*)&\ar[l]^-{\pi_x}  \Loc^\Xv_\Gv(\Sigma,\Sigma\setminus x)\ar[r]^-{q_x}& \Loc_\Gv}$$
where we recall that $\Loc^\Xv_\Gv(D)\simeq \Xv/\Gv$ and $\Loc^\Xv_\Gv(D,D^*)\simeq T[-1]\Xv/\Gv$.

We define the spectral $\Theta$-series as $\Ll_{\Xv,x}=\unshear \circ q_{x*}^{\shear} \pi_x^{!\shear}$
(in notation similar to \eqref{Whitsymbols}; in particular, $\unshear$ is the identification of 
  $QC^{!}(\Loc_{\Gv})^{\shear}$ with $QC^{!}(\Loc_{\Gv})$  from the trivialization of the $\GGm$ action). 
This admits a natural $\HECKE$-linearity from identifying the groupoid of Hecke modifications at $x$ on $\Loc_\Gv$ as the pullback from the groupoid $\pt/\Gv \times_{\Gv/\Gv} \pt/\Gv$ on local systems on the disc. Again, there is a normalized version $\Ll_{\Xv, x}^{\norm}$ by incorporating the same $\varepsilon$ and $\beta$ twists as in \eqref{LXnormdef}. 

Just as in \eqref{PXThetaX}
  there is a natural identification $$\Ll_\Xv\simeq \Ll_{\Xv,x}(\left[ i_*\omega_X \right]^{\shear} )$$
between the $L$-sheaf of $\Xv$ and the $\Theta$-series applied to the basic object $i_*\omega_X$ (the base change calculation is identical to the identification of the period and $\Theta$-series in Section~\ref{Theta section}, with the roles of $*$ and $!$ exchanged). 
Here the notation $[\dots]^{\shear}$ means the following: $i_* \omega$ {\em a priori} belongs to the unsheared category, but we may regard it in the sheared
category because the $\Gm$-equivariant objects are identified between usual and sheared categories, as in \eqref{dataeq}. 

The Hecke-linearity of the $L$-functor implies that the image of the basic object $i_*\omega_X$ carries an action of the enriched (or inner) endomorphisms of the basic object, which is just the algebra of functions $\bO_\Mv:=\cO(\Mv/\Gv)\in \QCshear(\fgxv/\Gv)$:

\begin{corollary}\label{local L observable action}
Fix $x\in \Sigma$ and consider $QC^!(\Loc_\Gv)$ as a module category for the Hecke category $\QCshear(\fgxv/\Gv)$ through modification of local systems at $x$. Then the $L$-sheaf $\Ll_\Xv\in QC^!(\Loc_\Gv)$ carries the structure of module over the algebra $\bO_{\Mv}\in \QCshear(\fgxv/\Gv)$.
\end{corollary}

We are ready, now, to state the one-point form of local-global  compatibility in the case when both sides are polarized.
Thus we assume given a dual hyperspherical pair $(G, M=T^*X)$ and $(\check{G}, \check{M}=T^* \check{X})$.
We further assume that the spectral eigencharacter $\check\eta:\Gv\to \Gm$ is trivial (see however Remark~\ref{epsilon gerbe remark}). In this case we do not need to introduce spectral epsilon factors to normalize the spectral side;
the normalization is simply a shift
${\Ll_{\Xv,x}^{\norm}:= \Ll_{\Xv,x}\langle -\beta_\Xv \rangle}$:

\begin{conjecture}[Local-Global Conjecture, 1-point version]\label{local-global conjecture}
The equivalence $\mathbb{L}_X$ of the local conjecture (Conjecture~\ref{local conjecture}) and the geometric Langlands correspondence are intertwined by the normalized spectral and automorphic $\Theta$-series:
$$
   \xymatrix{
 \SHV(\Bun_G^X(D,D^*)) \ar[d]^{\Theta_{X,x}^{\norm}}  \ar[r]^{\mathbb{L}_X} &\QC^!(\Loc_\Gv^{\Xv}(D,D^*)) \ar[d]^{\Ll_{\Xv,x}^{\norm}}\\
 \AUT(\Bun_G(\Sigma))  \ar[r]^{\mathsf{GL}} & \QC^!(\Loc_\Gv(\Sigma))} 
$$
Moreover applied to basic objects this identification recovers the identification of period and L-sheaves given by the global conjecture (Conjecture~\ref{GlobalGeometricConjecture}).
 \end{conjecture}
 
 Equivalently, the conjecture asserts that the identification of period and $L$-sheaves is linear with respect to the identification $\Planch_X\simeq \bO_\Mv$ of Conjecture~\ref{Plancherel algebra conjecture} and the actions $\Planch_X\actson \cP_X^{\norm}$ of Proposition~\ref{plancherel symmetry} and $\bO_\Mv\actson \Ll_\Xv^{\norm}$ of Corollary~\ref{local L observable action}. (Again note that in our special case with trivial eigencharacter normalization does not affect the local spectral category at all.)

\begin{remark}[Epsilon factors and normalized local-global compatibility]\label{epsilon gerbe remark} \index{epsilon factor}
In the general case when we are given a nontrivial eigencharacter $\check\eta:\Gv\to \Gm$, we need to define the normalized form 
$$\QCshear(\Mv/\Gv)^{\norm}=\QCshear(\Mv/\Gv)\otimes_{\Rep(\Gv)} \varepsilon_{1/2,D}^\vee(\Gv)$$
of the local spectral category by tensoring with a {\em half-epsilon gerbe}, as we briefly sketch. The half-epsilon gerbe on $D^*$, an invertible sheaf of categories over $\Loc_\Gv(D^*)$, is pulled back from a half-epsilon gerbe for $\Gm$ via $\check\eta:\Loc_\Gv(D^*)\to \Loc_\Gm(D^*)$. The abelian half-epsilon gerbe in turn is constructed (by analogy with  \S \ref{spinstructure-spectral}) out of the skyscraper sheaf of categories at $\{\cK^{1/2}\}\in \Pic(D^*)$ by applying local geometric class field theory. The local and global half-epsilon factors are related by a functor 
 $$\varepsilon_{1/2,\Sigma\setminus x}^\vee:\varepsilon_{1/2,D}^\vee \longrightarrow \QC^!(\Loc_\Gv(\Sigma))$$ much as the local category and $L$-sheaf are related by the $\Theta$-series, and this functor allows us to twist the spectral $\Theta$-series to its normalized form $$\Ll_{\Xv,x}^{\norm}:\QCshear(\Mv/\Gv)^{\norm} \longrightarrow QC^!(\mathrm{\Loc}_{\check{G}}).$$ 
 \end{remark}

  \section{Automorphic factorization} \label{automorphic-factorization}
 \label{factorization algebra}
 
The theory of factorization algebras is an algebraic counterpart of the theory of $E_2$-algebras in topology  first introduced by Beilinson and Drinfeld~\cite{BD, BDchiral} to capture the commutativity of unramified Hecke operators through a mechanism of colliding points on a curve. (We refer to \S \ref{En algebras} and \S \ref{factorization algebras section} for a brief overview of $E_n$-algebras, factorization algebras and factorization homology.) 

 In particular, the spherical Hecke category $\HECKE_G$ has not only the structure of monoidal category but that of monoidal {\em factorization category}. In the constructible setting this results~\cite{E3spherical} in an $E_3$-monoidal structure on $\HECKE_G$ --a  derived weakening of the notion of symmetric monoidal category.
The spectral form of the spherical Hecke category $\QCshear(\fgxv/\Gv)$ 
correspondingly has a well-understood $E_3$-structure in the constructible setting; however,  the de Rham form and the factorizable geometric Satake correspondence are unavailable at the time of writing (but expected in upcoming work~\cite{CampbellRaskin}). 

In this section we discuss, in the de Rham and \'etale settings, the construction of a factorizable form of the local category $\SHV(X_F/G_O)$ and of the Plancherel algebra $\Planch_X$ of \S \ref{PlancherelCoulomb}. In particular, we explain that the structure on $\Planch_X$ of associative algebra object of the Hecke category can be enhanced to a locally constant factorization associative algebra object of the factorizable Hecke category on any smooth curve $\Sigma$. In the \'etale setting over $\C$ this makes  $\Planch_X$ into an $SO(2)$-fixed $E_3$-algebra object in $\HECKE_G$, also known as an associative oriented disc algebra. 

These constructions are largely variants of familiar constructions in the literature, in particular the factorization structure of loop spaces of Kapranov and Vasserot~\cite{KapranovVasserot} and the relative Coulomb branch construction (ring object) of Braverman-Finkelberg-Nakajima~\cite{BFNring}, as seen through the formalism of factorization (or chiral) categories as developed by Raskin and Gaitsgory (see~\cite{raskinchiral} and~\cite{raskinCPS2} as well as ~\cite{gaitsgoryWhatActs}, which is closely related to the construction of this section in the Eisenstein case). 

A crucial output of the theory of factorization is the mechanism of {\em factorization homology} (see \S 
\ref{fact homology sect}) which ``integrates'' or globalizes factorization algebras over $\Sigma$.
 Factorization homology appears to play a role in geometric settings analogous to the construction of Euler products (as suggested for example in~\cite{gaitsgoryAB} and private communications by J. Francis and C. Barwick). This notion is much better behaved and understood in the topological setting of $E_n$-algebras, hence in the Betti spectral setting, to which \S \ref{spectral-factorization} and \S \ref{spectral geometric quantization} are restricted, but we continue to provide an overview on both sides.

The factorization homology of the Hecke category (in both automorphic and spectral forms) produces the {\em global Hecke category} $\bH_{\Sigma}$ through which the actions of spherical Hecke functors on the global category factors. Likewise we can ``globalize'' the Plancherel algebra to produce an object that we call the {\em RTF algebra}
$\RTF_{X,\Sigma}$, which is an associative algebra object in the global Hecke category. 

We sketch the idea behind a factorizable form of the $\Theta$-series construction, conditionally on the further development of sheaf theory in infinite type.
 The existence of this factorizable $\Theta$-series implies that    
 the one-point action of the Plancherel algebra on the period sheaf descends to the RTF algebra.
  Just as the Plancherel algebra encodes maps between Hecke transforms of the basic object in the local setting, the structure on $\RTF_{X,\Sigma}$ of algebra in the global Hecke category produces maps between arbitrary Hecke functors applied to the period sheaf. 
Thus the RTF algebra plays the role of a geometric counterpart (categorification) of the Relative Trace Formula, the self-pairing of the period functional with arbitrary insertions of Hecke operators, as we discuss in \S \ref{RTF discussion}. We leave as an open problem the detailed study of the RTF algebra and its relation to the more familiar forms of the relative trace formula. 
This chapter is included as motivation and for the benefit of a more complete conjectural picture that ties in the local and global conjectures. 
 We explain the main ideas but do not verify all the  technical details; we hope these
 can be examined in a more thorough treatment of the topic.
In the next section we present the spectral counterpart of this story in the Betti setting, which is not plagued by the same technical difficulties.

 The contents are as follows:
\begin{itemize}
\item[-] \S \ref{factorization categories section} presents the setup of factorizable loop spaces, factorization categories and factorization homology.
\item[-] \S \ref{factorizable Hecke}
 discusses factorizable forms of the spherical Hecke category $\HECKE$, the $X$-spherical category $\HECKE^X$ and the Plancherel algebra $\Planch_X$.
\item[-] \S \ref{factorizable Theta} 
 introduces the RTF algebra, the factorization homology of the Plancherel algebra, and sketches the idea of the factorizable $\Theta$-series, which endows the period sheaf with an action of the RTF algebra.

\end{itemize}
 
We will make use of some notions concerning sheaves of categories, ULA objects and rigid tensor categories which are collected in Appendix
 \S \ref{shvcats}.
In Appendix~\ref{fieldtheorysec} we explain how the structures discussed in this chapter arise naturally from the algebraic formalization of boundary conditions in topological quantum field theory. 
The constructions of this section work in the de Rham and constructible (in particular \'etale)  sheaf-theoretic settings as described in Appendix~\ref{sheaf theory}.

\subsection{Factorization Categories}\label{factorization categories section}
\index{factorization category}

\subsubsection{Factorization algebras from loop spaces} \label{factorization setup}
We first discuss the basic geometric source of factorization relevant to us, the factorizable version of arc and loop spaces as introduced in~\cite{KapranovVasserot} -- see~\cite[Section 2]{raskinCPS2} for an excellent overview.
We defer to Appendix~\ref{fieldtheorysec} for background material including generalities on factorization and the Ran space.

 For the purpose of this section, $X$ can be an arbitrary smooth affine $G$-variety.
Of course, when we discuss the global conjecture, we need to make the further restrictions on $X$ used by that conjecture.

For a point $x\in \Sigma$ let $O_x\subset F_x$ denote the complete local ring of $\Sigma$ at $x$ and its field of fractions. Given an affine scheme $X$ we have an ind-scheme of loops and a subscheme of arcs $$LX_x=X(F_x)\supset LX_{+,x}=X(O_x)$$ both of infinite type over $k$. 
As $x$ varies, these spaces assemble to ind-schemes $LX_\Sigma \to \Sigma$ (and likewise for arcs). 

More generally, by ~\cite[Proposition 3.5.2]{KapranovVasserot}, given any finite set $I$ the arc and loop constructions, applied to the formal completions of finite subsets of $\Sigma$, are representable by ind-schemes (of ind-finite, respectively ind-infinite type) $LX_{+,\Sigma^I}\subset LX_{\Sigma^I}$ over $\Sigma^I$. Moreover, for $X$ smooth these multipoint arc spaces are ``pro-smooth'': $LX_{+,\Sigma^I}\to \Sigma^I$ can be represented 
as the filtered inverse limit of smooth schemes under smooth affine morphisms over $\Sigma^I$, see~\cite[Example 4.2.5]{KapranovVasserot} and ~\cite[Lemma 2.5.1]{raskinCPS2}. 

We also need a hybrid of the loop and arc constructions: given a map $I\to J$ and a $J$-tuple $\{x_j\}_{j\in J}$ of points in $\Sigma$, we may consider the ind-scheme of loops into $X$ at the points $x_j$ which are required to be arcs (i.e., integral) at the $\{x_j\}_{j\in J\setminus I}$ (where $J\setminus I$ denotes the complement of the image). This defines an ind-scheme $LX_{\Sigma^J, +:\Sigma^{J\setminus I}}$ over $\Sigma^J$ (see also ~\cite[Section 2.10]{raskinCPS2}). 
As discussed in \S \ref{units discussion}, it's convenient to extend these assignments over the category of possibly empty finite sets.

These objects carry three additional structures, which make $\{LX_{\Sigma^I}\}$ (and its arc subscheme) into a factorization algebra over $\Sigma$ in the correspondence category of ind-schemes:
\begin{enumerate}
\item[$\bullet$] {\bf Ran's Condition:} For every surjection $I\twoheadrightarrow J$ we have an isomorphism $$LX_{\Sigma^I}\times_{\Sigma^I} \Sigma^J \simeq LX_{\Sigma^J}.$$
\item[$\bullet$] {\bf Factorization:} For every decomposition $I\simeq I_1\coprod I_2$ we have an isomorphism
$$LX_{\Sigma^I}|_{U_{I_1,I_2}}\simeq [LX_{\Sigma^{I_1}} \times LX_{\Sigma^{I_2}}]|_{U_{I_1,I_2}}$$
of the restrictions to the locus ${U_{I_1,I_2}}\subset \Sigma^I$ of disjoint $I_1$- and $I_2$-tuples.
\item[$\bullet$] {\bf Unitality:} For every injection $I\hookrightarrow J$ we have a correspondence 
$$\xymatrix{& LX_{\Sigma^J, +:\Sigma^{J\setminus I}} \ar[dr]\ar[dl]  &\\
\Sigma^J\times_{\Sigma^I} LX_{\Sigma^I}& & LX_{\Sigma^J}}$$ compatible with factorization data.
\end{enumerate}

 We can summarize the conditions as saying we have an ind-scheme $LX^{fact}$ over $\Ran_{\Sigma}$ which is multiplicative, as well as an extension of this structure over the unital Ran space (i.e., replacing surjections of finite sets by arbitrary maps of possibly empty finite sets).   

Replacing $X$ by an affine group-scheme $G$ over $k$, we have corresponding versions of the loop group $LG_x, LG_{\Sigma^I}$ and $LG^{fact}$ which are group ind-schemes over $k$, $\Sigma^I$ and $\Ran_{\Sigma}$ respectively, and the subgroups $LG_+$ of arcs. If $X\onacts G$ is a $G$-variety, we obtain factorizable versions of the actions on the loop spaces $LG\actson LX$ compatible with factorization structures.

\subsubsection{Factorization Categories}
We now pass from spaces to categories of sheaves. Factorization categories are generalizations of $E_2$-monoidal categories, which themselves are derived versions of braided tensor categories. The theory of factorization categories is developed in ~\cite{raskinchiral, gaitsgoryunital} in the de Rham setting (see also~\cite{butson2} for a useful overview and applications in a context very close to ours and~\cite{E3spherical} for the factorizable Hecke category).
We do not present the fully structured $\infty$-categorical definition here, for which we refer to the above references, but only a practical snapshot thereof.
The discussion below applies equally well in constructible sheaf theories but not in the Betti setting -- we crucially use the $!$-tensor product.

\begin{definition}\label{fact category def}
A {\em faithful unital factorization category} over $\Sigma$ consists of the following data 
\begin{enumerate}
\item[$\bullet$] For every finite set $I$, we are given a sheaf of categories $\cC_{\Sigma^I}$ over $\Sigma^I$.
\item[$\bullet$] {\bf Ran's Condition:} For every surjection $\alpha:I\twoheadrightarrow J$ we have an isomorphism $$\Delta_{\alpha}^!\cC_{\Sigma^I}\simeq \cC_{\Sigma^J}.$$
\item[$\bullet$] {\bf Factorization:} For every decomposition $I\simeq I_1\coprod I_2$ we have a full embedding
$$[\cC_{\Sigma^{I_1}} \boxtimes \cC_{\Sigma^{I_2}}]|_{U_{I_1,I_2}} \hookrightarrow \cC_{\Sigma^I}|_{U_{I_1,I_2}}$$
of the restrictions to the locus ${U_{I_1,I_2}}\subset \Sigma^I$ of disjoint $I_1$- and $I_2$-tuples.
\item[$\bullet$] {\bf Unitality:} For every injection $I\hookrightarrow J$ we have a morphism $$SHV(\Sigma^J)\otimes_{SHV(\Sigma^I)} \cC_{\Sigma^I}\longrightarrow \cC_{\Sigma^J}$$ of sheaves of categories compatible with factorization data.
\end{enumerate}
The unital structure on $\cC$ is said to be ULA if the unit morphism $$u_{\Sigma^I}:SHV(\Sigma^I)\to \cC_{\Sigma^I},$$ defined by the injection of the empty set to $I$, is ULA over $\Sigma^I$ for all $I$ (i.e., has a $\SHV(\Sigma^I)$-linear continuous right adjoint).
\end{definition}

Again as noted the full unital factorization structure is best expressed as an assignment $I\to \cC_{\Sigma^I}$ over the category of possibly empty finite sets, or a multiplicative sheaf over the unital Ran space. 

Define $\omega_{\cC_{\Sigma^I}}=u_{\Sigma^I}(\omega_{\Sigma^I})$. Thanks to the ULA property of the unital structure we can consider its inner endomorphisms $$\cA_{\Sigma^I}=\underline{\mathrm{End}}(\omega_{\cC_{\Sigma^I}})=u_{\Sigma^I}^R u_{\Sigma^I}(\omega_{\Sigma^I}),$$ where $u^R$ and $u$ are linear over sheaves on $\Sigma^I$ (so in particular the construction is compatible with restriction maps) This guarantees that the $\cA_{\Sigma^I}$ form a factorizable sheaf valued in associative algebras, and in the constructible setting we may further apply Lurie's results summarized in Corollary~\ref{Lurie fact corollary}: 

\begin{proposition}\label{fact category gives E3 algebra}
Let $\cF$ denote a faithful unital factorization category over $\Sigma$, with a ULA unit. 
\begin{itemize}
\item The internal endomorphisms of the unit object $$\cA_{\Sigma^I}=\underline{End}(\omega_{\cC_{\Sigma^I}})\in SHV(\Sigma^I)$$ form a factorization associative algebra on $\Sigma$. 
\item In the constructible setting over $\C$, if we further assume that $\cA$ is locally constant, then the $!$-fibers $\cA_x$ $(x\in \Sigma)$ form an associative $E_\Sigma$-algebra, or (for $\Sigma=\AA^1$) an $SO(2)$-fixed $E_3$-algebra. 
\end{itemize}
\end{proposition}

Recall from \S \ref{factLurie} that local constancy is the property of a factorization algebra (as a $!$-sheaf on $\Ran(\Sigma)$) that its $!$-restriction to the strata of $\Ran(\Sigma)$ (configurations of $I$ distinct points) are locally constant, together with a hypercompleteness assumption (which is automatically satisfied for bounded below cochain complexes such as we will encounter).

\subsubsection{Factorization homology}\label{independent version}
The first two items in Definition~\ref{fact category def} give the definition of a (!-){\em sheaf of categories} over the Ran space $\Ran_\Sigma$ (as in~\cite{gaitsgoryunital}), and we let 
$$\cC_{\Ran_\Sigma}:=\lim_{\leftarrow, \Delta_{\alpha}^!} \cC_{\Sigma^I}$$ denote its global sections. However, since we are using $!$-sheaves, $\cC_{\Ran_\Sigma}$ behaves like a {\em homology} theory for $\Sigma$ with coefficients in $\cC$: indeed by passing to left adjoints we can rewrite this limit as a colimit over $!$-pushforwards, which we consider the {\em factorization homology} of $\cC$:
$$\int_\Sigma \cC:= \lim_{\rightarrow, \Delta_{\alpha,!}} \cC_{\Sigma^I} \simeq \cC_{\Ran_\Sigma}.$$ This is parallel to the definition of 
factorization algebras in topology (see~\S\ref{factorization algebras section}, especially Remark~\ref{cosheaf remark}) as factorizable {\em cosheaves} on the Ran space of a manifold, and their factorization homology is defined as the cosheaf homology. 

The Ran space $Ran(\Sigma)$ is homologically contractible~\cite{BD,HA}, so that $C_*\Ran(\Sigma)\simeq \kk$. However the Ran space of course still carries a large category of constructible sheaves or $D$-modules $\SHV(\Ran_\Sigma)=\lim_{\rightarrow} \SHV(\Sigma^I)$, which is itself the factorization homology of the unit factorization category.  As a result the notion $\int_\Sigma \cC$ of factorization homology for a factorization category is too large: a unital structure on $\cC$ defines a functor $u:\SHV(\Ran_\Sigma)\to \int_\Sigma \cC$, so rather than a single unit (a pointing by $\Vect$) we have a $\SHV(\Ran_\Sigma)$-worth thereof. In particular the images of skyscrapers at distinct points of $\Sigma$ will not typically be isomorphic. 
 
Therefore it is useful to refine the notion of factorization homology for unital factorization categories to a unital or ``independent'' version~\cite{gaitsgoryunital} by erasing the contribution of $\SHV(\Ran_\Sigma)$. Namely the unital structure on $\cC$ endows $\int_\Sigma \cC$ with the structure of module for the monoidal category $(\SHV(\Ran_\Sigma),\ast)$ of sheaves on $Ran(\Sigma)$ with the convolution monoidal structure (arising from the algebra structure on the Ran space in the correspondence category).
This monoidal category is augmented to $\Vect$, and the unital factorization homology $$\int_\Sigma^u \cC= \int_\Sigma \cC\ot_{(\SHV(\Ran_\Sigma),\ast)} \Vect$$ is defined as the coinvariants of the naive factorization homology. This construction takes the place for factorization categories of the (Betti) factorization homology of $E_n$-algebras, though is not as well behaved (see e.g.~\cite{darioTFT}). 
By construction $\int_\Sigma^u \cC$ is pointed by $\Vect$, in particular the image of skyscrapers at distinct points of $\Sigma$ have been identified.

\subsection{The Factorizable Plancherel Algebra} \label{factorizable Hecke}
In this section we present the factorization categories associated to loop spaces and the factorization algebras extracted from them, leading to the factorizable form of the Plancherel algebra $\Planch_X$.

\subsubsection{Sheaves on loop spaces}

Consider the factorization space $LX$ over $\Sigma$. Applying the functor $\SHV$ to the cosheaf $LX_{\Sigma^I}$ over $\Sigma^I$ we obtain a sheaf of categories $\ul{\SHV}(LX)_{\Sigma^I}$ over $\Sigma^I$. Varying $I$ 
we have an assignment  $$I\mapsto \ul{\SHV}(LX)_{\Sigma^I}$$
from finite sets to sheaves of categories, which assemble to a sheaf of categories $\ul{\SHV}(LX)^{fact}$ over the Ran space of $\Sigma$.
This assignment further satisfies the faithful form of the factorization axiom above. In the $D$-module setting, these maps are equivalences by the symmetric monoidal property of the assignment $X\mapsto \D(X)$, producing the strong notion of factorization category as it appears in~\cite{raskinchiral}.

 \medskip
 
\begin{proposition}~\cite[2.10]{raskinCPS2},~\cite[6.3]{raskininfinite}
The categories $\{\ul{\SHV}(LX_{\Sigma^I})\}$ define a unital factorization category $\SHV(LX)^{fact}$ over $\Sigma$. The unital structure is given by the factorizable basic objects $\Phi_{X,\Sigma^I}=i_*\omega_{LX_{+,\Sigma^I}}\in \SHV(LX_{\Sigma^I})$ for varying $I$, which are ULA over $\Sigma^I$.
\end{proposition}

Note the ULA property of the unit follows from the pro-smoothness of the arc space and the preservation of the ULA condition under smooth pullback and proper pushforward.

It follows from the proposition that in the constructible setting, the internal endomorphisms $\ul{\mathrm{End}}(\Phi_{X,\Sigma})^{fact}=\{\ul{\mathrm{End}}(\Phi_{X,\Sigma^I})\in \SHV(\Sigma^I)\}$ form a locally constant factorization associative algebra on $\Sigma$. (Hypercompleteness is automatic since the sheaf is bounded below, and constructibility for the stratification by diagonals is evident). 
Therefore we may apply Lurie's results from \S \ref{factLurie}. The corresponding $E_3$ algebra over $\C$ is the simply the commutative algebra of cochains on $L_+X$ (equivalently on $X$). To get a more interesting algebra we need to invoke equivariance.

\subsubsection{Factorizable spherical category}
In order to incorporate  equivariance we first recall the factorizable spherical category (see~\cite{E3spherical}). See~\cite{E3spherical} and~\cite{butson2} in particular for the notion of factorization monoidal category (associative algebra object in factorization categories).

\begin{definition}  \label{factorizable Hecke defn} The factorizable spherical category $\HECKE^{fact}$ is the unital factorization monoidal category defined by the assignment
$$I\mapsto \HECKE_{\Sigma^I}= (\SHV(LG_{+,\Sigma^I}\backslash LG_{\Sigma^I} / LG_{+,\Sigma^I}),\ast) \in Alg(\SHV(\Sigma^I))$$ of the convolution category of $LG_+$-equivariant sheaves on the Beilinson-Drinfeld affine Grassmannian, with its natural unital factorization structure. 
\end{definition}

Note that (as discussed in \S \ref{sssderivedSatake},\S \ref{renormalization section}) we use the ind-finite (or ``renormalized'') form of the spherical category~\cite{ArinkinGaitsgory}.
  
  The factorization structure on the spherical category is naturally compatible with the action of $\HECKE$ on $SHV(\Bun_G(\Sigma))$ by Hecke modifications over varying points -- in other words, these actions assemble to an action of $\int_\Sigma\HECKE$ on $SHV(\Bun_G(\Sigma))$ (which is compatible in a suitable sense with the unital structure).

 \begin{remark} \label{local constancy is hard}
The factorizable spherical category satisfies a nontrivial local constancy property~\cite{E3spherical}, which results in $\HECKE$ carrying an $E_3$ (or in general associative $E_\Sigma$) monoidal structure in the constructible setting. This local constancy uses the ind-properness of the affine Grassmannian in a fundamental way. See~\ref{E3 spherical category} for the spectral origin of this $E_3$ structure. 
 \end{remark}

\begin{lemma} 
The monoidal category $\HECKE_{\Sigma^I}$ is rigid over $\Sigma^I$.
\end{lemma}

\begin{proof}
The rigidity over $\Sigma^I$ follows, as for the standard rigidity of the spherical category $\HECKE$ over a point, from the ind-properness of the affine Grassmannian, which results in the ind-proper convolution map having a continuous right adjoint. Note that the compactness of the unit in the spherical category is a feature of working with ind-finite (renormalized) sheaves, and fails eg in $\D(LG_+\backslash LG/LG_+)$.
\end{proof}

\subsubsection{Equivariant version.}

We now consider the factorizable version of the local category $\SHV(X_F/G_O)$, by passing to the Hecke-module category of $LG_{+,\Sigma^I}$-equivariant sheaves on $LX_{\Sigma^I}$.  The rigidity of the factorizable Hecke category guarantees (through Proposition~\ref{rigid ULA}) that the ULA property of the unit (the basic sheaf $\Phi_{X,\Sigma^I}$) upgrades to the Hecke-linear setting as well, and the following proposition is a formal consequence of the setup 
(though again we do not present the details here):

\begin{prop}\label{factorizable Sph^X}
Let $X$ be a smooth affine $G$-variety.
\begin{enumerate}
\item[$\bullet$]  The assignment $$I\mapsto \HECKE_{X,\Sigma^I}:=\SHV(LG_{+,\Sigma^I}\backslash LX_{\Sigma^I})$$ extends to define a $\HECKE^{fact}$-module in faithful unital factorization categories $\HECKE^{X,fact}$ over $\Sigma$, with ULA unit given by the equivariant basic object $\Phi_{X}$.

\item[$\bullet$] The internal endomorphisms $$ \{\Planch_{X,\Sigma^I} = \ul{\mathrm{End}}(\Phi_{X,\Sigma^I})\in \HECKE_{\Sigma^I}\}$$ form a locally constant factorization associative algebra $\Planch_{X}^{fact}\in Alg(\HECKE^{fact})$ in the factorizable spherical category, the {\em factorizable Plancherel algebra}. 

\item[$\bullet$] In the constructible setting over $\C$, this endows the Plancherel algebra $\Planch_{X,x}$ with the structure of associative $E_\Sigma$-algebra object in the Hecke category, or for $\Sigma=\AA^1$ with an $SO(2)$-fixed $E_3$-algebra structure. 
\end{enumerate}
\end{prop}
\index{Plancherel algebra, factorizable}

 \begin{remark}[Normalized version] 
 Just as in the one-point case explained in \S \ref{Theta section}, 
the construction of Proposition~\ref{factorizable Sph^X} has a natural ``normalized'' modification in our setting of spherical varieties. Namely we use the $\GGm$-action on $X$ to twist the stacks $LX_{\Sigma^I}$ of maps of punctured discs on $\Sigma$ into $X$ by $\cK^{1/2}$, and also shift the Hecke action by the degree $\deg_\eta$ as in \S \ref{normalized-local}.   \end{remark}

 \begin{remark}
The local constancy of the factorizable Plancherel algebra here is not a subtle geometric property, like the local constancy of the Hecke category itself~\cite{E3spherical}, but rather a formal consequence of the same local constancy:  any section of the factorizable Hecke category over powers of the curve is automatically locally constant on the strata. We view a factorization algebra in $\HECKE$ as a (lax monoidal) functor of factorization categories from the unit $\SHV(\Ran_\Sigma)$. The former is not itself locally constant but the inclusion of locally constant categories into sheaves of categories admits a left adjoint (as in~\cite[Lemma G.1.6]{AGKRRV1}) through which any morphism to a locally constant category factors. In particular stratum by stratum a factorization algebra in $\HECKE$ is given by a functor from local systems on the stratum.  
\end{remark}

\begin{problem} \label{localconstancyHeckeX}
For $X$ a smooth affine spherical variety, and working in the constructible setting over $\C$, is the factorizable form of the $X$-spherical category locally constant? In particular this would endow $\HECKE^X$ with an $SO(2)$-fixed $E_2$ structure (i.e., upgrade it to a balanced braided tensor category).  
\end{problem}

Note that as discussed in \S \ref{Hecke local conj},~\ref{Poisson local conj} this local constancy is in fact implied by the local conjecture, since by the affineness of the spectral category (\S \ref{Spectralaffineness}) the entire Hecke category is given by modules for the Plancherel algebra, hence inherits its locally constant factorization structure. 

As discussed in \S \ref{basic examples}, there is substantial recent progress in understanding the categories $\HECKE^X$ in many examples. In particular the work on the Gaiotto conjecture (~\cite{BravermanMirabolic},~\cite{BravermanOrthosymplectic}, ~\cite{BravermanOrthosymplecticII},~\cite{TravkinYangI} and ~\cite{TravkinYangII}) provides an explicit understanding of how to construct and describe monoidal structures on $\HECKE^X$ in a series of examples, which one hopes will lead to an understanding of the factorization structures on the automorphic category in general.

 \begin{remark} (The factorizable relative Grassmannian and the product on the Plancherel algebra:) 
  Let us sketch the factorization structure on the Plancherel algebra
in a more explicit way, in particular again making clear its essentially finite-dimensional content.
We also indicate explicitly how it recovers a product structure on $\Planch_X$, which we anticipate (but don't check) agrees with the one from ~\cite{BFN}.

In  \S \ref{GXgrassdef} we introduced the {\em relative Grassmannian}
$\Gr^X$ as the subvariety of $\Gr \times X_O$
consisting, informally, of pairs $(x,g)$ with $xg \in X_O$. 
 Now, both $X_O$ and $\Gr$ extend to factorization spaces;
and the relative Grassmannian extends to a factorization subspace of 
the product space.  

Let us explicate this in the case $\Sigma =\mathbb{P}^1, 
G=\GL_n,  X =\mathbb{A}^n$; the same discussion applies
to the general case with only notational changes. 

The fiber of $X_O$
and $\Gr$ over the point $(z_1, z_2) \in \Sigma^2$ are given, respectively, 
are defined as the projective and direct limits over $N$ of: 

\begin{itemize}
\item  Points $\mathbf{x}$ of $\mathbb{A}^n$ valued in the ring
 $\C[t]/f^N$, with $f=f_{z_1, z_2} = (t-z_1) (t-z_2)$; 
\item  projective $\C[t]$-submodules $\Lambda$ of $\left( \frac{  f^{-N} \C[t]}{f^N \C[t]} \right)^n$. 
\end{itemize}

Clearly, along the diagonal $z_1=z_2$, the fiber reduces to a similarly
defined family over $\Sigma$ itself. 
Moreover, the factorization subspace corresponding to $\Gr^X$
is defined by the condition that $\mathbf{x} \in \Lambda$; 
this is clearly compatible with specialization to the diagonal. 

Now let us spell out how this constructs  the product on $\Planch_X$. 
Recall that $\Planch_X^{(V)}$ denotes the $V$-multiplicity space of the Plancherel algebra. In what follows, the role of the factorizable spaces
can be replaced by the finite dimensional versions sketched above. 

There exists
sheaves $\mathcal{T}_{V,W}$ on the factorizable
affine Grassmannian over $\Sigma \times \Sigma$
whose fibers are $T_V \boxtimes T_W$ away from the diagonal
and $T_{V \otimes W}$ at the diagonal.
Working still over $\Sigma \times \Sigma$, take the $!$-pullback of this sheaf to the factorization version of $\Gr^X$,
and then the $*$-pushforward to $\Sigma \times \Sigma$. 
Recalling the computation \eqref{shiftystuff0} of the multiplicity spaces of the Plancherel algebra,
the resulting sheaf on $\Sigma \times \Sigma$ comes with canonical identifications:
$$ \mbox{off diagonal $!$-stalks} \simeq  \Planch_{X}^{(V)} \otimes \Planch_{X}^{(W)},
\mbox{diagonal $!$-stalks} \simeq  \Planch_X^{(V \otimes W)}.$$
Then a specialization map -- which, for $!$-stalks, goes from the ``nearby'' stalk to the ``special'' stalk --  
  gives rise to the product  $$ \Planch_{X}^{(V)} \otimes \Planch_{X}^{(W)} \rightarrow \Planch_{X}^{(V \otimes W)}$$
  \end{remark}
 
\subsection{Factorizable $\Theta$-series and the RTF algebra}\label{factorizable Theta}

 \index{RTF algebra} 
Proposition~\ref{factorizable Sph^X} defines a factorizable version $\Planch_X^{fact}$ of the Plancherel algebra, an algebra object in the factorizable Hecke category. We now pass to factorization homology, i.e., compactly supported sections of $\Planch_X^{fact}$ over $\Ran_\Sigma$,
an algebra object in the monoidal category $\bH_{G,\Sigma}=\int_\Sigma \HECKE^{fact}$ given by the colimit of $\HECKE_{\Sigma^I}$ over $!$-pushforwards (see~\cite{gaitsgoryunital} for a discussion of monoidal structure on factorization homology):

 \begin{definition} \label{RTFalgebradefinition}
 We define the {\em RTF algebra} to be the factorization homology of the  
 factorization associative algebra $\Planch_{X}^{fact}$, 
 \begin{equation} \label{RTFalgebradef}\RTF_{X,\Sigma}=\int_\Sigma \Planch_X^{fact}  \in Alg(\bH_{G,\Sigma})
\end{equation} 
 \end{definition}

The constructions of the previous sections show the way to assemble the one-point $\Theta$-series $\Theta_{X,x}:\HECKE^X\to SHV(\Bun_G(\Sigma))$ as the point $x$ varies into an ``ad\`elic geometric'' object, the Ran version of the $\Theta$-series. However to do so one must face the following difficulty: the $\Theta$-series is defined by $*$-pullback on $*$-sheaves, while the factorization structure was defined using $!$-pullbacks on sheaves of categories defined using $!$-tensor products. Thus in order to carry out this construction we must assume that
\begin{itemize}
\item the categories of $*$ and $!$ sheaves on $LG_{+,\Sigma^I}\backslash LX_{\Sigma^I}$ are identified for all $I$, compatibly with $!$-pushforward along diagonal maps, and 
\item the $!$-pushforwards along diagonal maps satisfy base change with $*$-pullbacks of $*$-sheaves. 
\end{itemize}
This is beyond the scope of the current work, and we leave the following as an open problem:

\begin{problem}\label{endpoint}
Fix $X$ a smooth affine spherical $G$-variety. Show the following:
\begin{enumerate}
\item The one-point $\Theta$-series $\Theta_{X,x}$ for any $x\in\Sigma$ factors through
a $\bH_{\Sigma}=\int_\Sigma \HECKE$-linear functor 
$$\xymatrix{ \HECKE^X_x \ar[rr]^-{\Theta_{X,x}} \ar[dr]_-{i_{x,!}} && \SHV(\Bun_G(\Sigma))\\
& \int_\Sigma \HECKE^X \ar[ur]_-{\int_\Sigma \Theta_X} & }$$

\item The period sheaf $\cP_X$ admits the structure of module for the RTF algebra $\RTF_{X,\Sigma}\in Alg(\bH_{\Sigma})$ extending the action of $\Planch_{X,x}$ for fixed $x\in \Sigma$.\end{enumerate}
\end{problem}

The $\bH_\Sigma$-linear functor $act_{\cP_X}$  from $\RTF_{X,\Sigma}$-modules to the global category $\SHV(\Bun_G(\Sigma))$ can be interpreted as a coherent mechanism of constructing maps $Hom_{Bun_G}(V\ast \cP_X, W\ast \cP_X)$ between arbitrary Hecke functors applied to the period sheaf. Unlike the local situation with the Plancherel algebra, in general this will not produce {\em all} such maps, however, but rather all maps ``of local origin''. 

\begin{remark}
The compatibility between the $\Planch_{X,x}$- and $\RTF_{X,\Sigma}$-actions on the period sheaf can be expressed as a commutative diagram of pointed categories
$$\xymatrix{ \Planch_{X,x}\module_{\HECKE_x} \ar[rr]^-{act_{\cP_X}} \ar[dr]_-{i_{x,!}} && \SHV(\Bun_G(\Sigma))\\
& \RTF_{X,\Sigma}\module_{\bH_\Sigma}  \ar[ur]_-{act_{\cP_X}} & }.$$ 
Namely the pushforward of $\Planch_{X,x}$ under the functor $\HECKE_x\to \bH_\Sigma$ (insertion of a point) maps to $\RTF_{X,\Sigma}$
by the universal property of factorization homology, which allows us to compare the (pointed) categories of modules for the two algebra objects. 
\end{remark}

\begin{remark}
The spherical Hecke action on $\SHV(\Bun_G(\Sigma))$ is unital, in the sense that it factors the action of $\SHV(\Ran_\Sigma)$ on $\bH_\Sigma$ through the functor $\Gamma_c:\SHV(\Ran_\Sigma)\to \Vect$. Likewise one can ask for the factorizable $\Theta$-series to be naturally unital, i.e., to factor coinvariants of $\SHV(\Ran_\Sigma)$ to define a $\int_\Sigma^u \HECKE$-linear functor
$$\Theta_X: \int_\Sigma^u \HECKE^X \longrightarrow SHV(\Bun_G(\Sigma))$$ sending the unit ($\Vect$-pointing) to the period sheaf.  
\end{remark}

\begin{proof}[Discussion]
We sketch the main idea behind factorizable $\Theta$-series, contingent on the sheaf-theoretic hypotheses above. Let $\Bun^X_G(\Sigma^I)\to \Sigma^I$ denote the moduli stack of $G$-bundles on $\Sigma$ with a section of the associated $X$-bundle on the complement of the universal $I$-tuple of points of $\Sigma$. This stack sits in a correspondence 
$$\xymatrix{LG_{+,\Sigma^I}\backslash LX_{\Sigma^I}&\ar[l]^-{\pi_I}  \Bun^X_G(\Sigma^I)\ar[r]^-{q_I}& \Bun_G(\Sigma)},$$
and the $\Theta$-series is defined for each finite set $I$ as 
$$\Theta_{X,\Sigma^I}=q_{I!}\pi_I^*: \HECKE^X_{\Sigma^I}\to SHV(\Bun_G(\Sigma)).$$
By our sheaf-theoretic hypothesis, the $\Theta$-series commutes with the diagonal pushforward functors $\Delta_{\Sigma^I,!}$ and so  the $\Theta_{X,\Sigma^I}$ assemble to a functor out of the colimit of the $\HECKE_{X,\Sigma^I}$, i.e., the factorization homology $\int_\Sigma \HECKE^X$. Moreover this structure extends over the unital Ran space by inserting arc spaces (and the corresponding objects, the basic sheaves $\Phi$) as we saw in the one-point case:
for every injection $I\hookrightarrow J$ we have a commutative diagram with Cartesian square
$$\xymatrix{ LG_{+,\Sigma^J}\backslash LX_{\Sigma^J, +:\Sigma^{J\setminus I}}. \ar[d]^-{i} & \ar[l]^-{\pi} \Bun^X_G(\Sigma^I)\ar[d]^-{i_x}\ar[dr]^-{q}&\\
 LG_{+,\Sigma^J}\backslash LX_{\Sigma^J}&\ar[l]^-{\pi_I}  \Bun^X_G(\Sigma^J)\ar[r]^-{q_I}& \Bun_G(\Sigma)}$$
 The $\HECKE$-linearity follows as in the one-point setting, for example by lifting all the diagrams to $G(K)_{\Sigma^I}$-equivariant diagrams by picking full level structures and then passing to the quotient by $G(O)_{\Sigma^I}$.

In the second part,  the functor $\Planch_{X,x}\module_{\HECKE_x}\to \RTF_{X,\Sigma}\module_{\bH_\Sigma}$ is given as follows (see Lemma~\ref{affine observables} for a related discussion on the spectral side). The inclusion at a point defines a composite lax monoidal functor $$\xymatrix{\HECKE_x\ar[r]^-{i_{x,!}}& \HECKE_{\Sigma}\ar[r]^-{\Delta_{1,!}}& \lim_{\rightarrow,\Delta_{I,!}}\HECKE_{\Sigma^I}= \bH_\Sigma}$$ from fixed one-point to varying one-point to global Hecke categories. The image of $\Planch_{X,x}$ under the composite thereby maps to the image of $\Planch_{X,\Sigma}$ under $\Delta_{1,!}$ and thence to the colimit $ \lim_{\rightarrow,\Delta_{I,!}}\Planch_{X,\Sigma^I}=\RTF_{X,\Sigma}$. 
   \end{proof}

\subsubsection{Relation to the Relative Trace Formula}   \label{RTF discussion}\index{Relative Trace Formula}
We briefly comment on the reasoning behind the nomenclature for the RTF algebra. We refer to~\cite{YunICM} for an overview of the geometric interpretation of relative trace formulas in the function field setting, and to~\cite{FrenkelNgo} (see a review in~\cite{Frenkeltrace}) for a study of a geometric version of the Arthur and Kuznetsov trace formulas very close in spirit to our current work. 
In particular see {\it op.cit.} for the interpretation of the spectral side of these trace formulas in terms of the spectral side of geometric Langlands; the spectral analog of the RTF algebra, the L-observables, is studied in \S \ref{spectral geometric quantization}.

The relative trace formula expresses the inner product of $\Theta$-series, i.e., the inner products of Hecke operators applied to period functionals. Let us approach this problem in our unramified geometric setting, replacing the inner product by a Hom pairing (cf. Lemma \ref{Homlemma}).

Given two spherical $G$-varieties $X,Y$ we can consider the Hom space $\Hom_{\SHV(\Bun_G(\Sigma))}(\cP_X,\cP_Y)$ between the corresponding period sheaves $\cP_X$ and $\cP_Y$;
a geometric incarnation of this is the fiber product $\Bun_G^X\times_{\Bun_G}\Bun_G^Y$. More generally we can insert Hecke functors $H_{V_i}$ at finitely many points $\{x_i\}$ of $\Sigma$ and evaluate $\Hom_{\SHV(\Bun_G(\Sigma))}(\bigotimes H_{V_i,x_i} \ast \cP_X,\cP_Y)$. This Hom space is related to the stack of pairs of $G$-bundles related by a Hecke modification of prescribed type and sections of associated $X$- and $Y$-bundles (whose point count is related to the relative trace formula as explained in~\cite[2.2.4]{YunICM}).

As explained in \S \ref{case5} (in the discussion of the closely related $L^2$ version of the global conjecture), and in parallel with the characterization of the Plancherel algebra in Remark~\ref{inner endos}, these twisted Hom spaces for varying Hecke functors assemble into a single object, the inner Hom 
between the period sheaves taking values in the global Hecke category. 
(Here we should first pass from $\cP_X, \cP_Y$ to their spectral projection to  make the
action of Hecke functors vary nicely; for simplicity we will not explicitly keep track of this in the notation.)
This is, in other words, an enrichment of the usual $\Hom$ -- one recovers the usual $\Hom$
by taking morphisms from the identity object of the global Hecke category;
we will denote it by $\Hom_{\bH_{\Sigma}}(\cP_X, \cP_Y)$
and call it the ``RTF space.''

The RTF algebra $\RTF_{X,\Sigma}$ 
comes with a morphism, in the global Hecke category,
  $$\RTF_{X,\Sigma}\longrightarrow 
\Hom_{\bH_\Sigma}(\cP_X,\cP_X).$$
We regard the left-hand side as  an approximation to the RTF space in the case $X=Y$. 
More precisely, $\RTF_{X,\Sigma}$ can be considered the part of the RTF space which is ``of local origin'', i.e., comes from integrating the local version of the RTF provided by the Plancherel algebra $$\Planch_X\stackrel{\sim}{\longrightarrow}\Hom_{\HECKE}(\Phi_X,\Phi_X).$$ 
Based on the spectral description studied in \S \ref{spectral geometric quantization}, we propose the following heuristic picture: when  
we localize the story at a Langlands parameter with a unique fixed point on $\Xv$, then the RTF algebra and full RTF space should agree. This is closely related to the assertion of Conjecture~\ref{L2conj}.
 In general, the former provides the ``$\Xv$-local part'' of the later: on one hand, the RTF decomposes as a sum over automorphic representations of periods squared. On the other hand, the RTF algebra geometrizes the part of this sum which is supported on the diagonal of $\Xv$ -- the sum of squares of contributions associated to individual fixed points on $\Xv$, rather than the square of the sum.

\begin{remark}[More general RTF] To put this construction in context, it is useful to consider the following perspective: the RTF algebra is associated not just to the pair $(X,X)$ of $G$-spherical varieties but to the identity map between them. More generally, given two spherical varieties $X,Y$ {\em together with} a $G$-equivariant correspondence $X\leftarrow Z \rightarrow Y$ we can ask to quantize $Z$ to an intertwiner, a $\HECKE_G$-linear functor between the categories $\HECKE^X$ and $\HECKE^Y$ compatible with factorization. The ``inner endomorphisms'' construction then produces a $(\Planch_X,\Planch_Y)$-bimodule  $\Planch_{Z}$. Passing to factorization homology we find a $(\RTF_{X,\Sigma},\RTF_{Y,\Sigma})$-bimodule $\RTF_{Z,\Sigma}$ together with a map $$\RTF_{Z,\Sigma}\longrightarrow RTF_{X,Y}(\Sigma)=\Hom_{\bH_\Sigma}(\cP_X,\cP_Y)$$ to the global RTF space, which can be considered a contribution to the RTF from the intertwiner $Z$. This is part of the richer story of the higher category of periods suggested by the interpretation as boundary conditions, which we explore in forthcoming work. 
\end{remark}

%% file: spectral-factorization.tex
\newcommand{\Flat}{\mathrm{Flat}}
\newcommand{\weak}{weak }
\newcommand{\hvx}{\check{\mathfrak{h}}^*}
 \newcommand{\For}{\mathrm{For}}
 
 \section{Local Spectral Quantization} \label{spectral-factorization}
 
   In this section and the next we adopt the point of view that the construction of $L$-sheaves for a hyperspherical variety $\Mv$ is a problem of quantization, relative to the stack $\Loc_\Gv(\Sigma)$ of Langlands parameters on a curve. 
   
   Namely, for each local system $\rho$, we aim to produce an associative algebra (``observables'') and a module for this algebra (``states''). 
The states are the fiber $\Ll_\Mv|_\rho$ of the $L$-sheaf at $\rho$, while the observables are the fiber at $\rho$ of a sheaf of algebras we introduce, the algebra of ``$L$-observables'' $\bO_{\Mv,\Sigma}$, which provides the spectral counterpart of the RTF algebra $\RTF_{M,\Sigma}$ from the previous chapter.  We view the $L$-observables and $L$-sheaf at $\rho$ as the output of deformation quantization and geometric quantization, respectively,  applied to a symplectic variety, which for $\Mv=T^*\Xv$ polarized is simply the cotangent bundle of the stack of $\rho$-twisted maps from $\Sigma$ to $\Xv$, i.e., the (homotopy) fixed points of the Galois group or fundamental group of $\Sigma$, acting on $\Xv$ through $\rho$. 

From the perspective of number theory, the $L$-observables and its action as endomorphisms of the $L$-function provide a geometric counterpart of the $L$-function of the symplectic variety $\Mv$ and its square-root provided by the $L$-function on $\Xv$. 
See also Remark~\ref{2-shifted quantization} for a discussion from the perspective of shifted symplectic geometry, 
and Appendix~\ref{fieldtheorysec} where this problem is placed in the context of the problem of quantizing hamiltonian $\Gv$-varieties to boundary conditions for a 4d TQFT $\cB_\Gv$.

 \begin{quote}
 {\em 
 We restrict ourselves entirely to the topological setting over $\C$, so that local systems are always meant in the Betti sense, and factorization algebras are meant in the locally constant sense as $E_n$-algebras, see Appendix~\ref{factorization algebras section}. For the local considerations of \S \ref{spectral-factorization} there is no distinction between the Betti and \'etale settings over $\C$; the global considerations of \S \ref{spectral geometric quantization} are confined to the Betti setting.}
\end{quote}

In this section we focus on the {\em local} and factorizable origin of this deformation quantization problem, which we encode  by introducing the notion of a {\em
 spectral deformation quantization} of $\Mv$: a rotation-invariant $E_3$ (or factorization associative) algebra $\bO_{\Mv}$ in the spectral Hecke category, the ``$L$-algebra'', which deforms the $\cO(\fgxv)$-algebra of functions $\cO(\Mv)$. This notion recovers in particular a deformation quantization in the usual sense of $\cO(\Mv)$ as a Poisson algebra.

We shall then describe explicit spectral deformation quantizations in the polarized and twisted polarized cases.
We can combine this discussion with the factorizable Plancherel algebra of the previous section to give a factorizable form of the local conjecture (contingent on a factorizable form of the Satake correspondence, cf.~\cite{CampbellRaskin}): 
namely, there should exist a
a factorizable identification $$\Planch_X\stackrel{?}{\simeq} \bO_\Mv.$$ 
In \S \ref{spectral geometric quantization} we will use the factorization structure (via the mechanism of factorization homology) to describe the global deformation quantization obtained as the factorization homology of the $L$-algebra, $\bO_{\Mv,\Sigma}=\int_\Sigma \bO_\Mv$ and its relation to geometric quantization, namely its action on the $L$-sheaf. From the arithmetic perspective, the structure described in this section can be viewed as a subtle extra structure on the $\Gv$-representation $\cO(\Mv)$ needed to construct a geometric form of the Euler product.
  
 In more detail:
 
\begin{itemize}

\item \S \ref{E3 spherical category} discusses the factorizable form of the spectral Hecke category;
\item \S \ref{spectral quantization} gives the formal definition of a spectral deformation quantization for a given $(\check{G}, \check{M})$
and discuss its relationship to the deformation of $\check{M}$ arising from loop rotation;

 \item \S \ref{SDQcotangent} constructs spectral deformation quantizations for polarized hamiltonian varieties $\Mv=T^*\Xv$, while 
 \item \S \ref{SDQcotangent2} extends that construction to the twisted polarized setting $\Mv=T^*_\Psi\Xv$, and
 \item \S \ref{local Whittaker induction} explains how the spectral quantization of twisted cotangent bundles produces an Arthur (or spectral Whittaker) induction functor on quantized Hamiltonian spaces.
  \end{itemize}

\subsection{The factorizable spectral Hecke category}\label{E3 spherical category}
 As mentioned we will work in the simpler Betti setting over $\C$, where the tools of topological field theory (specifically factorization homology and $E_n$-algebras) allow us to give a precise and tight relation between the local and global conjectures. We refer the reader to Appendix~\ref{fieldtheorysec}, specifically Section~\ref{En algebras}, for a brief review of $E_n$-algebras and their relation to Poisson geometry. In particular we take advantage of Lurie's identification of $E_n$-algebras with locally constant factorization algebras on $\R^n$ (see Section~\ref{factLurie}) to go back and forth between the two perspectives.

\begin{terminology} \label{TerminologyFramed}
For the rest of this chapter, we will often use the term {\em factorization algebra} to denote
\begin{quote}   an $\SO(2)$-equivariant
 locally constant factorization algebra $\mathcal{A}$ on  $\R^2$,
\end{quote}
where a factorization algebra is understood in the topological sense (\S \ref{factLurie}),
i.e., a cosheaf on the Ran space of $\R^2$ with multiplicative structure, or equivalently as an $E_2$-algebra.
\end{terminology}

  Such $\mathcal{A}$ then gives rise to a locally constant factorization algebra on any oriented surface $\Sigma$,
  in a fashion that is compatible with pullback under morphisms of such.  In practice this is how we shall think of them, i.e., 
 a ``$\SO(2)$-equivariant factorization algebra'' is a compatible collection of locally constant factorization algebras on all oriented surfaces.  
From the $E_n$ perspective, an $SO(2)$-equivariant factorization algebra is an $E_2$-algebra that is invariant (in the derived sense)
 for the natural action of the rotation group on $E_2$-algebras. These are identified with the notion of ``framed $E_2$-algebra,'' i.e., an algebra over the operad of framed little discs, and give rise to algebras over the colored operad of discs in $\Sigma$ for any oriented surface $\Sigma$ (see Section~\ref{factLurie}).

We can also speak of ``factorization categories,'' which
are factorization algebras, in the sense above, now taken to be valued in $\mathsf{DGCAT}_k$.
Symmetric monoidal categories give rise to factorization categories, see e.g.~\cite{raskinCPS2}. 
A natural example is the abelian Hecke category of equivariant perverse sheaves on the Grassmannian,
which can be understood as a factorization category in the sense above
whose stalk at any point of $\R^2$ is identified with the ``usual'' abelian Hecke category.
Crucially for us, the factorization picture extends to the entire derived Hecke category $\HECKE$, as has been carried out in detail by Nocera \cite{E3spherical}. Here $\HECKE$ is as in \S \ref{sssderivedSatake} but we drop explicit mention of the group $G$.
This gives the Hecke category $\HECKE$ an $SO(2)$-fixed $E_3$ structure, as we now discuss:
 
\subsubsection{The $E_3$ structure on the Hecke category}
The convolution structure on the Hecke category is compatible with the factorization $E_2$ structure (as first observed by Lurie), making $\HECKE$ into a monoidal factorization category, or equivalently an $E_3$-category. It is important to note that unlike its abelian version this is truly a noncommutative object, i.e., is not given by a symmetric monoidal (i.e., $E_\infty$) structure.

 This gives us the ability to talk about, e.g., a ``factorization associative
algebra'' inside $\HECKE$.   
Just as
we can think of an associative algebra as (the image of the unit under) a lax monoidal functor from $\mathrm{Vect}_k$,
we can define a factorization associative algebra object   as a a lax monoidal functor of monoidal factorization categories from $\mathrm{Vect}_k$
to the factorizable Hecke category. 
Passing to a stalk at a given point, we get an algebra object of $\HECKE$ in the usual sense,
and via derived Satake we also get  a $\Gv$-equivariant algebra
over $\fgxv[2]$. By an abuse of notation we will say ``$A \in \HECKE$ is a factorization associative algebra in $\HECKE$''
in this situation, i.e., we will use the same notation 
for the factorization object and its stalk.

On the spectral side, the $E_3$ structure on the Hecke category just described is visible in 
the description of $\HECKE$ (see~\cite{ArinkinGaitsgory}), 
which already arose in \S \ref{SSP1}, as ind-coherent sheaves on the stack
$$\textrm{Map}(S^2,B\Gv)\simeq \fgv[-1]/\Gv\simeq B\Gv\times_{\fgv/\Gv} B\Gv$$ of local systems on $S^2$.  
Here the fiber product description comes from the decomposition of $S^2$ into hemispheres. 
This mapping stack description makes $\Loc_\Gv(S^2)$ into an $SO(3)$-fixed $E_3$-algebra in the correspondence category, whence a corresponding monoidal structure on ind-coherent sheaves thanks to the formalism of~\cite{GR}, as explained in~\cite{ToenBrane}.

There is not to our knowledge a published proof that
the ``automorphic'' and ``spectral'' $E_3$ structures match. We will however not be making formal use of this.

\begin{remark}  \label{Xlocallyconstant}
A crucial point in \cite{E3spherical} is the {\em local constancy}.
This arises from the following fact: if we let $\mathcal{G}r$
be the affine Grassmannian -- i.e., the space of modifications of a $G$-bundle at a point --
and define similarly  $\mathcal{G}r(D)$ to be the space of modifications of a $G$-bundle
supported inside a disc, the map $$\mathcal{G}r \rightarrow \mathcal{G}r(D)$$
is a stratified homotopy equivalence with respect to natural stratifications on both sides.
As we have noted in Problem
\ref{localconstancyHeckeX}, it's an  important 
question to clarify when the corresponding
local constancy holds for the $X$-spherical category.
\end{remark}

\begin{remark}
 Note that we can regard
 stack $\fgxv[2]/\Gv$ as the shifted cotangent bundle $T^*[3]B\Gv$, which carries a natural 3-shifted symplectic structure. The $E_3$-monoidal structure on the Hecke category provides a deformation quantization~\cite{CPTVV} of this structure on $T^*[3]B\Gv$ (see also~\cite{secondaryproducts} for a discussion).
 \end{remark}

\subsection{Spectral Deformation Quantization.}\label{spectral quantization}
We define the notion of ``spectral deformation quantization.''

\subsubsection{}
Recall that the cohomology of an $E_3$-algebra forms a graded Poisson algebra of degree $-2$ (graded $P_3$-algebra).
Now the sheared algebra $\cO^\shear({\Mv})$ has a graded Poisson bracket of degree $-2$, i.e., is a graded $P_3$ algebra, and we would like to deform it into an $E_3$ algebra, or equivalently a locally constant factorization associative algebra  on the Euclidean plane $\mathbb{R}^2$. Moreover we would like to do this compatibly with two key structures:
\begin{enumerate}
\item[$\bullet$] action of changes of coordinates ($SO(2)$-action on the plane): we would like an associative factorization algebra defined on any oriented surface $\Sigma$ , and 
\item[$\bullet$] the Hamiltonian $\Gv$-action: we would like to quantize ${\Mv}$ compatibly with the $E_3$-structure on the spherical category via the sheared moment map $\mu:{\Mv}\to \fgxv/\Gv$.
\end{enumerate}

These desiderata are captured in the following definition:  
\begin{definition} \label{sdq}
A {\em spectral deformation quantization} of the Hamiltonian $\Gv$-variety ${\Mv}$ consists of the following:  
\begin{enumerate}
\item[$\bullet$] a  locally constant factorization associative algebra object $\bO_{{\Mv}}$ in the factorizable Hecke category $\HECKE$ (see
\S \ref{E3 spherical category})
 and 
\item[$\bullet$]  an identification of the cohomology of $\bO_{{\Mv}}$ with sheared functions on $\Mv$, $$H^\ast(\bO_{{\Mv}})\simeq \cO^\shear({\Mv})$$
as Poisson algebras in $\QCshear(\fgxv/\Gv)$. 
\end{enumerate}
In this situation, the {\em ${\Mv}$-Hecke category} $\HECKE^{\Mv}=\bO_{{\Mv}}\module_{\HECKE}$ denotes the corresponding {\em quantum Hamiltonian $\Gv$-space}, by which we connote a factorization $\HECKE$-module category.
\end{definition}

 Note that we can pass between the deformation $\bO_{\Mv}$ and the deformation 
$\HECKE^{\Mv}$ of its module category: to pass from the latter to the former, we take
endomorphisms of the unit object.

 \begin{example}
 (Spectral deformation quantization in the presence of the local conjecture):
 We already saw in \S \ref{PlancherelCoulomb} that, in the context of the local unramified conjecture, $\mathcal{O}(\check{M})^{\shear}$ -- considered
as an algebra object in the Hecke category -- arises as the Plancherel algebra $\Planch_X$, 
the endomorphisms of the unit object in the $X$-spherical category. The discussion of Section~\ref{automorphic-factorization} upgrades this construction to a locally constant associative factorization algebra. Hence the Plancherel algebra is naturally an $SO(2)$-invariant $E_3$ algebra, and gives (assuming the local conjecture) a spectral deformation quantization of $\mathcal{O}(\check{M})^{\shear}$. 
\end{example}

 \begin{remark}[Comparing two forms of quantization]
\label{Poissonloop2}
A spectral deformation quantization of $\check{M}$ gives rise in particular to a deformation of the commutative algebra $\cO^\shear(\Mv)$ to an $E_3$-algebra $\bO_{\Mv}$:

This is a special case of a general construction that can be used to degenerate an object to its cohomology: the Postnikov tower construction (i.e., the t-structure on the dg category $\Vect$ of chain complexes) upgrades $\bO_{\Mv}$ to a filtered object in $\Vect$ (cf.~\cite{HA}). We then apply the Rees construction (as in~\cite{moulinos}) to obtain a ${\mathbb A}^1/\Gm$-object, with fiber at zero (the associated graded) given by the cohomology, namely $\cO^\shear(\Mv)$.

  On the other hand, there is another deformation of $\cO^{\shear}(\Mv)$ of interest.
 Namely,  $SO(2)$-equivariance gives rise to a deformation of $\cO^{\shear}(\Mv)$ over $\kk[[u]] = H^*(B \SO(2), \kk)$,
by considering derived $\mathrm{SO}(2)$-invariants (i.e., $\mathrm{SO}(2)$-equivariant cohomology). This deformation moreover acquires the structure of associative $\kk[[u]]$-algebra. This can be seen (after passing to cohomology) by considering the $SO(2)$-invariants of the $E_3$-algebra $\bO_{\Mv}$, which has one associative multiplication preserved by $SO(2)$. 

There is a strong compatibility between these two deformations of $\cO^\shear(\Mv)$, namely they give rise to the {\em same} 2-shifted Poisson bracket:

\begin{quote}
{\em Claim:} the $SO(2)$-equivariant deformation defines a deformation quantization of $\cO^\shear(\Mv)$, i.e., the associated Poisson bracket is identified with the  Poisson bracket arising from the symplectic structure. 
 \end{quote}
 
  This justifies, in other words, the definition of the Poisson bracket given in 
   \S \ref{Poisson from loop}.
The quoted can be proved by an argument of Ben-Zvi and Neitzke (a version appears in~\cite[Section 6]{secondaryproducts}, see~\cite[Proposition 25.1.1]{butson} for a more precise general version) and amounts to 
 a computation in the $\SO(2)$-equivariant homology of $S^2$,
 which is indexing binary operations on the cohomology. 
 We briefly explain the relevant identity there. 
 Consider the cohomology of $B \SO(2)$ with 
 coefficients in the homology of $S^2$. This is supported in degrees $-2,0,2,\dots$;
and as a $k[u]$ module  (where $u$ is a generator for $H^2(B\SO(2))$) is given by
taking a free module on generators $L, R$ in degree $0$ (corresponding
to the inclusion of the two fixed points) and adjoining
the class $P$ arising from the fundamental class  of $S^2$
satisfying $Pu = L-R$.  
This identity $Pu=L-R$  corresponds, after translating, to the desired identity $u \{x,y\} = xy-yx$.  \end{remark}

Return now to Definition \ref{sdq}.
For a general (shifted) symplectic variety (or Hamiltonian $\Gv$-variety)
the question of constructing a spectral deformation quantization 
poses a problem of shifted deformation quantization, as studied (and solved using formality) in~\cite{CPTVV} -- or more rigidly, invoking the grading, a problem of {\em filtered} deformation quantization (as in~\cite{Losev-deformations}). 
In the coming sections (\S \ref{SDQcotangent}, \S \ref{SDQcotangent2}) we explicitly
construct a spectral deformation quantization in all polarized and twisted polarized cases. 
The general case of our hyperspherical varieties is thereby reduced to the vectorial case, i.e., to the spectral form of the Weil representation,
much as automorphic quantization can be achieved using the theory of the Weil representation
(\S~\ref{remark:unpol}).

\subsection{Spectral quantization for cotangents} \label{SDQcotangent}
We describe explicitly how to give a spectral deformation quantization when $M=T^* X$ is a cotangent bundle
and the $\Gm$ action on $X$ is trivial:

  Let $i:\Xv\to \Ll\Xv\simeq T[-1]\Xv$ denote the inclusion of constant loops (zero section). Thanks to the Koszul resolution we find an algebra isomorphism $\mathrm{End}(i_*\omega_\Xv)\simeq \cO^\shear({\Mv})$, and we use the resulting Koszul duality equivalence
\begin{equation}\label{Koszul on X} 
Hom(i_*\omega_\Xv,-):QC^!(\Ll \Xv)\longrightarrow \QCshear({\Mv}),
\end{equation} which naturally upgrades to a monoidal equivalence with respect to convolution on the source and tensor product on the target. 

On the other hand, the monoidal category $QC^!(\Ll \Xv)$ naturally upgrades to a factorization category, since $\Ll \Xv=Map(S^1,\Xv)$ carries a framed $E_2$-structure in the correspondence category of stacks, by considering maps from complements of little discs into $\Xv$. 

This discussion can be upgraded to be $\Gv$-equivariant. The action of Hecke modifications on $\Gv$-local systems on the disc induces an action on the stack $\Ll \Xv/\Gv$ of local systems with a section of the $\Xv$-bundle on the punctured disc, i.e., the basechange of the descent groupoid in the following diagram:  
$$\xymatrix{&\Ll\Xv/\Gv\ar[r]^-{\wt{i}}\ar[d]& \Ll(\Xv/\Gv)\ar[d]^-{\pi}\\
 (\pt/\Gv)^{S^2}\ar[r]<.5ex>\ar[r]<-.5ex> & pt/\Gv \ar[r]^-{i}& \Ll pt/\Gv}.$$

Hence $QC^!(\Ll \Xv/\Gv)$ forms a factorization algebra object in modules for the factorization monoidal category $\HECKE_\Gv$.

From this we can formally deduce that the internal endomorphisms of the unit  in $QC^{!}(\mathcal{L}\check{X}/\Gv)$ defines a factorization associative algebra object $\bO_{\Mv}$ in $\HECKE_\Gv$. 
It is easy to describe the underlying associative algebra object, since as a mere module category for a monoidal category we can identify 
the $\HECKE_\Gv$-module $QC^!(\Ll\Xv/\Gv)$  with the $\QCshear(\fgxv/\Gv)$-module
 $\QCshear({\Mv}/\Gv)$. The internal endomorphisms of the unit $\cO^\shear(\Mv)\in  \QCshear({\Mv}/\Gv)$ are identified with the image
 of $\cO^\shear(\Mv)$ under the moment map $\mu:{\Mv}\to \fgxv$, so it follows that we have an equivalence $\bO_{{\Mv}}\simeq \mu_*\cO^\shear({\Mv})$ 
of algebra objects in the Hecke category.

\subsection{Spectral quantization for twisted cotangents} \label{SDQcotangent2}
Next we show how to spectrally deformation quantize twisted cotangent bundles. The construction is closely related to our calculation of spectral Whittaker $L$-sheaves on $\mathbb{P}^1$ in \S\ref{Whittaker P1}, in particular the calculation of Equation~\ref{exp via shifted skyscraper}.

Let $\Gv\actson (\Xv,\Psi)$ be as in Section~\ref{case4}. 
The twisted cotangent bundle $T^*_\Psi\Xv$ is obtained by Hamiltonian reduction by $\Ga$ from the cotangent $T^*\Psi$ of the total space of $\Psi$, i.e., we have a $\Gv$-equivariant identification $$T^*_\Psi\Xv \simeq T^*\Psi\times_{\gax}^{\Ga} \pt_1.$$ We thus will start from the spectral quantization of $T^*\Psi$ as constructed above and apply a quantized form of Hamiltonian reduction, by quantizing $\pt_1$.

 To do so it will be important to keep track of gradings. First recall that $\GGm\actson \Ga$ with weight 2. Thus the category $$\HECKE_\Ga=\QCshear(\gax/\Ga) \simeq QC(\gax)\otimes QC(\Ga^\shear)$$ has a factor given by {\em unsheared} sheaves on the affine line, and in particular contains a skyscraper object $\cO_1$ at $1\in \gax$. 

We would like to endow $\cO_1\in \HECKE_\Ga$ with the structure of factorization associative algebra quantizing its evident commutative algebra structure with respect to the symmetric monoidal structure on $\QCshear(\gax/\Ga)$. For this we note that for any group $H$ the category $\HECKE_H$ is linear over the $SO(2)$-fixed $E_4$-algebra\footnote{This is in fact commutative, and is 
 physically the ring of local operators in the field theory, i.e., functions on the Coulomb branch $(\fh^*/H)^\shear$.}
$$\End(1_{\HECKE_{H}})\simeq \cO^\shear(\fh^*)^H$$  Thus we may specialize the entire Hecke category $\HECKE_\Ga$ over $1\in \gax$, and recover $\cO_1$ as the unit in the specialized category, whence its factorization associative algebra ($SO(2)$-fixed $E_3$) structure. 

Now we suppose $\Psi\to \Xv$ is a $\Gv\times \GGm$-equivariant affine bundle, where $\GGm$ acts on $\Ga$ as above. It follows that the category $\QCshear(T^*\Psi/\Gv\times \Ga)$ has the structure of quantum Hamiltonian $\Gv\times \Ga$-space. We now need to impose the moment map condition for $\Ga$ in a factorizable fashion.

\begin{prop}\label{quantum twisted cotangents}
The category $\QCshear(\Mv/\Gv)$ associated to a twisted cotangent bundle $\Mv=T^*_\Psi\Xv$ carries a canonical structure of factorization $\HECKE_\Gv$-module through its identification with the category $$\QCshear(\Mv/\Gv)\simeq \cO_1\module(QC^!(\Ll \Psi/\Gv\times \Ga))$$
of modules for the factorization associative algebra $\cO_1\in \HECKE_{\Ga}$ in the factorization $\HECKE_{\Gv}\otimes \HECKE_\Ga$-module category 
$$\QCshear(T^*\Psi/\Gv\times \Ga)\simeq QC^!(\Ll \Psi/\Gv\times \Ga)).$$
\end{prop}

The category $\QCshear(\Mv/\Gv)$ is pointed by the sheared structure sheaf of $\Mv/\Gv$. This corresponds under Koszul duality (compare with the calculation in \S\ref{Whittaker P1}) to the unit in the factorization category $\cO_1\module(QC^!(\Ll \Psi/\Gv\times \Ga))$. This unit is given by $\cO_1\ast i_*\omega_{\Xv/\Gv}$, the action of $\cO_1\in \HECKE_\Ga$ on the unit  of $QC^!(\Ll \Psi/\Gv\times \Ga)$, itself given as the pushforward of the dualizing sheaf under the inclusion of constant loops $$i:\Xv/\Gv\simeq \Psi/\Gv\times \Ga\hookrightarrow \Ll \Psi/\Gv\times \Ga.$$

\subsection{Arthur Induction: Local Case}\label{local Whittaker induction}
We now explain another perspective on the quantization of twisted cotangent bundles: it provides the local counterpart for the Arthur induction of $L$-sheaves and local, spectral counterpart of the Whittaker induction for hamiltonian actions (Sections~\ref{Whittaker induction},~\ref{case4}, and ~\ref{Whittaker functoriality}).

As a warm-up let us consider the case when $\mathfrak{sl}_2$ is trivial, so that $\Hv \subset \Gv$ --
although we could equally well work with a homomorphism $\Hv \rightarrow \Gv$. Classically, we define symplectic induction using the Hamiltonian bi-module $T^*\Gv$ for $\Hv$ and $\Gv$. 

Equivalently (taking the point of view of Hamiltonian spaces mentioned
in Remark \ref{fancy induction}) we have the Lagrangian correspondence between the coadjoint quotients, which fits in the following commutative diagram with pullback square:
$$\xymatrix{\fgxv/\Gv\ar[d]& \fgxv/\Hv\ar[r]\ar[l]\ar[d]&\fhxv/\Hv \ar[dl]\\
\pt/\Gv&\pt/\Hv\ar[l]&}$$ 
  We may ``quantize'' this diagram to produce a factorization algebra in bimodules $\HECKE_{\Gv\gets \Hv}:=\HECKE_\Gv\ot_{\Rep(\Gv)} \Rep(\Hv)$ for $\HECKE_\Gv$ and $\HECKE_\Hv$:
$$\xymatrix{\HECKE_\Gv\ar[r]& \HECKE_{\Gv\gets\Hv}&\ar[l]\HECKE_\Hv\\
\Rep(\Gv)\ar[u]\ar[r]&\Rep(\Hv)\ar[u]\ar[ur]&}$$ 
where as before  $\HECKE_{\Hv}$ denotes the local (spectral) Hecke category for $\Hv$,  
i.e., the shear of quasi-coherent sheaves on $\fhxv/\Hv$, but now considered as a factorization monoidal category. 
As a plain category the bimodule  
$\HECKE_{\Gv\gets \Hv}$ is simply $ \QCshear(\fgxv/\Hv)$
as expected, but we have now written it in a manifestly factorizable fashion: we are applying the functoriality of $QC^!$ for correspondences to the diagram

$$\xymatrix{\Loc_\Gv(S^2)\ar[d]& \ar[r]\ar[d]\ar[l]\Loc_{\Gv,\Hv}(S^2,D)&\ar[dl] \Loc_\Hv(S^2)\\
\Loc_\Gv(D)&\ar[l]\Loc_{\Hv}(D)&}$$ 

where $\Loc_{\Gv,\Hv}(S^2,D)$ denotes $\Gv$-local systems on $S^2$ with a reduction to $\Hv$ on one hemisphere.

\begin{remark} The factorization bimodule $\HECKE_{\Gv\gets \Hv}$ is the 3-shifted form of the familiar $(\D_X,\D_Y)$-bimodule $\D_{Y\gets X}$ used to define $D$-module functoriality, where $X=pt/\Hv\to Y=pt/\Gv$.
\end{remark}

We can now define quantized symplectic induction: given a factorization $\HECKE_\Hv$-module $\cS$ (the quantized analog of a hamiltonian $\Hv$-space $S$) we have $$ind_{\Hv}^\Gv \cS:= \HECKE_{\Gv\gets \Hv}\ot_{\HECKE_\Hv} S.$$ 

The analogous definition in the Whittaker case is now clear. Fix $\Hv\times SL_2\to \Gv$ with even $SL_2$.

\begin{definition}  
The functor of Arthur induction from quantized Hamiltonian $\Hv$-spaces to quantized Hamiltonian $\Gv$-spaces is given by
 $$\mathsf{AI}(\cS):= \left[ \HECKE_{\Gv \gets \Hv\times SL_2} \otimes_{\HECKE_\Hv} \cS\right]^{\varpi\unshear},$$
 where $$\HECKE_{\Gv \gets \Hv\times SL_2}:=\cO_1\module(QC^!(\frac{\Ll \Psi}{\Gv\times \Hv\times  \Ga}))$$ is the bimodule of Proposition \ref{quantum twisted cotangents} quantizing the Hamiltonian $\Gv\times \Hv$-space $T^*_\Psi \Gv/U$.
\end{definition}

 \section{Global Spectral Quantization}\label{spectral geometric quantization}
 In this section we study the global spectral quantization of hyperspherical varieties $\Mv$. 
 This consists of three ingredients\footnote{We are indebted to Pavel Safronov for teaching us this tripartite point of view on quantization in general and specifically of shifted symplectic geometry, cf.~\cite{SafronovQuantization}.},  which can all be explicitly constructed in the polarized and twisted polarized cases $\Mv=T^*_\Psi\Xv$:
 
  \begin{itemize}
 \item[$\bullet$] geometric quantization: the $L$-sheaf $\Ll_{\Xv}$, an object of $\QC^!(\Loc_\Gv(\Sigma))$;
 \item[$\bullet$] deformation quantization: the algebra of $L$-observables $\bO_{\Mv,\Sigma}$, an algebra in the {\em global Hecke category} $\bH_\Sigma=\int_\Sigma\HECKE$ which acts on $\QC^!(\Loc_\Gv(\Sigma))$;
 \item[$\bullet$] compatibility: the action $\bO_{\Mv,\Sigma}\actson \Ll_{\Xv}$ encoded in the spectral $\Theta$-series, a functor of $\bH_\Sigma$-modules
 $$\bO_{\Mv_\Sigma}\module\longrightarrow \QC^!(\Loc_\Gv(\Sigma)).$$
 \end{itemize}
 
To carry this out, we shall review and apply the theory developed by Beraldo (building on work of Arinkin and Gaitsgory), which describes the global Hecke category on the spectral side explicitly as a refinement of $QC(\Loc_\Gv)$ which captures not only support on $\Loc_\Gv$ (Langlands parameters) but also {\em singular support} (expected to be related to Arthur parameters).
We conclude with a foray into the geometric study of Arthur parameters, whose construction is closely related to the theory of $L$-sheaves for twisted polarizations.

The material of this chapter relates in particular to three prior sections of the paper:
 \begin{itemize}
 \item[(a)] In \S \ref{Lindep}, and in particular \S \ref{specWeil-start}, we discussed (in the case $\check{X} = \mbox{ a vector space}$) how the $L$-sheaf
 in the vectorial case could be seen as a representation of a certain algebra
deforming the  $L$-sheaf not of $\check{X}$ but of $\check{M}$.
 
 \item[(b)]  In  \S \ref{case5} we examined the algebra
 of (a) in more detail over the locus $\Loc_\Gv^\circ\subset\Loc_\Gv$ where
 there is a unique fixed point on $\check{M}$ and
 saw that the deformation of (a) is simply a Clifford algebra
 deforming an exterior algebra. 
 In our current terminology,
 this Clifford algebra is (an incarnation of)
 the algebra of $L$-observables, i.e.,
 a deformation quantization;
 and the existence of
 geometric quantization is the question of Morita trivialization of the family of Clifford algebras,
 which, we anticipate, is controlled by the anomaly (Remark \ref{anomalyspectral}).

  \item[(c)]  The corresponding phenomenon in number theory is expressed by
  \eqref{L^2 normalized equation}: the endomorphisms of the period sheaf
  corresponds to the square of the period, and at the numerical level,
  there is no distinction between $\bO_{\check{M}}$ and its
   undeformed version $L_{\check{M}}$.
   \end{itemize}
However, the discussion here refines all of these prior discussions in an important way,
as we now explain -- for
$\bO_{\check{M},\Sigma}$ has 
   a structure finer than ``sheaf of algebras over $\Loc$,''
and we will now seek to construct it together with this finer structure. 
 Indeed, 
the global Hecke category provides a quantization of the $1$-shifted cotangent bundle $T^*[1]\Loc_\Gv$, a higher analog of the sheaf of differential operators,
and thus $\bO_{\check{M}, \Sigma}$
will be equipped with a structure  loosely analogous to a $D$-module and the $L$-sheaf with that of a solution.  
  See \S \ref{cohmic} for an informal discussion.

 \begin{remark}[Beyond the polarized case] 
In general quantization problems are very hard without the data of a polarization. However as in several other points in this paper we can apply the rigid structure theory of hyperspherical varieties. We describe Arthur induction functors for both the local (\S \ref{local Whittaker induction}) and global (\S \ref{WhitArth}) quantization problems, which reduce the three quantization problems above to the linear case of a symplectic representation, i.e., to the spectral analog of the theory of the Weil representation.
We will discuss this case in more detail elsewhere.  \end{remark}

 \medskip

\begin{itemize}
 
 \item \S \ref{microlocalization} reviews the notions (relative flat connections, shifted differential operators and microlocalization
of coherent sheaves) needed for global spectral quantization, mainly following Beraldo;

\item \S \ref{L-observables} shows that a spectral deformation quantization gives a sheaf of algebras -- the {\em $L$-observables} $\bO_{\check{M},\Sigma}$ -- over $\Loc$, and 
as alluded to above, 
something more {\em microlocal} -- an algebra over the Hochschild cohomology of $\Loc$. This generalizes the
Clifford algebra encountered in \S \ref{Lindep}, and should be considered the {\em global} spectral deformation quantization. The compatibility with global geometric quantization is expressed by the condition that the $L$-sheaf be a module for $L$-observables, which we interpret as a solution to a ``categorified holonomic differential equation''.  

\item \S \ref{polarized spectral} verifies the compatibility between our global spectral deformation and geometric quantizations in the polarized case. We apply a result of~\cite{HoLi} to identify the $L$-observables in this case with ``relative differential operators'' along $\Loc^\Xv\to\Loc_\Gv$. 

\item \S \ref{Whittaker L-observables} describes the modifications needed to describe $L$-observables and spectral $\Theta$-series in the twisted polarized case.

\item \S \ref{WhitArth} introduces the notion of Geometric Arthur Parameters: we study global Arthur functoriality,  a geometric construction of nontempered Langlands parameters parallel to Arthur functoriality (\S \ref{Arthur functoriality}) which combines the singular support theory of~\cite{ArinkinGaitsgory} with an eigenproperty for Hecke operators that is sheared by an Arthur $SL_2$.

 \end{itemize}

 \begin{remark}[Shifted Symplectic Geometry Perspective]\label{2-shifted quantization}
 \index{shifted symplectic geometry}
The constructions in this section and the previous fit very naturally into the framework of shifted symplectic geometry~\cite{PTVV} and its origin as the semiclassical phase space geometry of the BV-AKSZ construction of quantum field theories~\cite{AKSZ}, see~\cite{CalaqueScheimbauer}. 
We refer to ~\cite{PantevVezzosi,SafronovQuantization} for discussions of shifted deformation and geometric quantization, respectively.

As discussed in Remark~\ref{fancy induction},~\ref{SS2} the equivariant moment map  $$\mu:\Mv/\Gv\to \fgxv/\Gv$$ gives the symplectic variety $\Mv$ the structure of 1-shifted Lagrangian in $\fgxv/\Gv=T^*[1]B\Gv$. We instead want to consider the sheared version $\Mv^{\shear}$ which is a 2-shifted symplectic stack in the sense of~\cite{PTVV}; for example, in the polarized case when $\GGm$ acts trivially on $\Xv$, $\Mv^\shear=T^*[2]\Xv$ is a 2-shifted cotangent bundle.

As we have noted elsewhere (see e.g.\ Remark
~\ref{coaffine},\ref{shearing0}) it is troublesome to work with $\Mv^{\shear}$ directly as a geometric object, and in this text we have generally
``simulated'' it by carrying out constructions on $\Mv$ and shearing them. 

Ignoring this issue for this discussion, the sheared moment map $\mu^\shear:\Mv^{\shear}/\Gv\to \fgxv[2]/\Gv$ defines a shifted Lagrangian in the 3-shifted symplectic stack $\fgxv[2]/\Gv=T^*[3]B\Gv$. 
In the AKSZ formalism, quantizing spaces of maps into $\Mv^\shear$ defines a 3d TQFT (Rozansky-Witten theory), while quantizing spaces of maps into $\Mv^\shear/\Gv\to \pt/\Gv$ (coupling maps to $\Mv$ to $\Gv$-local systems) defines a boundary condition for a 4d TQFT (the 4d B-model $\cB_\Gv$). The AKSZ description of Kapustin-Witten theory in terms of $T^*[3]B\Gv$ is due to Elliott and Yoo~\cite{Elliott:2015rja} and is the basis of recent work~\cite{HilburnYoo} on S-duality for boundary theories and the analysis~\cite{elliottgwilliamwilliams} of associated factorization algebras of observables. 

The quantization of $T^*[3]B\Gv$ is the $E_3$-monoidal spectral Hecke category, while the quantization of $\mu^\shear$ is given by the spectral deformation quantization $\bO_{\Mv}$ constructed in this section (and its category of modules, the factorizable local category $\QCshear(\Mv/\Gv)$), an $E_3$-algebra object in (and $E_2$-module category for) the Hecke category.

Now we evaluate these constructions on an oriented compact surface $\Sigma$. The mapping stack $\mathrm{Map}(\Sigma,\Mv^{\shear})$ carries a natural (unshifted) symplectic structure inherited by integration over $\Sigma$ from the 2-shifted symplectic structure on $\Mv^{\shear}$. In the polarized case $\Mv=T^*\Xv$, this construction recovers the cotangent bundle of the space of maps into $\Xv$. More generally, the moment map defines a shifted Lagrangian $\mu_\Sigma$, and a relative symplectic variety $q_\Mv$:

\begin{equation} \label{musigmadef}  
\xymatrix{\Map(\Sigma, \check{M}^{\shear}/\check{G})\ar[d] \ar[rr]^-{\mu_\Sigma} \ar[drr]^-{q_{\Mv}}&& \mathrm{Map}(\Sigma, \fgxv[2]/\Gv)\simeq T^*[1] \Loc_{\check{G}}\ar[d]\\ 
\mathrm{Map}(\Sigma, \Xv/\Gv)\simeq\Loc_\Gv^\Xv\ar[rr]^-{q_\Xv}&& \mathrm{Map}(\Sigma, B\Gv)\simeq\Loc_{\Gv}
}\end{equation}

 The shifted Lagrangian $\mu_\Sigma$ is the spectral counterpart (``BBB brane'') of the Gaiotto Lagrangian ``BAA brane''~\cite{GaiottoSduality,GinzburgRozenblyum,YuLi}, the Lagrangian given by the moment map of $M$:
$$\Map(\Sigma,M/G\times^{\GGm}\cK)\longrightarrow T^*\Bun_G,$$
whose quantization is the de Rham period sheaf of $M$. 
\index{Gaiotto Lagrangian}

Given a local system $\rho$, the fiber of $q_\Mv$ defines an (unshifted) symplectic variety whose geometric and deformation quantization give the $L$-sheaf and $L$-observables at $\rho$, respectively.  As $\rho$ varies, the $L$-sheaf assembles into an object of $\QC^!(\Loc_\Gv)$, the geometric quantization of  $T^*[1]\Loc_\Gv$, and the $L$-observables assemble into an algebra in the global Hecke category $\bH_\Sigma=\int_\Sigma \HECKE$, which provides the deformation quantization of $T^*[1]\Loc_\Gv$ as described below. 
 \end{remark}

 \subsection{Relative flat connections and coherent microlocalization}\label{microlocalization}
 In this section, to provide technical background for what follows 
  we review the theory of relative flat connections, shifted differential operators and microlocalization of coherent sheaves following Beraldo~\cite{darioH,dariocenterofH}.
This theory provides a categorified or {\em shifted}\footnote{Categorification and shifting are meant to be roughly synonymous -- $n$-shifted symplectic spaces are the classical phase spaces of $(n+1)$-dimensional quantum field theories, whose topological aspects are captured by $n$-categories.} analog of the quantization of the cotangent bundle by differential operators and the corresponding theory of singular support and wave front sets of distributions.

 \subsubsection{Coherent microlocalization and ``shifted'' differential operators} \label{cohmic}

Now, on any variety $Y$, the endomorphisms of any sheaf $\Ll$ are of course an $\mathcal{O}_Y$-algebra, 
but they also have a finer structure: There is a morphism of sheaves of algebras
\begin{equation} \label{hhe} HH^*(Y) \rightarrow \mathrm{Ext}^*(\Ll,\Ll),\end{equation} 
from Hochschild cohomology of $Y$ to the derived endomorphisms of $\Ll$.
For $Y$ smooth this arises from a canonical action (the ``Atiyah class'') of the tangent bundle $\mathcal{T}_Y$
by degree $1$ endomorphisms
 $\mathcal{T}_Y \longrightarrow
   \mathrm{Ext}^1(\Ll, \Ll) $
   arising by associating to a vector field $X$ the 
self-extension of $\mathcal{L}$ arising by 
  infinitesimally displacing $\mathcal{L}$ along $X$.  

Said differently, the self-ext algebra on the right 
can be regarded as the pushforward of a sheaf on the shifted cotangent bundle $T^*Y[1]$.
This    describes a microlocal geometry for the coherent sheaf $\mathcal{L}$ \footnote{ This is an analogue
of the following process in real analysis: for a function on a smooth manifold $M$, 
 there is a standard way (the ``Wigner distribution'') to lift 
 $|f|^2$
 to a distribution on the classical phase space $T^* M$.
 This distribution expresses how the
 mass of $f$ is distributed microlocally, and bears a close relation to the action of quantum observables on $f$, in particular the module over the algebra of differential operators on $M$ generated by $f$.  } 
 and is indeed a shifted version of the usual versions of microlocalization, e.g.\ 
  the ring of functions on $T^* Y[1]$,
 carries a Poisson structure with bracket of degree $-2$,  which is the Gerstenhaber bracket
 when identified with  $HH^{\ast}(Y)$.
 
When developed in a derived setting -- replacing the graded algebra $HH^*(Y)$ by the algebra of Hochschild cochains -- it yields quite a fine invariant: endomorphisms of a sheaf that commute with the Hochschild action are automatically locally constant. For example, if $\Ll_1, \Ll_2$ are line bundles, so that $\mathrm{End}(\Ll_i) = \mathcal{O}$, coincidence of the associated morphisms \eqref{hhe} at the derived level is equivalent to $\Ll_1 \otimes \Ll_2^{-1}$
admitting a flat structure.  See Remark~\ref{center of H remark} for a more precise formulation.

Return now to the general concern of this paper.  Suppose that $\check{M}$ is polarized. In  \S \ref{case5}
we have constructed a certain sheaf of algebras $\bO_{\Mv,\Sigma}$ over the ``nonpolar locus'' $\Loc_{\check{G}}^\circ$ ---
essentially a deformation of the $L$-sheaf of $\check{M}$ itself -- 
 acting on the $L$-sheaf, 
i.e., with a morphism $$ \bO_{\check{M},\Sigma} \longrightarrow \mathrm{End}(\mathcal{L}_{\check{X}}).$$
Now, according to our discussion above, we can seek to construct
$\bO_{\check{M},\Sigma}$ {\em not merely as an $\mathcal{O}$-algebra, but as an algebra equipped with a morphism from $HH^*(\Loc_{\check{G}})$, compatibly with its action on the $L$-sheaf.}   
  In the remainder of this section we will do this (in fact over the entire stack of local systems), if $\check{M}$ has a spectral deformation quantization.  
But the technical implementatiom of all this is difficult
because  $\Loc$ is not a variety.    We will now recall the work of Beraldo~\cite{darioH,dariocenterofH}, which constructs a suitable version
of ``the category of $HH^*(\Loc)$-modules.''

\subsubsection{Relative flat connections}\label{relative flat connections} 
Let $f:  Y \to Z$ be a morphism. Then we can speak of the ``completion $\hat{Z}_Y$ of $Z$ along $f$.'' 
This is represented by the functor sending a scheme $U$
to morphisms $U \rightarrow Z$ together with a lift $U_{\mathrm{red}} \rightarrow Y$
of the induced morphism from the reduced subscheme of $U$. 
This has a hybrid behavior: if $f$ is a closed immersion, this will recover
the formal completion of $Z$ along $Y$; if $Z$ is a point, it will
recover the de Rham stack of $Y$.  In general the completion $\widehat{Z}_Y$ is
 the fiber product $\widehat{Z}_Y=Y_{dR}\times_{Z_{dR}} Z$ of $Z$ with the de Rham space of $Y$.

The category of {\em relative flat connections} $\Flat(Y/Z)=\Flat(f:Y\to Z)$  was introduced by Arinkin-Gaitsgory and Beraldo~\cite{ArinkinGaitsgory2,dariocenterofH}
 with the notation $\mathrm{IndCoh}_0(\widehat{Z}_Y)$ and is a modification of ind-coherent sheaves on the completion of $Z$ at $Y$.\footnote{We only consider representable morphisms of {\em bounded}, locally of finite presentation, perfect stacks, avoiding many of the key technicalities of~\cite{dariocenterofH} and~\cite{darioH}}. It is defined as the pullback
$$\xymatrix{
\Flat(Y/Z)\ar[r]\ar[d]&QC(Y)\ar[d]^-{\Upsilon_Y}\\
QC^!(\widehat{Z}_Y)\ar[r]^-{\widehat{f}^!} & QC^!(Y)
}$$ where $\widehat{f}:Y\to \widehat{Z}_Y$ is the completion of $Z$ along $Y$.
In other words, relative flat connections on $f$ are ind-coherent sheaves on the completion of $f$ which are quasicoherent\footnote{Our use of the term relative flat connections instead of relative $D$-modules is meant to evoke this regularity condition -- or to emphasize the role of $\Flat(Y/Z)$ as relative {\em left} $D$-modules, cf. ~\cite{dariocenterofH}.}, i.e., in the image of the symmetric monoidal functor $$\Upsilon_Y= -\ot \omega_Y:QC(Y)\to QC^!(Y),$$ when pulled back to $Y$.

Let us first illustrate two extreme instances of this notion: if $Z$ is a point we recover $D$-modules on $Y$, while if $Y$ is smooth and $f$ a closed embedding we recover (ind-coherent) sheaves on the completion of $Z$ along $Y$.

\subsubsection{Relative flat connections as modules for relative differential operators}
In general~\cite[0.3.8]{dariocenterofH} relative flat connections $\Flat(Y/Z)$ are identified with quasi-coherent sheaves on $Y$ equipped with a ``compatible'' action of the relative tangent complex $T_{Y/Z}$,
i.e., the fiber of $T_Y \rightarrow f^* T_Z$. 
Here $T_{Y/Z}$ is  identified with $\mathcal{O}_Z$-linear derivations of $\mathcal{O}_Y$  and  as such  
can be regarded as a sheaf of differential graded Lie algebras on $Y$, acting nontrivially on the structure sheaf of $Y$. 
 
Equivalently, then, we can regard $\Flat(Y/Z)$ as modules in $QC(Y)$ for  
the {\em relative differential operators} for $Y/Z$:
\begin{equation} \label{RDOdef} \cDD(f):= \mbox{
the universal enveloping algebra of the Lie algebroid $T_{Y/Z}$.}\end{equation} Just as with the ordinary sheaf of differential operators this is not an algebra in $QC(Y)$ but a monad, i.e., an algebra object in endofunctors of $\QC(Y)$, representable (via the mechanism of integral transforms) by a sheaf on $Y\times Y$ supported near the diagonal. We regard $\cDD(f)$ (or rather its pushforward to $Z$) as a deformation quantization of the relative cotangent complex $T^*_{Y/Z}$, as a relative symplectic variety over $Z$.

 \begin{remark}[Enveloping algebras]\label{enveloping alg remark}
A general definition of universal enveloping algebra in the $QC^!$ setting is given in~\cite[Volume 2, Chapter 8, Section
4.2]{GR}, as the groupoid algebra of a formal groupoid. Namely given a Lie algebroid $T$ on $Y$ 
 one has an adjoint pair $(\Ind,\For)$
of functors, induction and the forgetful functor,
 between $\QC^!(Y)$ and $T$-modules in $QC^!(Y)$.
 We then define ${\mathcal U} T:=\For\circ \Ind$ 
 to be the resulting monad, i.e., an algebra object in endofunctors of $\QC^{!}(Y)$.   One can think of ${\mathcal U} T$ geometrically as $\pi^!\pi_*$ where $\pi:Y\to Y/\exp(T)$ is the quotient by the formal groupoid corresponding to $T$. 
 
  In our setting (cf.~\cite[0.3.8]{dariocenterofH})
 there is again an adjunction $(\Ind,\For)$ between $QC(Y)$ and $\Flat(Y/Z)$, and the quasicoherent sheaf underlying the algebra of relative differential operators is $\cDD(f)=\For(\Ind(\cO_Y))$. For example:
 \begin{itemize}
 \item[-]
For $Z=\pt$ and $Y$ smooth this recovers the usual induction / forgetful adjunction between $\cO$-modules and $D$-modules, i.e., the description of $\cDD$ as the underlying $\cO$-module of the induced $D$-module $\Ind(\cO)$. 
\item[-] For $Y \rightarrow Z$ a closed embedding,  again with $Y$ smooth, the functor $\For$  $!$-restricts a sheaf on a formal neighborhood of $Y$
to a sheaf on $Y$, whereas $\Ind$ is the pushforward. Correspondingly, $\For(\Ind(\cO))$ recovers $\Hom_{\cO_Z}(\cO_Y, \cO_Y)$. 

\end{itemize}

Note that since we are taking the enveloping algebra of the {\em relative} tangent complex of $Y/Z$, the algebra $\cDD(f)$ is naturally $\QC(Z)$-linear, i.e., its pushforward to $Z$ defines an algebra object in $QC(Z)$ -- a quasicoherent sheaf of algebras -- rather than just a monad.

We note for future reference that there is a ``de Rham pushforward'' functor $f_{*,dR}:\Flat(Y/Z)\to QC^!(Z)$ (again we're assuming $f$ representable), such that the composition $f_{*,dR}\circ \Ind\simeq f_*\circ \Upsilon$ is the ind-coherent pushforward.
\end{remark}
 
A little more on the two basic examples:
\begin{itemize}
\item When $f: Y \rightarrow Z$ is a smooth morphism of smooth varieties, $\cDD(f)$
consists of differential operators on $Y$ that are ``along the fibers,'' which is to say, they commute
with multiplication by functions on $Z$; this quantizes the usual relative cotangent bundle of $f$.
\item When $f: Y \rightarrow Z$ is an LCI immersion, then $\cDD(f)$ consists, as noted above, $Z$-linear endomorphisms of $\mathcal{O}_Y$.  This can also be seen from the  Koszul dual description of sheaves on the formal completion in terms of the monad $i^!i_*$.    
 In this case the relative cotangent bundle $T^*(Y/Z)$ has fibers (when considered as a sheaf on $Z$) given by odd symplectic vector spaces $N[1] \oplus N^*[-1]$ where $N$ is the normal bundle to $Y\hookrightarrow Z$, and its quantization $\cDD(f)$ should be considered as its Clifford algebra quantization.  \end{itemize}

\subsubsection{Shifted $\D$-modules}\label{shifted D-modules}
Beraldo defines the category $\bH(Y)$ of shifted differential operators as (suitably regular) integral transforms for ind-coherent sheaves supported near the diagonal:
 
\medskip

\begin{definition}~\cite{darioH,dariocenterofH}
The monoidal category $\bH(Y)$ of shifted differential operators on a quasi-smooth stack $Y$ is the convolution category of sheaves on the formal neighborhood of the diagonal of $Y$,
$$\bH(Y):=\Flat(\Delta:Y\to Y\times Y)$$
\end{definition}

By~\cite[Corollary 3.5.3]{dariocenterofH} $\bH(Y)$ is in fact a {\em rigid} monoidal category for $Y$ a perfect stack (such as all the stacks arising in this chapter).
Locally (for $Y$ affine) the category $\bH(Y)$ can be identified Koszul dually with the category of modules for the algebra of Hochschild cochains $HC^{\ast}(Y)$ on $Y$, the self-Ext of the structure sheaf of the diagonal, i.e., with integral transforms generated by the identity. Its monoidal structure comes in this realization from the $E_2$-algebra structure on Hochschild cochains. On the level of cohomology, this $E_2$ structure recovers the shifted Poisson ($P_2$) algebra of polyvector fields, i.e., functions on $T^*[1]Y$. Thus the algebra $HC^{\ast}(Y)$ and its module category $\bH(Y)$ define deformation quantizations of the shifted symplectic stack $T^*[1]Y$.

\begin{remark} 
Note the close parallel with Grothendieck's definition of differential operators as integral transforms for functions supported set-theoretically on the diagonal.
Indeed, for $Y$ smooth, the subtleties of $QC^!$ disappear and $\bH(Y)$-modules are precisely {\em crystals of categories} on $Y$, which are sheaves of categories with a flat connection (sheaves of categories on the de Rham space $Y_{dR}$) or equivalently $\D(Y)$-module categories. In general, to illustrate the rigidity of $\bH(Y)$-actions, observe that a $\bH(Y)$-linear functor $F:QC(Y)\to QC(Y)$ is determined by the object $F(\cO)$, which acquires the structure of $D$-module, and thus if $F$ is proper (preserves compact objects) $F(\cO)$ is in fact a local system. 
For $Y$ quasi-smooth, $\bH(Y)$-modules are singular generalizations of crystals of categories and admit a theory of singular support in $T^*[1]Y$ just as $D$-modules have singular support in $T^*Y$. From this perspective crystals of categories play the role of local systems -- they are characterized by having zero singular support. 
\end{remark}

The two natural $\bH(Y)$-modules $QC(Y)$ and $QC^!(Y)$ play the role of smooth functions and distributions, and in the latter case Beraldo's theory recovers the notion of singular support of (ind-)coherent sheaves~\cite{ArinkinGaitsgory}, which is zero for quasicoherent sheaves.
Namely given a coherent sheaf $\cF\in QC^!(Y)$ we may consider its internal endomorphisms $\ul{\End}_{\bH(Y)}(\cF)$ as an algebra object in $\bH(Y)$,
which we regard as a ``chain-level microlocalization'' of $\cF$, or as the system of differential equations satisfied by $\cF$. Passing to cohomology we find its classical limit, the classical (or cohomological) microlocalization, an algebra on $T^*[1]Y$ as before whose support recovers the singular support of $\cF$.
 There is a monoidal functor $QC(Y)\to \bH(Y)$ with a continuous and lax-monoidal right adjoint $\bH(Y)\to QC(Y)$
(see~\cite[0.3.4]{dariogluing}), which means that we can ``forget'' an algebra object in $\bH(Y)$ to a quasicoherent sheaf of algebras -- in particular the algebra $\ul{\End}_{\bH(Y)}(\cF)$ refines the usual sheaf inner hom $\ul{\End}_{QC(Y)}(\cF)$, lifting it to a microlocal object.

\subsubsection{Enhanced relative differential operators} \label{subsecERDO}
 The algebra of relative differential operators $\cDD(f)$ was defined above
-- after pushforward to $Z$ -- is a sheaf of algebras on $Z$. 
 However this pushforward has a richer microlocal structure -- it naturally defines an algebra in $\bH(Z)$, as we now describe.

For a morphism $f:Y\to Z$, the category $\Flat(Y\to Z)$ of relative flat connections is itself a {\em shifted $D$-module} on $Z$, i.e., it has the structure of module for shifted differential operators $\bH(Z)$. In Lie algebra language this action is given by the map of Lie algebras $f^*T_Z[-1]\to T_{Y/Z}$ coming from the relative tangent sequence of $f$. More formally, it is a consequence of applying the functoriality of $QC^!$ to the action of the formal groupoid given by the completion of the diagonal $\Delta:Z\to Z\times Z$ on that given by the completion of $f:Y\to Z$, see also Remark~\ref{functoriality of H}.) We can use this to define an ``enhanced'' algebra of relative differential operators, 
an algebra object in $\bH(Z)$ whose underlying quasicoherent sheaf is the pushforward to $Z$ of $\cDD(f)$:

\begin{definition}\label{enhanced RDO}
The enhanced relative differential operators along $f:Y\to Z$ (for $Z$ perfect) are the algebra object $\underline{D}(f)\in Alg(\bH(Z))$ given by internal endomorphisms of $\Ind(\cO_Y)\in \Flat(Y/Z)$.
\end{definition}

\begin{remark}\label{functoriality of H} The shifted $D$-module structure $\bH(Z)\actson \Flat(Y/Z)$ is part of the rich functoriality of the theory of shifted $\D$-modules developed in~\cite{darioH} (see~\cite[0.3.4]{dariogluing} for a useful summary). Namely $\Flat(Y/Z)=\bH(\pt\gets Y\to Z)$ is identified as the ``Gauss-Manin'' sheaf of categories on $Z$, the pushforward of the structure sheaf of categories $QC(Y)$ to $Z$ in the world of categories with flat connection,
$$\Flat(Y/Z)\simeq QC(Y)\ot_{\bH(Y)} \bH(Y\to Z)\actson \bH(Z).$$
\end{remark}

\begin{remark}[The center of $\bH$ and microlocalization]\label{center of H remark}
To continue the metaphor of coherent microlocalization from \S \ref{cohmic}, we would like to understand the greater rigidity imposed by enhancing a $QC(Z)$-module to an $\bH(Z)$-module. One way to measure this is to look at the functorial endomorphisms of modules for the two monoidal categories, i.e., their centers. The center of $QC(Z)$ is identified~\cite{BZFN} with sheaves $QC(\Ll Z)$ on the derived inertia, which contains $QC(Z)$ itself. On the other hand the main theorem of~\cite{dariocenterofH} identifies the center of $\bH(Z)$ with [a subtle derived version of] $\D$-modules on the inertia $\Ll Z$. In particular if $Z$ is a scheme this is identified with $\D$-modules on $Z$ itself, and the only coherent objects in the center are flat vector bundles (local systems). 
Compare with the distinction between center of $\cO(Z)$ for $Z$ affine, which are 
(classically) functions $\cO(Z)$ or (derivedly) distributions on $\Ll Z$, and the center of $\D(Z)$
which is (classically)
 locally constant functions or (derivedly)  the de Rham cohomology of $Z$. We interpret this loosely as limiting the ambiguity in prescribing a coherent sheaf by the $\bH$-module it generates to tensoring with flat vector bundles. 
\end{remark}

\subsection{L-observables and Hecke constraints}\label{L-observables}
Having recalled the necessary technical preliminaries,
will now define the algebra of $L$-observables.  Recall that in \S \ref{case5} we defined a quasicoherent sheaf of algebras $\bO_{\Mv}^\circ$ on an open $\Loc^\circ\subset \Loc$ in the vectorial case. In that setting the relative symplectic variety of \eqref{musigmadef}
is a symplectic odd vector bundle, equivalently the data of a quadratic vector bundle and $\bO_{\Mv}^\circ$ is its deformation quantization to a Clifford algebra.
Moreover Conjecture~\ref{L2conj} (combined with the global period conjecture) identifies this algebra with the inner endomorphisms of the L-sheaf, in the cases where it has been defined. The construction of this section, which depends on being given a spectral deformation quantization, defines a sheaf of algebras on all of $\Loc$, together with a microlocal refinement (algebra in $\bH(\Loc)$). As we will see, it still acts on the $L$-sheaf, though in general we do not expect it to give the full inner endomorphism algebra.  

Fix now a smooth projective curve $\Sigma$ over $\CC$, which we consider in its Betti realization as an oriented topological surface, and $\Loc_\Gv=\Loc_\Gv(\Sigma)$ the associated stack of Langlands parameters, a quasi-smooth algebraic stack.

Recall (\S \ref{fact homology sect}) that given a locally constant factorization category $\cC$ on a surface $\Sigma$,  the {\em factorization homology}~\cite{ayalafrancis} , or topological chiral homology in the terminology of \cite{HA}, is the universal category equipped $\int_\Sigma \cC$  with functors $i_x: \cC_x\to \int_\Sigma \cC$ for $x\in \Sigma$ compatible with the unital factorization structure. If $\cC=QC(Z)$ is the symmetric monoidal category of sheaves on a perfect stack, then $\int_\Sigma \cC\simeq QC(\Map(\Sigma,Z))$ calculates sheaves on the mapping stack.

We can apply this to globalize the Hecke action.
The global spectral category $QC^!(\Loc_\Gv)$ carries an action $\HECKE_x\actson QC^!(\Loc_\Gv)$ of the spherical Hecke category for every $x\in \Sigma$. These actions 
\index{global Hecke category}  \index{$\bH_{\Sigma}$} assemble together into the action of the {\em global Hecke category} $$\bH_{\Sigma}:=\int_{\Sigma} \HECKE,$$
which was identified by Beraldo~\cite{darioTFT} with shifted differential operators on $\Loc_\Gv$:

\medskip

\begin{theorem}~\cite{darioTFT}\label{dario fact homology}
There is an equivalence of monoidal categories between the global Hecke category $\bH_{\Sigma}$, i.e., the factorization homology of the spherical category, and shifted differential operators on $\Loc_{\Gv}$,
$$\bH_{\Sigma}\simeq \bH(\Loc_{\Gv}).$$ 
In particular $\bH_\Sigma$ is rigid, and the local Hecke actions factor through the global action
$$\xymatrix{\bigotimes_{x\in \Sigma}\HECKE \ar[rr]\ar[dr] &&\End(QC^!(\Loc_\Gv))\\ & \bH_{\Sigma} \ar[ur]&}$$ 
\end{theorem}

\subsubsection{L-observables} \label{LHecke}

We are now in a position to define the ``algebra of $L$-observables,'' 
and to define abstractly a notion of $L$-sheaf in terms of it:

\begin{definition} \label{Lobsdef}
Suppose given a spectral deformation quantization $\bO_{\check{M}}$.
The {\em L-observables sheaf} 
$$ \bO_{{\Mv},\Sigma}=\int_\Sigma \bO_\Mv \in Alg(\bH_{\Sigma})$$
is the algebra object of the global Hecke category $\bH_{\Sigma}$ that is defined by factorization homology of the factorization associative algebra $\bO_{\check{M}}$ in the Hecke category.
 \end{definition}

That the $L$-observables sheaf is indeed an algebra object follows from the lax monoidal functoriality of factorization homology~\cite{ayalafrancis} 
applied to the functor from $\Vect$ associated to $\bO_{\Mv}$. 
As noted in \S \ref{shifted D-modules} that we may consider such an algebra  object in $\bH_\Sigma=\bH(\Loc)$ as a ``microlocal refinement'' of an underlying quasicoherent sheaf of algebras on $\Loc$.

\begin{remark} \label{semiclassical remark}
One way of taking a semiclassical limit here 
is by replacing the factorization structures by ``naive'' commutative multiplications:

Above we use the ``sophisticated'' factorization structure on the spherical Hecke category. 
If we instead endow it with its naive symmetric monoidal structure $(\QCshear(\fgxv/\Gv),\otimes)$ coming
from tensor product of coherent sheaves,
the factorization homology recovers the shear of $QC(\Map(\Sigma,\fgxv/\Gv))$ (since shearing commutes with factorization homology).
  The latter category  -- the semiclassical limit of the global Hecke category -- is a decompleted form of the category of sheaves on the shifted
  cotangent bundle $T^*[1]\Loc_\Gv$, where the  shift comes from the orientation class of $\Sigma$.

Similarly, the ``semiclassical limit'' of the algebra of $L$-observables
can be obtained by using 
the naive structure of $\cO_{\Mv}$ as a commutative algebra over $\fgxv/\Gv$
in place of the spectral deformation quantization. 
In this limit, 
 the $L$-observable sheaf degenerates to the structure sheaf of the ``spectral Gaiotto Lagrangian'' from Remark~\ref{2-shifted quantization}, i.e., the
pushforward of the sheared structure sheaf under the integrated moment map $\mu_\Sigma$ of~\ref{musigmadef}.   The underlying quasicoherent sheaf on $\Loc$, i.e., the pushforward $q_{\Mv,*}\cO$, is a variant of the $L$-sheaf of $\Mv$.
\end{remark}

 \subsubsection{$L$-eigensheaves and shifted differential equations}
  The rigid monoidal category $\bH_\Sigma$ acts on the spectral category $QC^!(\Loc_\Gv)$ -- an action we can interpret as either by Hecke functors or by shifted differential operators. Therefore, for an algebra object $A\in \bH_\Sigma$,
  we can consider $A$-module objects $L$ in $\QC^!(\Loc)$. This notion can be expressed equivalently as giving an algebra map
  $A\to \underline{\End}_{\bH_\Sigma}(L)$ to internal endomorphisms of $L$, or as a pointed $\bH$-linear functor 
  $$A\module_{\bH_\Sigma}\longrightarrow \QC^!(\Loc), \hskip.3in A\mapsto L,$$ and refines the more familiar notion of module for the underlying quasicoherent sheaf of algebras.

 \begin{definition}\label{L-eigensheaf}  An {\em L-eigensheaf} for ${\Mv}$ is a module object 
 in $\QC^{!}(\Loc_{\Gv})$ under the $L$-observables $\bO_{{\Mv},\Sigma}$.
  Equivalently, an $L$-eigensheaf is given by an the data of $\bH_\Sigma$-linear morphism 
 $$\Ll_{\Mv,\Sigma}:\bO_{\Mv,\Sigma}\module\longrightarrow QC^!(\Loc_\Gv).$$ We call such a functor a {\em spectral $\Theta$-series} for $\Mv$ and the category $\bH_\Sigma^\Mv:=\bO_{\Mv,\Sigma}\module$ the {\em category of $L$-observables} for $\Mv$.
\end{definition}

\index{spectral $\Theta$-series}
\index{$L$-observables} \index{$\bH^{Mv}_{\Sigma}$ $L$-observables category}\index{$L$-eigensheaf}

 We are going to verify in \S \ref{polarized spectral} that this notion of $L$-eigensheaf is compatible with the
constructions appearing earlier in the paper (i.e., ``an $L$-sheaf is an $L$-eigensheaf.'') 
The terminology ``$L$-eigensheaf'' is meant to suggest, roughly speaking, that it satisfies all Hecke constraints encoded by the L-observables. This is a less flabby notion than it might appear, or than the notion of module for the underyling quasicoherent sheaf of algebras might suggest: 
 the ambiguity in an $L$-eigensheaf
 is, roughly speaking, that of twisting by a flat bundle. 
Recall that by Theorem~\ref{dario fact homology} the global Hecke category $\bH_\Sigma$ recovers the shifted differential operators on the stack of Langlands parameters. Thus the pointed $\bH_\Sigma$-module $\bH_\Sigma^\Mv$ is prescribing a categorified system of differential equations, and an $L$-eigensheaf is precisely a solution of this system inside $QC^!(\Loc_\Gv)$. This is a categorified analog of the problem of finding a distributional solution of an algebraic system of differential equations as a map from a cyclic $D$-module. Our expectation (cf. Remark~\ref{center of H remark}) is that the eigensheaf condition loosely speaking determines the $L$-sheaf up to tensoring with flat vector bundles; it would be useful to formulate this more precisely.

 \begin{remark}[Holonomicity]
 In fact the $L$-eigensheaf condition is the shifted analog of a {\em holonomic} $\D$-module, in the following sense: the semiclassical limit of the $L$-observables form a shifted Lagrangian $\mu_\Sigma$ (~\ref{musigmadef}) of Remark~\ref{semiclassical remark} in $T^*[1]\Loc_\Gv$. It is a very interesting problem to 
 establish analogues of the strong finiteness properties of holonomic differential equations.  \end{remark}
 
 \begin{remark}[Numerical analogue]
 The analogue of the $L$-eigensheaf property in the numerical relative Langlands program is that the period function transforms
 under each local Hecke algebra as the basic function on $X_F$ does. This is not a trivial constraint; for example,
it is frequently enough to force the period function only to pair nontrivially with forms that are functorially lifted from another group.
However, it is nonetheless a much less rigid property than the categorified version that appears above. One key reason for this
is that the above notion keeps track of the derived Hecke and factorization structure. 
 \end{remark}

\subsubsection{Affineness and categorical factorization homology}
We spell out a more abstract perspective we will need in \S \ref{polarized spectral} to interface with results of~\cite{HoLi}. Namely, the given definition of $\bH^\Mv_\Sigma$ is a ``shortcut'' made possible by the affineness of $\Mv$: the module category $$\bH^\Mv_\Sigma=\bO_{\Mv,\Sigma}\module(\bH(\Loc_\Gv))\in \bH(\Loc_\Gv)\module$$ is the global counterpart to the local Hecke-module category $$\HECKE^\Mv=\bO_\Mv\module(\HECKE_\Gv)\in \HECKE_\Gv\module.$$ Specifically, $\HECKE^\Mv$ is a factorization algebra object in $\HECKE_\Gv$-modules, hence its factorization homology defines a global module $\int_{\Sigma} \HECKE^{{\Mv}}\in \bH(\Loc_\Gv)\module$. 

\begin{lemma}\label{affine observables}
The $L$-observable category $\bH^\Mv_\Sigma$ is identified as $\bH(\Loc_\Gv)$-module category with the factorization homology $\int_{\Sigma} \HECKE^{{\Mv}}$.
\end{lemma}

\begin{proof} 
The factorization homology is defined a colimit over disc embeddings of the $\HECKE_\Gv$-module categories $\HECKE^\Mv$. This is calculated as a colimit of the induced $\bH_\Sigma=\int_\Sigma \HECKE_\Gv$-module categories $\HECKE_\Mv\ot_{\HECKE_\Gv} \bH_\Sigma$. Thus we are reduced to calculating a colimit of categories of modules for algebra objects in a fixed monoidal category. But the functor $R\mapsto R\module$ taking algebras to their pointed categories of modules preserves colimits (it is the left adjoint of taking inner endomorphisms of the pointing), so the colimit category is identified with modules for the colimit algebra.
\end{proof}

\subsection{L-sheaves and L-observables in the polarized case}\label{polarized spectral}
We now examine  the polarized case $\Mv=T^*\Xv$ for a $\Gv$-variety $\Xv$,
and, in particular, the compatibility of the constructions of this section with 
the rest of the paper. 
 
 Specifically we shall  verify that the $L$-sheaf as constructed previously -- 
i.e., 
 the push-forward $L_{\Xv}(\Sigma)=q_*\omega_{\Loc^\Xv}$ of the dualizing sheaf
under \begin{equation} \label{qdef} q:\Loc^\Xv\longrightarrow \Loc_{\Gv} \end{equation} -- 
is in fact a $L$-eigensheaf, with respect to the spectral deformation quantization constructed in  \S \ref{spectral quantization}.  In fact we interpret the results of~\cite{HoLi} in our case as giving a complete description of the $L$-observables and the spectral $\Theta$-series construction in terms of the theory of relative flat connection and relative differential operators from \S \ref{microlocalization}.
 
 This is all in accordance with the semiclassical picture described in the introduction to the section.
 In the case at hand, 
 the ``relative'' symplectic variety over $\Loc_G$ described in \eqref{musigmadef}
is the  relative cotangent bundle of $\Loc^{\Xv} \rightarrow \Loc_G$;
 as we discussed, the algebra of $L$-observables is a ``relative'' deformation quantization
 of this, and it is therefore reasonable to expect that it is described by relative differential operators.

\subsubsection{Spectral $\Theta$-series in the polarized case}
We first give a direct argument for the eigensheaf property of the $L$-sheaf in the polarized case by exhibiting it as the image of a spectral $\Theta$-series functor (as in Definition~\ref{L-eigensheaf}). In other words, we present the spectral counterpart of the argument we explained for Problem~\ref{endpoint} (without the accompanying technical difficulties in the automorphic setting).
\medskip

\begin{prop}[Polarized L-sheaf]\label{polarized L functor}
The L-sheaf $L_{\Xv}(\Sigma)\in QC^!(\Loc_\Gv)$, as defined in Section~\ref{Lsheaf},
carries the structure of $L$-eigensheaf attached to $\check{M}$, 
given the spectral deformation quantization described in \S \ref{spectral quantization}.
That is to say, the $L$-sheaf carries the structure of module for the L-observables $\bO_{\Mv,\Sigma}\in Alg(\bH_{\Sigma})$.
\end{prop}
 
  \begin{proof}[Prop.~\ref{polarized L functor}]

 In Section~\ref{enhance} (and more specifically \eqref{LXLoc}) we enhanced the L-sheaf to the image of the unit under the one-point $\HECKE$-linear spectral $\Theta$-series 
 $$L_{\Xv,x}:QC^!(\Ll \Xv/\Gv)\to QC^!(\Loc_\Gv).$$ The construction is easily seen to extend to a factorizable morphism of $\HECKE$-module, which automatically gives a global spectral $\Theta$-series, hence an $L$-eigensheaf. 
 
 More precisely, given an arbitrary collection of embedded discs $\iota_I: \coprod_I D_I\hookrightarrow \Sigma$ with boundary circles $\partial \iota_I: \coprod \partial D_I\hookrightarrow \Sigma$ we consider $\Gv$-local systems on $\Sigma$ with sections of the associated $\Xv$-bundle on the complement of the discs, with its natural action of Hecke modifications of local systems and forgetful map to $\Loc_\Gv$:
 
 $$\xymatrix{(\Xv/\Gv)^I\ar[d]^-{i} & \ar[l]^-{\pi} \Loc^\Xv\ar[d]^-{i_{\iota_I}}\ar[dr]^-{q}&\\
\Loc_\Gv(\iota_I,\Xv\mbox{ on }\partial \iota_I)&\ar[l]^-{\pi_{\iota_I}}  \Loc_\Gv(\Sigma,\Xv\mbox{ on }\Sigma\setminus \iota_I)\ar[r]^-{q_{\iota_I}}& \Loc_\Gv}$$
 
This correspondence defines the $I$-fold version of the $\Theta$-series, and is easily seen to be unital, sending the pushforward of the dualizing sheaf on the closed locus $\Loc^\Xv$ to the L-sheaf, and compatible with the factorizable $\HECKE$-action. 

It follows from the universal property of factorization homology (its construction as a colimit) that these multipoint $\Theta$-series descend to the factorization homology. In particular this endows the L-sheaf with a factorizable action of the local L-observables, hence an action of $\bO_{\Mv,\Sigma}$.
 \end{proof}
 
\subsubsection{$L$-observables and relative differential operators}
 
 We now turn to the identification of $L$-observables with relative differential operators on $q$. 
 This identification will be deduced from (a special case of) a recent result of Ho and Li~\cite{HoLi}, which calculates factorization homology (and more generally associated topological field theory structure) of a wide class of Hecke categories generalizing Beraldo's description of the global Hecke category~\cite{darioTFT}. 

First recall that Lemma~\ref{affine observables} allows us to identify the category $\bH^{\Mv}:=\bO_{\Mv,\Sigma}\module$ of modules for the factorization homology of $\bO_\Mv$ with the factorization homology of the local category $\HECKE^\Mv$, which is the object calculated by~\cite{HoLi}:

\medskip 

 \begin{theorem}~\label{HoLi theorem}\cite[Theorem 3.3.1]{HoLi}
 Let $\check{M}=T^* \Xv$ equipped with the spectral deformation quantization described earlier. 
 The category of $L$-observables is identified with relative flat connections, i.e.
\begin{equation} \label{HMVHL} \bH^{\Mv}\simeq \Flat(\Loc^\Xv/ \Loc_\Gv).\end{equation} 
\end{theorem}

To check compatibility of this identification with our constructions it is useful to explain the translation between our notation and that of Ho and Li. First recall that we are using $\Flat(f:Y\to Z)$ to denote what appears in~\cite{HoLi} (as well as~\cite{ArinkinGaitsgory2,dariocenterofH}) as $IndCoh_0(\widehat{Z}_Y)$. 
We consider a morphism ${\mathcal Y}\to \cZ$ of perfect stacks of locally finite presentation. 
Given a manifold with boundary $(\partial M, M)$ (or more generally a morphism of Betti spaces or anima $N\to M$) we consider the mapping stack $(\cY,\cZ)^{(M,\partial M)}:= \cY^{\partial M}\times_{\cZ^{\partial M}} \cZ^M$ together with the morphism of stacks
$$q_{(M,\partial M)}:\cY^{\partial N} \longrightarrow (\cY,\cZ)^{(M,\partial M)}.$$

Moreover to an open embedding $N\subset M$ of manifolds with boundary, Ho and Li attach the correspondence
$$\xymatrix{
\cY^N\ar[d] & \ar[l]\ar[r] \cY^M \ar[d]& \cY^M\ar[d]\\
(\cY,\cZ)^{(N,\partial N)} & \ar[l]\ar[r] (\cY,\cZ)^{(M,M\setminus N)} & (\cY,\cZ)^{(M,\partial M)}
}$$

We specialize this construction to our setting ${\mathcal Y}:=\Xv/\Gv\to \cZ:=pt/\Gv$ and to 2-manifolds as follows:

\begin{itemize}
\item For $(M,\partial M)=(D^2,S^1)$ we obtain
$$q_{(D^2,S^1)}: \Loc^\Xv(D)=\Xv/\Gv\longrightarrow \Loc_\Gv(D,\Xv\mbox{ on }S^1)=\Ll \Xv/\Gv.$$

\item For $(M,\partial M)=(\Sigma,\emptyset)$ we obtain
$$q_{(\Sigma,\emptyset)}: \Loc^\Xv(\Sigma)\longrightarrow \Loc_\Gv(\Sigma).$$

\item For a disc embedding $\iota_I: N=\coprod_I D_I\hookrightarrow M=\Sigma$ as in Proposition~\ref{polarized L functor} we obtain the correspondence used to define the spectral $\Theta$-series,

 $$\xymatrix{(\Xv/\Gv)^I\ar[d]& \ar[l] \Loc^\Xv\ar[d]\ar[r]^-{=}& Loc^\Xv\ar[d]\\
(\Ll \Xv/\Gv)^I &\ar[l]  \Loc_\Gv(\Sigma,\Xv\mbox{ on }\Sigma\setminus N)\ar[r] & \Loc_\Gv.}$$

\end{itemize}  

\cite{HoLi} then applies the functoriality of $IndCoh_0$ under correspondences from~\cite{darioH} to construct a functor from the category of manifolds with open embeddings to dg categories, which on objects attaches $M\mapsto \Flat(q_{(M,\partial M)})$. 
The main result of~\cite{HoLi} identifies the category attached to any $n$-manifold with the factorization homology of the $E_n$-category obtained by restricting the functor to embedding of discs $(D^n,\partial D^n)$.

When specialized to our setting, we find:

\begin{itemize}
\item For $(M,\partial M)=(D^2,S^1)$ we obtain an $E_2$-category ($H_2({\mathcal Y},\cZ)$ of {\it op. cit.} 3.2), which evaluates to the local category
$$\Flat(q_{(D^2,S^1)}: \Xv/\Gv\to \Ll \Xv/\Gv)\simeq \HECKE^{\Mv}=QC^!(\Ll \Xv/\Gv)$$
since $\Xv/\Gv\hookrightarrow \Ll \Xv/\Gv$ is a nil-isomorphism from a smooth stack, so the relative construction simply encodes $QC^!$ on the target.

\item For $(M,\partial M)=(\Sigma,\emptyset)$ we obtain the global ``Eisenstein'' category of {\it op. cit.} as the factorization homology of the Hecke category, which evaluates to the description of the global category in Theorem~\ref{HoLi theorem},
$$\bH^\Mv_\Sigma=\int_\Sigma \HECKE^\Mv\stackrel{\sim}{\longrightarrow} \Flat(q_{(\Sigma,\emptyset)}: \Loc^\Xv\to \Loc_\Gv).$$

\item The pointing of $\bH^\Mv_\Sigma$ is given by applying the functorality of $\Flat$ to the unit of $\HECKE^\Mv=QC^!(\Ll \Xv/\Gv)$, i.e., the pushforward of the dualizing sheaf on $\cY^D=\Xv/\Gv$. We find the pushforward of $\cO_{\Loc^{\Xv}}$ under $$\Flat(\Loc^\Xv/\Loc^{\Xv})=QC(\Loc^\Xv)\longrightarrow \Flat(\Loc^\Xv/\Loc^\Gv),$$ i.e., the sheaf of relative differential operators. 
\end{itemize}  

In particular we highlight the following additional consequences of the general construction of~\cite{HoLi}, which identify the L-observables, the global $\Theta$-series, the L-sheaf and its eigenproperty in terms of relative differential operators: 

\medskip

\begin{corollary}\label{HoLi corollary}
\begin{itemize}
\item The algebra of $L$-observables $\bO_{\Mv,\Sigma}\in \mathrm{Alg}(\bH_\Sigma)$ (as endomorphisms of the pointing in $\bH^\Mv_\Sigma$) is identified with the enhanced algebra of relative differential operators $\underline{\D}(q)$ along $q:\Loc^\Xv\to \Loc_\Gv$ (cf. Definition~\ref{enhanced RDO}).
\item The spectral $\Theta$-series $\bH^{\Mv}_\Sigma\to QC^!(\Loc_\Gv)$ 
of Definition \ref{L-eigensheaf} factors, with reference to the identification
\eqref{HMVHL}, 
 through the relative de Rham pushforward $q_{*,dR}:\Flat(\Loc^\Xv/\Loc_\Gv)\to QC^!(\Loc_\Gv)$ of Remark~\ref{enveloping alg remark}(and similarly for the factorizable version). 
 \item Under this identification the $L$-sheaf $q_*\omega_{\Loc^\Xv}=q_*\Upsilon(\cO_{\Loc^\Xv})$ is identified with the relative de Rham pushforward $q_{*,dR}$ of the induced $\D$-module of relative differential operators $\mathrm{Ind}(\cO_{\Loc^\Xv})\in \Flat(\Loc^\Xv/\Loc_\Gv)$ .
\end{itemize}
\end{corollary}

\subsection{The twisted polarized case}\label{Whittaker L-observables}
We briefly discuss the modifications necessary to carry out the construction of one-point and factorizable spectral $\Theta$-series and the description of $L$-observables in the case of twisted cotangent bundles $\Mv=T^*_\Psi\Xv$. 

We follow the setup and approach of \S \ref{SDQcotangent2}.
Namely we first replace the role of the $\Gv$-space $\Xv$ by that of the $\Gv\times \Ga$-space
$\Psi\to \Xv$ to obtain one-point and factorizable $\Theta$-series functors
$$\xymatrix{QC^!(\Ll \Psi/\rG\times \Ga)\ar[rr]\ar[dr]&&  QC^!(Loc_\Ga)\otimes QC^!(\Loc_\rG)\\
&\bH_{T^*\Psi}\ar[ur]&}$$
as in \S \ref{enhance} and \S \ref{polarized spectral}. Moreover the factorization homology $\bH_{T^*\Psi}$ is again identified with relative flat connections for $\Loc^{\Psi}$ over $\Loc_{\Ga}\times \Loc_\Gv$. The unit in $\HECKE^{\Psi}$ is $i_*\omega_{\Psi/\rG\times \Ga}$, and maps to the $L$-sheaf on $\Loc_\Ga\times \Loc_\Gv$. 

We now change perspective and consider this $L$-sheaf as representing the integral transform
$$\xymatrix{& \Loc_{\Gv\times \Ga}^\Psi=\Loc^\Xv_{\Gv} \ar[dl]\ar[dr]& \\
\Loc_\Ga&& \Loc_\Gv}$$
which takes the spectral exponential sheaf on $\Loc_\Ga$ to the spectral Whittaker $L$-sheaf as in \S \ref{case4}.
This construction is linear for the (one-point or factorizable) action of the Hecke category for $\Ga$.
Thus we may now pass everywhere to modules for the factorization algebra $\cO_1\in \HECKE_\Ga$ (the quantum version of imposing the moment map value $1$ for the Hamiltonian $\Ga$-action on $T^*\Psi$).
This operation takes $\HECKE^\Psi$ to the factorization $\HECKE_\Gv$-module $$\QCshear(T^*_\Psi\Xv/\Gv)
 \simeq \cO_1\module(QC^!(\Ll \Psi/\Gv\times \Ga))$$
as in Proposition~\ref{quantum twisted cotangents}.
By Lemma~\ref{affine observables} its factorization homology is described by modules for the factorization homology $\cO_{1,\Sigma}$ of $\cO_1\in \HECKE_\Ga$ in the factorization homology of $QC^!(\Ll \Psi/\Gv\times \Ga)$. 
Thus the $\Theta$-series functors above become
$$\xymatrix{\HECKE^\Mv\ar[rr]\ar[dr]&&  Hom(\cO_{1,\Sigma}\module(QC^!(Loc_\Ga)), QC^!(\Loc_\Gv))\\
&\bH_{\Mv}\ar[ur]&}$$
(where we have also used the Hecke-linear self-duality of $QC^!(\Loc)$ to turn the tensor product into a functor category).
Now observe that the spectral exponential sheaf on $\Loc_\Ga$ is naturally an $\cO_{1,\Sigma}$-module. Thus we may apply the functors produced by the spectral $\Theta$-series to the spectral exponential sheaf, producing the desired $\Theta$-series functors
$$\xymatrix{\HECKE^\Mv\ar[rr]\ar[dr]&&  QC^!(\Loc_\Gv)\\
&\bH_{\Mv}\ar[ur]&}$$

One can likewise carry out the quantum hamiltonian reduction by $\Ga$ to identify the $L$-observables and their modules with the natural {\em twisted} counterparts of relative differential operators and flat connections as in the previous section. We omit the details.

\subsection{Geometric Arthur Parameters}\label{WhitArth}
 \index{Arthur parameter}

In this section we discuss the geometric counterpart of the theory of Arthur parameters and its relation to the process of Arthur functoriality on the automorphic side, as discussed in \S \ref{Arthur functoriality} on the numerical level. 
While Arthur's conjectures propose the parametrization of nontempered automorphic forms, the geometric counterpart proposes the parametrization of automorphic sheaves with a {\em sheared} Hecke eigenproperty and nontrivial singular support. The shearing and singular support are both captured by the datum of an $SL_2$ homomorphism. Thus the  
basic situation for this section (as in \S \ref{Whittaker induction} and \S \ref{Whittaker functoriality}) is that we are given
a subgroup $\check{H} \subset \check{G}$ and a commuting $SL_2$:
$$\iota: \check{H} \times \SL_2 \rightarrow \check{G}.$$  
We restrict ourselves to even $SL_2$'s, i.e., we demand that the corresponding cocharacter $\varpi$ acts on $\fgv$ with only even weights (see Remark~\ref{odd SL2s} for a discussion of the odd case).
Recall from \S \ref{nontempered} that in the setting of Arthur's conjectures \cite{Arthur-unipotent-conjectures}
one considers Arthur parameters  $$ \phi_A:  \Gamma_F \times \SL_2 \rightarrow \check G(k),$$ where $\phi_A|\Gamma_F$ is a Langlands parameter into the centralizer $\Hv$ of the $\SL_2$ which is pure of weight zero, as defining Langlands parameters\begin{equation*}  \phi_L = 
 \phi_A \circ  \left( \mathrm{id} \times \left[\begin{array}{cc} \varpi^{1/2} & 0 \\ 0 & \varpi^{-1/2} \end{array} \right]\right):   \Gamma_F \rightarrow \check G(k),\end{equation*}
into $\Gv$ (where $\varpi$ is the cyclotomic character). 

The geometric counterpart of this is the construction of {\em sheared local systems} for $\Gv$ out of (usual) $\Hv$-local systems and the commuting $SL_2$;
informally, these will be $\Gv$-local systems whose associated vector bundles are given cohomological gradings (and weights) by the associated cocharacter $\iota|_\Gm$.
 
 In \S \ref{Whittaker functoriality} we constructed an Arthur (or spectral Whittaker) induction functor from sheaves on $\Loc_{\check{H}}$
to sheaves on $\Loc_{\check{G}}$, while in \S \ref{local Whittaker induction} we constructed a local counterpart from $\HECKE_\Hv$-modules to $\HECKE_\Gv$-modules. We now establish local-global compatibility and use it to check that Arthur induction satisfies a Hecke eigenproperty with eigenvalues given by sheared local systems. In other words, this process constructs objects with Hecke eigenvalues given by a derived local system, cohomologically regraded by the diagonal part $\iota|_\Gm$ of the Arthur $SL_2$, and with singular support given by the Arthur nilpotent. This is a global counterpart of the ``Arthur'' properties of the Hecke module $\QCshear(\Mv/\Gv)$ described in Section~\ref{spectral local Arthur}.

\begin{remark}\label{odd SL2s}[Odd $SL_2$'s]
The case of odd $SL_2$ triples is not covered by our current construction. The essential issue is that of defining the global spectral geometric quantization -- i.e., $L$-sheaf -- of a non-polarized symplectic representation, in this case the action of the $SL_2$ centralizer $\Hv$ on $W=\fu/\fu_+$. Again, this question reduces to the case of the spectral Weil representation \S \ref{specWeil-start},  which we intend to discuss in more detail elsewhere. 
\end{remark}

\begin{remark}[Functoriality and domain walls]\label{domain wall remark} \index{domain wall}\index{Nahm pole}
More broadly, Arthur functoriality is expected to have the structure of a domain wall or interface (\S \ref{defects section}) 
between the arithmetic TQFTs 
$$\mathsf{AI}: \cB_\Hv\longrightarrow \cB_\Gv$$ describing the theory of Langlands parameters for the groups $\Hv$ and $\Gv$. Such a domain wall encodes in particular maps $\cB_\Hv(\Xi)\to \cB_\Gv(\Xi)$ for arbitrary input geometries $\Xi$, compatible with all the relations and structures of the field theory. (Indeed such a domain wall structure is to be expected for the spectral quantization of any hyperspherical $\Gv\times \Hv$-variety; see \S \ref{groupcase} for a related discussion of periods vs. functoriality.)

This domain wall, in the physics context of supersymmetric gauge theory, was constructed in the work of Gaiotto and Witten~\cite{GaiottoWittenboundary, GaiottoWittenSduality, GaiottoSduality} (see also Remark~\ref{GaiottoWitten intro}). Namely to an $SL_2$ homomorphism into a group $\Gv$ they attach a maximally supersymmetric (1/2 BPS) boundary condition for ${\mathcal{N}=4}$ super-Yang-Mills theory for (the compact form of) $\Gv$, the {\em Nahm pole boundary condition} associated to the $SL_2$ data $\iota$. The Nahm pole boundary condition has flavor symmetry given by the (compact form of) the centralizer $\Hv$ of the $SL_2$, which enables it to be coupled to the $\Hv$ super-Yang-Mills theory and thereby upgraded to a domain wall between the theories, which may further be topologically twisted to define a domain wall between the TQFTs $\cB_\Hv$ and $\cB_\Gv$. 
\end{remark}

\subsubsection{Sheared local systems} \label{local shear 2}
\label{sheared Loc}

In  Example \ref{repGinnershear} we discussed
the equivalence $\Rep(\Gv )^{\varpi \shear} \simeq \Rep(\Gv)$, which does not respect the standard fiber functor.
As a result, the tensor category $\Rep(\Gv)$ has many nonisomorphic fiber functors
$$(-)^{\varpi \shear}:\Rep(\Gv)\to \Vect, \hskip.3in V\mapsto \underline{V}^{\varpi \shear}$$ taking a representation to the shear of the underlying vector space by the given cocharacter
$\varpi:\Gm\to \Gv$. Geometrically, a cocharacter defines a map $B\varpi: B\Gm\to B\Gv$ and such a map defines a sheared fiber functor 
$$\xymatrix{\Rep(\Gv)\ar[r] & \Rep(\Gm) \ar[r]^-{(-)^\shear} &\Rep(\Gm)\ar[r]^-{\textrm{forget}} & \Vect.}$$
We will be interested in these functors primarily as defining $\Rep(\Gv)$-module category structures on $\Vect$ -- i.e., potential ``eigenvalues'' for $\Rep(\Gv)$-actions. We denote $\Vect$ with this $\Rep(\Gv)$-module structure as  $\Vect^{\varpi\shear}$ (or $\Vect_\rho^{\varpi\shear}$ if we keep track of twisting by a $\Gv$-torsor).

This is a remarkable feature of Tannakian formalism in the derived setting that diverges from classical experience:
 seemingly Tannakian categories like $\Rep(\check{G})$ have infinitely many non-isomorphic fiber functors due to the phenomenon of shearing. Indeed Tannakian reconstruction in the derived setting~\cite{SAG,BhattHL,GermanTannaka} requires that one impose a connectivity (or t-exactness) hypothesis to circumvent such phenomena. In our current setting this bug is a crucial feature:
  the sheared actions of $\Rep(\Gv)$ on $\Vect$ arise naturally as the Hecke eigenvalues for Arthur parameters (see e.g.~\cite{FrenkelNgo,LafforgueLysenko}). 
Recently S. Slaoui and G. Stefanich have shown that for a large class of geometric stacks $X$ shearing is the {\em only} obstruction to Tannakian reconstruction, in that isomorphism classes of tensor functors $QC(X)\to \Vect$ are parametrized not just by points of $X$ but by points and inertial cocharacters $\underline{x}:\pt/\Gm\to X$ via the shearing construction $\cF\mapsto (\underline{x}^*\cF)^\shear$.

Recall that $\Gv$-local systems $\rho$ on a curve $\Sigma$ are identified with connective tensor functors $\underline{\rho}:\Rep(\Gv)\to \Loc(\Sigma)$ to the tensor category of local systems on $\Sigma$, through the assignment of the associated local system to a representation $V\in \Rep(\Gv)$, $\underline{\rho}(V)=\rho\times^\Gv V$.
   If $\rho$ is endowed with a $\Gm$ of automorphisms $\varpi:\Gm \rightarrow \Aut(\rho)$, 
   then so is the symmetric monoidal functor $\underline{\rho}$ (i.e., each of the associated local systems $\underline{\rho}(V)$ is functorially and multiplicatively assigned a $\Gm$ symmetry). Hence 
   $\underline{\rho}$ can be sheared to
give a new functor $\underline{\rho}^{\varpi \shear}:\Rep(\Gv)\to \Loc(\Sigma)$, a {\em sheared $\Gv$-local system}, whose associated 
local systems $\underline{\rho}^{\varpi \shear}(V)$ carry nontrivial cohomological gradings. 

From another point of view, we can fix a cocharacter $\varpi:\Gm\to \Gv$ and consider $\Gv$-local systems $\rho$ endowed with a reduction to the centralizer $\Gv^{\varpi}$ of $\varpi$. We may then shear $\rho$ by the induced $\Gm$-symmetry. The associated local systems $\underline{\rho}^{\varpi \shear}(V)$ are given by giving the local systems $\underline{\rho}(V)$ (for $\Gv$-representations $V$) a cohomological grading determined by the $\varpi$-grading of $V$.

 These sheared local systems are not  
  points of $\Loc_\Gv$ in the usual sense, but rather ``Tannakian points'', in that they give non-connective tensor functors from $QC(\Loc_\Gv)$
  to $\Vect$. Namely, the $\Gm$-symmetry of the local system $\rho$ defines a morphism $i_\rho:\pt/\Gm\to \Loc_\Gv$, and hence a tensor functor $$\xymatrix{QC(\Loc_\Gv)\ar[r]^-{i_\rho^*} & \Rep(\Gm) \ar[r]^-{(-)^\shear} &\Rep(\Gm)\ar[r]^-{\textrm{forget}} & \Vect.}$$

Sheared local systems are the geometric avatars of the Langlands parameters associated to Arthur parameters.
Given an $\Hv$-local system $\rho_\Hv$ and a homomorphism $\iota:\Hv\times SL_2\to \Gv$, the induced local system $\rho= \rho_\Hv\times^{\Hv} \Gv$ has in particular a reduction to the centralizer of the induced cocharacter $\varpi=\iota|_\Gm$ of $\Gv$
and we can form the sheared local system $\rho^{\varpi \shear}$, 
which will play an important role in what follows.

\subsubsection{Arthur induction}
Recall that in \S \ref{Whittaker functoriality} we used the $L$-sheaf for the Whittaker space $\Mv=T^*_\Psi \Gv/U$ as the kernel for a spectral Whittaker induction functor
 $$\mathsf{AI}: QC^!(\Loc_\Hv) \longrightarrow QC^!(\Loc_\Gv)$$ from sheaves on $\Loc_{\check{H}}$
to sheaves on $\Loc_{\check{G}}$. On the other hand in \S \ref{local Whittaker induction} the local spectral quantization of $\Mv$ provided the kernel $\HECKE_{\Gv\gets \Hv\times SL_2}$ for a local counterpart, taking $\HECKE_\Hv$-modules to $\HECKE_\Gv$-modules. We now establish local-global compatibility by applying the Whittaker version of spectral $\Theta$-series from \S \ref{Whittaker L-observables}:

\begin{proposition}\label{Whittaker local-global} Local and global forms of Arthur functoriality are compatible: 
\begin{itemize}
\item  (One-point) Fixing a point $x\in \Sigma$,
 $\mathsf{AI}$ lifts to a natural map $\mathsf{AI}_{\HECKE}$ of $\HECKE_\Gv$-modules
$$ \xymatrix{ QC^!(\Loc_\Hv)\ar[rr]^-{\mathsf{AI}}\ar[dr]_-{1\ot id} && QC^!(\Loc_\Gv)\\
&\HECKE_{\Gv\gets \Hv\times SL_2}\otimes_{\HECKE_\Hv} QC^!(\Loc_\Hv) \ar[ur]_-{\mathsf{AI}_\HECKE}& } $$
\item  (Factorizable) More generally, $\mathsf{AI}$ lifts to a map $\mathsf{AI}_{\bH}$ of $\bH_{\Gv}$-modules  
$$\bH_{\Gv\gets \Hv\times SL_2}\otimes_{\bH_\Hv} QC^!(\Loc_\Hv)\longrightarrow QC^!(\Loc_\Gv)$$
where $\bH_{\Gv\gets \Hv\times SL_2}$ is the factorization homology
of  $\HECKE_{\Gv\gets \Hv\times SL_2}$.
\end{itemize}
\end{proposition}
\proof
Using the self-duality of $\QC^!(\Loc_\Hv)$ and the tensor product identification $$\QC^!(\Loc_\Hv)\otimes \QC^!(\Loc_\Gv)\simeq \QC^!(\Loc_{\Hv\times \Gv})$$ the desired functors are represented by integral transform constructions
$$\HECKE_{\Gv\gets \Hv\times SL_2} \longrightarrow QC^!(\Loc_{\Hv\times\Gv})$$ and 
$$\bH_{\Gv\gets \Hv\times SL_2} \longrightarrow QC^!(\Loc_{\Hv\times\Gv})$$
linear for the one-point and global Hecke categories for $\Hv\times \Gv$, respectively.
Now observe that the source of these functors as the (one-point and factorizable) spectral quantizations of the hyperspherical $\Gv\times \Hv$-space
$\Mv=T^*_\Psi\Gv/U$. Thus the desired functors are provided by the (one-point and factorizable) spectral $\Theta$-series construction of \S \ref{Whittaker L-observables}. 
  \qed
 
 We will now apply Proposition~\ref{Whittaker local-global} to establish properties of the Arthur induction functor which are global analogues of the local Arthur properties from \S \ref{spectral local Arthur}. To do so we analyze $\HECKE_{\Gv\gets \Hv\times SL_2}$, a factorization algebra in $\HECKE_\Gv\otimes \HECKE_\Hv$-modules. Recall from Proposition~\ref{quantum twisted cotangents} that ignoring its factorization structure (fixing a point $x\in \Sigma$), i.e., as a plain module category for the monoidal category $\HECKE_\Gv\otimes\HECKE_\Hv\simeq \QCshear(\fgxv/\Gv)\otimes \QCshear(\fhxv/\Hv)$, we have an equivalence $$\HECKE_{\Gv\gets \Hv\times SL_2}\simeq \QCshear(T^*_\Psi (\Gv/U)/\Gv\times\Hv)$$ where the module structure comes from applying $\QCshear$ to the $\GGm$-equivaraint diagram
 $$\xymatrix{&T^*_\Psi( \Gv/U) /\Gv\times\Hv\ar[dl]\ar[dr]& \\
 \fgxv/\Gv & & \fhxv/\Hv.}$$
 This identification does not respect factorization structures but it does respect the $SO(2)$-action of changes of coordinates, so gives an equivalence that is locally constant on $\Sigma$.

The Slodowy slice description of  Example~\ref{Slodowy slice example}, 
provides an identification of $\Hv\times \Gv\times \GGm$-spaces 
$T^*_\Psi \Gv/U\simeq \fgv_e\times \Gv,$ where $\Hv\actson\fgv_e$ via the coadjoint action, compatibly with the moment map to $\fhxv$, though not with the moment map to $\fgxv$.
Thus we have an isomorphism
\begin{equation}\label{HeckeID}
\HECKE_{\Gv\gets \Hv\times SL_2}\simeq \QCshear(\fgv_e/\Hv)
\end{equation} 
as plain $\Rep(\Gv)\otimes \HECKE_\Hv$-modules. Explicitly the $\Rep(\Gv)$ action is given by
\begin{equation}\label{sheared fiber action} \Rep(\Gv) \simeq QC(\mathrm{pt}/\Gv)^{\shear} \rightarrow QC(\pt/\Hv) 
\rightarrow 
QC(\fgv_e/H)^{\shear}\end{equation}
where the first isomorphism comes from the inner structure of the $\varpi$-action on $\Gv$ 
(see Example \ref{repGinnershear} and \S \ref{sheared Loc}). 

From this one can deduce the following

 \begin{corollary} \label{sheared eigenvalues}
Arthur induction takes Langlands parameters to Arthur parameters:
\begin{enumerate}

\item The functor $\mathsf{AI}$ is naturally $\varpi$-sheared: for any $x\in \Sigma$ it intertwines the $\Rep(\Hv)$ and $\Rep(\Gv)$-actions 
via the sheared forgetful functor $(-)^{\varpi\shear}:Rep(\Gv)\to Rep(\Hv)$. Moreover the identification is locally constant in $\Sigma$.

\item The functor $\mathsf{AI}$ produces $f$-antitempered sheaves
$$\mathsf{AI}:QC^!(\Loc_\Hv)\longrightarrow QC^!_{f-anti}(\Loc_\Gv)\subset QC^!(\Loc_\Gv),$$
i.e., sheaves annihilated by the local Hecke action at any $x\in \Sigma$ of the sheared structure sheaf $\cO^\shear_{\cNN_{<f}}\in \QCshear(\fgxv/\Gv)$.

\item Furthermore, when restricted to $QC(\Loc_\Hv)$ the functor $\mathsf{AI}$ produces sheaves which are also $f$-tempered,
i.e., with singular support contained in $\cNN_{<f}$.

\end{enumerate}
\end{corollary}

 In particular for any smooth point $\{\rho_\Hv\}\in \Loc_\Hv$, the Whittaker induction of the skyscraper $\mathsf{AI}(\kk_{\rho_\Hv})$ is an $f$-tempered and $f$-antitempered sheared Hecke eigensheaf with eigenvalue the shear $(\rho_\Gv)^{\varpi\shear}$ of the induced $\Gv$-local system (cf. Section~\ref{sheared Loc}).

 This construction of geometric Arthur parameters suggests an explicit description of the graded pieces of the ``Arthur filtration'' of the spectral category $\QC^!_{\cNN}(\Loc_\Gv)$ by $f$-tempered sheaves for nilpotent orbits $f$ (as discussed for example in~\cite{LysenkoFourier}), whose numerical counterpart is discussed in \S \ref{Arthur functoriality}. Namely we conjecture
 
 \begin{conjecture}
The functor $\mathsf{AI}$, restricted to quasicoherent sheaves on $\Loc_\Hv$, generates the $f$-tempered graded piece of the Arthur filtration.

\end{conjecture}

  If the Arthur restriction is identified with the adjoint of Arthur induction (see Remark \S \ref{ArthurResRemark}) this would then produce a monadic descripition of the pieces of the Arthur filtration in terms of the bi-Whittaker reduction $U {}_\psi \lGIT T^*\Gv\GIT_\psi U$.

%% file: shearing.tex
  \section{Koszul Duality.}\label{shearingsec}
 
  We discuss some examples of shearing in  relation to Koszul duality and the construction of a spectral analog of the exponential sheaf.

\subsection{Koszul duality and sheaves on lines.}\label{shearing affine line}
We now examine in some detail various categories of modules for symmetric and exterior algebras, which arise throughout this work in many guises. As elsewhere, we work with coefficients of characteristic zero, and all constructions are derived, i.e., take place in the relevant dg or $\infty$-categories.  (We recommend~\cite{GKM} for a thorough discussion of the basic issues of Koszul duality in the language of triangulated categories, and~\cite{DrinfeldGaitsgory} for the $\infty$-categorical setting.)

The essential features are all visible in the case of the symmetric algebra $S= k[x_0]$ in a single variable of cohomological degree $0$, which we give $\GGm$-weight $-2$, and the Koszul dual exterior algebra $\Lambda=k[y_{-1}]$ on a generator in cohomological degree $-1$, which we give $\GGm$-weight $2$. (We keep track of cohomological degrees with subscripts and have $x,y$ denote dual coordinates, indicating $\Gm$-weights. All algebras appearing in this section will be formal, so we will
often identify a dg algebra with its cohomology ring.)
 
\begin{remark}[Even shearings]\label{even shearings} 
Note that the $\Gm$-weights appearing in this section -- and hence the associated shearings -- are all even, so that we do not encounter changes of parity and all vector spaces are considered even.
\end{remark}

\begin{table}
\centering
\caption{Modules under $S, S^{\shear}, \Lambda$}
  \begin{tabular}{|c|c|c|c|c|}
  \hline
 &  $S$-mod & $S^{\shear}$-mod & $\Lambda$-mod &   Ind-coherent$(\Lambda)$ \\
 \hline
  (a) & $S$ & $S^{\shear}$ & $k$  & $k$, see \eqref{kplus} \\
  (b) & $k_0$ & $k^{\shear}$ & $\Lambda$ & $\Lambda$ \\
  (c) & $k_1$ &   &   & See \S \ref{spectral exponential} \\
  (d) &  $S[x^{-1}]$ & $(S[x^{-1}])^{\shear}$ &   &\eqref{Tate module} \\ 
\hline
 \end{tabular}
 \end{table}
 \subsubsection{$S$-modules}
The basic category we start from is the usual module category for $S$, i.e., quasicoherent sheaves on the affine line $\AA^1=\Spec(S)$,
$$S\module\simeq QC(\AA^1).$$
We will keep track of four basic $S$-modules: 
\begin{itemize}
\item[(a)] the structure sheaf $k[x_0]$;
\item[(b)] the skyscraper $k_0$ at $0$, arising from $S \rightarrow k$ sending $x_0$ to $0$; 
  \item[(c)] the skyscraper  $k_1$ at $1$,arising from $S \rightarrow k$ sending $x_0$ to $1$.
\item[(d)] the structure sheaf of the punctured line $k[x_0, x_0\inv]$. 
\end{itemize}
The first three objects are represented by perfect complexes of $S$-modules, which categorically speaking form the small category of compact objects
 in the big category of arbitrary unbounded complexes of $S$-modules 
$$\mathrm{Perf}(S)=(S\module)^c \hookrightarrow S\module\simeq  \mathrm{Ind}(\mathrm{Perf}(S)).$$

 \subsubsection{$S^{\shear}$-modules} \label{S shear modules}
Now let us shear, according to the $\GGm$-action: we have $S^\shear\simeq k[x_2]$ where $u=x_2$ has cohomological degree $2$, which we identify with the $\Gm$-equivariant cohomology ring $H^\ast(B\Gm)$. 
Passing to categories of modules, we have $$S^\shear\module \simeq \QCshear(\AA^1).$$

As discussed in \S~\ref{shearcategory}, $\QCshear(\AA^1)$ is not equivalent to $QC(\AA^1)$ itself, though the corresponding categories of {\em graded} modules  are equivalent,
i.e., quasi-coherent sheaves on $\AA^1/\Gm$ is equivalent to its shear. 
  The category $S^\shear\module$ has a topological interpretation as a form of the $\Gm$-equivariant category of a point. More precisely it is the ind-finite' or ``renormalized'' form of $\Shv_\Gm(\pt)$, defined as the ind-category of $\Gm$-equivariant constructible sheaves on a point.  
 
 We note three basic objects in $\QCshear(\AA^1)$:
 \begin{itemize}
 \item[(a)]  the regular module $S^\shear=H^\ast_\Gm(\pt)$ itself;
 \item[(b)]  the augmentation $k_0^{\shear} =k=H^\ast_\Gm(\Gm)$;
 \item[(d)] the periodic module $k[x_2,x_2\inv]$ (the structure sheaf of the sheared punctured line, or $\Gm$-Tate cohomology of a point). 
 \end{itemize}
 The first two are compact objects, i.e., objects of the small category $\Perf(S^\shear)=\QCshear(\AA^1)^c$. These objects
 (endowed with evident $\Gm$-equivariant structures) are
 the shears of the correspondingly labelled objects (a), (b), (d) 
 inside $S$-modules. However, 
the ungraded $S$-module $k_1$ does not have an analogue in the sheared category.   
 
 \subsubsection{$\Lambda$-modules}
Let us consider the Koszul dual setting. Let $$\Lambda=k[y_{-1}]\simeq \Hom_{S^{\shear}}(k,k) $$ be the Koszul dual (homological) exterior algebra, 
 which is naturally identified with the homology $H_\ast(\Gm)$ with its Pontrjagin product (i.e., the ``topological group algebra'' of $\Gm$). 
 We can also consider $\Lambda$ as the ring of functions on the affine derived scheme $\Spec(\Lambda)\simeq \AA^1[-1]$, the shifted version of the {\em dual} affine line, so that $\Lambda\module\simeq QC(\AA^1[-1])$. Topologically, $\Lambda$-modules realize locally constant actions of $\Gm$, i.e., $\Gm$-equivariant sheaves on a point;
 indeed $\Lambda$-modules recover the standard notion of equivariant $D$-modules on a point,
 $$\Lambda\module\simeq \D(B\Gm).$$

  We note two basic bounded coherent (finite dimensional) objects in $\Lambda\module$: 
  \begin{itemize}
  \item[(a)] the augmentation $k$, and
   \item[(b)]
the regular module $\Lambda$.
\end{itemize}
 Again the small category of compact objects here is given by perfect complexes of $\Lambda$-modules. 
The first object is not perfect (=compact), while the latter is.
 In the realization $\Lambda\module\simeq QC(\AA^1[-1])$, $k$ corresponds to the skyscraper sheaf at the origin, a singular point.

 We have functors switching the augmentation and regular modules for $S^\shear$ and $\Lambda$ 
\begin{equation} \label{Koszul-naive}- \ot_{S^{\shear}}k: S^\shear\module \longleftrightarrow \Lambda\module: Hom_{\Lambda}(k,-)\end{equation}
$$ k \longleftrightarrow \Lambda[-1]$$
$$ S^{\shear} \longleftrightarrow k.$$
 i.e., they exchange equivariant and ordinary cohomology of $\Gm$-spaces~\cite{GKM}; 
  $\Lambda[-1]$ can be thought of as the dualizing sheaf of $\Lambda$, or as the cohomology of $\Gm$. 
 
 In particular the functors can not be inverse equivalences on the full (unbounded) categories of modules, since the noncompact augmentation of $\Lambda$ is taken to the regular module for $S^\shear$,
 and the latter is compact inside the category of $S^{\shear}$-modules.  In other words, they don't restrict to equivalences of the categories of perfect modules for $S^\shear$ and $\Lambda$, though they do identify perfect $S^\shear$-modules with bounded coherent $\Lambda$-modules.  

Another way to see the failure of \eqref{Koszul-naive} to be an equivalence is that the periodic module $k[x_2,x_2^{-1}]\in S^{\shear}\module$ vanishes under \eqref{Koszul-naive}, for it is sent to the acyclic complex
\begin{equation}\label{Tate module} \underline{\mathsf{Per}}: \cdots \rightarrow  \kk[y_{-1}] \stackrel{y_{-1}}{\rightarrow} \kk[y_{-1}]  \stackrel{y_{-1}}{\rightarrow} \kk[y_{-1}] \rightarrow \cdots 
\end{equation}
In this equation, the various copies of $\kk[y_{-1}]$ are placed in degrees that differ by $2$ from each other. 

\subsubsection{Ind-coherent modules for $\Lambda$}\label{indcoh for exterior algebra}

We can make Koszul duality \eqref{Koszul-naive} an equivalence by ``correcting'' either side. 
 
The solution we will adopt is  to correct the side of $\Lambda$-modules:  we enlarge $\Lambda\module=QC(\AA^1[-1])$ to the {\em ind-coherent category} of $\Lambda$,
i.e., the ind-category of the category of $\Lambda$-modules
with finite-dimensional cohomology. We will denote this category by $QC^!(\AA^1[-1])$.
This enlarged category of $\Lambda$-modules consists of formal colimits of finite dimensional (coherent) $\Lambda$-modules.
If we replace $\Lambda$ by a usual ring (rather than a dg-ring) the
ind-coherent  category can be equivalently
described as the (dg-enhanced) homotopy category of injective complexes of $\Lambda$-modules~\cite{Krause};
in such a category acyclic complexes of the general form \eqref{Tate module} need not be trivial. 

Since the category of $\Lambda$ modules is the ind-completion
of the category of perfect complexes, extending the inclusion of perfect
into coherent complexes gives 
\begin{equation} \label{Upsilondef} \Xi : QC \rightarrow QC^{!}.\end{equation} 
 We warn that $\Xi$ does not agree, when restricted to coherent sheaves,
to the tautological inclusion of coherent sheaves into $QC^!$ --
see examples (b), (b)' below. 

\begin{remark} \label{remark-colocalization}
We note that there is also a functor $$\Psi: QC^{!} \rightarrow QC$$ in the opposite direction (we follow the lettering of \cite{indcoh}),
which is right adjoint to $\Xi$, and 
  arises from the ind-completion of the inclusion of coherent sheaves into $QC$.
This functor is not fully faithful.  (In the case when $\Lambda$ is a usual Noetherian ring,
  this $\Psi$ is just the tautological functor from injective complexes to the derived category and is denoted $Q$ by Krause.)
\end{remark}

The passage from usual modules to the ind-coherent category has
the effect of (in fact is designed for) forcing the augmentation of $\Lambda$ to be a compact object (in fact a compact generator) like $S^\shear\in (S^{\shear})\module$, and indeed Koszul duality now extends to an equivalence
 $$S^\shear\module=\QCshear(\AA^1)  \longleftrightarrow QC^!(\AA^1[-1])=\Ind(\Lambda^\cdot\module_{f.d.}).$$

Some examples:

\begin{itemize}

\item[(b)]
$\kk$ as a $\Lambda$-module is coherent and thus 
gives an object of the ind-coherent category (via the inclusion 
of the coherent category to the ind-coherent). 
It will be convenient, for comparison with what follows
to represent $\kk$   by the coherent complex of free modules
\begin{equation}\label{kplus}  \kk: \qquad \cdots \rightarrow 0  \rightarrow  0 \rightarrow \kk[y_{-1}]  \stackrel{y_{-1}}{\rightarrow} \kk[y_{-1}] \stackrel{y_{-1}}{\rightarrow} \cdots 
\end{equation} 
(where the various copies of $\kk[y_{-1}]$ are generated in degrees $1, 3, 5, \dots$) and  under Koszul duality corresponds to   $S^{\shear}$.

\item[(b)']
The image of $k$ under the functor $\Xi$ of \eqref{Upsilondef}
can be thought of as   the {\em bounded above} complex of free $\Lambda$-modules:  
   \begin{equation}\label{kminus} \Xi(\kk) :  \qquad  \cdots \rightarrow   \kk[y_{-1}]   \stackrel{y_{-1}}{\rightarrow} \kk[y_{-1}]  \rightarrow 0  \rightarrow  0 \rightarrow \cdots
  \end{equation}
where the various copies of $\kk[y_{-1}]$ are generated in degrees $0,-2,-4,\dots$, and 
which we regard  as being in the ind-coherent category
by taking only finitely many terms of the above, and then taking a direct limit.  
   $\Xi(k)[1]$ corresponds under Koszul duality to $S^{\shear}[x_2^{-1}]/S^{\shear}$,
   which is a direct limit of torsion modules.

\item[(c)]  The infinite acyclic complex   encountered in~\ref{Tate module},
understood as the colimit of bounded-below truncations as above, 
now yields a nonzero object $\underline{\mathsf{Per}}$ of the ind-coherent category,
 which fits now into a sequence
 $k \rightarrow  \underline{\mathsf{Per}} \rightarrow \Xi(k)[1]$.
 This  $\underline{\mathsf{Per}}$ is the image of $S^{\shear}[x_2^{-1}]$ under Koszul duality.  Under
 the functor  $\Psi$ of Remark \ref{remark-colocalization} it is carried to the trivial object.

\end{itemize}

\subsubsection{Singular support}
The distinction between the categories  $$\Lambda\module\simeq QC(\AA^1[-1]), S^\shear\module\simeq QC^!(\AA^1[-1])$$ can also be captured by the theory of singular support of (ind-)coherent sheaves~\cite{BIK,ArinkinGaitsgory}.

Namely, the singular support of an object of $QC^{!}(\AA^1[-1])$ 
is a closed conical subset of $\AA^1$ -- i.e., either $\{0\}$ or all of $\AA^1$!
Explicitly
the singular support is zero if and only if
the Koszul dual module is torsion.  ``Quasicoherent sheaves are  ind-coherent sheaves with trivial singular support,'' in the sense that
$$ \Xi: QC(\AA^1[-1]) \rightarrow QC^{!}(\AA^1[-1])$$
 is fully faithful with essential image those objects with zero singular support.

\begin{remark}
Care is needed with these notions! As we saw above, the augmentation object $k$ in $\Lambda$-modules
can be considered in $QC^{!}$ in two different ways, giving
the objects called $k$ and $\Xi(k)$ above.

\begin{itemize}
\item[-] 
The singular support of $k$  is $\mathbb{A}^1$, but 
\item[-] The singular support of $\Xi(k)$ is $\{0\}$. 
\end{itemize}
 
 \end{remark}

\subsection{The spectral exponential sheaf}\label{spectral exponential}
We now define an exotic object on the affine line which plays the role of the exponential $D$-module or Artin-Schreier sheaf in the coherent setting. It can be considered a spectrally quantized form of the Hamiltonian $\Ga$-space
$$\pt_1 := \mbox{ a point with moment map value $1\in \gax$}.$$ In other words, we are seeking an algebraic
avatar of the exponential function, which defines a character of the Lie algebra $$\exp:\ga:=\mathrm{Lie}(\Ga)\to k$$ of the additive group -- but not of the additive group itself.
Correspondingly, we are going to construct this spectral exponential not 
 as a quasicoherent sheaf on $B\Ga$, but rather inside   a sheared version of that category.
Indeed, the category $QC(B\Ga)=\Rep(\Ga)$ of representations of the additive group is identified by Cartier duality with sheaves on the dual additive formal group, i.e., the formal completion of the dual Lie algebra $\gax$ at $0$, and thus doesn't have an object corresponding to the skyscraper $k_1$. However, we can perform the decompletion formally using Koszul duality and shearing.

To formalize this, we use the same setting as the previous section, with Koszul dual symmetric and exterior algebras $S=k[x_0]$ (where $x_0$ has $\Gm$ weight $-2$) and $\Lambda=k[y_{-1}]$ (where $y_1$ has $\Gm$ weight $2$). To make the link
with what we just said, we adopt the following point of view: 
\begin{itemize}
\item[-]  We view the $y$-line $\AA^1[-1]\simeq \Omega \Ga$ as the based loops in the additive group.\footnote{Equivalently, the shift by 2 of the classifying stack $B\Ga\simeq \AA^1[1]$ (a coaffine stack, cf.Remark~\ref{coaffine}).}
\item[-]  We view the $x$-line $\Spec(S)\simeq \gax$ as the dual Lie algebra of the additive group. 
\end{itemize}

As discussed above, 
Koszul duality identifies ind-coherent sheaves on the $y$-line
 $$QC^!(\AA^1[-1]) \simeq QC^{\shear}(\gax)$$ with the shear of ordinary quasicoherent sheaves on $x$-line. 
This is an equivalence of $\Gm$-categories, where, writing out the actions:
\begin{itemize}
\item  the shearing $\Gm$ is acting on $\AA^1$ by  the inverse square character,  and
 on the coordinate $y$ on $\AA^1$ by the square character;
 \item dually, the shearing $\Gm$ acts on $\gax$ by the square character, and on the coordinate $x$ on $\gax$ by the inverse square character. 
\end{itemize}

 Now observe that since the $\Gm$-weights of $x$ and $y$ are opposite, we find a sheared Koszul duality equivalenc 
 $$QC^!(\AA^1[-1])^{\shear, \textrm{negated action}} \simeq QC(\gax).$$
On the left,  
we are  now shearing by the   the {\em inverse} of the action just described, i.e.
the $\Gm$-action on $\AA^1[-1]$  is through squaring. We are led to the following definition:

\begin{definition} \label{sesdefn}  Let $\Gm$ act on $\AA^1$ through squaring. 
The spectral exponential sheaf  $$\exp\in QC^!(\AA^1[-1])^\shear$$ is the image of the skyscraper $k_1$ at 1 (evaluation module at $1\in \gax$) by the sheared Koszul duality equivalence  $QC(\gax)\simeq QC^!(\AA^1[-1])^\shear$. 
\end{definition}

Here are a couple of different ways to think about this construction:
\begin{itemize}
\item  Recall that by definition, the shear of a category has the same $\Gm$-equivariant objects as the original category, but with the enriched graded Hom-spaces modified by a shear, and from this data one formally reconstructs the whole category. 
Thus we start with the $\Gm$-equivariant skyscraper, the augmentation object $\kk \in QC^{!}(\AA[-1])$, which corresponds to $S^\shear\in S^\shear\module$.
Thus its endomorphism algebra is the symmetric algebra $S^\shear$. Now, the sheared category $\QC^{!}(\AA[-1])^{\shear}$ is built so that (i) it has an object by name $\kk^{\shear}$ but (ii) this object has endomorphisms the naive symmetric algebra $S=k[x_0]$ itself (now entirely in degree zero), so that (iii) we can form the triangle $\kk^{\shear} \stackrel{x_0-1}{\rightarrow} \kk^{\shear}$. This formal construction is exp.  

\item 
We can see $\exp$ as a deformation of $\mathcal{O}^{\shear}$, the shear of the structure sheaf $\mathcal{O} \in QC^{!}(\AA^1[-1])$, 
in the sense that there is a functor $QC(\AA^1) \rightarrow \QC^{!}(\AA[-1])^{\shear}$
carrying the skyscraper $k_0$ at $0$ to $\mathcal{O}$ and the skyscraper at $1$ to  $\exp$ (this is just a rephrasing of Definition~\ref{sesdefn}).
With reference to this functor, the ``deformation class''
in $\mathrm{Ext}^1(k_0, k_0)$ is carried to the degree $+1$ endomorphism of $\mathcal{O}^{\shear}$,
which comes from taking $y_{-1} : \mathcal{O} \rightarrow \mathcal{O}$, and shearing: 
  $y_1 = y_{-1}^{\shear}: \mathcal{O}^{\shear} \rightarrow \mathcal{O}^{\shear}$  has degree $1$.   We can think of 
  the pair of $\mathcal{O}$ and   the self-map
  $y_{-1}$ as a  kind of semiclassical limit of $\exp$. 
\end{itemize}

\begin{remark}
Note that the exponential sheaf $exp$ is naturally a character sheaf on the group-stack $\AA^1[-1]=\Omega\AA^1$, i.e., a (commutative) algebra object with respect to the convolution symmetric monoidal structure. This follows from the commutative algebra structure with respect to tensor product on the skyscraper sheaf under the symmetric monoidal functors of shearing and Koszul duality (which identifies the convolution symmetric monoidal structure with tensor product of modules for the symmetric algebra). A similar observation was made by Hilburn and Yoo~\cite{HilburnYoo}.
\end{remark}

%% file: sheaf-theory.tex
\newcommand{\Tr}{\mathrm{Tr}}
\newcommand{\cN}{\mathcal{N}}
\newcommand{\conjg}{\mathfrak{c}}
\newcommand{\Sing}{\mathrm{Sing}}

\section{Sheaf theory}
\label{sheaf theory}

In this appendix we survey the somewhat bewildering array of different categories of sheaves that we encounter in the paper. Definitive reference for many of the features we recall include~\cite{DrinfeldGaitsgory,GR,AGKRRV1}.   The contents of the section are as follows:

\begin{itemize}
\item \S \ref{HigherCatAppendix} recalls general features of the category theory we will use. 
\item \S \ref{coherent sheaf theories} discusses categories of algebraic sheaves (``coherent sheaf theory").
\item \S \ref{constructible sheaf theories} begins our discussion of categories of topological sheaves (``constructible sheaf theory").
\item \S \ref{sheaves on stacks} continues by discussing constructible sheaf theory on stacks.
\item \S \ref{renormalization section} discusses the ``finiteness versus safety'' distinction for constructible sheaf theory on stacks (also known as renormalization). 
\item \S \ref{infinite type} discusses constructible sheaves on infinite type objects. 
\item \S \ref{tensor products of sheaves} discusses duality structures on categories of sheaves.
\item \S \ref{shvcats} colllects some notions we will make use of concerning sheaves of categories, the ULA condition and rigid tensor categories.
\end{itemize}

\subsection{The format of sheaf theories: synopsis}\label{synopsis}

Let us first start with a synopsis of the overall format of sheaf theories.
 The (dg) categories we encounter come in two general flavors, ``small'' (consisting of sheaves with finiteness conditions, such as coherence or constructibility) and ``large'' (consisting of unbounded complexes and closed under arbitrary direct sums). The sheaf theories are divided into two archetypes: topological (A-side) (starting in \S \ref{constructible sheaf theories}), appearing on the automorphic side of the Langlands correspondence and algebraic (B-side), Section~\ref{coherent sheaf theories}, appearing on the spectral side. 
The topological theories are further divided in three types: de Rham, Betti and (intermediate between the two) \'etale. 

The general format of the construction of sheaf categories is as follows. One first defines sheaves on the basic building blocks -- 
for example, finite type schemes (see below for more precise discussion). This assembles into a contravariant functor from such schemes to categories under (suitably chosen) pullback. 
To incorporate pushforwards, base change and adjunctions, it is extremely convenient to use the formalism of correspondences as laid out in~\cite{GR} making sheaf categories a functor out of a category of correspondences of schemes.

Next we need to define sheaf theory on more general objects in algebraic geometry, such as stacks, infinite-type schemes, ind-schemes, and most generally prestacks (arbitrary ``functors of points'', i.e., functors from affines to simplicial sets). A unifying theme in sheaf theory is that we first define categories of sheaves on a class of objects which we take to be basic building blocks. We then define the category of sheaves in general as 
the limit of the categories of sheaves on building blocks under pullback -- i.e., as a {\em right Kan extension}.

So the question  is: what are  the basic building blocks? We will encounter three basic variants of this idea in the ``topological'' setting.

\begin{itemize}
\item[$\bullet$][Safe sheaves]
Here we take the basic building blocks to be finite type schemes. The resulting
categories built by right Kan extension will be called {\em safe sheaf categories}.
These are the sheaf categories most commonly used in the geometric Langlands program, eg~\cite{AGKRRV1}.
The corresponding compact objects in these categories on algebraic stacks are the safe sheaves of~\cite{DrinfeldGaitsgory} and hence we refer to objects of the large categories as ind-safe sheaves.

\item[$\bullet$][Finite sheaves and ``renormalization''].

There is an alternative which is often better adapted to {\em equivariant} settings:  take the basic building blocks to be finite type algebraic stacks 
and there we use the category of {\em ind-finite sheaves}, i.e., the ind-category of ``finite'' objects (constructible sheaves or coherent $D$-modules). We then use these sheaf categories as the basic building blocks to define the categories of ind-finite sheaves on arbitrary prestacks (still locally of finite type).  These are referred to as {\em renormalized} sheaf categories in the geometric Langlands literature.

The distinction between ind-finite sheaves and ind-safe sheaves on stacks 
arises from the fact that, for example, the constant sheaf on $BG$ 
is   certainly finite in a reasonable sense -- for example, working locally on covers --  but is 
not a compact object of the safe sheaf category.
This is an aspect of the fundamental theme of completion/decompletion in equivariant topology, and plays the role on the automorphic side that the distinction between perfect and coherent complexes (an aspect of singularity theory) plays on the spectral side.

\item[$\bullet$][Sheaves in infinite type] 
Finally, in local Langlands we encounter geometric objects of infinite type, such as the loop group $G_F$, the arc and loop schemes $X_O, X_F$ for $X$ an affine $G$-scheme and the quotient stack $X_F/G_O$. To define sheaf theory in this context requires an opposite procedure: (affine, or quasi-compact quasi-separated) schemes of infinite type are naturally constructed as {\em limits} of schemes of finite type. One defines sheaf categories on such objects as {\em colimits} of sheaf categories on finite type schemes under pullback, in other words as a {\em left} Kan extension. (Here again there's a bifurcation, whether we use $!$- or $*$-pullback, though the two theories give equivalent answers in the {\em placid} setting.) We again can choose how to treat equivariance for affine group schemes like $G_O$ (with safety or finiteness). Once (equivariant) sheaves have been defined on all schemes in this fashion we can use them as building blocks (under colimits) for general stacks, ind-schemes and prestacks and use the same (right Kan or limit) procedure as before to extend sheaf theory.
\end{itemize}

\subsection{Higher categories: small and large}  \label{HigherCatAppendix}
Let us begin by specifying the world in which we take categories of sheaves as living. 
Recall that we work over a fixed field $\kk$ of coefficients, of characteristic zero. We work primarily with $\kk$-linear differential graded (dg) categories, and often abuse notation to refer to such objects simply as {\em categories}. Moreover our dg categories 
will always be stable (or {\em pre-triangulated}, a property which implies that their homotopy category is triangulated). An equivalent notion is provided by the theory of stable $\kk$-linear $\infty$-categories, and we use the terms ``dg category'' and ``stable $\kk$-linear $\infty$-category'' interchangeably. 
Eventually, of course, it would be desirable to have a formulation of our conjectures in arbitrary characteristic, in which case the language of stable $\infty$-categories would be more suitable, but in our present setting we find it conceptually easier to speak of dg-categories.

 We will use the language of homotopical algebra developed by Lurie in~\cite{HTT,HA}, for which we refer to the exposition in~\cite[I.1.5-8]{GR}. 

We consider two main classes of dg-categories, informally referred to as ``small'' and ``large''. This refers to the size of the category in a set-theoretic sense: the categories that we call ``small''  are essentially small (i.e., their isomorphism classes are sets), while the large categories are not, and are typically obtained ind-completing small categories. The small categories we consider are closed under finite limits and colimits, while the large ones are closed under all (small) limits and colimits. We will distinguish them notationally using either lower case/upper case letters, or without/with an overline; \index{$\Shv$} \index{$\SHV$}
\index{$\overline{\mathcal{H}}$}
\begin{quote}
\begin{center}
example of notation: $\Shv$ or $\mathcal{S}$ (small) versus $\SHV$ or $\overline{\mathcal{S}}$ (big).
\end{center}
\end{quote}

\index{category over $\kk$}
\index{presentable category}
\index{$\infty$-category}
\index{monoidal $\infty$-category}

$\bullet$ The $\infty$-category $\mathsf{DGCat}_{\mathbf k}$ has as objects 
small, idempotent-complete dg-categories with morphisms given by exact functors. Such dg-categories are closed under finite (homotopy) limits and colimits. Examples include (dg-enhanced derived) categories of constructible sheaves, perfect complexes or bounded coherent complexes or finitely presented modules over a ring. 

$\bullet$ The $\infty$-category $\overline{\mathsf{DGCat}}_{\mathbf k}$ has as objects presentable dg-categories with morphisms given by colimit-preserving functors. Presentability encodes that these categories are cocomplete (closed under all small colimits - in particular infinite direct sums), but also accessible (are generated under suitably controlled colimits by a small set of objects -- a weaker notion than compact generation). 
They are automatically closed also under all small limits ~\cite[Corollary 5.5.2.4]{HTT}, and enjoy the $\infty$-categorical form of the adjoint functor theorem ~\cite[Corollary
5.5.2.9]{HTT}. Thus the large categories are much more suited for formal categorical operations. Examples include the (dg-enhanced) unbounded derived categories of quasicoherent sheaves or of all modules over a ring.

Both $\mathsf{DGCat}_{\mathbf k}$ and $\overline{\mathsf{DGCat}}_{\mathbf k}$ admit natural structures of {\em symmetric monoidal} $\infty$-categories under the Lurie tensor product. In particular we may perform higher algebra in these settings and speak of algebra objects (which are themselves monoidal dg categories and have their internal theory of algebras and modules), modules, tensor products and so on, with all notions taken in the homotopical ($\infty$-categorical) sense. 

One can formally add filtered colimits to a small category $\cC$ to obtain a (compactly-generated) presentable category $\Ind(\cC)$, and this defines a symmetric monoidal functor
$$Ind: \mathsf{DGCat}_{\mathbf k}\longrightarrow \overline{\mathsf{DGCat}}_{\mathbf k}.$$ Conversely we may pass from a presentable category $\D$ to its small category of compact objects $\D^c$. These operations define quasi-inverse equivalences between the category of small categories $\mathsf{DGCat}_{\mathbf k}$ and that of {\em compactly generated} presentable categories, with morphisms restricted to functors which preserve compact objects. This allows us to pass back and forth between the small and large settings when convenient, and indeed all large categories of interest to us will be compactly generated. 

\begin{remark}[Calculating colimits of large categories]\label{calculating colimits} 
An important technique for working with large categories is that a colimit in $\overline{\mathsf{DGCat}}_{\mathbf k}$ under functors which have colimit-preserving right adjoints (i.e., in the compactly generated setting, functors that preserve compact objects) can be identified with the limit in $\overline{\mathsf{DGCat}}_{\mathbf k}$ of the same categories under the right adjoints. This is ~\cite[Corollary 5.5.3.4]{HTT}, identifying the opposite of the $\infty$-category ${\mathcal Pr}^L$ of presentable categories under left adjoints (where $\overline{\mathsf{DGCat}}_{\mathbf k}$ is defined as $\Vect_{\mathbf k}$-modules) the $\infty$-category ${\mathcal Pr}^R$ of presentable categories under right adjoints.
\end{remark}
\index{$Ind$}

\subsubsection{Dual categories}\label{duality section}
Recall that a (presentable dg) category $\mathcal{C}$ is {\em dualizable} if there exists another presentable category $\cC^\vee$ and (colimit-preserving) functors 
$$u:\Vect\to \mathcal{C}\otimes \mathcal{C}^\vee$$ 
(the unit or coevaluation) and 
$$\epsilon:\mathcal{C}\otimes \mathcal{C}^\vee \to \Vect$$
(the counit or evaluation)
 satisfying a standard relation. (See~\cite[4.6]{HA} for the general $\infty$-categorical notion of duality.) Such duality data are uniquely determined in the higher categorical sense, i.e., up to contractible choices, if they exist. Compactly generated categories $\cC=\Ind(\cC^c)$ are automatically dualizable, and the dual category can be described explicitly as the ind-category of the opposite to the small category of compact objects, $$\cC^\vee\simeq \Ind((\cC^c)^{op}).$$ 

\index{dual category}

Self-duality for a category is an identification $\mathcal{C}^\vee\simeq \mathcal{C}$, which could come from a contravariant autoequivalence of the generating category of compact objects.
Self-duality is equivalent to giving unit and counit maps
$$U:\Vect\to \mathcal{C}\otimes \mathcal{C}, \hskip.3in E:\mathcal{C}\otimes \mathcal{C}\to \Vect$$ 
satisfying the same relation. Self-duality is additional data, but it suffices to specify the unit $U$ (which reproduces an isomorphism $\mathcal{C}^\vee\simeq \mathcal{C}$ by tensoring with $\cC^\vee$ and contracting by $\epsilon$). In the presence of self-duality we can convert bilinear forms to endomorphisms, i.e., we have an equivalence
$\cC\ot\cC\to End(\cC).$

\subsection{Coherent sheaf theories}\label{coherent sheaf theories}
On the spectral side of the Langlands correspondence we will make use of categories of coherent sheaves on derived schemes and stacks of finite type over $\kk$, a field of characteristic zero. These come in two main variants:  \index{$QC$} \index{$QC^{"!}$} \index{$\Perf$} \index{$\Coh$}
\begin{itemize}
\item[-] On the one hand, the big category of quasicoherent sheaves $QC(X)$ and its small version $\Perf(X)$ consisting of perfect complexes;
\item[-]   on the other hand the big category $QC^!(X)$ of ind-coherent sheaves and its small version $\Coh(X)$ consisting of coherent sheaves
 (bounded coherent complexes). 
 \end{itemize}
 Both have natural pullback and pushforward functors and symmetric monoidal structures. The former enjoys $f^*$ functoriality without any restrictions and plays the role of ``functions'', with tensor unit $\cO$, while the latter enjoys $f^!$ functoriality without any restrictions and plays the role of ``distributions'', with tensor unit the dualizing sheaf $\omega$. The book~\cite{GR} provides a comprehensive account of the general sheaf theory, while the articles~\cite{DrinfeldGaitsgory, BZFN} provide convenient starting points and the articles~\cite{indcoh, ArinkinGaitsgory} develop the key properties of ind-coherent sheaves and the notion of singular support, which interpolates between the two sheaf theories. For the convenience of the reader we provide a very brief synopsis.

\subsubsection{Quasicoherent sheaves}
\index{quasicoherent sheaves}
The big category $QC(X)$ of quasicoherent sheaves on any derived stack is defined as follows. We first define $QC$ as a functor on affine schemes by assigning $QC(S)=\cO(S)\module$ with pullback $f^*$ given by tensor product. We can then use right Kan extension to define $QC(X)$ for any stack (or prestack) $X$, i.e., $$QC(X)=\lim_{\leftarrow} QC(S)$$ where the limit is over the category of affines over $X$:
``a quasicoherent sheaf on $X$ is a (star-pullback-) compatible system of quasicoherent sheaves on affines mapping to $X$.'' 

 By construction we have pullback functors $f^*:QC(Y)\to QC(X)$ for any morphism $f:X\to Y$.
For $X$ a QCA (quasicompact with affine automorphism groups and finitely presented classical inertia) algebraic stack of finite type, the compact objects $QC(X)^c=\Perf(X)$ are given by the perfect complexes (locally representable by finite complexes of vector bundles). Moreover $QC(X)$ is dualizable, satisfies the K\"unneth formula $QC(X\times Y)\simeq QC(X)\otimes QC(Y)$ and (as a result) is canonically self-dual. For arbitrary \index{QCA}  QCA morphisms $f:X\to Y$ the functor $f^*$ has a colimit-preserving right adjoint $f_*$ satisfying base change and the projection formula. Under the more stringent condition (satisfied in most common stacks in characteristic zero, in particular every stack we'll encounter coherent sheaf theory on) that $X$ is perfect we also have that $$QC(X)\simeq \Ind(\Perf(X)),$$ i.e., $QC(X)$ is compactly generated by the perfect complexes.  Finally, we recall that $*$-tensor product of sheaves ($*$-pullback of external tensor products along the diagonal) endows $QC(X)$ with a symmetric monoidal structure, with unit the structure sheaf $\cO_X$ (the $*$-pullback of $\kk$ from a point), satisfying the projection formula. 

\subsubsection{Ind-coherent sheaves}\index{ind-coherent sheaves}\label{IndCoh section}
For a scheme of finite type $X$, we have the familiar small category $\Coh(X)$ of coherent sheaves (bounded coherent complexes), a full subcategory of $QC(X)$. We have an inclusion $\Perf(X)\subset \Coh(X)$ which is an equivalence precisely when $X$ is smooth. The corresponding large category $$QC^!(X)=\Ind(\Coh(X))$$ is the category of ind-coherent sheaves. By construction it comes with a unique colimit-preserving functor $\Psi_X:QC^!(X)\to QC(X)$ extending the inclusion $\Coh(X)\to QC(X)$. This functor is an equivalence on the bounded below subcategories with respect to standard t-structures, and can be reconstructed purely from the t-structure on $QC^!(X)$ as left completion. For $X$ a {\em bounded} (eventually coconnective -- i.e., the structure sheaf is supported in a finite number of cohomological degrees) derived scheme this functor is essentially surjective and exhibits $QC(X)$ as a colocalization of $QC^!(X)$. Moreover $QC^!(X)$ can be recovered from $QC(X)$ with its t-structure as {\em anti-completion}~\cite[C.5.5.]{SAG}. The ``difference'' between quasicoherent and ind-coherent sheaves, i.e., the kernel of $\Psi_X$, sits in cohomological degree $-\infty$, i.e., all the cohomologies of an object in the kernel of $\Psi_X$ vanish.

Ind-coherent sheaves are extended to stacks in parallel fashion to the definition of $QC$. 
First, ind-coherent sheaves on quasicompact schemes enjoy a continuous $!$-pullback functor and $*$-pushforward functor, which form an adjoint pair $(p_*,p^!)$ for proper morphisms and satisfy base change. In fact they naturally assemble together to form a functor out of the correspondence category of schemes~\cite[Part III]{GR}. As with $QC$ and $*$-pullbacks, we can then right Kan extend $QC^!$ (now with $!$-pullbacks) to all prestacks locally of finite type. The full correspondence formalism also extends to this setting, where we can $!$-pullback along any morphism and $*$-pushforward along any quasicompact schematic morphism. Moreover $(p_*,p^!)$ adjunction holds for any (ind-)proper morphism. 

The correspondence formalism automatically encodes a symmetric monoidal structure on $QC^!$ -- the $!$-tensor product $\otimes^!$, defined as $!$-pullback to the diagonal of the external product. Moreover this monoidal structure satisfies the projection formula (see~\cite[Section 2, Introduction to Part III]{GR}). An important further structure on $QC^!$ is that of module over $(QC,\otimes)$~\cite[II.6]{GR}. The functors $\Psi$ is naturally $QC$-linear, as is pushforward -- one of two ``mixed'' forms of the projection formula for the action of $QC$ on $QC^!$ (see~\cite[Proposition I.4.3.3.7]{GR} and \S \ref{Koszul dual via volume forms}). We will use the notation $\otimes$ for both the tensor structure on $QC$ and its action on $QC^!$ and the notation $\otimes^!$ for the tensor of $QC^!$.

For any QCA algebraic stack, $QC^!(X)=\Ind(\Coh(X))$ is compactly generated by coherent complexes, satisfies the K\"{u}nneth formula and is canonically self-dual~\cite{DrinfeldGaitsgory}. This self-duality provides a natural general formulation of Serre duality, and makes $!$-pullbacks and $*$-pushforwards dual functors. It carries a symmetric monoidal structure, the $!$-tensor product, with unit the dualizing sheaf $\omega_X$ (the $!$-pullback of $\kk$ from a point). It also carries the structure of $QC(X)$-module category compatible with pushforwards (via the projection formula). Acting on $\omega_X$ defines a functor $\Upsilon_X:QC(X)\to QC^!(X)$, which is dual to the colocalization $\Psi_X$ when the latter makes sense (for schemes or algebraic stacks). 

\subsubsection{Singular support}\label{singular support sec}
\index{singular support of coherent sheaves}
The theory of singular support of coherent sheaves on quasi-smooth\footnote{Quasi-smooth stacks are the ``derived lci'' stacks -- for our purposes they are the quotients by affine groups of schemes with tangent complex of amplitude $[0,1]$, i.e., locally isomorphic to the (derived) fiber of a map of affine spaces. } stacks~\cite{ArinkinGaitsgory} allows one to quantify and control the difference between perfect complexes and coherent sheaves, or equivalently between their ``large'' counterparts, quasicoherent and ind-coherent sheaves.

Informally speaking, the notion of singular support of a coherent sheaf $\cF\in \Coh(X)$ is a microlocal measure of singularity of sheaves (see Section~\ref{cohmic}): it records not only the points $x$ where $\cF$ is not perfect, but also the 1-shifted codirections $\xi\in H^{-1}(T^*_xX)$ where this failure occurs. 
In the context of the Koszul dual algebras $\Lambda = \mathcal O(\AA^1[-1])$, $S^\shear=H^*_{\Gm}(\pt)$, as in \S \ref{shearing affine line}, the singular support of a coherent $\Lambda$-module is the support of the corresponding $S^\shear$-module.

For a general quasi-smooth stack, the assignment $x\mapsto \xi\in H^{-1}(T^*_xX)$ forms a classical stack, the \emph{stack of singularities} $\Sing(X)$. We study local deformations of $X$ near $x$, governed by $H^1$ of the tangent complex at $x$, and record for which ones $\cF$ is obstructed from deforming, as evidenced by a corresponding class in $\Ext^2(\cF,\cF)$.  This defines a conical closed subset of $\Sing(X)$.

Conversely, given $\Lambda\subset \Sing(X)$ we can consider
 $\Coh_\Lambda(X)$, the category of coherent sheaves with singular support contained in $\Lambda$, and its ind-category $QC^!_\Lambda(X)$ which sits between
  $QC(X)=QC^!_{\{0\}}(X)$ and $QC^!(X)=QC^!_{\Sing(X)}(X)$.

\subsection{Topological sheaf theories on finite type schemes}\label{constructible sheaf theories}
On the automorphic side of the Langlands correspondence we will make use of variants of categories of constructible sheaves on schemes and stacks, as reviewed in~\cite[Appendix A]{GKRV} and~\cite[Appendix E,F,G]{AGKRRV1}. As a general convention we use $\Shv$ to denote small categories of sheaves of constructible nature, and $\SHV$ for corresponding large categories. Our notation is designed so that for finite type schemes $X$ (but not for general finite type stacks!) the large and small categories recover each other by passing to ind-objects  and to compact objects:
 \begin{equation}\label{ind sheaf equation} \SHV(X)=\Ind \Shv(X) \mbox{ and } \Shv(X)=\SHV(X)^c.
 \end{equation}
 We generally refer to objects of $\Shv(X)$ as {\em finite sheaves} and objects of $\SHV(X)$ as {\em ind-finite sheaves}.
 
We will consider three types of topological sheaf theories: de Rham, constructible and Betti. In each case, the sheaf theories have a microlocal aspect -- they admit a notion of singular support in the cotangent bundle. In fact, {\bf our convention is that ``Betti'' always refers to sheaves with Lagrangian singular support}, as discussed below. 

 \index{$\Shv$} \index{category of sheaves}

\subsubsection{Constructible/\'etale} \label{CST:ET}
 Let $\kk=\overline{\QQ_\ell}$ and $X$ be a scheme of finite type over an algebraically closed field $\FF$ with characteristic different than $\ell$. In this setting we take $\Shv^{et}(X)$ to refer to (the derived dg-category of) bounded constructible complexes of $\ell$-adic \'etale sheaves. (See~\cite{gaitsgorylurie} for an $\infty$-categorical treatment of $\ell$-adic \'etale sheaves.) 
Its ind-category is the large category $\SHV^{et}(X)=\Ind(\Shv^{et}(X))$ of ind-constructible \'etale sheaves. 

For a scheme of finite type $X/\C$ and any $\kk$, we can consider the dg derived category of constructible complexes of sheaves of $\kk$-vector spaces in the analytic topology. We will also refer to this category as the ``\'etale'' category of sheaves, $\Shv^{et}(X)$, although ``constructible'' is a more standard name, because the two can be treated simultanesously, in the context of the present paper. We let $\SHV^{et}(X)=\Ind(\Shv^{et}(X))$ the corresponding category of ind-constructible sheaves.

\'Etale sheaves (in either sense) enjoy the full six-functor formalism, in particular have adjoint pairs $(f^*,f_*)$ and $(f_!,f^!)$ of functors for arbitrary morphisms. 

In both settings there is a notion of singular support, which is a conical Lagrangian $\Lambda\subset T^*X$ and one can consider full subcategories of sheaves $\Shv_\Lambda(X)$ and $\SHV_\Lambda(X)$ with prescribed singular support.
 In the case of $\ell$-adic sheaves in positive characteristic, the notion of singular support is the one coming from Beilinson's definition in~\cite{BeilinsonSingularSupport}. 

\subsubsection{de Rham}  \label{CST:DR} Let $X$ denote a scheme of finite type over $\kk$. The de Rham model of sheaf theory is given by $D$-modules: $\SHV^{dR}(X)=\cDD(X)$ denotes the large category of all quasicoherent $D$-modules on $X$ (see~\cite{DrinfeldGaitsgory}). This category is most familiar (and realized as modules for the sheaf of differential operators) for $X$ smooth, but can be defined for any $X$ either by embedding $X$ as a closed subscheme of a smooth scheme or intrinsically as ind-coherent sheaves on the de Rham space $X_{dR}$. It is compactly generated, with compact objects $\Shv^{dR}(X)=\cDD^{coh}(X)$ forming the derived category of bounded coherent complexes of $D$-modules. 
Thus we have $\cDD(X)\simeq \Ind(\cDD^{coh}(X))$. Note that coherence is (as usual) taken in the sense of $D$-modules, so that e.g.\ $\cDD$ itself is coherent, though it is far from coherent as an $\OO$-module. 

$D$-modules have $f^!$ and $f_*$ functoriality in general, with $(f_*,f^!)$ adjoint for proper morphisms and $(f^!,f_*)$ adjoint up to a shift for smooth morphisms.

Among these compact objects we find holonomic $D$-modules, which give constructible sheaves under the ``solutions'' functor of the Riemann--Hilbert correspondence. (The RH functor restricts to an equivalence on holonomic $D$-modules with regular singularities.)
They enjoy the full six functor formalism, with the same formal properties as in the \'etale sheaf theories above. 

For example, fixing a conical Lagrangian\footnote{Though we will only need smooth schemes $X$, the theory of singular support for each of the constructible sheaf theories discussed extends naturally to singular $X$, see~\cite[E.6]{AGKRRV1}.} $\Lambda\subset T^*X$ we have the full subcategory $\cDD^{coh}_\Lambda(X)\subset \cDD^{coh}(X)$ of coherent $D$-modules with singular support (or characteristic variety) contained in $\Lambda$, which are in particular holonomic, and its ind-category which is a full subcategory $\cDD_\Lambda(X)\subset \cDD(X)$. 

Holonomic $D$-modules are far from generating all $D$-modules, and the existence of coherent $D$-modules such as $\cDD_X$ itself (with singular support all of $T^*X$) encodes useful phenomena such as algebraically varying families of connections.

\subsubsection{Betti}\label{Betti sheaves}\label{CST:B}
Finally for a scheme of finite type $X/\C$ and any $\kk$, we have the rather wild large category $\SHV^{all}(X)$ of {\em all} sheaves of $\kk$-vector spaces on the underlying topological space $X^{an}$ of $X$ in the (Hausdorff) complex analytic topology. This theory is far less familiar in algebraic geometry than the de Rham and \'etale variants; we merely summarize the main facts the reader can find in~\cite[Appendix G]{AGKRRV1}. This category is not compactly generated, and for us only
 plays a role similar to the role played by this paper to most readers:  a giant storage bin in which to find objects of interest. One can specify nice classes of sheaves by picking a stratification of $X$, or by the closely related method~\cite{KashiwaraSchapira} of fixing the allowed singular support of sheaves -- the codirections on $X$ outside of which we require sheaves to be locally constant. Namely, for any conical Lagrangian $\Lambda\subset T^*X$ there is a full subcategory $\SHV^{B}_\Lambda(X)\subset \SHV^{all}(X)$ of sheaves with singular support contained in $\Lambda$, which we refer to as Betti sheaves. Such sheaves are automatically locally constant along an associated stratification, whose union of conormals contains $\Lambda$. More generally we define Betti sheaves to be sheaves with singular support contained in some (conic, algebraic) Lagrangian $\Lambda$, i.e.,
\begin{equation} \label{Bettisheafdef} \SHV^{B}(X)=\lim_{\rightarrow, \Lambda} \SHV^{B}_\Lambda(X).\end{equation}
This is a compactly generated presentable dg category which we consider as a less wild storage bin that in particular contains all constructible sheaves on $X$.

\index{$\SHV^{B}$} \index{Betti sheaves} \index{$\SHV^{all}$}

An illustrative example comes by requiring zero singular support, i.e., $\Lambda$ is the zero section of $T^*X$. In this case we find $\SHV^{B}_{\{0\}}(X)=LOC(X)$, the large category of locally constant sheaves on $X$, representations of the fundamental $\infty$-groupoid of $X$ (a derived refinement of the familiar categories of representations of the fundamental group). For $X$ connected, a compact generator for this category can be given by choosing a point $x\in X$ and taking the pushforward of the constant sheaf under the path fibration $P_x\to X$ (in degree zero this is the ``universal cover'' local system). This sheaf is a locally constant replacement for the skyscraper sheaf $\kk_x$ -- it is obtained by applying to $\delta_x$ the left adjoint to the inclusion of locally constant sheaves into all sheaves. Note, however, that it is not a finite rank local system (i.e., it is not a constructible sheaf). 
\index{locally constant versus constructible}
 For example for $X\simeq T$ a torus we find $$LOC(X)=\kk[\pi_1(T)]\mbox{-modules},$$ and its compact objects are
 $$Loc(X) \simeq \mbox{finitely presented $\kk[\pi_1(T)]$-modules}. $$  This is in contrast to the small and large categories of locally constant sheaves in the \'etale (here meaning complex constructible) setting, which correspond to finite and locally finite $\kk[\pi_1(T)]$-modules, respectively.
 Cartier dually, finite rank local systems are given by coherent sheaves with finite support on the dual torus $T^\vee$;
  compact local systems $\Loc(X)$ correspond to all of $\Coh(T^\vee)$,
  and the large category  $\LOC(X)$ corresponds to $QC(T^{\vee})$. 
  In other words, studying
 $LOC(X)=\SHV_{\{0\}}^{B}(X)$ and its compact version $\Loc(X) =\Shv_{\{0\}}^{B}(X)$ allows us to consider local systems whose monodromy varies algebraically in families, while restricting to constructible objects only allows formal deformation of monodromies. 

In general the categories of Betti sheaves $\SHV^{B}_\Lambda(X)$  {\em are} compactly generated, and we let $\Shv^{B}_\Lambda(X)$ denote the corresponding category of compact objects $$\SHV^{B}_\Lambda(X)\simeq \Ind(\Shv^{B}_\Lambda(X)).$$ 
Compact generators can be given explicitly by enforcing the prescribed singular support on skyscrapers on strata (applying the left adjoint to the inclusion $\SHV^{B}_\Lambda(X)\subset \SHV^{all}(X)$,
compare discussion above in the case when $\Lambda$ is the zero-section.). As for local systems (or in the $D$-module setting), constructible sheaves with singular support $\Lambda$ give compact objects of $\Shv^B_\Lambda(X)$ (because of the finite type properties of the homotopy types of strata), but they are far from generating the category -- the compact objects will typically restrict to infinite rank local systems on strata.

\subsection{Topological sheaf theories on finite type stacks}\label{sheaves on stacks}
Next we would like to extend our different flavors of sheaf theories to more general finite type (pre)stacks (though {\bf we only consider Betti sheaves on {\em algebraic} stacks}).

First we can use the extension paradigm (as in \S\ref{synopsis}) in the de Rham or \'etale settings to define a large category of sheaves on any prestack locally of finite type, the {\em safe category} (or category of ind-safe sheaves) $\SHV_{s}(X)$, by right Kan extension over $!$-pullbacks:
 we set $$\SHV_{s}(X)=\varprojlim_{  f:Y\to X} (\SHV(Y),f^!)$$ as the limit over affine schemes of finite type $Y$ mapping to $X$.
 \footnote{Moreover this description guarantees that $\SHV_{s}(X)$ is dualizable, with dual given by the colimit
$$\SHV_{s,co}(X)\simeq \varinjlim_{ f:Y\to X} (\SHV(Y),f_*)$$ over $*$-pushforwards, the dual functors to $!$-pullbacks.}  In the next subsection, we will define a different ``large'' category, that we will denote by $\SHV(X)$, by changing the class of compact objects.

In the de Rham settting,  $\SHV_{s}^{dR}(X)=\cDD(X)$ is the standard big category of all $D$-modules on a stack.
 On the other hand, in the Betti setting of $\SHV^{all}$ we don't generally have a colimit-preserving $!$-pullback. 

There are two main motivations for the use of $!$- rather than $*$-pullback. One is practicality in de Rham setting, which only has $f^!$ and $f_*$ functoriality. A more substantial one is the desire to have a good sheaf theory for ind-schemes $X=\Cup_j X_j$, where the limit of categories over functors $i^!$ for closed embeddings is identified with the colimit $$\SHV_{s}(X)\simeq \varinjlim_j (\SHV(X_j), i_*)$$ over the left adjoint functors $i_*$ (by Remark~\ref{calculating colimits}). This allows one e.g.\ to show the categories are compactly generated with compact objects coming by extension from the finite type subschemes $X_j$. 

In the \'etale setting, where we always have a $(f_!,f^!)$-adjunction, we can apply the same argument in general to write $\SHV_{s}(X)$ as a colimit over $!$-pushforwards (which are defined for schematic morphisms). This establishes $\SHV_{s}(X)$ as compactly generated by $!$-pushforwards of constructible sheaves from affines. Also, since $f^!$ preserves constructibility, we can define the (small) category of \'etale sheaves on an arbitrary prestack $X$ as the limit
$$\Shv(X)=\varprojlim_{  f:Y\to X} (\Shv(Y),f^!)$$  
In other words, a sheaf on $X$ is a system of sheaves on affines over $X$ compatible under $!$-pullback.

\subsubsection{Algebraic stacks} Now let us narrow our focus and assume $X$ is an algebraic stack (so that $X$ has a smooth cover by an affine) with affine diagonal (so that the pullback of affines is affine). In this case we can replace the index category of all affines of finite type over $X$ with that of affines which are {\em smooth} over $X$ and with only smooth morphisms. For such morphisms $!$-pullbacks of Betti sheaves are colimit-preserving (and agree with $*$-pullbacks up to shift). Moreover the resulting definition of sheaves (in any of our flavors of sheaf theory) as a right Kan extension gives equivalent categories whether we use $!$- or $*$-pullbacks (since the two differ by shifts for smooth morphisms). See~\cite[Appendix G.7]{AGKRRV1} for Betti sheaves on stacks, including the good behavior (in particular compact generation) of the categories of Betti sheaves with fixed singular support.

Now observe that $!$-pullbacks for smooth morphisms preserve compact objects, since $f^!$ has a continuous right adjoint, a shift of $f_*$. (Note this is not the case for arbitrary morphisms in the de Rham setting, e.g., restriction of $D$-modules along $pt\hookrightarrow \AA^n$ takes $\cDD$ to an infinite dimensional vector space.)
As a result we can define the (standard) small category $\Shv^{dR}(X)=\cDD^{coh}(X)$ of coherent $D$-modules on an algebraic stack following the general format
$$\Shv(X)=\varprojlim_{  f:\Spec(R)\to X\mbox{ smooth}} (\Shv(Y),f^!).$$   
In other words, a sheaf on $X$ is a system of sheaves on affines smooth over $X$ compatible under $!$-pullbacks.
  
 \subsection{Finiteness, renormalization and safety}\label{renormalization section}
 \index{renormalized category}
We now come to a fundamental issue about sheaf theory on algebraic stacks which often goes by the (somewhat unfortunate) name renormalization (see in particular~\cite[Appendix F.5]{AGKRRV1}.) This is an instance of the issue of completion / decompletion commonplace in equivariant topology.

Namely, we have defined both large and small categories of sheaves $\SHV_{s}(X)$ and $\Shv(X)$ as a limit over smooth atlases (in all three settings). However, the finite objects (constructible sheaves, coherent $D$-modules or Betti sheaves that are compact on smooth covers) $\Shv(X)\subset \SHV_{s}(X)$ are {\em not} in general compact objects. For example, for $X=BG$, $\SHV_{s}(X)$ (in any of our sheaf theories) is identified with modules for $H_*(G)$.  In this category
 the constant sheaf on $BG$, which corresponds to the augmentation module for the exterior algebra $H_*(G)$, is not a compact object / perfect complex of modules. 
 More generally equivariant constructible sheaves are not compact in general. The compact objects 
 $$ \Shv_s(X):=\SHV_{s}(X)^c\subset \Shv(X)$$ form a full subcategory of finite sheaves, the {\em safe} sheaves introduced in~\cite{DrinfeldGaitsgory}.

 One can fix this by ``formally declaring'' finite objects to be compact, i.e., passing to the ind-category:
 
 \begin{definition} \label{renormalized sheaf definition} For a quasicompact algebraic stack $X$ we define the category of {\em ind-finite sheaves} as the ind-category
  $$\SHV(X):=\Ind(\Shv(X)),$$
  so that $\SHV(X)^c\simeq \Shv(X)$. 
 For $X$ an arbitrary algebraic stack locally of finite type we define the category $\SHV(X)$ of ind-finite sheaves by right Kan extension, i.e., as a limit over $\SHV(U)$ over quasicompact open substacks. 
 \end{definition}

The ind-finite sheaf category automatically comes with a functor 
\begin{equation} \label{unrendef} \operatorname{safe}:\SHV(X)\to \SHV_{s}(X) \end{equation} which is a colocalization, precisely analogous 
  to the functor
\begin{equation} \label{unrendef2}   
  QC^!(Y)=\Ind(\Coh(Y)) \to QC(Y)=\Ind(\Perf(Y))
\end{equation}
induced from the inclusion $\Coh(Y)\hookrightarrow QC(Y)$,   on a stack  
   of finite type. As in that setting, $\SHV(X)$ and $\SHV_{s}(X)$ differ only ``in cohomological degree $-\infty$'' with respect to the standard $t$-structure. (In fact one expects that the two categories can be formally obtained from each other by manipulations -- left completion and anti-completion, respectively -- of $t$-structures.) 
  
In the example $X=BG$, ind-sheaves $$\SHV(BG)\simeq H^*(BG)\module$$ recover the Koszul dual picture to $\SHV_{s}(BG)\simeq H_*(G)\module$, with the constant sheaf corresponding to the compact generator given by the regular module for $H^*(BG)$. (See \S \ref{shearing affine line} or~\cite{DrinfeldGaitsgory} for further discussion of this example). 
More generally, for quasicompact quotient stacks $X=Y/G$ the difference between the theories is captured in the support theory of (compact) objects as modules over $H^*(BG)$: ind-safe sheaves give torsion modules (supported at $0$) while the ind-finite sheaves have arbitrary support. 

Ind-finite sheaves are arguably the natural choice in the constructible world, where the focus is usually on the small sheaf categories and we just define the large categories formally by passing to inductive limits. On the other hand for $D$-modules we typically start from all $D$-modules and then impose finiteness conditions such as coherence, and there ind-safe sheaves recover the standard form of the large category (as defined e.g.\ in~\cite{BD}).

\subsection{Sheaf Theory in Infinite Type}\label{infinite type}

We now discuss some formal properties of constructible  sheaf theories on schemes and stacks of infinite type, following~\cite{raskininfinite}, see also~\cite[Section 4]{BKV}. 
Our primary application for this material in the main text is to sheaf theory on $X_F/G_O$ for $X$ a smooth affine spherical variety,
as in \S \ref{section-unramified-local}.

We will use the theory of $!$-{\em sheaves}, defined using the $!$-pullback functors.\label{infinite type scheme}
As discussed in
\S \ref{sheaves on stacks}, this
 is well-adapted to the de Rham and \'etale settings but not to the Betti setting of ``all'' sheaves, which don't admit continuous $!$-pullbacks for general maps.

Recall we have a (contravariant) functor $X\mapsto \SHV(X), f\mapsto f^!$ from the category of schemes of finite type. 
For $X$ a scheme of infinite type, we define $\SHV^!(X)$ as the left Kan extension of this functor, i.e., $\SHV^!(X)$ is the colimit of $\SHV(U)$ for finite type schemes  $X\to U$ under $X$, under $!$-pullbacks. Unfortunately, in general $f^!$ doesn't have a continuous right adjoint, so this colimit cannot be accessed concretely by rewriting it as a limit over right adjoints.

\index{shriek sheaves}
\index{! sheaves}
The resulting sheaf theory automatically comes with a colimit-preserving pullback functor $f^!:\SHV^!(Y)\to \SHV^!(X)$ for any map $f:X\to Y$, and a symmetric monoidal structure, the $!$-tensor product, from pullback along diagonal maps. Moreover, the formalism of~\cite{GR} can be used to enhance the functor $\SHV$ to a functor out of the correspondence category, as shown in the de Rham setting in~\cite[Section 3]{raskininfinite} (though the arguments apply in the \'etale setting as well). In other words, we also have $*$-pushforward functors -- which we only consider for proper morphisms -- satisfying base change. For proper morphisms $f$, we also have the $(f_*,f^!)$ adjunction. In general we don't have a form of Verdier duality, though we discuss in Section~\ref{appendix placid setting} the enhanced features of sheaf theory for {\em placid} ind-schemes of infinite type such as $G_F$ for $G$ a reductive group scheme and $X_F$ for $X$ a vector space. 

Given the definition of sheaf theory on arbitrary schemes, we now right-Kan-extend along $!$-pullback to define a sheaf theory $\SHV^!_{s}$ for arbitrary prestacks, equipped with $!$-pullbacks and left adjoint $*$-pushforwards for proper maps satisfying base change. We refer to the resulting objects as ind-safe sheaves.

\subsubsection{$*$-sheaves} \index{star sheaves} There is a ``dual'' theory of sheaves in infinite type, the $*$-sheaves~\cite{raskininfinite}.  
It is defined in terms of the dual functor $f_*$ of $f^!$. Namely, we have a (covariant) functor $X\mapsto \SHV(X), f\mapsto f_*$ from the category of schemes of finite type.
We then define the functor $\SHV^*, f_*$ on schemes of infinite type as the {\em right} Kan extension, i.e., $\SHV^*(X)$ is the limit of $\SHV(U)$ for finite type schemes under $X$, $X\to U$, under $*$-pushforwards. In other words, a $*$-sheaf is a system of sheaves on finite type approximations of $X$, compatible under $*$-pushforward. 
As noted in~\cite[Prop.3.19.1]{raskininfinite}, thanks to the (Verdier) duality between $f^!$ and $f_*$ in finite types it follows that  {\em if} $\SHV^!(X)$ is dualizable then its dual is given by $\SHV^*(X)$. Moreover, $*$-sheaves enjoy a correspondence formalism and a $\otimes$-action by $!$-sheaves satisfying a strong form of the projection formula.

\subsubsection{Ind-finite categories in infinite type}
The discussion above concerns large categories of sheaves on infinite-type schemes. If we wish to study small categories we run into the difficulty in the de Rham setting that $!$-pullback does not preserve coherent $D$-modules for morphisms that aren't smooth, for example, the inclusion of a point in a scheme; and that $*$-pushforward doesn't have a left adjoint $f^*$ in general. Thus we will now restrict our attention to constructible sheaves (either \'etale or Betti) or {\em holonomic} $D$-modules, all of which are preserved by $!$-pullbacks.
\footnote{Note that in our intended applications to $X_F/G_O$ for $X$ spherical {\em all} coherent $D$-modules are holonomic, so this restriction is harmless.}

In this constructible setting $!$-pullback preserves finiteness, so we can define (by left Kan extension again) a functor $\Shv^!$ of finite (i.e., constructible) sheaves with the full package of functoriality enjoyed by the large categories $\SHV^!_{s}$. Moreover, the $(f^*,f_*)$-adjunction lets us rewrite $\SHV^*(X)$ as the colimit of $\SHV(U)$ under $*$-pullbacks. As a result, $\SHV^*(X)$ for a quasicompact scheme is automatically compactly generated by $*$-pullbacks from finite type schemes (hence in particular dualizable). As a result, its dual $\SHV^!(X)\simeq Ind(\Shv^!(X))$ is also compactly generated (by $!$-pullbacks). 
This compact generation allows us to extend the $(f_*,f^!)$-adjunction and base change from proper maps to {\em ind-}proper maps, using the general extension machinery of~\cite[Theorem I.7.3.2.2]{GR} (i.e., by defining $f_*$ as the left adjoint of $f^!$)\footnote{We are indebted to Harold Williams for helpful remarks on sheaf theory in infinite type.}.

\subsubsection{Placid setting}\index{placid scheme} \label{appendix placid setting} We now recall the notion of placidity of a scheme $X$ in infinite type, which is a very strong form of the notion that the singularities of $X$ are finite dimensional (again following~\cite{raskininfinite}, as well as~\cite{DrinfeldInfinite}, see also the earlier~\cite[Definition 3.2.4]{KapranovVasserot}). Namely, a {\em placid presentation} of $X$ is an identification $X\simeq \varprojlim_{  i} U_i$ as a filtered inverse limit of finite type schemes under {\em smooth, affine} transition maps. We say $X$ is placid if it admits a placid presentation. For instance, if $X$ is pro-smooth, so that all the $U_i$ are themselves smooth, then $X$ is in particular placid.

On placid schemes we have a form of Verdier duality. Since for a smooth morphism $f^!$ forms a left adjoint of $f_*$ (i.e., agrees with $f^*$) {\em up to a shift}, placid schemes allow for a very tight relation between $!$- and $*$-sheaves, in which we absorb the (infinite!) shifts into the definition. Namely, for $X$ placid there's a canonical object 
$$\omega_X^{ren}\in \Shv^*(X),$$
the {\em renormalized dualizing sheaf}. It can be described as a suitable shift of the $*$-pullback of the dualizing sheaf of any of the $U_i$ in the placid presentation. In particular if $X$ is pro-smooth then the renormalized dualizing sheaf is simply the (unshifted!) constant sheaf $\omega_X^{ren}=\kk_X:=p^*\kk$ (for $p:X\to pt$).

The tensor action of $!$-sheaves on $*$-sheaves (defined in general) now results in an equivalence~\cite[Prop.4.8.1]{raskininfinite}
$$act_{\omega^{ren}}:\SHV^!(X)\to \SHV^*(X)$$ 
given by acting on the renormalized dualizing sheaf. In particular the equivalence identifies the dualizing sheaf $\omega_X=p^!\kk\in \SHV^!(X)$ with the renormalized dualizing sheaf $\omega_X^{ren}\in \SHV^*(X)$ (i.e., the constant sheaf in the pro-smooth setting). 

The notion of {\em placid morphism} of placid schemes is introduced in~\cite[Sec.4.10]{raskininfinite}, as a morphism which factors through smooth coverings on placid presentations. For such a morphism, one has a $*$-pullback functor and the equivalence of $!$- and $*$-sheaves intertwines the $!$- and $*$-pullback functors
~\cite[Prop.4.11.1]{raskininfinite} -- i.e., the equivalence absorbs the dimension shifts relating $!$- and $*$-pullback for smooth morphisms.

\subsection{Duality and Tensor products of sheaf categories}\label{tensor products of sheaves} We now collect some general facts about duality and tensor product theorems for categories of sheaves.

\subsubsection{Tensor products}
To translate between geometry and category theory it is often essential to know if the canonical tensor product functor $$\Shv(X)\otimes \Shv(Y)\to \Shv(X\times Y)$$ is an equivalence, in which case we say $\Shv$ satisfies the tensor product theorem in this setting.\footnote{A more structured version of this statement is the assertion that $\Shv$ forms a strict, rather than merely lax, symmetric monoidal functor out of a suitable category of stacks.}

Coherent sheaf categories and $D$-modules typically satisfy the tensor product theorem. 
 As explained in~\cite[Section 4.2]{DrinfeldGaitsgory}, such a tensor product theorem in any sheaf theory $\Shv$ follows formally from the combination of two statements: the tensor product result for {\em affine} schemes, and the compact generation, or more generally dualizability, of the categories of sheaves on the factors.
The dualizability of $QC$ and $\cDD$ on affines is automatic from their descriptions as categories of modules, and for $QC^!$ it is~\cite[Proposition 4.6.2]{indcoh}. Moreover $QC, QC^!$ are dualizable on QCA stacks (quasicompact with affine automorphism groups and finitely presented classical inertia)~\cite[4.2]{DrinfeldGaitsgory}\index{QCA}, as are $D$-modules on quasicompact stacks. Hence the tensor product theorems hold in these settings.

Compact generation (hence the tensor product theorem) of $\cDD(X)$ is established in~\cite{DrinfeldGaitsgorycompact} for a class of non-quasicompact stacks including the crucial case of moduli stacks $Bun_G(C)$ of bundles on curves. The relevant class are {\em truncatable} stacks $X$, which are those covered by open quasicompact substacks $U$ for which the inclusion $i_!$ (the would-be left adjoint of restriction) is well defined on all of $\cDD(U)$. In this case $\cDD(X)$ is compactly generated by such $!$-pushforwards of compact objects on opens.

In topology, tensor product theorems are far rarer. The enormous category $\SHV^{all}(X)$ of all sheaves on a locally compact Hausdorff topological space satisfies the tensor product theorem by~\cite[Theorem 7.3.3.9, Prop. 7.3.1.11]{HTT}. This stands in stark contrast to categories of (ind-)constructible sheaves of different flavors, which (despite being compactly generated) essentially never obey tensor product theorems: the K\"unneth formula implies that we have a full embedding $$\Shv(X)\otimes \Shv(Y)\hookrightarrow \Shv(X\times Y),$$ however the constant sheaf on the diagonal is rarely in the essential image (cannot be resolved by external powers of constructible sheaves on the factors). The Tensor Product Theorem of~\cite{AGKRRV1} for \'etale sheaves with nilpotent singular support on $\Bun_G$ is a striking exception, see Section~\ref{tensor and miraculous}

\subsubsection{Self-Duality} \label{self-duality}
Let us recall some results about self-duality for sheaf categories, see~\cite{AGKRRV2} or older references (e.g.~\cite{DrinfeldGaitsgory} or~\cite{gaitsgorykernels}). A self-duality $\Shv(X)^\vee\simeq \Shv(X)$ is uniquely specified by either its unit $u\in \Shv(X)\ot \Shv(X)$ or its counit $c:\Shv(X)\ot \Shv(X)\to Vect_\kk$. Thus it suffices to present a suitable functor out of $\Shv(X\times X)$ (which receives a functor from $\Shv(X)\ot \Shv(X)$), or, in the presence of tensor product theorems, to specify a suitable sheaf on $X\times X$.

Morally, self-duality comes from the diagonal: the diagonal correspondence $$\xymatrix{pt &\ar[l]^-{\pi} X\ar[r]^-{\Delta} & X\times X}$$ and its opposite present any space as self-dual in the correspondence category, so a suitably functorial linearization is canonically self dual. In other words, the natural candidates for units are versions of the constant sheaf on the diagonal, and for counits are versions of global sections of restriction to the diagonal. 

Depending on the sheaf theory, one finds two general flavors of self-dualities. In the better-known one, the unit is $\Delta_*\pi^! \kk=\Delta_*\omega_X$ the dualizing sheaf on the diagonal and the counit is $\pi_*\Delta^!(\cF\boxtimes \cG)=\Gamma(X,\cF\otimes^! \cG)$, the cohomology of the $!$-tensor product of sheaves. Serre duality for ind-coherent sheaves is of this form~\cite[4.4]{DrinfeldGaitsgory}, as is Verdier duality for $D$-modules or constructible sheaves on a quasicompact scheme. 
A similar statement holds for $D$-modules or constructible sheaves on quasicompact stacks, except that the global sections functor has to be replaced by a colimit-preserving version, the {\em renormalized} global sections. This is by definition the unique colimit-preserving functor agreeing with global sections on compact objects. (Concretely, for a quotient stack $X=Y/G$ this means we push forward to $BG$ and then take {\em homology} -- tensoring with the trivial sheaf -- rather than cohomology -- Hom from the trivial sheaf.)  

On the other hand, we have another of self-duality, with unit $\Delta_! \pi^*\kk$ the ($!$-extended) constant sheaf on the diagonal and counit $\pi_!\Delta^*(\cF\boxtimes \cG)=\Gamma_c(X,\cF\otimes^* \cG)$ the compactly supported cohomology of the $*$-tensor product. This describes the self-duality of the category of all sheaves $\SHV^{all}(X)$ on a locally compact Hausdorff space. In the setting of $D$-modules, a scheme, a quasicompact or truncatable stack is said to be {\em miraculous} if $\Delta_! \pi^*\kk$ is the unit of a self-duality. As the name suggests, this is quite rare, and encodes a form of homological smoothness. The miraculous duality for $\Bun_G$ ~\cite{DrinfeldGaitsgorycompact, gaitsgorystrange} is a striking exception, see Section~\ref{tensor and miraculous}.

\subsection{Sheaves of categories, ULA and rigidity}\label{shvcats}
We briefly review some categorical notions that will be needed in the next section.

 For the purposes of this section it will be important to work with the ``large'' versions of categories of sheaves.
 We will work in either de Rham or constructible (e.g., \'etale) sheaf theory, so that we have a symmetric monoidal category $(\SHV(M), \ot^!)$ of sheaves equipped with the $!$-tensor product.
 We will require only a coarse ``affinized'' version of the notion of sheaf of categories on a scheme $M$:
 
 \begin{definition} A {\em sheaf of categories} over $M$ is a $(\SHV(M),\otimes^!)$-module category $\cC_M\in \mathsf{DGCAT}_{\mathbf k}$.
 \end{definition}

This notion is well-adapted to $(f_!,f^!)$-functoriality rather than $(f^*,f_*)$; indeed a better name might be $!$-sheaf of categories. In particular, for a closed embedding $i:Z\hookrightarrow M$ let $\cC_Z=\cC_M\otimes_{\SHV(M)} \SHV(Z)$, the induced sheaf of categories over $Z$. Then the adjunction $(i_!,i^!)$ on sheaves induces an adjunction which we also denote $(i_!,i^!)$ between $\cC_Z$ and $\cC_M$.

\begin{remark}[Quasicoherent sheaves of categories] In the de Rham setting, as in~\cite{raskinchiral}, we can use Gaitsgory's 1-affineness theorem to identify sheaves of categories in this coarse sense with honest sheaves of categories $U\mapsto \cC(U)\in D(U)\module$ which are {\em quasicoherent}. 
\end{remark}

We now consider the tensor product of sheaves of categories (see \S \ref{tensor products of sheaves} for a discussion of tensor product theorems).
The assignment $\SHV$ defines a lax symmetric monoidal functor from stacks over $M$ to sheaves of categories over $M$. Concretely, for $Z\to M$ we have a sheaf of categories $\SHV(Z)\in \SHV(M)\module$ on $M$, and external product defines a functor
\begin{equation} \label{prodfunctor} \SHV(Z)\otimes \SHV(W)\longrightarrow \SHV(Z\times_M W), \end{equation} of $\SHV(M)$-modules. 

Unlike in the $D$-module setting, this functor \eqref{prodfunctor} fails to be an equivalence in the constructible world.  However, it does enjoy a weak variant of the K{\"u}nneth theorem, namely the comparison maps above are {\em fully faithful}. In other words, while sheaves on a product are not generated by external products of sheaves on the factors, the morphisms between external product sheaves are given as external products.

\subsubsection{The ULA condition} \index{ULA condition}
 We give a brief description of the ULA condition in its categorical formulation, see~\cite[Appendix B]{raskinCPS1} and~\cite[Appendix D]{GKRV} for related treatments and~\cite[Appendix A.2]{xinwenSatake},~\cite[IV.16]{reich} for more traditional treatments.

\begin{definition} Given a monoidal category $\cA$ and a module category $\cM$, an object $\cF\in \cM$ is said to be {\em universally locally acyclic} (or ULA) over $\cA$ if the functor $act_\cF:\cA\to \cM$ given by acting on $\cF$ has an $\cA$-linear colimit-preserving right adjoint. In this case the algebra object $act_{\cF}^R act_\cF(1_\cA)\in \cA$ is denoted by $\underline{End}(\cF)$, the {\em internal endomorphisms} of $\cF$ in $\cA$.
\end{definition}

The notion of ULA object derives from that of a ULA sheaf on a space $X$ with respect to a morphism $p:X\to Y$ (taking $\cA$ to be sheaves on $Y$ and $\cM$ to be sheaves on $X$); see \cite[Arcata, \S V]{SGA4.5} for this notion in its original algebro-geometric context. The notion of ULA object is preserved by colimit-preserving $\cA$-linear morphisms of $\cA$-module categories. In particular it follows from adjunctions and the projection formula that the ULA property for sheaves is preserved by smooth pullbacks and proper pushforwards of spaces over a fixed base $Y$ (see also~\cite[Appendix A.2]{xinwenSatake})

\subsubsection{Rigidity}\index{rigid tensor category} \label{rigid tensor categories}
We recall that the standard notion of ``rigidity' for a Tannakian category asserts that (in a small-category setting) objects have duals, 
which permits one to define internal $\Hom$.
We  now recall a corresponding notion of rigidity for a monoidal category in our setting,  and one of its main features, as exposed in~\cite[Section 1.9]{GR}, adapted to the setting of sheaves of categories. \footnote{In fact this formulation the notion and its main features are developed in a fashion readily adaptable to a general symmetric monoidal 2-category -- they do not refer to {\em objects} of the category, but only adjunction properties of morphisms.}
We work relative to some symmetric monoidal category ${\mathcal R}$ (e.g.\ in the $\kk$-linear setting we would take ${\mathcal R}=\Vect_\kk$, or for sheaves of categories over $M$ we would take $\SHV(M)$.)

\begin{definition} Let ${\mathcal R}$ denote a symmetric monoidal category.  A {\em rigid monoidal category over ${\mathcal R}$} is an algebra object $(\cC,\ast)\in ({\mathcal R}\module,\otimes_{\mathcal R})$ in ${\mathcal R}$-module categories for which 
\begin{enumerate}
\item the unit morphism ${\mathcal R}\to \cC$ has a colimit-preserving right adjoint, i.e., the unit $1_\cC\in \cC$ is ULA over ${\mathcal R}$; and 
\item the multiplication $\ast:\cC\ot\cC\to \cC$ has a colimit-preserving, $\cC$-linear right adjoint.
\end{enumerate} 
\end{definition}

\begin{proposition}\label{rigid ULA}
Fix $(\cC,\ast)$ rigid over ${\mathcal R}$ and $\cM$ any $\cC$-module category.
\begin{enumerate}
\item The action $act:\cC\ot\cM\to \cM$ has a colimit-preserving, $\cC$-linear right adjoint. 
\item For an object $\cF\in \cM$ over ${\mathcal R}$ (i.e., an ${\mathcal R}$-linear functor ${\mathcal R}\to\cM$), $\cF$ is ULA over ${\mathcal R}$ if and only if it is ULA over $\cC$.
\end{enumerate}
\end{proposition}

For example, if ${\mathcal R}$ is itself rigid then $\cF$ is ULA if and only if it is compact.

%% file: geometric-Langlands.tex
 \newcommand{\res}{r}

\section{The geometric Langlands correspondence}
\label{geometric Langlands}

In this section we briefly describe the different sheaf theoretic settings for the geometric Langlands correspondence. The conjectures all have the following general form:

\begin{quote}
  A full subcategory $\AUT_{\text{?`}}^{?}(\Bun_G(\Sigma))$ of sheaves on the stack of $G$-bundles is identified with a category of ind-coherent sheaves on a stack of Langlands parameters $\QC^!_{\text{?`}}(\Loc_\Gv^{?}(\Sigma))$, compatibly with actions of Hecke functors.
\end{quote}

Here $\Sigma$ is a smooth projective curve over a field $\FF$ and the categories are linear over a field $\kk$ of characteristic zero. We will give a very brief
sketch of the notation now, and then proceed to more detailed definitions.

\subsubsection{Cheat sheet}  \label{Cheatsheat}
 \begin{itemize}
\item 
$\SHV$ or $\Shv$ denotes all sheaves on $\Bun_G$,
defined by ``general purpose'' definitions as in \S 
\ref{sheaf theory}.  

\item $\AUT$ or $\Aut$ denotes the category of {\em automorphic sheaves},
the ``largest subcategory on which it is reasonable to study Hecke actions,''
see \S \ref{automorphic categories} for explicit definition and \S \ref{spectral action} for discussion.

\item There are adornments $?= dR, B, et$ and ${\text{?`}} = \mathcal N, s$ for the various categories, explained below. If we write  $\AUT$ or $\Aut$ without adornment, 
it means that one should take $? $ to be $dR, B, et$
according to the context of the usage, 
 and take ${\text{?`}}$ to be empty.
\end{itemize}

\subsubsection{? = Betti, de Rham or \'etale}  \label{BdRR}

There are three ``flavors'' of the geometric Langlands conjecture -- that is to say, three possibilities for the $?$
that appears above:  

\begin{itemize}
\item[$\bullet$]  de Rham~\cite{BD,ArinkinGaitsgory}, denoted by
$? = dR$; here  $\FF=\kk=\CC$.
\item[$\bullet$] Betti~\cite{BettiLanglands}, denoted by $?=B$; here
 $\FF=\CC$, $\kk$ arbitrary.

\item[$\bullet$] \'etale~\cite{AGKRRV1}, denoted by $?=et$; 
it makes sense in any  sheaf-theoretic context, and in particular both
 $\FF=\kk=\CC$, as well as $\FF$: of positive characteristic and $\kk$: $\ell$-adic, are admissible.
\end{itemize}

\subsubsection{$\text{?`}$ = $\mathcal{N}$ or $s$ -- nilpotency or safety conditions}
There is another parameter, the $\text{?`}$, that we can vary in formulating the geometric Langlands conjecture in each of its flavors, which has to do with how our sheaf theories treat singularities on the spectral side and stackiness on the automorphic side. 
Namely, on the spectral side we can allow all ind-coherent sheaves $\QC^!(\Loc_\Gv^{?}(\Sigma))$ or consider only sheaves with nilpotent singular support $\QC^!_\cNN(\Loc_\Gv^{?}(\Sigma))$as in~\cite{ArinkinGaitsgory}.
On the automorphic side, this corresponds to allowing all ind-finite sheaves (the ``renormalized'' automorphic category, which we denote simply $\AUT^{?}(\Bun_G)$) or restricting to the ``safe setting'' of ind-safe sheaves (which we denote $\AUT_{s}^{?}(\Bun_G)$), see \S \ref{sheaves on stacks}, \ref{renormalization section} for the definitions, and~\ref{safe Langlands} for further discussion. 

We emphasize that the definition of the automorphic categories $\Aut$, $\AUT$ \emph{already includes a ``nilpotent singular support'' condition in the Betti and \'etale settings}, see \S \ref{automorphic categories}. This condition is imposed by the spectral decomposition (see \S \ref{spectral action}), and is not related to the nilpotent support condition on the spectral side.

\begin{remark}[Extended groups and spin structures] We refer the reader to \S \ref{extended-group appendix} 
for a discussion of formulations in a way that does not depend on choices of spin structure. 
 \end{remark}
\subsection{Automorphic side} \label{Autside}

The automorphic categories in the de Rham, Betti and \'etale conjectures consist of different classes of sheaves on the same stack $\Bun_G(\Sigma)$ of $G$-bundles on the curve $\Sigma$.
We review the sheaf theory and then define the automorphic categories.

\subsubsection{The automorphic categories}\label{automorphic categories}

The sheaf theories on the automorphic side considered are ``constructible sheaf theories'' as reviewed in \S \ref{sheaf theory} following~\cite[Appendix A]{GKRV} and~\cite[Appendices E,F,G]{AGKRRV1}.
As we have mentioned, there is an important subtlety in formulating the geometric Langlands conjecture in each of its flavors, which on the automorphic side has to do with how our sheaf theories treat stackiness -- whether we allow all ind-finite sheaves (the ``renormalized'' sheaf category, which we denote simply $\SHV^{?}(\Bun_G)$) or restrict to the ``safe setting'' of ind-safe sheaves (which we denote $\SHV_{s}^{?}(\Bun_G)$). 
{\em Our default is to work with the larger ind-finite categories and to restrict to safety when necessary, see \S\ref{safe Langlands} below.}

Let us now describe $\AUT$ and $\SHV$ for each value of $?$, freely using the generalities
of \S \ref{sheaves on stacks} and \S \ref{renormalization section} to define the categories $\SHV$
on the $\FF$-stack $\Bun_G$. 

\begin{itemize}
\item
The de Rham automorphic category, when $\FF=\CC$, \[\AUT^{dR}(\Bun_G(\Sigma)):= \SHV^{\dR}(\Bun_G(\Sigma)) = \cDD(\Bun_G(\Sigma)),\] consists of all (ind-coherent) $D$-modules, i.e., there is no distinction between $\AUT$ and $\SHV$ in the de Rham setting. 

 \item The Betti automorphic category, when $\FF=\CC$,  
 \begin{multline*} \AUT^{B}(\Bun_G(\Sigma)):=\SHV^{B}_\cNN(\Bun_G(\Sigma)) \\
 \hookrightarrow \SHV^B(\Bun_G(\Sigma))  \mbox{ defined as in 
 \eqref{Bettisheafdef},}   \end{multline*}
  consists   
  of all (renormalized) sheaves of $\C$-vector spaces on the underlying topological stack, in the complex topology, whose singular support is contained in the global nilpotent cone $\cNN\subset T^*\Bun_G(\Sigma)$, the zero-fiber of the Hitchin map. Note that the global nilpotent cone is Lagrangian \cite{Faltings, Ginzburg-nilpotent}, forcing these sheaves to have cohomology that is locally constant along the strata of an associated stratification, explicitly described in \cite{BD} (see also \cite[\S D.3]{AGKRRV1}). 
  
\item The \'etale automorphic category  \begin{multline*} 
\AUT^{et}(\Bun_G(\Sigma)):=\SHV^{cons}_\cNN(\Bun_G(\Sigma)) \\  \hookrightarrow \SHV^{et}(\Bun_G(\Sigma)) =  \SHV^{cons}(\Bun_G(\Sigma)) \mbox{ as in 
\S \ref{CST:ET}} \end{multline*}
is the category of ind-constructible sheaves of $\kk$-vector spaces with nilpotent singular support,
which is a full subcategory of the category of ind-constructible sheaves, without singular support conditions. 
  The papers of Arinkin, Gaitsgory, Kazhdan, Raskin, Rozenblyum and Varshavsky develop the theory of the \'etale automorphic category of nilpotent sheaves\footnote{The main results are proved under the assumption on the characteristic of the ground field that there exists a non-degenerate G-equivariant pairing bilinear form on $\fg$, whose restriction
to the center of any Levi subalgebra remains non-degenerate.} $\SHV_\cNN(\Bun_G)$ 
(For example in the case $G=\Gm$, the nilpotent cone is simply the zero section and $\SHV_\cNN(\Bun_G)$ is the category of locally constant sheaves.) 
\end{itemize}

\begin{remark}
We recall from \S \ref{Betti sheaves} the substantial distance between Betti sheaves $\SHV^{B}_\cNN(\Bun_G(\Sigma))$ and nilpotent ind-constructible sheaves $\SHV^{et}_\cNN(\Bun_G(\Sigma))$: both are locally constant along the same stratification, but the compact objects in the former need not have finite rank cohomology sheaves. For instance in the case of $G=\Gm$, the nilpotent cone is the zero section and we are in the setting illustrated in {\it loc. cit.}, with Betti sheaves giving  \emph{arbitrary} locally constant sheaves on the Picard group of $\Sigma$, while \'etale sheaves correspond to \emph{locally finite} representations of the fundamental groups of components. 
\end{remark}

Why restrict to nilpotent singular support in the Betti and \'etale settings? There are several concrete answers based on convenience, matching with examples and experience in geometric representation theory (going back to Harish-Chandra's study of distributional characters and Lusztig's theory of character sheaves). A ``first-principles'' answer is provided by the results of~\cite{AGKRRV1} 
 on the spectral action (the ``converse to the Nadler-Yun theorem''), namely these categories are universally characterized by the requirement that Hecke functors depend in a locally constant way on points of $\Sigma$,  see \S \ref{spectral action} below.

Over $\FF=\kk=\CC$, the \'etale automorphic category maps naturally to both $D$-modules (landing in ind-coherent $D$-modules) by forgetting the singular support condition and applying the Riemann-Hilbert correspondence, and to the Betti automorphic category, by forgetting ind-constructibility. Crucially, as proved\footnote{Note that there are some assumptions behind that theorem which may not be satisfied when the characteristic of $\FF$ is ``small'' compared to $G$, see op.cit.\ \S 14.4.1.} in~\cite[Theorem 14.4.4]{AGKRRV1}, the \'etale category contains all Hecke eigensheaves in any of its ambient categories, and as such is a suitable ``core'' for the geometric Langlands correspondence.

\subsection{Spectral side} \label{spectral side section}

The spectral categories in the de Rham, Betti and \'etale conjectures are given by applying the same sheaf theory (of coherent nature)
but on different versions of the stack of local systems. 
 We first describe the different versions of the stack, and then discuss the sheaf theory. The stack of local systems is, in general, derived, and this has to be taken into account for the coherent theory.

\begin{itemize}
 \item The de Rham stack $\Loc_\Gv^{dR}(\Sigma)$ is the moduli stack of flat $\Gv$-connections on $\Sigma$, equivalently $\Gv$-torsors over the de Rham functor $\Sigma_{dR}$ or $\otimes$-functors $$\Rep(\Gv)\longrightarrow \cDD(\Sigma)$$ from representations of $\Gv$ into $D$-modules on $\Sigma$. That is to say, $R$-points of this stack are given by tensor functors to $R$-modules in $\cDD(\Sigma)$, $R\module\otimes \cDD(\Sigma)$.

\item  The Betti stack $\Loc_{\Gv}^{B}(\Sigma)$ is the moduli stack of ($R$-families of) locally constant $\Gv$-torsors on $\Sigma$, equivalently $\Gv$-torsors over the underlying homotopy type $\Sigma_{top}$ of $\Sigma$, representations into $\Gv$ of the fundamental $\infty$-groupoid of $\Sigma$ (or just the fundamental group when $\Sigma$ has positive genus), or $\otimes$-functors $$\Rep(\Gv)\longrightarrow \Loc(\Sigma)$$ from representations of $\Gv$ to locally constant sheaves on $\Sigma$.

\item The stack of local systems of restricted variation $\Loc^{et}_{\Gv}(\Sigma)$ introduced in~\cite{AGKRRV1} (which we denote with ``et'' because it matches the \'etale sheaf theory on the automorphic side), parametrizes $\otimes$-functors taking finite dimensional representations of  $\Gv$ to {\em locally finite} representations of the  ($\infty$-)fundamental group(oid). (This can be expressed in terms of big categories as $\otimes$-functors from $Rep(\Gv)$ to R-modules in ``quasi-lisse'' local systems, the t-completion of the ind-category of finite rank local systems on $\Sigma$.)
In particular this ensures that the semisimplification of the resulting local systems are (locally) constant in families, whence ``restricted variation.''
\end{itemize}

In fact~\cite[Theorem 4.8.4]{AGKRRV1} establishes that, over $\C$, $\Loc^{et}_{\Gv}(\Sigma)$ is the disjoint union of the formal completions of the Betti space $\Loc^B_{\check G}(\Sigma)$ over semi-simple local systems (which form the coarse moduli space or affinization of $\Loc$) -- indeed, when $\FF=\CC$, we have embeddings
$$\Loc_\Gv^{dR}(\Sigma) \hookleftarrow \Loc^{et}_{\Gv}(\Sigma) \hookrightarrow \Loc^B_{\Gv}(\Sigma)$$
so that restricted local systems form the ``common core'' for the de Rham and Betti spaces. The definition of $\Loc^{et}$ applies equally well in other sheaf theories, in particular for $\ell$-adic local systems, and the general structure of $\Loc^{et}$ is similar -- it is a disjoint union of functors that are relative algebraic stacks over formal affine schemes.
 
Next we discuss the category of sheaves. 
On the spectral side, all the forms of the geometric Langlands conjecture concern the same sheaf theory -- namely the category 
$$ QC^!(\Loc_{\Gv}^{?}(\Sigma))$$ of ind-coherent sheaves on the 
various versions of the stack of $\Gv$-local systems on $\Sigma$. In the ``safe'' version of the conjecture we further restrict to ind-coherent sheaves with {\em nilpotent singular support} (see \S \ref{singular support sec} for the notion of singular support).

The condition of nilpotent singular support  can be defined in all three settings as follows.
 All versions of the stack of local systems have tangent complex described as the cohomology of the associated adjoint local system, shifted by $1$. In particular the ($-1$)-st cohomology of the cotangent complex is identified with locally constant sections of the adjoint local system, and thus comes with a characteristic polynomial map to $\fgv\sslash \Gv$,with fiber over $\{0\}$ giving the spectral analog of the global nilpotent cone.  For example, for $\Gv=\Gm$, coherent sheaves on $\Loc_{\Gv}$ with nilpotent support are just perfect complexes.

\subsection{Unramified Geometric Langlands Conjecture} \label{uglc}

The geometric Langlands conjecture refines the spectral action, by seeking to precisely describe the automorphic category as a sheaf of categories over $\Loc_\Gv$. Recall that $\AUT^{?}(\Bun_G(\Sigma))$ refers to the full category of (ind-coherent) $D$-modules in the de Rham setting and the categories of (all or ind-constructible) sheaves with nilpotent singular support in the Betti or \'etale settings.

\medskip
 
\begin{conjecture}\label{GLC} 
 Let $?$ denote Betti, de Rham or \'etale setting.
\begin{enumerate}
\item[$\bullet$] (\cite{ArinkinGaitsgory},\cite{BettiLanglands},~\cite[Conjecture 21.2.7]{AGKRRV1})\label{restricted GLC}
There is a Hecke-equivariant equivalence of categories   
$$\AUT^{?}_{s}(\Bun_G(\Sigma))\simeq QC^!_{\cNN}(\Loc^{?}_\Gv(\Sigma))$$ between ind-safe automorphic sheaves on $\Bun_G$ and ind-coherent sheaves with nilpotent support on local systems, in each of the three settings.
\item[$\bullet$] More generally, there is a Hecke-equivariant equivalence of categories 
$$\AUT^{?}(\Bun_G(\Sigma)) \simeq QC^!(\Loc^{?}_\Gv(\Sigma))$$ 
between ind-finite automorphic sheaves on $\Bun_G$ and ind-coherent sheaves on local systems, in each of the three settings.
\end{enumerate}
\end{conjecture}

We will explain how the ind-safe version can be recovered from the ind-finite one in \S \ref{safe Langlands} below.

 \begin{remark}
Very recently,  a proof of the ind-safe version of the de Rham and Betti geometric Langlands correspondence has been announced by
 D. Arinkin, D. Beraldo, J. Campbell, L. Chen, D. Gaitsgory, J. Faergeman, K. Lin, S. Raskin and N. Rozenblyum. 
 See \cite{Gaitsgorypage} and references therein. 
 \end{remark}

\begin{remark} \label{GLCsuper}
In certain settings (see \S \ref{supersloppy} as well as~\cite[Remark 5.4.6]{BD}) it is desirable 
to also consider the super-version of the conjecture, which is as follows:  
\begin{itemize}
\item
On both sides of the conjecture, we consider sheaves of super-$\kk$-vector spaces.
\item Even sheaves on the automorphic side  correspond to sheaves on the spectral side
whose parity coincides with the action of $e^{2\rho}(-1)$.
\end{itemize}
Note that this does not affect the underlying categories or their module structure for the Hecke categories.
\end{remark}

\subsubsection{Safety, renormalization and singular support}\label{safe Langlands}  \index{renormalized category} 
As we have mentioned there is an important subtlety in formulating the geometric Langlands conjecture in each of its flavors, which on the automorphic side has to do with how our sheaf theories treat stackiness -- whether we allow all ind-finite sheaves (the ``renormalized'' sheaf category, which we denote simply $\SHV^{?}(\Bun_G)$) or restrict to the ``safe setting'' of ind-safe sheaves (which we denote $\SHV_{s}^{?}(\Bun_G)$). 

At the time of writing, all of the literature on (and evidence for) the geometric Langlands correspondence concerns the safe version. 
For the purposes of this paper, however, it is somewhat more natural to work with the stronger ind-finite conjecture (since $L$-sheaves don't naturally have nilpotent singular support). Thus, our default is not to impose equivariant support conditions automorphically or singular support conditions spectrally, while noting that all our statements have nilpotently projected/safe counterparts.

Recall that the category of ind-finite sheaves $\SHV(X)$ contains ind-safe sheaves $\SHV_s(X)$ as a full subcategory, as does $QC^!(Y)$ contain $QC^!_\Lambda(Y)$ for any singular support condition $\Lambda$. The pairs of categories are related by a colocalization and differ only in cohomological degree $-\infty$. Moreover, the difference between the two flavors can be measured by a support condition with respect to $Z=H^*(BG)$, with the smaller (ind-safe/nilpotent support) categories characterized by support at the origin. To see this, one can pick 
a point $x\in \Sigma$, obtaining actions of the Hecke category $\HECKE_G$ on both automorphic and spectral categories, hence an action of the endomorphism ring of the unit $$Z:=End(1_{\HECKE_G})\simeq H^*(BG) \simeq (\cO^\shear(\fgxv))^{\Gv}.$$
(This action can be described automorphically and spectrally in terms of the induced map $\Bun_G(\Sigma)\to BG$ and the presentation of $\Loc_\Gv(\Sigma)$ as a derived fiber of the stack of local systems with ramification allowed at $x$, respectively).
 Moreover we have:
 
 \medskip

\begin{proposition}\cite[Proposition 12.7.3]{ArinkinGaitsgory}
 The action of the (ind-finite) Hecke category $\overline{\mathcal{H}_G} \simeq \QC^{\shear}(\fgxv/\Gv)$ on $QC^!(\Loc_\Gv(\Sigma))$ preserves $QC^!_{\cNN}(\Loc_\Gv(\Sigma))$, and its restriction there factors through the colocalization functor $\QC^{\shear}(\fgxv/\Gv) \to \QC^{\shear}_{\cNN}(\fgxv/\Gv)$ of \eqref{unrendef2}. Vice versa, the action of the full subcategory $\QC^{\shear}_{\cNN}(\fgxv/\Gv)$ maps $QC^!(\Loc_\Gv(\Sigma))$ to $QC^!_{\cNN}(\Loc_\Gv(\Sigma))$.
\end{proposition}

In other words, $QC^!_{\cNN}(\Loc_\Gv(\Sigma))$, as a full subcategory of $QC^!(\Loc_\Gv(\Sigma))$, is characterized by its support as a $Z$-module, for any chosen point of $\Sigma$. 

A similar result holds on the automorphic side (see ~\cite[Remark 12.8.8]{ArinkinGaitsgory} for the support property for $D$-modules on $Bun_G$  and~\cite[F.5.6]{AGKRRV1} for the corresponding characterization of renormalization for quotient stacks\footnote{To apply this directly to $Bun_G$ we need to work on quasicompact opens and replace the reductive $G$ by $G_O/G_O^{(n)}$ for some congruence subgroup.} $Y/G$). Hence, the ind-finite version of the Geometric Langlands Conjecture \ref{GLC} strictly implies the ind-safe one.

From the perspective of topological field theory (cf. Appendix~\ref{fieldtheorysec}), the algebra $Z$ is the $E_4$ algebra of local operators in the theories $\cA_G\simeq \cB_\Gv$ (which however is in fact strictly commutative). The $Z$-linearity of the Langlands correspondence (and hence the role of nilpotent singular support on the spectral side) is interpreted in~\cite{ElliottYoosingular} as the dependence of  4d $\cN=4$ Yang-Mills theory on its Coulomb branch parameters.

 \begin{remark}[Nilpotent support and duality]
We reiterate that the appearances of nilpotent support on the two sides of the geometric Langlands correspondence do {\em not} correspond under duality.
Indeed, nilpotence on the spectral side corresponds -- automorphically -- to safety (the fact that ind-constructible sheaves or $D$-modules on $\Bun_G$ are torsion for the action of the ring $H^*(BG)$ of equivariant parameters), and this condition can be removed by ``renormalizing''.

By contrast, nilpotence on the automorphic side, in the Betti and \'etale settings, is forced on us by the requirement that the Hecke action on the spectral side factors through $QC(\Loc^?_{\check G})$, i.e., that Hecke functors vary locally constantly along the curve. We will discuss this spectral action in \S \ref{spectral action}.
\end{remark}

\begin{remark}[Projecting to nilpotent support]
Period sheaves on $\Bun_G$ are not themselves nilpotent in general, so don't naturally lie in the Betti or \'etale categories of automorphic sheaves. However, if we are interested in periods of automorphic forms, or geometrically in Hom pairings between eigensheaves and period sheaves, we are implicitly studying the period sheaf only as a functional on automorphic sheaves, or equivalently considering only its image under the spectral projector (also to be discussed in \S \ref{spectral action} below). 

Analogously, many natural sheaves on local systems, in particular the $L$-sheaves that we study in this paper, do not have nilpotent singular support. One could similarly apply a (much less dramatic) nilpotent projection functor to them, but it seems more natural not to do so and instead work with the larger ind-finite (``renormalized'') version of the geometric Langlands correspondence, which accommodates all ind-coherent sheaves on $\Loc_\Gv$. 

For further discussion, see \S \ref{nilpotent projection is good}.
\end{remark}

\begin{remark}[Compatibility with abelian duality and fluxes]\label{fluxes}
We record here a basic compatiblity between the Langlands correspondence and abelian duality for the center of $G$, known in physics as duality between electric and magnetic fluxes (see~\cite[\S 7.2]{KapustinWitten}). We don't know references in the mathematical literature. 

Namely, we can twist $G$-bundles by bundles for the center of $G$, giving rise to a translation action of $\Bun_{Z(G)}$ on $\Bun_G$ and hence on the automorphic categories (the action of magnetic fluxes). We will encounter this action only through the restricted action of $\Bun_{\Z/2}$ for a central involution $z:\Z/2\to Z(G)$, in the context of making statements independent of spin structures (see Remarks~\ref{Heisenberg comment},~\ref{Abelian double covers} and \S \ref{parity returns}). 

On the other hand we have a dual homomorphism $\Gv\to BZ(G)^{\vee}$, arising
from  a central extension of $\Gv$ by $Z(G)^\vee$;
 for semisimple groups this map comes from the identification $Z(G)^\vee=\pi_1(\Gv)$.
  Passing to stacks of local system we obtain a $\Bun_{Z(G)^\vee}$-torsor classifying lifts to this extension over $\Loc_\Gv$. We now appeal to the abelian duality (Poincar\'e-Pontrjagin or Weil pairing) 
$$\Bun_Z\times \Bun_{Z^\vee}\to B\Gm$$ to define associated line bundles on $\Loc_\Gv$ to $Z(G)$-torsors on $\Sigma$, i.e., a homomorphism
$$\Bun_{Z(G)}\to \Pic(\Loc_\Gv).$$ This defines a tensor product action (by electric fluxes) of $\Bun_{Z(G)}$ on $QC^!(\Loc_\Gv)$. 
Again, we will encounter this action only through the restricted action of $\Bun_{\Z/2}$ associated to a double cover $\Gv_z$ classified by $z^\vee:\Gv\to B\Z/2$.

The assertion is that these two actions are identified under the geometric Langlands correspondence -- indeed, they are identified unconditionally under the spectral action (i.e., the spectral action of  $\Bun_{\Z/2}\to QC(\Loc_\Gv)$ on the automorphic category agrees with the translation action). 
This is a consequence of the geometric Satake correspondence, specifically of its effect on translation by $Z(G)$.
\end{remark}

\subsection{Spectral action and the spectral projector}\label{spectral action}

The spectral action, or ``automorphic-to-spectral'' direction of the Langlands correspondence, establishes a sheafification or spectral decomposition of the automorphic side over the corresponding stack $\Loc_\Gv$ of Langlands parameters. Namely for every $x\in \Sigma$ we have an action of $\Rep(\Gv)$ by Hecke functors on sheaves on $Bun_G$, equipped with factorization structure (compatibility as the points vary and collide). In each of the three settings, a spectral action theorem asserts that this action factors through an action of quasicoherent sheaves on the stack of Langlands parameters: 

\begin{itemize}
\item de Rham: Gaitsgory's vanishing theorem~\cite[Theorem 4.5.2]{gaitsgoryoutline} asserts that the spherical Hecke action descends to an action (the {\it spectral action}) of quasi-coherent sheaves 
$QC(\Loc_{\check{G}}^{dR})$ on the de Rham space of $\Gv$-connections on 
$\AUT_s^{\dR}(\Bun_G)$. 
 
 \item Betti: Nadler and Yun~\cite{NadlerYunSpectral} (see also~\cite{GKRV}) proved that $QC(\Loc^B_\Gv)$ acts on $\AUT_s^{B}(\Bun_G)$. 
 
 \item  \'Etale: ~\cite{AGKRRV1} establishes an action of $QC(\Loc^{et}_{\Gv}(\Sigma))$ on the \'etale automorphic category $\AUT^{et}$. 
\end{itemize}

In each of these contexts, the spectral action sheafifies the automorphic category over the stack of Langlands parameters, identifying it as the global sections of a quasicoherent sheaf of categories (as studied in~\cite{1affine}) obtained by localization. Concretely, it means that for any two automorphic sheaves $\cF,\cG\in \Shv^{?}(\Bun_G)$ the Hom space $\Hom(\cF,\cG)$ localizes as a quasicoherent sheaf on $\Loc^{?}_\Gv$

In fact~\cite{AGKRRV1} establishes a much stronger form of the spectral action which characterizes the automorphic categories $\AUT\subset \SHV(\Bun_G)$. Namely, given any category $\cM$ with a factorizable action of $\Rep(\Gv)$ depending on points of $\Sigma$, they define the {\em spectrally decomposable} part of $\cM$
$$\iota:\cM^{spec}\hookrightarrow \cM.$$ This is the largest full subcategory of $\cM$ on which the $\Rep(\Gv)$ action is {\em locally constant} in $\Sigma$, i.e., factors through an action of $QC(\Loc^{?}_\Gv)$. Moreover, there is a canonical idempotent projector
\index{$\cM \mapsto \cM^{\spec}$} \index{spectral projector}  $$M\in \cM\mapsto \cM^{spec}\subset \cM^{spec},$$ the {\em Beilinson spectral projector}. This is a special case of a projector defined from $\cM$ to Hecke eigenobjects associated to any algebraic family of eigenvalues, here applied to the ``universal'' family $\Loc_\Gv$. The definition is a factorizable (or Ran-space) form of a general construction of projectors for modules over tensor categories, and is applicable in any of the sheaf theories.

In the de Rham setting, thanks to the spectral action encoded in Gaitsgory's vanishing theorem, we have 
$$\cDD(\Bun_G)^{spec}=\cDD(\Bun_G)$$ and the spectral projector is the identity. 

A major result of~\cite{AGKRRV1} establishes that in both the topological setting of all sheaves and the setting of ind-constructible sheaves, the spectrally decomposable parts of sheaves on $\Bun_G$ are precisely given by sheaves with nilpotent singular supports: 
\begin{equation} \label{agk0} \SHV^{et}_s(\Bun_G)^{spec}= \AUT^{et}_s(\Bun_G) (=\SHV^{et}_{s,\cNN}(\Bun_G)).\end{equation}
\begin{equation} \ \label{agk1} \SHV^B_s(\Bun_G)^{spec} =\AUT^B_s(\Bun_G) (=\SHV^B_{s,\cNN}(\Bun_G)).\end{equation}

This provides an intrinsic characterization of (and meaning for) nilpotent sheaves in terms of the Hecke action. 

\begin{remark} (Caveat about the spectral projector in safe versus ind-finite categories:) 
We anticipate that the same statement will be true also for the ind-finite categories, i.e., dropping
the subscript $s$.  In the text, we have allowed ourselves to use $P \mapsto P^{spec}$ in that setting without comment, on the assumption
that the corresponding results apply; it would be desirable to prove this. For the purposes of the main text, however, the statements can always
be ``projected'' into the safe category by the colocalization functor, as in \S \ref{safe Langlands}, so this sloppiness should not cause any essential problem. 
\end{remark}

The spectral projector $P\mapsto P^{spec}$ has very different properties in the Betti and \'etale settings, arising fundamentally from the distinct geometry of $\Loc$ in the two settings:
\begin{itemize}
\item[-] In the Betti setting $\Loc$ is closely modeled on affine schemes, and in particular its structure sheaf is a compact object. As a consequence it is proved in~\cite[Section 18]{AGKRRV1} that $(-)^{spec}$ provides a left adjoint to the inclusion of nilpotent sheaves in all sheaves (the {\em left nilpotent projection}), which exists in general from the theory of microlocalization (see~\cite[Section G.7]{AGKRRV1}). (The inclusion also has a continuous right adjoint.) 

\item[-] On the other hand, in the \'etale setting $\Loc$ is modeled on affine {\em formal} schemes, so that the structure sheaf is naturally a pro-object. There is correspondingly a pro-counterpart of the spectral projector~\cite[Section 17.1]{AGKRRV1} which is identified with the pro-left adjoint to the embedding of nilpotent sheaves. It is expected that the spectral projector $(-)^{spec}$ itself provides a {\em right} adjoint to the inclusion of nilpotent sheaves in ind-constructible sheaves. This is shown to be equivalent to~\cite[Conjecture 14.1.8]{AGKRRV1} 
that the subcategory $\SHV_\cNN(\Bun_G)\hookrightarrow \SHV(\Bun_G)$ is generated by compact objects that are compact in the ambient category, so that the right adjoint to the inclusion({\em right nilpotent projection}) is continuous. 
\end{itemize}

\subsection{Tensor product and self-duality}\label{tensor and miraculous}
We now discuss special algebraic properties of the automorphic sheaf categories, the tensor product and self-duality properties (see \S \ref{tensor products of sheaves} for a general discussion).

\subsubsection{Tensor product}
As one of the applications of the spectral projection, 
we have the following unexpected tensor product theorem for \'etale sheaves:

\begin{theorem}~\cite[Theorem 16.3.3]{AGKRRV1}\label{tensor product thm} For a pair of reductive groups $G,H$, there is an equivalence
$$\AUT^{et}_s(\Bun_{G\times H})\simeq \AUT^{et}_s(\Bun_G)\otimes \AUT^{et}_s(\Bun_H).$$ \end{theorem}

\begin{remark}
 The theorem is stated for the safe categories of automorphic sheaves. It would be desirable to establish the same for the categories of ind-finite sheaves. \ \end{remark}

\begin{remark}[Langlands for product groups]\label{Langlands for product groups}
The equivalence of Theorem~\ref{tensor product thm} and its (much easier) counterparts for $D$-modules on $\mathrm{Bun}$ and ind-coherent sheaves on $\Loc$ are all compatible with the spectral action of quasicoherent sheaves on the stacks $\Loc$ of Langlands parameters. (Indeed in the former case this action is used in the proof.)  Thus we can deduce the geometric Langlands correspondence  for $G\times H$ from those for $G$ and $H$ (in the de Rham or \'etale settings).\footnote{The tensor decomposition also respects Whittaker normalization, so we can expect a similar compatibility for the normalized geometric Langlands correspondence.}
In other words we expect a natural commutative diagram of equivalences (in both the ind-safe and ind-finite settings)
\begin{equation}\label{product Langlands}
\xymatrix{\AUT^{?}(\Bun_{G\times H})\ar[r]\ar[d]& QC^!_\cNN(\Loc_{\Gv\times \Hv})\ar[d]\\
 \AUT^{?}(\Bun_H)\otimes\AUT^{?}(\Bun_G) \ar[r] & QC^!_\cNN(\Loc_\Hv)\otimes QC^!_\cNN(\Loc_\Gv)}
 \end{equation}
\end{remark}

\subsubsection{Miraculous duality}
As recalled in \S \ref{tensor products of sheaves}, for quasicompact stacks, Verdier duality provides a canonical self-duality for both $D$-modules and ind-constructible sheaves, with unit $\Delta_*\omega$ and counit the (renormalized) global sections of the $!$-tensor product.\footnote{See
Appendix \ref{duality section} for a quick review of duality for categories.}
 For non-quasicompact but truncatable stacks such as $Bun_G$, Verdier duality fails to provide a self-duality~\cite{DrinfeldGaitsgorycompact}; rather it allows one to concretely identify the abstract dual of the category of sheaves (as the ``co-category'' $\SHV(X)_{co}$, the colimit of sheaf categories on quasicompact opens under $*$-pushforward). Nonetheless, all three flavors of the automorphic sheaf category are canonically self-dual by the {\em miraculous} duality, with unit the spectral projection of $\Delta_!\kk$ and counit the compactly supported cohomology of $*$-tensor product:

\begin{theorem}[Miraculous duality]\label{miraculous thm}
In either of the de Rham~\cite{gaitsgorystrange}, Betti~\cite[G.9.3]{AGKRRV1} and \'etale~\cite{AGKRRV2} settings, the object $\Delta_!\kk^{spec}$ is the unit for a self-duality, the {\em miraculous duality}, of the category of automorphic sheaves. \end{theorem}

See~\cite{drinfeldwang,wangbilinear} (especially~\cite[A.8-A.9]{drinfeldwang}) for a discussion of miraculous duality in relation to bilinear forms on automorphic forms in the classical setting of the Langlands correspondence. The corresponding duality can be described explicitly in terms of the {\em pseudo-identity functor} (see~\cite{gaitsgorykernels} for a detailed study in the de Rham setting and~\cite{AGKRRV2} for the subtler ``enhanced'' version in the \'etale setting). The pseudo-identity is 
the $!$-integral transform represented by $\Delta_!\kk$, $$\operatorname{PsId}_{Bun_G,!}:\cF\mapsto \pi_{2,\bullet}(\pi_1^!\cF\otimes^! \Delta_!\kk)$$ (where $\pi_{\bullet}$ refers to the renormalized pushforward, see \S \ref{self-duality}). For miraculous stacks this functor is an equivalence,\footnote{In the non-quasicompact setting we have the extra subtlety that this integral transform defines not an endofunctor but rather a functor $\SHV(X)_{co}\to \SHV(X)$ from the co-category of sheaves to the ordinary category of sheaves; this is the asserted equivalence.} and its composite with Verdier duality allows us to identify $\SHV$ with its dual.

\subsection{From geometric to arithmetic Langlands}\label{geometric to arithmetic Langlands}
The tensor product and self-duality theorems for the categories of nilpotent \'etale sheaves~\cite{AGKRRV1,AGKRRV2} open the way to defining and evaluating categorical traces of endofunctors. The main theorem of~\cite{AGKRRV3}, the ``Trace Conjecture,'' calculates the trace of Frobenius, asserting that automorphic functions are recovered precisely as the categorical trace of Frobenius on automorphic sheaves:

\medskip

\begin{theorem}~\cite{AGKRRV3}\label{trace conjecture}
The function-sheaf correspondence  
 induces an isomorphism 
$$\mbox{Frobenius trace }(\AUT^{et}(\Bun_G))\simeq  \mbox{compactly
supported functions on }\Bun_G(\F_q).$$
\end{theorem}

The trace of Frobenius acting on $\QC^!(\Loc^{et}_{\Gv})$ is identified with the space of algebraic distributions\index{algebraic distributions} (global sections of the dualizing sheaf) on the stack of arithmetic local systems (representations of the Weil group) on $\Sigma/\F_q$, $\Gamma(\Loc^{arith}_{\Gv}(\Sigma),\omega)$~\cite[24.6]{AGKRRV1}. 
Note that $\Loc^{arith}_{\Gv}$ is not quasi-smooth, and its dualizing complex is unbounded in positive cohomological degrees, while the structure sheaf is unbounded in the opposite direction -- in other words, distributions here behave very differently from functions, a distinction that is crucial for the geometric understanding of nontempered automorphic forms. 
(See e.g.\ Remark~\ref{AGKRRV Arthur} and \S \ref{derived volume forms}.)

It is also conjectured in {\it loc. cit.} that the same holds for sheaves with nilpotent singular support, $QC^!_\cNN(\Loc^{et}_{\Gv})$.  Thus Conjecture~\ref{restricted GLC}, in combination with Theorem~\ref{trace conjecture}, implies the following form of the arithmetic Langlands conjecture:

\begin{conjecture}\cite[Conjecture 24.8.6]{AGKRRV1} \label{geometric arithmetic Langlands}
There is an isomorphism 
$$\kk_c[\Bun_G(\F_q)]\simeq \Gamma(\Loc^{arith}_{\Gv}(\Sigma),\omega)$$
between unramified automorphic functions and algebraic distributions on the stack of arithmetic Langlands parameters (compatible with actions of unramified Hecke operators).
\end{conjecture}

 To normalize the isomorphism one must fix a spin structure on $\Sigma$ and a square-root of $q$; this dependence can be fixed by the use of extended groups as in  \S \ref{extended-group appendix}.

 \subsection{Extended groups and spin structures} \label{extended-group appendix} 
 
 \index{extended dual group} 
  The following discussion is not used in an essential way in the main text, but
 is referred to at several points in relation to making the statements more manifestly
 independent of choices of square roots.
 
 It is well-known in the classical Langlands program (see in particular~\cite{DelignetoSerre,BuzzardGee,Bernsteinhidden,xinwenSatake}) that many statements become cleaner in terms of an extended version of the Langlands dual group. The issue can be traced back for example to the need to choose a square root of the order $q$ of the residue field in order to identify the spherical Hecke algebra with the representation ring of the dual group. In other words the fundamental mechanism defining the spectral action, and hence pairing Langlands parameters with Hecke eigenvalues, needs to be modified to be made independent of choices. 
 
Another well-known phenomenon is the need to make some choices in order to normalize the Langlands correspondence. This usually appears in the form of picking Whittaker data in order to normalize Hecke eigenforms. In other words, given the spectral action we need to make choices to set up the correspondence between eigenforms and Langlands parameters.

Both of these issues can be corrected by replacing the groups $G$ and $\Gv$ by extended versions, which gives what appears to be the most symmetric 
  formulation of the Langlands correspondence.  
  Let ${}^C G$ be as defined in \S \ref{subsection-extended-group}:  the quotient of $G \times \GGm$ by the central element 
  $(e^{2 \check{\rho}}(-1), -1)$. We remind that $\GGm$ denotes the group $\Gm$, but we use different notation in order to refer to this distinguished instance of the group, and sometimes refer to it as the ``grading group''.

\begin{definition} 
\begin{enumerate}
\item The space of twisted $G$-bundles $\wt{Bun}_G(\Sigma)$ is the moduli stack of ${}^CG$-bundles $\cP$ equipped with an identification
$$\cP/G \simeq  \cK_\Sigma$$
of the associated $\GGm$-bundle (equivalently, line bundle) with the canonical bundle of $\Sigma$.
\item The space of twisted $\Gv$-local systems $\wt{Loc}_{\Gv}(\Sigma)$ is the moduli stack of ${}^C\Gv$-local systems $P$ equipped with an identification 
$$P/\Gv \simeq \varpi$$
of the associated rank one local system with the cyclotomic character.\footnote{The ``cyclotomic character'' is simply the (trivializable) orientation local system on the geometric curve, but when the curve is defined over $\FF_q$, we will also want to keep track of the action of Frobenius, in which case we think of this local system as the restriction of the cyclotomic character \eqref{cyclotomicdef} }
\end{enumerate}
\end{definition}

Thus, for example,  when dealing with $\SL_2$ on the automorphic side,
one deals instead with ``$\GL_2$-bundles with determinant the canonical class''
and when dealing with $\SL_2$ on the spectral side, one deals instead with
``$\GL_2$ local systems with cyclotomic determinant.''

The advantage of the moduli stack of twisted $G$-bundles is that it carries a canonical Whittaker sheaf
(i.e., without choosing a spin structure) 
and the advantage of the moduli stack of twisted $\Gv$-local systems is that the spectral action is defined 
without choices in a Frobenius-equivariant fashion (i.e., without choosing a square root of the 
cyclotomic character, which comes from the orientation sheaf). 
The twisted form of the geometric Langlands correspondence, Conjecture~\ref{GLC}, then predicts the following.

\medskip
\begin{quote} {\em Extended group formulation of the geometric Langlands equivalence:}    Let $?$ denote Betti, de Rham or \'etale setting. In each of the three settings, there is a equivalence of categories 
\begin{equation}\label{GLCextended}\AUT^{?}_{\invertQM}(\wt{\Bun}_G(\Sigma))\simeq QC^!_{\invertQM}(\wt{\Loc}^{?}_\Gv(\Sigma)).\end{equation}
\end{quote}

In this form of the conjecture, the various compatibilities characterizing the conjecture (e.g.\ the matching of Whittaker objects and structure sheaves,
or more generally the matching between period and $L$-sheaves) 
are well defined independently of choices of spin structure or square root of $q$ in the finite setting.
We can similarly define a twisted form of the arithmetic conjecture as formulated in Conjecture~\ref{geometric arithmetic Langlands} (see also \cite{BuzzardGee} over number fields).

\medskip

Consistently formulating the Langlands correspondence in this way would mean
that our treatment would depart from the literature, so we have not done so. 
However, it is convenient to highlight the choices involved in untwisting both sides, since they arise naturally throughout the Langlands program:

\begin{description}
\item[Automorphically] a spin structure, i.e., a choice of square root of the canonical bundle $\cK_\Sigma$ or the different, defines an isomorphism $\wt{\Bun}_G\simeq \Bun_G$. Different choices will differ by a translation on $\Bun_G$ by a two-torsion line bundle on $\Sigma$, acting via the central homomorphism $2 \check\rho(-1):\Z/2\to Z(G)$.

\item[Spectrally, in the finite context] a square root of the cyclotomic character $\varpi$, 
defines an isomorphism $\wt{\Loc}_\Gv\simeq \Loc_\Gv$.  Different choices 
will differ by a translation on $\Loc_{\check{G}}$ by a central $2$-torsion  element. 
Evaluating this square root at different points $x$ of the curve gives a choice of square root $\sqrt{q_x}$
of the size of the residue field $q_x$ at $x$, i.e., a compatible choice of square roots. In our text
we have chosen the choice $\sqrt{q_x} = \sqrt{q}^{\mathrm{deg} \ x}$ for a fixed choice of $\sqrt{q}$. 
\end{description}

%% file: fieldtheory.tex
\newcommand{\Bord}{\mathrm{Bord}}
\section{Algebraic Quantum Field Theory}\label{fieldtheorysec}

In this section we review some of the mathematical structures underlying quantum field theory and how these structures inform our view of the Langlands   correspondence and its relative version. The only technical aspects used in the bulk of the text concern $E_n$-algebras and factorization, and we direct the impatient reader to \S \ref{En algebras} and \S \ref{factorization algebras section} for a self-contained review. 

The interplay of quantum field theory and the Langlands correspondence dates back (at least) to work of Witten~\cite{WittenGrassmannians} and Beilinson-Drinfeld~\cite{BD} starting in the late 1980s, with a crucial turning point coming in the work of Kapustin and Witten~\cite{KapustinWitten}. To express these structures we formulate a notion of quantum field theory on algebraic curves valued in a symmetric monoidal higher category, which is an outgrowth of the algebraic approach to Segal's definition of conformal field theory~\cite{Segal:1987sk} (specifically a chiral CFT or modular functor) pioneered in~\cite{BeilinsonFeiginMazur} (see~\cite{bakalovkirillov,gaitsgoryCFT}), incorporating Beilinson and Drinfeld's theory of factorization algebras~\cite{BDchiral} and factorization homology, which provides a geometric counterpart to the theory of ad\`elic restricted products. An algebraic quantum field theory consists of an algebra of observables (a factorization algebra) and a module of states (a functional on its factorization homology). We drop the strong finiteness assumptions underlying modular functors (designed for rational CFTs or 3d topological field theory) and let our field theories take value in higher categories, so as to model aspects of 4-dimensional topological quantum field theories. We explain how this bare-bones definition is already sufficient to capture the field-theoretic aspects of the geometric Langlands program following the ideas of~\cite{KapustinWitten}. This formalization also allows us to define boundary theories, and we explain how the main structures appearing in this paper -- in particular Hamiltonian group actions and theta series -- fit naturally in this framework. We do not attempt to spell out all the higher categorical subtleties, but rather suggest parts of the quantum field theory intuition which can be made rigorous with current technology.

The contents of this appendix are as follows:
\begin{itemize}
\item \S \ref{TFT section} gives a motivational overview of mathematical structures underlying quantum field theory. We emphasize the roles of the theory of {\em states} (geometric quantization), as captured by the formalism of functors on bordism categories, and the theory of {\em observables} (deformation quantization), as captured by the theory of factorization algebras.
\item \S \ref{finite group TFT} illustrates some of the representation theoretic structures captured by TFT in the toy setting of finite group gauge theory.
\item \S \ref{En algebras} reviews the theories of $E_n$-algebras, and
\item \S \ref{factorization algebras section} reviews factorization algebras and factorization homology. 
\item \S \ref{algebraic field theory section} introduces a definition of algebraic quantum field theories, while
\item \S \ref{Langlands field theory section} summarizes how the Langlands correspondence fits into this formalism.
\item \S \ref{boundary AFT} discusses boundary conditions in algebraic quantum field theory, while
\item \S \ref{boundary Langlands} summarizes how relative Langlands duality can be viewed through this lens. 
 In particular we formulate our Meta-Conjecture~\ref{meta-conjecture} 
$$ \Theta_M\in \cA_G \longleftrightarrow \Ll_{\Mv}\in \cB_\Gv$$ 
which encapsulates much of the formal structure (although not the details) of this paper via boundary algebraic quantum field theories.
\end{itemize}

\begin{terminology}\label{category terminology} 
For a pointed $(\infty,n)$-category $\cC\ni 1_\cC$ we use the notation
$ \Ccirc=End(1_\cC)$ for the monoidal $(\infty,n-1)$--category of endomorphisms of the pointing (or unit), and $\Coo$ for the $(\infty,n-2)$ category of endomorphisms of the unit in $\Ccirc$.  In other words, $\Ccirc=\Omega \cC$ is the based loops in $\cC$ and $\Coo=\Omega^2\cC$ is the two-fold based loops (endomorphisms of the unit endomorphisms of the unit). 
\index{$\Ccirc$ endomorphisms of the unit}\index{$\Coo$ two-fold loops in a category}

The class of examples we have in mind are given by the notion of $\kk$-linear higher categories, constructed by iteratively  passing
to module categories: starting with
 $\Coo=Vect_\kk=\kk\module$ we take $\Ccirc=DGCat_\kk$, i.e., $Vect_\kk\module=\kk\module^2$ and $$\cC=DG2Cat_\kk:=DGCat_\kk\module=\kk\module^3$$ a category of $\kk$-linear 2-categories. (Here we are using the notation $\kk\module^n$ for the iterated higher categories of modules, taken from ~\cite{stefanichPr}.) While we only use this notion for motivation, to make precise sense of this one needs to address size issues -- in particular we think of $DGCat_\kk\in Pr^L$ as a presentable $\infty$-category but its category of modules is no longer presentable. These issues are addressed by the notion of n-presentable category developed in~\cite{stefanichPr}, where these objects are constructed not only as $(\infty,1)$- but as $(\infty,n)$-categories.  The structures we discuss do not involve non-invertible $n$-morphisms for $n>1$ so we stay in the world of $(\infty,1)$-categories, specifically enriched categories as developed in~\cite{gepnerhaugseng} (so one can for example replace $\kk\module^3$ by the category of $DGCat_\kk$-enriched categories). 

Given an object $\cM\in \cC$, we will refer to a morphism $M:1_\cC\to \cM$ from the unit of $\cC$ as an {\em object} of $\cM$, denoted $M\in \cM$. Dually we'll refer to a morphism $\Mv: \cM\to 1_\cC$ as a functional on $\cM$, which we think of as a generalized object. 
\index{$M\in \cM$ ``objects'' in an object $\cM\in \cC$}

 The ``object'' terminology is motivated by our motivating setting of higher categories of categories. For example if $\cC=DGCat_k$ then the unit $1_\cC$ is $Vect_k$ and morphisms $F:Vect_k\to \cM$ ($k$-linear colimit preserving functor) to a dg category $\cM$ are identified with objects $M=F(k)\in \cM$. Likewise if $\cC=\kk\module^3$ then ``objects'' of $\cC$ are $\kk$-linear 2-categories.
\end{terminology}

\subsection{States and Observables}\label{TFT section}\index{topological field theory}
The mathematical structure underlying quantum theories can be roughly broken up into three components: states, observables and correlation functions (the link between the first two). This trichotomy is familiar from mathematical approaches to quantum mechanics (1-dimensional quantum field theory) via quantization of a symplectic manifold $M$: one seeks to attach to $M$ a {\em geometric quantization} -- a Hilbert space $\cHH$ (states), a {\em deformation quantization} -- an associative algebra $A$ (observables), and a module structure $A\actson \cHH$, leading to a trace functional $\langle \cdots \rangle_{\cHH}:A^{\otimes n}\to \C$ (correlation functions). 
We are indebted to Pavel Safronov for teaching us this tripartite point of view on quantum theory, cf. in particular~\cite{SafronovQuantization}. We assume passing familiarity with $E_n$- and factorization algebras, which are then reviewed in the following two sections. 

\subsubsection{States} \index{states} In the setting of $n$-dimensional topological quantum field theories, the structure of states is captured by the Atiyah-Segal formalism of functors out of bordism categories, and the more general notion of {\em extended} topological field theories developed (among others) by Lawrence, Freed, Baez-Dolan, Costello, Hopkins and Lurie. (For simplicity we only consider the {\em oriented} version, so all manifolds will be oriented.) An n-dimensional TFT $\cZ$ is a representation of the symmetric monoidal $(\infty,n)$-category $\Bord_n$ of bordisms of $n$-manifolds, i.e., a symmetric monoidal functor 
$$\cZ:(\Bord_n,\coprod)\longrightarrow (\cC,\otimes).$$ Such a functor takes the empty $0$-manifold to the unit $1_\cC$, closed 1-manifolds (as self-bordisms of the empty manifold) to endomorphisms of $1_\cC$, etc., so that closed n-manifolds $M$ are taken to $n$-fold iterated based loopsin $\cC$  (endomorphisms of endomorphisms of$\dots$of$1_\cC$) . The target $\cC$ is typically taken to be $\C$-linear (a manifestation of the superposition principle), and in fact an $n$-fold  delooping of the complex numbers (in the sense that the iterated endomorphisms of the unit in $\cC$ are $\mathbb C$, cf.Terminology~\ref{category terminology}), so that $\cZ$ includes in 
particular the data of
\begin{itemize}
\item  partition functions $\cZ(N^n)\in \C$ for closed $n$-manifolds,
\item $\C$-vector spaces of states $\cZ(M^{n-1})$ for closed $n-1$-manifolds, 
\item vectors $\cZ(N)\in \cZ(\partial N)$ associated to $n$-manifolds with boundary,
\item $\C$-linear maps $\cZ(N):\cZ(\partial_{in} N)\to\cZ(\partial_{out}N)$ associated to bordisms,
\item $\C$-linear categories $\cZ(\Sigma^{n-2})$ for closed $n-2$-manifolds, 
\item $\dots$
\end{itemize}
Moreover, all of these assignments are multiplicative under disjoint union and locally constant over the classifying spaces of manifolds.

It is prohibitively hard to construct examples of this structure in its entirety. For example, $n$-dimensional TFTs coming out of gauge theory are essentially never defined on closed $n$-manifolds and at best ``top out'' at assigning vector spaces to $n-1$-manifolds. When one can assign partition functions to closed $n$-manifolds, it is typically through analytic and not algebraic or categorical means. Likewise it is rare that we know how to extend TFTs all the way down to a point (or even sometimes which higher delooping $\cC$ of $\C$ we should be working in). 

For instance (as discussed in \S\ref{TFTintro} and in more detail in\S \ref{Langlands field theory section}) we are primarily motivated by the Kapustin-Witten approach to the geometric Langlands program, which is formulated as an equivalence $$\cA_G\simeq \cB_{\Gv}$$ of four-dimensional oriented TFTs. In this case, neither side extends to 4-manifolds. The $\cB$- (or spectral) side is much better understood (in the Betti model), and in particular there are good candidates for extending the theory all the way down to a point, but on the $\cA-$ (or automorphic) side the TFT structure is much less evident and in particular we are not aware of any promising candidates for extending down to a point. 

Thus in practice it is convenient to work with a fragment of the full structure, not insisting on numbers in the top dimension or extending down to a point. Namely we might consider symmetric monoidal functors
\begin{equation}\label{TFT fragment}
\cZ:\Bord_{[k-1,k]} \to \cC
\end{equation}
 from the $(\infty,1)$-category of $k$-dimensional oriented bordisms of $k-1$-manifolds to some target symmetric monoidal $(\infty,1)$-category $\cC$. Such a functor assigns objects $\cZ(N)\in \cC$ to closed $k-1$-manifolds and $\cZ(M)\in \Ccirc$ to closed $k$-manifolds. 
 In the case $k=2$, the category  $\Bord_{[1,2]}$ can be described as having objects labeled by natural numbers, with morphisms $Hom(m,n)$ from $m$ to $n$ given by the classifying space of oriented surfaces with $m$ incoming and $n$ outgoing boundary components. 
 Moreover we can ``turn around'' bordisms to make all boundary components incoming, i.e., identify morphism spaces $Hom(m,n)\simeq Hom(m+n,0)$.

\subsubsection{Observables} \index{observables}
The theory of factorization algebras and factorization homology (which we review in \S \ref{factorization algebras section} below) arose in~\cite{BDchiral} and later~\cite{CostelloGwilliam} as a mathematical framework for the algebraic structure of observables (or equivalently symmetries) in quantum field theory (see~\cite{FBZ} for an introduction). In the Costello-Gwilliam formalism, we are given the assignment $U\mapsto \cF(U)$ of observables on open subsets of spacetime. In the Beilinson-Drinfeld formalism we are given instead an assignment $S\mapsto \cF_S$ for finite subsets $S\subset \Sigma$ of spacetime, which we think of as encoding observables on the formal neighborhood of $S$ in $\Sigma$. We then extract observables on open subsets as the output of factorization homology. 

Specifically, the theory of factorization homology arose as a formalization of the notion of conformal blocks in 2d conformal field theory, which capture the constraints (the {\em Ward identities}) satisfied by the partition function of a field theory imposed by identifying a given algebra of symmetries (formulated as a vertex, chiral or factorization algebra). In particular identifying a factorization algebra $\cF$ as observables in a given TFT $Z$ means the partition function $\cZ(M)\in \C$ on an $n$-manifold refines to a functional $$\langle - \rangle_{\cF,M}: \int_M \cF\to \C$$ on factorization homology of $\cF$ (recovering $\cZ(M)$ when applied to the unit). Analogously, the vector space of states $\cZ(N)\in Vect$ on an $n-1$-manifold is enhanced to a module for the associative algebra $\int_{N\times (0,1)} \cF$. 

Given a topological field theory $\cZ:\Bord_{[n,n-1]}\to \cC$ we can extract a factorization algebra $A=End_\cZ$ valued in $\cC$ of {\em local operators}, as the value of $\cZ$ on the sphere $S^{n-1}=\partial(D_n)$. 
By considering the bordisms obtained from embeddings of discs, the object $End_\cZ$ is endowed with the structure of oriented n-disc algebra (or framed $E_n$-algebra) in $\cC$, with unit given by a disc $\cZ(D_n)\in \cZ(S^{n-1})$, hence a locally constant factorization algebra on oriented $n$-manifolds. \index{local operators}

\begin{remark}[Theories of observables]
More generally a {\em theory of observables $A$ for $\cZ$} can be defined as a morphism of factorization algebras $A\to End_\cZ$. Thus $End_\cZ$ is the {\em final} factorization algebra of observables. This generalizes the description of a module $V$ for an associative algebra $A$ as classified by a homomorphism $A\to End(V)$. As in the case of modules for associative algebras, we think of the states $V$ and observables $A$ as independent variables linked by the morphism $A\to End(V)$ -- i.e., it is important to consider theories of observables different than the universal one $End_\cZ$ determined by states, and most examples of theories of observables in physics are not of this tautological form. 
\end{remark}

Now let us restrict attention to a fixed oriented $n$-manifold $M$. By considering discs embedded in $M$, the observables $A=End_\cZ$ define an $E_M$-algebra or factorization algebra $A_M$ on $M$, valued in $\cC$. The factorization homology 
$$\int_M A_M = \lim_\rightarrow \cZ(\partial D),$$  the colimit of the $E_M$ algebra over all disc embeddings, provides the space of global observables of the TFT on $M$. 

The relation between observables and states is captured by the data of {\em correlation functions} of local observables: by considering the complement of a union of discs in $M$ as a bordism we obtain
\begin{itemize}
\item[-] a morphism $$\langle - \rangle_{\cZ, D\to M}=\cZ(\Sigma\setminus D):\cZ(\partial D)=\bigotimes \cZ(\partial D_i)\longrightarrow 1_\cC$$ attached to oriented disc embeddings $D=\coprod D_i\hookrightarrow \Sigma$.
\end{itemize}
Thus if we label each boundary component by an object
$\cM_i\in \cZ(\partial D_i)$, we obtain an invariant
$$\langle \{\cM_i\}\rangle_{\cZ, D} \in \overline{\cC}.$$

The data of all the correlation functions on $M$ assembles into the data of a {\em state} on $A_M$: a single morphism
$$\langle - \rangle_{\cZ,M}: \int_M A_M \longrightarrow 1_\cC$$
out of the factorization homology of $A_M$ over $M$ (since a map out of this colimit amounts to a compatible collection of maps for arbitrary disc embeddings). For example the global state $\cZ(M)\in \Ccirc$ is obtained from the empty disc embedding, or equivalently by inserting the vacuum state on all discs. 

We will be interested in the case $n=2$ of oriented surfaces. In this case the data of correlation functions for varying surfaces $\Sigma$ captures a great part of the structure of the field theory. Namely the value $\cZ(N)$ on any closed oriented 1-manifold is given as a tensor product of copies of $A_\cZ$. Moreover $A_\cZ$ is automatically dualizable, in fact canonically self-dual, so we can describe the action of $Hom(n,m)$ on $A_\cZ$ in terms of the $m+n$-point correlation functions $Hom(m+n,0)\to Hom(A_{cZ}^{\otimes m+n},1_\cC)$. However, note that we have not encoded the {\em composition} of bordisms in this fashion.

\subsubsection{Defects}\label{defects section} \index{defect}
Much of the rich structure of a quantum field theory is provided by the notion of defects (and the corresponding operators or observables) of various dimensions.
As observed by Kapustin and Witten, unramified (spherical) Hecke operators in geometric Langlands arise naturally from considering the  't Hooft line defects in 4d Yang-Mills theory. These are generalizations of the Dirac monopole, and by a generalization of Gauss' law are labelled by states in a 2-sphere linking the codimension 3 line singularity.  Surface defects~\cite{GukovWitten} (codimension 2) capture ramification, local operators (codimension 4) capture the notion of singular support~\cite{ArinkinGaitsgory,ElliottYoosingular}, while this paper contends that boundary conditions or domain walls~\cite{GaiottoWittenboundary,GaiottoWittenSduality} (codimension 1) capture functoriality and periods.  (See~\cite{BettiLanglands} for an exposition of some of this structure in the setting of extended topological field theory.)

The functorial definition of TFTs naturally incorporates a notion of defect operators (or nonlocal observables) of various dimensions, most elegantly expressed via the Cobordism Hypothesis with Singularities~\cite{jacobTFT}. First, as we observed above, the local operators $\cZ(S^{n-1})$ form an $E_n$-algebra, which acts on the vector spaces of states $\cZ(N^{n-1})$ on $n-1$-manifolds. (Indeed this action descends to the factorization homology $\int_N \cZ(S^{n-1})$.)

Next come the {\em line defects}. These form the category $\cZ(S^{n-2})$, the value of the theory on the link of an embedded line. The collision (or ``operator product expansion'') of line operators endows this category with an $E_{n-1}$-monoidal structure, and the categories of states $\cZ(\Sigma)$ on closed $n-2$-manifolds form modules. Likewise, the {\em surface defects} are the $E_{n-2}$-monoidal 2-category $\cZ(S^{n-3})$ attached to the link of an embedded surface.

\index{domain wall}\index{interface}
Finally, the richest class of defects in topological field theory (and those most relevant to our work) is given by those of codimension 1, the  {\em interfaces} or {\em domain walls} between two field theories $\D:\cW\to \cZ$, which provide the natural notion of morphism of field theories. For a thorough study of domain walls and boundary theories, see~\cite{StewartThesis}. A domain wall can be defined as a functor out of the bordism category of {\em bipartite} manifolds -- manifolds with an embedded separating codimension one submanifold, and a marking of the components of the complement by the symbols $\cW$ and $\cZ$ (we'll suppress details of framing or orientation). Such a functor defines in particular two field theories $\cW,\cZ$ by considering only manifolds with the corresponding marking. Considering $\D$ on manifolds of the form $N\times I$ (separated into $\cW$ and $\cZ$ halves) gives rise to morphisms $\D(N):\cW(N)\to \cZ(N)$.

 An important special case of an interface is the notion of boundary condition (or better {\em boundary theory}) for a field theory, a morphism $\Theta:1\to \cZ$ from the trivial theory (taking any manifold to the unit) to $\cZ$. Variants of this notion have been formalized as {\em relative field theories}~\cite{FreedTeleman} and {\em twisted field theories}~\cite{StolzTeichner}, see ~\cite{JFScheimbauer} and~\cite{StewartThesis}. A boundary theory for a TFT $\cZ$ (or a {\em lax $\cZ$-twisted theory} in the terminology of~\cite{JFScheimbauer}) may be viewed a field theory of one dimension lower, valued in $\cZ$ -- indeed if $\cZ$ itself trivial then a boundary theory is simply a field theory valued in $\Ccirc$~\cite[Theorem 7.4]{JFScheimbauer}.

\begin{example}[Conformal field theories and modular functors]
The original motivation for spaces of conformal blocks (hence factorization homology) and modular functors was to express all the constraints (Ward identities) satisfied by the partition function of a 2d conformal field theory that come from knowing a chiral algebra of symmetries. This leads to the realization of conformal field theories as boundary theories for 3d topological field theories, the most famous example being the relation of the Wess-Zumino-Witten model to Chern-Simons theory. The framing anomaly of Chern-Simons theory itself is interpreted as making Chern-Simons a relative / twisted / boundary theory for a 4d invertible field theory.
\end{example}

\begin{remark}[Boundaries vs. domain walls]
Note that if we assume enough dualizability then a morphism (or interface) $\cW\to \cZ$ of field theories is identified with a boundary condition $1\to \cW^{\vee}\ot \cZ$. This reduces the notion of interfaces - but not of their composition - to that of boundary conditions. 
\end{remark}

\subsection{Extended example: finite group gauge theory.}\label{finite group TFT}
In order to understand the utility of the language of TFT for representation theory, it is invaluable to consider the most elementary examples of gauge theories, the finite-group versions of Yang-Mills theory (or Dijkgraaf-Witten theories with trivial action). See~\cite{FHLT,Teleman5,FreedHigher,freedmooreteleman} for more details. Thus we fix a finite group $G$ and discuss some aspects of $n$-dimensional TFTs $\cZ_G^n$ for $n=2,3,4$. 

These theories are all given by linearizing the same spaces of gauge fields $\Loc_G(M)$, attaching to manifolds the finite stacks (orbifolds) of $G$-local systems; since $G$ is a finite group,
these are the same as $G$-principal bundles.
In other words, these theories provide toy models for {\em both} the automorphic and the spectral theories $\cA_G$, $\cB_\Gv$ associated to reductive groups. The assignment $M\mapsto \Loc_G(M)$ defines a functor from the bordism category of manifolds to the correspondence category of orbifolds (which extends to higher categories by allowing correspondences of correspondences etc., see~\cite{haugseng,CalaqueScheimbauer}).

To linearize these spaces of fields we attach  to a finite orbifold $\cX$ 
\begin{itemize}
 \item point counts $\# \cX = \sum_{[\gamma]\in \cX} 1/|\Aut(\gamma)|$,
 \item spaces of functions $\C[\cX]\simeq \bigoplus_{[\gamma]\in \cX} (\C^{\Aut(\gamma)}\simeq \C)$, 
 \item categories of sheaves / vector bundles $\Vect[\cX]\simeq \bigoplus_{[\gamma]\in \cX} \Rep(\Aut(\gamma))$, \end{itemize} 
 with functoriality given by natural push-pull operations attached to bordisms.\footnote{We will suppress all discussion of duals and orientations, since everything we consider in the finite setting, e.g.\ functions on a finite set, is canonically self-dual.}

Thus in the 2d oriented TFT $\cZ^2_G$ we have 
\begin{itemize}
 \item $\cZ^2_G(\Sigma)=\# \Loc_G(\Sigma)$,
 \item $\cZ^2_G(S^1)=\C[\Loc_G(S^1)]=\C[G/G]$ gives class functions, 
 \item $\cZ^2_G(pt)=\Vect(\Loc_G(pt))=\Rep(G)$ the category of representations, and
 \item $\cZ^2_G([0,1])= \C[G] \in \cZ^2_G(\partial [0,1])=\Rep(G\times G)$ is the regular representation.
 \end{itemize}
 
 The $E_2$-structure of local operators for $\cZ^2_G$ is the commutative algebra structure on class functions, the center of the group algebra 
 $(\C[G],\ast)$. This is in fact a commutative Frobenius algebra, with trace given by the (outgoing) disc, i.e., (weighted) evaluation at the identity. Its spectrum is the dual $\widehat{G}$, the set of characters, equipped with (a rescaled) Plancherel measure. The linearity of $\cZ^2_G$ over its local operators (centrality of class functions) amounts to a ``Plancherel decomposition'' of $\cZ$ into a direct sum of theories over $\widehat{G}$. On the level of the numbers $\cZ^2_G(\Sigma)$ this recovers Mednykh's formula for the point-count of the $G$-character variety (see e.g.~\cite[Section 2.3]{HRV}).
 The category of line operators in the 3d theory is the braided ($E_2$) tensor category $\cZ^3_G(S^1)=\Vect(G/G)$, the Drinfeld center of $(\Vect(G),\ast)$. 
 
 Since $\Loc_G(S^2)=pt/G$, there are no nontrivial local operators in the 3d theory $\cZ^3_G(S^2)=\C$. However the 4d theory has an interesting symmetric monoidal category of line operators, $\cZ^4_G(S^2)=\Rep(G)$. Given a pointed surface $\Sigma\ni x$ we obtain an action 
 of $\Rep(G)$ on $\cZ^4_G(\Sigma)=\Vect(\Loc_G(\Sigma))$ by ``Hecke operators'': we linearize the bordism $$(\Sigma\times I )\setminus B_x:  \Sigma \coprod S^2 \to \Sigma$$
 given by removing a ball $B_x$ around the point $x\times \frac{1}{2}$. However since $G$-local systems on $S^2$ are trivial these operators become simply multiplication operators: $V\in \Rep(G)$ acts on $\cZ^4_G(\Sigma)$ by tensor product with the tautological vector bundle $V_x$ obtained by pullback along the evaluation map $ev_x:\Loc_G(\Sigma)\to pt/G$.
 
Given a finite $G$-set $X$, we obtain boundary theories $\Theta_X^n$ for each of the theories $\cZ^n_G$ by ``coupling the sigma model into $X$ to the gauge theory". First, given a manifold with boundary we consider the space of local systems with twisted maps of the boundary into $X$  
that is to say, with sections on the boundary of the associated $X$-bundle. In other words, we consider  the pullback
$$\xymatrix{ \Loc_G^X(N,\partial_N) \ar[r]\ar[d] & \Map(\partial_N, X/G)\ar[d]\\
\Loc_G(N)\ar[r]& \Loc_G(\partial_N).}$$
Linearizing these spaces defines an extension of $\cZ^n_G$ to a bordism category of manifolds with marked boundary. Concretely, for any closed manifold $M$ of dimension less than $n$,  put $N=M\times [0,1]$ and repeat the above reasoning with $\partial N$
replaced by its  component $M\times \{0\}$. 
By taking pushforward of the constant function (or vector bundle etc) along the map
$$\pi^X: \Loc_G^X(M\times [0,1], M\times \{0\}) \to \Loc_G(M).$$
we  obtain an invariant $\Theta^n_X(M)\in \cZ^n_G(M)$.

Let us call the left hand side above
 $\Loc^X_G(M) $ for short --  it is the orbifold of pairs of a local system $\rho\in \Loc_G(M)$ and a fixed point $x\in X^{\rho}$ on $X$, so pushforward along $\pi^X$ counts fixed points. For a $G$-orbit $X=G/H$ the map $\Loc^X_G \rightarrow \Loc_G$ becomes the induction map $\Loc_H(M)\to \Loc_G(M)$.
More generally, given two groups $G,H$ and a $G\times H$-space $X$ we obtain a domain wall or interface (cf. \S \ref{defects section}) $\Theta^n_X$ between $\cZ_G^n$ and $\cZ_H^n$. It is obtained by linearizing spaces of fields on bipartite manifolds, where one part carries a $G$-local system, another carries an $H$-local system and the interface carries a twisted map to $X$ (lift of the two local systems to a map to $X/G\times H$.)

Spelling out this structure leads to many familiar structures:
\begin{itemize}
\item $\Theta^2_X(pt)=\C[X]\in \Rep(G)$ is the associated representation, and
\item $\Theta^2_X(S^1)=\pi^X_*1\in \C[G/G]$ is the Atiyah-Bott formula for its character.
\item For the $G$-orbit $X=G/H$, we get the induced representation $\C[G/H]=\Ind_H^G(\C)$, and
\item the Frobenius character formula, $\chi_{\C[G/H]}=(H/H\to G/G)_* 1$.
\item The ``Neumann boundary condition'' $X=\pt$ produces the trivial representation, while 
\item the ``Dirichlet boundary condition'' $X=G$ produces the left-regular representation. 
\item For $H\subset G$ the $G\times H$-space $X=G$ defines a domain wall between $\cZ^2_H$ and $\cZ^2_G$,
which can be considered as an incarnation of either of the adjoint functors of ``induction'' or ``restriction.'' 
\end{itemize}

As a further example, if we take $G, T, U, B$ to be the points of 
a split reductive group over a finite field,  taking the action of the group $G \times T$ o $G/U$,  we obtain (finite field) parabolic induction as a domain wall. It is represented by the correspondence $$\Loc_G\leftarrow \Loc_B \rightarrow \Loc_T.$$

Given $G$-sets $X,Y$ we can consider $\cZ^2_G$ on $[0,1]$ with the two boundary components marked with $X,Y$. This produces the Hom space $\Hom_G(\C[Y],\C[X])$ as linearizing the space of fields $X\times^G Y$. For example, 
\begin{itemize}
\item Taking $Y=\pt$ produces the invariants $\C[X]^G$,
\item taking $Y=G$ produces the underlying vector space $\C[X]$, and
\item taking $Y=G/K$, $X=G/H$ produces the Mackey description of intertwiners as $\C[K \backslash G/H]$.
\end{itemize}

The sheaf $\Theta^3_X(S^1)\in \Vect(G/G)$ is the ``character sheaf" of the categorical representation 
$$\Theta^3_X(pt)=\Vect(G)\in G-cat=\cZ^3_G(pt).$$ Indeed if we formally take $X=G/B$ the flag variety for a reductive group this reproduces the Grothendieck-Springer sheaf, the pushforward of the constant sheaf on $B/B\simeq \wt{G}/G$ to $G/G$~\cite{character}.
The function $\Theta^3_X(T^2)
\in
 \C[\Loc_G(T^2)]$ is the {\em $2$-character} (or iterated trace)~\cite{KapranovGanter,secondary,campbellponto} of the categorical representation $\Theta^3_X(\pt)$. 
\footnote{More generally, a boundary theory $\Theta\in \cZ$ in an $n$-dimensional TFT determines invariants of the object $\Theta(pt)\in \cZ(pt)$, the ``higher characters'' $\Theta(N)\in \cZ(N)$ on $n-1$-dimensional manifolds, invariant under diffeomorphisms of $N$.}

Now let $\Xi$ be a $3$-manifold.
The function $\Theta^4_X(\Xi)=\pi^X_*1\in \C[\Loc_G(\Xi)]$ counting $X$-fixed points   is a toy model for both period and $L$-functions, while the vector bundle $\Theta^4_X(\Sigma)=\pi^X_*\underline{\C}\in Vect[\Loc_G(\Sigma)]$ (for a surface $\Sigma$)  is a toy model for both period and $L$-sheaves. For two $G$-sets $X,Y$ the value of $\cZ^4$ on $\Xi \times [0,1]$ with the two boundary components marked by $X$ and $Y$ agrees with the $L^2$-pairing of $\Theta^4_X(\Xi)$ and $\Theta^4_Y(\Xi)$ in $\C[Loc_G(\Xi)]$, and provides a toy model for the relative trace formula.

A boundary theory also gives richer structure when evaluated on manifolds with boundary. Namely given a TFT $\cZ$ and a boundary theory $\Theta$, we get for a manifold with boundary $(M,\partial_M)$ two objects $$\cZ(M)\in \cZ(\partial_M) \ni \Theta(\partial_M).$$ Moreover the boundary theory gives a canonical morphism between these objects: this arises from considering the manifold with corners $M\times [0,1]$ and evaluating the boundary theory on $M\times\{0\}$. Concretely in our example, this morphism is given by linearizing the correspondence
$$\xymatrix{ &\ar[dl] \Loc_G^X(M\times [0,1], M\times \{0\})\ar[dr] & \\
\Loc_G(M)\ar[dr] & &\ar[dl] \Map(\partial M, X/G) \\
& \Loc_G(\partial M). & }$$
For example:
\begin{itemize}
\item For $M=[0,1]$, 
$\Theta^2_X(\partial M)=\C[X]\otimes \C[X]\in \Rep(G \times G)$, and
\item the morphism $\Theta^2_X(\partial M)\to \cZ^2_G(M)=\C[G]$ is the matrix elements map.
\end{itemize}

 For $M=\Sigma \setminus \coprod D_i$ a surface minus discs around $x_i\in \Sigma$, we likewise get a morphism in $\bigotimes_i \cZ(S^1)$ between $\bigotimes \Theta(S^1)$ and $\cZ(\Sigma \setminus \coprod D_i)$. This is the toy model for the general construction of $\Theta$-series.

\subsection{$E_n$-algebras}\label{En algebras} \index{$E_n$ algebra}
We recall the notion of $E_n$-algebra, or algebra over the little $n$-discs operad (which we will only need for $n=1,2,3$), following the treatment in~\cite[Chapter 5]{HA}.\index{$E_n$ algebra}\index{disc algebra} Fix a symmetric monoidal $\infty$-category $\cC$. An $E_n$-algebra in $\cC$ is an object $A$ equipped with operations parameterized by the configuration space of disjoint discs in $\R^n$, $$\Conf_k(\R^n)\to Hom_{\cC}(A^{\ot k},A)$$ together with compositions corresponding to the embedding of discs inside larger discs. (Up to homotopy $\Conf_k(\R^n)$ is the configuration space of points in $\R^n$, and the compositions correspond to collisions of points.) 
In the case $n=1$, this structure is identified with that of an associative algebra object in $\cC$ (in the homotopical sense), also known as $A_\infty$-algebra. The case $n=2$ is most closely related to the geometry of configuration spaces on the affine line, or more general algebraic curves.

For $A$ an $E_n$ algebra in chain complexes, this structure amounts to maps from {\em chains} on the configuration spaces to $k$-ary operations on $A$. We may then pass to cohomology: $H^\ast(A)$ carries operations labelled by the homology of the same configuration spaces. For $n>1$ this produces a much simpler structure than the chain version: all the operations are generated by binary operations, i.e., by 
$$H_\ast(\Conf_2(\R^n))\simeq H_\ast(S^{n-1})\simeq k\oplus k[n-1],$$ and the degree $1-n$ class produces a Poisson bracket of degree $1-n$.  Equivalently, the unshearing $H^\ast(A)^\unshear$ (with respect to the cohomological-grading $\Gm$) is a Poisson algebra, equipped with a $\Gm$-action for which the bracket has weight $-2$. In our setting $A=\cO(\Mv)^\shear$ will arise as the shear of a graded cochain complex, and our aim is to produce an $E_3$-structure on $A$, whence we deduce a Poisson structure on $\cO(\Mv)$.

We also recall Lurie's form of Dunn additivity:

\medskip

\begin{theorem}~\cite{HA}\label{dunn additivity} For $k+l=n$, the structure of an $E_n$-algebra on $A\in \cC$ is equivalent to commuting structures of $E_k$ and $E_l$ algebra on $A$. 
\end{theorem}

We will need the case $1+2=3$: i.e., we will produce an $E_3$ algebra structure on $A$ out of compatible associative ($E_1$) and $E_2$ structures. 

\subsubsection{$E_n$-algebras on manifolds}
The group $O(n)$ acts on the little $n$-discs operad by changing the framing of $\R^n$, and hence on the collection of $E_n$-algebras. A {\em framed $E_n$-algebra} or {\em oriented little $n$-disc algebra} $A$ is a (homotopy) fixed point for the induced action of $SO(n)$. This equivariance for changes of coordinates means an oriented $n$-disc algebra $A$ defines a tensor functor out of the category of all oriented $n$-manifolds which are disjoint unions of discs with morphisms given by open embeddings~\cite{ayalafrancis}. We can also twist the structure of $E_n$-algebra by the tangent bundle of any oriented manifold $M$. The resulting structure, an $E_M$-algebra in the terminology of~\cite[Section 5.2]{HA}, can be thought of as a family of $E_n$-algebras $\{A_x\}_{x\in M}$ parametrized by points of $M$ but twisted by the tangent bundle of $M$ (i.e., the operations on $A_x$ are given by discs embedded in the tangent space $T_xM$), or (via~\cite[Theorem 5.2.4.9]{HA}) as a functor out of the category of disjoint unions of discs embedded in $M$. 

We will also encounter a hybrid notion between orientation and framing, namely $SO(2)$-fixed $E_3$ algebras. These objects give rise by Dunn-Lurie additivity to associative $E_\Sigma$-algebras for $\Sigma$ an oriented surface.

\subsubsection{Factorization homology}\label{En fact homology sect}\index{factorization homology} 
Factorization homology, also known as {\em topological chiral homology}, of an (oriented) $E_n$ algebra $A$ (cf.~\cite{ayalafrancis} and~\cite[Sections 5.3.2-5.3.4]{HA}) is a globalization procedure, which produces a homology theory defined on oriented $n$-manifolds. The factorization homology of $M$ with coefficients in $A$ is defined by considering $A$ as an $E_M$ algebra as above, i.e., a functor on disjoint unions of discs embedded in $M$, and then taking a colimit of this functor 
$$M\mapsto \int_M A := \lim_{\coprod D_i \to M} \bigotimes A(D_i)\in \cC.$$ Informally, we take the tensor product of copies of $A$ indexed by all discs embedded in $M$, and when a unions of $k$ discs in $M$ factors through a larger disc we factor through the corresponding $k$-ary operation $A^{\ot k}\to A$. 
A key theorem of~\cite{ayalafrancis} asserts that passing from $A$ to its factorization homology identifies oriented $n$-disc algebras (framed $E_n$ algebras) with homology theories defined only on oriented $n$-manifolds, i.e., tensor functors out of the category of manifolds under open embeddings satisfying excision.

\subsection{Factorization Algebras}\index{factorization algebra}\label{factorization algebras section}
We now give a highly impressionistic synopsis of the theory of factorization algebras. This theory is an extremely versatile generalization of $E_n$-algebras originating from the study of vertex algebras~\cite{BDchiral,FBZ}, and more generally the algebraic structure of observables in quantum field theory~\cite{CostelloGwilliam}. Just like $E_n$ algebras, factorization algebras have a local aspect, as describing multiplication operations parametrized by the collision of points, and a global aspect, as attaching measurements to open subsets of space covariantly functorial under embeddings (factorization homology).

We would like to consider factorization algebras in the setting of algebraic geometry and valued in symmetric monoidal $\infty$-categories $\cC$. This is motivated by the (local geometric) Langlands correspondence, in which case $\cC$ is a category of 2-categories, and our goal is to give a feeling for some of the key structures relevant to this setting rather than to give a detailed treatment. Indeed, such a theory is not currently available in the literature. An informal but more detailed discussion of this notion along similar lines is presented in~\cite{butson2}. 

In the Betti topological setting of manifolds, the theory of factorization algebras valued in an arbitrary symmetric monoidal $\infty$-category is developed\footnote{However, we're not aware of a reference that explicitly compares the formulation of unital factorization structures from~\cite{raskinchiral,gaitsgoryAB} with that in~\cite{HA}.} in ~\cite[Section 5.3]{HA} (see also~\cite{knudsen}), and crucially reduces to the theory of $E_n$-algebras under a local constancy hypothesis (see \S \ref{factLurie} below). 
 In the \'etale setting in positive characteristic, factorization (on the level of chain complexes) appears in the work of Gaitsgory and Lurie~\cite{GaitsgoryLurieTake1,gaitsgoryAB} on Weil's conjecture for function fields. In the de Rham setting there is ample literature starting with~\cite{BDchiral} for factorization valued in chain complexes, as well as the theory of factorization {\em categories} developed in ~\cite{raskinchiral, gaitsgoryunital} (see also \S \ref{factorization categories section}). In practice factorization categories built out of constructible sheaves only satisfy a weaker lax-monoidal version of the factorization axioms, which we suppress here (see \S \ref{lax remark}).

We will not specify precisely what properties we require for $\cC$-valued sheaf theory on schemes over $\FF$. (We also lump in without comment the topological theory of factorization on smooth manifolds from~\cite{HA}.) Such a theory includes at the minimum the data of a lax symmetric monoidal functor
$$SHV_\cC:Corr_{/\FF}\to \cC$$ from the correspondence category of schemes over $\FF$ to $\cC$,
where we denote the pullback functor for $f:X\to Y$ by $f^!:SHV_\cC(Y)\to SHV_\cC(X)$ and the pushforward by $f_*$.

\begin{remark}[Cosheaves]\label{cosheaf remark}
\index{cosheaves}
We emphasize that factorization algebras on a space $M$ are most naturally formulated as {\em co}sheaves on $M$, its powers and its Ran space of finite subsets. This is evident for example from their origin expressing observables in a field theory supported on a given patch of spacetime~\cite{CostelloGwilliam}, which are covariant under inclusion of opens. We recall that the category of $\cC$-valued cosheaves on a topological space $M$ is by definition the opposite category of $\cC^{op}$-valued sheaves,
$$cShv(M,\cC)=Shv(M,\cC^{op})^{op},$$ which amounts to covariant functors from opens on $M$ to $\cC$ taking open covers to colimit diagrams. 
 The fundamental measurement associated to a factorization algebra, its factorization homology, is naturally a {\em homology} invariant, given as global sections of a cosheaf, i.e., as a colimit, which receives maps from costalks, rather than cohomology, given as global sections of a sheaf, which maps to stalks. This cosheaf aspect is explicit in the topology literature, in particular in~\cite{HA}, and factorization homology is characterized axiomatically as a homology theory in~\cite{ayalafrancis}. However, thanks to the covariant form of Verdier duality~\cite[5.5.5]{HA} (see also the exposition in at~\cite[9.4]{GaitsgoryLurieTake1}), the theory of sheaves and cosheaves on locally compact Hausdorff  topological spaces is identified by the operation  
$$ \mbox{sheaf $\mathcal{F}$} \mapsto \mbox{cosheaf of compactly supported sections of $\mathcal{F}$}.$$
 This identification identifies the natural functoriality $(f_+,f^+)$ on cosheaves (given by the $(f^*,f_*)$ functoriality on sheaves valued in the opposite category) with the $(f_!,f^!)$ functoriality of $!$-sheaves. 
 
The cosheaf aspect of factorization algebras (or the covariant form of Verdier duality) is to our knowledge not discussed in the algebraic geometry literature. Instead, factorization algebras valued in vector spaces (or chain complexes) are formulated~\cite{gaitsgoryAB} as $!$-sheaves, and factorization homology is given as $!$-pushforward (compactly supported cohomology). One level up, factorization categories~\cite{raskinchiral,gaitsgoryunital} can be described as $!$-sheaves of categories with respect to $!$-tensor product (i.e., the value on a cover is given by a limit under $!$-pullback), but passing to left adjoints expresses them equally as cosheaves (the value on a cover is given as a colimit under $!$-pushforwards), and the global sections of a sheaf of categories is described either as a limit (cohomology) or colimit (homology). (See also \S \ref{shvcats}.)

Since we work primarily in the algebraic context we will default in \S \ref{automorphic-factorization} to the $!$-sheaf language and reserve cosheaves for the topological setting and informal discussion. 
\end{remark}

\subsubsection{The Ran Space and factorization algebras}\label{Ran space}
\index{$Ran(\Sigma)$ the Ran space}

Informallly, a factorization algebra $\cF$ over $\Sigma$ valued in $\cC$ attaches objects $\cF_S\in \cC$ to finite subsets $S\subset \Sigma$ in a fashion that varies well in families, takes disjoint unions to tensor products and is compatible with forgetting multiplicities, i.e., with diagonal maps (or collisions of points):
\begin{enumerate}
\item[$\bullet$] For every finite set $I$, we are given a $\cC$-valued $!$-sheaf $\cF_{\Sigma^I}\in SHV_\cC(\Sigma^I)$ on $\Sigma^I$. 
\item[$\bullet$] {\bf Ran's Condition:} For every surjection $\alpha:I\twoheadrightarrow J$ we have an isomorphism $$\Delta_{\alpha}^!\cF_{\Sigma^I}\simeq \cF_{\Sigma^J}.$$
\item[$\bullet$] {\bf Factorization:} For every decomposition $I\simeq I_1\coprod I_2$ we have an isomorphism
$$\cF_{\Sigma^I}|_{U_{I_1,I_2}}\simeq [\cF_{\Sigma^{I_1}} \boxtimes \cF_{\Sigma^{I_2}}]|_{U_{I_1,I_2}}$$
of the restrictions to the locus ${U_{I_1,I_2}}\subset \Sigma^I$ of disjoint $I_1$- and $I_2$-tuples.
\end{enumerate}

The definition of factorization algebra~\cite{BDchiral} is formulated in terms of the 
{\em Ran space} $Ran(\Sigma)$, the space of all finite subsets of $\Sigma$. In the algebraic setting this forms a prestack, i.e., functor from [derived] rings to [simplicial] sets, given as the colimit of $\Sigma^I$ over all diagonal maps $Ran(\Sigma)=\lim_\rightarrow \Sigma^I$. 

Given an open subset $U\subset Ran(\Sigma)$, its {\em support} is the open subset of $\Sigma$ defined as the union of finite subsets of $\Sigma$ parametrized by $U$ (see~\cite[5.5.4.3]{HA}). Open subsets $U,V\subset Ran(\Sigma)$ are said to be {\em independent} if their supports are disjoint. In this case $U\times V$ is naturally identified with an open subset $U\star V\subset Ran(\Sigma)$. 

The Ran space has a semigroup structure under union of subsets, but also a partially defined operation of {\em disjoint} union of finite sets. Disjoint union defines a correspondence $$\xymatrix{& [Ran(\Sigma)\times Ran(\Sigma)]^{disj}\ar[dl]\ar[dr]& \\
Ran(\Sigma)\times Ran(\Sigma)&& Ran(\Sigma)}$$ which makes $Ran(\Sigma)$ a (nonunital) commutative algebra in the correspondence category.
This induces a symmetric monoidal structure, {\em convolution}, on $!$-sheaves on $Ran(\Sigma)$.

A factorization algebra on $\Sigma$ is then defined as a $!$-sheaf (morally, cosheaf) $\cF$ on $Ran(\Sigma)$ which is multiplicative with respect to disjoint union~\cite{knudsen,raskinchiral}. The multiplicativity amounts to compatible isomorphisms $$\bigotimes_{I} \cF(U_i)\simeq \cF(\star_I U_i)$$ for independent subsets. Formally, $\cF$ has the structure of cocommutative coalgebra for the convolution symmetric monoidal structure, whose comultiplication map is an isomorphism on the disjoint locus. This can also be expressed in terms of the colored operad of discs in $\Sigma$~\cite[Section 5.3.4]{HA}.

\subsubsection{Units}\label{units discussion}
The notion of unital factorization algebra is developed in~\cite{raskinchiral,gaitsgoryAB} (one can also incorporate units in the operadic formulation in topology~\cite{HA}). Informally, a unital structure on a factorization algebra $\cF$ is an extension of the assignment $I\mapsto \cF_{\Sigma^I}\to \Sigma^I$, functorial for surjections of finite sets, to be functorial also over inclusions and the induced projections $\pi_{I\hookrightarrow J}: \Sigma^J\to \Sigma^I$:

\begin{enumerate}
\item[$\bullet$] {\bf Unitality:} For every injection $I\hookrightarrow J$ we have a morphism $$\Sigma^J\times_{\Sigma^I} \cF_{\Sigma^I}\longrightarrow \cF_{\Sigma^J}$$ compatible with factorization data.
\end{enumerate}

A unital structure allows one to extend the assignment $I\mapsto \cF_{\Sigma^I}\to \Sigma^I$ to the full category of possibly empty finite sets. This can be expressed elegantly as an extension of the multiplicative sheaf on the Ran space $Ran(\Sigma)$ to the {\em unital} Ran space~\cite{raskinchiral}, a unital algebra object in correspondences over $Ran(\Sigma)$.

\subsubsection{Factorization and $E_n$-algebras} \label{factLurie}
The crucial dictionary between factorization and $E_n$ structures is provided by a theorem of Lurie, identifying {\em locally constant} factorization algebras on an $n$-dimensional manifold $M$ with $E_M$ algebras. 
A factorization algebra on $M$ is said to be {\em locally constant} if $\cF$ is constructible with respect to the stratification of the Ran space of $M$ (or equivalently of the products $M^I$) by diagonals, in the sense of~\cite[Definition 5.5.11]{HA}. This entails that the $!$-restrictions to the strata are locally constant  (in the language of $!$-sheaves -- this corresponds to the $+$-restriction of cosheaves, opposite to the $*$-restriction to sheaves as in {\it op.cit.}), together with a hypercompleteness hypothesis. 

\medskip 

\begin{theorem}~\cite[Theorem 5.5.4.10]{HA}\label{factorization vs En}
There is an equivalence of $\infty$-categories between $\cC$-valued locally constant factorization algebras on $M$ and $E_M$-algebras in $\cC$.
\end{theorem}

The factorization algebra $\cF$ on $M$ attached to an $E_M$-algebra $A$ is characterized by its costalks at points of $Ran(M)$ -- i.e., finite subsets $S\subset M$ -- given by the tensor product $\bigotimes_{x\in S} A_x$ of $A$ over $S$. (In the reverse direct section we evaluate the factorization homology of $\cF$ -- see below -- over disc embeddings $U\subset M$.)

Let us spell out a special case (and its combination with Dunn-Lurie additivity Theorem~\ref{dunn additivity}):

\begin{corollary}\label{Lurie fact corollary}
\begin{enumerate}
\item There is an equivalence of $\infty$-categories between locally constant factorization algebras on $\AA^1_\C$ and $E_2$ algebras.
\item Likewise there is an equivalence of $\infty$-categories between locally constant factorization associative algebras on $\AA^1_\C$ and $E_3$ algebras.  
\end{enumerate}
\end{corollary}

Here a factorization associative algebra is a factorization algebra valued in associative algebras (or an $E_1$ object in the symmetric monoidal $\infty$-category of factorization algebras).
\index{factorization associative algebra}

\subsubsection{Factorization homology}\label{fact homology sect}\index{factorization homology}
Factorization homology was introduced by Beilinson and Drinfeld in the de Rham setting~\cite{BDchiral} as a way to capture correlation functions in conformal field theory, and in close analogy with ad\`elic constructions. It forms a refinement of a restricted tensor product $\bigotimes'_{x\in \Sigma} \cF_x$ of the values of a factorization algebra over points of a curve (pointed by the units $1_x\in \cF_x$), in which we impose local constancy in $x$ and multiplicativity (factorization or OPE) under collision of points.\index{restricted tensor product} \index{Euler product} Indeed factorization homology is a geometric counterpart of the notion of Euler product (as suggested by discussions with John Francis). This parallel is made explicit in~\cite[Sections 0.2.2, 14.1.7, 20.1.2]{gaitsgoryAB}, where the cohomological product formula for cohomology of $\Bun_G$ is interpreted as a categorified Euler product, and shown to recover the Euler product for the Tamagawa number upon taking trace of Frobenius.
 
From the factorizable cosheaf point of view, factorization homology is simply the homology (global co-sections) of $A$
$$\int_\Sigma A:= \Gamma_c(Ran(\Sigma),A), $$ the colimit of values of $A$ over opens of $Ran(\Sigma)$. 
 In particular the canonical map from any costalk $A_x$ of $A$ at $x\in \Sigma$ to the homology factors through the colimit of the diagram of tensor products $\otimes_{x\in S} A_x$ over finite subsets $S\subset\Sigma$ (defined via the unital structure of $A$), i.e., the {\em restricted tensor product} 
 $$\bigotimes_{x\in \Sigma}' A_x\longrightarrow \int_\Sigma A$$
 (where at all but finitely many points we insert the unit).
\index{restricted tensor product}

\subsubsection{Universal factorization algebras}\label{universal}
We have defined the notion of factorization algebra on a fixed smooth curve $\Sigma$ or manifold $M$, which are generalizations of the notion of $E_\Sigma$-algebra. However most factorization algebras one encounters are defined ``universally'' on arbitrary smooth curves: the notion of universal factorization algebra valued in vector spaces is studied in~\cite{BDchiral,gaitsgoryCFT,FBZ,Cliff1,Cliff2} and is equivalent to a factorization algebra on the disc, equivariant for the action of changes of coordinates, which in turn is identified with the notion of {\em quasiconformal vertex algebra}. This is in precise analogy with the passage from an $SO(n)$-fixed $E_n$-algebra to an $E_M$ algebra on any oriented manifold $M$. In physics this structure is typically expressed through the mechanism of {\em stress tensors}, providing an inner action of [the factorization algebra describing] changes of coordinates (see~\cite{CostelloGwilliam}). This richer version of universality (generalizing the notion of {\em conformal} vertex algebra, one endowed an inner action of the Virasoro vertex algebra) also encodes for example the [projectively] flat connection on factorization homology over the moduli of curves~\cite{BDchiral,FBZ}.

\index{universal factorization algebra}
\index{factorization algebra on curves}

For our purposes, we take the threadbare approach of defining a {\em factorization algebra on curves valued in $\cC$} as simply a factorization algebra defined on the universal smooth curve over the moduli stack of curves. We leave it as an open problem define universal factorization algebras in our current generality, either following~\cite{Cliff2} by attaching factorization algebras to curves functorially for arbitrary \'etale morphisms (as a suitable algebraic analog for the topological formulation as functors on manifolds with open inclusions), or using the notion of stress tensor.

\subsection{Algebraic quantum field theories}\label{algebraic field theory section}
We now imitate the structure of quantum field theory in the setting of algebraic curves, inspired by the work of Beilinson-Feigin-Mazur~\cite{BeilinsonFeiginMazur} and Beilinson-Drinfeld~\cite{BDchiral}. We begin with some motivating discussion; the eager reader may look ahead to Definition~\ref{AFT} and its explication. 

The original definition of algebraic quantum field theory~\cite{BeilinsonFeiginMazur} -- see~\cite{bakalovkirillov} for a detailed exposition -- captures algebraically the structure of $[1,2]$-dimensional part of a 3d topological field theory such as Chern-Simons theory (or equivalently the chiral part of a rational conformal field theory such as Wess-Zumino-Witten theory). Such a theory attaches a finite semisimple abelian category $\cZ_x$ to [the punctured disc around] a point of any smooth curve $\Sigma$, standing in for the circle. It also attaches a finite dimensional vector space $\langle \bigotimes \cM_i\rangle_{\cZ,S}$ to a curve marked with objects of the local categories $\cZ_{x_i}$ for a finite set $S=\{x_i\}\subset \Sigma$ (standing in for a 2-manifold with marked boundary). The category $\cZ_x$ acquires the structure of balanced braided tensor category, i.e., framed $E_2$-category, and (in modern language) gives rise to a locally constant factorization category on $\Sigma$ (an $E_\Sigma$-category in the Betti version), and the spaces of conformal blocks respect this factorization structure. One also has the crucial gluing law expressing the behavior of the spaces of conformal blocks under semistable degenerations of curves. In particular the invariant of any curve can be described in terms of $\cZ_x$ by parallel transporting to the boundary of the moduli space and using the gluing laws to reduce to genus zero, as in the proof of the Verlinde formula~\cite{faltingsverlinde}.

We are interested in algebraic models of {\em four}-dimensional topological field theories, structures that are expected to attach vector spaces to 3-manifolds, categories to 2-manifolds and a 2-category to $S^1$ (the counterpart to $\cZ_x$ above). Moreover these invariants have an inherently homotopic ($\infty$-categorical) nature, as for example derived categories of coherent sheaves on derived stacks of Langlands parameters. (This homotopical aspect is a universal feature of topological field theories arising -- as almost all examples do -- from the process of topological twisting of supersymmetric quantum field theories.)
In other words, the values on 3-manifolds are chain complexes over $\kk$ (objects of $\Coo=\kk\module=Vect_\kk$ following Terminology~\ref{category terminology}), on 2-manifolds $\kk$-linear dg categories (objects of $\Ccirc=DGCat_\kk=\kk-module^2$) and on 1-manifolds take value in an $(\infty,3)$-category $\cC=\kk\module^3$ of $\kk$-linear 2-categories.

Thankfully for us of low category number, we are only trying to model {\em fragments} of 4d TFTs (as in~\ref{TFT fragment}), which are symmetric monoidal functors 
$$\Bord_{[1,2]}\longrightarrow \cC$$
where in the source category there are only invertible bordisms between 2-manifolds. 
Hence the source here is ``only'' an $(\infty,1)$-category, i.e., has only invertible $n$-morphisms for $n>1$. 
Any such functor lands in the underlying $(\infty,1)$-category of $\cC$ (where we discard noninvertible higher morphisms), so we will continue not to need any higher category theory.  However since we do have invertible higher morphisms we may have access to the analogues of the invariants associated to certain 3-manifolds, namely mapping tori of diffeomorphisms. 

In fact we don't come close to modeling the full category $\Bord_{[1,2]}$, with objects $n\in \mathbb N$ and morphisms $Hom(m,n)$ given by classifying spaces of bordisms. Instead we only consider correlation functions, i.e., the morphism spaces $Hom(m,0)$ where all boundary components are incoming (configuration spaces of curves with $m$ marked points). As discussed in~\S \ref{TFT section} one can use duality to encode the action of all bordisms. However the {\em composition} of bordisms is largely missing from our definition: we only retain a local shadow of it in the factorization structure, expressing bordisms given by collections of discs inside larger discs.

\subsubsection{The definition}

 Fix a symmetric monoidal $\infty$-category $\cC$ and a $\cC$-valued sheaf theory as in \S \ref{factorization algebras section}.
  
\begin{definition}\label{AFT}

 An {\em algebraic quantum field theory} $\cZ$ consists of two sets of data as follows: 

\begin{enumerate}
\item[$\bullet$] [Local:] a unital factorization algebra $\cZ$ on smooth curves valued in $\cC$, and
\item[$\bullet$] [Global:] a state on $\cZ$ -- a natural transformation from the factorization homology of $\cZ$ (as a functor from curves to $\cC$) 
to the constant functor $1_\cC$, i.e., 
a functional
$$\langle - \rangle_{\cZ,\Sigma}: \int_\Sigma \cZ \longrightarrow 1_\cC$$
on the factorization homology of $\cZ$ functorial in isomorphisms of smooth curves.
\end{enumerate}

The {\em trivial field theory} valued in $\cC$ consists of the data of the unit factorization algebra $x\mapsto 1_\cC$ on every curve $\Sigma$, together with the identity map $\int_\Sigma 1_\cC\simeq 1_\cC$.

\end{definition}

The ``local'' data may be thought of alternatively as the value of the theory on 1-manifolds (objects) or as defining an algebra of observables, while the ``global'' data may thought of as the value on punctured 2-manifolds (morphisms) or as defining the states (a module over observables) and correlation functions.

\begin{remark}[Missing pieces]
We note some deficiencies in Definition~\ref{AFT}. First, recall from \S \ref{universal} that by a {\em factorization algebra on curves} we mean a factorization algebra on the universal curve. As noted in {\it loc.\ cit.}, one might instead ask for the observables to form a {\em universal} factorization algebra on smooth curves, functorial for arbitary \'etale morphisms, or better yet to define $\cZ$ as a factorization algebra on the disc endowed with a ``stress tensor'' or inner action of changes of coordinates. 

Another natural requirement is to ask for the state to factor through the {\em unital} factorization homology of $\cZ$ (cf. Remark~\ref{independent version}).

More substantially, we do not attempt to address the key structure of {\em composition of bordisms}. 
This can be expressed algebraically via the mechanism of gluing  --- the behavior of states under semistable degeneration of curves (the ``Verlinde formula''). A suitably general algebraic version of gluing inspired by~\cite{BeilinsonFeiginMazur} and especially~\cite{faltingsverlinde} is described in~\cite{murali} as a formal consequence of the data of extending the factorization algebra $\cZ$ from algebraic to rigid analytic curves.
\end{remark}

 \begin{remark}[Algebraic TFT]\label{algebraic TFT}
Definition~\ref{AFT} attempts to capture features of quantum field theories depending algebraically on an algebraic curve $\Sigma$. However, the quantum field theories relevant to the Langlands program are {\em topological}, and thus one could try to strengthen Definition~\ref{AFT} to capture some form of topological invariance. These algebraic TFTs come in three flavors: Betti, de Rham and \'etale. Informally, we ask for the theory $\Sigma \mapsto (\cZ_\Sigma,\langle - \rangle_{\cZ,\Sigma})$ to factor through the assignment $\Sigma\mapsto \Sigma_?$ ($?=B,dR,et$) of the Betti, de Rham or \'etale spaces associated to smooth curves. 

In a Betti TFT we ask for the factorization algebra $\cZ$ on curves and its state $\langle - \rangle_{\cZ,\Sigma}$ to descend to the moduli of the Betti spaces of curves, i.e., the classifying spaces of diffeomorphism groups or homotopy type of the moduli of curves. Thus a Betti TFT amounts to (a small amount of the structure of) a topological field theory in the traditional sense, a functor of (oriented) topological manifolds, with invariants $\cZ(\Sigma)$ forming local systems over the moduli of curves and carry actions of mapping class groups. 

A de Rham field theory is a higher analog of the structure of {\em topological conformal field theory} \cite{getzlerTCFT},\cite{SegalTFT}, i.e., a CFT with a homotopic trivialization of the stress tensor. Indeed this is a general feature of ``topological twists'' of supersymmetric quantum field theories, explored in~\cite{ElliottSafronov}, in which the action of changes of coordinates (and deformations of metrics) is made exact in a structured way (through the action of the supersymmetry algebra). This gives a weak form of topological invariance -- in particular the categories $\cZ(\Sigma)$ attached to curves carry a flat connection over the moduli of curves which is {\em not} integrable in general (i.e., parallel transport is not defined) and do not inherit actions of mapping class groups. (This is a major difference with the classical setting of rational CFT or 3d TFT as in~\cite{BeilinsonFeiginMazur,bakalovkirillov}.)

Likewise in the \'etale version for $\FF={\overline \FF}_q$, we might ask an \'etale field theory to be 
functorial for isomorphisms of the \'etale site of $\Sigma$. Let us point out one important feature of this. For a curve $\Sigma$ defined over $\FF_q$, the \'etale site of $\Sigma/{\overline \FF}_q$ carries a canonical Frobenius automorphism. Hence the invariant $\cZ(\Sigma)$ carries a Frobenius automorphism (the action of an ``arithmetic mapping class group element''), and if it is dualizable, we may evaluate its Frobenius trace and consider it as the invariant attached to the corresponding ``arithmetic 3-manifold'', the mapping torus of Frobenius. 
\end{remark}

\subsubsection{Global states}

We now spell out the key data encoded in Definition~\ref{AFT}, and in particular the motivating case where objects of $\cC=\kk\module^3$ are dg 2-categories, of $\Ccirc=\kk\module^2$ are dg categories and of $\Coo=\kk\module$ are vector spaces (following our Terminology~\ref{category terminology}).

The most basic data associated to an algebraic quantum field theory is the ``space of states'':
 
$\bullet$ the states on $\Sigma$ define an invariant $\cZ(\Sigma)=\langle 1_{\cZ_\Sigma} \rangle_{\cZ,\Sigma}\in \Ccirc$, i.e., a dg category. 

Explicitly, the unital structure of $\cZ$ in particular endows the factorization homology on any curve $\int_\Sigma \cZ$ with a pointing $1_{\cZ,\Sigma}$ (morphism from $1_\cC$), and $\cZ(\Sigma)$ is the value of the state $\langle - \rangle_{\cZ,\Sigma}$ on the pointing.

If we assume $\cZ(\Sigma)$ is dualizable then we further get that

$\bullet$ any automorphism $F\in Aut(\Sigma)$ defines an invariant $Tr(F, \cZ(\Sigma))\in \Coo$, i.e., a dg vector space, which we think of as the states $\cZ(M_{F,\Sigma})$ of the theory on the mapping torus of $F$ (a 3-manifold).

\subsubsection{Factorization algebra of observables / local data}
The local data of $\cZ$ on a fixed curve $\Sigma$ includes:

$\bullet$ an invariant $\cZ_x=\cZ(D^*_x)\in \cC$ (i.e., a 2-category) to every point $x\in \Sigma$, which we think of as the states of the field theory on the ``1-manifold'' $D^*_x$, the punctured disc at $x$. 

$\bullet$ The invariant extends multiplicatively to finite subsets $S\subset \Sigma$: $$\cZ_S\simeq \bigotimes_{x\in S} \cZ_x.$$
 
$\bullet$ We are given a unit: an object $1_{\cZ_x}=\cZ(D_x)\in \cZ_x$ ($x\in \Sigma$), the ``vacuum state'' at $x$, and

$\bullet$ we insert units to extend $\cZ_S$ to be functorial under all maps of finite sets over $\Sigma$.

The structure of factorization algebra glues together the objects $\cZ_x$ for varying $x$, and most importantly encodes the ``operator product expansion'' (OPE), i.e., behavior when points collide: namely the invariants $\cZ_S$ assemble to a factorizable cosheaf over the [unital version of the] Ran space $Ran(\Sigma)$ of all finite subsets of $\Sigma$. 

\subsubsection{States / correlation functions}
The global aspect of $\cZ$ is given by spaces of correlation functions:

$\bullet$ for every finite subset $S\subset \Sigma$ we are given a functional
$$\langle - \rangle_{\cZ,S}: \cZ_S \to 1_{\cC}.$$

Concretely, applying this morphism to objects $\cM_i\in \cZ_{x_i}$ (i.e., composing with morphisms $\cM_i:1_\cC\to \cZ_{x_i}$)  we obtain an invariant 
$$\langle \bigotimes \cM_i\rangle_{\cZ,S}\in \overline{\cC}.$$ These assignments are asked to be invariant under maps of finite sets, in particular insertion of the unit at new points. Thus if we insert the unit everywhere we obtain the global states $\cZ(\Sigma)$. 

We require that all of the assignments $\langle - \rangle_{\cZ,S}$ to assemble to a single state on $\cZ$, i.e., a morphism out of the global observables $\int_\Sigma \cZ$, i.e., the homology of the Ran space with coefficients in the cosheaf defined by $\cZ$.
This guarantees that the assignment of correlation functions respects the factorization structure on the $\cZ_S$, in particular varies well with $x$ and is compatible with collisions of points.

\subsubsection{Defects in algebraic quantum field theories}
Algebraic field theories following Definition~\ref{AFT} afford defect operators of various dimensions. We continue to think of an algebraic quantum field theory as the $[1,2]$-dimensional part of a 4d TFT, and our terminology reflects this choice: as common in physics we label defects by their dimension in space-time (although emphasizing {\em co}dimension is more natural when thinking of an abstract $[1,2]$-dimensional field theory). 

$\bullet$ Surface defects: the assignment $x\mapsto \cZ(D_x^*)$ itself plays the role of the 2-category of surface defects in a 4d TFT, while its factorization structure captures the operator product expansion ($E_2$-structure) on surface defects. They play the role of possible codimension 2 singularities (ramification data) on $\Sigma$. 

$\bullet$ Line defects: the ``unramified Hecke category'' $\HECKE_{\cZ}$ of $\cZ$ plays the role of the category $\cZ(S^2)$ of line defects. These are by definition endomorphisms of the unit observable $$1_{\cZ_x}=\cZ(D_x) \in \cZ_x =\cZ(D^*_x),$$ i.e., $$\HECKE_{\cZ}:=End(1_{\cZ_x})\in \Ccirc.$$ (The algebraic avatar of the 2-sphere here being two discs joined along the punctured disc.)
The factorization algebra structure on $\cZ_x$, with $1_{\cZ_x}$ as its unit, provides the Hecke category with the structure of an factorization associative algebra on $\Sigma$ valued in $\Ccirc$ -- in the Betti setting this amounts to the structure of $E_3$-algebra in $\Ccirc$, as expected from line defects in 4d TFT.

$\bullet$ Observables from line defects: for any choice of point $x\in \Sigma$ we have an action of $\HECKE_\cZ$ on $\cZ(\Sigma)$ by Hecke modifications. This is a consequence of the unitality: we identify $\cZ(\Sigma)$ with correlation functions with the unit inserted at $x$, whence an action of endomorphisms of the unit. More generally $\langle \bigotimes \cM_i\rangle_{\cZ,S}$ carries an action of $\HECKE_\cZ$ at all points away from the ramification set $S$. This action descends to the global Hecke category, given by the factorization homology 
\begin{equation}\label{unramified Hecke action}\bH_{\cZ,\Sigma\setminus\{x_i\}}:=\int_{\Sigma\setminus \{x_i\}} \HECKE_{\cZ} \hskip.1in \actson \hskip.1in \langle \bigotimes \cM_i\rangle_{\cZ,S}
\end{equation}
In particular for $S$ empty we obtain an action of $\bH_{\cZ,\Sigma}$ on $\cZ(\Sigma)$. 

$\bullet$ Local operators: the point defects in a 4d TFT are captured in the algebraic quantum field theory setting by endomorphisms of the unit line defect, i.e., endomorphisms of the unit in $\HECKE_{\cZ}$. These endomorphisms form an $E_2$ factorization algebra on $\Sigma$, i.e., the algebraic counterpart of the $E_4$ structure on local operators in four dimensions. 

$\bullet$ Domain walls and boundary theories: we postpone to Section~\ref{boundary AFT} the discussion of the algebraic avatars of the richest class of defects in 4d TFT.

\subsection{Langlands correspondence via field theory}\label{Langlands field theory section}
We now outline how the structure of algebraic quantum field theory is meant to apply in the setting of the Langlands correspondence, in parallel to the discussion of the previous section (see~\cite{BettiLanglands} for a more limited discussion along similar lines). Namely, we explain how the geometric Langlands correspondence on a smooth projective curve $\Sigma$ over an algebraically closed field might be understood as an equivalence $\cA_G\simeq \cB_{\Gv}$ between algebraic quantum field theories, valued in $\kk$-linear 2-categories $\cC=\kk\module^3$. Moreover in the \'etale setting when $\Sigma$ is defined over a finite field, we discuss how passing to Frobenius traces  following~\cite{AGKRRV3} -- and assuming that we have available the structure described in
Remark  \ref{algebraic TFT} -- would recover a form of the Langlands conjecture for function fields. 

 \begin{remark}[Categorical difficulties]
The full structure of $\cA_G$ and $\cB_{\Gv}$ as algebraic quantum field theories is far from being precisely formulated. For example we can ask: what is the precise target category $\cC$ and the associated sheaf theory $SHV_\cC$? This is part of the problem of formulating the local geometric Langlands conjecture (being developed by D. Arinkin, D. Gaitsgory, S. Raskin and others): what kind of objects are the  ``local 2-categories $\cA_G(D^*_x)$ and $\cB_{\Gv}(D^*_x)$''?
Just as ind-coherent sheaves are not sheaves in a naive sense (projecting to $QC$ kills some objects) it is not clear we should consider them as actual 2-categories (i.e., objects of $\kk\module^3$) or as more sophisticated objects constructed by a higher sheaf theory $2IndCoh$ and valued in the 3-category $3IndCoh(pt)$ (these objects are defined in forthcoming work of Stefanich~\cite{stefanichIndCoh}).  However their unramified parts form a ``sub-field theory'' that is fairly well understood and already contains a great amount of structure (in particular all that is directly relevant to this paper), and we hope our schematic overview of the general expectations is useful regardless. 
\end{remark}

\begin{remark}[Topological invariance]\label{topological Langlands}
The field theories $\cB_\Gv$ on the spectral side (and hence, by the Langlands conjectures, the automorphic theories $\cA_G$) in the Betti, de Rham and \'etale settings satisfy a strong ``topological invariance'' property: they are
built from the
corresponding stacks $\Loc^?_\Gv(\Sigma)$ of ($?=$B, dR or et) local systems, which are themselves functors of $\Sigma_?$. In other words, they form algebraic TFTs as in Remark~\ref{algebraic TFT}. In the Betti setting, this means we are in the traditional setting of topological field theory (which was precisely the motivation for the Betti Langlands conjecture~\cite{BettiLanglands}).
In the \'etale setting, we only note again that this form of topological invariance implies that the value of the theory on curves defined over finite fields inherit Frobenius automorphisms. 
\end{remark}
 
\medskip

$\bullet$ [Global geometric] The global unramified automorphic invariants $\cA_G(\Sigma)=\AUT^?(\Bun_G(\Sigma))$ are given by the [\'etale, de Rham or Betti form of the] dg category of automorphic sheaves on $\Bun_G(\Sigma)$ (as in \S \ref{geometric Langlands}). 
On the spectral side, $\cB_{\Gv}(\Sigma)=\QC^!(\Loc_\Gv^?(\Sigma))$ is given by ind-coherent sheaves  on the [\'etale, de Rham or Betti] stack of local systems on $\Sigma$. Thus the conjectured equivalence $\cA_G\simeq \B_{\Gv}$ on $\Sigma$ thus realizes the geometric Langlands conjecture.

\medskip

$\bullet$ [Global arithmetic] If we assume a form of ``\'etale topological invariance'' as in Remarks~\ref{algebraic TFT} and~\ref{topological Langlands}, we can see some of the structure of the arithmetic Langlands correspondence for function fields. The categories $\cA_G(\Sigma), \cB_\Gv(\Sigma)$ are dualizable, hence one can evaluate the trace of any automorphism of $\Sigma$.  For $\Sigma$ defined over $\FF_q$ we may consider the categorical trace of Frobenius, which defines the values of the field theories $\cA_G$ and $\cB_\Gv$ (in the \'etale setting of~\cite{AGKRRV1}) on the ``arithmetic 3-manifold'' $M_{F,\Sigma}$ corresponding to $\Sigma$, the mapping torus of Frobenius. 

The trace $Tr(F,\B_\Gv(\Sigma))$ is identified with volume forms on the stack $\Loc_\Gv^{arith}(\Sigma)$ of {\em arithmetic} restricted local systems (Frobenius fixed points on $Loc_\Gv$) (see~\cite[Section 24]{AGKRRV1}). On the other hand, the Trace Theorem~\cite{AGKRRV3} recovers the space of unramified automorphic forms -- compactly supported functions on $Bun_G(\Sigma)(\mathbb F_q)$ as the categorical trace of Frobenius on $\AUT^{res}(\Bun_G(\Sigma))$.
Thus the conjectural isomorphism of vector spaces $\cA_G(M_{F,\Sigma})\simeq \cB_{\Gv}(M_{F,\Sigma})$ recovers the unramified Langlands conjecture for function fields, as formulated in~\cite{AGKRRV1} following the work of V. Lafforgue~\cite{LafforgueICM}.

\medskip

 $\bullet$  [Local geometric]
 The local automorphic 2-category $\cA_G(D^*)$ is expected to be given by categorical representations of the algebraic loop group $G_F$, i.e., by a suitable version of module categories for the convolution monoidal category $SHV(G_F)$.
 
 The local spectral 2-category $\cB_{\Gv}(D^*)$ in the de Rham setting is expected to be Stefanich's $2IndCoh(\Loc_\Gv(D^*))$, a modification (in the spirit of the modification $QC\leadsto QC^!$)  of the 2-category of quasicoherent sheaves of categories on the stack of local systems on $D^*$. 
The conjectured equivalence $\cA_G\simeq \B_{\Gv}$ on $D^*$ is the local geometric Langlands conjecture.

\medskip

$\bullet$ [Local unramified] While these 2-categories are poorly understood, their {\em unramified} parts are quite familiar. The unit object (or vacuum)  $$\cA_G(D)\in \cA_G(D^*)$$  is given by the $G_F$-category of sheaves on the affine Grassmannian $Gr=G_F/G_O$. 
Its endomorphisms, the Hecke category of ('t Hooft) line operators $\HECKE_{\cA_G}$, form the spherical Hecke category $SHV(G_O\backslash G_F/G_O)$. The (factorization) 2-category of modules for $\HECKE_{\cA_G}$ defines the well-understood unramified part of the local 2-category. 

On the spectral side $\cB_{\Gv}(D)$ is given by the category $Rep(\Gv)\simeq QC(pt/\Gv)$ of representations of the dual group thought of as sheaves on the substack $pt/\Gv\to Loc_{\Gv}(D^*)$ of trivial local systems on the punctured disc. 
Its endomorphisms, the (Wilson) line operators $\HECKE_{\cB_{\Gv}}$ recover the derived spherical category $QC^!(Loc_{\Gv}(S^2))$, and the Langlands duality of field theories on the (factorization) monoidal categories of line operators $\HECKE_{\cA_G}\simeq \HECKE_{\cB_\Gv}$ predicts the derived Satake correspondence~\cite{BezFink} (and its factorizable form~\cite{CampbellRaskin}). 

We may also pass to the Frobenius trace on the geometric Satake correspondence and recover the classical Satake correspondence (as in~\cite{xinwentrace}). Namely, the Frobenius trace on $\HECKE_{\cA_G}=SHV(G_O\backslash G_F/G_O)$ produces the spherical Hecke algebra, endomorphisms of the unramified representation $k[Gr]$ while the Frobenius trace on $\HECKE_{\cB_\Gv}$ produces functions on Frobenius-twisted conjugacy classes in $\Gv$.

\medskip

$\bullet$ [Unramified Hecke action and shtukas] The Hecke categories $\HECKE_{\cA_G}\simeq \HECKE_{\cB_\Gv}$ act on the global states $\cA_G(\Sigma)\simeq \cB_\Gv(\Sigma)$ by modifications at any point $x\in \Sigma$(~\ref{unramified Hecke action}). Thus the duality of field theories predicts the Hecke-linearity of the geometric Langlands correspondence. Moreover this equivalence respects the factorization monoidal structure of the Hecke categories, so that we may integrate over $x\in \Sigma$ to obtain an equivalence as modules for the unramified global Hecke categories $\bH_\Sigma=\int_{\Sigma} \HECKE$. 
In the Betti setting, the global Hecke category $\HECKE_{\Sigma,\cB_{\Gv}}$ is Beraldo's global Hecke category $\bH_\Sigma$ with its action on $\cB_{\Gv}(\Sigma)=QC^!(Loc_{\Gv}(\Sigma))$. 
 In general the global Hecke action detects singular support of coherent sheaves on $\Loc_\Gv(\Sigma)$ (the ``sheafification'' of the automorphic category over $T^*[1]\Loc_\Gv(\Sigma)$).  

In the global arithmetic setting, we obtain likewise the action of the spherical Hecke algebra on unramified automorphic forms. 
More significantly, for any unramified Hecke functor $H_S$ acting at $S\subset \Sigma$ we can take the categorical trace of Frobenius composed with $H$, 
$$Tr(F\circ H, \cA_G(\Sigma)\simeq \cB_{\Gv}(\Sigma))$$ -- in other words consider the value of the field theory on the arithmetic 3-manifold $M_{F,\Sigma}$ with insertions of line defects along $S$. As explained in~\cite{AGKRRV1,AGKRRV3} this categorical trace recovers the Langlands conjecture for the cohomology of moduli of unramified shtukas on $\Sigma$ with prescribed legs at $S$.

\medskip

$\bullet$ [Local operators and singular support]  
The local operators, i.e., endomorphisms of the unit in $\HECKE_{\cZ}$, are identified on the automorphic side with the equivariant cohomology ring $H^*(BG)$ and spectrally with the shifted invariant polynomials $\cO(\fgxv[2])^\Gv$. The vanishing of these operators measures the difference between ind-finite and ind-safe (renormalized and unrenormalized) categories of sheaves on $\Bun_G$ on the automorphic side (cf. \S \ref{renormalization section}), and the difference between ind-coherent sheaves with or without nilpotent singular support on the spectral side \S \ref{spectral side section}.

\medskip

$\bullet$ [Ramified global]
The local Langlands 2-categories of surface defects $\cA_G(D^*)$ and $\cB_{\Gv}(D^*)$ control ramification: one may view the action of a surface defect at $x\in \Sigma$ as modifying the ramification data considered at $x$. 

Namely, given objects $\cM_i\in \cA_G(D^*)$ we can consider the {\em correlation functions}, the automorphic categories  $\langle \cM_i \rangle_{\cA_G,\Sigma}\in DGCat_\kk$. These give the categories of sheaves on the stacks of $G$-bundles on $\Sigma$ with given ramification data at points $x_i\in\Sigma$ -- for example taking $\cM_i=SHV(G_F/H)$ for a congruence subgroup $H$ corresponds to imposing $H$-level structure at $x$. 
On the other hand inserting the unit surface defect, i.e., the unramified $G$-category, 
captures the imposition of no ramification at $x$. 

 The colimit of all these assignments, the factorization homology $\int_\Sigma \cA_G$ (i.e., the global observables on $\Sigma$), is a version of a ``restricted tensor product'' of the local automorphic 2-categories
$$\int_\Sigma \cA_G\sim \bigotimes_{x\in \Sigma}' G_F\mbox{-cat},$$ and the ``correlation function functional''  $$\langle - \rangle_{\cA_G,\Sigma}:\int_\Sigma \cA_G\to \Ccirc,$$ plays the role of the functor of $G_F$-coinvariants, or its representing object, the geometric counterpart of the space of all ad\`elic automorphic forms. 

Spectrally, for local ramification data $\check{\cM}_i\in \cB_\Gv(D^*)$, the correlation functions  $\langle \cM_i \rangle_{\cB_\Gv,\Sigma}$ give the dg categories of ind-coherent sheaves on moduli of local systems with prescribed ramification. The duality $\cA_G\simeq \cB_\Gv$ on correlation functions gives the ramified geometric Langlands conjecture.

Thus the equivalence of algebraic quantum field theory encodes a very general and flexible form of the ramified Langlands conjecture in the setting of curves. 
On the one hand, we recover the action of unramified Hecke modifications as the action~\ref{unramified Hecke action} of  $\HECKE_{\cA_G,\Sigma\setminus S}=\int_{\Sigma \setminus S} \HECKE_{\cA_G}$ on spaces of correlation functions $\langle \cM_i \rangle_{\cA_G,\Sigma}$. On the other hand, 
passing to Frobenius traces is then expected to give the ramified Langlands conjecture for function fields.

\medskip

$\bullet$ [Local arithmetic] Finally, we note that the Frobenius trace on the local geometric Langlands conjecture is expected to realize the categorical local Langlands conjecture in the spirit of~\cite{xinwentrace, farguesscholze} (see also~\cite{BZCHN} for a discussion).
Namely, the trace of Frobenius on $\cB_\Gv(D^*)=2IndCoh(\Loc_\Gv(D^*))$ formally produces the category $QC^!(Loc_\Gv^{arith}(D^*)$ of ind-coherent sheaves on the stack of local arithmetic Langlands parameters.
Much more speculatively, as suggested in~\cite{dennisshtuka}, the trace of Frobenius on the 2-category $\cA_G(D^*)$  is expected to produce the category of sheaves on the Kottwitz space of $G$-isocrystals $G_F/_{\sigma} G_F$. In particular it should contain as a full subcategory the local Langlands category of smooth representations of $G_F$.

Thus the duality $\cA_G\simeq \cB_{\Gv}$ on the ``arithmetic 2-manifold'' $F$ is meant to produce a full embedding
$$Rep^{sm}(G_F)\hookrightarrow QC^!(Loc^{arith}_\Gv(D^*))$$ of smooth representations of $G_F$ into ind-coherent sheaves on the stack of local Langlands parameters.

\subsection{Boundaries in algebraic quantum field theory}\label{boundary AFT}

We now discuss the formulation of interfaces (or domain walls) and boundary theories in the setting of algebraic quantum field theory.

\begin{definition} 

Let $(\cW, \langle - \rangle_{\cW})$ and $(\cZ, \langle - \rangle_{\cZ})$ be algebraic quantum field theories on curves valued in $\cC$, as in Definition~\ref{AFT}.

\begin{enumerate}
\item A morphism or {\em interface} $\D:\cW\to \cZ$ between field theories is a lax morphism of unital factorization algebras on curves $\cW\to \cZ$ together with a specified natural transformation $$\langle - \rangle_{\D,\Sigma}: \langle - \rangle_{\cZ,\Sigma} \circ {\int_\Sigma \D}\longrightarrow\langle - \rangle_{\cW,\Sigma} $$ making the following diagram commute:
$$\xymatrix{\int_\Sigma \cW \ar[dr]_-{\langle - \rangle_{\cW,\Sigma}}\ar[rr]^-{\int_\Sigma \D}& & \int _\Sigma \cZ \ar[dl]^-{\langle - \rangle_{\cZ,\Sigma}} \\
&1_\cC&}$$
\item A {\em boundary theory} $\Theta\in \cZ$ for the field theory $\cZ$ is a morphism from the trivial theory $$\Theta: 1_\cC \to \cZ.$$ 
\end{enumerate}
\end{definition}

\begin{remark}[Lax structures] \label{lax remark}
The structure of boundary theory we describe here is [an algebraic version of] a fragment of the notion of lax boundary TFT from~\cite{JFScheimbauer}. The laxness is evident in two places: first, we ask for a morphism of field theories only to commute with correlation functions up to specified natural transformation rather than natural isomorphism.

Second, we ask for the map of factorization algebras itself to be lax. It's useful to recall one role of lax monoidal functors. If we think of an object $M\in C$ in a $k$-linear monoidal category $C$ as the image of the unit $k\in Vect_k$ under a functor from the unit category $\underline{M}:Vect_k\to C$ with $M=\underline{M}(k)$, then the structure of associative algebra object on $M$ (in particular the map $M\otimes M\to M$) corresponds to a lax monoidal structure on the functor $\underline{M}$ (so that we have a morphism $M(k)\otimes M(k)\to M(k\otimes k)=M(k)$ which is not required to be an isomorphism). 
\end{remark}

Let us spell out the notion of boundary field theory $\Theta\in \cZ$. It is very helpful to refer back to the case when $\cZ$ is trivial:

\begin{lemma} A boundary theory $\D:1_\cC\to 1_\cC$ for the trivial theory is equivalent to the data of a field theory $\D$ valued in $\Ccirc$. 
\end{lemma}

In general we find {\em relative} versions of all the structures of a field theory, valued in the ``bulk theory'' $\cZ$ rather than in $\Ccirc$. 
Alternatively, it is useful to refer to the next section where we identify these structures with (relatively) familiar objects in the relative Langlands program.

\subsubsection{Global state}
The unital structure of $\Theta:1_\cC\to \cZ$ applied to correlation functions produces 

$\bullet$ the states on $\Sigma$ define a ``$\Theta$-series'' functional $\Theta(\Sigma): \cZ(\Sigma)\to 1_\Ccirc$,
i.e., a [representing] object in the dg category $\cZ(\Sigma)$.

Given sufficient dualizability we can then pass to trace of any automorphism $F$ of $\Sigma$ to obtain
$$\Theta(M_{F,\Sigma}):\cZ(M_{F,\Sigma})\to 1_\Coo,$$ i.e., a linear functional on the vector space $cZ(M_{F,\Sigma})$.

\subsubsection{Local observables}
Considering the morphism $\Theta$ at points $x\in \Sigma$ produces the following algebraic structures:

$\bullet$ a unital {\em factorization algebra object} $\{\Theta_x\in  \cZ_x\}_{x\in \Sigma}$ of $\cZ$ -- i.e., a section of the factorization algebra compatible with factorization~\cite{raskinchiral}, and equipped with a unit morphism $\Phi:\cZ(D)\to \Theta$ from the unit of $\cZ$. 

$\bullet$ a unital {\em factorization associative algebra object} $\HECKE_\Theta\in \HECKE_\cZ$ in the Hecke category of $\cZ$. 

(The latter, viewed as a lax morphism $\HECKE_\Theta:1_\Ccirc\to \HECKE \cZ$ of unital factorization associative algebras, comes by applying the lax morphism of unital factorization algebras $\Theta:1_\cC\to \cZ$ to the endomorphisms of the unit.)

How should we understand these structures? Recall that from the point of view of 4d TFT, the [factorization or $E_2$] 2-category $\cZ$ encodes surface -- i.e., codimension 2 -- defects, while its [factorization associative or $E_3$] Hecke category $\HECKE_\cZ$ encodes line -- i.e., codimension 3 defects. 
The relative versions $\Theta$ and $\HECKE_\Theta$ describe codimension 2 (line) and 3 (point) defects in the boundary theory (and reduce to line and point defects, respectively in the $\Ccirc$-valued theory or ``3d TFT'' $\Theta$ in the case when $\cZ$ is trivial). In this language, $\Theta$ is the value of the 4d TFT on a cylinder $S^1 \times I$ with one end marked by the boundary condition $\Theta$, while $\HECKE_\Theta$ is the value on $S^2\times I$ with a similiarly marked boundary.

\subsubsection{Unramified local observables}
We may also consider the unramified part of $\Theta$, $$\Theta^{unr}:=Hom_{\cZ}(\cZ(D), \Theta).$$
This forms a unital factorization algebra in $\Ccirc$. Keeping track of linearity over endomorphisms of the vacuum $\cZ(D)$ we can further consider $\Theta^{unr}$ as a unital factiorization algebra in modules for the line operators (Hecke category) $\HECKE_\cZ$. Topologically this is the value of the theory on a 2-disc ending on the boundary, i.e., the link of a line in the boundary. Physically these operators form the codimension 2 defects in the boundary theory. The unit of $\Theta^{unr}$ is provided by the basic object $\Phi$, i.e., the unit of $\Theta$.
(In the presence of sufficient dualizability assumptions to recover the ULA and rigidity conditions as in \S \ref{automorphic-factorization} we recover
$$\HECKE_\Theta\simeq \underline{End}_{\HECKE_\cZ}(\Phi)$$ as its internal endomorphisms.)

\subsubsection{Boundary correlation functions}

We now consider the global data produced by a boundary condition. We can evaluate factorization homology of the factorization algebra object $\Theta\in \cZ$ to obtain an object $\int_\Sigma \Theta \in \int_\Sigma \cZ$, and then apply correlation functions $\langle - \rangle_\cZ$:

$\bullet$ a functional on correlation functions with values in $\Theta$, i.e., a morphism 
$$\langle - \rangle_{\Theta,\Sigma}:  \langle \int_\Sigma \Theta \rangle_{\cZ,\Sigma} \to 1_{\Ccirc}.$$ This amounts to a compatible family of morphisms
$$\langle \bigotimes_S \Theta_{x_i} \rangle\longrightarrow 1_{\Ccirc},$$
which for $S$ empty reproduces the ``$\Theta$-series'' functional $$\Theta(\Sigma): \cZ(\Sigma)\to 1_{\Ccirc}.$$

Again if $\cZ$ is trivial this amounts to the data of correlation functions in the $\Ccirc$-valued theory $\Theta$. 

Finally we observe that boundary observables on $\Sigma\setminus S$ act on boundary correlation functions, in the boundary version of~\ref{unramified Hecke action} (for simplicity we formulate this only in the everywhere unramified case):

$\bullet$ The boundary observable algebra $\bH_{\Theta,\Sigma}:=\int_{\Sigma} \HECKE_{\Theta}$, an algebra object in $\bH_{\cZ,\Sigma}$, acts on the functional $\Theta(\Sigma)$ on the $\bH_{\cZ,\Sigma}$-modules $\cZ(\Sigma)$.

\subsection{Relative Langlands duality via field theory}\label{boundary Langlands}
We now explain how to match the structures underlying relative Langlands duality with those of boundary conditions in field theory.
To a hyperspherical $G$-variety $M$ we would like to attach its {\em automorphic quantization}, a boundary theory  $\Theta_M\in \cA_G$ for the automorphic field theory. Likewise to a hyperspherical $\Gv$-variety $\Mv$ we would like to attach its {\em spectral quantization}, a boundary theory 
$\Ll_\Mv\in \cB_\Gv$ for the spectral field theory (see Definition~\ref{sdq} for a weaker notion of spectral deformation quantization). We recount below some of the expected values of these boundary theories. The relative Langlands duality studied in this paper is organized by the following principle:

\medskip

\begin{conjecture}[Meta-conjecture]\label{meta-conjecture}
The automorphic and spectral quantizations of dual hyperspherical varieties are identified under the conjectural Langlands correspondence of algebraic quantum field theories $\bL:\cA_G\simeq \cB_\Gv$: i.e., there is a commutative diagram
$$\xymatrix{1_\cC\ar[d]_-{\Theta_M} \ar[r]^-{\sim}& 1_\cC\ar[d]^-{\Ll_{\Mv}} \\  
\cA_G \ar[r]^-{\bL}_-{\sim}  &\cB_\Gv}$$ 
\end{conjecture}

\index{automorphic quantization}
\index{spectral quantization}

\medskip
Let us spell out how this meta-conjecture encodes various structures that we have encountered throughout the main text.

$\bullet$ [Local geometric] 
The local object $\Theta_{M}\in \cA_G(D^*)$ is given in the polarized case $M=T^*X$ by the category of sheaves $\SHV(X_{F})$ with its $G_F$-action. Its unramified part $\Theta_{M}^{unr}$ is the unramified local category $\HECKE^X_G=SHV(X_F/G_O)$, a factorization algebra in $\HECKE_{\cA_G}$-module categories. The unit, i.e., the basic sheaf $\Phi\in \Theta_{M}^{unr}$, is the constant sheaf on $X_O$.

The local object $\Ll_{\Mv}\in \cB_\Gv(D^*)$ is given in the polarized case $\Mv=T^*\Xv$ by the category $QC^!(\Loc_\Gv^\Xv(D^*))$ of ind-coherent sheaves on the stack of locally constant maps into $\Xv/\Gv$, considered relative to the stack $\Loc_\Gv(D^*)$ of Langlands parameters.
 Its unramified form $\Ll_{\Mv}^{unr}$ is identified with $QC^!(\Ll \Xv / \Gv)\simeq \QCshear(\Mv/\Gv)$, the unramified local spectral category. 
 
Conjecture~\ref{local conjecture} identifies $\Theta_{M}$ and $\Ll_{\Mv}$ as module categories for the line operators (Hecke category) $\HECKE_{\cA_G}\simeq \HECKE_{\cB_\Gv}$, and is expected to upgrade to an equivalence of factorization algebras in Hecke modules. 

The ramified local geometric duality is the (currently imprecisely formulated) expectation that $\Theta_{M}$ and $\Ll_{\Mv}$ are identified under the conjectural local geometric Langlands correspondence (in fact compatibly with factorization).

\medskip

$\bullet$ [Local observables] The boundary observables, the factorization associative algebra $\HECKE_{\Theta_M}\in \HECKE_{\cA_G}$, recover the {\em Plancherel algebra} or relative Coulomb branch algebra of \S \ref{PlancherelCoulomb}. Indeed the description of $\HECKE_\Theta$ in terms of the theory on $S^2\times I$ amounts to the description of the Plancherel algebra as homology of the relative Grassmannian, while the description as internal endomorphisms of the basic sheaf (the unit) $\Phi\in \Theta_{M}^{unr}$ in the Hecke category recovers the definition of the Plancherel algebra.\index{Plancherel algebra}

On the spectral side we find the local $\Ll$-observables $$\HECKE_{\Ll_{\Mv}}\simeq \bO_\Mv,$$ the spectral deformation quantization of functions on the Hamiltonian $\Gv$-space $\Mv$ (Definition~\ref{sdq}). This is a factorization associative algebra in the Hecke Category $\HECKE_{\cB_\Gv}$, quantizing the shear of the moment map $\Mv/\Gv\to \fgxv/\Gv$. In other words, we recover $\Mv$ and its Hamiltonian $\Gv$-action directly out of its spectral quantization $\Ll_\Mv$ by passing to cohomology of local observables. 

\

\medskip

$\bullet$ [Global geometric, unramified]
 The boundary condition $\Theta_M\in \cA_G$ defines a ``$\Theta$-series'' functional $$\Theta_M(\Sigma): \cA_G(\Sigma)=\AUT^?(\Sigma)\to 1_\Ccirc,$$
which is representable by an object in the automorphic category: the $M$-period sheaf $\cP_M$. Likewise $\Ll_\Mv\in \cB_\Gv$ determines an object
$$\Ll_\Mv(\Sigma)\in QC^!(\Loc^?_\Gv(\Sigma)),$$ the $L$-sheaf of $\Mv$. The duality between these recovers the global period conjecture.

\medskip

\index{RTF algebra}
$\bullet$ [Global observables, unramified] The global observables $\bH_{\Theta,\Sigma}=\int_\Sigma \HECKE_\Theta \in \bH_{\cA_G,\Sigma}$ recover the RTF algebra, an associative algebra object in the global Hecke category $\bH_\Sigma$. The RTF algebra acts by endomorphisms of the period sheaf $\Theta_M(\Sigma)\in \cA_G(\Sigma)$, compatibly with the $\bH_\Sigma$-action.

Likewise the factorization homology $\bH_{\Ll,\Sigma}\in \bH_{\cB_\Gv,\Sigma}$ of the local $\Ll$-observables gives the algebra object of $\Ll$-observables studied in \S \ref{spectral-factorization}, which acts as endomorphisms of the $\Ll$-sheaf $\Ll_\Mv(\Sigma)\in \cB_{\Gv}(\Sigma)$.

\medskip

$\bullet$ [Correlation functions and $\Theta$-series] 
 \index{$\Theta$-series}
For a finite subset $S\subset \Sigma$, we may consider the ``$M$-ramified'' global category $\langle \bigotimes_S \Theta_{M,x_i} \rangle_{\cA_G,\Sigma},$ a module for the global Hecke category away from $S$ $\bH_{\cA_G,\Sigma\setminus S}$. 
In the polarized case $M=T^*X$ this is the category of sheaves on the stack of $G$-bundles with $X$-level structure along $S$ (i.e., with a section of the associated $X$-bundle on the punctured neighborhood of $S$). 
The boundary state $\langle  - \rangle_{\Theta,S}$ is then [represented by] an object of this category, preserved by the S-RTF algebra $$\bH_{\Theta_M,\Sigma\setminus S}=\int_{\Sigma\setminus S} \HECKE_{\Theta_M}.$$
(In the case where $S$ is empty, we recover [the functor represented by] the $X$-period sheaf.) 

This object $\langle  - \rangle_{\Theta,S}$ is an avatar of the $\Theta$-series operation: assuming sufficient dualizability, such an object is equivalent to the data of a functor 
$$\bigotimes_{x_i\in S} \Theta_{M,x_i}\to \SHV(Bun_{G}(\Sigma,S))$$ respecting $\prod_S G(F_{x_i})$-actions
from the local $X$-category on $S$ to sheaves on the stack of bundles with full level structure along $S$.
(See \S \ref{automorphic-factorization} for the unramified version of $\Theta$-series.)
In the colimit over $S$, the factorization homology $$\int_\Sigma \Theta_{M} \in \int_\Sigma \cA_G$$ is a geometric analog of the $G(\AA_\Sigma)$-representation $L^2(X(\AA))$, and the state $\langle - \rangle_{\Theta}$ plays the role of the ad\`elic $\Theta$-series.

We have a parallel story on the spectral side: the category $\langle \bigotimes_S \Ll_{\Mv,x_i}\rangle_{\cB_\Gv,S}$ is given in the polarized case $\Mv=T^*\Xv$ by ind-coherent sheaves on the stack of local systems with $\Xv$-fixed points structure along $S$ (i.e., with a flat section of the associated $\Xv$-bundle around $S$). The state $\langle  - \rangle_{\Ll_{\Mv,x_i},S}$ is represented by an object of this category -- a ramified generalization of the $\Ll$-sheaf -- preserved by the action of the $\Ll$-observables away from $S$ 
$$\bH_{\Ll_\Mv,\Sigma\setminus S}=\int_{\Sigma\setminus S} \bO_{\Mv}.$$ The unramified part of this construction 
recovers the $\Ll$-functor construction discussed in \S \ref{spectral-factorization}, while for $S$ empty we recover the $\Ll$-sheaf $\Ll_{\Xv}\in QC^!(Loc_\Gv(\Sigma))$.

\medskip
 $\bullet$ [Arithmetic] We leave to future work to spell out the form of these conjectures obtained by taking Frobenius traces everywhere. 
At the top level, taking Frobenius trace of the period sheaves we obtain the period, as a linear functional on the global space of automorphic forms, and taking Frobenius trace of the $L$-sheaves we obtain a geometric avatar of the L-function, realized as a derived volume form on the stack of arithmetic local systems. However the necessary dualizability conditions to perform these trace only hold if we localize (on the automorphic side via the Hecke action) to the open locus $\Loc^\circ$, away from the poles of the L-function, where the map $\Loc^\Xv\to \Loc_\Gv$ is proper.

.

\begin{remark}[Interfaces and functoriality]
More generally, it is an interesting problem to understand Langlands functoriality and its generalizations in terms of interfaces of field theory. For example, a group homomorphism $\Hv\to \Gv$ induces natural morphisms on moduli of local systems and hence a natural interface $\cB_{\Hv}\to\cB_{\Gv}$. Applying the Langlands corresponds produces a mysterious conjectural morphism $\cA_H\to \cA_G$ of automorphic theories. But there are many more interfaces $\cB_{\Hv}\to\cB_{\Gv}$ -- in particular, those coming from spectral quantization of hamiltonian $\Hv\times \Gv$-spaces. In other words, we might think of the assignment $\Gv\to \cB_\Hv$ as a ``quantization'' functor out of a higher category of reductive groups and bihamiltonian actions, or equivalently (in the language of Remark~\ref{fancy induction}) shifted Lagrangian correspondences between the 3-shifted symplectic stacks $T^*[3]pt/\Gv=\fgx[2]/\Gv$. As remarked in \S \ref{groupcase}, in a formal sense, i.e., ignoring problems with duality, and suppressing the structure of composition, the study of functoriality can be subsumed in the study of boundary theories. For example the theta correspondence can be viewed either as an operation on spaces of automorphic forms or as the study of a particular period for a product group. 
\end{remark}

%% file: MiscellaneousComputations.tex
\newcommand{\Whit}{\mathrm{Whit}}
\newcommand{\GrXtilde}{\widetilde{\Gr^X_G}}
\newcommand{\Fpbar}{\overline{\mathbb{F}_p}}
\section{Some miscellaneous computations}
\label{Eisappendix}
We gather here various computations that were postponed from the main body of the
text, to avoid disturbing too far the flow of proof. 

\subsection{Unnormalized Eisenstein periods}  \label{eis period numerical comp}
We carry out the computation of the numerical period in the case of of minimal unramified Eisenstein series.
The computation is in fact completely straightforward, but it is delicate in the matter
of signs and shifts. The shift is interesting, and this is why we explicate it, although
we will be quite terse.   This computation
has been alluded to in the main text at various points (\S \ref{Eisgeom}, \S \ref{Eisperiodmaintext}).

  After setup, we derive in Lemma \ref{nwnEis} the usual formula for the constant term of the Eisenstein series, 
  adapted to our notation. In \S \ref{Eisperiodappendix} we then pass from
  this to the computation of the numerical period.

  \subsubsection{Notation}
  Let $X^*$ be the cocharacter group of the  torus quotient $A$ of a Borel subgroup $B \subset G$, and let
\begin{equation} \label{Xstartensor} \lambda \in X^* \otimes (\mbox{everywhere unramified characters of id{\`e}le class group})\end{equation}
  (tensor product of abelian groups). 
  Such a $\lambda$ defines a character $A(\mathbb{A}) \rightarrow \C^{\times}$. 
For such $\lambda$, and $\alpha$ a root of $G$, we can form  \index{$\lambda_{\alpha}$}
$$\lambda_{\alpha} := \langle \lambda, \alpha^{\vee} \rangle,$$
  which is now an everywhere unramified character of the idele class group. 
  
  In what follows we regard $2\rho$ both as an element of $X^*$ but also
  as an element of the above group  \eqref{Xstartensor} via tensoring with the
  norm $|\cdot|$  on the idele class group which sends
  each uniformizer $\varpi_v$ to $q_v^{-1}$.
  Let $\partial$ be an idele representing the diiferent, as in \eqref{diffdef}, i.e.
  an idele of everywhere even local valuation equal to $n_v$, where the $n_v$
  are the vanishing orders of a $1$-form.  As in \eqref{psidef} this normalizes
an additive character $\psi: \mathbb{A}/F \rightarrow \C^{\times}$,
which is locally trivial on $\partial_v^{-1}$ but nontrivial on any larger open compact subgroup of $F_v$. 
  As an example of our notation  $\langle 2 \rho, \alpha^{\vee} \rangle (\partial) =|\partial|= q^{-(2g-2)}$;
  In what follows we write 
  $$ D = |\partial|^{-1} = q^{2g-2}$$
  for the ``discriminant.''
  
  \subsubsection{Notation concerning completed zeta functions} 
  Let $\xi(-)$ be the completed $L$-function of the global field; we have the functional equation
  $\xi(0,\chi) = \chi(\partial)  q^{g-1} \xi(1, \chi^{-1})$.
  As in \eqref{Lnormdef} we put
 $$\xi^{\norm}(s,\chi) = \epsilon(s, \chi)^{-1/2} \xi(s,\chi)
 =\chi(\partial)^{-1/2} q^{(s-1/2)(g-1)} \xi(s, \chi).$$
so that  
\begin{equation} \label{01ch} \xi^{\norm}(0, \chi) = \chi(\partial)^{-1/2} q^{(1-g)/2} \xi(0,\chi),
\xi^{\norm}(1, \chi)= \chi(\partial)^{-1/2} q^{(g-1)/2} \xi(1, \chi).\end{equation}
If we omit $s$ it means that we take $s=0$. Thus $\xi^{\norm}(\chi') = \xi^{\norm}(\chi)$, with $\chi' := |\cdot| \chi^{-1}$. 
Suppose now that $\mathfrak{s}$ is any $T$-stable subspace
of $\mathfrak{g}/\mathfrak{t}$; then we write
\begin{equation} \label{xisubspacenotn} \xi(s,\mathfrak{s} ) = \prod \xi(\lambda_{\alpha}, s)\end{equation}
where the product is taken over roots  $\alpha$ occuring 
in $\mathfrak{s}$. Similarly for $\xi^{\norm}$ etc.

  \subsubsection{The Eisenstein series}
We briefly summarize the construction of Eisenstein series
and compute the  constant term after Whittaker normalization in Lemma \ref{nwnEis}. 

  We use the standard adelic uniformization of the rational points of $\Bun_G$.
For $\lambda$ in \eqref{Xstartensor} we let
   \begin{equation} \varphi_{\lambda}: G(\mathbb{A}) \rightarrow \C\end{equation}
   be the unique left $N(\mathbb{A})$-invariant and right $K=G(\widehat{\mathfrak{o}})$-invariant
   function whose value on   $A(\mathbb{A})$ is given by $\lambda+\rho$
   (here, we interpret $\rho$ as $2 \rho \otimes |\cdot|^{1/2}$
   inside the group \eqref{Xstartensor}).
    The standard Eisenstein series is  obtained  by  summing this over the $F$-points of $G/B$;
    this is convergent for $\lambda$ sufficiently ``positive,'' e.g.\ $\lambda = t \rho$ for $t$ positive
    and extends by meromorphic continuation to other $\lambda$. 
 
   Let us normalize Haar measure on $U(\mathbb{A})$ so that the mass of $U(F) \backslash U(\mathbb{A})$ is trivial.
   This measure  equals $D^{-\dim(U)/2}$ multiplied by the measure which assigns to each $U(\mathfrak{o}_v)$ the mass $1$. 
   
The constant term of $E_{\lambda}$ 
may be computed as usual using the Bruhat decomposition
of $B_F \backslash G_F$ into $U_F$-orbits of the form $w U_F$
with stabilizer $U^w_F = \Ad(w) U_F \cap U_F$. This gives:

      \begin{equation} \label{GWB} \int_{[U]}  E_{\lambda}(ug)  = \sum_{w} \int_{U^w_F \backslash U(\mathbb{A})}  \varphi_{\lambda}(wu g) \stackrel{GK}{=}
  \sum_{w \in W} \prod_{\alpha >0, w \alpha <0} \frac{D^{-1/2}  \xi(0, \lambda_{\alpha}) }{\xi(1, \lambda_{\alpha})} \varphi_{w \lambda}\end{equation}
(here the $D^{-1/2}$ arises from normalization of measure, as explained above, and ``GK'' stands for Gindinkin and Karpilevic).  
The Whittaker coefficient has a similar formula with only the long Weyl element contributing;  we obtain
 $$ E_{\lambda}^{\psi}(g) := \int_{[U]} E_{\lambda}(ug) \psi(u) du =  \prod_{\alpha} D^{-1/2}  \frac{1}{\xi(1, \lambda_{\alpha})} W_{\lambda},$$ with $W_{\lambda} 
 =\prod W_{\lambda, v}$, and
 $W_{\lambda,v}$ is given by $\xi_v(1, \lambda_{\alpha})\cdot   \int \varphi(w u g) \psi(u) du $,
 where the integral is taken over $U(F_v)$ with respect to the measure assigning mass one to
 $U(\mathfrak{o}_v)$.  
 
 We write
$  \psi(u) = \psi_0(a_{0,v}^{-1} u a_{0,v})$ 
with (cf. \eqref{a0def}): 
  $$a_{0,v} = e^{2\rho^{\vee}}(\pi_v^{-n_v/2}), a_0 = (a_{0,v})_v =e^{2\rho^{\vee}}(\partial^{-1/2}).$$
 so that $\psi_0$ is ``unramified,''  i.e.
 equal to one on $U(\mathfrak{o}_v)$ but not on  larger compact subgroups; 
 then  
  \begin{eqnarray} W_{\lambda,v}(g_v) = \xi_v(1, \lambda_{\alpha}) \int \varphi(w   a_{0,v} u a_{0,v}^{-1}g) \psi_0(u) d(\Ad(a_{0,v}) u)
  \\ =  
| \langle \lambda + \rho-2\rho, \rho^{\vee} \rangle  (\pi_v^{n_v}) |
  W^{0}(a_{0,v}^{-1} g),
\end{eqnarray}
where $W^0$ is the unramified local Whittaker function with $W^0(1)=1$,
explicitly given as $\int_{U_v} \varphi(w u g) \psi_0(u) du$. 
 Globally we get 
  $$E_{\lambda}^{\psi}(g) = \left(  \prod_{\alpha} D^{-1/2}  \xi(1, \lambda_{\alpha})^{-1}  \right) \langle \lambda -\rho, \rho^{\vee} \rangle(\partial) W^0(a_{0}^{-1} g).$$
By Example
 \ref{Whit explicit computation example} the
 normalized period
  $ \langle P_{\Whit}^{\norm}, E_{\lambda} \rangle$ thereby equals 
  the value of this at $g=a_0$, 
 multiplied by $q^{\beta/2}$
 where
 $  \beta = (g-1) [   \dim U- \langle 2 \rho, 2 \rho^{\vee} \rangle ]$.
Since $\langle \rho, \rho^{\vee} \rangle (\partial)
  =q^{-2(g-1) \langle \rho, \rho^{\vee} \rangle}$ we get
  $$   \langle P_{\Whit}^{\norm}, E_{\lambda} \rangle =q^{b_U/2} \left(  \prod_{\alpha} D^{-1/2}  \xi(1, \lambda_{\alpha})^{-1}  \right)
 \langle \lambda, \rho^{\vee} \rangle(\partial) $$
 with $b_U = (g-1) \dim U$, as usual. 
Taking into account that $D^{-\dim U/2} = q^{-(g-1) \dim U} = q^{-b_U}$ we get, in the shorthand of
\eqref{xisubspacenotn}:
\begin{equation} \label{Whit period Eis}   \langle P_{\Whit}^{\norm}, E_{\lambda} \rangle =q^{-b_U/2}  \xi(1, \mathfrak{n})^{-1}
 \langle \lambda, \rho^{\vee} \rangle(\partial) \end{equation}

  \begin{lemma} \label{nwnEis}
Write 
  $\mathcal{E}_{\lambda} := \frac{E_{\lambda}}{\langle P_{\Whit}^{\norm}, E_{\lambda} \rangle }$.
Then the constant term of $\mathcal{E}_{\lambda}$, for the Haar probability measure on $U(\mathbb{A})/U(F)$, is given by
\begin{equation} \label{Whiteisperiod} \mathcal{E}_{\lambda}^U =   \sum_{w \in W} \xi^{\norm}(0, w^{-1} \overline{\mathfrak{n}}). \varphi_{w \lambda},\end{equation}
where on the  right we again use the shorthand of \eqref{xisubspacenotn}.
 \end{lemma}
 
 \proof
  Dividing \eqref{GWB} by \eqref{Whit period Eis} we get
\begin{equation} \label{cwdef} 
  \mathcal{E}_{\lambda}^U =  q^{b_U/2} \sum_{w \in W} \langle - \lambda,   \rho^{\vee}  \rangle (\partial)  \cdot c_w    \varphi_{w\lambda }, 
\mbox{ where }c_w =  \prod_{\alpha > 0 }  \begin{cases}  D^{-1/2} \xi(0, \lambda_{\alpha}), & w \alpha < 0 \\     \xi(1, \lambda_{\alpha})    ,  & w \alpha > 0 \end{cases}.
\end{equation} 
We must check that 
\begin{equation} \label{cwsat}  \langle - \lambda, \rho^{\vee} \rangle (\partial) \cdot c_w =  q^{-b_U/2} \xi^{\norm}(0, w^{-1} \overline{\mathfrak{n}}).\end{equation}
 for which we rewrite 
 the definition \eqref{cwdef} via \eqref{01ch} and the functional equation. 
 $$ c_w =  \prod_{\alpha > 0 } D^{-1/4}   \begin{cases}  \lambda_{\alpha}(\partial)^{1/2}  \xi^{\norm}(0, \lambda_{\alpha}), & w \alpha < 0 \\   \lambda_{\alpha}(\partial)^{1/2}   \xi^{\norm}(0, -\lambda_{\alpha})    ,  & w \alpha > 0 \end{cases}.
 $$
Our result follows after observing that the set $\{\alpha: \alpha>0, w \alpha <0\} \cup  \{-\alpha: \alpha > 0,  w \alpha > 0 \}$
 is precisely the set of roots of the form $\{w^{-1} \beta: \beta < 0\}$, and that $D^{\dim(U)/4} = q^{b_U/2}$.

 \qed

  \subsubsection{The numerical Eisenstein period}\label{Eisperiodappendix}
  
 The space of interest is $X=U \backslash G$ as a $G \times T$
 space (i.e., via the right action $(g,t): U x \mapsto U t^{-1} x g$);
 the point stabilizer is $T^{\Delta} U \hookrightarrow G \times T$
 with $T^{\Delta}$ the diagonal copy of $T$:
\begin{equation} \label{Xprez2} X \simeq (T^{\Delta} U) \backslash G \times T.\end{equation}

 The putative dual space is $$\check{X} = \check{G}/\check{U}$$
   with action (expressed on the left given by $(\check{g}, \check{t}): x \check{U} \mapsto g x t \check{U}$,
   that is to say $\check{X} = G/T^{-\Delta} U$ in evident notation). 
  
    We take the $\GGm$ action on both $X$ and $\check{X}$ to be trivial. 
  The modular characters
 are given by $\eta_{\aut}: (g,t) \mapsto e^{2 \rho}(t)$
and the same for $\check{X}$. 
\begin{example}
Here is an example, to help with figuring signs. In the case $G=\SL_2$, we have $X=\mathbb{A}^2-\{0\}$
where $\Gm$, identified with $T$ via the positive coroot, acts by scaling
and $\SL_2$ acts by right multiplication, and $\check{X} = \mathbb{A}^2-\{0\}/\mu_2$,
where  $\Gm$, identified with $\check{T}$ via the posiitve root, 
acts via scaling and $\PGL_2$ acts by left multiplication.
\end{example}
  
   Now consider the automorphic form on $G \times T$ given by
   $\varphi = \mathcal{E}_{\lambda} \boxtimes (w \lambda)^{-1}$,
  for $\lambda$ as in \eqref{Xstartensor}, and some $w \in W$, the Weyl group for $G$. 
 We now  compute the ratio of periods \begin{equation} \label{GTUratio} \frac
 { \langle P^{\mathrm{norm}}_{X}, \varphi  \rangle}
 {\langle P^{\mathrm{norm}}_{\mathrm{Whit}}, \varphi   \rangle}
 \end{equation}
The denominator equals one by choice of $\mathcal{E}_{\lambda}$ as in Lemma \ref{nwnEis}, for the normalized Whittaker period is the same computed on $G$ or on $G \times T$. For 
 the numerator   we use  \eqref{Xprez2} as well as \eqref{PXform3};
  im the notation of the latter we have
  $$ \langle P^{\norm}_{G/U}, \varphi   \rangle
  = q^{-(b_T+b_U)/2} \int_{U_F T_F^{\Delta} \backslash (T \times G)_{\mathbb{A}}}  \Phi^0(g) |\eta_{\aut}(g)|^{1/2} \varphi(g) $$
  where $\Phi^0$ is the characteristic function of $U_{\mathbb{A}} T_{\mathbb{A}}^{\Delta} 
  (G \times T)(\widehat{\mathfrak{o}})$ and the measure on $T \times G$
  assigns mass $1$ to the maximal compact. 
    The modular character here is given by $\eta: (t,t) \mapsto e^{2\rho}(t)$,
   cf. discussion of \S \ref{Htrans}:
 the right adjoint action of $T^{\Delta}$ on the tangent space at the identity  coset of
 $(T^{\Delta} U) \backslash (T \times G)$ is given by
 $\mathrm{Ad}(t^{-1})$ on $\mathfrak{b}^-$. 
 Using Iwasawa decomposition and noting
 that the integrally normalized meausre $du$ is $q^{b_U}$ multiplied by 
 the probability measure, we rewrite the above as:
 $$ q^{(b_U-b_T)/2}  \int_{T_F\backslash T_{\mathbb{A}}} \mathcal{E}_{\lambda}^U(t) (w\lambda)^{-1} (t)
 |e^{-\rho}(t)| $$
 where $dt$ is now the Haar measure on $T$ corresponding to the standard measure on $\Bun_T$
 where each bundle is weighted by inverse-automorphism;
 the factor $e^{-\rho}$ arises from $e^{-2\rho}$ from the measure in Iwasawa decomposition,
 combined with the square root of the modular character.

  In fact, the right hand side of \eqref{GTUratio} will diverge, since it involves (in effect) integrating a character over $T$. 
  We treat this purely formally: we regard it as nonzero only when the character is trivial,
  in which cae it will be given by the $\zeta$-value computing the volume of $[T]$,
which we shall formally understand to $q^{(g-1) r} \xi(1)^r=q^{b_T} \xi(1)^r$ with $r=\dim(T)$
(this is only a formal expression, for $\xi$ has a pole at $1$).

  Computing formally  thus,  \eqref{GTUratio} equals 
\begin{equation} \label{GTUratio2}  q^{b_U/2}  q^{b_T/2}  \xi^{\norm}(0, w^{-1} \overline{\mathfrak{n}}) \cdot \xi(1)^{\mathrm{rank}}  = 
   q^{b_U/2} \xi^{\norm}(0, w^{-1} \overline{\mathfrak{b}}).\end{equation}

  We compare this with what we would expect starting from \eqref{PXshriekprediction} 
 (which is of course formulated only for cusp forms). The dual space is
 $\check{X} = \check{G}/\check{U}$  with left action given by  $(\check{g}, \check{t}): x  \check{U}   \mapsto   \check{g} x \check{t} \check{U}$;
the parameter of $\varphi$
 is given by $(\lambda, (w \lambda)^{-1})$, and fixes precisely $\check{x} = w^{-1} \check{U}$. 
 The adjoint action of $\pi_1$ on the tangent space $T_{\check{x}} \check{X}$
 corresponds to the adjoint action of $\lambda$ on $w^{-1} \overline{\mathfrak{b}}$. 
 Therefore, taking the ratio of \eqref{PXshriekprediction} for both $G/U$ and the Whittaker period
 we would be led instead to a prediction of $ \xi^{\norm}(0, w^{-1} \overline{\mathfrak{b}})$;
 that is to say, 
   \eqref{GTUratio} is in line with \eqref{PXshriekprediction} {\em except for the factor $q^{b_U/2}$}. 
This discrepancy -- which, geometrically, would manifest itself as a shift $\langle b_U \rangle$ --
was already mentioned in the text, and would be interesting to understand. 
The paper \cite{ChenVenkatesh} contains other examples which exhibit similar discrepancies.

 \subsection{Numerical derivation of the effect of twisting}
   \begin{example} \label{num twisting} 
   We revisit Lemma \ref{shiftpredlem} from a numerical standpoin, examining  the effect of twisting the choice of $K^{1/2}$. This computation serves solely as a sign-check on that Lemma.
   We put ourselves in the situation of Conjecture \ref{numconj}.    
  Returning to \eqref{PXform2} 
   \begin{equation}   P_X^{\norm} : g \in G(\mathbb{A}) \mapsto q^{-\beta_X/2}  |\eta_{\aut}(g)|^{1/2} \sum_{x \in X(F)} g \cdot (\partial^{1/2} \cdot \Phi(  x)),.\end{equation}
  and replace $X$ by $X' = X[\eta_{\spec}^{\vee,-1}]$. Then $\beta_{X'} =  \beta_X  + \tau$ with $\tau=(g-1) \langle \eta_{\aut}, \eta_{\spec} \rangle$ we get 
 the result is 
 $$P_{X'}=   \sum_{x} (g, \partial^{1/2}) \cdot' \Phi(x) = 
   \eta_{\spec}^{\vee}(\partial^{-1/2})  \left[ q^{\beta_X/2}  |\eta_{\aut}(g)|^{-1/2} P_X^{\norm} \right].$$
   
   Pairing this against $f$ with central character $\omega_f$ gives 
   can be computed in terms of  the same pairing for the normalized period:
   $$P_{X'}(f) =q^{\beta_X/2}   \omega_f(\eta_{\spec}^{\vee}(\partial^{1/2}))   P_X^{\norm}(f') $$
   where $f'$ is $f$ twisted by $|\eta_{\aut}|^{-1/2}$.   
   This has the effect of twisting the Galois parameter of $f$ through the composite
   $ \textrm{cyclotomic}^{1/2} \circ \eta_{\aut}^{\vee, -1}$.  
   Applying Conjecture \ref{numconj} we find
   \\\
 $$ P_{X'}(f) =  q^{(\beta_X -b_G)/2}  \omega_f( \eta_{\spec}^{\vee}(\partial^{1/2}))  L^{\norm}(0,   T^{\shear,'}) $$
where $T^{\shear}$ is the tangent space
sheared by the twisted action for $\check{X}' = \check{X}[\eta_{\aut}^{\vee, -1}]$. 

Now $ \omega_f( \eta_{\spec}^{\vee}(\partial^{1/2})) $
is the square root of the central $\epsilon$-factor for the Galois parameter
of $f$ acting on $\det(T)$. We apply \eqref{epsilon0k} which says  
$$ \epsilon(0, T^{\shear,'})=   \varepsilon(1/2, T) q^{b_G - \beta_{\check{X}'}} $$
to rewrite 
$$P_{X'}(f) = q^{\beta_{X}/2 + \beta_{\check{X}'}/2 -b_G} \ L^{\norm}(0,   T^{\shear,'}) $$
and using $\beta_{\check{X}'} = \beta_{\check{X}} - \tau$ this becomes
$$P_{X'}(f) =  q^{-b_G} q^{\beta_{X}/2 + \beta_{\check{X}'}/2 - \tau/2} L^{\norm}(0,   T^{\shear,'}) $$
 which can be compared with \eqref{LX0}, i.e.
 the conclusion is that $X'$ and $\check{X}'$ are dual,
 but one has to twist the period sheaf by $\beta_X+\beta_{\check{X}}-\tau$.

  \end{example}

\qed

\subsection{Proof of Proposition \ref{anomalousautomorphic}} \label{prooflemmametaplectic}

\proof (sketch)
Soul{\'e}'s etale Chern class \cite[II]{Soule} defines a morphism
\begin{equation} \label{c22} c_{22}: H_2(\Sp_{2r}(F), \Z/2) \rightarrow H^2(F, \Z/2)\end{equation}
which can be verified (comparing with the discussion of  \cite{Deligne} and \cite{Prasad})
to define the metaplectic extension of $\Sp_{2r}(F)$, that is to say, 
the unique topological double cover thereof.

Now consider the two representations of $H$ defined by
$$H \rightarrow \mathrm{Sp}_{2r} \mbox{ and }H \stackrel{\theta}{\rightarrow} \mathbb{G}_m
\hookrightarrow \SL_2, $$ where the last map  is the inclusion
of a maximal torus.  It follows from \eqref{Vtheta} that these two representations of $H$ have the
same second Chern class. Therefore  the two  maps
$$H_2 (H(F), \Z/2) \rightarrow H^2(F, \Z/2)$$
arising from pulling back \eqref{c22} via $H \rightarrow \mathrm{Sp}_V$
or $H \rightarrow \SL_2$ must coincide. 

 We claim that the latter pullback of $c_{22}$, via $H \rightarrow \SL_2$ is trivial.
 Indeed, this factors through $\mathbb{G}_m(F)$, and our claim follows 
 from the fact that 
 the metaplectic cover of $\SL_2(F)$ splits over $\mathbb{G}_m(F) = F^{\times}$. \footnote{
At the level of cocycles, this amounts to fact
that the Hilbert symbol $F^{\times} \times F^{\times} \rightarrow \pm 1$
is trivialized as $2$-cocycle by
 the $1$-cocycle 
$ \eta:  u \varpi^n \mapsto \bar{u}^n \cdot (-1,-1)^{n \choose 2}$
where we write, for short, $\bar{u}  =\pm 1$ according to whether $u$ is residually a square or not.
} Therefore, the pullback of the metaplectic cover under $H \rightarrow \mathrm{Sp}_{2r}(F)$ does too. 

\qed

 The construction above globalizes:  
 if $F$ is a global function field,  then the metaplectic cover of $\mathrm{Sp}(\mathbb{A}_F)$
 splits upon pullback to $H(\mathbb{A}_F)$ compatibly with its standard splitting on $\mathrm{Sp}(F)$.
 Moreover, if we push out the metaplectic cover by $\pm 1 \rightarrow S^1$, the result
 above remains valid for 
 $F=\R$ or an extension of $\mathbb{Q}_2$. 
 In relation to rationality issues we note at least the following:
 
 \begin{lemma} \label{anomaly rationality}
Suppose $\rho: H \rightarrow \mathrm{Sp}$
 is a homomorphism of split $\Z$-groups,
 and that after base change to $\C$ the Chern class condition of Proposition  \ref{anomalousautomorphic}
 is satisfied (for some character $\theta: H \rightarrow \Gm$, which
 we may as well suppose defined over $\Z$). 
 
  Then the metaplectic cover pulled back to $H(F)$ splits over any field $F$  
in which $\pm 2$ are both squares. 
 \end{lemma}
 
 A similar assertion holds for $\Z[1/N]$,  with identical proof, but now additionally requiring prime divisors of $N$ to be squares. 
\proof
We are going to verify that the condition of Proposition \ref{anomalousautomorphic}, i.e.
$$c_2(\rho) = c_1(\theta)^2$$
in absolute {\'e}tale cohomology {\em over $F$}. 

To do so we will first compute in {\'e}tale cohomolgy over $\Z[1/2]$ and first show that 
there the difference is decomposable:
\begin{equation} \label{Apollo1} c_2(\rho) - c_1(\theta)^2 \in H^1 \cup H^3,\end{equation}
from which we will readily deduce the result in the final paragraph.

Let $P$ be a Borel subgroup of $H$ over $\Z$. 
Note that  $ c_2(\rho) - c_1(\theta)^2 $
dies in the mod $2$ {\'e}tale cohomology of $BP_{\Z[1/2]}$
(i.e., the classifying space of $P$ as a group scheme over $\Z[1/2]$,
which can be constructed as a simplicial $\Z[1/2]$-scheme). 
To see this, we note that $P$ can be replaced by a maximal torus,
and the cohomology $BT_{\Z[1/2]}$ is a polynomial
ring over $H^*(\Z[1/2])$ generated by the Chern classes
of a basis of $\Hom(T, \Gm)$. Both $c_2(\rho)$ and $c_1(\theta)$
can be expressed as polynomials in these Chern classes 
(with coefficients in $\Z/2\Z$); the subring of such classes in $H^*(BT_{\Z[1/2]})$
maps isomorphically to $H^*(BT_{\C})$, and therefore
the coincidence of $c_2(\rho)$ and $c_1(\theta)^2$ in $BT_{\C}$
implies the same over $\Z[1/2]$.

Now consider the map in {\'e}tale cohomology 
induced by $BP_{\Z[1/2]} \rightarrow BH_{\Z[1/2]}$.
We will show its kernel is decomposable,  and in fact belongs
to the right hand side of \eqref{Apollo1}, concluding the proof. 
 By the Serre spectral sequence
(applied in $H$-equivariant cohomology with mod $2$ coefficients for the fibration $\pi: H/P \rightarrow \Spec \Z[1/2]$)
we get a spectral sequence converging to the {\'e}tale cohomology of $BP_{\Z[1/2]}$
and whose $E_2^{pq}$ term is given by 
$$H^p(BH_{\Z[1/2]}, R^q \pi_* \Z/2).$$
The $R^q \pi_*$ terms are vanishing only for $q$ even;
we will repeatedly use this without comment. The only possible obstruction
to the injectivity of the edge map $H^4(BH) \rightarrow H^4(BP)$
comes from the differential $d_3: E_3^{12} \rightarrow E_3^{40}$. 
But $E_3^{12}= E_2^{12} = H^1(BH_{\Z[1/2]}, R^2 \pi_* \Z/2)$
and the action of $\pi_1(BH_{\Z[1/2]}) = \pi_1(\Z[1/2])$ on $R^2 \pi_* \Z/2$
is trivial (in fact, this is the mod $2$ reduction of $R^2 \pi_* \Z_2$
which is a sum of copies of $\Z_2(1)$, indexed by Schubert cells of codimension $1$). 

Consequently, each class in $E_2^{12}$ is a product
of a class in $E_2^{10} = E_3^{10}$ and $E_2^{02}=E_3^{02}$.
 Since the differential $d_3$ is a derivation,
we conclude that the image of such a class under $d_3$
has the form $E_3^{10} \cup E_3^{30}$. 
However, $ E_3^{10}=E_2^{10}  =H^1(BH_{\Z[1/2]}, \Z/2)$ and $E_2^{30} = E_3^{30} = H^3(BH_{\Z[1/2]}, \Z/2)$, 
so the image of such a class under $d_3$ corresponds
to a class in $H^4(BH)$ that decomposes as $H^1 \cup H^3$,
i.e.
is decomposable as in \eqref{Apollo1}.
 The kernel of
$H^4(BH) \rightarrow H^4(BP)$ consists entirely of such classes, and so we are done with the proof of \eqref{Apollo1}.

Finally, the claimed result follows from \eqref{Apollo1}, because
the natural map  $H^1(\Z[1/2], \Z/2) \rightarrow H^1(BH_{\Z[1/2]})$ is an isomorphism; the source group
 coincides with the group of square classes in $\Z[1/2]^{\times}$, and this is killed
by passage to any $F$ as in the statement.
\qed  

\subsection{Proof of Proposition \ref{dumblemma}} \label{dumblemmaproof}

In the main text, we used the following perhaps intuitively obvious
manifestation of the rigidity of homomorphisms between reductive groups:

     \begin{proposition} \label{dumblemma}
    Suppose
    that  $G_1, G_2$ are split groups over $\Z$ and
    $\mathcal{F}$ 
a conjugacy class of complex homomorphism $ (G_1)_{\C} \rightarrow (G_2)_{\C}$
containing a homomorphism defined over $\Q$. 
 For any sufficiently large $p$, the following is valid:
    
\begin{quote}
There exists, up to conjugacy, only one homomorphism $f_q: (G_1)_{\Fpbar} \rightarrow (G_2)_{\Fpbar}$
 of the associated reductive groups over $\Fpbar$
 which fits into a diagram
 $$  \mathcal{F} \ni f_{\C} \leftarrow f_{\Q} \rightarrow f_{\Q_p} \leftarrow f_{\Z_p} \rightarrow f_{\Fpbar}$$
 Here the $f_R: (G_1)_{R} \rightarrow (G_2)_R$ denote homomorphisms of 
 group schemes over $R$, and arrows $f_R \rightarrow f_{R'}$
 denote base change by a ring homomorphism $R \rightarrow R'$. 
 \end{quote}

    \end{proposition}

\proof
    By assumption there is a $\Q$-homomorphism $f_{\Q}^0$ in the conjugacy class of $f$. 
   We will fix such an $f_{\Q}^0$ once and for all. 
  The constants $p_0, N$ coming out of the argument are going to depend on the choice of $f_{\Q}^0$. 
 
 Suppose given diagrams
\begin{equation} \label{ffprime} f_{\C} \leftarrow f_{\Q} \rightarrow f_{\Z_p}, f'_{\C} \leftarrow f_{\Q}' \rightarrow f'_{\Z_p}\end{equation}
where $f'_{\C}, f_{\C}$ are conjugate.
 Let $f_{\bar{p}}, f'_{\bar{p}}$ be the base-changes of $f_{\Z_p}, f'_{\Z_p}$
 along 
  $\Z_p \rightarrow \overline{\mathbb{F}_p}$.  We will
  first of all show that $f_{\bar{p}}, f'_{\bar{p}}$ are conjugate to one another so long as $p \geq p_0$.

   We will need a version of ``the scheme of homomorphisms from $G_1$ to $G_2$.''
 We fix 
 a $G_2 \times G_2$-stable subspace 
subspace $W_2 \subset \C[G_2]$ such that $W_2 \cap \Z[G_2]$
(integral lattice defined by the Chevalley form) generates $\Z[G_2]$. 
Let us fix a similarly stable $W_1 \subset \C[G_1]$ 
with the property that $f_{0,\Q}^* W_1 \subset W_2$.    The same inclusion is then holds
true for both $f_{\Q}$ and $f_{\Q}'$: 
both $f_{\Q}$ and $f_{\Q}'$ send $W_{2} \cap \Q[G_2]$
to $W_1 \cap \Q[G_1]$.

 The functor that sends a ring $A$
to homomorphisms of $A$-group schemes
$\theta: G_{1,A} \rightarrow G_{2,A}$ with the property that 
$\theta^* W_{2,A} \subset W_{1,A}$ is representable by a finite type $\Z$-scheme that we will call $Y$. \footnote{
In    SGA3, \cite[Corollaire 7.2.3,XXIV]{SGA3} there is a more systematic treatment of the $\Hom$-scheme.
However, the ``full'' $\Hom$-scheme is much bigger, because it includes, e.g., Frobenius morphisms in characteristic $p$. By imposing
the condition that $\theta^*$ map $W_2$ to $W_1$, we  eliminate the Frobenius morphism in almost all characteristics
thereby producing a finite type scheme.}
 Indeed, $\theta^*$ defines an element of the affine space $\Hom(W_{2,A}, W_{1,A})$,
 thus realizing the functor as a subfunctor of this affine space. 
The condition that $\theta$ define a homomorphism of group schemes amounts 
to imposing the condition that it extend to a Hopf algebra homomorphism,
which is seen to be the $A$-points of a $\Z$-subscheme.
Moreover, there is an evident action of $G_2$ on $Y$
by post-composing a homomorphism with conjugation on $G_2$.

Now, 
 $f_{\Q}, f_{\Q}'$ both define $\Q$-points of $Y$.
 They lie in the same irreducible component of $Y_{\C}$,
comprising   the $(G_2)_{\C}$-orbit of either one. 
 Let $Y_0$ denote its Zariski closure of this irreducible component inside $Y$
(that is to say, the  Zariski closure of the corresponding set of closed points in the underlying scheme $Y$).
We equip $Y_0$, {\em a priori} a closed subset, with
the reduced scheme structure.

 By general principles, our ``reference'' morphism $f^0_{\Q}$
 extends to a $\Z[1/N]$-point of $Y$ for some $N$;
 necessarily, the resulting morphism $\Spec \ \Z[1/N] \rightarrow Y$
 factors uniquely through $Y_0$. 
 To simplify typography write $\Z' :=\Z[1/N]$. 
thus, $f^0_{\Q}$ extends to a  $\Z'$-point $f^0_{\Z'}$ of $Y_0$. 
Acting on this $f^0_{\Z'}$ we get a morphism of schemes
over $\Z'$
\begin{equation} \label{orbitGG} (G_2)_{\Z'} \rightarrow (Y_0)_{\Z'}\end{equation}
whose image is, by Chevalley's theorem, a constructible subset,
and it contains all points on the generic fiber $(Y_0)_{\Q}$. 
 The complement of this constructible
 set is itself constructible, and so has constructible image
 in $\mathrm{Spec}\ \Z'$, 
  disjoint from the generic point of $\mathrm{Spec} \ \Z'$. Therefore, enlarging $N$,
  we can suppose that this orbit map \eqref{orbitGG} is surjective at the pointwise level.
  This implies, in particular, that 
 each fiber $Y_0 \times \overline{\mathbb{F}_p}$
  is a single orbit of $(G_2)_{\overline{\mathbb{F}_p}}$
for sufficiently large $p$.

Return now to \eqref{ffprime}. 
Observe that $f_{\Z_p}$ and $f_{\Z_p}'$
correspond to maps $\Spec \  \Z_p \rightarrow Y$; the condition
that $f_{\Z_p}^* W_{2,\Z_p} \subset W_{1, \Z_p}$ and its analogue for $f'$ follow from the statement over $\Q$.
 The image of these maps $\Spec  \ \Z_p \rightarrow Y$   lie inside $Y_0$, 
 and therefore the maps  themselves factor uniquely through $Y_0 \hookrightarrow Y$. 
In particular,  so long as $p$ does not divide $N$,  the maps $\Spec \overline{\F_p} \rightarrow Y$ classifying $f_{\bar{p}}, f'_{\bar{p}}$
have images belonging to the same   $G_2$-orbit, i.e.
 $f_{\bar{p}}, f'_{\bar{p}}$ are conjugate to one another by  $G_2(\overline{\mathbb{F}_p})$.

  \qed